\documentclass{amsbook}

\usepackage{amssymb}
\usepackage{amsfonts}
\usepackage{amsthm}
\usepackage{amsmath}
\usepackage{mathtools}
\usepackage{bbm}
\usepackage{bussproofs}
\usepackage{mathrsfs}
\usepackage{tikz, tikz-cd}
\usepackage[all,2cell]{xy}
\usepackage{epigraph}
\usepackage{ stmaryrd } 
\usepackage{todonotes}
\usepackage{imakeidx} 
\usepackage{hyperref}
\usepackage{prerex}

\tikzset{%
	symbol/.style={%
		draw=none,
		every to/.append style={%
			edge node={node [sloped, allow upside down, auto=false]{$#1$}}}
	}
}

\usetikzlibrary{matrix,arrows}

\newtheorem{Theorem}{Theorem}

\newtheorem{conjecture}[Theorem]{Conjecture}

\newtheorem{proposition}[Theorem]{Proposition}
\newtheorem{lemma}[Theorem]{Lemma}
\newtheorem{corollary}[Theorem]{Corollary}

\newtheorem{thmx}{Theorem}
\newtheorem{propx}[thmx]{Proposition}
\newtheorem{lemmx}[thmx]{Lemma}
\newtheorem{corx}[thmx]{Corollary}

\theoremstyle{definition}
\newtheorem{example}[Theorem]{Example}
\newtheorem{remark}[Theorem]{Remark}
\newtheorem{definition}[Theorem]{Definition}

\newtheorem{desideratum}[Theorem]{Desideratum}

\DeclareMathOperator{\Ascr}{\mathscr{A}}
\DeclareMathOperator{\Bscr}{\mathscr{B}}
\DeclareMathOperator{\Cscr}{\mathscr{C}}
\DeclareMathOperator{\Dscr}{\mathscr{D}}
\DeclareMathOperator{\Escr}{\mathscr{E}}
\DeclareMathOperator{\Fscr}{\mathscr{F}}

\DeclareMathOperator{\Iscr}{\mathscr{I}}

\DeclareMathOperator{\Kscr}{\mathscr{K}}

\DeclareMathOperator{\Tscr}{\mathscr{T}}

\DeclareMathOperator{\Vscr}{\mathscr{V}}

\DeclareMathOperator{\Abb}{\mathbb{A}}

\DeclareMathOperator{\Ebb}{\mathbb{E}}
\DeclareMathOperator{\Fbb}{\mathbb{F}}
\DeclareMathOperator{\Gbb}{\mathbb{G}}

\DeclareMathOperator{\Ibb}{\mathbb{I}}

\DeclareMathOperator{\Nbb}{\mathbb{N}}

\DeclareMathOperator{\Pbb}{\mathbb{P}}

\DeclareMathOperator{\Dcal}{\mathcal{D}}

\DeclareMathOperator{\Ocal}{\mathcal{O}}
\DeclareMathOperator{\Pcal}{\mathcal{P}}

\DeclareMathOperator{\Tcal}{\mathcal{T}}

\DeclareMathOperator{\Afrak}{\mathfrak{A}}
\DeclareMathOperator{\afrak}{\mathfrak{a}}

\DeclareMathOperator{\Cfrak}{\mathfrak{C}}

\DeclareMathOperator{\Dfrak}{\mathfrak{D}}

\DeclareMathOperator{\Efrak}{\mathfrak{E}}

\DeclareMathOperator{\mfrak}{\mathfrak{m}}

\DeclareMathOperator{\LAvg}{LAvg_{\mathnormal{G}}^{\mathnormal{H}}}
\DeclareMathOperator{\RAvg}{RAvg_{\mathnormal{G}}^{\mathnormal{H}}}

\DeclareMathOperator{\QCoh}{\mathbf{QCoh}}
\DeclareMathOperator{\Shv}{\mathbf{Shv}}

\DeclareMathOperator{\id}{id}

\DeclareMathOperator{\Hom}{Hom}
\DeclareMathOperator{\Ker}{Ker}

\DeclareMathOperator{\Set}{\mathbf{Set}}

\DeclareMathOperator{\R}{\mathbb{R}}
\DeclareMathOperator{\C}{\mathbb{C}}
\DeclareMathOperator{\Z}{\mathbb{Z}}
\DeclareMathOperator{\N}{\mathbb{N}}

\DeclareMathOperator{\KVect}{\mathbf{K-Vect}}

\DeclareMathOperator{\RMod}{\mathbf{R-Mod}}

\DeclareMathOperator{\Ring}{\mathbf{Ring}}

\DeclareMathOperator{\Cring}{\mathbf{Cring}}

\DeclareMathOperator{\Ab}{\mathbf{Ab}}

\DeclareMathOperator{\Q}{\mathbb{Q}}

\DeclareMathOperator*{\colim}{colim}

\DeclareMathOperator{\CalO}{\mathcal{O}}

\DeclareMathOperator{\Spec}{Spec}
\DeclareMathOperator{\Sch}{\mathbf{Sch}}
\DeclareMathOperator{\FResl}{\mathbf{FResl}}
\DeclareMathOperator{\Dom}{Dom}

\DeclareMathOperator{\Free}{Free}
\DeclareMathOperator{\Forget}{Foret}
\DeclareMathOperator{\Cat}{\mathbf{Cat}}

\DeclareMathOperator{\Top}{\mathbf{Top}}

\DeclareMathOperator{\Codom}{Codom}

\DeclareMathOperator{\Iso}{\mathcal{I}so}

\DeclareMathOperator{\op}{op}

\DeclareMathOperator{\GL}{GL}

\DeclareMathOperator{\Per}{\mathbf{Perv}}

\DeclareMathOperator{\incl}{incl}

\DeclareMathOperator{\trace}{trace}

\DeclareMathOperator{\Sf}{\mathbf{Sf}}
\DeclareMathOperator{\GVar}{\mathnormal{G}-\mathbf{Var}}

\DeclareMathOperator{\OXMod}{\mathcal{O}_{\mathnormal{X}}-\mathbf{Mod}}

\DeclareMathOperator{\Ch}{\mathbf{Ch}}

\DeclareMathOperator{\VSet}{\mathscr{V}-\mathbf{Set}}
\DeclareMathOperator{\VAb}{\mathscr{V}-\mathbf{Ab}}
\DeclareMathOperator{\fTopos}{\mathfrak{Topos}}

\DeclareMathOperator{\CAT}{\mathbf{CAT}}
\DeclareMathOperator{\AlgGrp}{\mathbf{SAlgGrp}}

\DeclareMathOperator{\fCat}{\mathfrak{Cat}}
\DeclareMathOperator{\fCAT}{\mathfrak{CAT}}
\DeclareMathOperator{\SfResl}{\mathbf{SfResl}}
\DeclareMathOperator{\quo}{\mathsf{quot}}

\DeclareMathOperator{\Loc}{\mathbf{Loc}}
\DeclareMathOperator{\DMod}{\mathnormal{D}-\mathbf{Mod}}
\DeclareMathOperator{\DGMod}{\mathcal{D}_{\mathnormal{G}}-\mathbf{Mod}}
\DeclareMathOperator{\eq}{eq}
\DeclareMathOperator{\Eq}{Eq}
\DeclareMathOperator{\Coeq}{Coeq}
\DeclareMathOperator{\coeq}{coeq}
\DeclareMathOperator{\true}{\mathsf{true}}

\DeclareMathOperator{\PC}{\mathbf{PC}}

\DeclareMathOperator{\PreEq}{\mathbf{PreEq}}

\DeclareMathOperator{\HVar}{\mathnormal{H}-\mathbf{Var}}
\DeclareMathOperator{\Var}{\mathbf{Var}}
\DeclareMathOperator{\Ad}{Ad}
\DeclareMathOperator{\inv}{inv}

\DeclareMathOperator{\Open}{\mathbf{Open}}

\let\emptyset\varnothing
\let\epsilon\varepsilon

\newcommand{\os}[2]{\overset{#1}{#2}}
\newcommand{\us}[2]{\underset{#1}{#2}}
\newcommand{\ds}[3]{\us{#2}{\os{#1}{#3}}}

\newcommand{\Dbeq}[2]{\mathnormal{D}^{\mathnormal{b}}_{{#1}}({#2})}
\newcommand{\Dbeqstack}[1]{\mathnormal{D}^{\mathnormal{b}}_{\text{eq}}(\underline{[\mathnormal{G} \backslash \mathnormal{X}]}_{\bullet};\overline{\Q}_{\ell})}

\newcommand{\DbQl}[1]{\mathnormal{D}_{\mathnormal{c}}^{\mathnormal{b}}({#1};\overline{\Q}_{\ell})}
\newcommand{\DbeqQl}[1]{\mathnormal{D}_{\mathnormal{G}}^{\mathnormal{b}}({#1};\overline{\Q}_{\ell})}

\DeclareMathOperator{\fX}{\quot{\overline{\mathnormal{f}}}{\mathnormal{X}}}

\DeclareMathOperator{\hGamma}{\quot{\overline{\mathnormal{h}}}{\Gamma}}

\DeclareMathOperator{\XGamma}{\quot{\mathnormal{X}}{\Gamma}}

\DeclareMathOperator{\XGammap}{\quot{\mathnormal{X}}{\Gamma^{\prime}}}

\DeclareMathOperator{\AGamma}{\quot{\mathnormal{A}}{\Gamma}}
\DeclareMathOperator{\AGammap}{\quot{\mathnormal{A}}{\Gamma^{\prime}}}
\DeclareMathOperator{\AGammapp}{\quot{\mathnormal{A}}{\Gamma^{\prime\prime}}}
\DeclareMathOperator{\BGamma}{\quot{\mathnormal{B}}{\Gamma}}
\DeclareMathOperator{\BGammap}{\quot{\mathnormal{B}}{\Gamma^{\prime}}}

\DeclareMathOperator{\opi}{\overline{\pi}}
\DeclareMathOperator{\of}{\overline{\mathnormal{f}}}

\DeclareMathOperator{\oalphaGamma}{\quot{\overline{\alpha}_{\mathnormal{X}}}{\Gamma}}
\DeclareMathOperator{\oalphaGammap}{\quot{\overline{\alpha}_{\mathnormal{X}}}{\Gamma^{\prime}}}
\DeclareMathOperator{\fGX}{\quot{\of}{\mathnormal{G \times X}}}
\DeclareMathOperator{\muGamma}{\mu_{\Gamma}}
\DeclareMathOperator{\muGammap}{\mu_{\Gamma^{\prime}}}
\DeclareMathOperator{\opiGamma}{\quot{\opi_2}{\Gamma}}
\DeclareMathOperator{\opiGammap}{\quot{\opi_2}{\Gamma^{\prime}}}
\DeclareMathOperator{\pGamma}{\overline{\mathnormal{p}}_{\Gamma}}
\DeclareMathOperator{\pGammap}{\overline{\mathnormal{p}}_{\Gamma^{\prime}}}
\DeclareMathOperator{\aGamma}{\overline{\mathnormal{a}}_{\Gamma}}
\DeclareMathOperator{\aGammap}{\overline{\mathnormal{a}}_{\Gamma^{\prime}}}
\DeclareMathOperator{\rhoGamma}{\quot{\rho}{\Gamma}}
\DeclareMathOperator{\Resl}{\mathbf{Resl}}

\DeclareMathOperator{\Bicat}{\mathsf{Bicat}}
\DeclareMathOperator{\One}{\mathbbm{1}}
\DeclareMathOperator{\cnst}{cnst}
\DeclareMathOperator{\El}{El}
\DeclareMathOperator{\LC}{\mathbf{LC}}

\let\emptyset\varnothing
\let\epsilon\varepsilon

\newcommand{\pullbackcorner}[1][dr]{\save*!/#1-1.2pc/#1:(-1,1)@^{|-}\restore}
\newcommand{\quot}[2]{\,_{#2}{#1}}

\newcommand{\odGamma}[2]{\quot{\overline{\mathnormal{d}}_{#1}^{#2}}{{\Gamma}}}
\newcommand{\ul}[1]{\underline{{#1}}}

\newcommand{\GHPreq}[2]{({#1},{#2})\mathfrak{PreEq}(\mathnormal{X})}

\UseTwocells

\numberwithin{Theorem}{section}
\numberwithin{equation}{section}

\makeindex[name = terminology, title = Index of Terminology, intoc] 
\makeindex[name = notation, title = Index of Notation, intoc] 

\title{Categories of Pseudocones and Equivariant Descent}
\author{Geoff Vooys}
\date{\today}

\begin{document}

\begin{abstract}
In this monograph we provide an in-depth and systematic study of pseudolimits of pseudofunctors $F:\mathscr{C}^{\op} \to \mathfrak{Cat}$ in the $2$-category of categories where $\mathscr{C}$ is a $1$-category and use this to give an explicit and careful study of the category theory used in representation theory, equivariant algebraic geometry, and equivariant algebraic topology and give a unifying language to study equivariant sheaves, equivariant perverse sheaves, and their equivariant derived categories. We show how to use the pseudocone construction $\mathsf{Bicat}(\mathscr{C}^{\op},\mathfrak{Cat})(\operatorname{cnst}(1),F)$ in order to derive categorical and homological properties of the pseudolimit of $F$. We explicitly show when the pseudolimit of $F$ is complete, cocomplete, enriched in models of a Lawvere theory, (braided) monoidal, regular, triangulated, admits $t$-structures, and more. We use these various structural results to give a new category-theoretic proof and construction of the equivariant standard and pervese $t$-structures and equivariant six functor formalism for the equivariant derived category $D_G^b(X)$ in both the geometric and topological cases as well as for $D_G^b(X;\overline{\Q}_{\ell})$ in the geometric case. We also show in what sense precise sense we can view the equivariant derived category in terms of localizations. After restricting to the case of group resolution categories, we show the existence of a natural isomorphism $\Theta:\alpha_X^{\ast} \Rightarrow \pi_2^{\ast}$ which satisfies a pseudofunctorial version of the cocycle condition $d_1^{\ast}\Theta = d_2^{\ast}\Theta \circ d_0^{\ast}\Theta$. We also use the pseudocone formalism to give an in-depth analysis of change of groups functors. We use the pseudocone formalism and $\Theta$ to develop a notion of equivariant trace with an eye towards the representation theory of $p$-adic groups.
\end{abstract}


\maketitle
\tableofcontents

\chapter{Introduction}\label{Section: Intro}

Equivariant descent and equivariant homological algebra are very important tools (and areas of study in their own right) in algebraic geometry, algebraic topology, orbifold theory, the Langlands Programme, representation theory, $D$-module theory, topological quantum field theories, (sheaf-theoretic) microlocal analysis, (graded affine) Hecke algebras, arithmetic geometry, certain aspects of category theory, and nearly anywhere it is important to not only have local information about some sort of global object of interest but also where it is crucial to know that local information respects some sort of local or global symmetries on the object of interest. There are many technical difficulties which arise when asking for objects to behave reasonably with the symmetries which they need to interact or respect. For instance, we have known since \cite{GIT} that to describe this local information/global symmetry situation sheaf-theoretically it is enough to ask for a sheaf $\Fscr$ on a space $X$ carrying a group action $\alpha_X:G \times X \to X$ together with an isomorphism $\theta:\alpha_X^{\ast}\Fscr \to \pi_2^{\ast}\Fscr$ of sheaves on $G \times X$ which satisfies the cocycle condition
\[
d_1^{\ast}\theta = d_0^{\ast}\theta \circ d_2^{\ast}\theta;
\] 
cf.\@ Definition \ref{Defn: GIT Cocycle} for more information and a description of the maps $d_0$ and $d_2$. When studying equivariant cohomology, however, it is woefully insufficient to ask for a similar condition for complexes in the derived category $D^b_c(X)$\footnote{Namely the corresponding category of objects $A$ of $D_c^b(X)$ equipped with and isomorphisms $\theta:\alpha_X^{\ast}A \to \pi_2^{\ast}A$ in $D_c^b(G \times X)$ which satisfies the given cocycle condition need not even be a triangulated category; cf.\@ \cite{Torsten}.} and perhaps even worse it need not be the case that the equivariant derived category is a derived category of equivariant sheaves (for instance it is known that $D^b_{\Gbb_m}(\Pbb_{\C}^1) \not\simeq D^b(\Shv_{\Gbb_m}(\Pbb_{\C}^{1})$ and that $D^b_{\Gbb_m}(\Abb^1_{\C}) \not\simeq D^b(\Shv_{\Gbb_m}(\Abb^1_{\C}))$). As a result, a more robust and careful approach needs to be taken in order to study equivariant homological algebra which requires a significant use of descent-theoretic language. Luckily for us in the present day, however, setting this up has already been done.

It has been well-known since the 1990's (cf.\@ \cite{BernLun} for the topological case and \cite{LusztigCuspidal2} for the geometric case, for instance, or \cite{PramodBook} for a textbook account) that in order to make an equivariant derived category which behaves in the ways we require\footnote{This is vague at the moment, but a large list can be found in \cite[Desiderata 1.0.1]{MyThesis}. For the scope of this introduction, however, the important things an equivariant derived category $\Dscr_G(X)$ needs to have are: $\Dscr_G(X)$ needs a forgetful functor $\Dscr_G(X) \to D_c^b(X)$; $\Dscr_G(X)$ needs to be triangulated and have $t$-structures with hearts equivalent to the usual categories of equivariant sheaves and equivariant perverse sheaves; $\Dscr_G(X)$ needs to have a six functor formalism; and each of the six functors, the perverse cohomology functors ${}^{p}H_G:\Dscr_G(X) \to \Per_G(X)$, and the standard cohomology functors $H_G:\Dscr_G(X) \to \Shv_G(X)$ need to commute up to natural isomorphism with the forgetful functors $\Dscr_G(X) \to D_c^b(X)$, $\Shv_G(X) \to \Shv(X)$, and $\Per_G(X) \to \Per(X)$. This is developed and shown completely and carefully in Chapter \ref{Section: Six Functor Formalism} below, but can also be found in \cite{PramodBook} and in \cite{MyThesis}. Note also that we have left it completely vague as to whether or not $G$ and $X$ are topological spaces or varieties; both situations, as well as the $\ell$-adic situation for the case of varieties, are established in Chapter \ref{Section: Six Functor Formalism} below.}, we need to work with categories of resolutions of the group action on our object and track the ways in which objects become isomorphic as resolutions become more-and-more acyclic. This information culminates in a technically imposing gadget whose constructions often did not have explicit presentations in the literature or whose proofs were lacking in sufficient detail to make them usable in many practical cases. However, despite these limitations equivariant derived categories have remained important in the literature and in the Langlands Programme all the same; in fact, in \cite{CFMMX} equivariant derived categories and their $\ell$-adic versions have been used to provide a generalization of Arthur packets to what are called ABV packets for the study of $p$-adic groups. 

Because of the need to perform more detailed and precise calculations involving the equivariant derived category, a generalization to a class of categories called equivariant categories which make the actual category theory used more clear was performed in my PhD thesis \cite{MyThesis}. This generalization has led to an application in abstract algebraic differential geometry (cf.\@ \cite{DoretteMe}) and promises further applications both inside and outside algebraic geometry; however, the approach in \cite{MyThesis} is somewhat unrefined due to a missing observation which is, in hindsight, obvious: the equivariant derived category $D_G^b(X)$, and more generally an equivariant category $F_G(X)$, on a variety is really ``just'' the category of pseudocones from the constant pseudofunctor with value $\One$ into $F$. Consequently the equivariant descent therein is confined only for the case of smooth algebraic groups acting upon varieties, while the topological case of equivariant categories and the more general level of many of the proofs given in \cite[Chapters 2 -- 5]{MyThesis} were left untouched due to the missing mention of pseudocones.

Pseudocones are an important concept in higher category theory which allow us to argue about limits and colimits in ($2$-or-higher)-categories in which a cone $C$ to a diagram $D:I \to \Cfrak$ need not commute on-the-nose but instead only be a cone up to (coherent) invertible $2$-cell. For instance, if we have a cospan $X \xrightarrow{f} Z \xleftarrow{g} Y$ in a $2$-category $\Cfrak$ then a pseudocone over said cospan can be seen as an object $W$ in $\Cfrak$ for which there are morphisms $s:W \to X$ and $t:W \to Y$ together with an invertible $2$-cell $\alpha:f \circ s \xRightarrow{\cong} g \circ t$:
\[
\begin{tikzcd}
W \ar[r, ""{name = U}]{}{s} \ar[d, swap]{}{t} & X \ar[d]{}{f} \\
Y \ar[r, swap, ""{name = D}]{}{g} & Z \ar[from = U, to = D, Rightarrow, shorten <= 4pt, shorten >= 4pt]{}{\alpha} \ar[from = U, to = D, Rightarrow, shorten <= 4pt, shorten >= 4pt, swap]{}{\cong}
\end{tikzcd}
\]
Alternatively, viewing a cospan as a pseudofunctor $C:\Lambda_2^2 \to \Cfrak$, a pseudocone with object value $W$ of $C$ is a pseudonatural transformation $\cnst(W) \Rightarrow C$:
\[
\begin{tikzcd}
W \ar[r, equals, ""{name = UL}]{}{} \ar[d, swap]{}{t} & W \ar[d]{}[description]{\omega} \ar[r, equals, ""{name = UR}]{}{} & W \ar[d]{}{s} \\
Y \ar[r, swap, ""{name = DL}]{}{g} & Z & X \ar[l, ""{name = DR}]{}{f} \ar[from = UL, to = DL, swap, Rightarrow, shorten <= 4pt, shorten >= 4pt]{}{\cong} \ar[from = UR, to = DR, Rightarrow, shorten <= 4pt, shorten >= 4pt]{}{\cong}
\end{tikzcd}
\]

When our ambient $2$-category is the $2$-category of categories, $\fCat$, this perspective on pseudocones is particularly well-suited to equivariant descent. If $W = \One$ is the terminal category then we are saying that a pseudocone
\[
\begin{tikzcd}
	\One \ar[r, ""{name = U}]{}{s} \ar[d, swap]{}{t} & \Cscr \ar[d]{}{f} \\
	\Dscr \ar[r, swap, ""{name = D}]{}{g} & \Escr \ar[from = U, to = D, Rightarrow, shorten <= 4pt, shorten >= 4pt]{}{\alpha} \ar[from = U, to = D, Rightarrow, shorten <= 4pt, shorten >= 4pt, swap]{}{\cong}
\end{tikzcd}
\]
is a choice of objects $s$ in $\Cscr$ and $t$ in $\Dscr$ together with an isomorphism $f(s) \xrightarrow{\alpha} g(t)$ in $\Escr$. When the categories $\Cscr,$ $\Dscr,$ and $\Escr$ are ``algebraic-geometric'' categories (such as the derived category of a resolution of a group action or the category of sheaves) then such an isomorphism $\alpha$ can record geometric information; other examples can be found in the Galois sheaves of \cite[Expos{\'e} XIII]{SGA7}. In many cases, these natural isomorphisms are indexed or determined in some way by a group action $G$ on $X$ and the categories and natural isomorphisms record how the isomorphisms $G$ determines are respected by the objects of the categories which appear in the cone. We are then naturally led to work with pseudocones and the category of pseudocones $\Bicat(\Cscr^{\op},\fCat)(\cnst(\One),F)$ in the study of equivariant descent and so can use these pseudocone categories to aid our understanding of equivariant descent.

The categories $\Bicat(\Cscr^{\op},\fCat)(\cnst(\One),F)$ are important in the literature for purely categorical reasons: they describe and give models for pseudolimits of shape $F$ in $\fCat$. This has been known at least since \cite{GiraudCN} when they were first described in terms of Cartesian sections, and have also been described in bicategorical terms in \cite{TwoDimCat}. In this monograph we are not as much interested in their abstract theory or existence, although we do verify that $\Bicat(\Cscr^{\op},\fCat)(\cnst(\One),F)$ is a pseudolimit below for the sake of a sanity check and because this is not explicitly verified in \cite{GiraudCN}, we are mainly interested in the properties that these categories have and how they can aid us in our study and performance of equivariant descent.

Because the equivariant derived categories are all pseudocone categories (cf.\@ Sections \ref{Subsection: Equivariant Derived Cat of Var} and \ref{Subsection: EDC of a Space}) and the categories of equivariant sheaves are all equivalent to pseudocone categories, the pseudocone formalism gives a unifying language for studying equivariant sheaves and equivariant derived categories at the same time. It also provides easily accessible formal tools and perspectives which clean up, clarify, and explicitly prove many constructions brushed under the rug in the literature. For instance, the adjunction
\[
\begin{tikzcd}
D_G^b(X;\overline{\Q}_{\ell}) \ar[rr, swap, ""{name = D}, bend right = 20]{}{Rh_{\ast}} & & D_G^b(Y;\overline{\Q}_{\ell}) \ar[ll, bend right = 20, swap, ""{name = U}]{}{h^{\ast}} \ar[from = U, to = D, symbol = \dashv]
\end{tikzcd}
\]
when there is a $G$-equivariant morphism $h:X \to Y$ exists not because of arcane constructions but instead because we have an adjunction between derived category pseudofunctors in the $2$-category $\Bicat(\Cscr^{\op},\fCat)$ and the corresponding adjoint between equivariant derived categories arises as the pseudoconification of said $2$-categorical adjunction (cf.\@ Theorem \ref{Thm: Section 3: Gamma-wise adjoints lift to equivariant adjoints} and Corollary \ref{Cor: Pseudocone Functors: Pushforward functors for equivariant maps for scheme sheaves} below). These types of observations beg an explicit and careful study of pseudocone categories and the properties they have abstractly, as any such property can be applicable to equivariant categories.

While the pseudocone categories are particularly helpful for establishing the category theory of the categories that appear when doing equivariant geometry, they are also useful for informing the study of equivariant descent. For instance, if we have either a smooth algebraic group $G$ acting on a variety $X$ or a topological group acting on a space $X$ then there is a corresponding natural isomorphism $\Theta:\alpha_X^{\ast} \xRightarrow{\cong} \pi_2^{\ast}:F_G(X) \to F_G(G \times X)$. Similarly, this high-level technology allows a careful description and construction of the Change of Groups functors $\varphi^{\sharp}:F_H(X) \to F_G(X)$ and how they behave in practice. Additionally, we can also use this formalism to define a notion of equivariant trace for an object in a certain class of equivariant categories $F_G(X)$ by leveraging the existence of the isomorphism $\Theta$ as well as developing a notion of stalks for a certain class of pseudocone categories.

A large guiding principle in the theme of this monograph, based in part in our development of the equivariance isomorphism $\Theta$, the Change of Groups functors (and the perspective in terms of a pseudofunctor; cf.\@ Corollary \ref{Cor: Finally Chofg is pseudofunctorial}), is that pseudocones inform and guide our study of equivariant descent. By using and embracing the pseudofunctors which define the categories of equivariant sheaves, equivariant derived categories, equivariant perverse sheaves, and the corresponding higher-categorical perspective that comes with carrying the pseudofunctors, pseudonatural transformations, and modifications around we not only unify the methodology used to study equivariant objects but we provide new avenues with which we can explore these objects and use them in practice.

\section{Main Results}\label{Subsection: Results}
In this section we give a highlight reel of results that appear in this monograph. These are what I think of as a guide towards some of the bigger results which indicate not only the power afforded to the pseudocone formalism in how it may be applied to equivariant descent, but also in how it is of independent categorical interest and study.

The first main results of this monograph lie in the fact that for a pseudofunctor $F:\Cscr^{\op} \to \fCat$, the pseudocone category $\PC(F) = \Bicat(\Cscr^{\op},\fCat)(\cnst(\One),F)$ is the pseudolimit of shape $F$ in $\fCat$. While this has been known since \cite{GiraudCN}, there it was shown through the connection to the category of Cartesian sections of the elements fibration $p:\El(F) \to \Cscr$ of $F$ and even still was given without proof. Our first two main results here are sanity checks that $F$ not only is indeed a pseudolimit of $F$ in $\fCat$ but also coincides with the description given in \cite{GiraudCN}. We also explicitly and carefully establish when these pseudolimits are enriched in the $\Vscr$-models of a Lawvere theory $\Tcal$ when $\Vscr$ is a Grothendieck universe.

\begin{propx}[{Cf.\@ Proposition \ref{Prop: Pseudocone Section: CSections are pseudolimits}}]
	Let $F:\Cscr^{\op} \to \fCat$ be a pseudofunctor and let $p:\El(F) \to \Cscr$ be the associated elements fibration. Then if $\mathbf{CSect}_B(p)$ is the category of Cartesian sections of $p$ and based natural transformations then $\PC(F)$ is equal to the category $\mathbf{CSect}_B(p)^{\op}$.
\end{propx}
\begin{thmx}[{Cf.\@ Theorem \ref{Thm: Section Pseudocones: PCF is the pseudolimit of F}}]
Let $F:\Cscr^{\op} \to \fCat$ be a pseudofunctor. Then $\PC(F)$ is the pseudolimit of $F$ in $\fCat$.
\end{thmx}
\begin{lemmx}[{Cf.\@ Lemma \ref{Lemma: Pseudocone Section: Lawvere Enrichment}}]
	Let $\Vscr$ be a Grothendieck universe and let $F:\Cscr^{\op} \to \fCat$ be a $\Vscr$-locally small pseudofunctor. Let $\Tcal$ be a Lawvere theory. If each category $F(X)$ is enriched in $\Tcal(\VSet)$ for all $X \in \Cscr_0$ then the category $\PC(F)$ is enriched in $\Tcal(\VSet)$ as well.
\end{lemmx}

Our next main results in the monograph regard various categorical structures and properties that the category of pseudocones $\PC(F)$ can be asked to satisfy. In particular, they examine (co)limits of specified shape in $\PC(F)$, (braided) monoidal structures on $\PC(F)$, when $\PC(F)$ is regular, and when $\PC(F)$ admits a subobject classifier.

\begin{thmx}[{Cf.\@ Theorem \ref{Thm: Section 2: Equivariant Cat has lims}}]
Let $F:\Cscr^{\op} \to \fCat$ be a pseudofunctor and let $I$ be an index category for which there is a diagram $d:I \to \PC(F)$. For each $X \in \Cscr_0$ let $d_{X}:I \to F(X)$ denote the diagram functor
\[
\begin{tikzcd}
	I \ar[drr, swap]{}{d_{X}} \ar[r]{}{d} & \PC(F) \ar[r]{}{\tilde{\imath}} & \prod\limits_{X \in \Cscr_0} F(X) \ar[d]{}{\pi_{X}} \\
	& & F(X)
\end{tikzcd}
\]
and assume that for each $X \in \Cscr_0$ the limit $\lim d_{X}$ exists in $F(X)$. Furthermore, assume that for every morphism $f \in \Cscr_1$ with $f:X \to Y$ there is an isomorphism $\theta_f$ witnessing the limit preservation
\[
F({f})\left(\lim_{\substack{\longleftarrow \\ i \in I}}d_{X}(i)\right) \cong \lim_{\substack{\longleftarrow \\ i \in I}} F({f})(d_{X}(i)).
\]
Then the diagram $d$ admits a limit in $\PC(F)$.
\end{thmx}
\begin{thmx}[{Cf.\@ Theorem \ref{Theorem: Section 2: Monoidal preequivariant pseudofunctor gives monoidal equivariant cat}}]
Let $F:\Cscr^{\op} \to \fCat$ be a monoidal pseudofunctor. Then $\PC(F)$ is a monoidal category.
\end{thmx}
\begin{propx}[{Cf.\@ Proposition \ref{Prop: Section 2: Equivariant cat is symmetric monoidal}}]
		Let $F:\Cscr^{\op} \to \fCat$ be a braided monoidal pseudofunctor. Then $\PC(F)$ is braided monoidal with braiding $\beta_{A,B}:A \otimes B \to B \otimes A$ given by
		\[
		\beta_{A,B} := \left\lbrace \quot{\beta_{A,B}}{X} \; | \; X \in \Cscr_0 \right\rbrace
		\]
		where $\quot{\beta_{A,B}}{X}$ is a shorthand for the braiding isomorphism $\quot{\beta_{\quot{A}{X},\quot{B}{X}}}{X}$ in $F(X)$. Additionally if each category $F(X)$ is symmetric monoidal then so is $\PC(F)$.
\end{propx}
\begin{propx}[{Cf.\@ Proposition \ref{Prop: Regularity of Equivariant Category}}]
Let $F:\Cscr^{\op} \to\fCat$ be a pseudofunctor such that each fibre category $F(X)$ is regular and for which each fibre functor $F(f)$ is finitely complete and preserves regular epimorphisms. Then $\PC(F)$ is regular.
\end{propx}
\begin{propx}[{Cf.\@ Proposition \ref{Prop: Section 2: Equivariant Cat has Subobject Classifiers}}]
Let $F$ be a pseudofunctor such that each fibre category $F(X)$ has a subobject classifier, $\quot{\Omega}{X}$, and for which each fibre functor $F(f)$ preserves subobject classifiers, terminal objects, and subobject classifying pullbacks. Then $\PC(F)$ admits a subobject classifier $\Omega$ whose object collection has the form
\[
\Omega = \lbrace \quot{\Omega}{X} \; | \; X \in \Cscr_0 \rbrace.
\]
\end{propx}

Our next major results concern the functors between pseudocone categories. The first few show that the psuedocone construction constitutes a strict $2$-functor $\PC(-):\Bicat(\Cscr^{\op},F) \to \fCat$, how this allows a clean description of some adjunctions between pseudocone categories, and how to build functors 
\[
\Bicat(\Cscr^{\op},\fCat)(\cnst(\One),F) \to \Bicat(\Dscr^{\op},\fCat)(\cnst(\One),E).\] 
While not listed in this summary section, much of the importance of these results lies in the ability to produce, define, and describe from purely categorical perspectives the corresponding pseudofunctor the equivariant push-pull and exceptional push-pull adjunctions $h^{\ast} \dashv Rh_{\ast}$ and $Rh_! \dashv h^!$ between equivariant derived categories.

\begin{thmx}[{Cf.\@ Theorem \ref{Thm: Functor Section: Psuedonatural trans are pseudocone functors}}]
Let $\Cscr$ be a category and let $F,E:\Cscr^{\op} \to \fCat$ be a pair of pseudofunctors for which there is a pseudonatural transformation $\alpha:F \to E$. Then there is an induced functor $\ul{\alpha}:\PC(F) \to \PC(E)$ induced by the vertical post-composition functor
\[
\ul{\alpha} = \alpha \circ (-):\Bicat(\Cscr^{\op},\fCat)(\cnst(\One),F) \to \Bicat(\Cscr^{\op},\fCat)(\cnst(\One),E).
\]
Explicitly, if $(A, T_A) \in \PC(F)_0$ then
\[
\ul{\alpha}(A) = \left\lbrace \quot{\alpha}{X}\left(\quot{A}{X}\right) \; : \; X \in \Cscr_0, \quot{A}{X} \in A\right\rbrace
\]
and
\[
T_{\ul{\alpha}A} := \left\lbrace \quot{\alpha}{X}\left(\tau_f^A\right) \circ \quot{\alpha}{f}^{-1}_{\quot{A}{Y}} \; : \; f \in \Cscr_1, f:X \to Y \right\rbrace
\]
where $\quot{\alpha}{f}$ is the natural isomorphism
\[
\quot{\alpha}{f}: \quot{\alpha}{X} \circ F(f) \xRightarrow{\cong} E(f) \circ \quot{\alpha}{Y}.
\]
Similarly, if $P = \lbrace \quot{\rho}{X} \; | \; X \in \Cscr_0 \rbrace$ then
\[
\ul{\alpha}P = \left\lbrace \quot{\alpha}{X}\left(\quot{\rho}{X}\right) \; : \; X \in \Cscr_0\right\rbrace.
\]
\end{thmx}
\begin{lemmx}[{Cf.\@ Lemma \ref{Lemma: Modifications give equivariant natural transformations}}]
		Let $F,E:\Cscr^{\op} \to \fCat$ be pseudofunctors and let $\alpha,\beta:F \Rightarrow E$ be pseudonatural transformations. If $\eta:\alpha \Rrightarrow \beta$ is a modification, then there exists a natural transformation $\underline{\eta}$ fitting into the $2$-cell
		\[
		\begin{tikzcd}
			\PC(F) \ar[r,bend left = 30, ""{name = U}]{}{\underline{\alpha}} \ar[r,swap, bend right = 30,""{name = L}]{}{\underline{\beta}} & \PC(E) \ar[from = U, to = L, Rightarrow,shorten >= 5pt, shorten <= 5pt]{}{\underline{\eta}}
		\end{tikzcd}
		\]
		given on objects $C = \lbrace \quot{C}{} \; | \; X \in \Cscr_0\rbrace$ of $\PC(F)$ by
		\[
		\underline{\eta}_C := \lbrace \quot{\eta}{X}_{\quot{C}{X}}:\quot{\alpha}{X}(\quot{C}{X}) \to \quot{\beta}{X}(\quot{C}{X}) \; | \; X \in \Cscr_0\rbrace.
		\]
\end{lemmx}
\begin{thmx}[{Cf.\@ Theorem \ref{Thm: Section 3: Gamma-wise adjoints lift to equivariant adjoints}}]
		Let $F,E:\Cscr^{\op} \to \fCat$ be pseudofunctors and let $\alpha:F \Rightarrow E$ and $\beta:E \to F$ be pseudonatural transformations. Then if $\quot{\alpha}{X} \dashv \quot{\beta}{X}$ in $\fCat$ for all $X \in \Cscr_0$, there is an adjunction
		\[
		\begin{tikzcd}
			F \ar[r, bend left = 30, ""{name = U}]{}{\alpha} & E \ar[l, bend left = 30, ""{name = L}]{}{{\beta}} \ar[from = U, to = L, symbol = \dashv]
		\end{tikzcd}
		\]
		in the $2$-category $\Bicat(\Cscr^{\op},\fCat)$. In particular, this gives rise to an adjunction of categories:
		\[
		\begin{tikzcd}
			\PC(F) \ar[r, bend left = 30, ""{name = U}]{}{\underline{\alpha}} & \PC(E) \ar[l, bend left = 30, ""{name = L}]{}{\underline{\beta}} \ar[from = U, to = L, symbol = \dashv]
		\end{tikzcd}
		\]
\end{thmx} 
\begin{propx}[{Cf.\@ Proposition \ref{Section: Pseudocone Functors: PC cat is monoidal closed}}]
		Let $F:\Cscr^{\op} \to \fCat$ be a monodial closed pseudofunctor for which each functor $F(f)$ is monoidal closed. Then if each category $F(X)$ is monoidal closed or symmetric monoidal closed, so is $\PC(F)$.
\end{propx}
\begin{thmx}[{Cf.\@ Theorem \ref{Thm: Pseudocone Functors: Pullback induced by fibre functors in pseudofucntor}}]
	Let $\Cscr \xrightarrow{\varphi} \Escr \xleftarrow{\psi} \Dscr$ be a cospan of categories together with a functor $\gamma:\Cscr \to \Dscr$. If there is a $2$-cell
	\[
	\begin{tikzcd}
		\Cscr \ar[rr, ""{name = U}]{}{\varphi} \ar[dr, swap]{}{\gamma} & & \Escr \\
		& \Dscr \ar[ur, swap]{}{\psi} \ar[from = U, to = 2-2, Rightarrow, shorten <= 4pt, shorten >= 4pt]{}{\alpha}
	\end{tikzcd}
	\]
	and a pseudofunctor $F:\Escr^{\op} \to \fCat$ then $E \circ \varphi^{\op}$ admits pseudocone translations (cf.\@ Definition \ref{Defn: Pseudocone section: Coone translations}) along $\gamma^{\op}$ to $F \circ \psi^{\op}$ of shape $h = (\quot{h}{X}, \quot{h}{f})$ where, for all $X \in \Cscr_0$,
	\[
	\quot{h}{X} := F(\alpha_X):F(\gamma X) \to F(X)
	\]
	and for all $f:X \to Y$ in $\Cscr$,
	\[
	\quot{h}{f} := \phi_{\varphi f,\alpha_Y}^{-1} \circ \phi_{\alpha_X,(\psi \circ \gamma) f}.
	\]
\end{thmx}

Our next main results concern the homological algebra of pseudocone categories and show that in natural situations the categories $\PC(F)$ are triangulated and admit $t$-structures. Additionally, we can use these to conclude that if the pseudocone categories have six functor formalisms, then there are relations between these six functors whenever local versions of the relations between the functors hold.
\begin{thmx}[{Cf.\@ Theorem \ref{Theorem: Section Triangle: Equivariant triangulation}}]
		Let $F:\Cscr^{\op} \to \fCat$ be a triangulated pseudofunctor (cf.\@ Definition \ref{Defn: Section Triangle: Triangulated Pseudofunctor}). Then the category $\PC(F)$ is a triangulated category with triangulation induced by saying that a triangle
		\[
		\xymatrix{
			A \ar[r]^-{P} & B \ar[r]^-{\Sigma} & C \ar[r]^-{\Phi} & A[1]
		}
		\]
		is distinguished if and only if for all $X \in \Cscr_0$ the triangle
		\[
		\xymatrix{
			\quot{A}{X} \ar[r]^-{\quot{\rho}{X}} & \quot{B}{X} \ar[r]^-{\quot{\sigma}{X}} & \quot{C}{X} \ar[r]^-{\quot{\varphi}{X}} & \quot{A}{X}[1]
		}
		\]
		is distinguished in $F(\XGamma)$.
\end{thmx}
\begin{propx}[{Cf.\@ Proposition \ref{Prop: Section Triangle: Truncated preeq gives truncation functors}}]
		Let $F:\Cscr^{\op} \to \fCat$ be a truncated pseudofunctor (cf.\@ Definition \ref{Defn: Truncated preq pseudofunctor}). Then there are truncation functors
		\[
		\tau^{\leq 0}:\PC(F) \to \PC(F), \qquad \tau^{\geq 0}:\PC(F) \to \PC(F).
		\]
		The truncation $\tau^{\leq 0}$ is defined on objects by
		\[
		\tau^{\leq 0}(\lbrace \quot{A}{X}\; | \; X \in \Cscr_0 \rbrace) = \lbrace \quot{\tau^{\leq 0}}{X}(\quot{A}{X}) \; | \; X \in \Cscr_0 \rbrace,
		\]
		with transition isomorphisms $\tau_f^{A\leq 0}$ given by the diagram
		\[
		\xymatrix{
			F(f)\left(\quot{\tau^{\leq 0}}{Y}\quot{A}{Y}\right) \ar[rr]^-{\theta_f^{\leq 0}} \ar[drr]_{\tau_f^{A \leq 0}} & & \quot{\tau^{\leq 0}}{X}\left(F(f)\quot{A}{Y}\right) \ar[d]^{\quot{\tau^{\leq 0}}{X}\left(\tau_f^A\right)} \\
			& & \quot{\tau^{\leq 0}}{X}\left(\quot{A}{X}\right)
		}
		\]
		and defined on morphisms by
		\[
		\tau^{\leq 0}\left(\lbrace \quot{\rho}{X} \; | \; X \in \Cscr_0 \rbrace\right) = \left\lbrace \quot{\tau^{\leq 0}}{X}\left(\quot{\rho}{X}\right) \; | \; X \in \Cscr_0 \right\rbrace.
		\]
		Similarly, the other truncation $\tau^{\geq 0}$ is defined on objects by
		\[
		\tau^{\geq 0}(\lbrace \quot{A}{X}\; | \; X \in \Cscr_0 \rbrace) = \left\lbrace \quot{\tau^{\geq 0}}{X}\left(\quot{A}{X}\right) \; | \; X \in \Cscr_0 \right\rbrace,
		\]
		with transition isomorphisms $\tau_f^{A\geq 0}$ given by the diagram
		\[
		\xymatrix{
			F(f)\left(\quot{\tau^{\geq 0}}{Y}\quot{A}{Y}\right) \ar[rr]^-{\theta_f^{\geq 0}} \ar[drr]_{\tau_f^{A\leq 0}} & & \quot{\tau^{\geq 0}}{X}\left(F(f)\quot{A}{Y}\right) \ar[d]^{\quot{\tau^{\geq 0}}{X}\left(\tau_f^A\right)} \\
			& & \quot{\tau^{\geq 0}}{X}\left(\quot{A}{X}\right)
		}
		\]
		and defined on morphisms by
		\[
		\tau^{\geq 0}\left(\lbrace \quot{\rho}{X} \; | \; X \in \Cscr_0 \rbrace\right) = \left\lbrace \quot{\tau^{\geq 0}}{X}\left(\quot{\rho}{X}\right) \; | \; X \in \Cscr_0 \right\rbrace.
		\]
\end{propx}
\begin{thmx}[{Cf.\@ Theorem \ref{Thm: Section Triangle: The 
			Six Functor Relations}}]
	Assume that we have $\Cscr,$ $\Dscr,$ $\Escr,$ $F$, $\alpha$, $\varphi,$ $\psi,$ $\rho,$ $(F \ast \alpha)_{\square},$ $(F \ast \alpha)^{\square},$ $(F \ast \alpha)_{\dagger},$ and $(F \ast \alpha)^{\dagger}$ be given as in the setup prior to Theorem \ref{Thm: Section Triangle: The Six Functor Relations} and assume that for all $X \in \Cscr_0$ and for all objects $A \in F(\rho X)_0$ and for all objects $B, C \in F((\psi \circ \varphi)X)_0$ we have natural isomorphisms
	\begin{align*}
		\quot{(F \ast \alpha)_{\dagger}}{X}\left(\quot{A}{X}\right) \us{F((\psi \circ \varphi)X)}{\otimes} \quot{B}{X} &\cong \quot{(F \ast \alpha)_{\dagger}}{X}\left(\quot{A}{X} \us{F(\rho X)}{\otimes} (F \ast \alpha)^{\square}(\quot{B}{X})\right)  \\
		\left[\quot{(F \ast \alpha)_{\dagger}}{X}(\quot{A}{X}), \quot{B}{X}\right]_{F((\psi \circ \varphi)(X))} &\cong \quot{(F \ast \alpha)_{\square}}{X}\left(\left[\quot{A}{X}, \quot{(F \ast \alpha)^{\dagger}}{X}(\quot{B}{X})\right]_{F(\rho X)}\right) \\
		\quot{(F \ast \alpha)^{\dagger}}{X}\left(\left[\quot{C}{X}, \quot{B}{X}\right]_{F((\psi \circ \varphi)X)}\right) &\cong \left[\quot{(F \ast \alpha)^{\square}}{X}(\quot{C}{X}), \quot{(F \ast \alpha)^{\dagger}}{X}(\quot{B}{X})\right]_{F(\rho X)}
	\end{align*}
	Then for any objects $A \in \PC(F \circ \rho^{\op})_0$ and $B, C \in \PC(F \circ \psi^{\op} \circ \varphi^{\op})_0$ there are natural isomorphisms:
	\begin{align*}
		(F \ast \alpha)_{\dagger}(A) \us{\PC(F \circ \psi^{\op} \circ \varphi^{\op})}{\otimes} B &\cong (F \ast \alpha)_{\dagger}\left(A \us{\PC(F \ast \rho^{\op})}{\otimes} (F \ast \alpha)^{\square}(B)\right) \\
		\left[(F \ast \alpha)_{\dagger}(A), B\right]_{\PC(F \circ \psi^{\op} \circ \varphi^{\op})} &\cong (F \ast \alpha)_{\square}\left(\left[A, (F \ast \alpha)^{\dagger}(B)\right]_{\PC(F \circ \rho^{\op})}\right) \\
		(F \ast \alpha)^{\dagger}\left(\left[C, B\right]_{\PC(F \circ \psi^{\op} \circ \varphi^{\op})}\right) &\cong \left[(F \ast \alpha)^{\square}(C), (F \ast \alpha)^{\dagger}(B)\right]_{\PC(F \circ \rho^{\op})}
	\end{align*}
\end{thmx}

After establishing our categorical results, we change gears and move to study equivariant descent for both topological spaces and for varieties. We introduce the notion of an equivariant category $F_G(X)$ in terms of pseudocone categories on a space/variety and then show that equivariant categories have a natural isomorphism 
\[
\Theta:\alpha_X^{\ast} \xRightarrow{\cong} \pi_2^{\ast}:F_G(X) \to F_G(G \times X)
\]
which satisfies the pseudofunctorial version of the GIT cocycle condition (cf.\@ Defition \ref{Defn: GIT Cocycle}). We also show that these equivariant categories admit Change of Groups functors $\varphi^{\sharp}:F_H(X) \to F_G(X)$ whenever there are morphisms $\varphi:G \to H$ of group objects which vary pseudofunctorially in the category of groups objects over $H$. 

\begin{thmx}[{Cf.\@ Theorems \ref{Theorem: Formal naive equivariance} and \ref{Thm: GIT Cocycle Condition Formal Equivariance}}]
There is a natural isomorphism $\Theta:\alpha_X^{\ast} \xrightarrow{\cong} \pi_2^{\ast}:F_G(X) \to F_G(G \times X)$ which satisfies the GIT cocycle condition.
\end{thmx}
\begin{thmx}[{Cf.\@ Theorem \ref{Thm: Section 3: Change of groups functor}}]
	Let $G, H,$ $\varphi:G \to H$, $X$, and $F$ be given as in either of the cases below:
\begin{itemize}
	\item $G$ and $H$ are smooth algebraic groups, $\varphi$ is a morphism of algebraic groups, $X$ is a left $H$-variety, and $F$ is an $H$-pre-equivariant pseudofunctor;
	\item $G$ and $H$ are topological groups, $\varphi$ is a moprhism of topological groups, $X$ is a left $H$-space, and $F$ is an $H$-pre-equivariant pseudofunctor.
\end{itemize} 
Then in both cases there is a functor $\varphi^{\sharp}:F_H(X) \to F_G(X)$\index[notation]{PhiSharp@$\varphi^{\sharp}$} given by sending an object $(A,T_A)$ to the pair $(A^{\prime},T_{A^{\prime}})$ and a morphism $P$ to $P^{\prime}$, where, in the geometric case, the pair $(A^{\prime}, T_{A^{\prime}})$ is induced by the equation
\[
A^{\prime} := \left\lbrace \overline{F}(h_{\Gamma})\left(\quot{A}{H \times^{G} \Gamma^{\prime}}\right) \; : \; \Gamma^{\prime} \in \Sf(G)_0, \quot{A}{H \times^{G} \Gamma^{\prime}} \in A \right\rbrace
\]
and on $P^{\prime}$ is induced by
\[
P^{\prime} := \left\lbrace \overline{F}(h_{\Gamma})\left(\quot{\rho}{H \times^{G} \Gamma^{\prime}}\right) \; : \; \Gamma^{\prime} \in \Sf(G)_0, \quot{\rho}{H \times^{G} \Gamma^{\prime}} \in P  \right\rbrace,
\]
where $h_{\Gamma}$ is the isomorphism $\XGamma \xrightarrow{\cong} \quot{X}{H \times^{G} \Gamma}$ of Lemma \ref{Lemma: Section 3: Induction space quotient is iso to the Sf(G) quotient}. The topological case is given mutatis mutandis.
\end{thmx}
\begin{propx}[{Cf.\@ Proposition \ref{Prop: Section 3.3: Chagnge of groups interacting with change of grups}}]
	Let $G_0, G_1, G_2, X, F, \varphi,$ and $\psi$ be given as in either case below:
\begin{itemize}
	\item $G_0, G_1, G_2$ are smooth algebraic groups\footnote{It is unfortunate that there is a clash of notation, but in this case $G_2$ is meant only as the third entry in a list of algebraic groups and \emph{not} as the exceptional Lie group over $\Spec K$.}, $X$ is a left $G_2$-variety, $F$ is a $G_2$-pre-equivariant pseudofunctor on $X$, and both $\varphi:G_0 \to G_1$ and $\psi:G_1 \to G_2$ are morphisms of algebraic groups.
	\item $G_0, G_1, G_2$ are topological groups, $X$ is a left $G_2$-space, $F$ is a $G_2$-pre-equivariant pseudofunctor on $X$, and both $\varphi:G_0 \to G_1$ and $\psi:G_1 \to G_2$ are morphisms of topological groups.
\end{itemize}
Then there is an invertible $2$-cell
\[
\begin{tikzcd}
	F_{G_2}(X) \ar[dr, swap]{}{(\psi \circ \varphi)^{\sharp}} \ar[rr]{}{\psi^{\sharp}} & {} &F_{G_1}(X) \ar[dl, swap]{}{\varphi^{\sharp}} \\
	& F_{G_0}(X) \ar[from = 2-2, to = 1-2, Rightarrow, shorten >= 4pt, shorten <= 4pt]{}{\alpha}
\end{tikzcd}
\]
in $\fCat$.
\end{propx}
\begin{corx}[{Cf.\@ Corollary \ref{Cor: Finally Chofg is pseudofunctorial}}]
	Let $H$ be a smooth algebraic group (respectively a topological group), let $X$ be a left $H$-variety (respectively a left $H$-space), and let $(F,\overline{F})$ be an $H$-pre-equivariant psuedofunctor on $X$. Then there is a pseudofunctor
\[
F_{-}(X):\left(\mathbf{SAlgGrp} \downarrow H\right)^{\op} \to \fCat
\]
where:
\begin{itemize}
	\item For each smooth algebraic group morphism $\varphi:G \to H$, $F_{-}(X)(\varphi) := F_G(X)$ for the $(G,H)$-pre-equivariant pseudofunctor
	\[
	\left(\left(\overline{F}\left(G \backslash (H \times^{G} ((-) \times X))\right), \overline{F}\right), (F,\overline{F}), \id_{\overline{F}}\right);
	\]
	\item For each morphism $\varphi:G \to G^{\prime}$ in $\mathbf{SAlgGrp} \downarrow H$, the fibre functor is given by the Change of Groups functor
	\[
	\varphi^{\sharp}:F_{G^{\prime}}(X) \to F_G(X).
	\]
	\item For each composable pair of morphisms $\varphi_{01}:G_0 \to G_1$ and $\varphi_{12}:G_1 \to G_2$ in $\mathbf{SAlgGrp} \downarrow H$, the compositor $\phi_{\varphi_{01}, \varphi_{12}}$ is defined by
	\[
	\phi_{\varphi_{01}, \varphi_{12}} := \alpha_{012}^{-1},
	\]
	where $\alpha_{012}$ is the natural isomorphism of Proposition \ref{Prop: Section 3.3: Chagnge of groups interacting with change of grups}.
\end{itemize}
Respectively, there is also a pseudofunctor $F_{-}(X):\left(\mathbf{TopGrp} \downarrow H\right)^{\op} \to \fCat$ whose assignment is given similarly.
\end{corx}

The last main results of the monograph we present here provide the final components of our pseudocone-categorical proof of the complete six functor formalism for the equivariant derived categories $D_G^b(X)$ for topological spaces and for varieties. We also provide one final major result showing the additivity of equivariant trace on distinguished triangles provided that the pseudofunctors themselves readily admit notions of trace for arbitrary objects.
\begin{thmx}[{Cf.\@ Theorems \ref{Thm: Section SFF: The big results on the EDC for schemes}, \ref{Thm: Section SFF: The big results on the EDC for spaces}}]
		Let $G$ be a smooth algebraic group and let $h:X \to Y$ be a $G$-equivariant morphism of left $G$-varieties. Then the following hold for both pseudofunctors $D_c^b(-;\overline{\Q}_{\ell}):\Var_{/K}^{\op} \to \fCat$ and $D_c^b(-):\Var_{/K}^{\op} \to \fCat$.
		\begin{enumerate}
			\item The categories $D_G^b(X)$ and $D^b_G(X;\overline{\Q}_{\ell})$ are symmetric monoidal closed categories.
			\item The categories $D_G^b(X)$ and $D_G^b(X;\overline{\Q}_{\ell})$ are triangulated categories which have standard and perverse $t$-structures for which
			\[
			D_G^b(X;\overline{\Q}_{\ell})^{\heartsuit_{\operatorname{stand}}} \simeq \Shv_G^{\operatorname{na{\ddot{\imath}}ve}}(X;\overline{\Q}_{\ell})
			\]
			\[
			D_G^b(X;\overline{\Q}_{\ell})^{\heartsuit_{\operatorname{per}}} \simeq \Per_G^{\operatorname{na{\ddot{\imath}}ve}}(X;\overline{\Q}_{\ell})
			\]
			and similarly for $D_G^b(X)$.
			\item There are equivariant pullback and equivariant exceptional pullback functors $h^{\ast}:D_G^b(Y;\overline{\Q}_{\ell}) \to D_G^b(X;\overline{\Q}_{\ell})$ and $h^{!}:D_G^b(Y;\overline{\Q}_{\ell}) \to D_G^b(X;\overline{\Q}_{\ell})$ as well as equivariant pushforward and equivariant proper pushforward functors $Rh_{\ast}:D_G^b(X;\overline{\Q}_{\ell}) \to D_G^b(Y;\overline{\Q}_{\ell})$ and $Rh_{!}:D_G^b(X;\overline{\Q}_{\ell}) \to D_G^b(Y;\overline{\Q}_{\ell})$. The same also holds for $D_G^b(X)$ and $D_G^b(Y)$.
			\item There are adjunctions $h^{\ast} \dashv Rh_{\ast}$ and $Rh_{!} \dashv h^{!}$.
			\item For any $G$-equivariant open immersion $j:U \to X$ and corresponding complementary $G$-equivariant closed immersion $i:V \to X$ there is an open/closed distinguished triangle
			\[
			\begin{tikzcd}
				Rj_!\left(j^!(-)\right) \ar[r]{}{\epsilon^{j}_{-}} & \id_{D_G^b(X;\overline{\Q}_{\ell})}(-) \ar[r]{}{\eta_{-}^{i}} & Ri_{\ast}\left(i^{\ast}(-)\right) \ar[r] & \left(Rj_!\left(j^!(-)\right)\right)[1]
			\end{tikzcd}
			\]
			in $D_G^b(X;\overline{\Q}_{\ell})$ and the same is true of $D_G^b(X)$.
			\item The standard and perverse $t$-structure cohomology functors on $D_G^b(X;\overline{\Q}_{\ell})$ commute up to isomorphism with the forgetful functors and the same is true of the standard and perverse $t$-structures on $D_G^b(X)$.
			\item The six functors $(-)\otimes(-),$ $R[-,-]$, $h^{\ast}$, $Rh_{\ast}$, $Rh_{!}$, and $h^{!}$ all commute up to isomorphism with the corresponding forgetful functors.
			\item There are natural isomorphisms
			\begin{align*}
				Rh_{!}(A) \us{Y}{\otimes} B &\cong Rh_{!}\left(A \us{X}{\otimes} h^{\ast}(B)\right), \\
				\left[Rh_{!}(A), B\right]_{Y} &\cong Rh_{\ast}\left(\left[A, h^{!}(B)\right]_{X}\right), \\
				h^{!}\left([C,B]_{Y}\right) &\cong \left[h^{\ast}(C), h^{!}(B)\right]_{X}
			\end{align*}
			where $A \in D_G^b(X; \overline{\Q}_{\ell})_0$ and $B, C \in D_G^b(Y;\overline{\Q}_{\ell})_0$. The same also holds for $D_G^b(X)$ and $D_G^b(Y)$.
	\end{enumerate}
The same all hold when $X$ is a topological space, $G$ is a topological group which admits $n$-acyclic free $G$-spaces which are manifolds, and for the equivariant derived category $D_G^b(X)$ save that we have the additional exact triangle
\[
\begin{tikzcd}
	Ri_!\left(i^!(-)\right) \ar[r]{}{\epsilon^i_{-}} & \id_{D^b_G(X)}(-) \ar[r]{}{\eta^j_{-}} & Rj_{\ast}\left(j^{\ast}(-)\right) \ar[r] & \left(Ri_!\left(i^!(-)\right)\right)[1]
\end{tikzcd}
\]
when $i:Z \to X$ is an equivariant closed inclusion and $j:U \to X$ is the corresponding complementary equivariant open inclusion.
\end{thmx}
\begin{thmx}[{Cf.\@ Theorem \ref{Thm: Section SFF: Additive trace}}]
Let $F$ be a triangulated trace-class pre-equivariant pseudofunctor on $X$ and $G \times X$ for which:
			\begin{itemize}
				\item The category $\Cscr$ of values the objects $A$ of $\overline{F}(\quot{X}{G})$, $\overline{F}(X)$, $\overline{F}(\quot{(G \times X)}{G})$ and $\overline{F}(G \times X)$ take is triangulated and arises as the ``derived category'' of a ``homotopy category'' in the sense of \cite[Section 5]{MayTrace}.
				\item The operation of taking stalks gives triangulated exact functors $(-)_{(g,x)}:F_G(G \times X) \to \Cscr$ and $(-)_{x}:F_G(X) \to \Cscr$.
				\item Each functor 
				\[
				\overline{F}\left(\overline{\id_{\Gamma} \times \alpha_X}\right):\overline{F}(\XGamma) \to \overline{F}\left(\quot{(G \times X)}{\Gamma}\right)
				\]
				(cf.\@ the proof of Lemma \ref{Lemma: Section Change of Space: Action map pullback}) is additive. 
			\end{itemize}
			Then equivariant trace is additive on distinguished triangles in the sense that if
			\[
			\begin{tikzcd}
				A \ar[r]{}{P} & B \ar[r]{}{\Phi} & C \ar[r]{}{\Psi} & A[1]
			\end{tikzcd}
			\]
			 is a distinguished triangle for objects $A, B, C$ in $F_G(X)$ then
			\[
			\trace_g(B,x) = \trace_g(A,x) + \trace_g(C,x).
			\]
\end{thmx}

\section{Structure of the Monograph}\label{Subsection: Structure}
In this short subsection we detail the basic structure of the monograph and where results can, generically speaking, be found. In Chapters \ref{Section: Category of Pseudocones}, \ref{Section: Limits and Colimits in Pseudocones}, \ref{Section: Functors of Pseudocones}, and \ref{Section: Pseudocone Algebra} we perform a systematic study of the category of pseudocones $\PC(F)$ of a pseudofunctor $F:\Cscr^{\op} \to \fCat$ while in Chapters \ref{Section: Equivariant Descent}, \ref{Chapter: Chofg}, and \ref{Section: Six Functor Formalism} we use our work from the earlier chapters and apply it to studying equivariant descent for both a smooth algebraic group acting on a variety and when a topological group (under mild assumptions) acts on a topological space. 

In Chapter \ref{Section: Category of Pseudocones} we begin our study of pseudocone categories $\PC(F)$ by first by introducing the equivariant derived category $D_G^b(X)$ of a variety $X$ as an example in Section \ref{Subsection: Equivariant Derived Cat of Var} (as well as some of the equivariant resolution theory required, minus technical details deferred to Chapter \ref{Section: Equivariant Descent}). We also examine the equivariant derived category of a topological space $D_G^b(X)$ in Section \ref{Subsection: EDC of a Space} as a second motivating example in order to provide a juxtaposition which illustrates how pseudocones arise in equivariant geometry and topology. In Section \ref{Subsection: Preliminary Results} we introduce the category of pseudocones, $\PC(F)$, in full detail and then establish some basic properties of said category (such as the facts that $\PC(F) \simeq \operatorname{PseudoLim}(F)$ in $\fCat$ and $\PC(F) = \mathbf{CSect}_B(\Cscr)^{\op}$) and different ways we can view it. In Section \ref{Subsection: Basic Properties} we define what it means for a pseudofunctor to be locally small with respect to a Grothendieck universe and use this notion to establish a basic property of $\PC(F)$ by determining when it is enriched in the $\Vscr$-models of Lawvere theories $\Tscr$.

Our study of $\PC(F)$ as a category continues in Chapter \ref{Section: Limits and Colimits in Pseudocones} where we study the various properties of $\PC(F)$ which are largely determined by limits, colimits, and other categorical arguments that are ``locally'' given by properties of the fibre categories $F(X)$ and the fibre functors $F(f)$. This process starts in Section \ref{Subsection: Limits and Colimits in PCF} where we study and determine the most natural\footnote{In the sense of ``naturally occurring in practice.''} situations when $\PC(F)$ has limits and colimits of specific shape, and hence also determine the most natural situations where $\PC(F)$ is complete, cocomplete, and even additive or Abelian. These observations allow us to recognize in purely category-theoretic terms that many commonly occurring categories of equivariant sheaves (such as the categories of equivariant local systems, equivariant $D$-modules, equivariant perverse sheaves, equivariant torsion sheaves, equivariant sheaves of modules, equivariant $\ell$-adic sheaves, etc.) are all Abelian. In Section \ref{Subsection: Monoidal Structures on PCF} we establish when $\PC(F)$ is a monoidal category, braided monoidal category, and symmetric monoidal category in terms of the local information recorded by the fibre categories $F(X)$ and the fibre translation functors $F(f)$. After discussing monoidal structures on $\PC(F)$, in Section \ref{Subsection: PCF Regualr} we discuss when the category $\PC(F)$ is regular and when it admits a subobject classifier in terms of the local categories $F(X)$ and the fibre translation functors $F(f)$.

After studying the categorical properties of pseudocone categories we change gears and study functors $\PC(F) \to \PC(E)$ in Chapter \ref{Section: Functors of Pseudocones}. These, of course, are functors of the form
\[
\Bicat(\Cscr^{\op},\fCat)(\cnst(\One),F) \to \Bicat(\Dscr^{\op},\fCat)(\cnst(\One),E).
\]
In Section \ref{Subsection: Change of Fibre} we study those functors which are given by post-composition by a pseudonatural transformation $F \Rightarrow E$ and use them to study adjoints in the $2$-category $\Bicat(\Cscr^{\op},\fCat)$. We then use adjoints in $\Bicat(\Cscr^{\op},\fCat)$ to produce adjoints between pseudocone categories in $\fCat$; in particular we classify when the categories $\PC(F)$ are monoidal closed based on local information determined by the fibre categories $F(X)$ and fibre functors $F(f)$ in order to give clean, rigourous and purely categorical proofs of the existence of the closed symmetric monoidal structures on the equivariant derived categories. We also use the language formulated here to discuss localizations of pseudocone categories in terms of the information encoded by the pseudofunctor in order to establish in what sense the equivariant derived category may be seen as a localization\footnote{Because equivariant derived categories are not the derived categories of equivariant sheaves, it begs the question of in what sense can we see them as a localization of the equivariant category of sheaves. Our results show that $D_G^b(X)$ is, in essence, determined by a pseudocone localization of $\Ch^b_G(X)$. See Proposition \ref{Prop: Section Functors: EDC is pc localization} for details
}. In Section \ref{Subsection: Change of Domain} we show how to formulate certain functors
\[
\Bicat(\Cscr^{\op},\fCat)(\cnst(\One),F) \to \Bicat(\Dscr^{\op},\fCat)(\cnst(\One),E)
\]
in terms of the functors studied in Section \ref{Subsection: Change of Fibre} and use them to construct the equivariant pullback, pushforward, exceptional pushforward, and exceptional pullback purely in terms of descent data. We then use the pseudocone formalism developed so far to give purely categorical and rigourous proofs of the adjunctions $h^{\ast} \dashv Rh_{\ast}$ and $Rh_! \dashv h^!$ in the topological and scheme-theoretic cases.

After studying functors between pseudocone categories, in Chapter \ref{Section: Pseudocone Algebra} we study the homological algebra of the categories $\PC(F)$. Our main focus in Section \ref{Section: Triangles are go} is twofold: first on establishing that when we have pseudofunctors $F:\Cscr^{\op} \to \fCat$ for which the fibre categories $F(X)$ are triangulated and the fibre functors $F(f)$ are also triangulated, then the categories $\PC(F)$ are triangulated categories; second, we focus on when we can deduce the existence of $t$-structures on the categories $\PC(F)$ in terms of $t$-structures on the local categories $F(X)$ and $t$-exact fibre functors $F(f)$. We use these observations to discuss the hearts of these $t$-structure and use them to illustrate the standard and perverse $t$-structures on $D_G^b(X)$ and $D_G^b(X;\overline{\Q}_{\ell})$ as examples. In Section \ref{Section: SFF for Pseudocones} we study when pseudocone categories $\PC(F)$ have six functor formalisms as induced by six functor formalisms on the local categories $F(X)$ and then use these to give as an example the relations between the functors that the six functor formalisms on the equivariant derived category satsify. We are also very careful to illustrate where the Axiom of Choice is important in giving the triangulated and truncation-structures on the categories $\PC(F)$. The issue which requires the Axiom of Choice is fundamentally a fact that triangulation and truncation structures are homotopical and higher-categorical structures, and while we do not pursue it in this monograph I conjecture that a good $\infty$-categorical extension of what we discuss here is a way to get around the ``unnatural'' usage of the Axiom of Choice required. Careful remarks about said generalizations are given in Section \ref{Section: Triangles are go}.

We change gears from pure category theory in Chapter \ref{Section: Equivariant Descent} and begin an exploration into how categories of pseudocones $\PC(F)$ aid us in studying and performing equivariant descent. In Section \ref{Section: Smooth Free Geometry} we establish the geometric facts that allow us to use the category of smooth free resolutions of $X$ as a home for the equivariant derived category of a variety $X$. Additionally, in Section \ref{Section: Smooth Free Geometry} we also introduce the notion of a pre-equivariant pseudofunctor $F = (F,\overline{F})$ on a variety or topological space\footnote{The definition used coincides with that used in \cite{DoretteMe}, save that it is extended to also handle the case of topological spaces. It differs from the definition I first introduced in the thesis \cite{MyThesis} mainly in that it is no longer formally problematic.} and consequently introduce the notion of an equivariant category $F_G(X) = \PC(F)$ associated to the pre-equivariant pseudofunctor. We then introduce the notion of the forgetful functor $F_G(X) \to \overline{F}(X)$ which generalizes the forgetful functors $\Shv_G(X) \to \Shv(X),$ $\Loc_G(X) \to \Loc(X)$, $\Shv_G(X;\overline{\Q}_{\ell}) \to \Shv(X;\overline{\Q}_{\ell})$, $D_G^b(X;\overline{\Q}_{\ell}) \to D_c^b(X;\overline{\Q}_{\ell})$, and so on. In Section \ref{Subsection: The Equivariance Isos} we show that the equivariant categories $F_G(X)$ come equipped with functors $\pi_2^{\ast}, \alpha_{X}^{\ast}:F_G(X) \to F_G(G \times X)$ and prove that  there is a natural isomorphism $\Theta:\alpha_X^{\ast} \xRightarrow{\cong} \pi_2^{\ast}$ which we call the equivariance isomorphism that satisfies a cocycle condition which specializes to the cocycle condition described in GIT (cf.\@ Definition \ref{Defn: GIT Cocycle}, \cite[Section 0.2]{BernLun}, \cite[Section 1.3]{GIT}) when $F$ is a pre-equivariant pseudofunctor of sheaves, perverse sheaves, local systems, and other pseudofunctors of similar flavours.

After establishing the existence of the equivariance natural isomorphism $\Theta$, in Chapter \ref{Chapter: Chofg} we study the Change of Groups functors $\varphi^{\sharp}:F_H(X) \to F_G(X)$ when we have a morphism $\varphi:H \to G$ of smooth algebraic groups and a pre-equivariant pseudofunctor $F = (F,\overline{F})$ with respect to the $H$-action on $X$. Section \ref{Section: Section 3: Change of Groups} mainly concerns the existence of the Change of Groups functors, how they arise in terms of descent information and manipulation of bibundles through the use of induction spaces, and give a criterion for determining when Change of Groups functors commute up to natural isomorphism with other functors induced by the pseudocone formalism. After handling these situations we also show in this section that the Change of Groups functors have compositor isomorphisms $(\psi \circ \varphi)^{\sharp} \cong \varphi^{\sharp} \circ \psi^{\sharp}$ and unit functors $1_G^{\sharp}:F_G(X) \to F_{1}(X)$ which are naturally isomorphic to the forgetful functor $F_G(X) \to \overline{F}(X)$. In Section \ref{Section: High Level Perspective Change of Groups} we provide a higher-level perspective on the Change of Groups functors by proving that they arise as a the translation functors of a pseudofunctor
\[
F_{-}(X)^{\op}:(\mathbf{SAlgGrp}\downarrow H)^{\op} \to \fCat
\]
as well as the transition functors of a pseudonatural transformation:
\[
\begin{tikzcd}
\GHPreq{-}{H} \ar[rr, bend left = 20, ""{name = U}]{}{(-)_{H}(X)} \ar[rr, bend right = 20, swap, ""{name = D}]{}{(-)_{\square}(X)} & & \fCat \ar[from = U, to = D, Rightarrow, shorten <= 4pt, shorten >=4pt]{}{(-)^{\sharp}}
\end{tikzcd}
\]
For details see Proposition \ref{Prop: Pseudonat Chofg perspective}. In Section \ref{Section: Quot and Ind Equiv} we show that there are induction equivalences $F_H(H \times^G X) \simeq F_G(X)$ whenever $H$ is a topological group and $G$ is a closed subgroup and whenever $H$ is a smooth algebraic group and $G$ is a smooth algebraic subgroup (cf.\@ Propositions \ref{Prop: Section IndQuot Equiv: Ind equiv for top} and \ref{Prop: Section IndQuot Equiv: Ind equiv for var}). We also show that whenever $X$ is a free $H$-object there are quotient equivalences $F_H(X) \simeq F_{H/G}(G \backslash X)$ whenever either $G \trianglelefteq H$ is a closed normal subgroup of topological group $H$ or when $G \trianglelefteq H$ is a smooth normal algebraic subgroup of the algebraic group $H$; see Propositions \ref{Prop: Section IndQuot Equiv: Quotient equiv for top} and \ref{Prop: Section IndQuot Equiv: Quotient equiv for variety} for details.  Afterwards, in Section \ref{Subsection: Averaging Functors} we use the induction equivalences to study left and right averaging functors, which are left and right adjoints to the Change of Groups $\operatorname{incl}^{\sharp}:F_H(X) \to F_G(X)$ whenever $\operatorname{incl}:G \to H$ is the inclusion of a closed subgroup $G \leq H$ or smooth algebraic subgroup $G \leq H$. 

In Chapter \ref{Section: Six Functor Formalism}, the final chapter of this monograph, we classify the equivariant categories $F_G(X)$ on a topological space or a variety which carry an open/closed distinguished triangle as well as introduce the basic theory for equivariant stalks and trace. Section \ref{Section: SFF} is focused on establishing the open/closed distinguished triangles for equivariant categories and consequently requires developing a certain amount of equivariant descent techniques for equivariant immersions before determining which equivariant categories have open/closed distinguished triangles and then proceeding to give the open/closed distinguished triangles for the categories $D_G^b(X)$ and $D_G^b(X;\overline{\Q}_{\ell})$ when $X$ is a variety and $G$ is a smooth algebraic group and $D_G^b(X)$ when $X$ is a topological space and $G$ is a topological group as examples. We also establish from purely categorical arguments the fact that the standard cohomology $H_G$, perverse cohomology ${}^{p}H_G$, and each of the six functors for the equivriant derived categories $D_G^b(X)$ and $D_G^b(X;\overline{\Q}_{\ell})$ commute up to natural isomorphism with the forgetful functors. 

After working through the equivariant six functor formalism, in Section \ref{Section: Trace and Stalks} we introduce a notion of stalk and trace for equivariant categories $F_G(X)$ with a focus on applications to representation theory. This necessitates a discussion on the different perspectives we can take to define a stalk for the equivariant category and a justification of the choice we make for our definition. After said discussion, we develop some of the basic theory of equivariant stalks and how they interact with the various pseudocone categories $F_G(X)$ and functors between them. As a consequence of having said stalks we introduce a notion of equivariant trace (in the sense of having a trace of an object $A$ of an equivariant category $F_G(X)$ at points $g$ of the group $G$ and $x$ of the object $X$ where the equation $gx = x$ holds) based on taking stalks of the equivariance $\Theta:\alpha_X^{\ast}(A) \to \pi_2^{\ast}(A)$ at the point $(g,x)$ of $G \times X$. After exploring the basic properties of this definition we spend the remainder of the chapter showing that equivariant trace is additive on distinguished triangles, stable under immersions, and continuous with respect to pullback when the pseudofunctor $\overline{F}$ is strict. 

Below we give a leitfaden detailing chapter dependencies. Dashed arrows indicate light dependencies while filled arrows indicate fundamental dependencies.
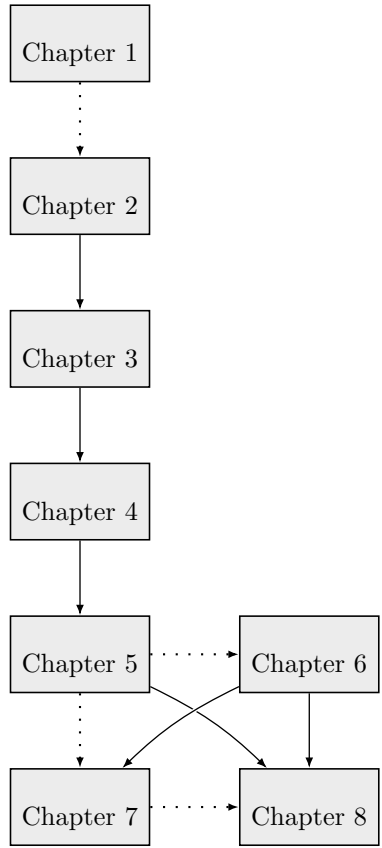
\begin{center}
\begin{chart}
	\begin{scope}[every rectangle node/.style={
			sharp corners,
			line width=0.6pt,
			fill=gray!15
		}]
		\reqhalfcourse 30,60:{}{Chapter 1}{}
		\reqhalfcourse 30,50:{}{Chapter 2}{}
		\reqhalfcourse 30,40:{}{Chapter 3}{}
		\reqhalfcourse 30,30:{}{Chapter 4}{}
		\reqhalfcourse 30,20:{}{Chapter 5}{}
		\reqhalfcourse 45,20:{}{Chapter 6}{}
		\reqhalfcourse 30,10:{}{Chapter 7}{}
		\reqhalfcourse 45,10:{}{Chapter 8}{}
	\end{scope}
	
	\prereq 30,50,30,40:
	\prereq 30,40,30,30:
	\prereq 30,30,30,20:
	\prereq 30,20,45,10:
	\prereq 45,20,45,10:
	\prereq 45,20,30,10:
	\coreq 30,60,30,50:
	\coreq 30,20,45,20:
	\coreq 30,20,30,10:
	\coreq 30,10,45,10:
\end{chart}
\end{center}

\newpage

\section{Notation and Other Conventions}\label{Subsection: Notation and Conventions}
In this short section we simply codify and make explicit some conventions used throughout this paper.
\begin{remark}\label{Remark: Notations and conventions}
	Here are some names and notations for specific (bi-or-2-or-1)-categories and how we will write them.
	\begin{itemize}
		\item We write $\Set$\index[notation]{Set@$\Set$} to denote the category whose objects are (small) sets and whose morphisms are total functions.
		\item Generically $1$-categories are written with script fonts (e.g., $\Cscr, \Dscr, \Ascr$) although some exceptions, such as named categories, are usually written in bold font (such as $\Set, \Ring, \Cring, \Ab,$ etc.).
		\item We write $\Ab$\index[notation]{Ab@$\Ab$} for the category of Abelian groups and, more generally, $\Ab(\Cscr)$\index[notation]{Ab(C)@$\Ab(\Cscr)$} to denote the category of Abelian group objects in a category $\Cscr$.
		\item We write $\Ring$\index[notation]{Ring@$\Ring$} for the category of rings with unit-preserving morphisms and $\Cring$\index[notation]{Cring@$\Cring$} for the category of commutative rings with unit-preserving morphisms. As in the case of $\Ab$, if $\Cscr$ is a category we write $\Ring(\Cscr)$\index[notation]{Ring(C)@$\Ring(\Cscr)$} and $\Cring(\Cscr)$\index[notation]{Cring(C)@$\Cring(\Cscr)$} for the category of ring objects and commutative ring objects, respectively.
		\item We write $\One$ for the terminal category, i.e., ``the'' category with exactly one object and one morphism; note that we also assume some fixed representative $\One = (\lbrace \ast \rbrace, \lbrace \id_{\ast} \rbrace)$ but do not care to specify the name of the element.
		\item If $(\Cscr,J)$ is a site we write $\Shv(\Cscr,J)$\index[terminology]{Shv(CJ)@$\Shv(\Cscr,J)$} to denote the category of sheaves of sets with respect to the Grothendieck topology $J$. Similarly, $\Ab(\Cscr,J)$\index[notation]{Ab(CJ)@$\Ab(\Cscr,J)$} (respectively) denotes the category of $J$-sheaves of Abelian groups on $\Cscr$.
		\item If $S$ is a scheme we write $\Sch_{/S}$ for the category of $S$-schemes.
		\item If $\Cscr$ is a category with a terminal object, generically we write $\top$\index[notation]{top!$\top$} for the terminal object in $\Cscr$. Similarly, if $\Cscr$ is a category with an initial object, generically we write $\bot$\index[notation]{bot!$\bot$} for the initial object. A notable exception is when $\Cscr$ has a zero object; in that case we write $0$ for the zero object (and hence for the terminal object) in $\Cscr$.
		\item We write $\Cat$\index[notation]{Cat@$\Cat$} for the $1$-category of small categories and functors and $\fCat$\index[notation]{CatbutFrak@$\fCat$} for the $2$-category of small categories, functors, and natural transformations. If we need the meta-categorical versions, we write $\CAT$\index[notation]{CAT@$\CAT$} and $\fCAT$\index[notation]{CATbutFrak@$\fCAT$} for the meta-category and meta-2-category of all locally small categories, respectively.
		\item Generically higher-categorical versions of named $1$-categories are denoted in fraktur fonts; e.g.\@ $\fCat$ for the $2$-category of categories and $\fTopos$ for the bicategory of toposes. The basic idea is script fonts denote categories with trivial $2$-and-higher cells and data while fraktur fonts denote categorical objects with potentially nontrivial $2$-and-higher cells and data.
		\item We assume that all pseudofunctors are normalized\index[terminology]{Pseudofunctor! Normalized} in the sense that $F(\id_X) = \id_{FX}$. By \cite[Expos{\'e} VI.9 Pages 180, 181]{SGA1} this may be done without loss of generality.
		\item If $\Cfrak$ is a bicategory and $f:X \to Y$ is a $1$-morphism we write $\iota_f$ for the identity $2$-cell on $f$.
		\item If $\Cscr$ is a category (or $2$-category or bicategory) we write $\Cscr_0$ for the class of objects in $\Cscr$, $\Cscr_1$ for the class of $1$-morphisms in $\Cscr$, and $\Cscr_2$ for the class of $2$-morphisms in $\Cscr$.
		\item If $\Cscr$ is a $1$-category we write $\Cscr(X,Y)$ for $X, Y \in \Cscr_0$ to denote the hom-set $\Hom_{\Cscr}(X,Y) = \lbrace f \in \Cscr_1 \; | \; \Dom(f) = X, \Codom(f) = Y \rbrace$.
		\item If $\Cscr$ is a $2$-category or a bicategory we write $\Cscr(X,Y)$ for $X, Y \in \Cscr_0$ to denote the category whose objects are $1$-morphisms $f:X \to Y$ and whose morphisms are $2$-morphisms with source and target $1$-morphisms with domain and codomain $X$ and $Y$.
		\item If $\Cfrak$ and $\Dfrak$ are bicategories, we write $\Bicat(\Cfrak, \Dfrak)$\index[notation]{Bicat@$\Bicat(-,-)$} to denote the bicategory whose objects are pseudofunctors $F:\Cfrak \to \Dfrak$, morphisms are pseudonatural transformations, and whose $2$-morphisms are modifications.
		\item In this paper all rings are assumed to have (multiplicative) identities and (additive) negatives.
		\item If $X$ is a sheafed space (a topological space equipped with a structure sheaf) we denote the underlying space of $X$ by $\lvert X \rvert$ and the structure sheaf by $\Ocal_X$.
		\item In this paper a variety (over a field $K$) is a reduced separated scheme of finite type over $K$, i.e., a scheme $X$ equipped with a structure map $\nu_X:X \to \Spec K$ such that $X$ is reduced and $f$ is separated and finite type. Note the our choice to work with reduced separated schemes of finite type for varieties is so that our varieties coincide with the definition given in, say, \cite{Vakil} and \cite{LusztigCuspidal2}. None of the results we use in this paper depend scheme-theoretically in any major way on the fact that our varieties are reduced; simply being separated and finite type (and quasi-projective; see below) is enough.
		\item In this paper we assume that all varieties are quasi-projective, i.e., there is an open immersion $j:X \to Y$ for some projective variety $Y$. In particular, every variety is a locally closed subscheme of some projective space $\Pbb_K^n$.
		\item If $A$ is a commutative ring we write $\Sch_{/A}$ as a shorthand for the slice category $\Sch_{/\Spec A}$. Similarly, the pullback of a cospan $X \rightarrow \Spec A \leftarrow Y$ is written $X \times_A Y$ in place of the more verbose $X \times_{\Spec A} Y$.
		\item If $K$ is a field we write $\Var_{/K}$ as a shorthand for the category $\Var_{/\Spec K}$ of varieties over $K$.
		\item If $K$ is a field we write $\mathbf{SAlgGrp}_{/K}$ for the category of smooth algebraic groups over $K$. This the category of group objects in $\Var_{/K}$ for which the structure map $G \to \Spec K$ is smooth.
		\item If $F:\Cscr^{\op} \to \fCat$ is a pseudofunctor and $X \xrightarrow{f} Y \xrightarrow{g} Z$ are morphisms in $\Cscr$ then we write $\phi_{f,g}$ for the compositor $2$-cell:
		\[
		\begin{tikzcd}
			& F(Y) \ar[dr]{}{F(f)} \\
			F(Z) \ar[ur]{}{F(g)} \ar[rr, swap, ""{name = D}]{}{F(g \circ f)} & & F(X) \ar[from = 1-2, to = D, Rightarrow, shorten <= 4pt, shorten >= 4pt]{}{\phi_{f,g}}
		\end{tikzcd}
		\]
	\end{itemize}	
\end{remark}

\section{Acknowledgments}
I would like to thank Clifton Cunningham for guidance in writing the material in this monograph in the original thesis version \cite{MyThesis} and for many helpful conversations and suggestions afterwards; James Steele for various conversations, suggestions, and interest in this work; the Vonganish group in general for helpful discussions and interest in this work and in particular both Andrew Fiori and Kristaps Balodis for asking thought-provoking questions; Kristine Bauer and Adam Topaz for an in-depth reading of an early version of this work while in its thesis form; Berndt Brenken for many interesting conversations, helpful edits, and suggestions on an early version of the text; Dorette Pronk for significant guidance, suggestions, advice, and interest in this work; Theo Johnson-Freyd for helpful conversations and questions; Bob Par{\'e} for noticing some very helpful categorical facts; the ATCAT seminar for allowing me to present some of the results in this monograph; and Deni Salja for helpful discussions regarding various aspects of the Grothendieck construction. Some of this work was supported while the author was an AARMS postdoctoral scholar at Dalhousie University.

\chapter{The Category of Pseudocones}\label{Section: Category of Pseudocones}
In this section we will get to know the main character in our paper: the category of pseudocones (with apex category $\mathbbm{1}$). While we will get to know some of the basic properties of these categories and their interaction with fibrations (by way of the Grothendieck construction), we first will see our motivation behind the pseuodocone perspective we take on equivariant geometry and topology. Note that in Subsection \ref{Subsection: Equivariant Derived Cat of Var} we cover algebro-geometric motivation while in Subsection \ref{Subsection: EDC of a Space} we cover algebro-topological motivation. Both these examples are relevant to representation theory (either directly to the theory of $p$-adic groups in the style of \cite{CFMMX} or in the topological and $D$-module theoretic flavour of \cite{BernLun}). However, category-theorists or other readers who needs no additional motivation to study the categories $\PC(F)$ of pseudocones can safely skip to Section \ref{Subsection: Preliminary Results} for the first introduction to our main character.

\section{The Equivariant Derived Category of a Variety}\label{Subsection: Equivariant Derived Cat of Var}
In this subsection we motivate our pseudocone construction by the study of the equivariant derived category of $\ell$-adic sheaves on a variety $X$ over a field $K$. The formulation we give of this category was developed by \cite{LusztigCuspidal2} and used variously in \cite{MirkVilDuality}, \cite{RepGradedHecke}, \cite{CFMMX}, \cite{RepGradedHecke}, and \cite{MyThesis}. We follow the definition as given in \cite{MyThesis} which itself mimics the definition of Lusztig in \cite{LusztigCuspidal2} and defines the equivariant derived category in terms of descent through certain resolutions and their quotients. To describe this we need to make some quick definitions before we can describe the equivariant derived category $D_G^b(X; \overline{\Q}_{\ell})$. 

In this section we let $K$ be a field and $G$ be a smooth algebraic group, i.e., a group object in the category $\Var_{/K}$ of $K$-varieties for which the structure map $G \to \Spec K$ is smooth. We also let $X$ be a (quasi-projective; cf.\@ Remark \ref{Remark: Notations and conventions}) left $G$-variety, i.e., a $K$-variety $X$ equipped with a morphism $\alpha_X:G \times X \to X$ of $K$-varieties such that the diagrams
\[
\begin{tikzcd}
	G\times (G \times X) \ar[d, swap]{}{\cong} \ar[rr]{}{\id_G \times \alpha_X} & & G \times X \ar[d]{}{\alpha_X} \\
	(G \times G) \times X \ar[dr, swap]{}{\mu_G \times \id_X} & & X \\
	& G \times X \ar[ur, swap]{}{\alpha_X}
\end{tikzcd}\qquad \begin{tikzcd}
	X \ar[r]{}{\cong} \ar[d, equals] & \Spec K \times X \ar[d]{}{1_G \times \id_X} \\
	X & G \times X \ar[l]{}{\alpha_X}
\end{tikzcd}
\]
commute. Equivariant morphisms of $G$-varieties are maps $f:X \to Y$ of left $G$-varieties for which the diagram
\[
\begin{tikzcd}
	G \times X \ar[r]{}{\alpha_X} \ar[d, swap]{}{\id_G \times f} & X \ar[d]{}{f} \\
	G \times Y \ar[r, swap]{}{\alpha_Y} & Y
\end{tikzcd}
\]
commutes. In order to describe the kinds of resolutions of the group action of $G$ on $X$ we will need to know about smooth free resolutions, which are in essence the {\'e}tale principal $G$-varieties.

\begin{definition}[{\cite[Section 1.9]{LusztigCuspidal2}}]\label{Defn: Pseudocone Section: Smooth free varieties}
	A left $G$-variety $\Gamma$ is a smooth free $G$-variety\index[terminology]{Smooth Free $G$-Variety} if:
	\begin{itemize}
		\item $\Gamma$ is smooth, pure dimensional, and admits a geometric quotient $q:\Gamma \to G \backslash \Gamma$ (for the definition of a geometric quotient see \cite[Definition 0.6]{GIT}).
		\item The geometric quotient 
		\[
		\quo:\Gamma \to G \backslash \Gamma
		\] 
		is ({\'e}tale) locally isotrivial with fibre $G$, i.e., there is an {\'e}tale cover $\lbrace \varphi_i:U_i \to G \backslash \Gamma \; | \; i \in I \rbrace$ of $G \backslash \Gamma$ such that in each pullback diagram
		\[
		\begin{tikzcd}
			\quo^{-1}(U_i) \ar[r] \ar[d] & \Gamma \ar[d]{}{q} \\
			U_i \ar[r, swap]{}{\varphi_i} & G \backslash \Gamma
		\end{tikzcd}
		\]
		there is an isomorphism $\quo^{-1}(U_i) \cong G \times U_i$, the corresponding map $G \times U_i \to U_i$ acts via projection $\pi_2:G \times U_i \to U_i$, and the {\'e}tale map $\varphi_i:U_i \to G \backslash \Gamma$ can be taken to be finite {\'e}tale.
	\end{itemize}
\end{definition}
\begin{remark}
Another way of saying that the morphism $\quo:\Gamma \to G \backslash \Gamma$ is ({\'e}tale) locally isotrivial is by saying that the morphism $\quo$ is {\'e}tale locally trivializable by a cover of finite {\'e}tale opens. We will generally argue with isotrivial morphisms with the perspective of {\'e}tale locally trivializable morphisms with a special condition on the corresponding {\'e}tale cover in this monograph.
\end{remark}
\begin{remark}
The careful reader may be wondering we we have chosen to ask for our smooth free varieties $\Gamma$ to have the quotient map $\quo:\Gamma \to G \backslash \Gamma$ be {\'e}tale locally trivializable as opposed to something like Zariski locally trivializable, fppf locally trivializable, fpqc locally trivializable, or something else. Our main reasons for working with {\'e}tale locally trivializable quotient maps $\quo:\Gamma \to G \backslash \Gamma$ is basically due to \cite[Exercise 6.1.1]{PramodBook}: there are affine smooth algebraic groups $G$ with smooth free $G$-varieties $\Gamma$ for which the quotient map $\quo:\Gamma \to G \backslash \Gamma$ is not Zariski locally trivializable. Asking for the the quotient to be {\'e}tale locally trivializable guarantees that the action $\alpha_X:G \times X \to X$ is free in the sense that the map
\[
\langle \alpha_X, \pi_2 \rangle: G \times X \to X \times X
\]
is a closed immersion. Finally asking for the quotient to be locally isotrivial is required so that we can apply Proposition \ref{Prop: Section 2.1: The quotient prop}. When $G$ is affine this poses no issue; by Remark \ref{Remark: Section Descent: Affine SfG Vars} $\Sf(G)$-varieties are exactly those pure dimensional smooth free $G$-varieties which admit a geometric quotient.
\end{remark}
\begin{definition}[{\cite[Section 1.9]{LusztigCuspidal2}}]\label{Defn: Pseudocone Section: Smooth free morphisms}
	A morphism of smooth free $G$-varieties\index[terminology]{Smooth Free $G$-Variety! Morphism} is a smooth $G$-equivariant morphism $f:\Gamma \to \Gamma^{\prime}$ where $\Gamma$ and $\Gamma^{\prime}$ are smooth free $G$-varieties and $f$ is a morphism of constant fibre dimension. Explicitly, for every $\gamma_0, \gamma_1 \in \lvert \Gamma \rvert$ we have $\dim f^{-1}(\gamma_0) = \dim f^{-1}(\gamma_1)$.
\end{definition}
This leads us to the definition of a category of smooth free $G$-varieties and the category of resolutions of a $G$-space $X$ by using smooth free varieties.
\begin{definition}\label{Defn: Pseudocone Section: Sf(G)}
	Define the category of smooth free $G$-varieties $\Sf(G)$ as follows:
	\begin{itemize}
		\item Objects: Smooth free $G$-varieties $\Gamma$ (cf.\@ Definition \ref{Defn: Pseudocone Section: Smooth free varieties}).
		\item Morphisms: Morphisms of smooth free $G$-varieties (cf.\@ Definition \ref{Defn: Pseudocone Section: Smooth free morphisms}).
		\item Composition and Identities: As in $\Var_{/K}$.
	\end{itemize}
\end{definition}
The category of smooth free $G$-resolutions of a $G$-variety is then the category ``$\Sf(G) \times X$.'' We'll justify this as a good category with which to do equivariant descent shortly, but we'll get to know the category of resolutions itself first.
\begin{definition}\label{Defn: Pseudocone Section: SfResl}
	Let $X$ be a left $G$-variety. We define the category $\SfResl_G(X)$ of smooth free $G$-resolutions of $X$ as follows.
	\begin{itemize}
		\item Objects: Maps $\pi_2:\Gamma \times X \to X$ for $\Gamma \in \Sf(G)_0$.
		\item Morphisms: Maps $f:\Gamma \times X \to \Gamma^{\prime} \times X$ of the form $f = h \times \id_X$ for $h \in \Sf(G)(\Gamma, \Gamma^{\prime})$.
		\item Composition and Identities: As in $\Var_{/K}$.
	\end{itemize}
\end{definition}
\begin{remark}\label{Remark: Pseudocone Section: SfResl is iso to Sf}
	From inspection of Definitions \ref{Defn: Pseudocone Section: Sf(G)} and \ref{Defn: Pseudocone Section: SfResl} there is an isomorphism of categories $\Sf(G) \cong \SfResl_G(X)$ for any $G$-variety $X$. This is a feature and not a bug; the presence of this isomorphism is crucial for defining both pullback and pushforward functors $h_{\ast}:D_G^b(X; \overline{\Q}_{\ell}) \to D_G^b(Y; \overline{\Q}_{\ell})$ and $h^{\ast}:D_G^b(Y;\overline{\Q}_{\ell}) \to D_G^b(X;\overline{\Q}_{\ell})$ when there is an equivariant morphism $h:X \to Y$ of $G$-varieties (cf.\@ Corollaries \ref{Cor: Pseudocone Functors: Existence of equivariant pullback  for schemes} and \ref{Cor: Pseudocone Functors: Pushforward functors for equivariant maps for scheme sheaves}).
\end{remark}

We will study the categories $\Sf(G)$ and $\SfResl_G(X)$ in more detail later in Chapter \ref{Section: Equivariant Descent}, but for the moment all we need to know about the category $\SfResl_G(X)$ and classical facts about the equivariant derived category are the following geometric facts:
\begin{enumerate}
	\item For every smooth free $G$-variety $\Gamma \in \Sf(G)_0$ and every $G$-variety $X$, the variety $\Gamma \times X$ admits a $G$-quotient variety $G \backslash (\Gamma \times X)$ which is functorial in $\SfResl_G(X)$ (cf.\@ Proposition \ref{Prop: Section 2.1: The quotient prop}). If $f \times \id_X:\Gamma \times X \to \Gamma^{\prime} \times X$ is a morphism in $\SfResl_G(X)$ then we write $\of:G \backslash (\Gamma \times X) \to G \backslash (\Gamma^{\prime} \times X)$
	\item The operation of taking quotients by $G$ gives a quotient functor $\quo:\SfResl_G(X) \to \Var_{/K}$.
	\item If $X$ is a free $G$-variety then there is an equivalence of categories 
	\[
	D_G^b(X;\overline{\Q}_{\ell}) \simeq D_c^b(G \backslash X; \overline{\Q}_{\ell}).
	\]
\end{enumerate}
These three facts (in particular Facts (2) and (3) above) allow us to describe the equivariant derived category by descending through the categories $D_c^b(G \backslash P; \overline{\Q}_{\ell})$ where $P \to X$ is a smooth free resolution of the $G$-action on $X$. This approach was first suggested in the topological case by Bernstein and Lunts in \cite{BernLun} (which we explore in Subsection \ref{Subsection: EDC of a Space} below) and later used in the geometric setting by Lusztig in \cite{LusztigCuspidal2}. Let us describe Lusztig's description of $D_G^b(X)$, extended to the equivariant derived category of $\ell$-adic sheaves.
\begin{definition}[{cf.\@ \cite[Section 1.10]{LusztigCuspidal2}}]
	The equivariant derived category $D_G^b(X; \overline{\Q}_{\ell})$ is described as follows:
	\begin{itemize}
		\item Objects: Pairs $(A, T_A)$ where $A$ is a collection
		\[
		A = \left\lbrace \quot{A}{\Gamma} \; | \; \quot{A}{\Gamma} \in D_c^b(G \backslash (\Gamma \times X))_0, \Gamma \in \Sf(G)_0 \right\rbrace
		\]
		(one object of each $D_c^b(G \backslash (\Gamma \times X); \overline{\Q}_{\ell})$ for each variety $\Gamma \times X$ in the category $\SfResl_G(X)$) and where $T_A$ is a collection of transition isomorphisms
		\[
		T_A = \left\lbrace \tau_{f}^{A}:\of^{\ast}(\AGammap) \to \AGamma \; | \; f:\Gamma \to \Gamma^{\prime}, f \in \Sf(G)_1 \right\rbrace
		\]
		subject to the following cocycle condition. Given a pair of composable morphisms $\Gamma \xrightarrow{f} \Gamma^{\prime} \xrightarrow{g} \Gamma^{\prime\prime}$ in $\Sf(G)$, the diagram
		\[
		\begin{tikzcd}
			\of^{\ast}\left(\overline{g}^{\ast}\left(\AGammapp\right)\right) \ar[rr]{}{\of^{\ast}(\tau_g^A)} \ar[d, swap]{}{\phi_{f, g}} & & \of^{\ast}\left(\AGammap\right) \ar[d]{}{\tau_f^A} \\
			(\overline{g} \circ \overline{f})^{\ast}\left(\AGammapp\right) \ar[rr, swap]{}{\tau_{g \circ f}^{A}} & & \AGamma
		\end{tikzcd}
		\]
		commutes where $\phi_{f,g}$ is a compositor natural isomorphism 
		\[
		\of^{\ast} \circ \overline{g}^{\ast} \Rightarrow (\overline{g} \circ \of)^{\ast}.
		\]
		\item Morphisms: A morphism $P:(A, T_A) \to (B, T_B)$ is a collection of maps
		\[
		P = \left\lbrace \quot{\rho}{\Gamma}:\AGamma \to \BGamma \; | \; \Gamma \in \Sf(G)_0, \AGamma \in A, \BGamma \in B, \quot{\rho}{\Gamma}:\AGamma \to \BGamma \right\rbrace,
		\]
		such that for any map $f:\Gamma \to \Gamma^{\prime}$ in $\Sf(G)$ the diagram
		\[
		\begin{tikzcd}
			\of^{\ast}\left(\AGammap\right) \ar[d, swap]{}{\tau_f^A} \ar[rr]{}{\of^{\ast}(\quot{\rho}{\Gamma^{\prime}})} & & \of^{\ast}\left(\BGammap\right) \ar[d]{}{\tau_f^B} \\
			\AGamma \ar[rr, swap]{}{\quot{\rho}{\Gamma}} & & \BGamma
		\end{tikzcd}
		\]
		commutes.
		\item Composition and Identities: Induced by the equations
		\[
		\left\lbrace \quot{\varphi}{\Gamma} \; | \; \Gamma \in \Sf(G)_0 \right\rbrace \circ \left\lbrace \quot{\rho}{\Gamma} \; | \; \Gamma \in \Sf(G)_0 \right\rbrace = \left\lbrace \quot{\varphi}{\Gamma} \circ \quot{\rho}{\Gamma} \right\rbrace
		\]
		and
		\[
		\id_{(A, T_A)} = \left\lbrace \id_{\AGamma} \; | \; \Gamma \in \Sf(G)_0 \right\rbrace.
		\]
	\end{itemize}
\end{definition}
\begin{remark}
	Strictly speaking the definition here differs from what is given by Lusztig in \cite{LusztigCuspidal2} and follows \cite{MyThesis} instead. The only difference is the presence of the compositor natural isomorphism $\phi_{f,g}:\of^{\ast} \circ \overline{g}^{\ast} \Rightarrow (\overline{g} \circ \of)^{\ast}$ in the cocycle condition; it is absent in \cite{LusztigCuspidal2}, but we make the convention that our pseudofunctors are not strict and so cannot omit it.
\end{remark}
In the definition above we see that we have a family of categories, namely the $D_c^b(G \backslash (\Gamma \times X);\overline{\Q}_{\ell})$, varying pseudofunctorially in $\SfResl_G(X)^{\op}$ with transition functors $\of^{\ast}:D_c^b(G \backslash (\Gamma \times X); \overline{\Q}_{\ell}) \to D_c^b(G \backslash (\Gamma \times X); \overline{\Q}_{\ell})$ and compositors $\phi_{f,g}$. Objects in $D_G^b(X;\overline{\Q}_{\ell})$ are given by specifying one $\AGamma$ for each $D_c^b(G \backslash (\Gamma \times X); \overline{\Q}_{\ell})$ and isomorphisms which translate along the morphisms in $\SfResl_G(X)^{\op}$ in a pseudonatural way. Because choosing the single objects $\AGamma$ correspond to choosing a functor 
\[
\AGamma:\mathbbm{1} \to D_c^{b}(G \backsim (\Gamma \times X);\overline{\Q}_{\ell})
\] 
and choosing the isomorphisms $\tau_f^A:\of^{\ast}(\AGammap) \to \AGamma$ with the cocycle condition corresponds to a family of natural transformations
\[
(\tau_{f}^{A})^{-1}:\AGamma \xRightarrow{\cong} (\of^{\ast} \circ \AGammap):\One \to D_c^b(G \backslash (\Gamma \times X); \overline{\Q}_{\ell})
\] 
which vary pseudonaturally in $\SfResl_G(X)^{\op}$. Consequently to describe an object $(A, T_A)$ is to describe a pseudocone of shape $D_{c}^{b}(G \backslash (\Gamma \times (-)); \overline{\Q}_{\ell})$ or, equivalently, a pseudonatural transformation from the constant pseudofunctor at $\One$ to the pseudofunctor $D_c^b(G \backslash ((-) \times X);\overline{\Q}_{\ell})$. In an analogous analysis, morphisms can be shown to correspond to modifications. Similar generalizations of this construction studied in \cite{MyThesis} (such as the equivariant derived category of perverse sheaves and the equivariant categories of: sheaves, Abelian sheaves, sheaves of modules, local systems, perverse sheaves, torsion modules, $D$-modules, etc.) all have the same basic structure: giving objects in these categories correspond to pseudocones and giving morphisms correspond to modifications. As such we end up finding that the categories of pseudocones with apex $\One$ generalizes and recaptures equivariant descent (in the sense of the study of equivariant categories in the style of \cite{MyThesis} and \cite{DoretteMe}) for varieties.

\begin{remark}\label{Remark: Pseudocone Section: Using SFG instead of Resls}
	The astute or experienced reader may wonder about our choice to work only with the smooth free $G$-resolutions of $X$ as opposed to the category of \emph{all} possible $G$-resolutions of $X$. Our choice is primarily motivated by following \cite{LusztigCuspidal2}, \cite{CFMMX}, and \cite{MyThesis} but the biggest reason to work with $\SfResl_G(X)$ in place of the category of all resolutions of $X$ is because $\SfResl_G(X)$ has particularly well-behaved objects and makes performing calculations simpler. 
	
	The good news is that this choice to work with $\SfResl_G(X)$ does not have any drawbacks in the sense that get equivalent equivariant derived categories. Define the category $\Resl_G(X)$\index[notation]{Reslg@$\Resl_G(X)$} as follows:
	\begin{itemize}
		\item Objects: Smooth principal $G$-varieties $P$ (smooth $G$-varieties with locally isotrivial quotient map $P \to G \backslash P$ with fibre $G$) equipped with smooth $G$-equivariant morphisms $p:P \to X$. 
		\item Morphisms: Smooth $G$-equivariant morphisms $f:P \to Q$ which give rise to a morphism
		\[
		\begin{tikzcd}
			P \ar[rr]{}{f} \ar[dr] & & Q \ar[dl] \\
			 & X
		\end{tikzcd}
		\] 
		in $\Var_{/X}$.
		\item Composition and Identities: As in $\Var_{/X}$.
	\end{itemize}
	Then $\SfResl_G(X)$ is a subcategory of $\Resl_G(X)$, and $\Resl_G(X)$ turns out to be the ``ultimate home'' for doing equivariant descent. Recall (cf.\@ \cite[Definition 6.1.18]{PramodBook}) that if $n \in \N$, a morphism $f:X \to Y$ of varieties is $n$-acyclic if and only if for every map $g:Z \to Y$ the base change map $\tilde{f}:X \times_Y Z \to Z$ has the property that for each perverse sheaf $\Fscr \in \Per(Z)_0$ the map
	\[
	\Fscr \to {}^{p}\tau^{\leq n}\left({}^{p}R\tilde{f}_{\ast}\left({}^{p}\tilde{f}^{\ast}\Fscr\right)\right)
	\]
	is an isomorphism. Note that there is a notion of $n$-acyclic map for the standard $t$-structure on $D_G^b(X;\overline{\Q}_{\ell})$, as opposed to the perverse $t$-structure, and \cite[Exercise 6.4.1]{PramodBook} provides a comparison between the notions and shows that they give equivalent categories. When $G$ is an affine smooth algebraic group there are $n$-acyclic $G$-resolutions $P_n \to \Spec K$ for all $n \in \N$ and these resolutions are all varieties in $\Sf(G)$; cf.\@ \cite[Proposition 6.1.23]{PramodBook} or \cite{MyThesis} for a proof of this fact. An important technical lemma (cf.\@ \cite[Lemma 6.4.8]{PramodBook}) allows us to determine that whenever $\Cscr \to \Resl_G(X)$ is a subcategory with the properties that:
	\begin{enumerate}
		\item The category $\Cscr$ admits finite products;
		\item For any $n \in \N$ there is an $n$-acyclic resolution $P \to X$ in $\Cscr$;
		\item The trivial resolution $\pi_2:G \times X \to X$ is an object of $\Cscr$;
	\end{enumerate}
	then, in the language of Definition \ref{Defn: Pseudocone Section: Pseudocone Cat} below, there is an equivalence of categories
	\[
	\PC\left(D_c^b\big(G \backslash (-)\big):\Cscr^{\op} \to \fCat\right) \simeq \PC\left(D_c^b\big(G \backslash (-)\big):\Resl_G(X)^{\op} \to \fCat\right)
	\]
	which extends to the $\ell$-adic sheaf case as well. In particular, because $\SfResl_G(X)$ satisfies all three conditions (the only nontrivial part is that $\Sf(G)$ contains $n$-acyclic resolutions for all $n \in \N$, but this given by \cite[Proposition 7.1.6]{MyThesis}) we get the equivalences of categories
	\begin{align*}
		D_G^b(X) &= \PC\left(D_c^b\big(G \backslash (-)\big):\Cscr^{\op} \to \fCat\right) \\
		&\simeq \PC\left(D_c^b\big(G \backslash (-)\big):\Resl_G(X)^{\op} \to \fCat\right).
	\end{align*}
	and
	\begin{align*}
		D_G^b(X; \overline{\Q}_{\ell}) &= \PC\left(D_c^b\big(G \backslash (-); \overline{\Q}_{\ell}\big):\Cscr^{\op} \to \fCat\right) \\
		&\simeq \PC\left(D_c^b\big(G \backslash (-); \overline{\Q}_{\ell}\big):\Resl_G(X)^{\op} \to \fCat \right).
	\end{align*}
\end{remark}
\newpage

\section{The Equivariant Derived Category of a Topological Space}\label{Subsection: EDC of a Space}

The case of the equivariant derived category of a topological space (as formulated by Bernstein and Lunts in \cite{BernLun}) is remarkably similar to the equivariant derived category of a variety, save for the fact that it has a more fibrational flavour and nature than the algebraic incarnation. Because of this, however, we will take a lighter touch in its exposition. Let $G$ be a topological group (not necessarily Hausdorff) and let $X$ be a left $G$-space (also not necessarily Hausdorff). Describing the category $D_G^b(X)$ is, as in the algebraic case, done by descending through the derived categories of $G$-resolutions of $X$. The definition of $D_G^b(X)$ itself, however, is more involved and not particularly suitable for our purposes\footnote{There are significant issues involving choice in the definition of $D_G^b(X)$ given in \cite{BernLun} which we would like to avoid.}, so we will give a more concrete categorical description using fibrations. In either case of giving a definition for $D_G^b(X)$, however, we need to discuss $G$-resolutions of $X$.

\begin{definition}[{\cite[Definition 2.1.1]{BernLun}}]\label{Defn: Free Gspace}
	Let $X$ be a left $G$-space. We say that $X$ is free if the quotient map $q:X \to G \backslash X$ is a locally trivializable fibration with fibre $G$. That is, there is an open cover $\lbrace U_i \to G \backslash X \; | \; i \in I \rbrace$  such that in the pullback diagrams
	\[
	\begin{tikzcd}
		q^{-1}(U_i) \ar[r]{}{p_{1,i}} \ar[d, swap]{}{p_{2,i}} & X \ar[d]{}{q} \\
		U_i \ar[r, swap]{}{} & G \backslash X
	\end{tikzcd}
	\]
	we have $q^{-1}(U_i) \cong G \times U_i$ and the corresponding diagram
	\[
	\begin{tikzcd}
	q^{-1}(U_i) \ar[r]{}{\cong} \ar[dr, swap]{}{p_{2,i}} & G \times U_i \ar[d]{}{\pi_2} \\
	& U_i
	\end{tikzcd}
	\]
	commutes.
\end{definition}
\begin{definition}
	Let $G$ be a topological group and let $X$ be a left $G$-space. Define the category $\FResl_G(X)$\index[notation]{FResl@$\FResl_G(X)$} of free $G$-resolutions of $X$ as follows:
	\begin{itemize}
		\item Objects: $G$-equivariant morphisms $p:P \to X$ where $P$ is a free $G$-space.
		\item Morphisms: $G$-equivariant morphisms $f:P \to Q$ such that the triangles
		\[
		\begin{tikzcd}
			P \ar[rr]{}{f} \ar[dr, swap]{}{p} & & Q \ar[dl]{}{q} \\
			& X
		\end{tikzcd}
		\]
		commute.
		\item Composition and Identities: As in $\Top_{/X}$.
	\end{itemize}
\end{definition}

As we alluded earlier, defining the category $D_G^b(X)$ directly in the style of \cite{BernLun} is not straightforward: the category $D_G^b(X)$ is defined, in a more modern language than what is used in\cite[Section 2.2]{BernLun}, as a certain pseudolimit which itself defines each of the objects in the diagram only up to equivalence and with a large shadow of the Axiom of Choice looming over the entire construction; while this is fixed to some degree by appending the term ``pseudo'' to everything, the actual working definition of an object in $D_G^b(X)$ under this definition is obtuse and difficult to work with by hand as a result. Instead, however, there is a fibrational construction given in \cite[Section 2.4]{BernLun} which does suit our needs and is a category not only with equivalent properties but also with better objects (which are given entirely by descent-theoretic information). Consider a pseudofunctor
\[
D^b(G \backslash (-)):\FResl_G(X)^{\op} \to \fCat
\]
induced by the assignments on objects and morphisms
\[
P \mapsto D^b(G \backslash P), \qquad f \mapsto f^{\ast}.
\]
Then it essentially follows from \cite[Section 2.4.3]{BernLun} that there is an equivalence $D_G^b(X) \simeq \mathbf{CSect}_B(p)^{\op}$ where $\mathbf{CSect}_B(p)$ is the category of based Cartesian sections of the elements fibration $p:\operatorname{El}(D_c^b(G \backslash (-))) \to \FResl_G(X)$ associated to our pseudofunctor. However, by Proposition \ref{Prop: Pseudocone Section: CSections are pseudolimits} below there is an equality of categories between the category of pseudocones of shape $D_c^b(G \backslash (-))$ with apex $\One$ and the category $\mathbf{CSect}_B(p)^{\op}$. Thus by studying pseudocones in the topological case as well we can learn about various properties of the equivariant derived category of $X$.

\begin{remark}
	It is worth noting that this begins to paint a picture that equivariant category theory, equivariant homological algebra, and equivariant geometry comes in two parts: first, in the part indexed and recorded by the categories of pseudocones (with apex $\One$) and second in the theory of resolutions we use to arrive at a good notion of equivariant descent. This philosophy has been leveraged in \cite{DoretteMe} to construct equivariant tangent structures on equivariant categories over varieties (and gives an answer for when you have tangents to an object which behave suitably with the group action) and the proofs there carry over mutatis mutandis to the topological case (and hence also to the case where a Lie group acts on a smooth manifold) as well by using this pseudocone technology; cf.\@ \cite[Section 5]{DoretteMe}.
\end{remark}
\newpage

\section{Preliminary Results and Relations with Fibrations}\label{Subsection: Preliminary Results}
Let us proceed to introduce, in broad strokes, the categories we will be studying in this monograph. Fix a category $\Cscr$ and a pseudofunctor $F:\Cscr^{\op} \to \fCat$. As stated earlier, we wish to study pseudocones of shape $F$ with apex $\One$. This consists of the following data:
\begin{itemize}
	\item For every $X \in \Cscr_0$ a functor $\quot{\alpha}{X}:\One \to F(X)$ such that if $f:X \to Y$ is any morphism in $\Cscr$, there is an invertible $2$-cell $\quot{\alpha}{f}$ in $\fCat$
	\[
	\begin{tikzcd}
		\One \ar[r, equals, ""{name = U}] \ar[d, swap]{}{\quot{\alpha}{Y}} & \One \ar[d]{}{\quot{\alpha}{X}} \\
		F(Y) \ar[r, swap, ""{name = D}]{}{F(f)} & F(X) \ar[from = U, to = D, Rightarrow, shorten <= 4pt, shorten >= 4pt]{}{\quot{\alpha}{f}} \ar[from = U, to = D, Rightarrow, shorten <= 4pt, shorten >= 4pt, swap]{}{\cong}
	\end{tikzcd}
	\]
	such that for any composable morphisms $X \xrightarrow{f} Y \xrightarrow{g} Z$ the pasting diagram
	\[
	\begin{tikzcd}
		\One \ar[r, equals, ""{name = UL}] \ar[d, swap]{}{\quot{\alpha}{Z}} & \One \ar[d]{}[description]{\quot{\alpha}{Y}}  \ar[r, equals, ""{name = UR}]{}{} & \One \ar[d]{}{\quot{\alpha}{Z}} \\
		F(Z) \ar[rr, swap, bend right = 40, ""{name = DD}]{}{F(g \circ f)} \ar[r, swap, ""{name = DL}]{}{F(g)} & F(Y) \ar[r, swap, ""{name = DR}]{}{F(f)} & F(X) \ar[from = 2-2, to = DD, Rightarrow, shorten <= 4pt, shorten >= 4pt]{}{\phi_{f,g}} \ar[from = UL, to = DL, Rightarrow, shorten <= 4pt, shorten >= 4pt]{}{\quot{\alpha}{g}} \ar[from = UR, to = DR, Rightarrow, shorten <= 4pt, shorten >= 4pt]{}{\quot{\alpha}{f}}
	\end{tikzcd}
	\]
	is equal to the pasting diagram
	\[
	\begin{tikzcd}
		& \One \ar[dr, equals] \\
		\One \ar[ur, equals] \ar[rr, equals, ""{name = U}] \ar[d, swap]{}{\quot{\alpha}{Z}} & & \One \ar[d]{}{\quot{\alpha}{X}} \\
		F(Z) \ar[rr, swap, ""{name = D}]{}{F(g \circ f)} & & F(X) \ar[from = U, to = D, Rightarrow, shorten <= 4pt, shorten >= 4pt]{}{\quot{\alpha}{f}} \ar[from = 1-2, to = U, Rightarrow, shorten <= 4pt, shorten >= 4pt]{}{\iota_{\id}}
	\end{tikzcd}
	\]
\end{itemize}
However, in writing down what it means to be a pseudocone with apex $\One$ we should note that this is exactly the information required of a pseudonatural transformation of pseudofunctors $\Cscr^{\op} \to \fCat$; to formalize this categorically we simply need to write the apex category of the pseudocone, $\One$, as a pseudofunctor $\Cscr^{\op} \to \fCat$ in some way. We can do this immediately, however, by using a constant pseudofunctor at $\One$.
\begin{definition}\index[terminology]{Constant Pseudofunctor}
	Let $\Cfrak$ and $\Dfrak$ be bicategories and let $Y \in \Dfrak_0$ be an object. We define the constant pseudofunctor at $Y$, $\cnst(Y):\Cfrak \to \Dfrak$\index[notation]{cnst@$\cnst(-)$} to be the psuedofunctor given by the assignments:
	\begin{itemize}
		\item For any object $X$ of $\Cfrak$, $\cnst(Y)(X) := Y$.
		\item For any $1$-morphism $f:X \to X^{\prime}$ of $\Cfrak$, $\cnst(Y)(f) := \id_Y$.
		\item For any $2$-morphism $\alpha:f \to g$ of $\Cfrak$, $\cnst(Y)(\alpha) := \iota_{\id}$.
	\end{itemize}
\end{definition}
\begin{remark}
	Note that whenever we discuss such a constant pseudofunctor, the bicategory in which the object $Y$ resides is considered part of the defining data. However, if $\Cscr$ is a category we will always assume if nothing is stated explicitly that $\cnst(\Cscr)$ is a pseudofunctor with values in $\fCat$.
\end{remark}

From the definition of the constant pseudofunctor $\cnst(\One):\Cscr^{\op} \to \fCat$ we see that to give a pseudocone of $F:\Cscr^{\op} \to \fCat$ with apex $\One$ it suffices to give a pseudonatural transformation $\alpha:\cnst(\One) \Rightarrow F:\Cscr^{\op} \to \fCat$. Consequently we use this as our tool to define the category of pseudocones; any distinction between the category of pseudocones and the hom-category of psuedonatural transformations and modifications from $\cnst(\One)$ to $F$ will be only semantic in nature and not structural in any practical way.

\begin{definition}\label{Defn: Pseudocone Section: Pseudocone Cat}\index[terminology]{Pseudocone Category}
	Let $\Cscr$ be a $1$-category, regarded as a (strict) bicategory with trivial $2$-cells, and let $F:\Cscr^{\op} \to \fCat$. We define the category of pseudocones of shape $F$ to be the category
	\[
	\PC(F) := \Bicat(\Cscr^{\op},\fCat)(\operatorname{cnst}(\mathbbm{1}),F).
	\]\index[notation]{PCF@$\PC(F)$}
\end{definition}
\begin{remark}\label{Remark: Pseudocone Section: Pseudocone explicit} 
	Let us unwrap the definition of $\PC(F)$ in order to give an explicit description of the objects and morphisms of $\PC(F)$, as working with $\PC(F)$ is much more straightforward with this knowledge. To this end recall that to give a functor $\quot{\alpha}{X}:\One \to F(X)$ is to pick out an object of $F(X)$ (and its identity morphism, of course). Similarly, if $f:X \to Y$ is a morphism in $\Cscr$ then the $2$-cell 
	\[
	\quot{\alpha}{f}:(\quot{\alpha}{X} \circ \id_{\One})(\ast) \xrightarrow{\cong} (F(f) \circ \quot{\alpha}{Y})(\ast)
	\]
	corresponds to giving an isomorphism $\quot{\alpha}{X}(\ast) \cong F(f)(\quot{\alpha}{Y}(\ast))$; we write 
	\[
	\tau_f:F(f)(\quot{\alpha}{Y}(\ast)) \xrightarrow{\cong} \quot{\alpha}{X}
	\] 
	for the corresponding inverse isomorphism and call these the transition isomorphisms (of the pseudocone) --- note that since the inverse of an isomorphism is unique, determining the $\quot{\alpha}{f}$ is equivalent to determining the $\tau_f$. Asking the pasting diagrams to coincide is equivalent to asking that for a composable pair $X \xrightarrow{f} Y \xrightarrow{g} Z$ of morphisms in $\Cscr$, the diagram
	\[
	\begin{tikzcd}
		F(f)\left(F(g)\quot{\alpha}{Z}(\ast)\right) \ar[d, swap]{}{\phi_{f,g}^{\quot{\alpha}{Z}(\ast)}} \ar[rr]{}{F(f)(\tau_g)} & & F(f)(\quot{\alpha}{Y}(\ast)) \ar[d]{}{\tau_f} \\
		F(g \circ f)\quot{\alpha}{Z}(\ast) \ar[rr, swap]{}{\tau_{g \circ f}} & & \quot{\alpha}{X}(\ast)
	\end{tikzcd}
	\]
	must commute. That is, the pseudonaturality condition gives rise to the cocycle condition 
	\[
	\tau_{g \circ f} \circ \phi_{f,g}^{\quot{\alpha}{Z}(\ast)} = \tau_{f} \circ F(f)(\tau_g).
	\]
	Similarly, if we have pseudocones $\alpha, \beta:\cnst(\One) \Rightarrow F$ a modification $\rho:\alpha \Rightarrow \beta$ is given by a family of natural transformations $\quot{\rho}{X}:\quot{\alpha}{X} \Rightarrow \quot{\beta}{X}$, one for each $X \in \Cscr_0$, such that for all morphisms $f:X \to Y$ in $\Cscr$ the diagrams of functors and natural transformations
	\[
	\begin{tikzcd}
		\quot{\alpha}{X} \ar[d, swap]{}{\quot{\alpha}{f}} \ar[rr]{}{\quot{\rho}{X}} & & \quot{\beta}{X} \ar[d]{}{\quot{\beta}{f}} \\
		F(f) \circ \quot{\alpha}{Y} \ar[rr, swap]{}{F(f) \ast \quot{\rho}{Y}} & & F(f) \circ \quot{\beta}{Y}
	\end{tikzcd}
	\]
	commute. However, this implies that we can replace the isomorphisms $\quot{\alpha}{f}$ with the corresponding transition isomorphisms to deduce that the diagrams
	\[
	\begin{tikzcd}
		F(f) \circ \quot{\alpha}{Y} \ar[rr]{}{F(f) \ast \quot{\rho}{Y}} \ar[d, swap]{}{\tau_f} & & F(f) \circ \quot{\beta}{Y} \ar[d]{}{\sigma_f} \\
		\quot{\alpha}{X} \ar[rr, swap]{}{\quot{\rho}{X}} & & \quot{\beta}{X}
	\end{tikzcd}
	\]
	also must commute. This allows us to give the following explicit description for $\PC(F)$:
	\begin{itemize}
		\item Objects: Pairs $(A, T_A)$ where 
		\[
		A = \lbrace \quot{A}{X} \; | \; X \in \Cscr_0, \quot{A}{X} \in F(X)_0 \rbrace
		\] 
		is a collection of objects and transition isomorphisms
		\[
		T_A = \lbrace \tau_f^{A}:F(F)(\quot{A}{Y}) \xrightarrow{\cong} \quot{A}{X} \; | \; f:X \to Y, f \in \Cscr_1 \rbrace
		\]
		which satisfy the cocycle condition
		\[
		\tau_{g \circ f}^A \circ \phi_{f,g} = \tau_f^A \circ F(f)\left(\tau_g^A\right)
		\]
		for any composable morphisms $X \xrightarrow{f} Y \xrightarrow{g} Z$ in $\Cscr$.
		\item Morphisms: A morphism $P:(A,T_A) \to (B, T_B)$ is a collection of morphisms 
		\[
		P = \lbrace \quot{\rho}{X}:\quot{A}{X} \to \quot{B}{X} \; | \; X \in \Cscr_0 \rbrace
		\]
		such that for any morphisms $f:X \to Y$ the diagram
		\[
		\begin{tikzcd}
			F(F)(\quot{A}{X}) \ar[rr]{}{F(f)(\quot{\rho}{X})} \ar[d, swap]{}{\tau_f^A} & & F(f)(\quot{B}{Y}) \ar[d]{}{\tau_f^B} \\
			\quot{A}{X} \ar[rr, swap]{}{\quot{\rho}{X}}	& & \quot{B}{X}
		\end{tikzcd}
		\]
		commutes.
	\end{itemize}
\end{remark}
\begin{remark}\label{Remark: Pseudocone Section: Direction of morphisms versus cones}
	It may seem strange at first glance to use the transition isomorphisms $\tau_f^A$ in place of the $\quot{\alpha}{f}$ in our explict description of $\PC(F)$. This is done primarily so that we can leverage our geometric intuition in more straightforward ways as well as make more direct comparisons with \cite{BernLun}, \cite{LusztigCuspidal2}, and other works which use the equivariant derived category; cf.\@ also Propositions \ref{Prop: Pseudocone Section: Sections are lax limit} and \ref{Prop: Pseudocone Section: CSections are pseudolimits} below. The main reasoning behind this discrepancy between the transition isomorphism/pseudocone direction (namely in defining $\tau_f^A = \quot{\alpha}{f}^{-1}$) is because the equivariant derived categories were developed with morphisms going in the geometric direction while the definition of pseudonatural transformations go in the algebraic direction. We have made the decision to stick with the geometric direction of morphisms so as to make the comparisons and applications to geometric and topological situations that motivate us more straightforward and immediate. 
\end{remark}
\begin{remark}
It is of general categorical interest in to drop the ``pseudo'' assumption in our pseudocone definitions and instead work with lax cones. If we do this we cannot use the trick where we chose $\tau_f = \quot{\alpha}{f}^{-1}$ as, of course, the witness morphisms in lax cones need not be invertible. Instead we simply must choose a convention when working with a category $\LC(F)$ of lax cones of $F$. In this monograph we choose to work with our witness morphisms working in the same direction as the pseudocone category (and so move in the geometric direction and not in the algebraic direction).
\end{remark}

In light of the remark above, we give an explicit definition of the category of lax cones of a pseudofunctor $F:\Cscr^{\op} \to \fCat$. While we will not work seriously with the category of lax cones in this monograph, it is relevant in locally to this section: the description of lax limits of $F$ and sections of the elements fibration rely on the category of lax cones.

\begin{definition}\index[terminology]{Lax Cone Category}
Let $F:\Cscr^{\op} \to \fCat$ be a pseudofunctor. The category of lax cones with apex $\One$ is defined to be\index[notation]{LCF@$\LC(F)$}
\[
\LC(F) := \Bicat_{\operatorname{lax}}\left(\Cscr^{\op},\fCat\right)(\cnst(\One),F)
\]
the category of lax natural transformations from $\cnst(\One)$ to $F$.
\end{definition}
\begin{remark}
As  with $\PC(F)$, we can give a more explicit description of $\LC(F)$. It is given as follows:
\begin{itemize}
	\item Objects: Pairs $(A, T_A)$ where 
	\[
	A = \lbrace \quot{A}{X} \; | \; X \in \Cscr_0, \quot{A}{X} \in F(X)_0 \rbrace
	\] 
	is a collection of objects and transition morphisms
	\[
	T_A = \lbrace \tau_f^{A}:F(F)(\quot{A}{Y}) \xrightarrow{} \quot{A}{X} \; | \; f:X \to Y, f \in \Cscr_1 \rbrace
	\]
	which satisfy the cocycle condition
	\[
	\tau_{g \circ f}^A \circ \phi_{f,g} = \tau_f^A \circ F(f)\left(\tau_g^A\right)
	\]
	for any composable morphisms $X \xrightarrow{f} Y \xrightarrow{g} Z$ in $\Cscr$.
	\item Morphisms: A morphism $P:(A,T_A) \to (B, T_B)$ is a collection of morphisms 
	\[
	P = \lbrace \quot{\rho}{X}:\quot{A}{X} \to \quot{B}{X} \; | \; X \in \Cscr_0 \rbrace
	\]
	such that for any morphisms $f:X \to Y$ the diagram
	\[
	\begin{tikzcd}
		F(F)(\quot{A}{X}) \ar[rr]{}{F(f)(\quot{\rho}{X})} \ar[d, swap]{}{\tau_f^A} & & F(f)(\quot{B}{Y}) \ar[d]{}{\tau_f^B} \\
		\quot{A}{X} \ar[rr, swap]{}{\quot{\rho}{X}}	& & \quot{B}{X}
	\end{tikzcd}
	\]
	commutes.
\end{itemize}
\end{remark}

For those familiar with the Grothendieck construction of a pseudofunctor $F:\Cscr^{\op} \to \fCat$, the categories $\LC(F)$ and $\PC(F)$ should look familiar by virtue of a comparison with the categories of elements of $F$; in Propositions \ref{Prop: Pseudocone Section: Sections are lax limit} and \ref{Prop: Pseudocone Section: CSections are pseudolimits} we make this comparison precise. Let us recall this construction; for an alternative presentation, see \cite{Vistoli}. Fix a pseudofunctor $F:\Cscr^{\op} \to \fCat$. The Grothendieck construction gives an equivalence (and in fact a $2$-equivalence of $2$-categories) between pseudofunctors $\Cscr^{\op} \to \fCat$ and cloven fibrations $\Escr \to \Cscr$. This allows us to view the pseudofunctor $F$ as a fibration $p:\El(F) \to \Cscr$ where $\El(F)$ is the category of elements of $F$. The category $\El(F)$\index[notation]{ElF@$\El(F)$} of elements of $F$\index[terminology]{Category of Elements} is defined as follows:
\begin{itemize}
	\item Objects: Pairs $(X, A)$ where $X \in \Cscr_0$ and $A \in F(X)_0$.
	\item Morphisms: A map $f:(X,A) \to (Y,B)$ is a pair $(f_0, f_1)$ where $f_0:X \to Y$ is a morphism in $\Cscr$ and $f_1:F(f_0)(B) \to A$ is a morphism in $F(X)$.
	\item Composition: Given morphisms $(X,A) \xrightarrow{(f_0, f_1)} (Y, B) \xrightarrow{(g_0, g_1)} (Z,C)$ the composite $(g_0, g_1) \circ (f_0, f_1)$ is defined by setting
	\[
	(g_0, g_1) \circ (f_0, f_1) := (g_0 \circ f_0, f_1 \circ F(f_0)(g_1) \circ \phi_{f_0, g_0}^{-1});
	\]
	note that the second morphism in the pair types correctly because it takes the form
	\[
	F(g_0 \circ f_0)(C) \xrightarrow{\phi_{f_0,g_0}^{-1}} F(f_0)\left(F(g_0)(C)\right) \xrightarrow{F(f_0)(g_1)} F(f_0)(B) \xrightarrow{f_1} A.
	\]
	\item Identities: The identity on $(X,A)$ is given by $(\id_{X},\id_A)$.
\end{itemize}
Note that for the identities to have the form expressed above we used that our psuedofunctors have been assumed to be normalized. If we instead have that $F(\id_{X}) \ne \id_{FX}$ we need to redefine $\id_A$ and use the unitor isomorphism $\phi_A:F(\id_X)(A) \to A$ in its place. 

The category $\El(F)$ comes naturally equipped with a projection $p:\El(F) \to \Cscr$ given by, simply, first projection. This functor is a fibration and gives one half of the equivalence between cloven fibrations and psuedofunctors. The other equivalence, which we will not use here, determines a pseudofunctor $F$ from a fibration $p:\Escr \to \Cscr$ by taking $F(X) := p^{-1}(\id_X)$ and then defining the translation functors $F(f)$ by taking Cartesian lifts which lie within a given cleavage; cf.\@ \cite{Vistoli} for details.

What is of interest to us lies in giving sections to the fibration $p:\El(F) \to \Cscr$. To give such a section $s:\Cscr \to \El(F)$ we must declare the following:
\begin{itemize}
	\item For each $X \in \Cscr_0$, an object $sX \in F(X)_0$. Such an object gives our object pair $s(X) := (X, sX)$ in $\El(F)$.
	\item For each morphism $f:X \to Y$ in $\Cscr$ we must determine a morphism $\tau_f:F(f)(sY) \to sX$. This gives our morphism $s(f) := (f, \tau_f)$.
	\item Asking for $s$ to be a functor asks that for all composable pairs of morphisms $X \xrightarrow{f} Y \xrightarrow{g} Z$ in $\Cscr$ the diagram
	\[
	\begin{tikzcd}
		F(f)\left(F(g)(sZ)\right) \ar[rr]{}{F(f)(\tau_g)} \ar[d, swap]{}{\phi_{f,g}} & & F(f)(sY) \ar[d]{}{\tau_f} \\
		F(g \circ f)(sZ) \ar[rr, swap]{}{\tau_{g \circ f}} & & sX
	\end{tikzcd}
	\]
	must commute. Similarly, we must also have that $\tau_{\id_X} = \id_{sX}$ for all $X \in \Cscr_0$.
\end{itemize}
However, such information is exactly that of a lax cone of shape $F$ with apex $\One$ by an earlier discussion of ours. 

Describing morphisms between sections is more difficult. We have to decide whether we work with \emph{all} possible natural transformations of sections or only with those which behave trivially over the base category $\Cscr$. In keeping with the fibrational theme, we elect to work with those sections which are invariant over the base in the following sense.
\begin{definition}
	Let $p:\Escr \to \Bscr$ be a fibration and let $s,t:\Bscr \to \Escr$ be sections of $p$. If $\alpha:s \Rightarrow t$ is a natural transformation such that the pasting diagram
	\[
	\begin{tikzcd}
		\Bscr \ar[rr, bend left = 20, ""{name = U}]{}{s} \ar[rr, bend right = 20, swap, ""{name = D}]{}{t} & & \Escr \ar[rr]{}{p} & & \Bscr \ar[from = U, to = D, Rightarrow, shorten <= 4pt, shorten >= 4pt]{}{\alpha}
	\end{tikzcd}
	\]
	is equal to the identity $2$-cell
	\[
	\begin{tikzcd}
		\Bscr \ar[rr, equals, bend left = 20, ""{name = U}]{}{} \ar[rr, equals, bend right = 20, swap, ""{name = D}]{}{} & & \Bscr \ar[from = U, to = D, Rightarrow, shorten <= 4pt, shorten >= 4pt]{}{\iota_{\id}}
	\end{tikzcd}
	\]
	then we say that $\alpha$ is a based natural transformation\index[terminology]{Based Natural Transformation} between sections. We write $\mathbf{Sect}_B(p)$\index[notation]{SectBp@$\mathbf{Sect}_B(p)$} for the category of sections of $p$ and based natural transformations.
\end{definition} 
Let us return to the case where we have a pseudofunctor $F:\Cscr^{\op} \to \fCat$ and two sections $s, t:\Cscr \to \El(F)$ of the fibration $p:\El(F) \to \Cscr$. When we have a based natural transformation
\[
\begin{tikzcd}
	\Cscr \ar[rr, bend left = 20, ""{name = U}]{}{s} \ar[rr, bend right = 20, swap, ""{name = D}]{}{t} & & \El(F) \ar[from = U, to = D, Rightarrow, shorten <= 4pt, shorten >= 4pt]{}{\alpha}
\end{tikzcd}
\]
we find from the condition $p \ast \alpha = \iota_{\id_{\Cscr}}$ that when we write each map 
\[
\alpha_X:(X,sX) \to (X,tX)
\]
as a pair $\alpha_X = (f,\rho_X)$ then $f = \id_X$ and $\rho_X$ is a morphism of the form
\[
\rho_X:tX \to sX;
\]
note that this follows from the fact that our pseudofunctors are normalized in the sense that $F(\id_X) = \id_{FX}$ so $F(\id_X)(tX) = \id_{FX}(tX) = tX$. Consequently, for any morphism $f:X \to Y$ in $\Cscr$ the commutativity of the naturality square
\[
\begin{tikzcd}
	(X,sX) \ar[d, swap]{}{(f, \sigma_f)} \ar[rr]{}{(\id_X,\rho_X)} & & (X, tX) \ar[d]{}{(f,\tau_f)} \\
	(Y, sY) \ar[rr, swap]{}{(\id_Y, \rho_Y)} & & (Y, tY)
\end{tikzcd}
\]
tells us that the diagram
\[
\begin{tikzcd}
	F(f \circ \id_X)(tY) \ar[dr, bend right = 20, equals]{}{} \ar[r, equals]{}{} & F(f)(tY) \ar[d, swap]{}{F(f)(\rho_Y)} \ar[r]{}{\tau_f} & tX \ar[d]{}{\rho_X} \\
	& F(f)(sY) \ar[r, swap]{}{\sigma_f} & sX
\end{tikzcd}
\]
commutes. However, this is exactly the condition that the collection of maps $P = \lbrace \rho_X:tX \to sX \; | \; X \in \Cscr_0 \rbrace$ is a morphism 
\[
P:(\lbrace tX \; | \; X \in \Cscr_0\rbrace, \lbrace \tau_f \; | \; f \in \Cscr_1\rbrace) \to (\lbrace sX \; | \; X \in \Cscr_0 \rbrace, \lbrace \sigma_f \; | \; f \in \Cscr_1 \rbrace).
\]
In this way we see that every based transformation of sections of $p$ is an $\op$-morphisms of lax cones. Perhaps unsurprisingly this construction works both ways and allows us to derive the following proposition.

\begin{proposition}\label{Prop: Pseudocone Section: Sections are lax limit}
	Let $F:\Cscr^{\op} \to \fCat$ be a pseudofunctor and let $p:\El(F) \to \Cscr$ be the associated elements fibration. Then $\mathbf{Sect}_B(p)^{\op} = \LC(F)$.
\end{proposition}
\begin{proof}
	The discussion prior to the statement of the proposition gives the subobject inclusions $\mathbf{Sect}_B(p)_0 \subseteq \LC(F)_0$ and $\mathbf{Sect}_B(p)_1^{\op} \subseteq \LC(F)_1$. However, the reverse inclusion $\LC(F)_0 \subseteq \mathbf{Sect}_B(p)_0$ is simply an unwinding of the definition using the objects $\quot{A}{X} := sX$ and morphisms $\tau_f^A := sf$. The inclusion $\LC(F)_1 \subseteq \mathbf{Sect}_B(p)^{\op}_1$ is shown similarly.
\end{proof}

We can capture the category $\PC(F)$ by using what we call Cartesian sections, i.e., sections $s:\Cscr \to \El(F)$ of the projection $p:\El(F) \to \Cscr$ for which each morphism $s(f) = (f,\tau_f)$ is a Cartesian arrow over $f$. In order classify these sections we need to know the Cartesian arrows over $f$. However, by construction each Cartesian arrow in $\El(F)$ over $f:X \to Y$ in $\Cscr$ is isomorphic to an arrow of the form $(f,\id_{F(f)A}):(X,F(f)A) \to (Y,A)$ for some $A \in F(Y)_0$.  Consequently the Cartesian arrows are of the form $(f,\tau_f):(X,B) \to (Y,A)$ where $\tau_f:F(f)A \to B$ is an isomorphism. This gives the characterization of $\PC(F)$ below.

\begin{proposition}\label{Prop: Pseudocone Section: CSections are pseudolimits}
	Let $F:\Cscr^{\op} \to \fCat$ be a pseudofunctor and let $p:\El(F) \to \Cscr$ be the associated elements fibration. Then if $\mathbf{CSect}_B(p)$ is the category of Cartesian sections of $p$ and based natural transformations then $\PC(F)$ is equal to the category $\mathbf{CSect}_B(p)^{\op}$.
\end{proposition}
\begin{proof}
	This follows mutatis mutandis to Proposition \ref{Prop: Pseudocone Section: Sections are lax limit} and the observation categorizing the Cartesian arrows above any given arrow in $\El(F)$.
\end{proof}

For various examples of the categories of pseudocones in algebraic geometry, see \cite[Example 3.2.7]{MyThesis}; we will recall some here, however. Various topological analogues of these can be given, and we present some below (in addition to other examples).
\begin{example}\label{Example: Pseudocone Section: Big list of examples for great justice}
	Let us consider the following various examples of $\PC(F)$ for various pseudofunctors $F$.
	\begin{enumerate}
		\item Let $K$ be a field, let $G$ be a smooth algebraic group over $K$, and let $X$ be a left $G$-variety. Then the category of ($\ell$-adic) equivariant perverse sheaves has a formulation in terms of a category of the form $\Per_G(X) = \PC(F)$\index[notation]{PerG@$\Per_G(X)$} by taking $F:\SfResl_G(X)^{\op} \to \fCat$ to be a pseudofunctor with object and morphism assignments
		\[
		\Gamma \times X \mapsto \Per(G \backslash (\Gamma \times X); \overline{\Q}_{\ell})
		\]
		and
		\[
		f \times \id_X \mapsto {}^{p}\of^{\ast}:\Per(G \backslash (\Gamma^{\prime} \times X); \overline{\Q}_{\ell}) \to \Per(G \backslash (\Gamma \times X); \overline{\Q}_{\ell}),
		\]
		respectively. That this is equivalent to the usual incarnation of the category of equivariant perverse sheaves is shown in \cite{PramodBook}, \cite{MyThesis}, and also in Corollary \ref{Cor: Section Triangle: Psedofunctor of perv is equiv perv}.
		\item Let $G$ be a smooth algebraic group over $K$ and let $X$ be a left $G$-variety. Then the category of equivariant local systems on $X$ is given by $\PC(L)$ where $L$ is a pseudofunctor $L:\SfResl_G(X)^{\op} \to \fCat$ with object and morphism assignments
		\[
		\Gamma \times X \mapsto \Loc(G \backslash (\Gamma \times X); \overline{\Q}_{\ell})
		\]
		and
		\[
		f \times \id_X \mapsto \of^{\ast}:\Loc(G \backslash (\Gamma^{\prime} \times X)) \to \Loc(G \backslash (\Gamma \times X)),
		\]
		respectively.
		\item Let $G$ be a smooth algebraic group over $K$ and let $X$ be a left $G$-variety. Then the equivariant derived category of perverse sheaves on $X$, $D_G^b(\Per(X))$\index[notation]{DGPer@$D_G^b(\Per(X))$} is determined by $\PC(F)$ for the pseudofunctor 
		\[
		F:\SfResl_G(X)^{\op} \to \fCat
		\] 
		where $F$ has object and morphism assignments
		\[
		\Gamma \times X \mapsto D^b(\Per(G \backslash (\Gamma \times X)))
		\]
		and
		\[
		f \times \id_X \mapsto L\left({}^{p}\of^{\ast}\right):D^b(\Per(G \backslash (\Gamma^{\prime} \times X))) \to D^b(\Per(G \backslash (\Gamma \times X))),
		\]
		respectively.
		\item Let $G$ be a topological group and let $X$ be a left $G$-space. The category of equivariant perverse sheaves on $X$ is determined by $\PC(F)$ for a pseudofunctor $F:\FResl_G(X)^{\op} \to \fCat$ where $F$ has object and morphism assignments
		\[
		P \mapsto \Per(G \backslash P)
		\]
		and
		\[
		f \mapsto {}^{p}\of^{\ast}:\Per(G \backslash Q) \to \Per(G \backslash P)
		\]
		respectively.
		\item Let $G$ be a smooth algebraic group and let $X$ be a left $G$-variety. We can then describe the category $D$-$\mathbf{Mod}_G(X)$\index[notation]{DModG@$D$-$\mathbf{Mod}_G(X)$} of $G$-equivariant $D$-modules as $\PC(D)$ where $D:\SfResl_G(X)^{\op} \to \fCat$ is a pseudofunctor with object and morphism assignments
		\[
		D(\Gamma \times X) := \DMod(G \backslash (\Gamma \times X))
		\]
		and
		\[
		D(f \times \id_X) := \of^{\ast}:\DMod(G \backslash (\Gamma^{\prime} \times X)) \to \DMod(G \backslash (\Gamma \times X))
		\]
		respectively. Note that in this case the functor $\of^{\ast}$ exists by \cite[Section 4.1]{benzvi2015character} because each morphism $\of$ is smooth by \@ Proposition \ref{Prop: Section 2.1: The quotient prop}.
		\item Let $G$ be a Lie group, regarded as a topological group, and let $X$ be a smooth manifold (regarded as a topological space) with a smooth $G$-action. Then the category of $G$-equivariant sheaves on $X$ is given by $\PC(F)$ where $F:\FResl_G(X)^{\op} \to \fCat$ is a pseudofunctor with object and morphism assignments
		\[
		P \mapsto \Shv(G \backslash P)
		\]
		and
		\[
		f \mapsto \of^{\ast}:\Shv(G \backslash Q) \to \Shv(G \backslash P),
		\]
		respectively.
		\item Let $(\Cscr,J)$ be a site and consider a pseudofunctor $S:\Cscr^{\op} \to \fCat$ defined on objects and morphisms by
		\[
		S(X) := \Shv(\Cscr_{/X}, J_{/X})
		\]
		and
		\[
		f \mapsto f^{\ast}:\Shv(\Cscr_{/Y},J_{/Y}) \to \Shv(\Cscr_{/X}, J_{/X}),
		\]
		respectively. Note that $f^{\ast}$ is the functor where, if $f:X \to Y$ is the morphism then for any sheaf $\Fscr \in \Shv(\Cscr_{/X},J_{/Y})$ the sheaf $f^{\ast}\Fscr$ acts on objects $u:U \to X$ by
		\[
		\Fscr(U \xrightarrow{u} X \xrightarrow{f} Y)
		\]
		and on $X$-morphisms $g:U \to V$ by evaluating $\Fscr$ on the morphism $g$ regarded as a commuting triangle:
		\[
		\begin{tikzcd}
			U \ar[rr]{}{g} \ar[dr, swap]{}{f \circ u} & & V \ar[dl]{}{f \circ v} \\
			& Y
		\end{tikzcd}
		\]
		Then $\PC(S)$ gives a category anti-equivalent to the global sections of the $J$-stack $\El(S) \to \Cscr$ which are effective on descent.
	\end{enumerate}
\end{example}

We close this section introducing the pseudocone categories $\PC(F)$ by showing that when $\Cscr$ has a terminal object there is an equivalence of categories $\PC(F) \simeq F(\top_{\Cscr})$ and use this to give a quick proof of the fact that for principal $G$-fibrations $X$ there is an equivalence of categories $D_G^b(X) \simeq D^b(G \backslash X)$ implied by the pseudocone formalism.

\begin{proposition}\label{Prop: Pseudocone Section: Terminal object in C gives us pseudocones as global sections}
	Let $\Cscr$ be a category with a terminal object and let $F:\Cscr^{\op} \to \fCat$ be a pseudofunctor. Then there is an equivalence of categories
	\[
	F(\top_{\Cscr}) \simeq \PC(F).
	\]
\end{proposition}
\begin{proof}
	We show this by providing two functors $L:F(\top_{\Cscr}) \to \PC(F)$ and $R:\PC(F) \to F(\top_{\Cscr})$ and then showing that $R \circ L \cong \id_{F(\top)}$ and $L \circ R \cong \id_{\PC(F)}$. Following alphabetic order, we first define the functor $L:F(\top_{\Cscr}) \to \PC(F)$ as follows. Given an object $A$ of $F(\top)$ we define the collection
	\[
	LA := \left\lbrace F(!_X)(A) \; | \; X \in \Cscr_0 \right\rbrace
	\]
	where $!_X:X \to \top$ is the unique map from $X$ to the terminal object in $\Cscr$. Then for any morphism $f:X \to Y$ we have the commuting diagram
	\[
	\begin{tikzcd}
		X \ar[rr]{}{f} \ar[dr, swap]{}{!_X} & & Y \ar[dl]{}{!_Y} \\
		& \top
	\end{tikzcd}
	\]
	in $\Cscr$. Consequently we define our transition isomorphisms by setting
	\[
	\tau_f^{LA} := \phi_{f, !_Y}^{A}:F(f)\left(F(!_Y)(A)\right) \to F(!_Y \circ f)(A) = F(!_X)(A).
	\]
	That the collection $T_{LA} = \lbrace \tau_{f}^{LA} \; | \; f \in \Cscr_1 \rbrace$ satisfies the cocycle condition follows from the fact that the $\phi_{-,-}$ are the compositor isomorphisms for the pseudofunctor $F$. The definition on morphisms is given similarly, i.e., for any morphism $\rho:A \to B$ in $F(\top)$ the map $L\rho$ is given by
	\[
	L\rho = \left\lbrace F(!_X)(\rho) \; | \; X \in \Cscr_0 \right\rbrace.
	\]
	That this is a functor is an immediate verification.
	
	Let us now define the functor $R:\PC(F) \to F(\top_{\Cscr})$. For objects $A = (A, T_A)$ define $RA$ by
	\[
	RA = R(\lbrace \quot{A}{X} \; | \; X \in \Cscr_0 \rbrace, \lbrace \tau_f^A \; | f \in \Cscr_1\rbrace) := \quot{A}{\top}
	\]
	and for morphisms $P = \lbrace \quot{\rho}{X} \; | \; X \in \Cscr_0 \rbrace$ we define $RP := \quot{\rho}{\top}$. This is even more straightforward than $R$ to verify to be a functor.
	
	To check that $R \circ L = \id_{F(\top)}$ we see that for any object $A \in F(\top)_0$,
	\begin{align*}
		(R \circ L)(A) &= R\left(\lbrace F(!_X)(A) \; | \; X \in \Cscr_0 \rbrace\right) = F(!_{\top})(A) = F(\id_{\top})(A) = \id_{F(\top)}(A) \\
		&= A
	\end{align*}
	and similarly for morphisms $\rho$. Thus $R \circ L = \id_{F(\top)}$.
	
	We now check that $L \circ R \cong \id_{\PC(F)}$. For this fix an object 
	\[
	\left(\left\lbrace \quot{A}{X} \; | \; X \in \Cscr_0 \right\rbrace, \left\lbrace \tau_f^A \; | \; f \in \Cscr_1 \right\rbrace\right)
	\]
	and note that we have
	\[
	\left(L \circ R\right)A = \left(\left\lbrace F(!_X)(\quot{A}{\top}) \; | \; X \in \Cscr_0 \right\rbrace, \left\lbrace \phi_{f,!_{Y}}^{\quot{A}{\top}} \; | \; f \in \Cscr_1, f:X \to Y\right\rbrace\right).
	\]
	We now claim that for there is an isomorphism $Z_A:(L \circ R)(A) \xrightarrow{\cong} A$ which is natural in $A$. For this we need to define $Z_A = \lbrace \quot{\zeta_A}{X} \; | \; X \in \Cscr_0 \rbrace$ where for each $X \in \Cscr_0$,
	\[
	\quot{\zeta_A}{X}:F(!_X)(\quot{A}{\top}) \to \quot{A}{X}.
	\]
	However, we simply set $\quot{\zeta_A}{X} := \tau_{!_X}^{A}$. To see that this is a morphism in the category $\PC(F)$ we need to show that for any morphism $f$ of $\Cscr$, say $f:X \to Y$, the diagram
	\[
	\begin{tikzcd}
		F(f)\left(\quot{((L \circ R)(A))}{Y}\right) \ar[rr]{}{F(f)(\quot{\zeta_A}{Y})} \ar[d, swap]{}{\tau_f^{(L \circ R)(A)}} & & F(f)\left(\quot{A}{Y}\right) \ar[d]{}{\tau_f^A} \\
		\quot{(L \circ R)(A)}{X} \ar[rr, swap]{}{\quot{\zeta_A}{X}} & & \quot{A}{X}
	\end{tikzcd}
	\]
	commutes. However, as this diagram is exactly the diagram
	\[
	\begin{tikzcd}
		F(f)\left(F(!_Y)(\quot{A}{\top})\right) \ar[rr]{}{F(f)(\tau_{!_Y}^{A})} \ar[d, swap]{}{\phi_{f,!_Y}^{\quot{A}{\top}}} & & F(f)(\quot{A}{Y}) \ar[d]{}{\tau_f^A} \\
		F(!_X)(\quot{A}{X}) \ar[rr, swap]{}{\tau_{!_X}^{A}} & & \quot{A}{X}
	\end{tikzcd}
	\]
	the commutativity of the diagram follows from the cocycle condition the transition isomorphisms of $A$ satisfy and the fact that $!_X = !_Y \circ f$. Thus $Z_A$ is an isomorphism in $\PC(F)$. Establishing that it is natural comes down to checking that if $P = \lbrace \quot{\rho}{X} \; | \; X \in \Cscr_0 \rbrace$ then the diagram
	\[
	\begin{tikzcd}
		(L \circ R)(A) \ar[r]{}{Z_A} \ar[d, swap]{}{(L \circ R)(P)} & A \ar[d]{}{P} \\
		(L \circ R)(B) \ar[r, swap]{}{Z_B} & B
	\end{tikzcd}
	\]
	commutes. But since this comes down to check the relevant diagram for each $X \in \Cscr_0$ and the $X$-local components have the form
	\[
	\begin{tikzcd}
		F(!_X)(\quot{A}{\top}) \ar[r]{}{\tau_{!_X}^{A}} \ar[d, swap]{}{F(!_X)(\quot{\rho}{\top})} & \quot{A}{X} \ar[d]{}{\quot{\rho}{X}} \\
		F(!_X)(\quot{B}{\top}) \ar[r, swap]{}{\tau_{!_X}^{B}} & \quot{B}{X}
	\end{tikzcd}
	\]
	and so the naturality of $Z$ is equivalent to the fact that morphisms in $\PC(F)$ are required to commute appropriately with the transition isomorphisms. It thus follows that $L \circ R \cong \id_{\PC(F)}$ via the natural isomorphism $Z$ and so, combined with the fact that $R \circ L = \id_{F(\top)}$, this gives the equivalence
	\[
	\PC(F) \simeq F(\top_{\Cscr}).
	\]
\end{proof}
We now use the proposition above to give short and completely categorical proofs of the facts that $D_G^b(X) \simeq D^b(G \backslash X)$ whenever either $G$ is a topological group $X$ is a free $G$-space or when $G$ is a smooth algebraic group and $X$ is an $\Sf(G)$-object.
\begin{corollary}\label{Cor: Pseudocone Section: Equivalence for derived cat free space and equiv cat}
	Let $G$ be a topological group and assume that $G$ acts freely on $X$. Then there is an equivalence of categories
	\[
	D_G^b(X) \simeq D^b(G \backslash X).
	\]
\end{corollary}
\begin{proof}
	Because $G$ acts freely on $X$, it follows that the identity morphism $\id_X:X \to X$ is a free $G$-resolution of $X$. Because $\FResl_G(X)$ is a subcategory of $\Top_{/X}$ it follows that $\id_X$ is a terminal object in $\FResl_G(X)$. Applying Proposition \ref{Prop: Pseudocone Section: Terminal object in C gives us pseudocones as global sections} to the pseudofunctor $D^b(G \backslash (-))$ of Subsection \ref{Subsection: EDC of a Space} gives the equivalence
	\[
	D^b(G \backslash X) = D^b(G \backslash (\id_X)) \simeq \PC(D^b(G \backslash (-))) = D_G^b(X).
	\]
\end{proof}
\begin{corollary}\label{Cor: Pseudocone Section: SfG variety quotient is edc}
	Let $G$ be an affine smooth algebraic group and let $X$ be a principal $G$-variety. Then $D_G^b(X) \simeq D_c^b(G \backslash X)$ and $D_G^b(X; \overline{\Q}_{\ell}) \simeq D_c^b(G \backslash X; \overline{\Q}_{\ell})$.
\end{corollary}
\begin{proof}
	The assumptions in the statement of the corollary imply that $\id_{X}:X \to X$ is a terminal object in $\Resl_G(X)$. Applying Proposition \ref{Prop: Pseudocone Section: Terminal object in C gives us pseudocones as global sections} and Remark \ref{Remark: Pseudocone Section: Using SFG instead of Resls} we then get that
	\begin{align*}
		D_c^b(G \backslash X) &= D_c^b(G \backslash \top_{\Resl_G(X)}) \\
		&\simeq \PC\left(D_c^b\big(G \backslash (-)\big):\Resl_G(X)^{\op} \to \fCat \right) \\
		&\simeq \PC\left(D_c^b\big(G \backslash (-)\big):\Cscr^{\op} \to \fCat\right) = D_G^b(X).
	\end{align*}
	The $\ell$-adic case follows similarly.
\end{proof}

We now show that the category $\PC(F)$ is the pseudolimit of the diagram $F:\Cscr^{\op} \to \fCat$ and that $\LC(F)$ is the lax limit of the diagram of shape $F$. Note that this gives a sanity check of the descriptions of lax limits and pseudolimits in \cite{GiraudCN}.
\begin{Theorem}\label{Thm: Section Pseudocones: PCF is the pseudolimit of F}
	Let $F:\Cscr^{\op} \to \fCat$ be a pseudofunctor. Then $\PC(F)$ is the pseudolimit of $F$ in $\fCat$.
\end{Theorem}
\begin{proof}
	We first show that $\PC(F)$ is a category with a pseudonatural transformation $\cnst(\PC(F)) \to F$, as that is what it means first to be a pseudocone over $F$. For this we first define the functors $p_X:\PC(F) \to F(X)$ for all $X \in \Cscr_0$  as the composite
	\[
	\begin{tikzcd}
		\PC(F) \ar[r]{}{p} & \prod\limits_{X \in \Cscr_0} F(X) \ar[r]{}{\pi_X} & F(X)
	\end{tikzcd}
	\]
	where $\pi_X$ is the usual projection out of the product and $p:\PC(F) \to \prod_{X \in \Cscr_0} F(X)$ is the functor which on objects $(A,T_A)$ discards the transition isomorphisms and sends $A$ to $(\quot{A}{X})_{X \in \Cscr_0}$ and on morphisms sends $P$ to $(\quot{\rho}{X})_{X \in \Cscr_0}$. Now define the natural isomorphism $\alpha_f:F(f) \circ p_Y \Rightarrow p_X$ by defining, for all objects $(A,T_A) \in \PC(F)_0$,
	\[
	\alpha_f^{A} := \tau_f^{A}.
	\]
	That this is a natural transformation is a straightforward check and follows immediately from the fact that the morphisms $P:(A, T_A) \to (B, T_B)$ must commute suitably with the transition morphisms in $T_A$ and $T_B$. Moreover, that these vary pseudonaturally in $\Cscr^{\op}$ also follows immediately from the cocycle condition
	\[
	\tau_{g \circ f}^{A} \circ \phi_{f,g} = \tau_{f}^{A} \circ F(f)(\tau_g^A)
	\]
	and hence gives that $\PC(F)$ is a pseudocone over $F$. We call this pseudocone $\alpha:\cnst(\PC(F)) \to F$.
	
	We now show that any pseudocone over $F$ in $\fCat$ factors through $\alpha$. To this end assume that we have a category $\Dscr$ equipped with functors $G_X:\Dscr \to F(X)$ for all $X \in \Cscr_0$ and natural isomorphisms
	\[
	\begin{tikzcd}
		& \Dscr \ar[dr, ""{name = R}]{}{G_Y} \ar[dl, swap, ""{name = L}]{}{G_X} \\
		F(X) & & F(Y) \ar[ll]{}{F(f)} \ar[from = R, to = L, Rightarrow, shorten <= 8pt, shorten >= 8pt]{}{\beta_f} \ar[from = R, to = L, swap, Rightarrow, shorten <= 8pt, shorten >= 8pt]{}{\cong}
	\end{tikzcd}
	\]
	for all $f \in \Cscr_1$ which vary pseudonaturally in $\Cscr^{\op}$. Now define a functor $G:\Dscr \to \PC(F)$ as follows. For an object $Z \in \Dscr_0$, define the collection
	\[
	GZ := \left\lbrace G_X(Z) \; | \; X \in \Cscr_0 \right\rbrace
	\]
	and define the transition isomorphisms $T_{GZ}$ by
	\[
	T_{GZ} := \left\lbrace \beta_f^{Z} \; | \; f \in \Cscr_1 \right\rbrace
	\]
	so that $\tau_f^{GZ} = \beta_f^{Z}$. To see that this indeed are transition isomorphisms, we simply use that $(G_X,\beta_f)$ describes a pseudonatural transformation $\cnst(\Dscr) \to F$ in $\Bicat(\Cscr^{\op},\fCat)$ and so we get that the equation
	\[
	\tau_{g \circ f}^{GZ} \circ \phi_{f,g} = \beta_{g \circ f}^{Z} \circ \phi_{f,g} = \beta_f^{Z} \circ F(f)\left(\beta_g^{Z}\right) = \tau_{f}^{GZ} \circ F(f)\left(\tau_g^{GZ}\right)
	\]
	holds. For morphisms $\rho:Z \to W$ in $\Dscr$ we similarly define $G\rho$ by
	\[
	G\rho := \left\lbrace G_X(\rho) \; | \; X \in \Cscr_0\right\rbrace,
	\]
	which is trivially checked to be a morphism in $\PC(F)$ by the naturality of the $\beta_f$. That this is a functor is then immediate to verify and follows from the fact that the $G_X$ are all functors.
	
	We now note that for any $X \in \Cscr_0$ and any $Z \in \Dscr_0$ we have that
	\[
	(p_X \circ G)(Z) = p_X\left(\left\lbrace G_X(Z)\; | \; X \in \Cscr_0 \right\rbrace, T_{GZ}\right) = G_X(Z)
	\]
	and similarly for morphisms. Thus $p \circ G = G_X$ for each $X \in \Cscr_0$. 
	
	Consider the pasting diagram
	\[
	\begin{tikzcd}
		& \Dscr \ar[ddr, bend left = 20]{}{G_Y} \ar[ddl, swap, bend right = 20]{}{G_X} \ar[d]{}[description]{G} & \\
		& \PC(F) \ar[dr, ""{name = R}]{}{p_Y} \ar[dl, swap, ""{name = L}]{}{p_X} \\
		F(X) & & F(Y) \ar[ll]{}{F(f)} \ar[from = R, to = L, Rightarrow, shorten <= 8pt, shorten >= 8pt]{}{\alpha_f}
	\end{tikzcd}
	\]
	where the empty cells denote the identity transformation. We claim this pastes to the $2$-cell describing $\beta_f$. To see this we simply note that since $\alpha_f$ picks out the transition isomorphism $\tau_f^A$ of an object in $\PC(F)$, for any $Z \in \Dscr_0$ we get
	\[
	\alpha_f^{GZ} = \tau_{f}^{GZ} = \beta_f^{Z}
	\]
	which proves the claim. Thus it follows that the pseudocone $\Dscr$ factors through $\PC(F)$.
	
	We now show the final remaining property of $\PC(F)$ being a pseudolimit. Assume that we have two pseudocones
	\[
	\begin{tikzcd}
		\Dscr \ar[r, shift left = 1]{}{G} \ar[r, shift right = 1, swap]{}{G^{\prime}} & \PC(F)
	\end{tikzcd}
	\]
	for which the $2$-cells
	\[
	\begin{tikzcd}
		& \Dscr \ar[dr, ""{name = R}]{}{p_Y \circ G} \ar[dl, swap, ""{name = L}]{}{p_X \circ G} \\
		F(X) & & F(Y) \ar[ll]{}{F(f)} \ar[from = R, to = L, Rightarrow, shorten <= 8pt, shorten >= 8pt]{}{\alpha_f \ast G}
	\end{tikzcd}
	\]
	and
	\[
	\begin{tikzcd}
		& \Dscr \ar[dr, ""{name = R}]{}{p_Y \circ G^{\prime}} \ar[dl, swap, ""{name = L}]{}{p_X \circ G^{\prime}} \\
		F(X) & & F(Y) \ar[ll]{}{F(f)} \ar[from = R, to = L, Rightarrow, shorten <= 8pt, shorten >= 8pt]{}{\alpha_f \ast G^{\prime}}
	\end{tikzcd}
	\]
	coincide. However, this implies that for every $Z \in \Dscr_0$ and $X \in \Cscr_0$ we have that
	\[
	\quot{GZ}{X} = p_X(G(Z)) = (p_X \circ G)(Z) = (p_X \circ G^{\prime})(Z) = p_X(G^{\prime}(Z)) = \quot{G^{\prime}Z}{X}
	\]
	and for every $f \in \Cscr_1$
	\[
	\tau_{f}^{GZ} = \alpha_f \ast G = \alpha_f \ast G^{\prime} = \tau_f^{G^{\prime}Z}.
	\]
	Thus we have that $GZ = (GZ, T_{GZ}) = (G^{\prime}Z, T_{G^{\prime}Z}) = G^{\prime}Z$ for all $Z \in \Cscr_0$. Checking that $G\rho = G^{\prime}\rho$ for all $\rho \in \Dscr_1$ is similar. It thus follows that $G = G^{\prime}$ and hence that there is a unique morphism $G \Rightarrow G^{\prime}$: the identity modification. From all this it follows that $\PC(F)$ is the pseudolimit of $F$ in $\fCat$.
\end{proof}
\begin{proposition}
	Let $F:\Cscr^{\op} \to \fCat$ be a pseudofunctor. Then $\LC(F)$ is the lax limit of $F$ in $\fCat$.
\end{proposition}
\begin{proof}
This follows mutatis mutandis to the proof of Theorem \ref{Thm: Section Pseudocones: PCF is the pseudolimit of F}.
\end{proof}
\newpage

\section{Basic Properties of the Category of Pseudocones}\label{Subsection: Basic Properties}
Let us close this introductory chapter by first introducing some terminology (so that we can talk succinctly about pseudofunctors that universally take values in additive categories and $\Vscr$-locally small categories, among other such things) before proceeding to develop some basic theory of pseudocone categories.

\begin{definition}\label{Defn: Additive, triangulated, and locally small pseudofunctors}
	Fix a pseudofunctor $F:\Cscr^{\op} \to \fCat$. If for all $X \in \Cscr_0$ and $f \in \Cscr_1$ the category $F(X)$ is additive and each functor $F(f)$ is additive, we call $F$ an additive pseudofunctor.\index[terminology]{Pseudofunctor! Additive}
\end{definition}
\begin{definition}\label{Defn: Pseudocone Section: Locally small pseudofunctors}
	Fix a pseudofunctor $F:\Cscr^{\op} \to \fCat$ and a Grothendieck universe $\Vscr$. If $\Cscr$ is $\Vscr$-small and for all $X \in \Cscr_0$ the category $F(X)$ is a $\Vscr$-locally small category we call $F$ a $\Vscr$-locally small pseudofunctor.\index[terminology]{Pseudofunctor! Locally small}
\end{definition}

We now show a quick smallness result; as such fix, for the moment, a Grothendieck universe $\Vscr$. With applications to additive and triangulated categories in mind, we need to know that $\PC(F)$ is a category in the same universe as $\Vscr$ whenever we have a $\Vscr$-locally small pseudofunctor $F:\Cscr^{\op} \to \fCat$ lest we have to jump universes just to define addition on hom-sets.

\begin{lemma}\label{Lemma: Section 2: Pseudofunctor values in Vlocally small means Equiv cat Vlocally small}
	Let $\Vscr$ be a Grothendieck universe and let $F:\Cscr^{\op} \to \fCat$ be a $\Vscr$-locally small pseudofunctor. Then $\PC(F)$ is $\Vscr$-locally small.
\end{lemma}
\begin{proof}
	Assume that $\Vscr$ has associated strongly inaccessible cardinal $\kappa$. Fix two arbitrary objects $(A,T_A), (B,T_B) \in \PC(F)_0$; note that it suffices to show that $\PC(F)(A,B)$ is a $\Vscr$-set, i.e., that $\lvert \PC(F)(A,B)\rvert < \kappa$.
	
	Note that by assumption $\lvert \Cscr_0 \rvert, \lvert \Cscr_1 \rvert < \kappa$. For what follows, consider the topos $\VSet$ of $\Vscr$-sets and let $\Omega \in \VSet_0$ be the corresponding subobject classifier. Because $\kappa$ is strongly inaccessible, the power objects $\Pcal_{\Vscr}(X) := [X,\Omega]$ in $\VSet_0$ satisfy $\lvert \Pcal_{\Vscr}(X)\rvert < \kappa$.
	
	Consider now that since $P$ is defined as the collection
	\[
	P = \lbrace \quot{\rho}{X}:\quot{A}{X} \to \quot{B}{X} \; |\; X \in \Cscr_0, \quot{A}{X} \in A, \quot{B}{X} \in B \rbrace
	\]
	we have that
	\[
	P \subseteq \bigcup_{X \in \Sf(G)_0} F(X)(\quot{A}{X},\quot{B}{X}).
	\]
	Because each of the categories $F(X)$ is locally $\Vscr$-small, each hom-collection \\$F(X)(\quot{A}{X},\quot{B}{X})$ is a $\Vscr$-set. Thus, since $\Cscr_0$ is a $\Vscr$-set, we have that the union
	\[
	Y := \bigcup_{X \in \Cscr_0} F(X)(\quot{A}{X},\quot{B}{X})
	\]
	is a $\Vscr$-set as well.
	
	Now note that for each $P \in \PC(F)(A,B)$, we have that $\lbrace P \rbrace \subseteq \Pcal_{\Vscr}(Y)$. Thus it follows that
	\[
	\Pcal_{\Vscr}(Y) \supseteq \bigcup_{P \in \PC(F)(A,B)}\lbrace P \rbrace = \PC(F)(A,B)
	\]
	and hence we have that $\PC(F)(A,B)$ is a $\Vscr$-set. This proves the lemma.
\end{proof}

We now conclude with one last lemma which shows that if $F:\Cscr^{\op} \to \fCat$ is a $\Vscr$-locally small pseudofunctor such that each category $F(X)$ is enriched in $\Tcal(\VSet)$ for a (potentially infinitary) Lawvere theory $\Tcal$ and each functor $F(f)$ is a $\Tcal(\VSet)$-functor then $\PC(F)$ is enriched in $\Tcal(\VSet)$ as well. The basic argument simply reduces to the fact that the various operations we need to be defined and coherences we need to satisfy are defined $X$-locally, they are defined globally as well. While we present the lemma below at a relatively large level of generality, we only need it explicitly to discuss when $\PC(F)$ is $\VAb$ enriched. As such we refer the reader unfamiliar with (infinitary) Lawvere theories to \cite[Appendix A]{brandenburg2021large} and also to \cite{BenThesis} for an introduction to the theory of limit sketches. All we need for what follows is that Lawvere theories $\Tcal$ (in the universe $\Vscr$) are categories $\Tcal$ with all $\Vscr$-indexed products and a distinguished object $L \in \Tcal_0$ such that for any object $K$ of $\Tcal$ there is a $\Vscr$ set $I$ and an isomorphism $K \cong L^I$ in $\Tcal$. Models of the theory $\Tcal$ (in a category $\Cscr$) are product-preserving functors $M:\Tcal \to \Cscr$ and morphisms of models are natural transformations between product-preserving functors.

\begin{remark}
	Because the object assignment of a model $M:\Tcal \to \Cscr$ of a Lawvere theory is largely determined by its value $M(L)$ we will usually refer to the model by the object $M(L)$ itself. In this vein we also will abuse notation and refer to a morphism of models $M \Rightarrow N$ as $\underline{\varphi}$ in reference to the morphism $\varphi = \varphi_{L}:M(L) \to N(L)$, as to a large degree this determines the natural transformation.
\end{remark}

\begin{lemma}\label{Lemma: Pseudocone Section: Lawvere Enrichment}
	Let $\Vscr$ be a Grothendieck universe and let $F:\Cscr^{\op} \to \fCat$ be a $\Vscr$-locally small pseudofunctor. Let $\Tcal$ be a Lawvere theory. If each category $F(X)$ is enriched in $\Tcal(\VSet)$ for all $X \in \Cscr_0$ then the category $\PC(F)$ is enriched in $\Tcal(\VSet)$ as well.
\end{lemma}
\begin{proof}
	Recall that in order to show that $\PC(F)$ is enriched in $\Tcal$ we need to show that each hom-set $\PC(F)(A,B)$ is a model of $\Tcal$ in $\VSet$ for all $A, B \in \PC(F)_0$. For this it suffices to show that there is a product-preserving functor $M:\Tcal \to \VSet$ for which $\Tcal(L) = \PC(F)(A,B)$ when $L$ is the distinguished object of the Lawvere theory. It then further suffices to show that if $S$ is any $\Vscr$-set then for all $s \in S$ the induced $S$-tuple
	\[
	\left(M(\pi_s):M(L^S) \to M(L)\right)_{s \in S}
	\]
	is a product diagram in $\VSet$.
	
	To define $M$ we use the full skeletal subcategory of $\Tcal$ consisting of powers $L^I$ for all $\Vscr$-sets $I$ for objects and taking our morphisms to be those morphisms between $L^I$ and $L^J$ for all $\Vscr$-sets $I$ and $J$ and then translate through the isomorphisms of objects $K \cong L^I$ in order to extend $M$ to $\Tcal$. For this we define $M$ as follows on powers of $L$
	\[
	M(L)^{I} := \PC(F)(A,B)^{I}
	\]
	and on morphisms $f:L^I \to L^J$ in $\Tcal$ by the assignment, for $(P_i)_{i \in I} \in \PC(F)(A,B)^I$
	\[
	M(f)(P_i)_{i \in I} := \left\lbrace \quot{M}{X}(f)(\quot{\rho_i}{X})_{i \in I} \; | \; X \in \Cscr_0 \right\rbrace
	\]
	where $\quot{M}{X}$ is the model functor 
	\[
	\quot{M}{X}:\Tcal \to \VSet
	\] 
	given by $\quot{M}{X}(L) = F(X)(\quot{A}{X}, \quot{B}{X})$ specified by the $\Tcal(\VSet)$-enrichment structure on $F(X)$. That this is a functor follows from the fact that each model $\quot{M}{X}$ is a functor. To see that $M$ preserves products fix a $\Vscr$-set $S$ and consider the diagram
	\[
	\left(M(\pi_s):M(L^S) \to M(L)\right)_{s \in S}
	\]
	and note that each morphism $M(\pi_s)$ takes the form
	\[
	M(\pi_s) = \left\lbrace \quot{M}{X}(\pi_s) \; | \; X \in \Cscr_0 \right\rbrace.
	\]
	Because each $\quot{M}{X}(\pi_s)$ is the projection from a product it follows from a straightforward verification\footnote{Alternatively, the reader who is willing to accept the poor form of forward-referencing a result can use Theorem \ref{Thm: Section 2: Equivariant Cat has lims} below.} that $M(\pi_s)$ is a projection morphism as well. Thus $M$ is product-preserving and $\PC(F)(A,B)$ is a model of $\Tcal(\VSet)$.
	
	We now must show that for any morphism $P:A \to B$ in $\PC(F)$, pre-and-post-composition by $P$ is a morphism of models of $\Tcal$. The pre-compositional argument is dual to the pre-compositional one, so we present only the post-compositional argument. Let $C$ be an object in $\PC(F)$ and recall that we must prove that the post-composition morphism $P_{\ast}:\PC(F)(C,A) \to \PC(F)(C,B)$ induces a morphism of models. 
	
	For this we take the same approach as above and only show the naturality of $\underline{P}_{\ast}$ along the skeletal subcategory of $\Tcal$ generated by the $L^I$. Let $f:L^I \to L^J$ be a morphism in $\Tcal$ and consider the induced morphisms of $\Vscr$-sets 
	\[
	\PC(F)(Z,A)(f):\PC(Z,A)^{I} \to \PC(Z,A)^{J}
	\] 
	and 
	\[
	\PC(F)(Z,B)(f):\PC(F)(Z,B)^{I} \to \PC(F)(Z,B)^{J}.
	\]
	By construction the maps $\PC(F)(Z,A)(f)$ all have the form 
	\[
	\PC(F)(Z,A)(f) = \left\lbrace F(X)(\quot{Z}{X}, \quot{A}{X})(f) \; | \; X \in \Cscr_0 \right\rbrace.
	\]
	Now write $P = \lbrace \quot{\rho}{X} \; | \; X \in \Cscr_0 \rbrace$ and note that for all $\Vscr$-sets $S$
	\[
	(\ul{P_{\ast}})_{L^J} = \left\lbrace\left(\ul{\quot{\rho}{X}_{\ast}}\right)_{L^S} \; | \; X \in \Cscr_0 \right\rbrace.
	\] 
	To see that $\ul{P}_{\ast}$ is natural we use that each map $\ul{\quot{\rho}{X}_{\ast}}$ is natural (as each category $F(X)$ is $\Tcal(\Vscr)$-enriched) and calculate
	\begin{align*}
		&\left(\ul{P_{\ast}}\right)_{L^J} \circ \PC(F)(Z,A)(f)\\ 
		&= \left\lbrace \left(\ul{\quot{\rho}{X}_{\ast}}\right)_{L^J} \; | \; X \in \Cscr_0 \right\rbrace \circ \left\lbrace F(X)(\quot{Z}{X}, \quot{A}{X})(f) \; | \; X \in \Cscr_0 \right\rbrace \\
		&= \left\lbrace \left(\ul{\quot{\rho}{X}_{\ast}}\right)_{L^J} \circ F(X)(\quot{Z}{X}, \quot{A}{X})(f) \; | \; X \in \Cscr_0\right\rbrace \\
		&= \left\lbrace F(X)(\quot{Z}{X}, \quot{B}{X})(f) \circ (\ul{\quot{\rho_{\ast}}{X}})_{L^I} \; | \; X \in \Cscr_0 \right\rbrace \\
		&= \left\lbrace F(X)(\quot{Z}{X}, \quot{B}{X})(f) \; | \; X \in \Cscr_0 \right\rbrace \circ \left\lbrace (\ul{\quot{\rho_{\ast}}{X}})_{L^I} \; | \; X \in \Cscr_0 \right\rbrace \\
		&= \PC(F)(Z,B) \circ \left(\ul{P}_{\ast}\right)_{L^I}.
	\end{align*}
	It thus follows that $\ul{P_{\ast}}$ is natural and so, after proceeding with a similar argument to show that $\ul{P^{\ast}}:\PC(F)(B,C) \to \PC(F)(A,C)$ is a natural transformation between models, we get that $\PC(F)$ is enriched over $\Tcal(\VSet)$.
	
\end{proof}
\begin{corollary}\label{Cor: Pseudocone Section: Ab Enrichment}
	Assume that $F:\Cscr^{\op} \to \fCat$ is a $\Vscr$-locally small pseudofunctor and assume that for all $X \in \Cscr_0$ the category $F(X)$ is enriched in $\VAb$. Then $\PC(F)$ is enriched in $\VAb$ as well.
\end{corollary}
\begin{proof}
	Simply apply Lemma \ref{Lemma: Pseudocone Section: Lawvere Enrichment} to the observation that $\VAb = \Ab(\VSet)$ where $\Ab$ denotes the Lawvere theory of Abelian groups.
\end{proof}
\newpage

\chapter{Limits, Colimits, and Structures in the Category of Pseudocones}\label{Section: Limits and Colimits in Pseudocones}

Some of the most important aspects when working with a category lie in its limits, colimits, and the structures it has available. For instance one of the reasons that toposes are often touted as a world in which to do logic and geometry lie in the fact that structurally toposes have not only limits and colimits, but also ample structure which is rich enough to be the world of not only higher order intuitionistic logic but also for sheaf-theoretic geometry and topology (in the sense that the categories $\Shv(\Cscr,J)$ of sets are all special kinds of toposes). Because in practice we need to know about these various properties and how they interact with descent (especially if one is interested in studying equivariant cohomology or working with objects like equivariant perverse sheaves, equivariant derived hom, or other such gadgets and functors) it becomes important to examine when the categories $\PC(F)$ admit limits, colimits, and other structures. Because of our motivation with the categories $D_G^b(X;\overline{\Q}_{\ell}), \Shv_G(X;\overline{\Q}_{\ell}),$ $\Per_G(X; \overline{\Q}_{\ell})$, and other categorical-geometric gadgets (such as the equivariant tangent categories studied in \cite{DoretteMe}) our main approach is to study limits, colimits, and structures in $\PC(F)$ which arise from limits, colimits, and structures in the fibre categories $F(X)$ and which of those limits, colimits, and structures are preserved by the translation functors $F(f)$.

\section{Limits and Colimits for Pseudocones}\label{Subsection: Limits and Colimits in PCF}
Our first goal in our development of pseudocone categories is to show that if each category $F(X)$ admits either a zero object or finite products and if each fibre functor $F(f)$ preserves these zero objects and products then $\PC(F)$. These, together with Corollary \ref{Cor: Pseudocone Section: Ab Enrichment}, will allow us to conclude that the category $\PC(F)$ is additive whenever each category $F(X)$ is and each of the functors $F({f})$ are additive. The first result we present below is admittedly a special case of Theorems \ref{Thm: Section 2: Equivariant Cat has lims} and \ref{Theorem: Section 2: Equivariant Cat has colims}, but we present it in isolation to give a straightforward example of how these proofs work before diving into more general cases.

\begin{lemma}\label{Lemma: Section 2: Equivariant Cat has zero objects}
	Let $F:\Cscr^{\op} \to \fCat$ and assume that for all $X \in \Cscr_0$ and for all $f \in \Cscr_1$ each category $F(X)$ has a zero object $\quot{0}{X}$ and that each functor $F({f})$ preserves each corresponding zero object. Then the object
	\[
	0 = \lbrace \quot{0}{X} \; | \; X \in \Cscr_0\rbrace
	\]
	is a zero object in $\PC(F)$.
\end{lemma}
\begin{proof}
	Begin by defining $(0,T_0)$ by setting the collection $0$ as above and defining
	\[
	T_0 := \lbrace \tau_f^0:F({f})(\quot{0}{Y}) \xrightarrow{\cong}  \quot{0}{X} \; | \; f \in \Cscr_1 \rbrace
	\]
	as the collection of unique ismorphisms witnessing the fact that the $F(f)$ preserve the zero object. 
	
	We will first verify that $0$ is terminal in $\PC(F)$; that it is initial follows by duality and will be omitted. First let $(A,T_A)$ be an aribtrary object in $\PC(F)$. Define the map $!_{A}:(A,T_A) \to (0,T_0)$ by taking the $X$-local collection of the terminal maps $!$ on each $\quot{A}{X}$; explicitly
	\[
	!_{A} := \lbrace !_{\quot{A}{X}}:\quot{A}{X} \to \quot{0}{X} \; | \; X \in \Cscr_0 \rbrace.
	\]
	That $!_{A}$ is a morphism in $\PC(F)$ is clear, as the diagrams
	\[
	\xymatrix{
		F(f)(\quot{A}{Y}) \ar[r]^-{\tau_f^A} \ar[d]_{!_{F(f)(\quot{A}{Y})}} & \quot{A}{X} \ar[d]^{!_{\quot{A}{X}}} \\
		F(f)(\quot{0}{Y}) \ar[r]_{\tau_f^0} & \quot{0}{X}
	}
	\]
	all commute for all $f \in \Cscr_1$ and for all $X \in \Cscr_0$. That this map $!_{A}$ is unique follows from the fact that for any other morphism $P:A \to 0$, we would have by each $\quot{0}{X}$ being terminal that $\quot{\rho}{X} = !_{\quot{A}{X}}$ and hence that $!_{A} = P$. This proves the lemma.
\end{proof}
By reusing the proof of Lemma \ref{Lemma: Section 2: Equivariant Cat has zero objects} we get the following.
\begin{lemma}
	If each category $F(X)$ has a terminal (or initial) object and each functor $F({f})$ preserves it, then $\PC(F)$ has a terminal (or initial) object.
\end{lemma}

For what follows below we will prove when it is possible to deduce that the equivariant category has limits over diagrams of specific shape based on when limits of the same shape exist in each fibre category $F(X)$. This will allow us to deduce, say in the cases of equivariant quasi-coherent sheaves, perverse sheaves, local systems, $\Ocal_X$-modules, and other examples where the fibre functors $F({f})$ are exact/left exact, that the category $\PC(F)$ has all finite limits (or limits of specific type; cf.\@ Corollaries \ref{Cor: Section 2: Limits of shape I in equiv cat} and \ref{Cor: Equivariant cat is finitely complete or cocomplete}).
\begin{Theorem}\label{Thm: Section 2: Equivariant Cat has lims}
	Let $F:\Cscr^{\op} \to \fCat$ be a pseudofunctor and let $I$ be an index category for which there is a diagram $d:I \to \PC(F)$. For each $X \in \Cscr_0$ let $d_{X}:I \to F(X)$ denote the diagram functor
	\[
	\xymatrix{
		I \ar[drr]_{d_{X}} \ar[r]^-{d} & \PC(F) \ar[r]^-{\tilde{\imath}} & \prod\limits_{X \in \Cscr_0} F(X) \ar[d]^{\pi_{X}} \\
		& & F(X)
	}
	\]
	and assume that for each $X \in \Cscr_0$ the limit $\lim d_{X}$ exists in $F(X)$. Furthermore, assume that for every morphism $f \in \Cscr_1$ with $f:X \to Y$ there is an isomorphism $\theta_f$ witnessing the limit preservation
	\[
	F({f})\left(\lim_{\substack{\longleftarrow \\ i \in I}}d_{X}(i)\right) \cong \lim_{\substack{\longleftarrow \\ i \in I}} F({f})(d_{X}(i)).
	\]
	Then the diagram $d$ admits a limit in $\PC(F)$.
\end{Theorem}
\begin{proof}
	Let $X \xrightarrow{f} Y \xrightarrow{g} Z$ be a pair of composable arrows in $\Cscr$. Recall the natural isomorphism $\phi_{f,g}:F(f) \circ F(g) \Rightarrow F(g\circ f)$; we will abuse notation and write
	\[
	\phi_{f,g}^{A_i}:F(f)\big(F(g)d_{Z}(i)\big) \xrightarrow{\cong} F(g \circ f)\left(d_{Z}(i)\right)
	\]
	in place of the more cumbersome $\phi_{f,g}^{d_{Z}(i)}$. Define the collection of objects
	\[
	A := \left\lbrace \lim_{\substack{\longleftarrow \\ i \in I}} d_{X}(i) \; : \; X \in \Cscr_0 \right\rbrace.
	\]
	To see how to define the transition isomorphisms $\tau_f^A$ for all $f \in \Cscr_1$, recall that for all $f:X \to Y$ there are natural isomorphisms 
	\[
	\theta_f: F(f)\left(\lim(d_{Y})\right) \xrightarrow{\cong} \lim(F(f) \circ d_{Y}).
	\] 
	Then for all $i \in I_0$ we have commuting diagrams
	\[
	\xymatrix{
		\lim(F(f)\circ d_{Y}) \ar[d]_{\pi_{i}^{\prime}} \ar@{-->}[rr]^-{\exists!\lim(\tau_f^{A_i})} & & \lim(d_{X}) \ar[d]^{\pi_i} \\
		F(f)(d_{Y}(i)) \ar[rr]_-{\tau_f^{A_i}} & & d_{X}(i)
	}
	\]
	and, because each $\tau_f^{A_i}$ is an isomoprhism, so is $\lim(\tau_f^{A_i})$. We thus define the transition isomorphism $\tau_f^{A}:F(f)\left(\quot{A}{Y}\right) \xrightarrow{\cong} \quot{A}{X}$ via the composition:
	\[
	\xymatrix{
		F(f)(\lim(d_{Y})) \ar[r]^-{\theta_f} \ar[dr]_{\tau_f^A} & \lim(F(f) \circ d_{Y}) \ar[d]^{\lim(\tau_f^{A_i})} \\
		& \lim(d_{X})
	}
	\]
	Note that because $\theta_f$ and $\lim(\tau_f^{A_i})$ are isomorphisms so is the composite $\tau_f^A$. 
	
	We now need to verify the cocycle condition
	\[
	\tau_f^A \circ F(f)\tau_g^A = \tau_{g \circ f}^{A} \circ \phi_{f,g}^{\lim A_i}
	\]
	for the composable pair of morphisms $f$ and $g$. With this in mind we will decompose the isomorphism $\theta_{g \circ f}$ in order to get an algebraic identity that will allow us to prove the cocycle condition. Consider that by definition
	\[
	\tau_{g \circ f}^{A} = \lim(\tau_{g \circ f}^{A_i}) \circ \theta_{g \circ f}
	\]
	and by the uniqueness of the ismorphism $\theta_{g \circ f}$ (as it is induced from the universal property of a limit) we find that the diagram (where we recall that $d_{Z}(i) = \quot{A_i}{Z}$ and $\lim \quot{A_i}{Z} = \quot{A}{Z}$)
	\[
	\xymatrix{
		F(g \circ f)(\lim(\quot{A_i}{Z})) \ar[rr]^-{\theta_{g \circ f}}_-{\cong} \ar[d]_{\left(\phi_{f,g}^{\lim A_i}\right)^{-1}}^{\cong} & & \lim(F(g \circ f)(\quot{A_i}{Z})) \\
		F(f)\big(F(g)(\lim(\quot{A_i}{Z}))\big) \ar[d]_{F(f)\theta_g}^{\cong} \\
		F(f)\big(\lim F(g)(\quot{A_i}{Z})\big)\ar[d]_{F(f)\big(\lim \tau_g^{A_i}\big)}^{\cong} & & \lim F(f)\big(F(g)(\quot{A_i}{Z}) \big) \ar[uu]^{\cong}_{\lim \phi_{f,g}^{A_i}} \\
		F(f)\big(\lim (\quot{A_i}{Y})\big) \ar[rr]_-{\theta_f}^-{\cong} & & \lim\big( F(f)(\quot{A_i}{Y})\big) \ar[u]^{\cong}_{\big(\lim(F(f)\tau_g^{A_i})\big)^{-1}}
	}
	\]
	commutes in $F(X)$.
	We now consider the morphism $\psi:\lim\big(F(g \circ f)\quot{A_i}{Z}\big) \to \lim\quot{A_i}{X}$ defined by the composition:
	\[
	\xymatrix{
		\lim\big(F(g \circ f)(\quot{A_i}{Z})\big) \ar[dd]_{\psi} \ar[rr]^-{\theta_{g \circ f}^{-1}} & & F(g \circ f)\big(\lim (\quot{A_i}{Z}) \big) \ar[d]^{\big(\phi_{f,g}^{\lim A_i}\big)^{-1}} \\ 
		& & F(f)(F(g)\lim (\quot{A_i}{Z}) \ar[d]^{F(f)\tau_{g}^{A}})\\
		\lim (\quot{A_i}{X}) & & F(f)\big(\lim (\quot{A_i}{Y})\big) \ar[ll]^-{\tau_f^A}
	}
	\]
	Substituting our calculation of $\theta^{-1}_{g \circ f}$ above into the definition of $\psi$ yields
	\begin{align*}
		\psi&=\tau_f^A \circ F(f)\tau_g^A \circ \big(\phi_{f,g}^{\lim A_i}\big)^{-1} \circ \theta_{g \circ f}^{-1} \\
		&= \tau_f^A \circ F(f)\tau_g^A \circ F(f)\theta_g^{-1} \circ F(f)(\lim \tau_g^{A_i})^{-1} \circ \theta_f^{-1} \\
		&\circ \lim(F(f)\tau_g^{A_i}) \circ \big(\lim \phi_{f,g}^{A_i}\big)^{-1} \\
		&= \tau_f^A \circ F(f)\tau_g^A \circ F(f)\big(\theta_g^{-1} \circ (\lim \tau_g^{A_i})^{-1}\big) \circ \theta_f^{-1}\\
		& \circ \lim(F(f)\tau_g^{A_i}) \circ \big(\lim \phi_{f,g}^{A_i}\big)^{-1} \\
		&= \tau_f^A \circ F(f)\tau_g^A \circ F(f)\big(\tau_g^{A}\big)^{-1} \circ \theta_{f}^{-1} \circ \lim(F(f) \tau_g^{A_i}) \circ \big(\lim \phi_{f,g}^{A_i}\big)^{-1} \\
		&= \tau_f^A \circ \theta_f^{-1} \circ \lim(F(f)\tau_g^{A_i}) \circ \big(\lim \phi_{f,g}^{A_i}\big)^{-1} \\
		&= \lim(\tau_f^{A_i}) \circ \theta_f \circ \theta_f^{-1} \circ \lim(F(f)\tau_g^{A_i}) \circ \big(\lim \phi_{f,g}^{A_i}\big)^{-1} \\
		&= \lim(\tau_f^{A_i}) \circ \lim(F(f)\tau_g^{A_i}) \circ \big(\lim \phi_{f,g}^{A_i}\big)^{-1} \\
		&= \lim\left(\tau_f^{A_i} \circ F(f)\tau_g^{A_i} \circ (\phi_{f,g}^{A_i})^{-1} \right) \\
		&=\lim\left(\tau_{g \circ f}^{A_i} \circ \phi_{f,g}^{A_i} \circ (\phi_{f,g}^{A_i})^{-1} \right) \\
		&= \lim\tau_{g \circ f}^{A_i},
	\end{align*}
	where we canceled the composite $\big(\phi_{f,g}^{\lim A_i}\big)^{-1} \circ \phi_{f,g}^{\lim A_i}$ without mention in the second equality and where we used the universal property of the limit to give the equation
	\[
	\lim(\tau_f^{A_i}) \circ \lim(F(f)\tau_g^{A_i}) \circ \big(\lim \phi_{f,g}^{A_i}\big)^{-1} = \lim\big(\tau_f^{A_i} \circ F(f)\tau_g^{A_i} \circ \phi_{f,g}^{A_i}\big).
	\]
	This gives that $\tau_f^A \circ F(f)\tau_g^A \circ \big(\phi_{f,g}^{\lim A_i}\big)^{-1} \circ \theta_{g \circ f}^{-1} = \lim\tau_{g \circ f}^{A_i}$, which in turn implies that
	\[
	\tau_f^A \circ F(f)\tau_g^A = \lim\left(\tau_{g \circ f}^{A_i}\right) \circ \theta_{g \circ f} \circ \phi_{f,g}^{\lim A_i} = \tau_{g \circ f}^{A} \circ \phi_{f,g}^{A}.
	\]
	Therefore it follows that $T_{A}$ satisfies the cocycle condition and hence the pair $(A,T_A)$ is an object in $\PC(F)$.
	
	We now show that $(A,T_A)$ is the object vertex of a cone over the $(A_i, T_{A_i})$, i.e., that there are morphisms $P_i:A \to A_i$ for all $i \in I_0$ that make the diagrams
	\[
	\xymatrix{
		& A \ar[dr]^{P_j} \ar[dl]_{P_i} & \\
		A_i \ar[rr]_{d(\alpha)} & & A_j
	}
	\]
	commute whenever there is a morphism $\alpha \in I(i,j)$. For any $i \in I_0$ we define the collection $P:A \to A_i$ by
	\[
	P := \lbrace \quot{\rho_i}{X}:\quot{A}{X} \to \quot{A_i}{X} \; | \; X \in \Cscr_0 \rbrace
	\]
	where each $\quot{\rho}{X} \in F(X)(\quot{A}{X},\quot{A_i}{X})$ is the corresponding cone map at $i$. That this makes $A$ into a cone over the $A_i$ is clear if $P$ is a morphism in $\PC(F)$, as for each $X \in \Cscr_0$ the corresponding cone identity holds. Thus all we need to do is verify that $P$ is a morphism in $\PC(F)$, i.e., we need to show that for all $f \in \Cscr_1$, say with $f:X \to Y$, the diagram
	\[
	\xymatrix{
		F(f)(\quot{A}{Y}) \ar[rr]^-{F(f)\quot{\rho_i}{Y}} \ar[d]_{\tau_f^A}& & F(f)(\quot{A_i}{Y}) \ar[d]^{\tau_f^{A_i}} \\
		\quot{A}{X} \ar[rr]_-{\quot{\rho_i}{X}} & & \quot{A_i}{X}
	}
	\]
	commutes. For this first recall that $\tau_f^A = \lim(\tau_f^{A_i}) \circ \theta_f$ and that the morphism $\lim(\tau_f^{A_i})$ is defined from the commuting diagram
	\[
	\xymatrix{
		\lim\big(F(f)(\quot{A_i}{Y})\big) \ar[rr]^-{\lim\tau_f^{A_i}} \ar[d]_{\pi_i} & & \lim(\quot{A_i}{X}) \ar[d]^{\quot{\rho_i}{X}} \\
		F(f)(\quot{A_i}{Y}) \ar[rr]_-{\tau_f^{A_i}} & & \quot{A_i}{X}
	}
	\]
	in $F(\quot{X}{X})$, where $\pi_i:\lim\big(F(f)\quot{A_i}{Y}\big) \to F(f)\quot{A_i}{Y}$ is the limit map. Furthermore, from the isomorphism $\theta_f$ we also derive the commuting triangle
	\[
	\xymatrix{
		F(f)\big(\lim (\quot{A_i}{Y})\big) \ar[drr]_{F(f)\quot{\rho_i}{Y}} \ar[rr]^-{\theta_f} & & \lim\big(F(f)(\quot{A_i}{Y})\big) \ar[d]^{\pi_i} \\
		& & F(f)(\quot{A_i}{Y})
	}
	\]
	in $F(X)$. These together give the commuting diagram
	\[
	\xymatrix{
		F(f)\big(\lim (\quot{A_i}{Y})\big) \ar[drr]_{F(f)\quot{\rho_i}{Y}} \ar[rr]^-{\theta_f} & &  \lim\big(F(f)(\quot{A_i}{Y})\big) \ar[rr]^-{\lim\tau_f^{A_i}} \ar[d]_{\pi_i} & & \lim(\quot{A_i}{X}) \ar[d]^{\quot{\rho_i}{X}} \\
		& & F(f)(\quot{A_i}{Y}) \ar[rr]_-{\tau_f^{A_i}} & & \quot{A_i}{X}
	}
	\]
	which gives the equality
	\[
	\tau_f^{A_i} \circ F(f)\quo{\rho_i}{Y} = \tau_f^{A_i} \circ \pi_i \circ \theta_f = \quot{\rho_i}{X} \circ \lim\big(\tau_f^{A_i}\big) \circ \theta_f = \quot{\rho_i}{X} \circ \tau_f^A.
	\]
	Thus the diagram
	\[
	\xymatrix{
		F(f)(\quot{A}{Y}) \ar[rr]^-{F(f)\quot{\rho_i}{Y}} \ar[d]_{\tau_f^A}& & F(f)(\quot{A_i}{Y}) \ar[d]^{\tau_f^{A_i}} \\
		\quot{A}{X} \ar[rr]_-{\quot{\rho_i}{X}} & & \quot{A_i}{X}
	}
	\]
	commutes and so $P_i$ is a morphism in $\PC(F)$. This shows that $A$ is a cone over the $A_i$.
	
	Finally we verify that $A$ is a universal cone over $\PC(F)$. For this let $(C,T_C)$ be a cone over the $A_i$ with morphisms $\Psi_i:C \to A_i$ for which 
	\[
	\Psi_i = \lbrace \quot{\psi_i}{X}:\quot{C}{X} \to \quot{A_i}{X} \; | \; X \in \Cscr_0\rbrace.
	\] 
	Because each $\quot{A}{X}$ is the limit in $F(\quot{X}{X})$ of the $\quot{A_i}{X}$, for all $X \in \Cscr_0$ we produce unique morphisms $\quot{\zeta}{X}:\quot{C}{X} \to \quot{A}{X}$ making the diagram
	\[
	\xymatrix{
		& \quot{C}{X} \ar@{-->}[d]^{\quot{\zeta}{X}} \ar@/^/[ddr]^{\quot{\psi_j}{X}} \ar@/_/[ddl]_{\quot{\psi_i}{X}} & \\
		& \quot{A}{X} \ar[dr]_{\quot{\rho_j}{X}} \ar[dl]^{\quot{\rho_i}{X}} \\
		\quot{A_i}{X} \ar[rr]_-{d(\alpha)} & & \quot{A_j}{X}
	}
	\]
	commute for all $i,j \in I_0$ and for all $\alpha \in I(i,j)$. We now will show that the collection
	\[
	\zeta := \lbrace \quot{\zeta}{X} \; | \; X \in \Cscr_0 \rbrace
	\]
	is a morphism in $\PC(F)(C,A)$. Note that this will complete the proof of the proposition, as any other cone map $\Phi:C \to A$ will be component-wise equal to $\zeta$ by the uniqueness of the $\quot{\zeta}{X}$.
	
	To see that $\zeta$ is a morphism we must verify that for any $f \in \Cscr_1$ the diagram
	\[
	\xymatrix{
		F(f)(\quot{C}{Y}) \ar[rr]^-{F(f)\quot{\zeta}{Y}} \ar[d]_{\tau_f^C} & & F(f)(\quot{A}{Y}) \ar[d]^{\tau_f^A} \\
		\quot{C}{X} \ar[rr]_-{\quot{\zeta}{X}} & & \quot{A}{X}
	}
	\]
	commutes where $\Dom f = X$ and $\Codom f = Y$. To this end, let $i,j \in I_0$ and let $\alpha \in I(i,j)$. Then, because $C$ is a cone to the $A_i$ and because the assignment $d:I \to \PC(F)$ is a functor (so in particular $d(\alpha)$ is a morphism in $\PC(F)$), we find that the diagram
	\[
	\xymatrix{
		& F(f)(\quot{C}{Y}) \ar@/^/[dr]^{F(f)\quot{\psi_j}{Y}} \ar@/_/[dl]_{F(f)\quot{\psi_i}{Y}} & \\
		F(f)(\quot{A_i}{Y}) \ar[rr]_-{F(f)\big(d_{Y}(\alpha)\big)} \ar[d]_{\tau_f^{A_i}} & & F(f)(\quot{A_j}{Y}) \ar[d]^{\tau_f^{A_j}} \\
		\quot{A_i}{X} \ar[rr]_-{d_{X}(\alpha)} & & \quot{A_j}{X}
	}
	\]
	commutes in $F(\quot{X}{X})$. Thus there exists a unique morphism $\xi_f$ making the diagram
	\[
	\xymatrix{
		& F(f)(\quot{C}{Y}) \ar@{-->}[d]^{\xi_f} \ar@/^/[ddr]^{\tau_f^{A_j} \circ F(f)\quot{\psi_j}{Y}} \ar@/_/[ddl]_{\tau_f^{A_i} \circ F(f)\quot{\psi_i}{X}^{\prime}} & \\
		& \quot{A}{X} \ar[dr]_{\quot{\rho_j}{X}} \ar[dl]^{\quot{\rho_i}{X}} \\
		\quot{A_i}{X} \ar[rr]_-{d(\alpha)} & & \quot{A_j}{X}
	}
	\]
	commute. Using the above diagram together with the identity 
	\[
	\tau_f^{A_i} \circ F(f)\quot{\rho_i}{Y} = \quot{\rho_i}{X} \circ \tau_f^A
	\] 
	gives that on one hand
	\begin{align*}
		\quot{\rho_i}{X} \circ \xi_f &= \tau_f^{A_i} \circ F(f)\quot{\psi_i}{Y} = \tau_f^{A_i} \circ F(f)\big(\quot{\rho_i}{Y} \circ \quot{\zeta}{Y}\big) \\
		&= \tau_f^{A_i} \circ F(f)\quot{\rho_i}{Y} \circ F(f)\quot{\zeta}{Y}  \\
		&=\quot{\rho_i}{X} \circ \tau_f^A \circ F(f)\quot{\zeta}{Y};
	\end{align*}
	following this same argument mutatis mutandis we derive also that 
	\[
	\quot{\rho_j}{X} \circ \xi_f = \quot{\rho_j}{X} \circ \tau_f^A \circ F(f)\quot{\zeta}{Y}.
	\] 
	On the other hand, using that $\Psi_i$ is an $\PC(F)$-morphism gives that
	\[
	\quot{\rho_i}{X} \circ \quot{\zeta}{X} \circ \tau_f^C = \quot{\psi_i}{X} \circ \tau_f^C = \tau_f^{A_i} \circ F(f)\quot{\psi_i}{Y} = \quot{\rho_i}{X} \circ \xi_f.
	\]
	As before, we also derive that $\quot{\rho_j}{X} \circ \quot{\zeta}{X} \circ \tau_f^C = \quot{\rho_j}{X} \circ \xi_f$. Putting all these identities together gives us that both $\quot{\zeta}{X} \circ \tau_f^C$ and $\tau_f^A \circ F(f)\quot{\zeta}{Y}$ produce fillers for the universal cone diagram in the same way that $\xi_f$ does; consequently, by the universal property of $\quot{A}{X}$ as the limit over the $\quot{A_i}{X}$, we derive that
	\[
	\tau_f^{A} \circ F(f)\quot{\zeta}{Y} = \xi_f = \quot{\zeta}{X} \circ \tau_f^C.
	\]
	This gives the commutativity of the diagram
	\[
	\xymatrix{
		F(f)(\quot{C}{Y}) \ar[rr]^-{F(f)\quot{\zeta}{Y}} \ar[d]_{\tau_f^C} & & F(f)(\quot{A}{Y}) \ar[d]^{\tau_f^A} \\
		\quot{C}{X} \ar[rr]_-{\quot{\zeta}{X}} & & \quot{A}{X}
	}
	\]
	and hence shows that $(A,T_A)$ is a universal cone over the $(A_i, T_{A_i})$. Because this shows that $d$ has a limit in $\PC(F)$, this completes the proof of the theorem.
\end{proof}
\begin{Theorem}\label{Theorem: Section 2: Equivariant Cat has colims}
	Let $F:\Cscr^{\op} \to \fCat$ be a pseudofunctor and let $I$ be an index category for which there is a diagram $d:I \to \PC(F)$. For each $X \in \Cscr_0$ let $d_{X}:I \to F(X)$ denote the diagram functor
	\[
	\xymatrix{
		I \ar[drr]_{d_{X}} \ar[r]^-{d} & \PC(F) \ar[r]^-{\tilde{\imath}} & \prod\limits_{X \in \Cscr_0} F(X) \ar[d]^{\pi_{X}} \\
		& & F(X)
	}
	\]
	and assume that for each $X \in \Cscr_0$ the colimit $\colim d_{X}$ exists in $F(X)$. Furthermore, assume that for every morphism $f \in \Cscr_1$ with $f:X \to Y$ there is an isomorphism $\theta_f$ witnessing the limit preservation
	\[
	F(f)\left(\lim_{\substack{\longrightarrow \\ i \in I}}d_{X}(i)\right) \cong \lim_{\substack{\longrightarrow \\ i \in I}} F(f)(d_{X}(i)).
	\]
	Then the diagram $d$ admits a colimit in $\PC(F)$.
\end{Theorem}
\begin{proof}
	Dualize the proof of Theorem \ref{Thm: Section 2: Equivariant Cat has lims}.
\end{proof}

\begin{corollary}\label{Cor: Section 2: Limits of shape I in equiv cat}
	Let $F:\Cscr^{\op} \to \fCat$ be a pseudofunctor. If $I$ is an index category and each category $F(X)$ has all limits of shape $I$ and if each functor $F(f)$ preserves these limits, then $\PC(F)$ has all limits of shape $I$ as well. Dually, if each category $F(X)$ has all colimits of shape $I$ and if each functor $F(f)$ preserves these colimits, then $\PC(F)$ has all colimits of shape $I$ as well.
\end{corollary}
\begin{proof}
	Use Theorem \ref{Thm: Section 2: Equivariant Cat has lims} for every diagram $d:I \to \Cscr$, as such a diagram gives rise to a compatable family of diagrams $d_{X}:I \to \Cscr$ by way of projection to the $X$-th component.
\end{proof}
\begin{corollary}\label{Cor: Equivariant cat is finitely complete or cocomplete}
	If each category $F(X)$ is (finitely) complete and each functor $F(f)$ is (finitely) continuous then $\PC(F)$ is (finitely) complete. Dually, $\PC(F)$ is (finitely) cocomplete if each category $F(X)$ is (finitely) cocomplete and if each functor $F(f)$ is finitely cocontinuous. 
\end{corollary}
\begin{proof}
	Use Corollary \ref{Cor: Section 2: Limits of shape I in equiv cat} with the categories $I = \emptyset$ to give terminal objects, $I = \coprod_{s \in S} \Delta^0$ (where $S$ is an indexing $\Vscr$-set) to give products indexed by $S$, and $I = \Lambda_2^2$ to give pullbacks. Invoking now the classical categorical result that a category $\Cscr$ admits all limits if and only if it admits products, terminal objects, and pullbacks we find that $\PC(F)$ is complete.
	%
	%
\end{proof}
\begin{remark}
	The finiteness assumptions above are used exclusively when determining which indexing sets $S$ we are allowed to use. If one only knows of the existence of finite limits or colimits, we restrict our attention to sets of finite cardinality.
\end{remark}

\begin{corollary}\label{Cor: Section 2: Additive pseudofunctor gives additive equivariat cat}
	If the pseudofunctor $F:\Cscr^{\op} \to \fCat$ is additive and locally $\Vscr$-small then the category $\PC(F)$ is a locally $\Vscr$-small additive category. 
\end{corollary}
\begin{proof}
	Corollary \ref{Cor: Section 2: Limits of shape I in equiv cat} shows that $\PC(F)$ admits finite products and a zero object. Local $\Vscr$-smallness is Lemma \ref{Lemma: Section 2: Pseudofunctor values in Vlocally small means Equiv cat Vlocally small}, while the $\VAb$ enrichment is Corollary \ref{Cor: Pseudocone Section: Ab Enrichment} as $\VAb = \Ab(\VSet)$. Finally, appealing to a classical result in homological algebra (cf.\@ \cite[Proposition II.9.1]{HiltonStammHA}, for instance) gives that $\PC(F)$ is additive.
\end{proof}
\begin{corollary}\label{Cor: Section 2: Additive complete and cocomplete equivariant cats}
	If each of the categories $F(X)$ are complete and cocomplete as well as additive, and if the functors $F(f)$ preserve limits and colimits, then $\PC(F)$ is additive, complete, and cocomplete as well.
\end{corollary}
\begin{proof}
	Corollary \ref{Cor: Section 2: Additive pseudofunctor gives additive equivariat cat} gives the additive structure on $\PC(F)$ while Corollary \ref{Cor: Equivariant cat is finitely complete or cocomplete} gives the completion and cocompletion.
\end{proof}

We now present some lemmas that allow us to deduce some cases when the pseudocone categories $\PC(F)$ are Abelian. Of particular use here is that this gives proofs that the categories of equivariant ${D}$-modules, local systems, and quasi-coherent ({\'e}tale) sheaves on a variety $X$ are Abelian (cf.\@ Corollary \ref{Cor: Section 2: Additive equivriant derived, local system, and qcoh cats} below). Similar results for when $X$ is a topological space or a smoooth manifold hold as well.

While illuminating the structure of Abelian pseudocone categories, we also present a proposition about descending regular categories through these pseudocone constructions. While this will not be fundamentally important for us, it will show the way in which descent arguments aid in developing a notion of pseudocone-based categorical logic.

We proceed by recalling that a category $\Cscr$ is Abelian if and only if it is additive, (finitely) complete, (finitely) cocomplete, and has the property that every monomorphism is a kernel and every epimorphism is a cokernel. We have already seen in Corollary \ref{Cor: Section 2: Additive complete and cocomplete equivariant cats} that under mild assumptions the category $\PC(F)$ is additive, complete, and cocomplete. Thus all we need to show are the regularity conditions, i.e., the fact that every monomorphism is a kernel and that every epimorphism is a cokernel.
\begin{lemma}\label{Lemma: Section 2: KErnels descend}
	Assume that $F$ is a pre-equivariant pseudofunctor, that each functor $F(f)$ is continuous and cocontinuous, and that each category $F(X)$ admits equalizers and cokernel pairs of morphisms. Then if every monic in each fibre category $F(X)$ is a kernel, every monic in $\PC(F)$ is a kernel.
\end{lemma}
\begin{proof}
	Begin by recalling that a morphism $\Psi$ is a kernel of a morphism $\Phi:A \to B$ if and only if there is an isomorphism $\Psi \cong \eq(I_0, I_1)$, where $I_0$ and $I_1$ are the cokernel pair of $\Phi$ that appear in the pushout below:
	\[
	\begin{tikzcd}
		\Dom \Phi \ar[rr]{}{\Phi} \ar[d, swap]{}{\Phi} & & \Codom \Phi \ar[d]{}{I_1} \\
		\Codom \Phi \ar[rr, swap]{}{I_0} & & \Codom \Phi \coprod_{\Dom \Phi} \Codom \Phi
	\end{tikzcd}
	\]
	
	To proceed we let $P = \lbrace \quot{\rho}{X} \; | \; X \in \Cscr_0\rbrace:A \to B$ be a monomorphism in $\PC(F)$ and we let the pair $(I_0,I_1)$ be the cokernel pair of $P$ arising from the pushout
	\[
	\begin{tikzcd}
		A \ar[r]{}{P} \ar[d, swap]{}{P} & B \ar[d]{}{I_1} \\
		B \ar[r, swap]{}{I_0} & B \coprod_A B
	\end{tikzcd}
	\]
	in $\PC(F)$. From Theorem \ref{Theorem: Section 2: Equivariant Cat has colims} we know that the pushout of $B \leftarrow A \rightarrow B$ is given by
	\[
	B \coprod_A B = \left\lbrace \quot{B}{X} \coprod_{\quot{A}{X}} \quot{B}{X} \; : \; X \in \Cscr_0 \right\rbrace
	\]
	so each $I_k$ has the form
	\[
	I_k = \left\lbrace \quot{i_k}{X}:\quot{B}{X} \to \quot{B}{X} \coprod_{\quot{A}{X}} \quot{B}{X} \; : \; X \in \Cscr_0 \right\rbrace,
	\]
	where $k \in \lbrace 0, 1 \rbrace$ and where each pair $(\quot{i_0}{X},\quot{i_1}{X})$ is the cokernel pair of $\quot{\rho}{X}$ in $F(X)$. Using Theorem \ref{Thm: Section 2: Equivariant Cat has lims} we also find that the morphism $P$, as a monic in $\PC(F)$, has the property that each $\quot{\rho}{X}$ is monic in $F(X)$. Because each monomorphism in each fibre category $F(X)$ is a kernel, there is a unique isomorphism $\quot{\theta}{X}:X \xrightarrow{\cong} \Eq(\quot{i_0}{X},\quot{i_1}{X})$ making the diagram
	\[
	\begin{tikzcd}
		\quot{A}{X} \ar[d, dashed, swap]{}{\exists!\quot{\theta}{X}} \ar[rr]{}{\quot{\rho}{X}} & & \quot{B}{X} \ar[d, equals] \\
		\Eq(\quot{i_0}{X},\quot{i_1}{X}) \ar[rr, swap]{}{\eq(\quot{i_0}{X}, \quot{i_1}{X})} & & \quot{B}{X}
	\end{tikzcd}
	\]
	commute in $F(X)$; note that $\Eq(\quot{i_0}{X},\quot{i_1}{X})$ is the equalizer object of the $\quot{i_k}{X}$ and $\eq(\quot{i_0}{X},\quot{i_1}{X})$ is the canonical equalizer morphism. Using Theorem \ref{Thm: Section 2: Equivariant Cat has lims} again, we have that the collection
	\[
	\Eq(I_0,I_1) := \left\lbrace \Eq(\quot{i_0}{X},\quot{i_1}{X}) \; | \; X \in \Cscr_0 \right\rbrace
	\]
	determines the object assignment of the equalizer $\Eq(I_0,I_1)$ in $\PC(F)$ and the collection
	\[
	\eq(I_0,I_1) := \left\lbrace \eq(\quot{i_0}{X},\quot{i_1}{X}) \;|\; X \in \Cscr_0\right\rbrace
	\]
	determines the equalizer morphism $\eq(I_0,I_1):\Eq(I_0,I_1) \to B$ in $\PC(F)$. Thus we will be done if we can show that the $\Vscr$-set
	\[
	\Theta := \left\lbrace \quot{\theta}{X} \; | \; X \in \Cscr_0 \right\rbrace
	\]
	is a morphism $\Theta:A \to \Eq(I_0, I_1)$ in $\PC(F)$, as then $\Theta$ is the unique morphism making the square
	\[
	\begin{tikzcd}
		A \ar[r]{}{P} \ar[d, dashed, swap]{}{\exists!\Theta} & B \ar[d, equals] \\
		\Eq(I_0,I_1) \ar[r, swap]{}{\eq(I_0,I_1)} & B
	\end{tikzcd}
	\]
	commute.
	
	For this fix a morphism $f \in \Cscr_1$ and let $f:X \to Y$. We now must show that $\quot{\theta}{X} \circ \tau_f^A = \tau_f^{\Eq(I_0,I_1)} \circ F(f)\eq(\quot{i_0}{Y},\quot{i_1}{Y})$. For this we calculate that
	\begin{align*}
		\eq(\quot{i_0}{X},\quot{i_1}{X}) \circ \quot{\theta}{X} \circ \tau_f^A &= \quot{\rho}{X} \circ \tau_f^A = \tau_f^B \circ F(f)\quot{\rho}{Y} \\
		&= \tau_f^B \circ F(f)\big(\eq(\quot{i_0}{Y},\quot{i_1}{Y}) \circ \quot{\theta}{Y}\big) \\
		&= \tau_f^B \circ F(f)\eq(\quot{i_0}{Y},\quot{i_1}{Y}) \circ F(f)\quot{\theta}{X{^\prime}} \\
		&= \eq(\quot{i_0}{X},\quot{i_1}{X}) \circ \tau_f^{\Eq(I_0,I_1)} \circ F(f)\quot{\theta}{Y}.
	\end{align*}
	Using that $\eq(\quot{i_0}{X},\quot{i_1}{X})$ is monic then allows us to deduce that
	\[
	\quot{\theta}{X} \circ \tau_f^A = \tau_f^{\Eq(I_0,I_1)} \circ F(f)\quot{\theta}{X},
	\]
	which proves that $\Theta \in \PC(F)(A,\Eq(I_0,I_1))$.
\end{proof}
By dualizing the above lemma, we also have the result below for cokernels and epimorphisms.
\begin{lemma}\label{Lemma: Section 2: Cokernels descend}
	Assume that $F$ is a pre-equivariant pseudofunctor, that each functor $F(f)$ is continuous and cocontinuous, and that each category $F(X)$ admits coequalizers and kernels pairs of morphisms. Then if every epimorphism in each fibre category $F(X)$ is a cokernel, every epimorphism in $\PC(F)$ is a cokernel.
\end{lemma}
This leads us to the fact that Abelian categories descend through the equivariance, so if each fibre category $F(X)$ is Abelian and if each fibre functor $F(f)$ is exact then $\PC(F)$ is Abelian.
\begin{proposition}\label{Prop: Section 2: Abelian Descends}
	Let $F:\Cscr^{\op} \to \fCat$ be a pseudofunctor such that each category $F(X)$ is Abelian and for which each functor $F(f)$ is exact. Then the category $\PC(F)$ is Abelian as well.
\end{proposition}
\begin{proof}
	By Corollary \ref{Cor: Section 2: Additive complete and cocomplete equivariant cats} it follows that the category $\PC(F)$ is complete, cocomplete, and additive. To see that it has the necessary regularity property, simply use Lemmas \ref{Lemma: Section 2: KErnels descend} and \ref{Lemma: Section 2: Cokernels descend} to derive that every monic and epic in $\PC(F)$ is a kernel or cokernel, respectively.
\end{proof}

As a corollary, we have that the categories of equivariant $D$-modules, local systems, and quasicoherent sheaves are all Abelian.
\begin{corollary}\label{Cor: Section 2: Additive equivriant derived, local system, and qcoh cats}
	Each of the categories $D_G^b(X)$, $\DGMod(X)$, $\Loc_G(X)$, $\Per_G(X)$, and $\QCoh_G(X)$ are all additive. Moreover, the categories $\DGMod(X)$, $\Loc_G(X)$, $\Per_G(X)$, and $\QCoh_G(X)$ are Abelian.
\end{corollary}
\begin{proof}
	That all four categories are additive follows from Corollary \ref{Cor: Section 2: Additive pseudofunctor gives additive equivariat cat}. For each $\Gamma \in \Sf(G)_0$ the quotient map 
	\[
	\quo_{\Gamma}:\Gamma \times X \to G\backslash (\Gamma \times X)
	\] 
	is smooth by Proposition \ref{Prop: Section 2.1: Smooth quotient map}. Similarly, for any $f:\Gamma \to \Gamma^{\prime}$ in $\Sf(G)$ the induced map between quotient varieties
	\[
	\of:G \backslash (\Gamma \times X) \to G \backslash (\Gamma^{\prime} \times X)
	\]
	is also smooth by Proposition \ref{Prop: Section 2.1: Smooth quotient map}. This then implies that  each functor $\overline{f}^{\ast}$ is exact on the level of {\'e}tale sheaves and both exists and is exact at the level of $D$-modules. Thus we can apply Corollary \ref{Cor: Section 2: Additive complete and cocomplete equivariant cats} to derive that the categories $\DGMod(X),$\\ $\Loc_G(X)$, $\Per_G(X)$, and $\QCoh_G(X)$ are all finitely complete and finitely cocomplete. Finally, Proposition \ref{Prop: Section 2: Abelian Descends} shows that the categories $\DGMod(X),$ $\Loc_G(X),$ $\Per_G(X)$ and $\QCoh_G(X)$ are all Abelian categories, as each of the fibre categories $\DMod(G \backslash (\Gamma \times X)),$ $\Loc(G \backslash (\Gamma \times X))$, $\Per(G \backslash (\Gamma \times X))$, and $\QCoh(G \backslash (\Gamma \times X))$ are Abelian for all $\Gamma \times X \in \SfResl(G)_0$ and have smooth fibre functors $\of^{\ast}$ for all $f \times \id_X \in \Sf(G)_1$.
\end{proof}
The same proof may be used to show that the category $\Shv_G(X,R)$ of ({\'e}tale) sheaves of (left) modules over a fixed ring $R$ (or modules over a fixed sheaf of rings $\CalO_X$) is Abelian; similarly, the category $\Shv_G(X;\overline{\Q}_{\ell})$ of $\ell$-adic sheaves is Abelian. We record these results in the corollary below.
\begin{corollary}
	The category $\OXMod_G$\index[notation]{OXMod@$\OXMod_G$} of equivariant $\Ocal_X$-modules and the category $\Shv_G(X;R)$\index[notation]{SheavesofRmodules@$\Shv_G(X,R)$}  of equivariant sheaves of $R$-modules are Abelian. Moreover, the categories $\Shv_G(X;\overline{\Q}_{\ell})$ of equivariant $\ell$-adic sheaves and $\Per_G(X;\overline{\Q}_{\ell})$ of equivariant $\ell$-adic perverse sheaves are also Abelian.
\end{corollary}
\newpage

\section{Monoidal Structures on Pseudocone Categories}\label{Subsection: Monoidal Structures on PCF}
We now begin providing cases when the equivariant categories $\PC(F)$ are monoidal. This will allow us to not only give a foundation for an equivariant monoidal theory, but also give us tools to study equivariant derived monoidal theory and equivariant dg-categories over the geometric objects about which we care. The main use for these tchniques in this paper are in homological algebraic settings, as it gives us tools with which to construct the six-functor-formalism of the equivariant derived category and other pseudocone-based descent-theoretic settings. I expect us to be able to use this in noncommutative algebraic geometry (following the foundations and methodology of \cite{Keller}, \cite{Orlov_2018}, \cite{toën2014derived}, \cite{ToenDAGDQ}, among others, which define and work with noncommutative schemes in dg-categorical language) to give an equivariant notion of noncommutative (derived) algebraic geometry. Additionally, the braided monoidal category structure on $\PC(F)$ could be useful in equivariant higher algebra or in a descent-based approach to TQFTs, fusion categories, tensor categories and other such structures useful in higher-categorical physics; for various references on tensor categories, TQFTs, fusion categories, and such see for example \cite{Arun}, \cite{Delaney2019fusion}, \cite{DeligneTanak}, \cite{Chris2Noah}, and \cite{TheoClaudia}. For the moment, however, we will simply provide the main situation that allows us to deduce when the category $\PC(F)$ is monoidal in terms of having a pseudofunctor with monoidal categories for fibre categories and (strong) monoidal functors for fibre functors.
\begin{definition}\label{Defn: Section 2: Monoidal Preequivariant Pseudofunctor}
	Let $F:\Cscr^{\op} \to \fCat$ be a pseudofunctor. We say that $F$ is {monoidal}\index[terminology]{Pseudofunctor! Monoidal} if the following hold:
	\begin{itemize}
		\item For all $X \in \Cscr_0$, each category $F(X)$ is a monoidal category \[
		\left(F(X), \os{X}{\otimes}, \quot{I}{X}, \quot{\alpha}{X}, \quot{\lambda}{X}, \quot{\rho}{X}\right);
		\]
		\item For all $f:X \to Y$ in $\Cscr_1$, $f$ is monoidal in the sense that there are natural isomorphisms
		\[
		\theta_f^{A,B}:F(f)\left(A \os{Y}{\otimes} B\right) \xrightarrow{\cong} F(f)A \os{X}{\otimes} F(f)B
		\]
		which are pseudofunctorial in the sense that for all composable pairs $X \xrightarrow{f} Y \xrightarrow{g} Z$ in $\Cscr$, the pasting diagram
		\[
		\begin{tikzcd}
			F(\quot{X}{Z}) \times F(\quot{X}{Z}) \ar[d]{}{F(g) \times F(g)} \ar[dd, bend right = 80, swap, ""{name = LL}]{}{F(g \circ f) \times F(g\circ f)} \ar[rrr, ""{name = UU}]{}{\os{Z}{\otimes}} & & & F(\quot{X}{Z}) \ar[dd, bend left = 80, ""{name = RR}]{}{F(g \circ f)} \ar[d,swap]{}{F(g)}\\
			F(\quot{X}{Y}) \times F(\quot{X}{Y}) \ar[d]{}{F(f) \times F(f)} \ar[rrr, ""{name = M}]{}[description]{\os{Y}{\otimes}} & & & F(\quot{X}{Y}) \ar[d]{}{F(f)} \\
			F(X) \times F(X) \ar[rrr, swap, ""{name = BB}]{}{\os{X}{\otimes}} & & & F(X)
			\ar[from = 2-4, to = RR, Rightarrow, shorten <=2pt, shorten >= 2pt, swap]{}{\phi_{f,g}}
			\ar[from = LL, to = 2-1, swap, Rightarrow, shorten <= 2pt, shorten >= 2pt]{}{\phi_{f,g}^{-1} \times \phi_{f,g}^{-1}}
			\ar[from = UU, to = M, Rightarrow, shorten <= 4pt, shorten >= 4pt]{}{\theta_g}
			\ar[from = M, to = BB, Rightarrow, shorten <= 8pt, shorten >= 4pt]{}{\theta_f}
		\end{tikzcd}
		\]
		is equal to the $2$-cell:
		\[
		\begin{tikzcd}
			F(\quot{X}{Z}) \times F(\quot{X}{Z}) \ar[r, ""{name = U}]{}{\os{Z}{\otimes}} \ar[d, swap]{}{F(g\circ f) \times F(g \circ f)} & F(\quot{X}{Z}) \ar[d]{}{F(g \circ f)} \\
			F(X) \times F(X) \ar[r, swap, ""{name = L}]{}{\os{X}{\otimes}} & F(X) \ar[from = U, to = L, Rightarrow, shorten >= 4pt, shorten <= 4pt]{}{\theta_{g \circ f}}
		\end{tikzcd}
		\]
		\item The functors $F(f)$ preserve the units pseudofunctorially in the sense that there for all $f \in \Cscr_1$, there is an isomorphism
		\[
		\sigma_f:F(f)\quot{I}{Y} \xrightarrow{\cong} \quot{I}{X}
		\]
		such that for all composable pairs $X \xrightarrow{f} Y \xrightarrow{g} Z$ in $\Cscr$, we have that
		\[
		\sigma_f \circ F(f)\sigma_g = \sigma_{g \circ f} \circ \phi_{f,g}.
		\]
	\end{itemize}
\end{definition}
\begin{remark}
	The functors we have called monoidal in this paper are frequently called {strong monoidal}\index[terminology]{Monoidal functor}\index[terminology]{Strong monoidal functor|see {Monoidal Functor}} in the category-theoretic literature (cf.\@ \cite{AguiMahajan.}, \cite{KellyStreet}, \cite{Leinster}, \cite{Yetter}, for instance). We have chosen to drop the adjective ``strong'' for a few reasons: First, we will need the isomorphisms for each use of a monoidal functor in this paper; second, we do not want to confuse strong monoidal functors with monoidal functors that have tensorial strengths; and third, having monoidal functors be implicitly lax is not something that our descent formalism will allow. We need to be working with tensorial transition isomorphisms, not simply transition maps; as such, we will need the full strength of a monoidal functor which preserves the tensor and units up to natural isomorphism.
\end{remark}
\begin{Theorem}\label{Theorem: Section 2: Monoidal preequivariant pseudofunctor gives monoidal equivariant cat}
	Let $F:\Cscr^{\op} \to \fCat$ be a monoidal pseudofunctor. Then $\PC(F)$ is a monoidal category.
\end{Theorem}
\begin{proof}
	First note that if we can define the functor $\otimes,$ the object $I$, and the natural isomorphisms $\lambda, \rho,$ and $\alpha$ which $X$-wise agree with the functors, objects, and transformations in the fibre category
	\[
	\left(F(X), \os{X}{\otimes}, \quot{I}{X}, \quot{\lambda}{X}, \quot{\rho}{X}, \quot{\alpha}{X}\right)
	\]
	we will be done, as all triangle and pentagonal identities will hold $X$-wise, which is how we we compute these identities. We begin by defining the monoidal functor $\otimes:\PC(F) \times \PC(F) \to \PC(F)$. Define the functor on objects by setting, for $(A,T_A)$ and $(B,T_B)$ objects in $\PC(F)$,
	\[
	A \otimes B := \left\lbrace \quot{A}{X} \os{X}{\otimes} \quot{B}{X} \; \big| \; X \in \Cscr_0, \quot{A}{X} \in A, \quot{B}{X} \in B \right\rbrace.
	\]
	To define the transition isomorphisms on $A \otimes B$ first fix an $f \in \Cscr_1$ and write $f:X \to Y$. Then we define the morphism $\tau_f^{A \otimes B}$ via the composite
	\[
	\xymatrix{
		F(f)\left(\quot{A}{Y} \os{Y}{\otimes} \quot{B}{Y} \right) \ar[rr]^-{\theta_f^{A,B}} \ar[drr]_{\tau_f^{A \otimes B}} & & F(f)(\quot{A}{Y}) \os{X}{\otimes} F(f)(\quot{B}{Y}) \ar[d]^{\tau_f^A \otimes \tau_f^B} \\
		& & \quot{A}{X} \os{X}{\otimes} \quot{B}{X}
	}
	\]
	where $\theta_f^{A,B}$ is used as an abuse of notation for
	\[
	\theta_f^{A,B} := \theta_f^{\quot{A}{Y},\quot{B}{Y}}.
	\]
	We then define the transition isomorphisms via the definition
	\[
	T_{A \otimes B} := \left\lbrace \tau_f^{A \otimes B} \; | \; f \in \Cscr_1 \right\rbrace. 
	\]
	
	To verify that the transition isomorphisms satisfy the cocycle condition, let $X \xrightarrow{f} Y \xrightarrow{g} Z$ be a pair of composable morphisms in $\Cscr$. Now, using the pseudofunctoriality of the $\otimes$ functors together with the naturality of the $\theta_f$'s gives that
	\begin{align*}
		\tau_{g \circ f}^{A \otimes B} &= \left(\tau_{g \circ f}^{A} \os{X}{\otimes} \tau_{g \circ f}^{B}\right) \circ \theta_{g \circ f}^{A,B} \\
		&= \left(\big(\tau_f^A \circ F(f)\tau_g^A\big) \os{X}{\otimes}\big(\tau_f^B \circ F(f)\tau_g^B \big)\right) \circ \phi_{f,g}^{-1} \circ \theta_{g \circ f}^{A,B} \\
		&= \left(\big(\tau_f^A \circ F(f)\tau_g^A\big) \os{X}{\otimes}\big(\tau_f^B \circ F(f)\tau_g^B \big)\right) \\
		&\circ \phi_{f,g}^{-1} \circ \phi_{f,g} \circ \theta_{f}^{A,B} \circ F(f)\theta_g^{A,B} \circ \phi_{f,g}^{-1} \\
		&=\left(\big(\tau_f^A \circ F(f)\tau_g^A\big) \os{X}{\otimes}\big(\tau_f^B \circ F(f)\tau_g^B \big)\right) \circ \theta_{f}^{A,B} \circ F(f)\theta_g^{A,B} \circ \phi_{f,g}^{-1} \\
		&= \left(\tau_f^A \os{X}{\otimes} \tau_f^B\right) \circ \left(F(f)\tau_g^A \os{X}{\otimes} F(f)\tau_g^B\right) \circ \theta_{f}^{A,B} \circ F(f)\theta_{g}^{A,B} \circ \phi_{f,g}^{-1} \\
		&= \left(\tau_f^A \os{X}{\otimes} \tau_f^B\right) \circ \theta_f^{A,B} \circ F(f)\left(\tau_g^A \os{Y}{\otimes} \tau_g^B\right) \circ F(f)\theta_{g}^{A,B} \circ \phi_{f,g}^{-1} \\
		&= \tau_f^{A \otimes B} \circ F(f)\tau_{g}^{A \otimes B} \circ \phi_{f,g}^{-1}.
	\end{align*}
	Thus we have that $\tau_{g \circ f}^{A \otimes B} \circ \phi_{f,g} = \tau_f^{A\otimes B} \circ F(f)\tau_g^{A \otimes B}$ and so 
	\[
	\big((A,T_A),(B,T_B)\big) \mapsto (A \otimes B, T_{A \otimes B})
	\] 
	determines the object assignment of the $\otimes$ functor.
	
	To define $\otimes$ on morphisms, fix a morphism $\Phi \in \PC(F)(A,B)$ and a morphism $\Psi \in \PC(F)(A^{\prime},B^{\prime})$. We then define the morphism \[
	\Phi \otimes \Psi:A \otimes A^{\prime} \to B \otimes B^{\prime}
	\] 
	by
	\[
	\Phi \otimes \Psi := \left\lbrace \quot{\varphi}{X} \os{X}{\otimes} \quot{\psi}{X} \; \big| \; X \in \Cscr_0, \quot{\varphi}{X} \in \Phi, \quot{\psi}{X} \in \Psi \right\rbrace.
	\]
	This is a morphism because for any map $f:X \to Y$ in $\Cscr$, we have the equalities $\quot{\varphi}{X} \circ \tau_f^A = \tau_f^{A^{\prime}} \circ F(f)\quot{\varphi}{Y}$ and $\quot{\psi}{X} \circ \tau_f^{B} = \tau_f^{B^{\prime}} \circ F(f)\quot{\psi}{Y}$ so it is immediate that
	\begin{align*}
		\quot{(\Phi \otimes \Psi)}{X} \circ \tau_f^{A \otimes B} &= \left(\quot{\varphi}{X} \os{X}{\otimes} \quot{\psi}{X}\right) \circ \left(\tau_f^A \os{X}{\otimes} \tau_f^B\right) \circ \theta_f^{A,B} \\
		&= \left(\tau_f^{A^{\prime}} \os{X}{\otimes} \tau_f^{B^{\prime}}\right) \circ \left(F(f)\quot{\varphi}{Y} \os{X}{\otimes} F(f)\quot{\psi}{Y}\right) \\
		&= \left(\tau_f^{A^{\prime}} \os{X}{\otimes} \tau_f^{B^{\prime}}\right) \circ \theta_f^{A^{\prime},B^{\prime}} \circ F(f)\left(\quot{\varphi}{Y} \os{Y}{\otimes} \quot{\psi}{Y}\right) \\
		&= \tau_{f}^{A^{\prime} \otimes B^{\prime}} \circ F(f)(\quot{(\Phi \otimes \Psi)}{Y}).
	\end{align*}
	It also follows by the fact that each of the $\os{X}{\otimes}$ are functors that $\otimes$ preserves compositions and identities, giving that $\otimes$ is a functor as well.
	
	Let us show that the object $I = \lbrace \quot{I}{X} \; | \; X \in \Cscr_0 \rbrace$ determines an object in $\PC(F)$. Define $T_I := \lbrace \sigma_f:F(f)\quot{I}{Y} \to \quot{I}{X} \; | \; f \in \Cscr_1 \rbrace$ and note that by assumption for any $X \xrightarrow{f} Y \xrightarrow{g} Y$ in $\Cscr$, we have that
	\[
	\tau_f^{I} \circ F(f)\tau_g^I = \sigma_f \circ F(f)\sigma_g = \sigma_{g \circ f} \circ \phi_{f,g}= \tau_{g \circ f}^{I} \circ \phi_{f,g}.
	\]
	Thus $(I,T_I)$ is an object in $\PC(F)$. 
	
	Next we prove that there are natural isomorphisms $\lambda:I \otimes \id_{\PC(F)} \to \id_{\PC(F)}$ and $\rho:\id_{\PC(F)} \otimes I \to \id_{\PC(F)}$;  we will only show that $\lambda$ exists, as the existence of $\rho$ will follow mutatis mutandis.
	
	Fix an object $A \in \PC(F)_0$ and define $\lambda_A: I \otimes A \to A$ by
	\[
	\lambda_A := \left\lbrace \quot{\lambda_{\quot{A}{X}}}{X}:\quot{I}{X} \os{X}{\otimes} \quot{A}{X} \xrightarrow{\cong} \quot{A}{X} \; | \; X \in \Cscr_0 \right\rbrace.
	\]
	It remains to verify that $\lambda_A$ is a morphism in $\PC(F)$, so let $f: X \to Y$ be a morphism in $\Cscr$. Observe that since the tensor functors vary pseudofunctorially over $\Cscr$ through the preservation isomorphisms, we get that the invertible $2$-cell
	\[
	\begin{tikzcd}
		F(f)\big(\quot{I}{Y} \os{Y}{\otimes} \quot{A}{Y}\big) \ar[r, ""{name = U}]{}{\tau_f^I \os{X}{\otimes} \tau_f^A} \ar[d,swap]{}{F(f)\quot{\lambda}{Y}} & \quot{A}{X} \os{X}{\otimes} \quot{I}{X} \ar[d]{}{\quot{\lambda_{\quot{A}{X}}}{X}} \\
		F(f)(\quot{A}{Y}) \ar[r, swap, ""{name = L}]{}{\tau_f^A} & \quot{A}{X} \ar[from = U, to = L, Rightarrow, shorten >= 4pt, shorten <= 4pt]{}{\theta_f}
	\end{tikzcd}
	\]
	induces the commuting diagram:
	\[
	\xymatrix{
		F(f)\big(\quot{I}{Y} \os{Y}{\otimes} \quot{A}{Y}\big)\ar[d]_{F(f)\quot{\lambda_{\quot{A}{Y}}}{Y}} \ar[r]^-{\theta_{f}^{I,A}} & F(f)(\quot{I}{Y}) \os{X}{\otimes} F(f)(\quot{A}{Y}) \ar[r]^-{\tau_f^I \os{X}{\otimes} \tau_f^A} & \quot{A}{X} \os{X}{\otimes} \quot{I}{X} \ar[d]^{\quot{\lambda_{\quot{A}{X}}}{X}} \\
		F(f)(\quot{A}{Y}) \ar[rr]_-{\tau_f^A} & & \quot{A}{X}
	}
	\]
	This gives that
	\[
	\tau_f^A \circ F(f)\left(\quot{\lambda_{\quot{A}{Y}}}{Y}\right) = \quot{\lambda_{\quot{A}{X}}}{X} \circ \left(\tau_f^I \os{X}{\otimes} \tau_f^A\right) \circ \theta_f^{I,A} = \quot{\lambda_{\quot{A}{X}}}{X} \circ \tau_f^{I \otimes A},
	\]
	which proves that $\lambda_A$ is an $\PC(F)$ morphism (and in fact an isomorphism because each $\quot{\lambda}{X}$ is an isomorphism). That $\lambda$ is natural is immediate from the natruality of each $\quot{\lambda}{X}$, so we indeed get that $\lambda:I \otimes \id_{\PC(F)} \to \id_{\PC(F)}$ is a natural isomorphism. Similarly, we also have that $\rho:\id_{\PC(F)} \otimes I \to \id_{\PC(F)}$ is a natural isomorphism.
	
	We now need only to verify the existence of the associator(s). For this fix three objects $A, B, C  \in \PC(F)_0$ and define $\alpha_{A,B,C}:(A \otimes B) \otimes C \to A \otimes (B \otimes C)$ via
	\begin{align*}
		&\alpha_{A,B,C} \\
		&:= \left\lbrace \quot{\alpha_{A,B,C}}{X}:\left(\quot{A}{X} \os{X}{\otimes} \quot{B}{X}\right) \os{X}{\otimes} \quot{C}{X} \xrightarrow{\cong} \quot{A}{X} \os{X}{\otimes} \left(\quot{B}{X} \os{X}{\otimes} \quot{C}{X}\right) \; : \; X \in \Cscr_0\right\rbrace.
	\end{align*}
	We now prove that this is an $\PC(F)$ morphism. Fix an $f \in \Cscr(X, Y)$ and write $f:X \to Y$. Then we compute that
	\begin{align*}
		\tau_f^{(A \otimes B) \otimes C} &= \left(\tau_f^{A \otimes B} \os{X}{\otimes} \tau_f^C\right) \circ \theta_f^{A \otimes B,C} \\
		&= \left(\left(\tau_f^A \os{X}{\otimes} \tau_f^B\right) \os{X}{\otimes} \tau_f^C\right)  \circ \left( \theta_f^{A,B}\os{X}{\otimes}\id_{F(f)\quot{C}{Y}}\right) \circ \theta_f^{A \otimes B,C}
	\end{align*}
	while on the other hand
	\begin{align*}
		\tau_f^{A \otimes (B \otimes C)} &= \left(\tau_f^A \os{X}{\otimes} \tau_f^{B \otimes C}\right) \circ \theta_{f}^{A,B \otimes C} \\
		&= \left(\tau_f^A \os{X}{\otimes} \left(\tau_f^B \os{X}{\otimes} \tau_f^C\right)\right) \circ \left(\id_{F(f)\quot{A}{Y}} \os{X}{\otimes} \theta_{f}^{B,C}\right) \circ \theta_{f}^{A,B \otimes C}.
	\end{align*}
	Furthermore, using the pseudofunctoriality of the tensor functors and the naturality of the associators in each fibre category gives that the diagram
	\[
	\begin{tikzcd}
		F(\of)\left(\left(\quot{A}{Y} \os{Y}{\otimes} \quot{B}{Y}\right)\os{Y}{\otimes} \quot{C}{Y}\right) \ar[ddddd, bend left = 90]{}{F(\of)\quot{\alpha_{A,B,C}}{Y}} \ar[d]{}{\theta_f^{A \otimes B, C}} \\
		F(\of)\left(\quot{A}{Y} \os{Y}{\otimes}\quot{B}{Y}\right)\os{X}{\otimes} F(\of)(\quot{C}{Y}) \ar[d]{}{\theta_f^{A,B} \os{X}{\otimes} \id_{F(\of)\quot{C}{Y}}} \\
		\left(F(\of)(\quot{A}{Y}) \os{X}{\otimes} F(\of)(\quot{B}{Y})\right) \os{X}{\otimes} F(\of)(\quot{C}{Y}) \ar[d]{}{\quot{\alpha_{F(\of)A,F(\of)B,F(\of)C}}{X}} \\
		F(\of)(\quot{A}{Y}) \os{X}{\otimes}\left( F(\of)(\quot{B}{Y}) \os{X}{\otimes} F(\of)(\quot{C}{Y}) \right) \ar[d]{}{\id_{F(\of)\quot{A}{Y}} \os{X}{\otimes} \left(\theta_f^{B,C}\right)^{-1}} \\
		F(\of)(\quot{A}{Y}) \os{X}{\otimes} F(\of)\left(\quot{B}{Y} \os{Y}{\otimes} \quot{C}{Y}\right) \ar[d]{}{\left(\theta_f^{A,B \otimes C}\right)^{-1}} \\
		F(\of)\left(\quot{A}{Y} \os{Y}{\otimes} \left(\quot{B}{Y} \os{Y}{\otimes} \quot{C}{Y}\right)\right)
	\end{tikzcd}
	\]
	commutes. Substituting the induced identities above yields that the composite morphism $\tau_f^{A \otimes (B \otimes C)} \circ F(f)(\quot{\alpha_{A,B,C}}{Y})$ is equal to 
	\begin{align*}
		& \left(\tau_f^A \os{X}{\otimes} \left(\tau_f^B \os{X}{\otimes} \tau_f^C\right)\right) \circ \left(\id_{F(f)\quot{A}{Y}} \os{X}{\otimes} \theta_{f}^{B,C}\right) \circ \theta_{f}^{A,B \otimes C} \circ F(f)\quot{\alpha_{A,B,C}}{Y} \\
		&= \left(\tau_f^A \os{X}{\otimes} \left(\tau_f^B \os{X}{\otimes} \tau_f^C\right)\right) \circ \left(\id_{F(f)\quot{A}{Y}} \os{X}{\otimes} \theta_{f}^{B,C}\right) \circ \theta_{f}^{A,B \otimes C} \\
		&\circ \left(\theta_f^{A,B\otimes C}\right)^{-1} \circ\left(\id_{F(f)\quot{A}{Y}} \os{X}{\otimes} \left(\theta_f^{B,C}\right)^{-1}\right) \\
		&\quad\circ \quot{\alpha_{F(f)A,F(f)B,F(f)C}}{X} \circ \left(\theta_f^{A,B} \os{X}{\otimes} \id_{F(f)\quot{C}{Y}}\right) \circ \theta_f^{A \otimes B, C} \\
		&= \left(\tau_f^A \os{X}{\otimes} \left(\tau_f^B \os{X}{\otimes} \tau_f^C\right)\right) \circ \quot{\alpha_{F(f)A,F(f)B,F(f)C}}{X} \\
		&\circ \left(\theta_f^{A,B} \os{X}{\otimes} \id_{F(f)\quot{C}{Y}}\right) \circ \theta_f^{A \otimes B, C} \\
		&=\quot{\alpha_{A,B,C}}{X} \circ \left(\left(\tau_f^A \os{X}{\otimes} \tau_f^B\right) \os{X}{\otimes} \tau_f^C\right) \circ \left(\theta_f^{A,B} \os{X}{\otimes} \id_{F(f)\quot{C}{Y}}\right) \circ \theta_f^{A \otimes B, C} \\
		&= \quot{\alpha_{A,B,C}}{X} \circ \tau_f^{(A \otimes B) \otimes C},
	\end{align*}
	which in turn proves that $\tau_f^{A \otimes (B \otimes C)} \circ F(f)\quot{\alpha_{A,B,C}}{Y} =\quot{\alpha_{A,B,C}}{X} \circ \tau_f^{(A \otimes B) \otimes C}.$ Thus $\alpha_{A,B,C}$ is a morphism in $\PC(F)$ (and in fact an isomorphism because each $\quot{\alpha}{X}$ is an isomorphism). Moreover, the naturality of $\alpha$ is immediate from the fact that each $\quot{\alpha}{X}$ is a natural isomorphism.
\end{proof}

Another important character in the monoidal categorical and representation-theoretic worlds is that of braided monoidal categories. These are monoidal categories $\Cscr$ equipped with a braiding natural isomorphism
\[
\begin{tikzcd}
	\Cscr \times \Cscr \ar[rr, bend left = 40, ""{name = U}]{}{\otimes} \ar[rr, bend right = 40, swap, ""{name = D}]{}{\otimes \circ s} & & \Cscr \ar[from = U, to = D, Rightarrow, shorten <= 4pt, shorten >= 4pt]{}{\beta}
\end{tikzcd}
\] 
where $s:\Cscr \times \Cscr \to \Cscr \times \Cscr$ is the switching isomorphism generated by $(A,B) \mapsto (B,A)$ on objects and $(f,g) \mapsto (g,f)$ on morphisms. This natural isomorphism gives a monoidally coherent way of interchanging and commuting the objects through the tensor ``up to some constant/relation'' and also a categorification of noncommutative algebra. Because of the use of braided monoidal categories in representation theory, the study of TQFTs, noncommutative algebra, Hopf algebras, algebraic topology, and higher algebra we are interested in when these structures can be translated to $\PC(F)$ and hence when they arise on equivariant categories (in the sense of those categories defined in \cite{DoretteMe}, \cite{MyThesis}, and below). I expect these to be of use especially when studying and working in the context of equivariant higher algebra and the study of generalized equivariant cohomology theories.

\begin{definition}
	A braided monoidal category\index[terminology]{Monoidal Category! Braided} is a monoidal category $(\Cscr, \otimes, I)$ together with a natural isomorphism
	\[
	\begin{tikzcd}
		\Cscr \times \Cscr \ar[rr, bend left = 40, ""{name = U}]{}{\otimes} \ar[rr, bend right = 40, swap, ""{name = D}]{}{\otimes \circ s} & & \Cscr \ar[from = U, to = D, Rightarrow, shorten <= 4pt, shorten >= 4pt]{}{\beta}
	\end{tikzcd}
	\] 
	such that for all objects $A, B, C \in \Cscr_0$ the diagrams
	\[
	\begin{tikzcd}
		(A \otimes B) \otimes C \ar[rr]{}{\alpha_{A, B, C}} \ar[d, swap]{}{\beta_{A,B} \otimes \id_C} & & A \otimes (B \otimes C) \ar[rr]{}{\beta_{A, B \otimes C}} & &  (B \otimes C) \otimes A \ar[d]{}{\alpha_{B, C, A}} \\
		(B \otimes A) \otimes C \ar[rr, swap]{}{\alpha_{B, A, C}} & & B \otimes (A \otimes C) \ar[rr, swap]{}{\id_Y \otimes \beta_{A,C}} & & B \otimes (C \otimes A)
	\end{tikzcd}
	\]
	and
	\[
	\begin{tikzcd}
		A \otimes (B \otimes C) \ar[d, swap]{}{\id_A \otimes \beta_{B,C}} \ar[rr]{}{\alpha_{A,B,C}^{-1}} & & (A \otimes B) \otimes C \ar[rr]{}{\beta_{A \otimes B, C}} & & C \otimes (A \otimes B) \ar[d]{}{\alpha_{C, A, Y}^{-1}} \\
		A \otimes (C \otimes B) \ar[rr, swap]{}{\alpha_{A,C,B}^{-1}} & & (A \otimes C) \otimes B \ar[rr, swap]{}{\beta_{A, C} \otimes \id_B} & & (C \otimes A) \otimes B
	\end{tikzcd}
	\]
	commute. Additionally, if $\beta_{B,A} \circ \beta_{A,B} = \id_{A \otimes B}$ for all $A, B \in \Cscr_0$ we say that $\Cscr$ is a symmetric monoidal category.
\end{definition}
\begin{definition}
	Let $(\Cscr, \otimes, I, \quot{\beta}{\Cscr})$ and $(\Dscr, \boxtimes, J, \quot{\beta}{\Dscr})$ be braided monoidal categories. A monoidal functor $F:\Dscr \to \Cscr$ is braided monoidal\index[terminology]{Monoidal Functor! Braided} if for all $A, B \in \Dscr_0$ the diagram
	\[
	\begin{tikzcd}
		F\left(A \boxtimes B\right) \ar[rr]{}{\theta_{A, B}} \ar[d, swap]{}{F\left(\quot{\beta_{A,B}}{\Dscr}\right)} & & F(A) \otimes F(B) \ar[d]{}{\quot{\beta_{A,B}}{\Cscr}} \\
		F\left(B \boxtimes A\right) \ar[rr, swap]{}{\theta_{B,A}} & & F(B) \otimes F(A)
	\end{tikzcd}
	\]
	commutes.
\end{definition}
\begin{definition}
	Let $F:\Cscr^{\op} \to \fCat$ be a pseudofunctor.\index[terminology]{Pseudofunctor! Braided Monoidal} We say that $F$ is a braided monoidal pseudofunctor if $F$ is a monoidal psuedofunctor and additionally every category $F(X)$ is braided monoidal while every functor $F(f)$ is also braided monoidal.
\end{definition}
We now show that when the pseudofunctor $F$ is braided monoidal so is the category $\PC(F)$.

\begin{proposition}\label{Prop: Section 2: Equivariant cat is symmetric monoidal}
	Let $F:\Cscr^{\op} \to \fCat$ be a braided monoidal pseudofunctor. Then $\PC(F)$ is braided monoidal with braiding $\beta_{A,B}:A \otimes B \to B \otimes A$ given by
	\[
	\beta_{A,B} := \left\lbrace \quot{\beta_{A,B}}{X} \; | \; X \in \Cscr_0 \right\rbrace
	\]
	where $\quot{\beta_{A,B}}{X}$ is a shorthand for the braiding isomorphism $\quot{\beta_{\quot{A}{X},\quot{B}{X}}}{X}$ in $F(X)$. Additionally if each category $F(X)$ is symmetric monoidal then so is $\PC(F)$.
\end{proposition}
\begin{proof}
	We know from Theorem \ref{Theorem: Section 2: Monoidal preequivariant pseudofunctor gives monoidal equivariant cat} that $\PC(F)$ is monoidal; we now only need to furnish $\PC(F)$ with the braiding structure suggested in the statement of the proposition and prove that it gives the desired natural isomorphism :
	\[
	\begin{tikzcd}
		\PC(F) \times \PC(F) \ar[rr, bend left = 40, ""{name = U}]{}{\otimes} \ar[rr, bend right = 40, swap, ""{name = D}]{}{\otimes \circ s} & & \PC(F) \ar[from = U, to = D, Rightarrow, shorten <= 4pt, shorten >= 4pt]{}{\beta}
	\end{tikzcd}
	\]
	Note that the switching functor $s:\PC(F) \times \PC(F) \to \PC(F) \times \PC(F)$ is determined by $X$-locally switching the objects and morphisms.
	
	We prove first that $\beta_{A,B}$ is a morphism for any $A, B \in\PC(F)_0$. For this let $f:X \to Y$ be a morphism in $\Cscr$ and recall that
	\[
	\tau_f^{A \otimes B} = \left(\tau_f^A \os{X}{\otimes} \tau_f^B\right) \circ \tau_f^{A,B}.
	\]
	We now must prove that
	\[
	\begin{tikzcd}
		F(f)\left(\quot{A}{Y} \os{Y}{\otimes} \quot{B}{Y}\right) \ar[rr]{}{F(f)\left(\quot{\beta_{A,B}}{Y}\right)} \ar[d, swap]{}{\tau_f^{A \otimes B}} & & F(f)\left(\quot{B}{Y} \os{Y}{\otimes} \quot{A}{Y}\right) \ar[d]{}{\tau_f^{B \otimes A}} \\
		\quot{A}{X} \os{X}{\otimes} \quot{B}{X} \ar[rr, swap]{}{\quot{\beta_{A,B}}{X}} & & \quot{B}{X} \os{X}{\otimes} \quot{A}{X}
	\end{tikzcd}
	\]
	commutes in order to conclude that $\beta_{A,B}$ is a morphism. For this we simply calculate that
	\begin{align*}
		\tau_f^{B \otimes A} \circ F(f)\left(\quot{\beta_{A,B}}{Y}\right) &= \left(\tau_f^B \os{X}{\otimes} \tau_f^A\right) \circ \theta_f^{B,A} \circ F(f)\left(\quot{\beta_{A,B}}{Y}\right) \\
		&= \left(\tau_f^B \os{X}{\otimes} \tau_f^A\right) \circ \quot{\beta_{F(f)A,F(f)B}}{X} \circ \theta_{f}^{A,B} \\
		&= \quot{\beta_{A,B}}{X} \circ \left(\tau_f^A \os{X}{\otimes} \tau_f^B\right) \circ \theta_f^{A,B} \\
		&= \quot{\beta_{A,B}}{X} \circ \tau_{f}^{A \otimes B},
	\end{align*}
	as was desired.
	
	We now prove that $\beta$ is a natural transformation. Let $(A,B)$ and $(C,D)$ be objects in $\PC(F) \times \PC(F)$ and let $(P,\Phi):(A,B) \to (C,D)$ be a morphism. Then we calculate that for any $X \in \Cscr_0$,
	\begin{align*}
		\quot{\left(\Phi \otimes P\right)}{X} \circ \quot{\beta_{A,B}}{X} &= \left(\quot{\phi}{X} \otimes \quot{\rho}{X}\right)	\circ \quot{\beta_{A,B}}{X} = \quot{\beta_{C,D}}{X} \circ \left(\quot{\rho}{X} \otimes \quot{\varphi}{X}\right) \\
		&=\quot{\beta_{C,D}}{X} \circ \quot{\left(P \otimes \Phi\right)}{X}
	\end{align*}
	by the fact that each $\quot{\beta_{-,-}}{X}$ is a natural transformation. Thus the diagram
	\[
	\begin{tikzcd}
		A \otimes B \ar[r]{}{\beta_{A,B}} \ar[d, swap]{}{P \otimes \Phi} & B \otimes A \ar[d]{}{\Phi \otimes P} \\
		C \otimes D \ar[r,swap]{}{\beta_{C,D}} & D \otimes C
	\end{tikzcd}
	\]
	commutes in $\PC(F)$ and so $\beta$ is a natural isomorphism. That it satisfies the Hexagon Axioms is straightforward by using that each $X$-local braiding isomorphism satisfies the corresponding braiding isomorphism. Thus $\PC(F)$ is a braided monoidal category. Finally, if every braiding gives a symmetric monoidal structure on $F(X)$ we then calculate that
	\[
	\beta_{B,A} \circ \beta_{A,B} = \left\lbrace \quot{\beta_{B,A}}{X} \circ \quot{\beta_{A,B}}{X} \; | \; X \in \Cscr_0 \right\rbrace = \left\lbrace \quot{\id_{A \otimes B}}{X} \; | \; X \in \Cscr_0 \right\rbrace = \id_{A \otimes B}
	\]
	which shows that in this case $\PC(F)$ is also symmetric monoidal.
\end{proof}
\begin{example}\label{Example: Tensor Functor on EDC}
	There are two standard examples of monoidal pre-equivariant pseudofunctors that we should keep in mind. The first is that of the pre-equivariant pseudofunctor $F:\SfResl_G(X)^{\op} \to \fCat$ given on objects by $F(\XGamma) = \DbQl{\XGamma}$ and on morphisms by $F(\of) = \of^{\ast}$. Then each category $F(\XGamma)$ is monoidal with tensor functor
	\[
	\otimes_{X} := (-) \overset{L}{\otimes}_{\Gamma} (-).
	\]
	Note that the functor $(-)\otimes_{\Gamma}^{L}(-):\DbQl{\XGamma} \times \DbQl{\XGamma} \to \DbQl{\XGamma}$ is the (left) derived tensor functor of $\ell$-adic sheaves on $\XGamma$.
	By Theorem \ref{Theorem: Section 2: Monoidal preequivariant pseudofunctor gives monoidal equivariant cat} and Proposition \ref{Prop: Section 2: Equivariant cat is symmetric monoidal} this produces a symmetric monoidal functor $(-)\otimes^{L}(-)$ on $\DbeqQl{X}$.
	
	The second example to keep in mind is the pre-equivariant pseudofunctor $F:\SfResl_G(X)^{\op} \to \fCat$ given on objects by
	\[
	F(\XGamma) = \Shv(\XGamma, \text{{\'e}t})
	\]
	and on morphisms by
	\[
	F(\of) = \of^{\ast}:\Shv(\XGammap,\text{{\'e}t}) \to \Shv(\XGamma,\text{{\'e}t}).
	\]
	The monoidal functor in this case is defined to be the one induced by
	\[
	\otimes_{\Gamma} := (-) \times_{\Gamma} (-)
	\]
	This then induces an alternative construction of the product functor on the category $\Shv_G(X,\text{{\'e}t})$ by emphasizing its monoidal properties.
\end{example}
\newpage

\section{Regularity and Subobject Classifiers for Pseudocones}\label{Subsection: PCF Regualr}
We now move to studying some cases when equivariant categories are regular categories and when equivariant categories have subobject classifiers.(cf.\@ Definition \ref{Defn: Subobject Classifier} below). We will defer the study of adjunctions to the next chapter where functors between pseudocone categories are studied in significantly higher depth than we have seen so far.

\begin{proposition}\label{Prop: Regularity of Equivariant Category}
Let $F:\Cscr^{\op} \to\fCat$ be a pseudofunctor such that each fibre category $F(X)$ is regular and for which each fibre functor $F(f)$ is finitely complete and preserves regular epimorphisms. Then $\PC(F)$ is regular.
\end{proposition}
\begin{proof}
	Recall that a category $\Cscr$ is regular if and only if $\Cscr$ is finitely complete, admits all coequalizers of kernel pairs, and has the property that regular epimorphisms are stable under base change, i.e., regular epimorphisms are stable under pullback. We now proceed to show that $\PC(F)$ is regular by first showing that it is finitely complete; however, for this we simply apply Corollary \ref{Cor: Equivariant cat is finitely complete or cocomplete}.
	
	To see that $\PC(F)$ admits coequalizers of kernel pairs, let $P \in \PC(F)(A,B)$ with $P = \lbrace \quot{\rho}{X} \; | \; X \in \Cscr_0\rbrace$. Now take the pullback
	\[
	\xymatrix{
		A \times_B A \pullbackcorner \ar[r]^-{\Pi_0} \ar[d]_{\Pi_1} & A \ar[d]^{P} \\
		A \ar[r]_{P} & B
	}
	\]
	and note that by Theorem \ref{Thm: Section 2: Equivariant Cat has lims} the pullback takes the form
	\[
	A \times_B A := \lbrace \quot{A}{X} \times_{\quot{B}{X}} \quot{A}{X} \; | \; X \in \Cscr_0 \rbrace
	\]
	with the induced transition isomorphisms, while for $k \in \lbrace 0, 1\rbrace$ the morphisms $\Pi_k$ take the form
	\[
	\Pi_k := \lbrace \quot{\pi_k}{X}:\quot{A}{X} \times_{\quot{B}{X}} \quot{A}{X} \to \quot{A}{X} \; | \; X \in \Cscr_0 \rbrace
	\]
	where each $\quot{\pi_k}{X}$ is the corresponding first or second pullback map in the fibre category $F(X)$. Now define the $\Vscr$-set
	\[
	C := \lbrace \Coeq(\quot{\pi_0}{X},\quot{\pi_1}{X}) \; | \; X \in \Cscr_0 \rbrace.
	\]
	By assumption, each coequalizer exists in each category $F(X)$ and each functor $F(f)$ preserves these coequalizers. Thus by Theorem \ref{Theorem: Section 2: Equivariant Cat has colims} $(C,T_C)$ is an object in $\PC(F)$ (with induced transition isomorphisms $T_C$) and satisfies $C \cong \Coeq(\Pi_0,\Pi_1)$ with $\coeq(\Pi_0,\Pi_1):A \to C$ given by the set 
	\[\coeq(\Pi_0, \Pi_1) = \lbrace \coeq(\quot{\pi_0}{X}, \quot{\pi_1}{X}) \; | \; X \in \Cscr_0 \rbrace.
	\] 
	Thus $\PC(F)$ admits coequalizers of kernel pairs.
	
	We now verify that regular epimorphisms are stable under base change. Assume that $\Sigma:A \to C$ is a regular epimorphism in $\PC(F)$ and let $P:B \to C$ be an arbitrary morphism. Take the pullback
	\[
	\xymatrix{
		A \times_C B \ar[r]^-{\Phi} \ar[d]_{Z} \pullbackcorner & A \ar[d]^{\Sigma} \\
		B \ar[r]_{P} & C
	}
	\]
	where $Z := \lbrace \quot{\zeta}{X} \; | \; X \in \Cscr_0\rbrace$. Now let $(\Pi_0, \Pi_1)$ be the kernel pair of $Z$ and let $\Coeq(\Pi_0, \Pi_1)$ be the coequalizer of the kernel pair; note that again by Theorem \ref{Theorem: Section 2: Equivariant Cat has colims}, we have that 
	\[
	\Coeq(\Pi_0, \Pi_1) = \lbrace \Coeq(\quot{\pi_0}{X},\quot{\pi_1}{X}) \; | \; X \in \Cscr_0 \rbrace
	\] 
	and 
	\[
	\coeq(\Pi_0, \Pi_1) = \lbrace \coeq(\quot{\pi_0}{X}, \quot{\pi_1}{X}) \; | \; X \in \Cscr_0 \rbrace.
	\]
	Because each category $F(X)$ is regular, we have that $\quot{\zeta}{X} \cong \coeq(\quot{\pi_0}{X},\quot{\pi_1}{X})$ via a unique isomorphism $\quot{\theta}{X}$ fitting into the diagram
	\[
	\xymatrix{
		\quot{A}{X} \times_{\quot{C}{X}} \quot{B}{X} \ar@{=}[d] \ar[rr]^-{\coeq(\quot{\pi_0}{X},\quot{\pi_1}{X})} & & \Coeq(\quot{\pi_0}{X}, \quot{\pi_1}{X}) \ar@{-->}[d]^{\exists!\quot{\theta}{X}}_{\cong} \\
		\quot{A}{X} \times_{\quot{C}{X}} \ar[rr]_-{\quot{\zeta}{X}} & & \quot{B}{X}
	}
	\]
	in $F(X)$. We need to check now that $\Theta = \lbrace \quot{\theta}{X} \; | \; X \in \Cscr_0 \rbrace$ is an $\PC(F)$ morphism, as it is automatically an isomorphism if it is an $\PC(F)$-morphism. For this we must prove that for all $f \in \Cscr_1$ (say with $\Dom f = X$ and $\Codom f = Y$),
	\[
	\tau_f^B \circ F(f)\quot{\theta}{Y} \overset{?}{=} \quot{\theta}{X} \circ \tau_f^{\Coeq(\Pi_0, \Pi_1)}.
	\]
	For this we recall that $Z:A \times_B A \to B$ and that $\coeq(\Pi_0, \Pi_1):A \times_B A \to \Coeq(\Pi_0, \Pi_1)$ are morphisms in $\PC(F)$. We then calculate that on one hand
	\begin{align*}
		\quot{\zeta}{X} \circ \tau_f^{A\times_C B} &= \quot{\theta}{X} \circ \coeq(\quot{\pi_0}{X}, \quot{\pi_1}{X}) \circ \tau_f^{A \times_C B} \\
		&= \quot{\theta}{X} \circ \tau_f^{\Coeq(\Pi_0,\Pi_1)} \circ F(f)\coeq(\quot{\pi_0}{Y}, \quot{\pi_1}{Y})
	\end{align*}
	while on the other hand
	\begin{align*}
		\quot{\zeta}{X} \circ \tau_f^{A\times_C B} &= \tau_f^B \circ F(f)\quot{\zeta}{Y} = \tau_f^B \circ F(f)\big(\quot{\theta}{Y} \circ \coeq(\quot{\pi_0}{X},\quot{\pi_1}{X})\big) \\
		&= \tau_f^B \circ F(f)\quot{\theta}{Y} \circ F(f)\coeq(\quot{\pi_1}{Y},\quot{\pi_1}{Y}).
	\end{align*}
	Because $F(f)$ preserves regular epimorphisms, it follows that $F(f)\coeq(\quot{\pi_0}{Y},\quot{\pi_1}{Y})$ is epic. Thus, after canceling $F(f)\coeq(\quot{\pi_0}{Y},\quot{\pi_1}{Y})$ we get that
	\[
	\tau_f^B \circ F(f)\quot{\theta}{Y} = \quot{\theta}{X} \circ \tau_f^{\Coeq(\Pi_0,\Pi_1)}.
	\]
	This shows that $\Theta$ is a morphism in $\PC(F)$ and hence proves that $Z$ is a regular epimorphism. This completes the proof that $\PC(F)$ is regular.
\end{proof}

Our next and final proposition for this section lies further in this direction of pseudocone-based categorical logic. We will show below that under mild assumptions that if all the fibre categories have subobject classifiers and if the fibre functors preserve these subobject classifiers, terminal objects, and subobject classification pullbacks, then $\PC(F)$ has a subobject classifier as well. 
\begin{definition}[{\cite[Chapter 1.3]{MacLaneMoerdijk}}]\label{Defn: Subobject Classifier}
	A category $\Cscr$ has a subobject classifier\index[terminology]{Subobject Classifier} $\Omega$ if $\Cscr$ has a terminal object $\top$ and there is an object $\Omega \in \Cscr_0$ with a morphism $\true:\top \to \Omega$ such that for any monic $\mu:A \to B$ in $\Cscr$, a unique morphism $\chi_{\mu}:B \to \Omega$ making the diagram
	\[
	\xymatrix{
		A \ar@{-->}[r]^{\exists!} \ar[d]_-{\mu} & \top \ar[d]^{\true} \\
		B \ar[r]_{\exists!\chi_{\mu}} & \Omega
	}
	\]
	into a pullback diagram.
\end{definition}
\begin{example}
	Let $\Cscr = \Set$ be the category of sets. Then the set $\Omega = \lbrace 0, 1 \rbrace$ is the subobject classier. The map $\true:\lbrace \ast \rbrace \to \Omega$ is given by $\ast \mapsto 1$ and for any monic $\mu:A \to B$ in $\Set$ the morphism $\chi_{\mu}:B \to \Omega$ is given by
	\[
	\chi_{\mu}(b) := \begin{cases}
		1 & \exists\,a \in A.\, \mu(a) = b; \\
		0 & \text{else}.
	\end{cases}
	\]
\end{example}
\begin{example}
Let $X$ be a topological space. Then in the category $\Shv(X)$ the subobject classifier is the sheaf $\Omega$ with object assignment
\[
\Omega(V) = \Top(V,\lbrace \eta, s\rbrace)
\]
where $\lbrace \eta, s \rbrace$ is the Sierpinski space with $\eta$ the open point and $s$ the closed point.
\end{example}
\begin{remark}
	Let $\Cscr$ and $\Dscr$ be categories with subobject classifiers $\Omega_{\Cscr}$ and $\Omega_{\Dscr}$, respectively. We say that a functor $G:\Cscr \to \Dscr$ {preserves subobject classifing pullbacks}\index[terminology]{Subobject Classifying Pullbacks} when it preserves subobject classifiers, terminal objects, and if for a monic $\mu \in \Cscr(A,B)$ which is preserved by $G$ the diagram\index[notation]{true@$\true$}
	\[
	\xymatrix{
		GA \ar[d]_{G\mu} \ar@{-->}[r]^{G!_{A}} \ar[d]_{G\mu} & G\top_{\Cscr} \ar[d]^{G\true_{\Cscr}} \\
		GB \ar[r]_{G\chi(\mu)} & G\Omega_{\Cscr}
	}
	\]
	is isomorphic to the pullback diagram
	\[
	\xymatrix{
		GA \ar[d]_{G\mu} \pullbackcorner\ar[r]^{!_{GA}} & \top_{\Dscr} \ar[d]^{\true_{\Dscr}} \\
		GB \ar[r]_{\chi(G\mu)} & \Omega_{\Dscr}
	}
	\]
	in $\Dscr$. That is, there are isomorphisms $G\top_{\Cscr} \to \top_{\Dscr}$ and $G\Omega_{\Cscr} \to \Omega_{\Dscr}$ for which the diagram
	\[
	\begin{tikzcd}
		& GA \ar[rr]{}{!_{GA}} \ar[dd,swap, near end]{}{G\mu} & & \top_{\Dscr} \ar[dd]{}{\true_{\Dscr}} \\
		GA \ar[equals, ur] \ar[dd, swap]{}{G\mu} \ar[rr, near end, crossing over]{}{G!_{A}} & & G\top_{\Cscr} \ar[ur]{}{\cong}  \\
		& GB  \ar[rr,swap, near start]{}{\chi(G\mu)} & & \Omega_{\Dscr} \\
		GB \ar[ur, equals] \ar[rr,swap]{}{G\chi(\mu)} & & G\Omega_{\Cscr} \ar[ur, swap]{}{\cong} \ar[from = 2-3, to = 4-3, near start, crossing over]{}{G\true_{\Cscr}}
	\end{tikzcd}
	\]
	commutes in $\Dscr$.
\end{remark}
\begin{proposition}\label{Prop: Section 2: Equivariant Cat has Subobject Classifiers}
	Let $F$ be a pseudofunctor such that each fibre category $F(X)$ has a subobject classifier, $\quot{\Omega}{X}$, and for which each fibre functor $F(f)$ preserves subobject classifiers, terminal objects, and subobject classifying pullbacks. Then $\PC(F)$ admits a subobject classifier $\Omega$ whose object collection has the form
	\[
	\Omega = \lbrace \quot{\Omega}{X} \; | \; X \in \Cscr_0 \rbrace.
	\]
\end{proposition}
\begin{proof}
	We begin by defining our subobject classifier. Because the statement of the proposition gave $\Omega = \lbrace \quot{\Omega}{X} \; | \; X \in \Cscr_0 \rbrace$, we only need to give the transition isomorphisms $T_{\Omega}$. For each $f \in \Cscr_1$, say with $\Dom f = X$ and $\Codom f = Y$, let $\theta_f:F(f)(\quot{\Omega}{Y}) \xrightarrow{\cong} \quot{\Omega}{X}$ be the witness isomorphism induced by the fact that each $F(f)$ preserves subobject classifiers. We then define the collection of transition isomorphisms by declaring that for all $f \in \Cscr_1$, $\tau_f^{\Omega} := \theta_f.$
	
	Let us now verify the cocycle condition on the transition isomorphisms. Fix a composable pair of arrows $X \xrightarrow{f} Y \xrightarrow{g} Z$ in $\Cscr$. We now verify that
	\[
	\tau_f^{\Omega} \circ F(f)\tau_g^{\Omega} = \tau_{g \circ f}^{\Omega} \circ \phi_{f,g}.
	\]
	For this begin by noting that the map $\tau_{g \circ f}^{\Omega}$ is the unique isomorphism between $F(g \circ f)(\quot{\Omega}{Z})$ and $\quot{\Omega}{X}$; similarly, since $F(f)\tau_g^{\Omega}$ is the unique isomorphism between $F(f)\big(F(g)(\quot{\Omega}{Z})\big)$ and  $F(f)(\quot{\Omega}{Y})$, while $\tau_f^{\Omega}$ is the unique isomorphism between $F(f)(\quot{\Omega}{Y})$ and $\quot{\Omega}{X}$, we will be done if we know that $\phi_{f,g}^{\quot{\Omega}{Z}}$ preserves subobject classifiers. However, this is immediate, as $\phi_{f,g}^{\quot{\Omega}{Z}}$ is a natural isomorphism and both $F(f)(F(g)(\quot{\Omega}{Z}))$ and $F(g \circ f)(\quot{\Omega}{Z})$ are subobject classifiers. Thus the diagram
	\[
	\xymatrix{
		F(f)\big(F(g)(\quot{\Omega}{Z})\big) \ar[rr]^-{F(f)\tau_g^{\Omega}} \ar[d]_{\phi_{f,g}^{\quot{\Omega}{Z}}} & & F(f)(\quot{\Omega}{Y}) \ar[d]^{\tau_f^{\Omega}} \\
		F(g \circ f)(\quot{\Omega}{Z}) \ar[rr]_-{\tau_{g \circ f}^{\Omega}} & & \quot{\Omega}{X}
	}
	\]
	commutes, proving that $(\Omega,T_{\Omega})$ is an object in $\PC(F)$.
	
	We now construct the map $\true:\top \to \Omega$ in $\PC(F)$. Define
	\[
	\true := \lbrace \quot{\true}{X}:\quot{\top}{X} \to \quot{\Omega}{X} \; | \; X \in \Cscr_0\rbrace;
	\]
	to prove that this is an $\PC(F)$ morphism, we must verify that for all $f \in \Cscr_1$,
	\[
	\quot{\true}{X} \circ \tau_{f}^{\top} \overset{?}{=} \tau_{f}^{\Omega} \circ F(f)\quot{\true}{Y}.
	\]
	For this, however, note that because each functor preserves object classification pullbacks, given any monic $\mu:C \to D$ in $\PC(F)$ we have that the cube
	\[
	\begin{tikzcd}
		& F(f)(\quot{C}{Y}) \ar[rr]{}{!_{F(f)C}} \ar[dd,swap, near end]{}{F(f)\mu} & & \quot{\top}{X} \ar[dd]{}{\quot{\true}{X}} \\
		F(f)(\quot{C}{Y}) \ar[equals, ur] \ar[dd, swap]{}{F(f)\mu} \ar[rr, near end, crossing over]{}{F(f)!_{A}} & & F(f)(\quot{\top}{Y}) \ar[ur]{}{\cong} \\
		& F(f)(\quot{D}{Y})  \ar[rr,swap, near start]{}{\chi(F(f)\mu)} & & \quot{\Omega}{X} \\
		F(f)(\quot{D}{Y}) \ar[ur, equals] \ar[rr,swap]{}{F(f)\chi(\mu)} & & F(f)(\quot{\Omega}{Y}) \ar[ur, swap]{}{\cong}  \ar[from = 2-3, to = 4-3, near start, crossing over]{}{F(f)\quot{\true}{Y}}
	\end{tikzcd}
	\]
	commutes in $F(X)$. Taking the right-most face shows that the square
	\[
	\xymatrix{
		F(f)(\quot{\top}{Y}) \ar[d]_{F(f)\quot{\true}{Y}} \ar[r]^-{\cong} & \quot{\top}{X} \ar[d]^{\quot{\true}{X}} \\
		F(f)(\quot{\Omega}{Y}) \ar[r]_-{\cong} & \quot{\Omega}{X}
	}
	\]
	commutes; because the horizontal isomorphisms in the square above are the unique isomorphisms between the objects, i.e., the witnessing isomorphisms, we have that the diagram takes the form
	\[
	\xymatrix{
		F(f)(\quot{\top}{Y}) \ar[d]_{F(f)\quot{\true}{Y}} \ar[r]^-{\tau_f^{\top}} & \quot{\top}{X} \ar[d]^{\quot{\true}{X}} \\
		F(f)(\quot{\Omega}{Y}) \ar[r]_-{\tau_f^{\Omega}} & \quot{\Omega}{X}
	}
	\]
	instead. However, this is exactly the verification that $\true$ is an $\PC(F)$ morphism.
	
	We now construct the classifying map of an arbitrary monic in $\PC(F)$. Begin by letting $M = \lbrace \quot{\mu}{X} \; | \; X \in \Cscr_0\rbrace$ be a monic $M \in \PC(F)(A,B)$; note that each map $\quot{\mu}{X}$ is monic in $F(X)$ by Theorem \ref{Thm: Section 2: Equivariant Cat has lims}. Then in each fibre category $F(X)$ we have that there are (pullback) diagrams
	\[
	\xymatrix{
		\quot{A}{X} \ar[r] \ar[d]_{\quot{\mu}{X}} & \quot{\top}{X} \ar[d]^{\quot{\true}{X}} \\
		\quot{B}{X} \ar[r]_-{\chi(\quot{\mu}{X})} & \quot{\Omega}{X}
	}
	\]
	We thus define $\chi(M):B \to \Omega$ via
	\[
	\chi(M) := \lbrace \chi(\quot{\mu}{X}) \; | \; X \in \Cscr_0, \quot{\mu}{X} \in M \rbrace;
	\]
	note that we now must show that for all $f \in \Cscr_1$ with $f:X \to Y$,
	\[
	\tau_f^{\Omega} \circ F(f)\chi(\quot{\mu}{Y}) \overset{?}{=} \chi(\quot{\mu}{X}) \circ \tau_f^{B}.
	\]
	For this we note that since classifying maps are uniquely determined, it suffices to prove that the diagram
	\[
	\xymatrix{
		\quot{A}{X} \ar[d]_{\quot{\mu}{X}} \ar@{-->}[rrrr]^-{\exists!} & & & & \quot{\top}{X} \ar[d]^{\quot{\true}{X}} \\
		\BGamma \ar[rrrr]_-{\tau_f^{\Omega} \circ F(\of)\big(\chi(\quot{\mu}{Y})\big) \circ \big(\tau_f^B\big)^{-1}} & & & & \quot{\Omega}{X}
	}
	\]
	is a pullback diagram. Because the diagram is isomorphic to the standard pullback with bottom edge given by the morphism $\chi(\quot{\mu}{X})$, showing that the diagram is a pullback is immediate once we know that it commutes. For this we calculate that
	\begin{align*}
		\tau_f^{\Omega} \circ F(\of)\left(\chi(\quot{\mu}{Y})\right) \circ \left(\tau_f^B\right)^{-1} \circ \quot{\mu}{X} &= \tau_f^{\Omega} \circ F(\of)\left(\chi(\quot{\mu}{Y})\right) \circ F(\of)\left(\quot{\mu}{Y}\right) \circ \left(\tau_f^A\right)^{-1} \\
		&= \tau_f^{\Omega} \circ F(\of)\left(\chi(\quot{\mu}{Y}) \circ \quot{\mu}{Y}\right) \circ \left(\tau_f^A\right)^{-1} \\
		&= \tau_f^{\Omega} \circ F(\of)\left(\quot{\true}{Y} \circ !_{\quot{A}{Y}}\right) \circ \left(\tau_f^A\right)^{-1} \\
		&= \tau_f^{\Omega} \circ F(\of)\left(\quot{\true}{Y}\right) \circ F(\of)\left(!_{\quot{A}{Y}}\right) \circ \left(\tau_f^A\right)^{-1} \\
		&=\quot{\true}{X} \circ \tau_f^{\top} \circ F(\of)\left(!_{\quot{A}{Y}}\right) \circ \left(\tau_f^A\right)^{-1} \\
		&= \quot{\true}{X} \circ !_{\quot{A}{X}} \circ \tau_f^A \circ \left(\tau_f^A\right)^{-1} \\
		&= \quot{\true}{X} \circ !_{\quot{A}{X}},
	\end{align*}
	which proves the diagram commutes. Using the pullback observation we made above it follows that
	\[
	\chi(\quot{\mu}{X}) = \tau_f^{\Omega} \circ F(\of)\left(\chi(\quot{\mu}{Y})\right) \circ \left(\tau_f^B\right)^{-1}.
	\]
	so $\chi(M)$ is indeed an  $\PC(F)$-morphism. Furthermore, we note immediately that the diagram
	\[
	\xymatrix{
		A \ar[d]_{M} \ar[r]^{!_A} & \top \ar[d]^{\true} \\
		B \ar[r]_{\chi(M)} & \Omega
	}
	\]
	commutes in $\PC(F)$, as for all $X \in \Cscr_0$ the diagrams
	\[
	\xymatrix{
		\quot{A}{X} \pullbackcorner \ar[r]^{!_{\quot{A}{X}}} \ar[d]_{\quot{\mu}{X}} & \quot{!}{X} \ar[d]^{\quot{\true}{X}} \\
		\quot{B}{X} \ar[r]_{\chi(\quot{\mu}{X})} & \quot{\Omega}{X}
	}
	\]
	commute. We also note that because all these diagrams are pullbacks, we get that
	\[
	\xymatrix{
		A \pullbackcorner \ar[d]_{M} \ar[r]^{!_A} & \top \ar[d]^{\true} \\
		B \ar[r]_{\chi(M)} & \Omega
	}
	\]
	is a pullback in $\PC(F)$. It also follows that the map $\chi(M)$ is unique, as any other map will be equal to $\chi(M)$ as it will be $X$-wise equal to $\chi(\quot{\mu}{X})$. This proves that $\Omega$ is a subobject classifier.
\end{proof}
\newpage

\chapter{Functors Between Categories of Pseudocones}\label{Section: Functors of Pseudocones}

As category theorists and mathematicians who make use of category theory, we know that a good knowledge of objects is important but only at most half of the story: to know a category we need to know not just about the objects but also about the $1$-and-higher morphisms. At this point, in Chapters \ref{Section: Category of Pseudocones} and \ref{Section: Limits and Colimits in Pseudocones} we have just got to know the main characters of our study: the pseudocone categories
\[
\PC(F) = \Bicat(\Cscr^{\op},\fCat)(\cnst(\One),F).
\]
However we also need to learn and study the ways our main characters interact with and relate to each other; we need to study the  morphisms and natural transformations between these categories. While much of the structure can be inherited from the fact that $\Bicat$ is a tricategory\footnote{As we are working generically with $\Bicat$ as determined by bicategories, pseudofunctors, pseudonatural transformations, and modifications.} (so $\Bicat(\Cscr^{\op},\fCat)$ is in particular a bicategory which, because of the strictness of $\Cscr$ and $\fCat$, is actually a $2$-category), we obtain nice descriptions of morphisms and transformations between pseudocone categories. We will also discuss how to translate between pseudocone categories $\PC(F)$ and $\PC(E)$ both when we have pseudofunctors defined on a common domain $E,F:\Cscr^{\op} \to \fCat$ and when they are instead defined on different domains $F:\Cscr^{\op} \to \fCat$ and $E:\Dscr^{\op} \to \fCat$.

Much of the reason for studying functors and transformations between pseudocone categories comes not just from categorical curiosity, but also from a desire to apply these techniques to equivariant geometric/topological situations and contexts. We would like to apply functors between equivariant categories to various flavours such as:
\begin{enumerate}
	\item Functors which change fibre categories like the inclusions $\Shv_G(X) \to D_G^b(X)$ or $\Loc_G(X) \to \Per_G(X)$ for a smooth algebraic group $G$ and a left $G$-variety $X$.
	\item Functors which change spaces between fibres like $\Shv_G(X) \to \Shv_G(Y)$, $\Shv_G(X) \to \Shv_G(Y)$, or $D_G^b(Y) \to D_G^b(X)$ when $G$ is a smooth algebraic group, $X$ and $Y$ are left $G$-varieties, and there is a $G$-equivariant morphism $h:X \to Y$.
	\item Combining the two changes of domains and fibres at the same time.
	\item Functors which change groups like $\Shv_G(X) \to \Shv_H(X)$ or $D_G^b(X) \to D_H^b(X)$ when there is a morphism $\varphi:H \to G$ of smooth algebraic groups.
	\item The topological analogues of each of the situations above.
\end{enumerate}
We will see below that Flavours (1) and (2) above can be captured by entirely the pseudocone formalism. Flavour (4), however, is ultimately a geometric and topological phenomenon and requires a use of induction spaces in both the variety-theoretic and topological contexts and so is deferred to Chapter \ref{Chapter: Chofg}.

\section{Pseudocone Functors: Change of Fibre}\label{Section: Section 3: Change of Fibre}\label{Subsection: Change of Fibre}

We begin this subsection by studying the most straightforward type of functor between pseudocone categories: the Change of Fibre functors. These will allow us to compare different pseudocone categories over the same base category and lift adjoints to the level of pseudocone categories in terms of having fibrewise adjunctions (cf.\@ Theorem \ref{Thm: Section 3: Gamma-wise adjoints lift to equivariant adjoints}). We will see that this is arguably a straightforward technical corollary by virtue of the discussion below, but the reason we have called it a theorem is due to the power of its applications. In fact, if we are being somewhat reductionist, much of this subsection can be rephrased as showing that the pseudoconification map $\PC:\Bicat(\Cscr^{\op},\fCat) \to \fCat$ induced by $F \mapsto \PC(F)$ on objects is a strict $2$-functor.

Perhaps unsurprisingly at this point, we will need to know some things about the tricategory $\Bicat$ and its various properties --- for explicit and more complete details about this tricategory see \cite[Chapter 11]{TwoDimCat}; what we describe here can largely be seen as a summary of the definitions at the beginning of \cite[Chapter 11.3]{TwoDimCat}. So far we have only used the fact that for any bicategories $\Cfrak, \Dfrak$ and for any two pseudofunctors $F, G:\Cfrak \to \Dfrak$ we have a category of transformations and modifications $\Bicat(\Cfrak, \Dfrak)(F,G)$ because $\Bicat(\Cfrak, \Dfrak)$ is a bicategory. The structure of $\Bicat$ is furnished so that the composition of $1$-morphisms in the bicategory $\Bicat(\Cfrak, \Dfrak)$ is the vertical composition of the pseudonatural transformations while the composition pseudofunctors
\[
\ast:\Bicat(\Cfrak,\Dfrak) \times \Bicat(\Dfrak, \Afrak) \to \Bicat(\Cfrak,\Afrak)
\]
describe the horizontal composition of pseudonatural transformations and pseudofunctors. The horizontal composition of pseudofunctors is the usual composition
\[
\ast(F,G) := G \circ F 
\]
but the real magic comes in at the level of $1$-and-$2$-morphisms. For any pseudonatural transformations 
\[
\begin{tikzcd}
	\Cfrak \ar[rr, bend left = 30, ""{name = U}]{}{F} \ar[rr, bend right = 30, swap, ""{name = D}]{}{G} & & \Dfrak \ar[from = U, to = D, Rightarrow, shorten <= 4pt, shorten >= 4pt]{}{\alpha}
\end{tikzcd}
\]
and
\[
\begin{tikzcd}
	\Dfrak \ar[rr, bend left = 30, ""{name = U}]{}{H} \ar[rr, bend right = 30, swap, ""{name = D}]{}{K} & & \Afrak \ar[from = U, to = D, Rightarrow, shorten >= 4pt, shorten <= 4pt]{}{\beta}
\end{tikzcd}
\]
the horizontal composition $\beta \ast \alpha:H \circ F \Rightarrow K \circ G$ is determined by
\[
\beta \ast \alpha = \ast(\alpha,\beta);
\]
note that $\beta \ast \alpha$ is defined via the whiskering equation\index[terminology]{Whiskering}
\[
\beta \ast \alpha := (K \ast \alpha) \circ (\beta \ast F).
\]
We can perform a similar construction with modifications. If we have the diagram below of bicategories, pseudofunctors, pseudonatural transformations, and modifications
\[
\begin{tikzcd}
	\Cfrak \ar[rrr, bend left = 40, ""{name = ULeft}]{}{F} \ar[rrr, bend right = 40, swap, ""{name = DLeft}]{}{E} & & & \Dfrak \ar[rrr, bend left = 40, ""{name = URight}]{}{G} \ar[rrr, bend right = 40, swap, ""{name = DRight}]{}{H} & & & \Efrak \ar[from = ULeft, to = DLeft, Rightarrow, shorten <= 4pt, shorten >= 4pt, bend left = 40, ""{name = LeftLeft}]{}{\gamma} \ar[from = ULeft, to = DLeft, Rightarrow, shorten <= 4pt, shorten >= 4pt, bend right = 40, swap, ""{name = LeftRight}]{}{\alpha} \ar[from = LeftRight, to = LeftLeft, symbol = \underset{\rho}{\Rrightarrow}] \ar[from = URight, to = DRight, Rightarrow, shorten <= 4pt, shorten >= 4pt, bend left = 40, ""{name = RightLeft}]{}{\delta} \ar[from = URight, to = DRight, Rightarrow, shorten <= 4pt, shorten >= 4pt, bend right = 40, swap, ""{name = RightRight}]{}{\beta} \ar[from = RightRight, to = RightLeft, symbol = \underset{\varphi}{\Rrightarrow}]
\end{tikzcd}
\]
then we define the modification
\[
\begin{tikzcd}
	\Cfrak \ar[rrr, bend left = 40, ""{name = Up}]{}{G \circ F} \ar[rrr, bend right = 40, swap, ""{name = Down}]{}{H \circ E} & & & \Efrak \ar[from = Up, to = Down, Rightarrow, shorten <= 4pt, shorten >= 4pt, bend left = 40, ""{name = Left}]{}{\delta \ast \gamma} \ar[from = Up, to = Down, Rightarrow, shorten <= 4pt, shorten >= 4pt, bend right = 40, swap, ""{name = Right}]{}{\beta \ast \alpha} \ar[from = Right, to = Left, symbol = \underset{\varphi \ast \rho}{\Rrightarrow}]
\end{tikzcd}
\]
as follows: for any object $X \in \Cfrak_0$ we have the induced diagram
\[
\begin{tikzcd}
	(G \circ F)(X) \ar[rr, bend left = 30, ""{name = ULeft}]{}{\beta_{FX}} \ar[rr, bend right = 30, swap, ""{name = DLeft}]{}{\delta_{FX}} & & (H \circ F)(X) \ar[rr, bend left = 30, ""{name = URight}]{}{H(\alpha_X)} \ar[rr, bend right = 30, swap, ""{name = DRight}]{}{H(\gamma_X)} & & (H \circ E)(X)  \ar[from = ULeft, to = DLeft, Rightarrow, shorten <= 4pt, shorten >= 4pt, swap]{}{\varphi_{FX}}  \ar[from = URight, to = DRight, Rightarrow, shorten <= 4pt, shorten >= 4pt]{}{H(\rho_X)} 
\end{tikzcd}
\]
and so we take the resulting $2$-morphism component $(\varphi \ast \rho)_X$ to be the horizontal composite
\[
(\varphi \ast \rho)_X := H(\rho_X) \ast \varphi_{FX}.
\]
We will see below that these give us a clean construction of the functors $\PC(F) \to \PC(E)$ for $F,E:\Cscr^{\op} \to \fCat$ and natural transformations between them. Later we will also investigate the $2$-morphisms in $\Bicat(\Cfrak,\Dfrak)$ and use those modifications to build adjunctions and natural transformations between pseudocone categories.

Let us proceed to actually apply these bicategorical observations to our pseudocone categories. Below we show that if we have a $2$-cell
\[
\begin{tikzcd}
	\Cscr^{\op} \ar[rr, bend left = 20, ""{name = U}]{}{F} \ar[rr, bend right = 20, swap, ""{name = D}]{}{E} & & \fCat \ar[from = U, to = D, Rightarrow, shorten <= 4pt, shorten >= 4pt]{}{\alpha}
\end{tikzcd}
\]
then post-composition by $\alpha$ at the $2$-cell level and post-whiskering by $\alpha$ at the $3$-cell level gives rise to a functor $\ul{\alpha}:\PC(F) \to \PC(E)$. While this is immediate by construction of the bicategory $\Bicat(\Cscr^{\op},\fCat)$, because of our interest in applications to algebraic geometry and algebraic topology and using the pseudocone formalism for practical calculations, we want to be more explicit than just having the abstract construction at hand. Consequently we also give an explicit description of how this transforms objects and morphisms as well as their transition isomorphisms; this will be especially useful when studying when the functor $\ul{\alpha}$ preserves limits and colimits.

\begin{Theorem}\label{Thm: Functor Section: Psuedonatural trans are pseudocone functors}
	Let $\Cscr$ be a category and let $F,E:\Cscr^{\op} \to \fCat$ be a pair of pseudofunctors for which there is a pseudonatural transformation $\alpha:F \to E$. Then there is an induced functor $\ul{\alpha}:\PC(F) \to \PC(E)$ induced by the vertical post-composition functor
	\[
	\ul{\alpha} = \alpha \circ (-):\Bicat(\Cscr^{\op},\fCat)(\cnst(\One),F) \to \Bicat(\Cscr^{\op},\fCat)(\cnst(\One),E).
	\]
	Explicitly, if $(A, T_A) \in \PC(F)_0$ then
	\[
	\ul{\alpha}(A) = \left\lbrace \quot{\alpha}{X}\left(\quot{A}{X}\right) \; : \; X \in \Cscr_0, \quot{A}{X} \in A\right\rbrace
	\]
	and
	\[
	T_{\ul{\alpha}A} := \left\lbrace \quot{\alpha}{X}\left(\tau_f^A\right) \circ \quot{\alpha}{f}^{-1}_{\quot{A}{Y}} \; : \; f \in \Cscr_1, f:X \to Y \right\rbrace
	\]
	where $\quot{\alpha}{f}$ is the natural isomorphism
	\[
	\quot{\alpha}{f}: \quot{\alpha}{X} \circ F(f) \xRightarrow{\cong} E(f) \circ \quot{\alpha}{Y}.
	\]
	Similarly, if $P = \lbrace \quot{\rho}{X} \; | \; X \in \Cscr_0 \rbrace$ then
	\[
	\ul{\alpha}P = \left\lbrace \quot{\alpha}{X}\left(\quot{\rho}{X}\right) \; : \; X \in \Cscr_0\right\rbrace.
	\]
\end{Theorem}
\begin{proof}
	The statement that the post-composition functor gives the induced functor $\PC(F) \to \PC(E)$ is semantics and omitted. The meat of what is to be shown here is that the objects/morphisms carry the explicit forms claimed above.
	
	Because on one hand the vertical composition 
	\[
	\ul{\alpha}A = \alpha \circ (A,T_A) = (\ul{\alpha}(A), T_{\ul{\alpha}A})
	\] 
	has object component
	\[
	\ul{\alpha}(A) = \left\lbrace \quot{\alpha}{X}\left(\quot{A}{X}\right) \; : \; X \in \Cscr_0 \right\rbrace
	\]	
	we only need to determine the transition isomorphisms $T_{\ul{\alpha}A}$. By the convention discussed in Remark \ref{Remark: Pseudocone Section: Direction of morphisms versus cones} and the definition of the vertical composite $\alpha \circ (A,T_A)$ we get that
	\[
	\tau_f^{\ul{\alpha}A} = \left(\quot{(\alpha \circ (A,T_A))}{f}\right)^{-1} = \quot{\alpha}{X}\left(\quot{(A,T_A)}{f}\right)^{-1} \circ \quot{\alpha}{f}^{-1} = \quot{\alpha}{X}(\tau_f^A) \circ \left(\quot{\alpha}{f}\right)^{-1}.
	\]
	The check that this indeed gives a pseudonatural transformation (and that we did not, in fact, make an error in defining the transition isomorphisms) follows from the fact that the pseudonaturality of $\alpha$ implies that if we have any composable pair of morphisms $X \xrightarrow{f} Y \xrightarrow{h} Z$
	\[
	E(f)\left(\quot{\alpha}{Y}\left(\tau_{g}^{A}\right)\right) = \quot{\alpha}{f}_{\quot{A}{Y}} \circ \quot{\alpha}{X}\left(F(f)\quot{\tau_g^A}{Y}\right) \circ \left(\quot{\alpha}{f}_{F(g)\quot{A}{Z}}\right)^{-1}.
	\]
	This allows us to derive that
	\begin{align*}
		\tau_{f}^{\underline{\alpha}A} \circ E(X)\tau_g^{\underline{\alpha}A} &= \quot{\alpha}{X}(\tau_f^A) \circ \left(\quot{\alpha_{\quot{A}{Y}}}{f}\right)^{-1} \circ E(X)\big(\quot{\alpha}{Y}(\tau_g^A) \circ \quot{\alpha_{\quot{A}{Z}}}{g}^{-1}\big) \\
		&=\quot{\alpha}{X}(\tau_f^A) \circ \left(\quot{\alpha_{\quot{A}{Y}}}{f}\right)^{-1} \circ E(X)\big(\quot{\alpha}{Y}(\tau_g^A)\big) \circ E(X)\big(\quot{\alpha_{\quot{A}{Z}}}{g}\big)^{-1} \\
		&= \quot{\alpha}{X}(\tau_f^A) \circ \left(\quot{\alpha_{\quot{A}{Y}}}{f}\right)^{-1} \circ \left(\quot{\alpha_{\quot{A}{Y}}}{f}\right) \circ \quot{\alpha}{X}\big(F(X)\tau_g^A\big) \\
		&\circ \left(\quot{\alpha_{F(g)\quot{A}{Z}}}{f}\right)^{-1} \circ E(X)\big(\quot{\alpha_{\quot{A}{Z}}}{g}\big)^{-1} \\
		&= \quot{\alpha}{X}(\tau_f^A) \circ \quot{\alpha}{X}(F(X)\tau_g^{A}) \circ \left(\quot{\alpha_{F(g)\quot{A}{Z}}}{f}\right)^{-1} \circ E(X)\big(\quot{\alpha_{\quot{A}{Z}}}{g}\big)^{-1} \\
		&= \quot{\alpha}{X}\big(\tau_f^A \circ F(X)\tau_g^A\big) \circ \left(\quot{\alpha_{F(g)\quot{A}{Z}}}{f}\right)^{-1} \circ E(X)\big(\quot{\alpha_{\quot{A}{Z}}}{g}\big)^{-1} \\
		&= \quot{\alpha}{X}(\tau_{g \circ f}^A \circ \quot{\phi_{f,g}}{F})\circ \left(\quot{\alpha_{F(g)\quot{A}{Z}}}{f}\right)^{-1} \circ E(X)\big(\quot{\alpha_{\quot{A}{Z}}}{g}\big)^{-1} \\
		&= \quot{\alpha}{X}(\tau_{g \circ f}^{A}) \circ \quot{\alpha}{X}(\quot{\phi_{f,g}}{F})\circ \left(\quot{\alpha_{F(g)\quot{A}{Z}}}{f}\right)^{-1} \circ E(X)\big(\quot{\alpha_{\quot{A}{Z}}}{g}\big)^{-1}
	\end{align*}
	We again use the pseudonaturality of $\alpha$ to obtain the commuting diagram
	\[
	\begin{tikzcd}
		E(X)\big(E(g)\big(\quot{\alpha}{Z}(\quot{A}{Z})\big)\big) \ar[rr]{}{E(X)\big(\quot{\alpha_{\quot{A}{Z}}}{g}\big)^{-1}} \ar[dd, swap]{}{\quot{\phi_{f,g}}{E}} & & E(X)\big(\quot{\alpha}{Y}(F(g)(\quot{A}{Z}))\big) \ar[d]{}{\big(\quot{\alpha_{F(g)\quot{A}{Z}}}{f}\big)^{-1}} \\
		& & \ar[d]{}{\quot{\alpha}{X}(\quot{\phi_{f,g}}{X})} \quot{\alpha}{X}\big(F(X)(F(g)(\quot{A}{Z}))\big) \\
		E(g \circ f)\big(\quot{\alpha}{Z}(\quot{A}{Z})\big)\ar[rr, swap]{}{\quot{\alpha_{\quot{A}{Z}}}{g \circ f}^{-1}} & & \quot{\alpha}{X}\big(F(g \circ X)(\quot{A}{Z})\big)
	\end{tikzcd}
	\]
	and substitute the induced equation to get
	\begin{align*}
		&\quot{\alpha}{X}(\tau_{g \circ f}^{A}) \circ \quot{\alpha}{X}(\quot{\phi_{f,g}}{F})\circ \left(\quot{\alpha_{F(g)\quot{A}{Z}}}{f}\right)^{-1} \circ E(X)\big(\quot{\alpha_{\quot{A}{Z}}}{g}\big)^{-1} \\&= \quot{\alpha}{X}(\tau_{g \circ f}^A) \circ \big(\quot{\alpha_{\quot{A}{Z}}}{g \circ f}\big)^{-1} \circ \quot{\phi_{f,g}}{E} = \tau_{g \circ f}^{\underline{\alpha}A} \circ \quot{\phi_{f,g}}{E}.
	\end{align*}
	This shows that $\tau_{g \circ f}^{\ul{\alpha}A} \circ \quot{\phi_{f,g}}{E} = \tau_{f}^{\ul{\alpha}A} \circ E(f)(\tau_g^{\ul{\alpha}A})$ and hence that $\ul{\alpha}(A)$ has the form claimed in the statement of the theorem.
	
	The statement regarding the description of morphisms follows from untangling the combinatorial description of the whiskering of the $2$-morphism $P$ by the $1$-morphism $\alpha$ in $\Bicat(\Cscr^{\op},\fCat)$. Explicitly we find
	\[
	\ul{\alpha}(P) = \alpha \circ P = \alpha \ast P
	\]
	and so the description
	\[
	\alpha \ast P = \ul{\alpha}P = \left\lbrace \quot{\alpha}{X}\left(\quot{\rho}{X}\right) \; : \; X \in \Cscr_0\right\rbrace
	\]
	follows. That this is indeed a modification (and hence a morphism in $\PC(E)$) is either a straightforward check using the identity
	\[
	E(f)\left(\quot{\alpha}{Y}\left(\quot{\rho}{Y}\right)\right) = \quot{\alpha}{f}_{\quot{B}{Y}} \circ \quot{\alpha}{X}\left(F(f)\quot{\rho}{Y}\right) \circ \left(\quot{\alpha}{f}_{\quot{A}{Y}}\right)^{-1}
	\]
	which allows us to deduce
	\begin{align*}
		&\tau_f^{\underline{\alpha}B} \circ E(X)\big(\quot{\alpha}{Y}(\quot{\rho}{Y})\big)\\
		&= \tau_f^{\underline{\alpha}B} \circ \quot{\alpha_{\quot{B}{Y}}}{f} \circ \quot{\alpha}{X}\big(F(X)\quot{\rho}{Y}\big)  \circ \left(\quot{\alpha_{\quot{A}{Y}}}{f}\right)^{-1} \\
		&= \quot{\alpha}{X}(\tau_f^B) \circ \left(\quot{\alpha_{\quot{B}{Y}}}{f}\right)^{-1} \circ \quot{\alpha_{\quot{B}{Y}}}{f} \circ \quot{\alpha}{X}\big(F(X)\quot{\rho}{Y}\big)  \circ \left(\quot{\alpha_{\quot{A}{Y}}}{f}\right)^{-1} \\
		&= \quot{\alpha}{X}(\tau_f^B) \circ \quot{\alpha}{X}\big(F(X)\quot{\rho}{Y}\big) \circ \left(\quot{\alpha_{\quot{A}{Y}}}{f}\right)^{-1} \\
		&= \quot{\alpha}{X}\big(\tau_f^B \circ F(X)\quot{\rho}{Y}\big) \circ \left(\quot{\alpha_{\quot{A}{Y}}}{f}\right)^{-1} \\
		&= \quot{\alpha}{X}(\quot{\rho}{X} \circ \tau_f^A) \circ \left(\quot{\alpha_{\quot{A}{Y}}}{f}\right)^{-1} = \quot{\alpha}{X}(\quot{\rho}{X}) \circ \quot{\alpha}{X}(\tau_f^A) \circ \left(\quot{\alpha_{\quot{A}{Y}}}{f}\right)^{-1} \\
		&= \quot{\alpha}{X}(\quot{\rho}{X}) \circ \tau_{f}^{\underline{\alpha}A}
	\end{align*}
	which is what was to be shown.
\end{proof}

We'll use functors of this format enough to justify naming them.
\begin{definition}\label{Defn: Functors Section: Fibre Functor}
	Let $\Cscr$ be a category and let $F, E:\Cscr^{\op} \to \fCat$ be pseudofunctors. A Change of Fibre\index[terminology]{Change of Fibre Functor} functor $G:\PC(F) \to \PC(E)$ is a functor such that $G$ is naturally isomorphic to a post-composition functor $\ul{\alpha}:\PC(F) \to \PC(E)$ for a pseudonatural transformation $\alpha:F \Rightarrow E$.
\end{definition}

We now move on to prove the first key technical result about certain cases when change of fibre functors preserve limits of fixed shape. These will let us prove that under natural\footnote{In the sense of ``occurring in (mathematical) nature.''} conditions (and in particular the conditions we will be considering in this paper) the functors induced by pseudofunctors either preserve all limits of a specific shape (cf.\@ Corollary \ref{Cor: Section 3: (Co)limit preserving pseudonaturals give (co)limit preserving functors}) or preserve all limits or colimits (cf.\@ Corollary \ref{Cor: Section 3: (Co)Continuous pseudonaturals give (co)continuous functors}). We use these in turn to establish when certain pseudocone functors are additive (cf.\@ Corollary \ref{Cor: Section 3: Additive equivraint change of fibre functors}) which is important in practice for doing homological algebra with pseudocone categories.

\begin{proposition}\label{Prop: Section 3: Change of fibre equivariant functor preserves limit of a specific shape}
	Let $F, E:\Cscr^{\op} \to \fCat$ be pseudofunctors such that for a fixed index category $I$ the hypotheses of Theorem \ref{Thm: Section 2: Equivariant Cat has lims} are satisfied for $F$ and $E$ for diagram functors $d_F:I \to\PC(F), d_E:I \to \PC(E)$, and for all $C \in \Cscr_0$
	\[
	\begin{tikzcd}
		I \ar[drr, swap]{}{\quot{d_F}{X}} \ar[r]{}{d_F} & \PC(F) \ar[r]{}{\hat{\imath}_F} & \prod\limits_{X \in \Cscr_0} F(X) \ar[d]{}{\pi_{X}^{F}} \\
		& & F(X)
	\end{tikzcd}
	\] 
	and:
	\[
	\begin{tikzcd}
		I \ar[drr, swap]{}{\quot{d_E}{X}} \ar[r]{}{d_E} & \PC(E) \ar[r]{}{\hat{\imath}_E} & \prod\limits_{X \in \Cscr_0} E(X) \ar[d]{}{\pi_{X}^{E}} \\
		& & E(X)
	\end{tikzcd}
	\]
	Assume furthermore that $\alpha:F \Rightarrow E$ is a pseudonatural transformation such that each functor $\quot{\alpha}{X}$ preserves the limit over each diagram $\quot{d_F}{X}$, i.e., the limit of $\quot{\alpha}{X} \circ \quot{d_F}{X}$ exists in $E(X)$ and is isomorphic to $\quot{\alpha}{X}$ applied to the limit of $\quot{d_F}{X}$. Then $\underline{\alpha}$ preserves the limit of the induced diagram $d_F$ in $\PC(E)$.
\end{proposition}
\begin{proof}
	Appealing to Theorem \ref{Thm: Section 2: Equivariant Cat has lims} gives that $\lim(d_F)$ exists in $\PC(F)$ and is given by the $X$-wise limit of the diagrams $\quot{d_F}{X}$; appealing to Theorem \ref{Thm: Section 2: Equivariant Cat has lims} again gives that the limit of $\underline{\alpha} \circ d_F$ exists in $\PC(E)$ and is given $X$-locally as
	\[
	\lim(\underline{\alpha} \circ d_F) = \lbrace \lim(\quot{\alpha}{X} \circ \quot{d_F}{X}) \; | \; X \in \Cscr_0 \rbrace.
	\]
	As such, to prove that $\underline{\alpha}$ preserves the limit of $d_F$, we only need to verify that there is an isomorphism
	\[
	\underline{\alpha}\left(\lim_{\substack{\longleftarrow \\ i \in I_0}} d_F(i)\right) \cong \lim_{\substack{\longleftarrow \\ i \in I_0}}\left(\underline{\alpha} \circ d_F\right)(i)
	\]
	in $\PC(E)$.
	
	To prove the existence of the isomorphism claimed above, for all $i \in I_0$ let
	\[
	A_i := d_F(i) = \lbrace \quot{A_i}{X} \; | \; X \in \Cscr_0\rbrace = \lbrace \quot{d_F(i)}{X} \; | \; X \in \Cscr_0\rbrace
	\]
	and let
	\[
	A := \lim_{\substack{\longleftarrow \\ i \in I_0}} A_i = \lim_{\substack{\longleftarrow \\ i \in I_0}} d_F(i) = \left\lbrace \lim_{\substack{\longleftarrow \\ i \in I_0}} \quot{A_i}{X} \; | \; X \in \Cscr_0\right\rbrace.
	\]
	To construct the isomorphism
	\[
	\underline{\alpha}\left(\lim_{\substack{\longleftarrow \\ i \in I_0}} d_F(i)\right) \cong \lim_{\substack{\longleftarrow \\ i \in I_0}}\left(\underline{\alpha} \circ d_F\right)(i)
	\]
	we begin by defining the witness isomorphisms induced by the limit preservation of the $\quot{\alpha}{X}$ be given by
	\[
	\quot{\theta_{\quot{A_i}{X}}}{X}:\quot{\alpha}{X}\left(\lim_{\substack{\longleftarrow \\ i \in I_0}}\quot{A_i}{X}\right) \to \lim_{\substack{\longleftarrow \\ i \in I_0}}\quot{\alpha}{X}\left(\quot{A_i}{X}\right).
	\]
	and then define the collection
	\[
	\Theta :=\lbrace \quot{\theta_{\quot{A_i}{X}}}{X} \; | \; X \in \Cscr_0\rbrace.
	\]
	Because each of the $\quot{\theta}{X}$ and $\tau_f^{-}$ maps are isomorphisms, we only need to verify that $\Theta$ is a $\PC(E)$ morphism to prove the proposition, i.e., we need to check that the diagram
	\[
	\begin{tikzcd}
		E(f)\big(\quot{\alpha}{Y}(\lim(\quot{A_i}{Y}))\big) \ar[d, swap]{}{\tau_f^{\underline{\alpha}A}} \ar[rr]{}{E(f)\quot{\theta_{\quot{A_i}{Y}}}{Y}} & & E(f)\big(\lim (\quot{\alpha}{Y}(\quot{A_i}{Y}))\big) \ar[d]{}{\tau_f^{\lim \underline{\alpha}A_i}}  \\
		\quot{\alpha}{X}\big(\lim(\quot{A_i}{X})\big) \ar[rr]{}{\quot{\theta_{\quot{A_i}{X}}}{X}} & & \lim(\quot{\alpha}{X}(\quot{A_i}{X}))
	\end{tikzcd}
	\]
	commutes.
	
	Recall first that $\quot{\theta_f}{E}$ is the witness isomorphism
	\[
	\quot{\theta_f}{E}:E(f)\big(\lim(\quot{\alpha}{Y} \circ \quot{d_F}{Y})\big) \xrightarrow{\cong} \lim\big(E(f) \circ \quot{\alpha}{Y} \circ \quot{d_F}{Y}\big),
	\]
	and similarly for $\quot{\theta_f}{F}:F(f)\big(\lim (\quot{d_F}{X})\big) \xrightarrow{\cong} \lim\big(F(f)\circ \quot{d_F}{Y}\big)$. By Theorem \ref{Thm: Functor Section: Psuedonatural trans are pseudocone functors} we have that
	\[
	\quot{\theta_{\quot{A_i}{X}}}{X} \circ \tau_f^{\underline{\alpha}A} = \quot{\theta_{\quot{A_i}{X}}}{X} \circ \quot{\alpha}{X}(\tau_f^{A}) \circ \left(\quot{\alpha_{\quot{A}{Y}}}{f}\right)^{-1}.
	\]
	
	Note that by Theorem \ref{Thm: Section 2: Equivariant Cat has lims}, we also have
	\[
	\tau_{f}^{\lim \underline{\alpha}A_i} = \lim\big(\tau_f^{\underline{\alpha}A_i}\big) \circ \quot{\theta_f}{E}
	\]
	so it follows that
	\begin{align*}
		\tau_{f}^{\lim \underline{\alpha}A} \circ E(f)\quot{\theta_{\quot{A_i}{Y}}}{Y} &= \lim\big(\quot{\alpha}{X}(\tau_f^{A_i})\big) \circ \quot{\theta_f}{E} \circ E(f)\quot{\theta_{\quot{A_i}{Y}}}{Y}.
	\end{align*}
	From the uniqueness of morphisms (isomorphisms in this case) induced by the universal property of a limit, we derive that the diagram
	\[
	\begin{tikzcd}
		E(f)\big(\quot{\alpha}{Y}(\lim (\quot{A_i}{Y}))\big) \ar[d, swap]{}{E(f)\quot{\theta_{\quot{A_i}{Y}}}{Y}}  \ar[rr]{}{\left(\quot{\alpha_{\quot{A}{Y}}}{f}\right)^{-1}} & & \quot{\alpha}{X}\big(F(f)(\lim(\quot{A_i}{Y}))\big) \ar[d]{}{\quot{\alpha}{X}\big(\quot{\theta_f}{F}\big)} \\
		E(f)\big(\lim(\quot{\alpha}{Y}(\quot{A_i}{Y}))\big)  & & \quot{\alpha}{X}\big(\lim (F(f)(\quot{A_i}{Y}))\big) \ar[d]{}{\quot{\theta_{F(f)\quot{A_i}{Y}}}{X}}\\
		\lim E(f)\big(\quot{\alpha}{Y}(\quot{A_i}{Y})\big)     \ar[u]{}{\left(\quot{\theta_f}{E}\right)^{-1}} & & \lim \left(\quot{\alpha}{X}\big(F(f)(\quot{A_i}{Y})\big)\right) \ar[ll]{}{\lim\left(\quot{\alpha_{\quot{A_i}{Y}}}{f}\right)}
	\end{tikzcd}
	\]
	commutes in $E(X)$. Substituting the induced identity into the equation above gives that
	\begin{align*}
		&\tau_{f}^{\lim \underline{\alpha}A} \circ E(f)\quot{\theta_{\quot{A_i}{Y}}}{Y} \\
		&= \lim\big(\tau_f^{\underline{\alpha}A_i}\big) \circ \quot{\theta_f}{E} \circ E(f)\quot{\theta_{\quot{A}{Y}}}{Y} \\
		&=\lim(\tau_f^{\underline{\alpha}A_i}) \circ  \quot{\theta_f}{E} \circ \left(\quot{\theta_f}{E}\right)^{-1} \circ \lim\left(\quot{\alpha_{\quot{A_i}{Y}}}{f} \right) \circ \quot{\theta_{F(f)\quot{A_i}{Y}}}{X} \\
		&\circ \quot{\alpha}{X}\big(\quot{\theta_f}{F}\big) \circ \left(\quot{\alpha_{\quot{A}{Y}}}{f}\right)^{-1} \\
		&= \lim(\tau_f^{\underline{\alpha}A_i}) \circ \lim\left(\quot{\alpha_{\quot{A_i}{Y}}}{f} \right) \circ \quot{\theta_{F(f)\quot{A_i}{Y}}}{X} \circ \quot{\alpha}{X}\big(\quot{\theta_f}{F}\big) \circ \left(\quot{\alpha_{\quot{A}{Y}}}{f}\right)^{-1} \\
		&= \lim\left(\tau_f^{\underline{\alpha}A_i} \circ \quot{\alpha_{\quot{A_i}{Y}}}{f} \right) \circ \quot{\theta_{F(f)\quot{A_i}{Y}}}{X} \circ \quot{\alpha}{X}\big(\quot{\theta_f}{F}\big) \circ \left(\quot{\alpha_{\quot{A}{Y}}}{f}\right)^{-1} \\
		&= \lim\left(\quot{\alpha}{X}(\tau_f^{A_i}) \circ \quot{\alpha_{\quot{A_i}{Y}}}{f} \circ \left(\alpha_{\quot{A_i}{Y}}{f}\right)^{-1}\right) \circ \quot{\theta_{F(f)\quot{A_i}{Y}}}{X} \\
		&\circ \quot{\alpha}{X}\big(\quot{\theta_f}{F}\big) \circ \left(\quot{\alpha_{\quot{A}{Y}}}{f}\right)^{-1} \\
		&= \lim\left(\quot{\alpha}{X}(\tau_f^{A_i})\right) \circ \quot{\theta_{F(f)\quot{A_i}{Y}}}{X} \circ \quot{\alpha}{X}\big(\quot{\theta_f}{F}\big) \circ \left(\quot{\alpha_{\quot{A}{Y}}}{f}\right)^{-1}.
	\end{align*}
	We now use the uniqueness of limit cones to give the commutativity of the diagram
	\[
	\begin{tikzcd}
		\quot{\alpha}{X}\big(\lim(F(f)(\quot{A_i}{X}))\big) \ar[rr]{}{\quot{\theta_{F(f)\quot{A_i}{Y}}}{X}} \ar[d, swap]{}{\quot{\alpha}{X}(\lim\tau_f^{A_i})} & &\lim \quot{\alpha}{X}\big(F(f)(\quot{A_i}{Y})\big) \ar[d]{}{\lim \quot{\alpha}{X}(\tau_f^{A_i})} \\
		\quot{\alpha}{X}\big(\lim (\quot{A_i}{X})\big) \ar[rr]{}{\quot{\theta_{\quot{A_i}{X}}}{X}} & & \lim\big( \quot{\alpha}{X}\big(\quot{A_i}{X}\big)\big)
	\end{tikzcd}
	\]
	in $E({X})$. This allows us to further derive that
	\begin{align*}
		&\lim\left(\quot{\alpha}{X}(\tau_f^{A_i})\right) \circ \quot{\theta_{F(f)\quot{A_i}{Y}}}{X} \circ \quot{\alpha}{X}\big(\quot{\theta_f}{F}\big) \circ \left(\quot{\alpha_{\quot{A}{Y}}}{f}\right)^{-1}\\
		&=\quot{\theta_{\quot{A_i}{X}}}{X} \circ \quot{\alpha}{X}\left(\lim \tau_f^{A_i}\right) \circ \quot{\alpha}{X}\big(\quot{\theta_f}{F}\big) \circ \left(\quot{\alpha_{\quot{A}{Y}}}{f}\right)^{-1} \\
		&= \quot{\theta_{\quot{A_i}{X}}}{X} \circ \quot{\alpha}{X}\left(\lim\big(\tau_f^{A_i}\big) \circ \quot{\theta_f}{F}\right) \circ \left(\quot{\alpha_{\quot{A}{Y}}}{f}\right)^{-1} \\
		&= \quot{\theta_{\quot{A_i}{X}}}{X} \circ \quot{\alpha}{X}\left(\tau_f^A\right) \circ \left(\quot{\alpha_{\quot{A}{Y}}}{f}\right)^{-1} \\
		&= \quot{\theta_{\quot{A_i}{X}}}{X} \circ \tau_f^{\underline{\alpha}A}.
	\end{align*}
	Thus it follows that $\tau_{f}^{\lim \underline{\alpha}A} \circ E(f)\quot{\theta_{\quot{A_i}{Y}}}{Y} = \quot{\theta_{\quot{A_i}{X}}}{X} \circ \tau_f^{\underline{\alpha}A}$, which in turn proves that $\Theta$ is a morphism in $\PC(E)$.
\end{proof}

\begin{corollary}\label{Cor: Section 3: (Co)limit preserving pseudonaturals give (co)limit preserving functors}
	Assume that $F, E:\Cscr^{\op} \to \fCat$ are pseudofunctors such that each fibre category $F(X)$ and $E(X)$ admit all limits of shape $I$ and that each fibre functor $F(f)$ and $E(f)$ preserves these limits. Assume furthermore that $\alpha:F \Rightarrow E$ is a pseudonatural transformation whose fibre-transition functors preserve limits of shape $I$. Then $\underline{\alpha}:\PC(F) \to \PC(E)$ preserves all limits of shape $I$ in $\PC(F)$. Dually, if $F, E:\Cscr^{\op} \to \fCat$ are pseudofunctors such that each fibre category $F(X)$ and $E(X)$ admit all colimits of shape $I$ and that each fibre functor $F(f)$ and $E(f)$ preserves these colimits. Assume furthermore that $\alpha$ is a pseudonatural transformation whose fibre-transition functors preserve colimits of shape $I$. Then $\underline{\alpha}:\PC(F) \to \PC(E)$ preserves all colimits of shape $I$ in $\PC(F)$.
\end{corollary}
\begin{proof}
	Simply use Proposition \ref{Prop: Section 3: Change of fibre equivariant functor preserves limit of a specific shape} for every limit to every diagram over the index category $I$ to give the result. To prove this for colimits, dualize Proposition \ref{Prop: Section 3: Change of fibre equivariant functor preserves limit of a specific shape} and proceed mutatis mutandis.
\end{proof}
\begin{corollary}\label{Cor: Section 3: (Co)Continuous pseudonaturals give (co)continuous functors}
	Assume that $F, E:\Cscr^{\op} \to \fCat$ are pseudofunctors such that each fibre category $F(X)$ and $E(X)$ is complete and for which each fibre functor $F(f)$ and $E(f)$ is continuous, i.e., preserves limits. If $\alpha:F \Rightarrow E$ is a pseudonatural transformation whose fibre-transition functors are all continuous, then $\underline{\alpha}:\PC(F) \to \PC(E)$ is continuous. Dually, if $F, E:\Cscr^{\op} \to \fCat$ are pseudofunctors such that each fibre category $F(X)$ and $E(X)$ is cocomplete and for which each fibre functor $F(f)$ and $E(f)$ is cocontinuous and if $\alpha:F \Rightarrow E$ is a pseudonatural transformation whose fibre-transition functors are all cocontinuous, then $\underline{\alpha}:\PC(F) \to \PC(E)$ is cocontinuous
\end{corollary}
\begin{corollary}\label{Cor: Section 3: Additive equivraint change of fibre functors}
	Assume that $F$ and $E$ are additive pre-equivariant pseudofunctors and that each fibre-transition functor $\quot{\alpha}{X}:F(X) \to E(X)$ is an additive functor. Then the functor $\underline{\alpha}:\PC(F) \to \PC(E)$ is an additive functor of additive categories.
\end{corollary}

We proceed with a lonely proposition that shows that monoidal pseudonatural transformations give rise to monoidal equivariant functors. However, we first recall what it means to have a natural transformation be monoidal.
\begin{definition}\index[terminology]{Monoidal Natural Transformation}
	Let $(\Cscr,\otimes,I)$ and $(\Dscr, \boxtimes,J)$ be monoidal categories with monoidal functors $F,G:\Cscr \to \Dscr$. A natural transformation $\alpha:F \Rightarrow G$ is {monoidal} if and only if the following hold: 
	\begin{itemize}
		\item Given the unit preservation isomorphisms $\gamma:F(I) \xrightarrow{\cong} J$ and $\widetilde{\gamma}:F(I) \to J$, the diagram
		\[
		\begin{tikzcd}
			F(I) \ar[dr, swap]{}{\gamma} \ar[r]{}{\alpha_I} & G(I) \ar[d]{}{\widetilde{\gamma}} \\
			& J
		\end{tikzcd}
		\]
		commutes in $\Dscr$;
		\item For any objects $A,B \in \Cscr_0$, given the tensor preservation isomorphisms
		\[
		\theta_{F}^{A,B}:F(A \otimes B) \xrightarrow{\cong} FA \boxtimes FB
		\]
		and
		\[
		\theta_{G}^{A,B}:G(A \otimes B) \xrightarrow{\cong} GA \boxtimes GB,
		\]
		the diagram
		\[
		\begin{tikzcd}
			F(A \otimes B) \ar[d,swap]{}{\alpha_{A \otimes B}} \ar[r]{}{\theta_F^{A,B}} & FA \boxtimes FB \ar[d]{}{\alpha_A \boxtimes \alpha_B} \\
			G(A \otimes B) \ar[r, swap]{}{\theta_{G}^{A,B}} & GA \boxtimes GB
		\end{tikzcd}
		\]
		commutes in $\Dscr$.
	\end{itemize} 
\end{definition}
%

\begin{proposition}\label{Prop: Section 3: Monoidal pseudonats give rise to monoidal equivariant functors}
	Let $F,E:\Cscr^{\op} \to \fCat$ be monoidal pseudofunctors over $X$ and let $\alpha:F \Rightarrow E$ be a pseudonatural transformation such that each functor $\quot{\alpha}{X}$ is monoidal and each natural transformation $\quot{\alpha}{f}$ is a monoidal transformation. Then the functor $\underline{\alpha}:\PC(F) \to \PC(E)$ is monoidal.
\end{proposition}
\begin{proof}
	For any $X \in \Cscr_0$ we let $\quot{I_F}{X}, \quot{I_E}{X}$ be the tensorial units in the categories $F(X)$ and $E(X)$, respectively; write
	\[
	\quot{\rho}{X}:\quot{\alpha}{X}(\quot{I_F}{X}) \xrightarrow{\cong} \quot{I_E}{X}
	\]
	for the unit preservation isomorphisms induced by the monoidal functors $\quot{\alpha}{X}$;  let $f \in \Sf(G)_1$ be arbitrary; and as in Theorem \ref{Theorem: Section 2: Monoidal preequivariant pseudofunctor gives monoidal equivariant cat}, let
	\[
	\sigma_f:F(f)(\quot{I_F}{Y}) \to \quot{I_F}{X}
	\]
	and
	\[
	\sigma_f^{\prime}:E(f)(\quot{I_E}{Y}) \to \quot{I_E}{X}
	\]
	be the unit preservation isomorphisms induced by the monoidal functors $F(f)$ and $E(f)$. Similarly, let
	\[
	\quot{\theta_f}{F}^{A,B}:F(f)\left(\quot{A}{Y} \ds{Y}{F}{\otimes} \quot{B}{Y} \right) \xrightarrow{\cong} \left(F(f)(\quot{A}{Y})\right) \ds{X}{F}{\otimes} \left(F(f)(\quot{B}{Y})\right)
	\]
	and
	\[
	\quot{\theta_f}{E}^{A,B}:\quot{\theta_f}{E}^{A,B}:E(f)\left(\quot{A}{Y} \ds{Y}{E}{\otimes} \quot{B}{Y} \right) \xrightarrow{\cong} \left(E(f)(\quot{A}{Y})\right) \ds{X}{E}{\otimes} \left(E(f)(\quot{B}{Y})\right)
	\]
	be the tensor preservation isomorphisms induced by the monoidal functors $F(f)$ and $E(f)$, respectively, where we write $\otimes_E^X$ and $\otimes_F^X$ as the tensor functors of $E(X)$ and $F(X)$, respectively. Finally let
	\[
	\quot{\theta_{X}}{\alpha}^{A,B}:\quot{\alpha}{X}\left(\quot{A}{X} \ds{X}{F}{\otimes} \quot{B}{X}\right) \xrightarrow{\cong} \quot{\alpha}{X}\left(\quot{A}{X}\right) \ds{X}{E}{\otimes} \quot{\alpha}{X}\left(\quot{B}{X}\right)
	\]
	be the tensor preservation isomorphism induced by $\quot{\alpha}{X}$.
	
	We first prove that $\underline{\alpha}$ preserves the unit object. Let $I_F := \lbrace \quot{I_F}{X} \; | \; X \in \Cscr_0 \rbrace$ be the object set of the unit object $(I_F, T_{I_F})$ in $F_G(X)$; note that by Theorem \ref{Theorem: Section 2: Monoidal preequivariant pseudofunctor gives monoidal equivariant cat}, $T_{I_F} = \lbrace \sigma_f \; | \; f \in \Sf(G)_1 \rbrace$. Similarly, let $I_E := \lbrace \quot{I_E}{X} \; | \; X \in \Cscr_0 \rbrace$ be the object set of the unit object $(I_E,T_{I_E})$ in $E_G(X)$; as before, $T_{I_E} = \lbrace \sigma_{f}^{\prime} \; | \; f \in \Sf(G)_1 \rbrace$. Define the collection of morphisms $P$ by
	\[
	P := \lbrace \quot{\rho}{X}:\quot{\alpha}{X}(\quot{I_F}{X}) \to \quot{I_E}{X} \; | \; X \in \Cscr_0 \rbrace.
	\]
	Because each of the $\quot{\rho}{X}$ is an isomorphism, if we can prove that $P$ is a $E_G(X)$-morphism, we will have proved that $\underline{\alpha}$ preserves unit objects up to isomorphism. To this end fix an $f \in \Sf(G)_1$ and write $f:X \to Y$. We must show that the diagram
	\[
	\xymatrix{
		E(f)\big(\quot{\alpha}{Y}(\quot{I_F}{Y})\big) \ar[rr]^-{E(f\quot{\rho}{Y})} \ar[d]_{\tau_f^{\underline{\alpha}I_F}} & & E(f)(\quot{I_E}{Y}) \ar[d]^{\tau_f^{I_E}} & \\
		\quot{\alpha}{X}(\quot{I_F}{X}) \ar[rr]_-{\quot{\rho}{X}} & & \quot{I_E}{X}
	}
	\]
	commutes. For this observe that the unit preservation isomorphism 
	\[
	\gamma: \quot{\alpha}{X}\big(F(f)(\quot{I_F}{Y})\big) \to \quot{I_E}{X}
	\] 
	has the form
	\[
	\gamma = \quot{\rho}{X} \circ \quot{\alpha}{X}(\sigma_f) = \quot{\rho}{X} \circ \quot{\alpha}{X}(\tau_f^{I_F})
	\]
	while the unit preservation isomorphism $\widetilde{\gamma}:E(f)\big(\quot{\alpha}{Y}(\quot{I_F}{Y})\big) \to \quot{I_E}{X}$ takes the form
	\[
	\widetilde{\gamma} = \sigma_f^{\prime} \circ E(f)\quot{\rho}{Y} = \tau_f^{I_E} \circ E(f)\quot{\rho}{Y}.
	\]
	Using that $\quot{\alpha}{f}$ is a monoidal natural transformation gives that the diagram
	\[
	\xymatrix{
		\quot{\alpha}{X}\big(F(f)(\quot{I_F}{Y})\big) \ar[rr]^-{\gamma} \ar[d]_{\quot{\alpha_{\quot{I_F}{Y}}}{f}} & & \quot{I_E}{X} \\
		E(f)\big(\quot{\alpha}{Y}(\quot{I_F}{Y})\big) \ar[urr]_{\widetilde{\gamma}}
	}
	\]
	commutes. This gives, however, that
	\begin{align*}
		\quot{\rho}{X} \circ \tau_{f}^{\underline{\alpha}I_F} &= \quot{\rho}{X} \circ \quot{\alpha}{X}(\tau_f^{I_F}) \circ \left(\quot{\alpha_{\quot{I_F}{Y}}}{f}\right)^{-1}= \gamma \circ \left(\quot{\alpha_{\quot{I_F}{Y}}}{f}\right)^{-1} \\
		&= \widetilde{\gamma} \circ \quot{\alpha_{\quot{I_F}{Y}}}{f} \circ \left(\quot{\alpha_{\quot{I_F}{Y}}}{f}\right)^{-1} = \widetilde{\gamma} 
		= \tau_f^{I_E} \circ E(f)\quot{\rho}{Y}.
	\end{align*}
	Thus $P$ is a morphism in $E_G(X)$ and hence $\underline{\alpha}$ preserves the unit object.
	We now prove that $\underline{\alpha}$ preserves the tensor functors up to natural isomorphism. For this fix objects $A, B \in F_G(X)_0$ and consider the isomorphisms, for all $X \in \Cscr_0$,
	\[
	\quot{\theta_{X}}{\alpha}^{A,B}:\quot{\alpha}{X}\left(\quot{A}{X} \ds{X}{F}{\otimes} \quot{B}{X}\right) \xrightarrow{\cong} \quot{\alpha}{X}\left(\quot{A}{X}\right) \ds{X}{E}{\otimes} \quot{\alpha}{X}\left(\quot{B}{X}\right).
	\]
	Define the collection $\Theta_{A,B}:\underline{\alpha}(A \otimes_F B) \to \underline{\alpha}A  \otimes_E \underline{\alpha}B$ by
	\[
	\Theta_{A,B}:= \left\lbrace \quot{\theta_{X}^{A,B}}{\alpha} \; : \; X \in \Cscr_0 \right\rbrace.
	\]
	Note that because each of the $\quot{\theta_{X}}{\alpha}$ are natural isomorphisms, if we can prove that $\Theta_{A,B}$ is a $E_G(X)$-morphism, then it will automatically induce a natural isomorphism $\underline{\alpha} \circ \big((-)\otimes_F(-)\big) \Rightarrow \underline{\alpha}(-) \otimes_E \underline{\alpha}(-)$. Thus to complete the proof of the proposition it suffices to prove that the diagram
	\[
	\xymatrix{
		E(f)\left(\quot{\alpha}{Y}\left(\quot{A}{Y} \ds{Y}{F}{\otimes} \quot{B}{Y} \right)\right) \ar[rr]^-{E(f)\quot{\theta_{Y}}{\alpha}^{A,B}} \ar[d]_{\tau_{f}^{\underline{\alpha}(A \otimes_F B)}} & & E(f)\left(\quot{\alpha}{Y}(\quot{A}{Y}) \ds{Y}{E}{\otimes} \quot{\alpha}{Y}(\quot{B}{Y})\right) \ar[d]^{\tau_{f}^{\underline{\alpha}A \otimes_E \underline{\alpha}B}} \\
		\quot{\alpha}{X}\left(\quot{A}{X} \ds{X}{F}{\otimes} \quot{B}{X}\right) \ar[rr]_-{\quot{\theta_{X}}{\alpha}^{A,B}} & & \quot{\alpha}{X}(\quot{A}{X}) \ds{X}{E}{\otimes} \quot{\alpha}{X}(\quot{B}{X})
	}
	\]
	commutes. For this we first compute that by Theorems \ref{Theorem: Section 2: Monoidal preequivariant pseudofunctor gives monoidal equivariant cat} and \ref{Thm: Functor Section: Psuedonatural trans are pseudocone functors} we have one hand
	\begin{align*}
		\tau_f^{\underline{\alpha}(A \otimes_F B)} &= \quot{\alpha}{X}(\tau_f^{A \otimes_F B}) \circ \left(\quot{\alpha_{\quot{A}{Y} \otimes \quot{B}{Y}}}{f}\right)^{-1} \\
		&= \quot{\alpha}{X}\left(\tau_f^A \ds{X}{F}{\otimes} \tau_f^B \right) \circ \quot{\alpha}{X}(\quot{\theta_{f}}{F}^{A,B}) \circ \left(\quot{\alpha_{\quot{A}{Y} \otimes \quot{B}{Y}}}{f}\right)^{-1},
	\end{align*}
	while on the other hand
	\begin{align*}
		\tau_f^{\underline{\alpha}A \otimes_E \underline{\alpha}B} &=  \left(\tau_f^{\underline{\alpha}A} \ds{X}{E}{\otimes} \tau_f^{\underline{\alpha}B}\right) \circ \quot{\theta_f}{E}^{\underline{\alpha}A,\underline{\alpha}B} \\ 
		&= \left(\left(\quot{\alpha}{X}(\tau_f^A) \circ \left(\quot{\alpha_{\quot{A}{Y}}}{f}\right)^{-1}\right) \ds{X}{E}{\otimes} \left(\quot{\alpha}{X}(\tau_f^B) \circ \left(\quot{\alpha_{\quot{B}{Y}}}{f}\right)^{-1}\right)\right) \circ \quot{\theta_f}{E}^{\underline{\alpha}A,\underline{\alpha}B} \\
		&= \left(\quot{\alpha}{X}(\tau_f^A) \ds{X}{E}{\otimes} \quot{\alpha}{X}(\tau_f^{B}) \right) \circ \left(\quot{\alpha_{\quot{A}{Y}}}{f} \ds{X}{E}{\otimes} \quot{\alpha_{\quot{B}{Y}}}{f}\right)^{-1} \circ \quot{\theta_f}{E}^{\underline{\alpha}A,\underline{\alpha}B}.
	\end{align*}
	Now observe that the tensor preservation isomorphism
	\[
	\quot{\alpha}{X}\left(F(f)\left(\quot{A}{Y} \ds{Y}{F}{\otimes} \quot{B}{Y}\right)\right) \xrightarrow[\quot{\theta_f}{\alpha \ast F}^{A,B}]{\cong} \left(\quot{\alpha}{X}\big(F(f)(\quot{A}{Y})\big) \right) \ds{X}{E}{\otimes} \left(\quot{\alpha}{X}\big(F(f)(\quot{B}{Y})\big)\right)
	\]
	has the form
	\[
	\quot{\theta_f}{\alpha \ast F}^{F(f)A,F(f)B} = \quot{\theta_{X}}{\alpha}^{A,B} \circ \quot{\alpha}{X}\big(\quot{\theta_f}{F}^{A,B}\big)
	\]
	while the tensor preservation isomorphism
	\[
	E(f)\left(\quot{\alpha}{Y}\left(\quot{A}{Y} \ds{Y}{F}{\otimes} \quot{B}{Y}\right)\right) \xrightarrow[\quot{\theta_{f}}{E \ast \alpha}^{A,B}]{\cong} \left(E(f)\big(\quot{\alpha}{Y}(\quot{A}{Y})\big)\right) \ds{X}{E}{\otimes} \left(E(f)\big(\quot{\alpha}{Y}(\quot{B}{Y})\big)\right)
	\]
	has the form
	\[
	\quot{\theta_{f}}{E \ast \alpha}^{A,B} = \quot{\theta_f}{E}^{\underline{\alpha}A,\underline{\alpha}B} \circ E(f)\big(\quot{\theta_{Y}}{\alpha}^{A,B}\big),
	\]
	where in both cases $\alpha \ast F$ and $E \ast \alpha$ denote the whiskering of $\alpha$ by $F$ and $E$  by $\alpha$, respectively. Note also that since $\quot{\alpha}{f}$ is a monoidal natural isomorphism, so is $\quot{\alpha}{f}^{-1}$. The monoidal transformation $\quot{\alpha}{f}^{-1}$ also satisfies the identity
	\[
	\left(\quot{\alpha_{\quot{A}{Y}}}{f} \ds{X}{E}{\otimes} \quot{\alpha_{\quot{B}{Y}}}{f}\right)^{-1} \circ \quot{\theta_f}{E \ast \alpha}^{A,B} = \quot{\theta_f}{\alpha \ast F}^{A,B} \circ \left(\quot{\alpha_{\quot{A}{Y} \otimes \quot{B}{Y}}}{f}\right)^{-1}.
	\]
	Using these computations and identities in various ways we compute that
	\begin{align*}
		&\tau_f^{\underline{\alpha}A \otimes_E \underline{\alpha}B} \circ E(f)\quot{\theta_{Y}}{\alpha}^{A,B} \\
		&= \left(\tau_f^{\underline{\alpha}A} \ds{X}{E}{\otimes} \tau_f^{\underline{\alpha}B}\right) \circ \quot{\theta_f}{E}^{\underline{\alpha}A,\underline{\alpha}B} \circ E(f)\quot{\theta_{Y}}{\alpha}^{A,B} = \left(\tau_f^{\underline{\alpha}A} \ds{X}{E}{\otimes} \tau_f^{\underline{\alpha}B}\right) \circ \quot{\theta_f}{E \ast \alpha}^{A,B} \\
		&= \left(\quot{\alpha}{X}(\tau_f^A) \ds{X}{E}{\otimes} \quot{\alpha}{X}(\tau_f^{B}) \right) \circ \left(\quot{\alpha_{\quot{A}{Y}}}{f} \ds{X}{E}{\otimes} \quot{\alpha_{\quot{B}{Y}}}{f}\right)^{-1} \circ \quot{\theta_f}{E \ast \alpha}^{A,B} \\
		&= \left(\quot{\alpha}{X}(\tau_f^A) \ds{X}{E}{\otimes} \quot{\alpha}{X}(\tau_f^{B}) \right) \circ \quot{\theta_f}{\alpha \ast F}^{A,B} \circ \left(\quot{\alpha_{\quot{A}{Y} \otimes \quot{B}{Y}}}{f}\right)^{-1} \\
		&= \left(\quot{\alpha}{X}(\tau_f^A) \ds{X}{E}{\otimes} \quot{\alpha}{X}(\tau_f^{B}) \right)\circ \quot{\theta_{X}}{\alpha}^{F(f)A,F(f)B} \circ \quot{\alpha}{X}\big(\quot{\theta_f}{F}^{A,B}\big) \circ \left(\quot{\alpha_{\quot{A}{Y} \otimes \quot{B}{Y}}}{f}\right)^{-1} \\
		&= \quot{\theta_{X}}{\alpha}^{A,B} \circ \quot{\alpha}{X}\left(\tau_f^A \ds{X}{F}{\otimes} \tau_f^B\right) \circ \quot{\alpha}{X}\big(\quot{\theta_f}{F}^{A,B}\big) \circ \left(\quot{\alpha_{\quot{A}{Y} \otimes \quot{B}{Y}}}{f}\right)^{-1} \\
		&= \quot{\theta_{X}}{\alpha}^{A,B} \circ \tau_{f}^{\underline{\alpha}(A \otimes_F B)}.
	\end{align*}
	This gives that $\tau_f^{\underline{\alpha}A \otimes_E \underline{\alpha}B} \circ E(f)\quot{\theta_{Y}}{\alpha}^{A,B} = \quot{\theta_{X}}{\alpha}^{A,B} \circ \tau_{f}^{\underline{\alpha}(A \otimes_F B)}$ which proves that $\Theta_{A,B}$ is a morphism in $\PC(E)$.
\end{proof}
\begin{corollary}\label{Cor: Pseudocone Functors: Braided monoidal pseudonats are braided monoidal functors}
	Let $F,E:\Cscr^{\op} \to \fCat$ be braided monoidal pseudofunctors over $X$ and let $\alpha:F \Rightarrow E$ be a pseudonatural transformation such that each functor $\quot{\alpha}{X}$ is braided monoidal and each natural transformation $\quot{\alpha}{f}$ is a monoidal transformation. Then the functor $\underline{\alpha}:\PC(F) \to \PC(E)$ is braided monoidal.
\end{corollary}
\begin{proof}
	By virtue of Proposition \ref{Prop: Section 3: Monoidal pseudonats give rise to monoidal equivariant functors} we only need to show that $\ul{\alpha}$ is braided. For this fix objects $A, B \in \PC(F)_0$ and note that the assumptions on everything in sight imply that for all $X \in \Cscr_0$ the diagram
	\[
	\begin{tikzcd}
		\quot{\alpha}{X}\left(\quot{A}{X} \us{F(X)}{\otimes} \quot{B}{X}\right) \ar[rr]{}{\quot{\theta_{A,B}}{X}} \ar[d, swap]{}{\quot{\alpha}{X}\left(\quot{\beta_{A,B}}{F(X)}\right)} & & \quot{\alpha}{X}\left(\quot{A}{X}\right) \us{E(X)}{\otimes} \quot{\alpha}{X}\left(\quot{B}{X}\right) \ar[d]{}{\quot{\beta_{\ul{\alpha}(A),\ul{\alpha}{B}}}{E(X)}} \\
		\quot{\alpha}{X}\left(\quot{B}{X} \us{F(X)}{\otimes} \quot{A}{X}\right) \ar[rr, swap]{}{\quot{\theta_{B,A}}{X}} & & \quot{\alpha}{X}\left(\quot{B}{X}\right) \us{E(X)}{\otimes} \quot{\alpha}{X}\left(\quot{A}{X}\right)
	\end{tikzcd}
	\]
	commutes. However, this in turn implies that the diagram
	\[
	\begin{tikzcd}
		\ul{\alpha}\left(A \us{\PC(F)}{\otimes} B\right) \ar[rr]{}{\Theta_{A,B}} \ar[d, swap]{}{\ul{\alpha}(\quot{\beta_{A,B}}{F}} & & \ul{\alpha}(A) \us{\PC(E)}{\otimes} \ul{\alpha}(B) \ar[d]{}{\quot{\beta_{\ul{\alpha}(A),\ul{\alpha}(B)}}{E}} \\
		\ul{\alpha}\left(B \us{\PC(F)}{\otimes} A\right) \ar[rr, swap]{}{\Theta_{B,A}} & & \ul{\alpha}(B) \us{\PC(E)}{\otimes} \ul{\alpha}\left(A \right)
	\end{tikzcd}
	\]
	commutes in $\PC(E)$ and hence shows that $\ul{\alpha}$ is braided.
\end{proof}

We now show that if we have a modification between pseudonatural transformations
\[
\eta:\alpha\Rrightarrow\beta:F \Rightarrow E:\Cscr^{\op}\to \Cfrak
\]
then there is an induced natural transformation
\[
\ul{\eta}:\ul{\alpha} \Rightarrow \ul{\beta}:\PC(F) \to \PC(E).
\]
For a reminder on the information that goes into a modification of this form, see the definition below.
\begin{definition}\index[terminology]{Modification}
	Let $F,E:\Cscr^{\op} \to \fCat$ be pseudofunctors for some $1$-category $\Cscr$ and let $\alpha,\beta:F \Rightarrow E$ be pseudonatural transformations. A modification $\eta:\alpha \Rrightarrow \beta$ is an assignment $\eta:\Cscr_0 \to \fCat_2$ where each natural transformation $\eta_A$ fits into a $2$-cell
	\[
	\begin{tikzcd}
		F(A) \ar[r,bend left = 30, ""{name = U}]{}{\alpha_A} \ar[r,swap, bend right = 30,""{name = L}]{}{\beta_A} & E(A) \ar[from = U, to = L, Rightarrow,shorten >= 5pt, shorten <= 5pt]{}{\eta_A}
	\end{tikzcd}
	\]
	in $\fCat$ such that for all morphisms $g:A \to B$ in $\Cscr$, the diagram of functors/natural transformations
	\[
	\xymatrix{
		\alpha_A \circ F(g) \ar[rr]^-{\eta_A \ast \iota_{F(g)}} \ar[d]_{\alpha_g} & & \beta_A \circ F(g) \ar[d]^{\beta_g} \\
		E(g) \circ \alpha_B \ar[rr]_-{\iota_{E(g)}\ast\eta_B} & & E(g) \circ \beta_B
	}
	\] 
	commutes (where $\ast$ denotes horizontal composition of natural transformations).
\end{definition}
\begin{lemma}\label{Lemma: Modifications give equivariant natural transformations}
	Let $F,E:\Cscr^{\op} \to \fCat$ be pseudofunctors and let $\alpha,\beta:F \Rightarrow E$ be pseudonatural transformations. If $\eta:\alpha \Rrightarrow \beta$ is a modification, then there exists a natural transformation $\underline{\eta}$ fitting into the $2$-cell
	\[
	\begin{tikzcd}
		\PC(F) \ar[r,bend left = 30, ""{name = U}]{}{\underline{\alpha}} \ar[r,swap, bend right = 30,""{name = L}]{}{\underline{\beta}} & \PC(E) \ar[from = U, to = L, Rightarrow,shorten >= 5pt, shorten <= 5pt]{}{\underline{\eta}}
	\end{tikzcd}
	\]
	given on objects $C = \lbrace \quot{C}{} \; | \; X \in \Cscr_0\rbrace$ of $\PC(F)$ by
	\[
	\underline{\eta}_C := \lbrace \quot{\eta}{X}_{\quot{C}{X}}:\quot{\alpha}{X}(\quot{C}{X}) \to \quot{\beta}{X}(\quot{C}{X}) \; | \; X \in \Cscr_0\rbrace.
	\]
\end{lemma}
\begin{proof}
	Before we prove the naturality of $\underline{\eta}$, we must prove that for all objects $A$ of $\PC(F)$, the membership $\underline{\eta}_A \in \PC(E)(\underline{\alpha}A, \underline{\beta}A)$. That is, for all $f \in \Sf(G)_1$ (implicitly with $f:X \to Y$) we have that
	\[
	\tau_f^{\underline{\beta}A} \circ E(f)\left(\quot{\eta_{\quot{A}{Y}}}{Y}\right) = \quot{\eta_{\quot{A}{X}}}{X} \circ \tau_f^{\underline{\alpha}A}.
	\]
	To this end recall that by Theorem \ref{Thm: Functor Section: Psuedonatural trans are pseudocone functors}
	\[
	\tau_{f}^{\underline{\alpha}A} = \quot{\alpha}{X}(\tau_f^A) \circ \left(\quot{\alpha_{\quot{A}{Y}}}{f}\right)^{-1}, \qquad \tau_f^{\underline{\beta}A} = \quot{\beta}{X}(\tau_f^A) \circ \left(\quot{\beta_{\quot{A}{Y}}}{f}\right)^{-1},
	\]
	and as usual we derive that the diagram
	\[
	\begin{tikzcd}
		F(f)\big(\quot{\alpha}{Y}(\quot{A}{Y})\big) \ar[rr]{}{E(f)\quot{\eta_{\quot{A}{Y}}}{Y}} \ar[d, swap]{}{\left(\quot{\alpha_{\quot{A}{Y}}}{f}\right)^{-1}} & & E(f)\big(\quot{\beta}{Y}(\quot{A}{Y})\big) \\
		\quot{\alpha}{X}\big(F(f)(\quot{A}{Y})\big) \ar[rr, swap]{}{\quot{\eta_{F(f)\quot{A}{Y}}}{X}} & & \quot{\beta}{X}\big(F(f)(\quot{A}{Y})\big) \ar[u, swap]{}{\quot{\beta_{\quot{A}{Y}}}{f}}
	\end{tikzcd}
	\]
	commutes. Using this and the naturality of $\quot{\eta}{X}$ gives us that
	\begin{align*}
		\tau_f^{\underline{\beta}f} \circ E(f)\quot{\eta_{\quot{A}{Y}}}{Y} &=\quot{\beta}{X}(\tau_f^A) \circ \left(\quot{\beta_{\quot{A}{Y}}}{f}\right)^{-1} \circ E(f)\quot{\eta_{\quot{A}{Y}}}{Y} \\
		&= \quot{\beta}{X}(\tau_f^A) \circ \left(\quot{\beta_{\quot{A}{Y}}}{f}\right)^{-1} \circ \quot{\beta_{\quot{A}{Y}}}{f} \circ \quot{\eta_{F(f)\quot{A}{Y}}}{X} \circ \left(\quot{\alpha_{\quot{A}{Y}}}{f}\right)^{-1} \\
		&= \quot{\beta}{X}(\tau_f^A) \circ \quot{\eta_{F(f)\quot{A}{Y}}}{X} \circ \left(\quot{\alpha_{\quot{A}{Y}}}{f}\right)^{-1} \\
		&= \quot{\eta_{\quot{A}{X}}}{X} \circ \quot{\alpha}{X}(\tau_f^{A}) \circ \left(\quot{\alpha_{\quot{A}{Y}}}{f}\right)^{-1} = \quot{\eta_{\quot{A}{X}}}{X} \circ \tau_f^{\underline{\alpha}A},
	\end{align*}
	which proves that $\underline{\eta}_A:\underline{\alpha}A \to \underline{\beta}A$ is a $\PC(E)$-morphism.
	
	We now verify the naturality of $\underline{\eta}$. For this let $P = \lbrace \quot{\rho}{X} \; | \; X \in \Cscr_0\rbrace \in \PC(F)(A,B)$. We must prove that $\underline{\eta}_B \circ \underline{\alpha}P = \underline{\beta}P \circ \underline{\eta}_A$.  To see this note for all $X \in \Cscr_0$ we get that
	\[
	\quot{\beta}{X}\left(\quot{\rho}{X} \circ \quot{\eta_{\quot{A}{X}}}{X}\right) = \quot{\eta_{\quot{B}{X}}}{X} \circ \quot{\alpha}{X}(\quot{\rho}{X})
	\]
	which shows the commutativity of
	\[
	\begin{tikzcd}
		\underline{\alpha}A \ar[r]{}{\underline{\alpha}P} \ar[d, swap]{}{\underline{\eta}_A} & \underline{\alpha}B \ar[d]{}{\underline{\eta}_B} \\
		\underline{\beta}A \ar[r, swap]{}{\underline{\beta}P} & \underline{\beta}B
	\end{tikzcd}
	\]
	in $\PC(E)$. Thus $\underline{\eta}$ is a natural transformation between $\underline{\alpha}$ and $\underline{\beta}$.
\end{proof}

\begin{example}
	Consider the pseudocone categories $\PC(F) = D^b_G(X;\overline{\Q}_{\ell}) = \PC(E)$ induced by the pseudofunctor $D^b_c(-; \overline{\Q}_{\ell}) \circ \quo^{\op}:\SfResl_G(X)^{\op} \to \fCat$ and let $\alpha = [1] \circ [2]:D^b(-) \Rightarrow D^b(-)$ and $\beta = [2] \circ [1]:D^b(-) \Rightarrow D^b(-)$ be the degree $1 + 2$ and degree $2 + 1$ shift pseudonatural transformations. Then the natural isomorphisms
	\[
	\quot{[1]}{\Gamma} \circ \quot{[2]}{\Gamma} \xrightarrow[\eta_\Gamma]{\cong} \quot{[2]}{\Gamma} \circ \quot{1}{\Gamma}
	\]
	assemble to a modification $\eta:[2] \circ [1] \Rightarrow [1] \circ [2]$ which gives the natural isomorphism witnessing that $[2] \circ [1] \cong [1] \circ [2]$.
\end{example}

We now show that the various constructions we performed on $\Bicat(\Cscr^{\op},\fCat)$ are strictly $2$-functorial in the sense that they give rise to a strict $2$-functor $\PC:\Bicat(\Cscr^{\op},\fCat) \to \fCat$ which we call pseudconification\index[terminology]{Pseudoconification}. This, combined with the well-known Lemma \ref{Lemma: Section 3: 2-functors preserve adjoints} below (which we reprove here for the sake of completeness) allows us to lift adjoints in the bicategory $\Bicat(\Cscr^{\op},\fCat)$ to $\fCat$, which is our next main goal.

\begin{lemma}\label{Lemma: Section 3: The embedding of pre-equivariant functors to equivariant categories is a strict 2-functor}
	There is a strict $2$-functor 
	\[
	\PC:\Bicat(\Cscr^{\op},\fCat) \to \fCat
	\] 
	which sends pseudofunctors $F$ to their pseudocone categories $\PC(F)$, sends pseudonatural transformations $\alpha:F \to E$ to their induced functor $\underline{\alpha} =: \PC(\alpha)$, and sends modifications $\eta$ to their induced natural transformations $\underline{\eta} =: \PC(\eta)$.
\end{lemma}
\begin{proof}
	Because the $2$-functor is claimed to be strict, we have $\phi_{f,g} = \id_{\PC(f,g)}$ for all composable morphisms $f$ and $g$. Consequently we must verify that $\PC$ strictly preserves identity transformations, composition, and vertical composition of modifications.
	
	We fist note that the facts identity pseudonatural transformations and identity modifications get sent to the corresponding identity functors and identity natural transformations, i.e., that
	\[
	\begin{tikzcd}
		\Cscr^{\op} \ar[rr, bend left = 20, ""{name = U}]{}{F} \ar[rr, bend right = 20, swap, ""{name = D}]{}{F} & & \fCat \ar[from = U, to = D, Rightarrow, shorten <= 4pt, shorten >= 4pt]{}{\iota_F}
	\end{tikzcd}
	\]
	gets sent to
	\[
	\id_{\PC(F)}:\PC(F) \to \PC(F)
	\]
	and
	\[
	\begin{tikzcd}
		\Cscr^{\op} \ar[rrr, bend left = 40, ""{name = Up}]{}{F} \ar[rrr, bend right = 40, swap, ""{name = Down}]{}{F} & & & \fCat \ar[from = Up, to = Down, Rightarrow, shorten <= 4pt, shorten >= 4pt, bend left = 40, ""{name = Left}]{}{\alpha} \ar[from = Up, to = Down, Rightarrow, shorten <= 4pt, shorten >= 4pt, bend right = 40, swap, ""{name = Right}]{}{\alpha} \ar[from = Right, to = Left, symbol = \underset{\iota_{\alpha}}{\Rrightarrow}]
	\end{tikzcd}
	\]
	gets sent to
	\[
	\begin{tikzcd}
		\PC(F) \ar[rr, bend left = 30, ""{name = U}]{}{\ul{\alpha}} \ar[rr, bend right = 30, swap, ""{name = D}]{}{\ul{\alpha}} & & \fCat \ar[from = U, to = D, Rightarrow, shorten <= 4pt, shorten >= 4pt]{}{\iota_{\ul{\alpha}}}
	\end{tikzcd}
	\]
	for all relevant $F$ and $\alpha$ are both immediate from construction. To see that (vertical) composition of pseudonatural transformations is preserved, let $\alpha:F \to E$ and $\beta:E \to L$ be pseudonatural transformations. Then the composite $\PC(\beta \circ \alpha)$ is given on objects $(A,T_A)$ and satisfies the equations
	\begin{align*}
		\PC(\beta \circ \alpha)A &= \lbrace \quot{(\beta \circ \alpha)}{X}\quot{A}{X} \; | \; X \in \Cscr_0\rbrace = \lbrace (\quot{\beta}{X} \circ \quot{\alpha}{X})\quot{A}{X} \; | \; X \in \Cscr_0 \rbrace \\
		&= (\PC({\beta}) \circ \PC({\alpha}))A.
	\end{align*}
	We also now need to prove that $\tau_{f}^{(\PC(\beta \circ \alpha))A} = \tau_f^{\PC(\beta)(\PC(\alpha)A)}$ for all morphisms $f \in \Cscr_1$. For this let $f:X \to Y$. We now compute that on one hand
	\[
	\tau_{f}^{(\PC(\beta \circ \alpha))A} = \quot{\beta \circ \alpha}{X}(\tau_f^{A}) \circ \left(\quot{\alpha_{\quot{A}{Y}}}{f}\right)^{-1} = \quot{\beta}{X}\left(\quot{\alpha}{X}\big(\tau_f^A\big)\right) \circ \left(\quot{\alpha_{\quot{A}{Y}}}{f}\right)^{-1}.
	\]
	On the other hand we get from the pseudonaturality of $\alpha$ and $\beta$ that
	\begin{align*}
		\tau_f^{\PC(\beta)(\PC(\alpha)A)} &= \quot{\beta}{X}\left(\tau_f^{\PC(\alpha)A}\right) \circ \left(\quot{\beta_{\quot{\PC(\alpha)A}{Y}}}{f}\right)^{-1}\\
		& = \quot{\beta}{X}\left(\quot{\alpha}{X}\big(\tau_f^A\big) \circ \left(\quot{\alpha_{\quot{A}{Y}}}{f}\right)^{-1}\right) \circ \left(\quot{\beta_{\quot{\alpha}{Y}(\quot{A}{Y})}}{f }\right)^{-1}  \\
		&= \quot{\beta}{X}\left(\quot{\alpha}{X}\big(\tau_f^A\big)\right) \circ \quot{\beta}{X}\left(\quot{\alpha_{\quot{A}{Y}}}{f}\right)^{-1} \circ \left(\quot{\beta_{\quot{\alpha}{Y}(\quot{A}{Y})}}{f }\right)^{-1} \\
		&= \quot{\beta}{X}\left(\quot{\alpha}{X}\big(\tau_f^A\big)\right) \circ \left(\quot{\alpha_{\quot{A}{Y}}}{f}\right)^{-1} \\
		&= \tau_{f}^{(\PC(\beta \circ \alpha))A}
	\end{align*}
	which proves that $\tau_f^{\PC(\beta)(\PC(\alpha)A)} = \tau_{f}^{(\PC(\beta \circ \alpha))A}$. Thus it follows that
	\begin{align*}
		\PC(\beta \circ \alpha)(A,T_A) &= (\PC(\beta \circ \alpha))(A,T_A) = ((\PC(\beta \circ \alpha))A, T_{\PC(\beta \circ \alpha)A}) \\
		&= ((\PC(\beta) \circ \PC(\alpha))A, T_{(\PC(\beta) \circ \PC(\alpha))A}) \\
		&= (\PC(\beta) \circ \PC(\alpha))(A,T_A) = (\PC(\beta) \circ \PC(\alpha))(A,T_A),
	\end{align*}
	and similarly $\PC(\beta \circ \alpha)(P) = \PC(\beta)\big(\PC(\alpha)P\big)$ on morphisms $P$. Thus $\PC$ preserves $1$-compositions strictly. 
	
	To see that $2$-morphisms work out as expected, let $F,E:\Cscr^{\op} \to \fCat$ be pseudofunctors with pseudonatural transformations $\alpha,\beta,\gamma:F \Rightarrow E$ and modifications $\eta:\alpha \Rrightarrow \beta$ and $\epsilon:\beta \Rrightarrow \gamma$ as in the diagram
	\[
	\begin{tikzcd}
		F \ar[rr, bend left = 60, ""{name = U}]{}{\alpha} \ar[rr, ""{name = M}]{}[description]{\beta} \ar[rr, swap, bend right = 60, ""{name = L}]{}{\gamma} & & E \ar[from = U, to = M, Rightarrow, shorten <= 4pt, shorten >= 4pt]{}{\eta} \ar[from = M, to = L, Rightarrow, shorten <= 4pt, shorten >= 4pt]{}{\epsilon}
	\end{tikzcd}
	\]
	in the bicategory $\Bicat(\Cscr^{\op},\fCat)$. We now must check that
	\begin{align*}
		\PC(\epsilon \circ \eta) = \PC\epsilon \circ \PC\eta.
	\end{align*}
	However, this is routine as  for all $X \in \Cscr_0$ and for all $A \in \PC(F)_0$
	\begin{align*}
		\quot{(\PC(\epsilon \circ \eta)_A)}{X} &= \quot{(\underline{\epsilon \circ \eta})_A}{X} = \quot{(\epsilon \circ \eta)_{\quot{A}{X}}}{X} = \quot{\epsilon}{X}_{\quot{A}{X}} \circ \quot{\eta_{\quot{A}{X}}}{X} = \quot{\underline{\epsilon}_A}{X} \circ \quot{\underline{\eta}_A}{X} \\
		&= \quot{\PC(\epsilon)_A}{X} \circ \quot{\PC(\eta)_A}{X}
	\end{align*}
	so we conclude that $\PC(\epsilon \circ \eta) = \PC(\epsilon) \circ \PC(\eta)$, as was desired. The verification for horizontal composition is similar and omitted. Finally, the verification of the last remaining identity is trivial because the compositors $\phi_{\alpha,\beta}:\PC(\beta) \circ \PC(\alpha) \Rightarrow \PC(\beta \circ \alpha)$ are all identity natural transformations. Thus $\PC$ is a strict $2$-functor.
\end{proof}

We will now use modifications to describe adjoints that live in the $2$-category $\Bicat(\Cscr^{\op},\fCat)$ in terms of equational data. One benefit of this is that because adjoints in these $2$-categories are given by equational data, they are preserved by any strict $2$-functor. In particular, the strict $2$-functor $\PC$ of Lemma \ref{Lemma: Section 3: The embedding of pre-equivariant functors to equivariant categories is a strict 2-functor} lifts adjoints in $\Bicat(\Cscr^{\op},\fCat)$ to adjoints in $\fCat$. We follow the notion of adjunctions defined in a bicategory as given in \cite{RiehlVertityElementss}; this definition is equivalent to the definition of adjoints given in \cite{FreeAdj}, which defines a minimal $2$-category of the syntactic data required of an adjoint pair of functors and then gives an adjunction as a strict $2$-functor from this free-living adjunction category into a $2$-category $\Cfrak$.
\begin{definition}[{\cite[Defition 2.1.1]{RiehlVertityElementss}}]\index[terminology]{Adjunctions!In a $2$-Category}
	An adjunction in a bicategory $\Cfrak$ is a pair of objects $A$ and $B$ together with $1$-morphisms $f:A \to B$ and $g:B \to A$ such that there exist $2$-morphisms
	\[
	\begin{tikzcd}
		A \ar[r, bend left = 30, ""{name = UL}]{}{\id_A} \ar[bend right = 30, swap, r, ""{name = BL}]{}{g \circ f} & A \ar[from = UL, to = BL, Rightarrow, shorten >= 2pt, shorten <= 2pt]{}{\eta} & B \ar[r, bend left = 30, ""{name = UR}]{}{f \circ g} \ar[r, swap, bend right = 30, ""{name = BR}]{}{\id_B} & B \ar[from = UR, to = BR, Rightarrow, shorten <= 2pt, shorten >= 2pt]{}{\epsilon}
	\end{tikzcd}
	\]
	which satisfy the triangle identities below. The pasting diagram
	\[
	\begin{tikzcd}
		& B \ar[rr, equals] \ar[dr]{}{g} & {} & B \\
		A \ar[ur]{}{f} \ar[rr, equals] & {} \ar[u, Rightarrow, shorten >=5pt, shorten <= 5pt]{}{\eta} & A \ar[ur, swap]{}{f} \ar[u, swap, Rightarrow, shorten >=2pt, shorten <= 2pt]{}{\epsilon}
	\end{tikzcd}
	\]
	is equal to the $2$-cell
	\[
	\begin{tikzcd}
		A \ar[r, bend left = 30, ""{name = U}]{}{f} \ar[r, swap, bend right = 30, ""{name = B}]{}{f} & B \ar[from = U, to = B, Rightarrow, shorten >=2pt, shorten <= 2pt]{}{\iota_f}
	\end{tikzcd}
	\]
	and the pasting diagram
	\[
	\begin{tikzcd}
		& A \ar[rr, equals] \ar[dr]{}{f} \ar[d, Rightarrow, shorten >=5pt, shorten <= 5pt]{}{\epsilon} & {} \ar[d, swap, Rightarrow, shorten >=2pt, shorten <= 2pt]{}{\eta} & A \\
		B \ar[rr, equals] \ar[ur]{}{g} & {} & B \ar[ur, swap]{}{g} 
	\end{tikzcd}
	\]
	is equal to the $2$-cell:
	\[
	\begin{tikzcd}
		B \ar[r, bend left = 30, ""{name = U}]{}{g} \ar[r, swap, bend right = 30, ""{name = B}]{}{g} & A \ar[from = U, to = B, Rightarrow, shorten >=2pt, shorten <= 2pt]{}{\iota_g}
	\end{tikzcd}
	\]
	We write $f \dashv g$ to show that the pair $(A,B, f, g, \eta, \epsilon)$ is an adjunction with $f$ the left adjoint of $g$ (and $g$ the right adjoint of $f$).
\end{definition}
\begin{lemma}[{\cite[Lemma 2.1.3]{RiehlVertityElementss}}]\label{Lemma: Section 3: 2-functors preserve adjoints}
	Let $\Cfrak$ and $\Dfrak$ be bicategories. If $F:\Cfrak \to \Dfrak$ is a strict $2$-functor, then $F$ preserves any adjoints that exist in $\Cfrak$ on the nose.
\end{lemma}
\begin{proof}
	Because $F$ is strict, it preserves all equational data on the nose. Thus it sends the triangle identities in $\Cfrak$ on $(A,B,f,g,\eta,\epsilon)$ to the triangle identities in $\Dfrak$ on $(FA,FB,Ff,Fg,F\eta,F\epsilon)$.
\end{proof}
\begin{corollary}\label{Cor: Section 3: adjoints in PreEqSfResl get sent to adjoints in GEqCat}
	The $2$-functor
	\[
	\PC:\Bicat(\Cscr^{\op},\fCat) \to \fCat
	\]
	sends adjunctions between in $\Bicat(\Cscr^{\op},\fCat)$ to adjunctions between pseudocone categories.
\end{corollary}
\begin{proof}
	Apply Lemma \ref{Lemma: Section 3: 2-functors preserve adjoints} to the strict $2$-functor $\PC$ of Lemma \ref{Lemma: Section 3: The embedding of pre-equivariant functors to equivariant categories is a strict 2-functor}.
\end{proof}

This gives us our first main strategy and source of examples for adjunctions between pseudocone categories: Those that come from modifications. However, we would also like to know when we can take Change of Fibre functors $\PC(\alpha):\PC(F) \to \PC(E)$ and $\PC(\beta):\PC(E) \to \PC(F)$ which have adjoints $\quot{\alpha}{X} \dashv \quot{\beta}{X}:F(X) \to E(X)$, is there an adjunction $\PC(\alpha) \dashv \PC(\beta)$? We prove that, perhaps surprisingly\footnote{The surprising part is that if we have pseudonatural transformations $\alpha, \beta$ for which for all $X$ in $\Cscr_0$ there is an adjunction $\quot{\alpha}{X} \dashv \quot{\beta}{X}$ then this is enough to determine the unit and counit modifications $\eta$ and $\epsilon$ which witness $\alpha \dashv \beta$ in $\Bicat(\Cscr^{\op},\fCat)$.}, this is the case below.

\begin{Theorem}\label{Thm: Section 3: Gamma-wise adjoints lift to equivariant adjoints}
	Let $F,E:\Cscr^{\op} \to \fCat$ be pseudofunctors and let $\alpha:F \Rightarrow E$ and $\beta:E \to F$ be pseudonatural transformations. Then if $\quot{\alpha}{X} \dashv \quot{\beta}{X}$ in $\fCat$ for all $X \in \Cscr_0$, there is an adjunction
	\[
	\begin{tikzcd}
		F \ar[r, bend left = 30, ""{name = U}]{}{\alpha} & E \ar[l, bend left = 30, ""{name = L}]{}{{\beta}} \ar[from = U, to = L, symbol = \dashv]
	\end{tikzcd}
	\]
	in the $2$-category $\Bicat(\Cscr^{\op},\fCat)$. In particular, this gives rise to an adjunction of categories:
	\[
	\begin{tikzcd}
		\PC(F) \ar[r, bend left = 30, ""{name = U}]{}{\PC(\alpha)} & \PC(E) \ar[l, bend left = 30, ""{name = L}]{}{\PC(\beta)} \ar[from = U, to = L, symbol = \dashv]
	\end{tikzcd}
	\]
\end{Theorem}
\begin{proof}
	For all $X \in \Cscr_0$, let $\quot{\eta}{X}$ be the unit of adjunction
	\[
	\quot{\eta}{X}:\id_{F(X)} \Rightarrow \quot{\beta}{X} \circ \quot{\alpha}{X}
	\]
	and let $\quot{\epsilon}{X}$ be the counit of adjuction
	\[
	\quot{\epsilon}{X}:\quot{\alpha}{X} \circ \quot{\beta}{X} \Rightarrow \id_{E(X)}
	\]
	We begin by showing that the collection $\eta = \lbrace \quot{\eta}{X} \; | \; X \in \Cscr_0\rbrace$ forms a modification $\eta:\iota_F \to \beta \circ \alpha$; proving that $\epsilon := \lbrace \alpha \circ \beta \to \iota_E \; | \; X \in \Cscr_0 \rbrace$ is a modification is dual and will be omitted.
	
	Let $f \in \Cscr_1$ be arbitrary and write $f:X \to Y$. We now must show that the diagram
	\[
	\begin{tikzcd}
		F(f)\quot{A}{Y} \ar[rr]{}{\quot{\eta_{F(f)\quot{A}{Y}}}{X}} & & \big(\quot{\beta}{X} \circ \quot{\alpha}{X}\big)\left(F(f)\quot{A}{Y}\right) \ar[d]{}{\quot{(\beta \circ \alpha)_{\quot{A}{Y}}}{f}} \\
		F(f)\quot{A}{Y} \ar[u, equals] \ar[rr, swap]{}{F(f)\quot{\eta_{\quot{A}{Y}}}{Y}} & & F(f)\big((\quot{\beta}{Y} \circ \quot{\alpha}{Y})\quot{A}{Y}\big)
	\end{tikzcd}
	\]
	commutes for all $\quot{A}{Y}$ in $F(\quot{X}{Y})_0$. Consider the diagram:
	\[
	\begin{tikzcd}
		F(f)\quot{A}{Y} \ar[drrr] \ar[d, swap]{}{F(f)\quot{\eta_{\quot{A}{Y}}}{Y}} \ar[rrr]{}{\quot{\eta_{F(f)\quot{A}{Y}}}{X}} & & &  \big(\quot{\beta}{X} \circ \quot{\alpha}{X}\big)\left(F(f)\quot{A}{Y}\right) \ar[d,dashed]{}{\exists!\big(\quot{\beta}{X} \circ \quot{\alpha}{X}\big)(\rho)} \\  
		F(f)\big((\quot{\beta}{Y} \circ \quot{\alpha}{Y})\quot{A}{Y}\big)\ar[rrr, swap]{}{\big(\quot{(\beta \circ \alpha)_{\quot{A}{Y}}}{f}\big)^{-1}} & & & \big(\quot{\beta}{X} \circ \quot{\alpha}{X}\big)\left(F(f)\quot{A}{Y}\right)
	\end{tikzcd}
	\]
	By the adjunction $\quot{\alpha}{X} \dashv \quot{\beta}{X}$ there exists a unique map $\rho:F(f)\quot{A}{Y} \to F(f)\quot{A}{Y}$ making the above diagram commute. It is routine to check using the pseudonaturality of $\alpha$ and $\beta$, as well as the fact that units of adjunction are given uniquely up to unique isomorphism, that the morphism $\rho$ must be the identity on $F(f)\quot{A}{Y}$, giving that the diagram
	\[
	\begin{tikzcd}
		F(f) \ar[drr] \ar[d, swap]{}{F(f)\quot{\eta_{\quot{A}{Y}}}{Y}} \ar[rr]{}{\quot{\eta_{F(f)\quot{A}{Y}}}{X}} & &  \big(\quot{\beta}{X} \circ \quot{\alpha}{X}\big)\left(F(f)\quot{A}{Y}\right) \ar[d, equals] \\  
		F(f)\big((\quot{\beta}{Y} \circ \quot{\alpha}{Y})\quot{A}{Y}\big)\ar[rr, swap]{}{\big(\quot{(\beta \circ \alpha)_{\quot{A}{Y}}}{f}\big)^{-1}} & & \big(\quot{\beta}{X} \circ \quot{\alpha}{X}\big)\left(F(f)\quot{A}{Y}\right)
	\end{tikzcd}
	\]
	commutes. By post-composing the induced equation by $\quot{(\beta \circ \alpha)}{f}$ at $F(f)\quot{A}{Y}$ we get 
	\begin{align*}
		\quot{(\beta \circ \alpha)_{\quot{A}{Y}}}{f} \circ \quot{\eta_{F(f)\quot{A}{X}}}{X} &= \quot{(\beta \circ \alpha)_{\quot{A}{Y}}}{f} \circ \left(\quot{(\beta \circ \alpha)_{\quot{A}{Y}}}{f}\right)^{-1} \circ F(f)\left(\quot{\eta_{\quot{A}{Y}}}{Y}\right) \\
		&= F(f)\left(\quot{\eta_{\quot{A}{Y}}}{Y}\right),
	\end{align*}
	which shows that the diagram
	\[
	\xymatrix{
		F(f)\quot{A}{Y} \ar[rr]^-{\quot{\eta_{F(f)\quot{A}{Y}}}{Y}} & & \big(\quot{\beta}{X} \circ \quot{\alpha}{X}\big)F(f)\quot{A}{Y} \ar[d]^{\quot{(\beta \circ \alpha)_{\quot{A}{Y}}}{f}} \\
		F(f)\quot{A}{Y} \ar@{=}[u] \ar[rr]_-{F(f)\quot{\eta_{\quot{A}{Y}}}{Y}} & & F(f)\big((\quot{\beta}{Y} \circ \quot{\alpha}{Y})\quot{A}{Y}\big)
	}
	\]
	commutes and hence that $\eta$ is a modification. That $\epsilon$ is a modification is given dually.
	
	For the triangle identities, consider the pasting diagram in the bicategory \\ $\Bicat(\Cscr^{\op},\fCat)$:
	\[
	\begin{tikzcd}
		& E \ar[rr, equals] \ar[dr]{}{\beta} & {} & B \\
		F \ar[ur]{}{\alpha} \ar[rr, equals] & {} \ar[u, Rightarrow, shorten >=5pt, shorten <= 5pt]{}{\eta} & A \ar[ur, swap]{}{\alpha} \ar[u, swap, Rightarrow, shorten >=2pt, shorten <= 2pt]{}{\epsilon}
	\end{tikzcd}
	\]
	To prove that this glues (or untwists, if you prefer) to the correct diagram, note that for all $X \in \Cscr_0$ there is an equivalence of pasting diagrams
	\[
	\begin{tikzcd}
		& E(X) \ar[rr, equals] \ar[ddr]{}{\quot{\beta}{X}} & {} & E(X) & & E(X)\\ 
		{}& {}& {}& {} & = \\
		F(X) \ar[uur]{}{\quot{\alpha}{X}} \ar[rr, equals] & {} \ar[uu, Rightarrow, shorten >=5pt, shorten <= 5pt]{}{\quot{\eta}{X}} & A \ar[uur, swap]{}{\quot{\alpha}{X}} \ar[uu, swap, Rightarrow, shorten >=2pt, shorten <= 2pt]{}{\quot{\epsilon}{X}} & & & F(X) \ar[uu, bend left = 30, ""{name = L}]{}{\quot{\alpha}{X}} \ar[uu, bend right = 30, swap, ""{name = R}]{}{\quot{\alpha}{X}} \ar[Rightarrow, from = L, to = R, shorten <= 2pt, shorten >=2pt]{}{\iota_{\quot{\alpha}{X}}}
	\end{tikzcd}
	\]
	which in turn shows that the pasting diagram in $\Bicat(\Cscr^{\op},\fCat)$ is equal to the $2$-cell:
	\[
	\begin{tikzcd}
		F \ar[r, bend left = 30, ""{name = U}]{}{\alpha} \ar[swap, bend right = 30, r, ""{name = L}]{}{\alpha} & E \ar[Rightarrow, from = U, to = L, shorten <= 2pt, shorten >= 2pt]{}{\iota_{\alpha}}
	\end{tikzcd}
	\]
	The other triangle identity is given dually and omitted. Thus we get that the $6$-tuple $(F,E,\alpha,\beta,\eta,\epsilon)$ describes an adjunction in $\Bicat(\Cscr^{\op},\fCat)$, and hence by Corollary \ref{Cor: Section 3: adjoints in PreEqSfResl get sent to adjoints in GEqCat} in $\fCat$ as well.
\end{proof}
\begin{corollary}\label{Cor: Section 3: Fibre-wise equivalences are equivalences of equivariant categories}
	Let $\PC(\alpha) \dashv \PC(\beta):\PC(F) \dashv \PC(E)$ be psuedocone functors such that the adjoints are induced as in Theorem \ref{Thm: Section 3: Gamma-wise adjoints lift to equivariant adjoints}. Then if the adjoints $\quot{\alpha}{X} \dashv \quot{\beta}{X}$ are all equivalences of categories, so is $\PC(\alpha) \dashv \PC(\beta)$.
\end{corollary}
\begin{proof}
	By Theorem \ref{Thm: Section 3: Gamma-wise adjoints lift to equivariant adjoints}, we have that the adjoint equivalences $\quot{\alpha}{X} \dashv \quot{\beta}{X}$ lift/descend to an adjunction $\PC(\alpha) \dashv \PC(\beta):\PC(F) \to \PC(E)$. From the $X$-local natural isomorphisms $\quot{\alpha}{X} \circ \quot{\beta}{X} \cong \id_{E(X)}$ and $\quot{\beta}{X} \circ \quot{\alpha}{X} \cong \id_{F(X)}$ induced by the units and counits of adjunction, we get isomorphisms $\PC(\beta) \circ \PC(\alpha) \cong \id_{\PC(F)}$ and $\PC(\alpha) \circ \PC(\beta) \cong \id_{\PC(E)}$ induced by the unit and counit of adjunction.
\end{proof}

Let us now proceed to examine the fact that Theorem \ref{Thm: Section 3: Gamma-wise adjoints lift to equivariant adjoints} provides us a clean way to discuss when pseudocone categories are monoidal closed categories, toposes, and other such things because it allows us to deduce the existence of adjoints based on the existence of $X$-local adjoints. In particular, this gives us the necessary tools to provide a categorical foundation for the various adjunctions on $\Bicat(\SfResl_G(X)^{\op},\fCat)$ and/or $\Bicat(\FResl_H(Y)^{\op},\fCat)$ we need to perform and study equivariant algebraic geometry or equivariant algebraic topology; note in this case $G$ is a smooth algebraic group and $X$ is a left $G$-variety while $H$ is a topological group and $Y$ is a left $H$-space.

\begin{definition}
	Let $F:(\Cscr, \otimes, I) \to (\Dscr, \boxtimes, J)$ be a functor between\\ monoidal closed categories. We say that $F$ is monoidal closed\index[terminology]{Monoidal Closed Functor} if $F$ is simultaneously monoidal and closed. Additionally we say that $F$ is symmetric monoidal closed if it is simultaneously symmetric monoidal and closed.
\end{definition}
\begin{proposition}\label{Section: Pseudocone Functors: PC cat is monoidal closed}
	Let $F:\Cscr^{\op} \to \fCat$ be a monodial closed pseudofunctor for which each functor $F(f)$ is monoidal closed. Then if each category $F(X)$ is monoidal closed or symmetric monoidal closed, so is $\PC(F)$.
\end{proposition}
\begin{proof}
	By Theorem \ref{Theorem: Section 2: Monoidal preequivariant pseudofunctor gives monoidal equivariant cat} and Proposition \ref{Prop: Section 2: Equivariant cat is symmetric monoidal} it follows that in either case we have that $\PC(F)$ is monoidal or symmetric monoidal, repsectively. We now prove this lemma by providing the case in which $\PC(F)$ is only know to be monoidal; the symmetric monoidal case follows mutatis mutandis. From Theorem \ref{Thm: Section 3: Gamma-wise adjoints lift to equivariant adjoints} and the fact that we know there are adjoints
	\[
	\quot{A}{X}\os{X}{\otimes} (-) \dashv  [\quot{A}{X},-]_{X}
	\]
	for all $X \in \Cscr_0$, it suffices to prove that for any object $A \in \PC(F)_0$, the functors $A \otimes (-):\PC(F) \to \PC(F)$ and $[A,-]:\PC(F) \to \PC(F)$ arise from pseudonatural transformations as in Theorem \ref{Thm: Functor Section: Psuedonatural trans are pseudocone functors}. To this end let
	\[
	\theta_f^{A,B}:F(f)\left(\quot{A}{Y} \os{Y}{\otimes} \quot{B}{Y}\right) \xrightarrow{\cong} F(f)(\quot{A}{Y}) \os{X}{\otimes} F(f)(\quot{B}{Y})
	\]
	and
	\[
	\xi_{f}^{A,B}:F(f)[\quot{A}{Y},\quot{B}{Y}]_{Y} \xrightarrow{\cong} [F(f)(\quot{A}{Y}),F(f)(\quot{B}{Y})]_{X}
	\]
	be the preservation isomorphisms. 
	
	Fix an object $A \in \PC(F)_0$. We then define pseudonatrual transformation $A \otimes (-):F \to F$ as follows: For all $X \in \Cscr_0$, define the functor
	\[
	\quot{(A \otimes (-))}{X}:F(X) \to F(X)
	\]
	by
	\[
	\quot{(A \otimes (-))}{X} := \quot{A}{X} \os{X}{\otimes} (-).
	\]
	Fix a morphism $f:X \to Y$ in $\Cscr$. Now observe that on one hand we have that for all $\quot{B}{Y} \in F(\quot{X}{Y})_0$,
	\[
	\big(\quot{(A \otimes (-))}{X} \circ F(f)\big)(\quot{B}{Y}) = \quot{A}{X} \os{X}{\otimes} F(f)(\quot{B}{Y})
	\]
	while on the other hand we have
	\[
	\big(F(f) \circ \quot{(A \otimes (-))}{Y}\big)(\quot{B}{Y}) = F(f)\left(\quot{A}{Y} \os{Y}{\otimes} \quot{B}{Y}\right).
	\]
	There is then an isomorphism
	\[
	\xymatrix{
		\quot{A}{X} \os{X}{\otimes} F(f)(\quot{B}{Y}) \ar[rrr]^-{\left(\tau_f^A\right)^{-1}\os{X}{\otimes} \id_{F(f)\quot{B}{Y}}} \ar[drrr] & & & F(f)(\quot{A}{Y}) \os{X}{\otimes} F(f)(\quot{B}{Y}) \ar[d]^{\left(\theta_f^{A,B}\right)^{-1}} \\
		&	 & & F(f) \big(\quot{A}{Y} \os{Y}{\otimes} \quot{B}{Y}\big)
	}
	\]
	which is natural in $\quot{B}{Y}$ by the bifunctoriality of $\otimes$ and the natruality of $\theta_f$ (and hence its inverse). This determines our natural transformation $\quot{(A \otimes (-))}{f}$ component of the pseudonatural transformation $A \otimes (-)$; its object assignment is given by
	\[
	\quot{(A \otimes (-))_{\quot{B}{Y}}}{f} := \left(\theta_f^{A,B}\right)^{-1} \circ \left(\tau_f^A\right)^{-1} \os{X}{\otimes} \id_{F(f)(\quot{B}{Y})}.
	\]
	A routine check shows that the collection 
	\[
	A \otimes (-) = \lbrace \quot{(A \otimes (-))}{X}, \quot{(A \otimes(-))}{f} \; | \; X \in \Cscr_0, f \in \Cscr_1\rbrace 
	\]
	determines a pseudonatrual transformation which we identify (through Lemma \ref{Lemma: Section 3: The embedding of pre-equivariant functors to equivariant categories is a strict 2-functor}) with the functor $A \otimes (-)$.
	
	The fact that $[A,-]$ gives rise to a pseudonatural transformation is dual to the argument above. More explicitly, the object functors $\quot{[A,-]}{X}$ are given by
	\[
	\quot{[A,-]}{X}(\quot{B}{X}) := [\quot{A}{X}, \quot{B}{X}]_{X}
	\]
	and the natural isomorphism $\quot{[A,-]}{f}:\quot{[A,-]}{X} \circ F(f) \Rightarrow F(f) \circ \quot{[A,-]}{Y}$ are given by
	\[
	\quot{[A,-]_{-}}{f} := \left(\xi_f^{A,-}\right)^{-1} \circ [\tau_f^A,\id_{F(f)(-)}].
	\]
	As before, we identify the pseudonatural transformation $[A,-]$ with the functor $[A,-]:\PC(F) \to \PC(F)$. Now appealing to Theorem \ref{Thm: Section 3: Gamma-wise adjoints lift to equivariant adjoints} gives the adjunction $A \otimes (-) \dashv [A,-]$, as desired. 
\end{proof}
\begin{corollary}\label{Cor: Pseudocone Functors: Monoidal Closed Derived Cat}
	If $G$ is a smooth algebraic group and if $X$ is a left $G$-variety then $D_G^b(X;\overline{\Q}_{\ell})$ and $D_G^b(X)$ are symmetric monoidal closed categories.
\end{corollary}
\begin{proof}
	Because each morphism $\overline{f \times \id_X}:G \backslash (\Gamma \times X) \to G \backslash (\Gamma^{\prime} \times X)$ is smooth it follows that each functor $(\overline{f \times \id_X})^{\ast}$ is monoidal closed (that each functor is monoidal is a standard result in {\'e}tale cohomology while the closed aspect follows from \cite[Equation 4.2.5.1]{BBD}). Consequently it follows that $D_G^b(X;\overline{\Q}_{\ell})$ is symmetric monoidal closed.
\end{proof}
\begin{corollary}\label{Cor: Pseudocone Functors: Monoidal Closed Derived Cat Topological}
	If $G$ is a topological group which admits $n$-acyclic free $G$-spaces which are manifolds and if $X$ is a left $G$-space then $D_G^b(X)$ is a symmetric monoidal closed category. 
\end{corollary}
\begin{proof}
	By the assumptions on $G$ it follows that for any $n \in \N$ if $M$ is a specified $n$-acyclic free $G$-space which is a manifold then the second projection $\pi_2:M \times X \to X$ is a smooth $n$-acyclic resolution of $X$. From here taking the subcategory $\mathbf{SResl}_G(X)$ of $G$-resolutions of $X$ generated by the projections $M \times X \to X$, their pullbacks against each other, and the smooth maps between them, we get that there is an equivalence 
	\[
	D_G^b(X) \simeq \PC(D^b_c(G \backslash (-)):\mathbf{SResl}_G(X)^{\op} \to \fCat)
	\]
	by \cite[Section 2.4.4]{BernLun}. However, because every resolution $p^{\ast}:P \to X$ in $\mathbf{SResl}_G(X)$ is smooth and every morphism $f:P \to Q$ between such resolution is smooth, the pullback functors $p^{\ast}$ and $f^{\ast}$ all commute with the tensor functors $(-) \otimes_{D^b(G \backslash P)} (-)$ and $[-,-]_{D^b(G \backslash P)}$ by \cite[Section 1.8, Appendix A]{BernLun}. As such it follows that we induce pseudonatural transformations
	\[
	\begin{tikzcd}
		\mathbf{SResl}_G(X)^{\op} \ar[rr, bend left = 30, ""{name = U}]{}{D^b(G \backslash (-))} \ar[rr, bend right = 30, swap, ""{name = D}]{}{D^b(G \backslash (-))} & & \fCat \ar[from = U, to = D, Rightarrow, shorten <= 4pt, shorten >= 4pt]{}[description]{(-) \otimes A}
	\end{tikzcd}\qquad
	\begin{tikzcd}
		\mathbf{SResl}_G(X)^{\op} \ar[rr, bend left = 30, ""{name = U}]{}{D^b(G \backslash (-))} \ar[rr, bend right = 30, swap, ""{name = D}]{}{D^b(G \backslash (-))} & & \fCat \ar[from = U, to = D, Rightarrow, shorten <= 4pt, shorten >= 4pt]{}[description]{[A,-]}
	\end{tikzcd}
	\]
	for all objects $A \in \PC(D^b(G \backslash (-)))_0$. We now apply Theorem \ref{Thm: Functor Section: Psuedonatural trans are pseudocone functors} to get our tensor and internal hom functors $(-) \otimes A, [A,-]:D_G^b(X) \to D_G^b(X)$ and then apply Theorem \ref{Thm: Section 3: Gamma-wise adjoints lift to equivariant adjoints} to deduce the desired adjunction $(-) \otimes A \dashv [A,-]$.
\end{proof}
\begin{remark}\label{Remark: Pseudocone Functors: Why the Assume on G?}
	This is an alternative, arguably less direct, approach to showing that $D_G^b(X)$ is symmetric monoidal when $G$ is a topological group, but the point of our proof is to illustrate that ultimately it is a consequence of the pseudocone formalism that the symmetric monoidal structure exists in the first place. Moreover, by \cite[Section 3.1]{BernLun} this holds whenever $G$ is either a linear Lie group (a closed subgroup of $\GL_n(\R)$ for some $n \in \N$) or whenever $G$ is a Lie group with a finite number of connected components. Note that in both cases this implies that the category $D_G^b(X)$ is symmetric monoidal whenever $G$ is a finite group.
\end{remark}
\begin{corollary}\label{Cor: Section 3: Equivariant Cartesian Closed Cats}
	Let $F:\Cscr^{\op} \to \fCat$ be a pseudofunctor such that each category $F(X)$ has finite products and each functor $F(f)$ preserves these products. Assume further that each category $F(X)$ is Cartesian closed and that each fibre functor $F(f)$ is a Cartesian closed functor. Then $\PC(F)$ is Cartesian closed as well.
\end{corollary}
\begin{proof}
	This is simply Proposition \ref{Section: Pseudocone Functors: PC cat is monoidal closed} applied to a Cartesian monoidal closed structure.
\end{proof}

This leads us to examine when $\PC(F)$ is a topos as determined by pseudofunctors $F:\Cscr^{\op} \to \fCat$ whose fibre categories are toposes and whose fibre functors are inverse images of geometric morphisms.
\begin{proposition}\label{Proposition: Section 3: Equivariant Toposes}
	Let $F:\Cscr^{\op} \to \fCat$ be a pseudofunctor such that each category $F(X)$ is a topos and each morphism $F(f)$ is the pullback of a geometric morphism which preserves the subobject classifier and internal hom functors. Then $\PC(F)$ is a topos.
\end{proposition}
\begin{proof}
	Recall that a topos is a category $\Cscr$ which has finite limits, is Cartesian closed, and has a subobject classifier (cf.\@ \cite{MacLaneMoerdijk}). By Proposition \ref{Prop: Section 2: Equivariant Cat has Subobject Classifiers}, $\PC(F)$ has a subobject classifier and by Corollary \ref{Cor: Equivariant cat is finitely complete or cocomplete} $\PC(F)$ is finitely complete. Finally by Corollary \ref{Cor: Section 3: Equivariant Cartesian Closed Cats} we have that $\PC(F)$ is Cartesian closed. Therefore $\PC(F)$ is a topos.
\end{proof}
\begin{example}
	An example of such a topos is one whose fibre functors $F(f)$ are all the inverse image functors of an atomic geometric morphism. In particular, if each functor $F(f) = f^{\ast}:F(Y) \to F(X)$ is an inverse image of an {\'e}tale geometric morphism $f:F(X) \to F(Y)$, then $\PC(F)$ is a topos.
\end{example}

We now close this section with a topic of interest for equivariant homotopy theory (and in particular for a consideration of equivariant homological algebra): pseudocone localization. Our primary motivation for this is to determine in what sense, if any, the equivariant derived category $D^b_G(X)$ arises as a localization of the category of bounded chain complexes $\Ch^b_G(X)$ at some suitable class of quasi-isomorphisms\footnote{Or as close as we can get, anyway.}. While it is known that the equivariant derived categories are not derived categories (in the sense that generically it is the case that $D_G^b(X) \not\simeq Q^{-1}\Shv_G^b(X)$) we \emph{do} want to know if it is the case that there is a class of morphisms in a pseudocone category which allow us to see $D_G^b(X)$ as a ``pseudolocalization\footnote{This term is left vague and undefined primarily to give a sense of what we are looking for.}'' as determined by some pseudofunctor and special class of morphisms cut out by the pseudofunctor. More precisely, this is akin to asking if $F:\Cscr^{\op} \to \fCat$ is a pseudofunctor which takes values in localizations of categories, in what (necessarily weaker) sense can we see $\PC(F)$ as a localization? The moderate technical conditions we will impose show when we can essentially localize the fibre categories of a pseudofunctor over the data specified by the categories $F(X)$ and the fibre functors $F(f)$.

Begin by letting $F:\Cscr^{\op} \to \fCat$ be a pseudofunctor and, for all $X \in \Cscr_0$, let $\quot{S}{X} \subseteq \PC(F)_1$. Our goal now is to study when we can descend the data of the localizations $\quot{S}{X}^{-1}F(X)$ through the fibre functors $F(f)$, i.e., when there are functors $(S^{-1}F)(f):(\quot{S}{Y})^{-1}F(Y) \to (\quot{S}{X})^{-1}F(X)$ which allow us to define a suitable localization $\lambda_S:\PC(F) \to S^{-1}\PC(F)$ (which we will call a pseudocone localization). The evident condition we use to produce the pseudofunctor $S^{-1}F$ is that we need $F(f)$ to carry every $\quot{S}{Y}$-morphism to a $\quot{S}{X}$-morphism; it turns out this is the only structure constraint we need in order to derive the following proposition, which shows that taking localizations in this way gives rise to a  pseudofunctor $S^{-1}F$ whose fibre categories are all the localizations of the fibre categories of $F$. Subsequently, we show that there is a pseudonatural transformation $\lambda:F \implies S^{-1}F$ which realizes the category $S^{-1}\PC(F)$ as a pseudocone localization. Finally, after discussing pseudocone localization and the property it has (cf.\@ Propositions \ref{Prop: Section 3: Equivariant Localization Pseudofunctor}, \ref{Prop: Section 3: Equivariant localization pseudonat} and Theorem \ref{Theorem: Section 3: Equivariant localization is a 2-colimit}), we will compare this with the na{\"i}ve localization. In particular, we show that the pseudocone localization is not a literal localization (cf.\@ Example \ref{Example: Section 3: Equivariant localizaiton is not the localization of equivariant category}) but does have an analogous universal property when restricted to pseudocone categories and functors and then provide an explicit description of in what precise sense can we see the equivariant derived category as a localization of $\Ch^b_G(X)$.

\begin{proposition}\label{Prop: Section 3: Equivariant Localization Pseudofunctor}
	Let $F:\Cscr^{\op} \to \fCat$ be a  pseudofunctor and for all $X \in \Cscr_0$ assume that there is a subclass $\quot{S}{X} \subseteq F(X)_1$ such that for all $f \in \Cscr_1$, $F(f)(\quot{S}{Y}) \subseteq \quot{S}{X}$ where $f:X \to Y$. Then if $S \subseteq \PC(F)_1$ is the subclass of morphisms
	\[
	S := \lbrace P \in \PC(F)_1 \; | \; \forall\,X \in \Cscr_0.\,\forall\,\quot{\rho}{X} \in P.\, \quot{\rho}{X} \in \quot{S}{X} \rbrace
	\]
	there is a  pseudofunctor $S^{-1}F:\Cscr^{\op} \to \fCat$\index[notation]{SinverseF@$S^{-1}F$} which sends objects $X $ to $(\quot{S}{X})^{-1}F(X)$ and sends morphisms $f \times \id_X:X  \to Y $ to the unique functor $(S^{-1}F)(f)$ fitting in the commuting square:
	\[
	\xymatrix{
		F(Y) \ar[rr]^{F(f)} \ar[d]_{\lambda_{\quot{S}{Y}}} & & F(X) \ar[d]^{\lambda_{\quot{S}{X}}} \\
		(\quot{S}{Y})^{-1}F(Y) \ar@{-->}[rr]_{\exists!(S^{-1}F)(f)} & & (\quot{S}{X})^{-1}F(X)
	}
	\]
\end{proposition}
\begin{proof}
	We first note that the categories $(\quot{S}{X})^{-1}F(X)$ all exist by \cite[Section 1.1]{GabrielZisman}. The functor $(S^{-1}F)(f)$ exists by first observing that for any $\varphi \in \quot{S}{Y}$ since $F(f)(\varphi) \in \quot{S}{X}$, $(\lambda_{\quot{S}{X}} \circ F(f))(\varphi)$ is an isomorphism in $(\quot{S}{X})^{-1}F(X)$; applying the universal property of the localization to the diagram
	\[
	\xymatrix{
		F(Y) \ar[rr]^-{\lambda_{\quot{S}{X}} \circ F(f)} \ar[dr]_{\lambda_{\quot{S}{Y}}} & & (\quot{S}{X})^{-1}F(X) \\
		& (\quot{S}{Y})^{-1}F(Y) \ar@{-->}[ur]_{\exists!(SF)^{-1}(f)}
	}
	\]
	gives the existence of $(S^{-1}F)(f)$.
	
	Note that since the diagram
	\[
	\xymatrix{
		F(X) \ar[d] \ar@{=}[r] & F(X) \ar[d] \\
		(\quot{S}{X})^{-1}F(X) \ar@{=}[r] & (\quot{S}{X})^{-1}F(X)
	}
	\]
	commutes, we have that $(S^{-1}F)(\overline{\id}_{X}) = \id_{(\quot{S}{X})^{-1}F(X)}$. Thus we only need to establish the existence of the compositor natural isomorphisms, for $X \xrightarrow{f} Y \xrightarrow{g} Z$ in $\Cscr$,
	\[
	\phi_{f,g}:(S^{-1}F)(f) \circ (S^{-1}F)(g) \Rightarrow (S^{-1}F)(g \circ f)
	\]
	and their identities to prove that $S^{-1}F$ eixsts. For this fix $f$ and $g$ as above and consider that because $F$ is a pseudofunctor we have natural transformations $\phi_{f,g}$ fitting into a pasting diagram:
	\[
	\begin{tikzcd}
		& F(Y) \ar[d, Rightarrow, shorten >= 4pt, shorten <= 4pt]{}{\phi_{f,g}} \ar[dr]{}{F(f)} \\
		F(Z) \ar[ur]{}{F(g)} \ar[rr, swap]{}{F(g \circ f)} & {} & F(X)
	\end{tikzcd}
	\]
	We now recall (cf.\@ \cite{BorceuxCatAlg1}, \cite{GabrielZisman}, or \cite{MacLaneCWM}) that from the construction of the categories $(\quot{S}{X})^{-1}F(X)$ we have that for all $Y \in \Cscr_0$,
	\[
	(\quot{S}{Y})^{-1}F({Y})_0 = F({Y})_0.
	\]
	In particular, it follows from the universal property of the localizations and the constructions of the $\lambda$ functors that the object assignments of each functor $(S^{-1}F)(f)$, $(S^{-1}F)(g)$, and $(S^{-1}F)(g \circ f)$ must be the same as the object assignments of $F(f)$, $F(g)$, and $F(g \circ f)$, respectively. Then for any $A \in (\quot{S}{Z})^{-1}F(Z)_0$, we have
	\[
	\big((S^{-1}F)(f) \circ (S^{-1}F)(g)\big)A = \big(F(f) \circ F(g)\big)A
	\] 
	(and similarly for $F(g \circ f)$). Thus any natural transformation 
	\[
	\alpha:(S^{-1}F)(f) \circ (S^{-1}F)(g) \Rightarrow (S^{-1}F)(g \circ f)
	\]
	must be a collection of $(\quot{S}{X})^{-1}F(X)$-morphisms
	\[
	\alpha = \lbrace \alpha_A:(F(f) \circ F(g))A \to F(g \circ f)A \; | \; A \in (\quot{S}{Z})^{-1}F(Z) \rbrace.
	\]
	
	We now prove that $\phi_{f,g}$ extends to a compositor natural isomorphism 
	\[
	S^{-1}\phi_{f,g}:(S^{-1}F)(f) \circ (S^{-1}F)(g) \Rightarrow (S^{-1}F)(g \circ f)
	\] 
	where
	\[
	S^{-1}\phi_{f,g}:= \lbrace \phi_{f,g}^{A}:(F(f) \circ F(g))A \to F(g \circ f)A\; | \; A \in (\quot{S}{Z})^{-1}F(Z)_0 \rbrace.
	\]
	That this is an isomorphism for each $A$ is trivial from the fact that $\phi_{f,g}$ is; as such, we only need to prove the naturality holds, i.e., that for any zig-zag describing a morphism $\alpha:A \to B$ in $(\quot{S}{Z})^{-1}F(Z)$, the diagram
	\[
	\xymatrix{
		\big((S^{-1}F)(f) \circ (S^{-1}F)(g)\big)A \ar[rr]^-{\phi_{f,g}^{A}} \ar[d]_{\big((S^{-1}F)(f) \circ (S^{-1}F)(g)\big)\alpha} & & (S^{-1}F)(g \circ f)A \ar[d]^{(S^{-1}F)(g \circ f)\alpha} \\
		\big((S^{-1}F)(f) \circ (S^{-1}F)(g)\big)B \ar[rr]_-{\phi_{f,g}^{B}} & & (S^{-1}F)(g \circ f)B
	}
	\]
	commutes.
	
	To prove naturality of the square above, let $\alpha$ be represented by some zig-zag of the form
	\[
	\xymatrix{
		A = A_0  & & A_1 &\cdots & A_n  & & A_{n+1} = B\\
		& C_0 \ar[ur]_{\rho_0} \ar[ul]^{\psi_0} &  & \cdots \ar[ur] \ar[ul] & & C_n \ar[ur]_{\rho_n} \ar[ul]^{\psi_n}
	}
	\]
	where $\rho_i \in F(Z)_1$ and where $\psi_i \in \quot{S}{Z}$ for all $ 0 \leq i \leq n$ and for $n \in \N$. Note that not only is the argument we will present independent of the representation of $\alpha$, but proving the other cases (where $\alpha$ starts with a $\rho_i$, where $\alpha$ ends with a $\psi_i^{-1}$, etc.) follows mutatis mutandis and will be omitted.
	
By construction, since
	\[
	\alpha = \prod_{i=0}^{n} \rho_i \circ \psi_i^{-1}
	\]
	from the universal properties defining all the objects we have that
	\begin{align*}
		&\big((S^{-1}F)(f) \circ (S^{-1}F)(g)\big)\alpha \\
		&= \big((S^{-1}F)(f) \circ (S^{-1}F)(g)\big)\left(\prod_{i=0}^{n} \rho_i \circ \psi_i^{-1}\right) \\
		& = \prod_{i=0}^{n} \big((S^{-1}F)(f) \circ (S^{-1}F)(g)\big)\rho_i \circ \big((S^{-1}F)(f) \circ (S^{-1}F)(g)\big)\psi_i^{-1} \\
		&= \prod_{i=0}^{n} \big(F(f) \circ F(g)\big)\rho_i \circ \big(F(f) \circ F(g)\big)(\psi_i)^{-1}
	\end{align*}
	and similarly
	\[
	(S^{-1}F)(g \circ f)\alpha = \prod_{i=0}^{n} F(g \circ f)\rho_i \circ F(g \circ f)(\psi_i)^{-1}.
	\]
	Now note that from the naturality of $\phi_{f,g}$ the diagrams
	\[
	\xymatrix{
		\big(F(f) \circ F(g)\big)C_i \ar[rr]^-{\phi_{f,g}^{C_i}} \ar[d]_{\big(F(f) \circ F(g)\big)\rho_i} & & F(g \circ f)C_i \ar[d]^{F(g \circ f)\rho_i} \\
		\big(F(f) \circ F(g)\big)A_{i+1} \ar[rr]_-{\phi_{f,g}^{A_{i+1}}} & & F(g \circ f)A_{i+1}
	}
	\]
	commute in both $F(X)$ and $(\quot{S}{X})^{-1}F(X)$ for all $0 \leq i \leq n$. Similarly, the diagrams
	\[
	\xymatrix{
		\big(F(f) \circ F(g)\big)C_i \ar[rr]^-{\phi_{f,g}^{C_i}} \ar[d]_{\big(F(f) \circ F(g)\big)\psi_i} & & F(g \circ f)C_i \ar[d]^{F(g \circ f)\psi_i} \\
		\big(F(f) \circ F(g)\big)A_{i} \ar[rr]_-{\phi_{f,g}^{A_{i}}} & & F(g \circ f)A_{i}
	}
	\]
	commute in both $F(X)$ and $(\quot{S}{X})^{-1}F(X)$ for all $0 \leq i \leq n$; however, upon inverting both vertical morphisms, we find that the diagrams
	\[
	\xymatrix{
		\big(F(f) \circ F(g)\big)C_i \ar[rr]^-{\phi_{f,g}^{C_i}} & & F(g \circ f)C_i  \\
		\big(F(f) \circ F(g)\big)A_{i} \ar[rr]_-{\phi_{f,g}^{A_{i}}} \ar[u]^{\big(F(f) \circ F(g)\big)(\psi_i)^{-1}} & & F(g \circ f)A_{i} \ar[u]_{F(g \circ f)(\psi_i)^{-1}}
	}
	\]
	must also commute in $(\quot{S}{X})^{-1}F(X)$ for all $0 \leq i \leq n$. Now observe that
	\begin{align*}
		&(S^{-1}F)(g \circ f)\alpha \circ \phi_{f,g}^{A} \\
		&= \left(\prod_{i=0}^{n} F(g \circ f)\rho_i \circ F(g \circ f)(\psi_i)^{-1}\right) \circ \phi_{f,g}^{A} \\
		&= \left(\prod_{i=1}^{n} F(g \circ f)\rho_i \circ F(g \circ f)(\psi_i)^{-1}\right) \circ F(g \circ f)\rho_0 \circ F(g \circ f)(\psi_0)^{-1} \circ \phi_{f,g}^{A_0} \\
		&= \left(\prod_{i=1}^{n} F(g \circ f)\rho_i \circ F(g \circ f)(\psi_i)^{-1}\right) \circ F(g \circ f)\rho_0 \circ \phi_{f,g}^{C_0} \\
		&\circ \big(F(f) \circ F(g)\big)(\psi_0)^{-1} \\
		&= \left(\prod_{i=1}^{n} F(g \circ f)\rho_i \circ F(g \circ f)(\psi_i)^{-1}\right) \circ \phi_{f,g}^{A_1}\circ \big(F(f) \circ F(g)\big)\rho_0 \\
		&\circ \big(F(f \circ g)\big)(\psi_0)^{-1} \\ 
		&= \phi_{f,g}^{A_{n+1}} \circ \prod_{i=0}^{n}\big(F(f) \circ F(g)\big)\rho_i \circ \big(F(f) \circ F(g)\big)(\psi_i)^{-1} \\
		&= \phi_{f,g}^{A_{n+1}} \circ \big(F(f) \circ F(g)\big)\alpha = \phi_{f,g}^{B} \circ \big((S^{-1}F)(f) \circ (S^{-1}F)(g)\big)\alpha.
	\end{align*}
	Thus $S^{-1}\phi_{f,g}$ is a natural isomorphism. That it satisfies all the required relations of a compositor follows in a straightforward but tedious calculation using the fact that $\phi_{f,g}$ is a compositor, proving that $S^{-1}F$ is a  pseudofunctor.
\end{proof}
\begin{proposition}\label{Prop: Section 3: Equivariant localization pseudonat}
	Assume $F$ and $S$ are as in the statement of Proposition \ref{Prop: Section 3: Equivariant Localization Pseudofunctor}. There is a pseudonatural transformation $\lambda_S:F \Rightarrow S^{-1}F$\index[notation]{LambdaS@$\lambda_S$} whose object functor $\quot{\lambda}{X}$ is given by
	\[
	\quot{\lambda}{X} := \lambda_{\quot{S}{X}}:F(X) \to (\quot{S}{X})^{-1}F(X).
	\]
\end{proposition}
\begin{proof}
	Since we have already defined the $\quot{\lambda}{X}$, it suffices to define the natural isomorphisms $\quot{\lambda}{f}$. For this recall that the morphism $(S^{-1}F)(f)$ is defined via the universal property of the localization making the triangle
	\[
	\xymatrix{
		F(Y) \ar[rr]^{\quot{\lambda}{X} \circ F(f)} \ar[dr]_{\quot{\lambda}{Y}} & & (\quot{S}{X})^{-1}F(X) \\
		& (\quot{S}{Y})^{-1}F(Y) \ar@{-->}[ur]_{\exists!(S^{-1}F)(f)}
	}
	\]
	commute. However, this says that the diagram
	\[
	\xymatrix{
		F(Y) \ar[rr]^-{F(f)} \ar[d]_{\quot{\lambda}{Y}} & & F(X) \ar[d]^{\quot{\lambda}{X}} \\
		(\quot{S}{Y})^{-1}F(Y) \ar[rr]_-{(S^{-1}F)(f)} & & (\quot{S}{X})^{-1}F(X)
	}
	\] 
	commutes on the nose, so we define $\quot{\lambda}{f}$ to simply be the equality natural transformation witnessing $\quot{\lambda}{X} \circ F(f) = (S^{-1}F)(f) \circ \quot{\lambda}{Y}$. The required identity on pasting diagrams then holds immediately based on the fact that the $\quot{\lambda}{f}$ are equalities and from how we defined our compositors in the proof of Proposition \ref{Prop: Section 3: Equivariant Localization Pseudofunctor}.
\end{proof}
\begin{Theorem}\label{Theorem: Section 3: Equivariant localization is a 2-colimit}
	Let $F$ and $E$ be  pseudofunctors on $X$ and let 
	\[
	S = \lbrace \quot{S}{X} \; | \; \quot{S}{X} \subseteq F(X)_1, X \in \Cscr_0 \rbrace
	\] 
	be a collection of morphism classes as in Proposition \ref{Prop: Section 3: Equivariant Localization Pseudofunctor}. Assume there is a pseudonatural transformation $\alpha:F \to E$ such that for all $X \in \Cscr$ the functors $\quot{\alpha}{X}$ make 
	\[
	\quot{\alpha}{X}(\quot{S}{X}) \subseteq \Iso(E(X)) 
	\]
	hold. Then there exists a unique pseudonatural transformation $\eta:S^{-1}F \to E$ for which the induced diagram
	\[
	\xymatrix{
		\PC(F) \ar[rr]^{\PC(\alpha)} \ar[dr]_{\underline{\lambda}_S} & & \PC(E) \\
		& \PC(S^{-1}F) \ar@{-->}[ur]_{\exists!\,\underline{\eta}}
	}
	\]
	\index[notation]{SinverseFGX@$S^{-1}\PC(F)$}commutes in $\Cat$.
\end{Theorem}
\begin{proof}
	We first construct the functor component of the pseudonatural transformation $\eta$. For this fix a $X \in \Cscr_0$ and define the functor $\quot{\eta}{X}:(S^{-1}F)(X) \to E(X)$ as the unique functor filling the commuting diagram:
	\[
	\xymatrix{
		F(X) \ar[rr]^{\quot{\alpha}{X}} \ar[dr]_{\quot{\lambda}{X}} & & E(X) \\
		& (S^{-1}F)(X) \ar@{-->}[ur]_{\exists!\quot{\eta}{X}}
	}
	\]
	To define the natural isomorphisms $\quot{\eta}{f}$, fix a morpshism $f \in \Cscr_1$ and write $f:X \to Y$. We observe that by assumption we have a $2$-cell
	\[
	\begin{tikzcd}
		F(Y) \ar[r, bend left = 30, ""{name = U}]{}{\quot{\alpha}{X} \circ F(f)} \ar[r, bend right = 30, swap, ""{name = L}]{}{E(f) \circ \quot{\alpha}{Y}} & E(X) \ar[Rightarrow, from = U, to = L, shorten <= 4pt, shorten >= 4pt]{}{\quot{\alpha}{f}}
	\end{tikzcd}
	\]
	On one hand, note that the upper arrow in the $2$-cell factors in the strictly commuting square
	\[
	\xymatrix{
		F(X) \ar[r]^-{\quot{\lambda}{X}} & (S^{-1}F)(X) \ar[d]^{\quot{\eta}{X}} \\
		F(Y) \ar[u]^{F(f)} \ar[r]_-{\quot{\alpha}{X} \circ F(f)} & E(X)
	}
	\]
	while on the other hand the bottom arrow in the $2$-cell factors in the strictly commuting square:
	\[
	\xymatrix{
		F(Y) \ar[r]^-{E(f) \circ \quot{\lambda}{Y}} \ar[d]_{\quot{\lambda}{Y}} & E(X) \\
		(S^{-1}F)(Y) \ar[r]_-{\quot{\eta}{Y}} & E(Y) \ar[u]_{E(f)}
	}
	\]
	Combining these gives a pasting diagram:
	\[
	\begin{tikzcd}
		F(X) \ar[r]{}{\quot{\lambda}{X}} & (S^{-1}F)(X) \ar[d]{}{\quot{\eta}{X}^{\prime}} \\
		F(Y) \ar[u]{}{F(f)} \ar[r, bend left = 30, ""{name = U}]{}{\quot{\alpha}{X} \circ F(f)} \ar[r, swap, bend right = 30, ""{name = L}]{}{E(f) \circ \quot{\lambda}{Y}} \ar[d, swap]{}{\quot{\lambda}{Y}} & E(X) \\
		(S^{-1}F)(Y) \ar[r, swap]{}{\quot{\eta}{Y}} & E(Y) \ar[u, swap]{}{E(f)} \ar[from = U, to = L, Rightarrow, shorten <= 4pt, shorten >= 4pt]{}{\quot{\alpha}{f}}
	\end{tikzcd}
	\]
	It is routine to check that since the square
	\[
	\xymatrix{
		F(Y) \ar[r]^-{\quot{\lambda}{X}} \ar[d]_{F(f)} & (S^{-1}F)(Y) \ar[d]^{(S^{-1}F)(f)} \\
		F(X) \ar[r]_-{\quot{\lambda}{X}} & (S^{-1}F)(X)
	}
	\]
	commutes on the nose and the categories $(S^{-1}F)(X)$ and $(S^{-1}F)(Y)$ arise as localizations, there is a unique (from the fact that the localization arises as a $(2,1)$-colimit) natural isomorphism 
	\[
	\quot{\eta}{f}:\quot{\eta}{X} \circ (S^{-1}F)(f) \to E(f) \circ \quot{\eta}{Y}
	\]
	which fits into and factorizes the above pasting diagram as:
	\[
	\begin{tikzcd}
		F(X) \ar[r]{}{\quot{\lambda}{X}} & (S^{-1}F)(X) \ar[d]{}{\quot{\eta}{X}^{\prime}} \\
		F(Y) \ar[u]{}{F(f)}  \ar[d, swap]{}{\quot{\lambda}{Y}} & E(X) \\
		(S^{-1}F)(Y) \ar[uur, ""{name = M}]{}[description]{(S^{-1}F)(f)} \ar[r, swap]{}{\quot{\eta}{Y}} & E(Y) \ar[u, swap]{}{E(f)} \ar[from = M, to = 3-2, Rightarrow, swap, shorten <= 12pt, shorten >= 4pt]{}{\quot{\eta}{f}}
	\end{tikzcd}
	\]
	However, this is exactly the $2$-cell
	\[
	\begin{tikzcd}
		(S^{-1}F)(Y) \ar[r, bend left = 30, ""{name = U}]{}{\quot{\eta}{X} \circ (S^{-1}F)(f)} \ar[r, swap, bend right = 30, ""{name = L}]{}{E(f) \circ \quot{\eta}{Y}} & E(X) \ar[from = U, to = L, Rightarrow, shorten <= 4pt, shorten >= 4pt]{}{\quot{\eta}{f}}
	\end{tikzcd}
	\]
	that we will prove gives the morphism assignment of a pseudonatural transformation.
	
	We now prove the pseudonaturality of $\eta$. For this fix a composable pair $X \xrightarrow{f} Y \xrightarrow{g} Z$ in $\Cscr$. We begin by considering the pasting diagram:
	\[
	\begin{tikzcd}
		& (S^{-1}F)(Y) \ar[d, Rightarrow, shorten <= 4pt, shorten >= 4pt]{}{\quot{\phi_{f,g}}{S^{-1}F}} \ar[dr]{}{(S^{-1}F)(f)} \\
		(S^{-1}F)(Z) \ar[ur]{}{(S^{-1}F)(g)} \ar[d, swap]{}{\quot{\eta}{Z}} \ar[rr, ""{name = U}]{}[description]{(S^{-1}F)(g \circ f)} & {} & (S^{-1}F)(X) \ar[d]{}{\quot{\eta}{X}} \\
		E(Z) \ar[rr, swap, ""{name = L}]{}{E(g \circ f)} & & E(X) \ar[from = U, to = L, Rightarrow, shorten <= 4pt, shorten >= 4pt]{}{\quot{\eta}{g \circ f}}
	\end{tikzcd}
	\]
	A routine, but extremely tedious, check using that $\quot{\phi_{f,g}}{S^{-1}F}^{A} = \quot{\phi_{f,g}}{F}^A$ for all $A \in (S^{-1}F)(Z)_0 = F(Z)_0$, the fact that $\quot{\lambda}{f}$ is an identity transformation, and the universal property that induced the $\quot{\eta}{f}$ allows us to translate the above diagram through the pseudonatruality of $\alpha$ (much like how we proved $\quot{\phi_{f,g}}{F}^{A} = \quot{\phi_{f,g}}{S^{-1}F}^{A}$), and then move back to the localizations to derive that the above pasting diagram is equivalent to the pasting diagram below:
	\[
	\begin{tikzcd}
		(S^{-1}F)(Z)  \ar[rr, ""{name = LeftT}]{}{(S^{-1}F)(g)} \ar[d, swap]{}{\quot{\eta}{Z}} & & (S^{-1}F)(Y) \ar[rr, ""{name = MidT}]{}{(S^{-1}F)(f)} \ar[d]{}{\quot{\eta}{Y}} & & (S^{-1}F)(X) \ar[d]{}{\quot{\eta}{X}} \\
		E(Z) \ar[rr, swap, ""{name = LeftB}]{}{E(g)}, \ar[rrrr, bend right = 30, swap, ""{name = BB}]{}{E(g \circ f)} & &  E(Y) \ar[rr, swap, ""{name = MidB}]{}{E(f)} & & E(X) \ar[from = LeftT, to = LeftB, Rightarrow, shorten >= 4pt, shorten <= 4pt]{}{\quot{\eta}{g}} \ar[from = MidT, to = MidB, Rightarrow, shorten >= 4pt, shorten <= 4pt]{}{\quot{\eta}{f}} \ar[from = 2-3, to = BB, Rightarrow, shorten <= 4pt, shorten >= 4pt]{}{\quot{\phi_{f,g}}{E}}
	\end{tikzcd}
	\]
	Thus we have that $\eta$ is a pseudonatural transformation. Finally, we note that $\eta$ is even the unique pseudonatural transformation making
	\[
	\xymatrix{
		F \ar[rr]^{\alpha} \ar[dr]_{\lambda} & & E \\
		& S^{-1}F \ar@{-->}[ur]_{\exists!\eta}
	}
	\] 
	commute, as both $\quot{\eta}{X}$ and $\quot{\eta}{f}$ are induced by the universal properites of $1$-categories. This concludes the proof of the theorem.
\end{proof}
\begin{definition}\label{Defn: Equivariant Localization}\index[terminology]{Pseudocone Localization}
	A pseudocone localization of a category $\PC(F)$ is a pair $(\PC(S^{-1}F),\underline{\lambda}_S)$ where $S \subseteq \PC(F)_1$ is a collection of morphisms as in Theorem \ref{Theorem: Section 3: Equivariant localization is a 2-colimit}.
\end{definition}

\begin{example}\label{Example: Section 3: Equivariant localizaiton is not the localization of equivariant category}
	This example shows that if we na{\"i}vely localize an equivariant category at a class of morphisms, we can lose the equivariant nature of the category. This means in particular that equivariant localization is a delicate process; when localizing and preserving equivariant data, we need to make sure that we localize in a way that respects the descent data we have at hand.
	
	Fix a variety $X$ over $\Spec K$, let $G = \Spec K$ be the trivial group, and let $R$ be a ring for which $a \in R$ is not a right zero divisor and $\afrak = Ra$ is a left ideal for which $\RMod(R/\afrak, R) = 0$ and $R/\afrak$ is a nonzero simple $R$-module (a specific example is $R = \Z_p$ and $a = p$). Now consider the  pseudofunctor $F:\Cscr^{\op} \to \fCat$ given by $F(X) = \Ch(\RMod)$ and $F(f) = \id_{\Ch(\RMod)}$. Consider now the complexes $A^{\bullet}$ and $B^{\bullet}$ in $\Ch(\RMod)$, where $A^{\bullet}$ is defined as
	\[
	\xymatrix{
		\cdots \ar[r] & 0 \ar[r] & R \ar[r]^{\rho_a} & R \ar[r] & 0 \ar[r] & \cdots
	}
	\]
	where $\rho_a$ is the right multiplication by $a$ map and $R$ appears in degrees $0$ and $1$, and $B^{\bullet}$ is defined by
	\[
	\xymatrix{
		\cdots \ar[r] & 0 \ar[r] & R/\afrak \ar[r] & 0 \ar[r] & \cdots
	}
	\]
	where $R/\afrak$ appears in degree $0$. It is routine to check that
	\[
	H^{n}(A^{\bullet}) = \begin{cases}
		0 & \text{if}\, n \ne 0; \\
		R/\afrak & \text{if}\, n = 0;
	\end{cases}
	\]
	and that
	\[
	H^{n}(B^{\bullet}) = \begin{cases}
		0 & \text{if}\, n \ne 0; \\
		R/\afrak & \text{if}\, n = 0.
	\end{cases}
	\]
	Because $R/\afrak$ is nonzero and simple, any nonzero morphism $\varphi^{\bullet}:A^{\bullet} \to B^{\bullet}$ extends to a quasi-isomorphism $H^{\ast}(\phi):H^{\ast}(A^{\bullet}) \xrightarrow{\cong} H^{\ast}(B^{\bullet})$. We take the morphism $\varphi^{\bullet}:A^{\bullet} \to B^{\bullet}$
	\[
	\varphi^{n} = \begin{cases}
		0 & \text{if}\, n \ne 0; \\
		\pi_{\afrak} & \text{if}\, n = 0
	\end{cases}
	\]
	as our nonzero map, which is a quasi-isomorphism $A^{\bullet} \to B^{\bullet}$.
	
	Now define $\quot{Q}{X} \subseteq F(X)_1$ by $\quot{Q}{X} := \lbrace f \in F(X)_1 \; | \; f \, \text{is\, a\, quasi-isomorphism}\rbrace$ and set $Q := \lbrace P \in \PC(F) \; | \; \forall\, X \in \Cscr_0 \forall\, \quot{\rho}{X} \in P, \quot{\rho}{X} \in \rbrace$. Write $(Q^{-1}F)_G(X)$ for the equivariant localization implied by Theorem \ref{Theorem: Section 3: Equivariant localization is a 2-colimit} and Definition \ref{Defn: Equivariant Localization} and write $Q^{-1}(\PC(F))$ for the localization of the category $\PC(F)$ at $Q$.
	
	Our goal now is to show that $(Q^{-1}F)_G(X) \ne Q^{-1}(\PC(F))$. To show this we need only show that $(Q^{-1}F)_G(X)_0 \ne Q^{-1}(\PC(F))_0$. To this end consider the object $(A,T_A) \in Q^{-1}\PC(F)$ given by saying that the $\quot{A}{X} \in A$ take the form
	\[
	\quot{A}{X} = \begin{cases}
		A^{\bullet} & \text{if}\, X \ne G; \\
		B^{\bullet} & \text{if}\, X = G;
	\end{cases}
	\]
	and defining the $\tau_f^A \in T_A$ by
	\[
	\tau_f^A = \begin{cases}
		\varphi^{\bullet}:A^{\bullet} \to B^{\bullet} & \text{if}\,\exists\,X \in \Cscr_0. f = \nu_{X}:X \to \Spec K; \\
		\id_{A^{\bullet}}:A^{\bullet} \to A^{\bullet} & \text{else}.
	\end{cases}
	\]
	It is a trivial verification that $(A,T_A) \in (Q^{-1}F)_G(X)_0$. However, note that $(A,T_A) \notin \PC(F)_0$ because $\varphi^{\bullet}$ is a quasi-isomorphism (and hence an isomorphism in $Q^{-1}\Ch(\RMod)$ which is not an isomorphism in $\Ch(\RMod)$). As such we have $\PC(Q^{-1}F)_G(X)_0 \ne Q^{-1}(\PC(F))_0$ and hence $\PC(Q^{-1}F) \ne Q^{-1}(\PC(F))$.
\end{example}
\begin{remark}
	It remains to be characterized when, if the pseudocone localization $\PC(S^{-1}F)$ is induced by categories where each localization  $(\quot{S}{X})^{-1}F(X)$ is induced by a calculus of fractions, if $\PC(S^{-1}F)$ is itself an Ore localization of $\PC(F)$ in any way. Determining this seems difficult in practice.
\end{remark}
\begin{proposition}\label{Prop: Section Functors: EDC is pc localization}
Let $G$ be a smooth algebraic group and let $X$ be a left $G$-variety or let $G$ be a topological group and let $X$ be a left $G$-space. Then in either case we have that
\[
D_G^b(X) \simeq \PC\left(\operatorname{QIso}^{-1}\Ch^b\big(G \backslash (-)\big)\right).
\]
\end{proposition}
\begin{proof}
Recall that in either case we have that $D_G^b(X) = \PC\left(D_c^b\big(G \backslash (-)\big)\right)$ by Section \ref{Subsection: Equivariant Derived Cat of Var} in the geometric case and $D_G^b(X) = \PC\left(D^b\big(G \backslash (-)\big)\right)$ in the topological case by Section \ref{Subsection: EDC of a Space}. Now recall that in both the topological and geometric cases we have unique isomorphisms
\[
D_c^b\left(\quot{X}{\Gamma}\right) \cong \operatorname{QIso}^{-1}\Ch_c^b\left(\quot{X}{\Gamma}\right)
\]
which vary pseudonaturally in $\SfResl_G(X)^{\op}$ due to the universal property of localizations of categories\footnote{Depending on how you define the derived category, these isomorphisms could be actually be equalities; however, our approach here is to be representation-agnostic.}. But it then follows the corresponding pseudofunctors $D_c^b\left(G \backslash (-)\right)$ and $\operatorname{QIso}^{-1}\Ch_c^b\left(G \backslash (-)\right)$ are pseudonaturally equivalent and hence by Corollary \ref{Cor: Section 3: Fibre-wise equivalences are equivalences of equivariant categories} we have that $D_G^b(X) \simeq \PC\left(\operatorname{QIso}^{-1}\Ch_c^b(G \backslash(-))\right)$ in the geometric case. The topological case follows mutatis mutandis.
\end{proof}

\newpage

\section{Pseudocone Functors: Change of Domain}\label{Section: Section 3: Change of Space}\label{Subsection: Change of Domain}
We now discuss when we can give functors which change the domain of a pseudocone category, i.e., functors $\alpha:\PC(F) \to \PC(E)$ where $F$ and $E$ are pseudofunctors in a cospan
\[
\begin{tikzcd}
	\Cscr^{\op} \ar[r]{}{F} & \fCat & \Dscr^{\op} \ar[l, swap]{}{E}
\end{tikzcd}
\]
for some $1$-categories $\Cscr, \Dscr$. The section prior to this and the strict $2$-functor $\PC:\Bicat(\Cscr^{\op},\fCat) \to \fCat$ suggest that we should look for some sort of pseudonatural transformation, or something close to one, between $E$ and $F$ in order to define a functor $\PC(F) \to \PC(E)$. We will see that this can be done, save with some technical alterations to account for the fact that $E$ and $F$ have distinct domains and so need not be directly comparable automatically. This section is less abstract and general than its predecessor, but is arguably more important than because it allows us to define the functors between equivariant categories which we think of as ``Change of Space'' equivariant functors (cf.\@ Corollaries \ref{Cor: Pseudocone Functors: Existence of equivariant pullback  for schemes}, \ref{Cor: Pseudocone Functors: Existence of equivariant pullback  for spaces}, \ref{Cor: Pseudocone Functors: Existence of equivariant pullback  for schemes}, \ref{Cor: Pseudocone Functors: Existence of equivariant pullback  for spaces}, for instance) in the geometric and topological situation. For instance, the existence of the functors
\[
h^{\ast}:D_G^b(Y) \to D_G^b(X), \qquad Rh_{\ast}:D_G^b(X) \to D_G^b(Y)
\]
and
\[
h^{!}:D_G^b(Y) \to D_G^b(X), \qquad Rh_!:D_G^b(X) \to D_G^b(Y)
\]
when there is a $G$-equivariant morphism of varieties $h:X \to Y$ (or a $G$-equivariant morphism of topological spaces) follow from the arguments in this section (cf.\@ Corollaries \ref{Cor: Pseudocone Functors: Existence of equivariant pullback  for schemes}, \ref{Cor: Pseudocone Functors: Pushforward functors for equivariant maps for scheme sheaves}, \ref{Cor: Pseudocone Functors: Exceptional Pushforward functors for equivariant maps for scheme sheaves}, \ref{Cor: Pseudocone Functors: Topological EDC has equivariant pushpull}, \ref{Cor: Pseudocone Functors: Topological EDC has equivariant exceptional pushpull}). The definition and conditions we described have appeared briefly and with minimal detail in \cite{DoretteMe}, so we present a systematic discussion and explain them in detail here.

To define what are essentially pseudonatural transformations from $F$ to $E$, even in spite of the fact that the domains of the pseudofunctors differ, we use a slight parametrization trick to move our pseudocones around. The basic idea is that if we have a functor $\gamma:\Cscr \to \Dscr$ we can view $\gamma$ in some sense as a path from $\Cscr$ to $\Dscr$ which can allow us to pick out a part of $\fCat$ in which the image of the pseudofunctors $F:\Cscr^{\op} \to \fCat$ and $E:\Dscr^{\op} \to \fCat$ interact and have a pseudonatural transformation moving between them. While this is ultimately a pseudonatural transformation between $F$ and the induced pre-composite pseudofunctor $E \circ \gamma^{\op}:\Cscr^{\op} \to \fCat$, I am of the opinion that the intuition just described is not just helpful but also that the results in this section are useful enough for practical situations that they warrant explicit statements.

\begin{definition}\label{Defn: Pseudocone section: Coone translations}
	Let $\Cscr, \Dscr$ be categories with $\gamma:\Cscr \to \Dscr$ a functor and $F:\Cscr^{\op} \to \fCat$ and $E:\Dscr^{\op} \to \fCat$ pseudofunctors. We say that $F$ admits pseudocone translations\index[terminology]{Pseudocone Translations} of shape $h = (\quot{h}{X}, \quot{h}{f})$ along $\gamma$ to $E$ if $h$ is a pseudonatural transformation:
	\[
	\begin{tikzcd}
		\Cscr^{\op} \ar[rr, bend left = 30, ""{name = U}]{}{F} \ar[rr, bend right = 30, swap, ""{name = D}]{}{E \circ \gamma^{\op}} & & \fCat \ar[from = U, to = D, Rightarrow, shorten <= 4pt, shorten >= 4pt]{}{h}
	\end{tikzcd}
	\]
	Similarly, we say that $E$ has an optranslation\index[terminology]{Pseudocone Optranslation} along $\gamma$ to $F$ of shape $h = (\quot{h}{X}, \quot{h}{f})$ if there is a pseudonatural transformation:
	\[
	\begin{tikzcd}
		\Cscr^{\op} \ar[rr, bend right = 30, swap, ""{name = U}]{}{F} \ar[rr, bend left = 30,""{name = D}]{}{E \circ \gamma^{\op}} & & \fCat \ar[from = D, to = U, Rightarrow, shorten <= 4pt, shorten >= 4pt]{}{h}
	\end{tikzcd}
	\]
	When no confusion will arise we will often be brief and simply say that $F$ has a translation to $E$ (respectively $E$ has a translation to $F$) and leave $h$ and $\gamma$ implicit.
\end{definition}
\begin{remark}
	Let $\Cscr, \Dscr$ be categories with $\gamma:\Cscr \to \Dscr$ a functor and $F:\Cscr^{\op} \to \fCat$ and $E:\Dscr^{\op} \to \fCat$ pseudofunctors. We now unwrap the definition of what it means to have translations along $\gamma$ so as to have an explicit description of the information at hand. That $F$ admits pseudocone translations of shape $h = (\quot{h}{X}, \quot{h}{f})$ along $\gamma$ to $E$ means the following hold:
	\begin{itemize}
		\item For every $X \in \Cscr_0$ there is a functor $\quot{h}{X}:F(X) \to E(\gamma X)$;
		\item For every $f \in \Cscr_1, f:X \to Y$, there is a natural isomorphism:
		\[
		\begin{tikzcd}
			F(Y) \ar[r, ""{name = U}]{}{F(f)} \ar[d, swap]{}{\quot{h}{Y}} & F(X) \ar[d]{}{\quot{h}{X}} \\
			E(\gamma Y) \ar[r, swap, ""{name = D}]{}{E(\gamma f)} & E(\gamma X) \ar[from = U, to = D, Rightarrow, shorten <= 4pt, shorten >= 4pt]{}{\quot{h}{f}} \ar[from = U, to = D, Rightarrow, shorten <= 4pt, shorten >= 4pt, swap]{}{\cong}
		\end{tikzcd}
		\]
		\item For any composable arrows $X \xrightarrow{f} Y \xrightarrow{g} Z$ the pasting diagram
		\[
		\begin{tikzcd}
			F(Z) \ar[r, ""{name = UpLeft}]{}{F(g)} \ar[d, swap]{}{\quot{h}{Z}} & F(Y) \ar[r, ""{name = UpRight}]{}{F(f)} \ar[d]{}[description]{\quot{h}{Y}} & F(X) \ar[d]{}{\quot{h}{X}} \\
			E(\gamma Z) \ar[r, swap, ""{name = DownLeft}]{}{E(\gamma g)} \ar[rr, swap, bend right =40, ""{name = Downest}]{}{E(\gamma(g \circ f))} & E(\gamma Y) \ar[r, swap, ""{name = DownRight}]{}{E(\gamma f)} & E(\gamma X) \ar[from = UpLeft, to = DownLeft, Rightarrow, shorten <= 4pt, shorten >= 4pt]{}{\quot{h}{g}} \ar[from = UpRight, to = DownRight, Rightarrow, shorten <= 4pt, shorten >= 4pt, swap]{}{\quot{h}{f}} \ar[from = 2-2, to = Downest, Rightarrow, shorten <= 4pt, shorten >= 4pt]{}{\quot{\phi_{f,g}}{E}}
		\end{tikzcd}
		\]
		is equal to the pasting diagram:
		\[
		\begin{tikzcd}
			& F(Y) \ar[dr]{}{F(f)} \\
			F(Z) \ar[ur]{}{F(g)} \ar[d, swap]{}{\quot{h}{Z}} \ar[rr, ""{name = U}]{}[description]{F(g \circ f)} & & F(X) \ar[d]{}{\quot{h}{X}}	 \\
			E(\gamma Z) \ar[rr, swap, ""{name = D}]{}{E(\gamma(g \circ f))} & & E(\gamma X) \ar[from = U, to = D, Rightarrow, shorten >= 4pt, shorten <= 4pt]{}{\quot{h}{g \circ f}} \ar[from = 1-2, to = U, Rightarrow, shorten <= 4pt, shorten >= 4pt]{}{\quot{\phi_{f,g}}{F}}
		\end{tikzcd}
		\]
	\end{itemize}
	An explicit description of what it means for $E$ to have optranslations along $\gamma$ to $F$ of shape $h$ is omitted.
\end{remark}
\begin{definition}\label{Defn: Pseudocone Functors: Descent Modification}
	Let $\gamma:\Cscr \to \Dscr$ be a functor such that if $F:\Cscr^{\op} \to \fCat$ and $E:\Dscr^{\op} \to \fCat$ are pseudofunctors then $F$ admits pseudocone translations of shape $\alpha = (\quot{\alpha}{X}, \quot{\alpha}{f})$ and $\beta = (\quot{\beta}{X}, \quot{\beta}{f})$ along $\gamma$ to $E$ or $E$ admits optranslations of shape $\alpha, \beta$ along $\gamma$. Then we say that there is a descent modification\index[terminology]{Descent Modification} $\rho:\alpha \Rrightarrow \beta$ if $\rho$ is a modification fitting in the diagram:
	\[
	\begin{tikzcd}
		\Cscr^{\op} \ar[rrr, bend left = 40, ""{name = Up}]{}{F} \ar[rrr, bend right = 40, swap, ""{name = Down}]{}{E \circ \gamma^{\op}} & & & \fCat \ar[from = Up, to = Down, Rightarrow, shorten <= 4pt, shorten >= 4pt, bend left = 40, ""{name = Left}]{}{\beta} \ar[from = Up, to = Down, Rightarrow, shorten <= 4pt, shorten >= 4pt, bend right = 40, swap, ""{name = Right}]{}{\alpha} \ar[from = Right, to = Left, symbol = \underset{\rho}{\Rrightarrow}]
	\end{tikzcd}
	\]
	or
	\[
		\begin{tikzcd}
		\Cscr^{\op} \ar[rrr, bend left = 40, ""{name = Up}]{}{E \circ \gamma^{\op}} \ar[rrr, bend right = 40, swap, ""{name = Down}]{}{F} & & & \fCat \ar[from = Up, to = Down, Rightarrow, shorten <= 4pt, shorten >= 4pt, bend left = 40, ""{name = Left}]{}{\beta} \ar[from = Up, to = Down, Rightarrow, shorten <= 4pt, shorten >= 4pt, bend right = 40, swap, ""{name = Right}]{}{\alpha} \ar[from = Right, to = Left, symbol = \underset{\rho}{\Rrightarrow}]
	\end{tikzcd}
	\]
\end{definition}
\begin{remark}
	As before, let us unwrap what this means explicitly. Fix categories, functors, and pseudofunctors as in Definition \ref{Defn: Pseudocone Functors: Descent Modification} when $F$ admits pseudocone translations of shape $\alpha, \beta$ along $\gamma$ to $E$. Then $\rho$ is determined by a collection of natural transformations
	\[
	\rho  = \left\lbrace \quot{\rho}{X}:\quot{\alpha}{X} \Rightarrow \quot{\beta}{X} \; | \; X \in \Cscr_0 \right\rbrace
	\]
	such that for any morphism $f:X \to Y$ in $\Cscr$ the diagram of functors and natural transformations
	\[
	\begin{tikzcd}
		\quot{\alpha}{X} \circ F(f) \ar[d, swap]{}{\quot{\rho}{X} \ast F(f)} \ar[r]{}{\quot{\alpha}{f}} & E(\gamma f) \circ \quot{\alpha}{Y} \ar[d]{}{E(\gamma f) \ast \quot{\rho}{Y}} \\
		\quot{\beta}{X} \circ F(f) \ar[r, swap]{}{\quot{\beta}{f}} & E(\gamma f) \circ \quot{\beta}{Y}
	\end{tikzcd}
	\]
	commutes. The case in which $E$ admits pseudocone optranslations of shape $\alpha, \beta$ along $\gamma$ to $F$ is similar but omitted.
\end{remark}

Because of the descriptions above we immediately get that if $F$ admits cone translations to $E$, then there is a functor $\PC(F) \to \PC(E)$ as an application of Theorem \ref{Thm: Functor Section: Psuedonatural trans are pseudocone functors}; similarly if there are descent modifications from $h$ to $k$ then there is a natural transformation between. While this is a special case of the material we presented in the subsection prior, because of its application to the various equivariant pushforward/pullback, cohomology, and various $t$-structure functors we present (cf. Chapter \ref{Section: Triangles are go}) in detail here for future use.

\begin{proposition}\label{Prop: Pseudocone Functors: TRanslation}
	Let $F:\Cscr{\op} \to \fCat$ and $E:\Dscr^{\op} \to \fCat$ be pseudofunctors and let $\gamma:\Cscr\to \Dscr$ be a functor. Then if $F$ admits pseudocone translations of shape $h$ along $\gamma$ to $E$ there is a functor $\ul{h}:\PC(F) \to \PC(E \circ \gamma^{\op})$ and if $E$ admits optranslations of shape $h$ along $\gamma$ to $F$ then there is a functor $\ul{h}:\PC(E \circ \gamma^{\op}) \to \PC(F)$. Similarly if $F$ admits translations of shape $\alpha$ and $\beta$ along $\gamma$ to $E$ and there is a descent modification $\rho:\alpha \Rrightarrow \beta$ then there is a natural transformation $\ul{\rho}:\ul{\alpha} \Rightarrow \ul{\beta}$ (respecitvely if $E$ admits optranslations of shape $\alpha, \beta$ along $\gamma$ to $F$ and there is a descent modification $\rho:\alpha \Rrightarrow \beta$ then there is a natural transformation $\ul{\rho}:\ul{\alpha} \Rightarrow \ul{\beta}$).
\end{proposition}
\begin{proof}
	Apply Lemma \ref{Lemma: Modifications give equivariant natural transformations}.
\end{proof}

We now show a useful proposition which determines when we can get pseudocone translation functors as induced from a natural transformation between subcategory inclusions and a pseudofunctor defined on the larger category. Note that as a corollary we can deduce the existence of the equivariant pullback functors
\[
h^{\ast}:D_G^b(Y) \to D_G^b(X), \qquad h^{\ast}:D_G^b(Y; \overline{\Q}_{\ell}) \to D_G^b(X;\overline{\Q}_{\ell})
\]
whenever $h:X \to Y$ is a $G$-equivariant\footnote{The astute reader will note that I have not specified whether or not $G$ is a topological group or an algebraic group or whether or not $X$ and $Y$ are schemes or spaces. In truth it does not matter provided the formalisms are all chosen appropriately: if $G$ is an algebraic group then $X$ and $Y$ are schemes and if $G$ is a topological group then $X$ and $Y$ are topological spaces.} morphism.

\begin{Theorem}\label{Thm: Pseudocone Functors: Pullback induced by fibre functors in pseudofucntor}
	Let $\Cscr \xrightarrow{\varphi} \Escr \xleftarrow{\psi} \Dscr$ be a cospan of categories together with a functor $\gamma:\Cscr \to \Dscr$. If there is a $2$-cell
	\[
	\begin{tikzcd}
		\Cscr \ar[rr, ""{name = U}]{}{\varphi} \ar[dr, swap]{}{\gamma} & & \Escr \\
		& \Dscr \ar[ur, swap]{}{\psi} \ar[from = U, to = 2-2, Rightarrow, shorten <= 4pt, shorten >= 4pt]{}{\alpha}
	\end{tikzcd}
	\]
	and a pseudofunctor $F:\Escr^{\op} \to \fCat$ then $E \circ \varphi^{\op}$ admits pseudocone translations along $\gamma^{\op}$ to $F \circ \psi^{\op}$ of shape $h = (\quot{h}{X}, \quot{h}{f})$ where, for all $X \in \Cscr_0$,
	\[
	\quot{h}{X} := F(\alpha_X):F(\gamma X) \to F(X)
	\]
	and for all $f:X \to Y$ in $\Cscr$,
	\[
	\quot{h}{f} := \phi_{\varphi f,\alpha_Y}^{-1} \circ \phi_{\alpha_X,(\psi \circ \gamma) f}.
	\]
\end{Theorem}
\begin{proof}
	Note that upon reversing all arrows in each underlying category in the $2$-cell
	\[
	\begin{tikzcd}
		\Cscr \ar[rr, ""{name = U}]{}{\varphi} \ar[dr, swap]{}{\gamma} & & \Escr \\
		& \Dscr \ar[ur, swap]{}{\psi} \ar[from = U, to = 2-2, Rightarrow, shorten <= 4pt, shorten >= 4pt]{}{\alpha}
	\end{tikzcd}
	\]
	we get a corresponding $2$-cell:
	\[
	\begin{tikzcd}
		& \Dscr^{\op} \ar[dr]{}{\psi^{\op}} \\
		\Cscr^{\op} \ar[ur]{}{\gamma^{\op}} \ar[rr, swap, ""{name = D}]{}{\varphi^{\op}} & & \Escr^{\op} \ar[to = D, from = 1-2, Rightarrow, shorten <= 4pt, shorten >= 4pt]{}{\alpha^{\op}}
	\end{tikzcd}
	\]
	Taking the pseudofunctor $F:\Escr^{\op} \to \fCat$ and applying the post-composition pseudofunctor
	\[
	F \circ (-):\Bicat(\Cscr^{\op},\Escr^{\op}) \to \Bicat(\Cscr^{\op},\fCat)
	\]
	then sends the diagram above to a $2$-cell:
	\[
	\begin{tikzcd}
		\Cscr^{\op} \ar[rr, bend left = 40, ""{name = U}]{}{F \circ \psi^{\op} \circ \gamma^{\op}} \ar[rr, bend right = 40, swap, ""{name = D}]{}{F \circ \varphi^{\op}} & & \fCat \ar[from = U, to = D, Rightarrow, shorten <= 4pt, shorten >= 4pt]{}{F \ast \alpha^{\op}}
	\end{tikzcd}
	\]
	The object functors of the pseudonatural transformation $F \ast \alpha^{\op}$ are
	\[
	\quot{(F \ast \alpha^{\op})}{X} := F(\alpha_X)
	\]
	and the morphism natural transformations are
	\[
	\quot{(F \ast \alpha^{\op})}{f} := \phi_{\alpha_X, (\psi \circ \gamma)f}^{-1} \circ \phi_{\varphi f, \alpha_Y}
	\]
	for objects $X$ in $\Cscr$ and morphisms $f:X \to Y$ in $\Cscr$. Taking $h := F \ast \alpha^{\op}$ then proves the theorem.
\end{proof}
We now present two corollaries/examples of the theorem above in order to show how the equivariant pullback functors exist as a formal consequence of our pseudocone formalism.
\begin{corollary}\label{Cor: Pseudocone Functors: Existence of equivariant pullback  for schemes}
	Let $G$ be a smooth algebraic group and let $X, Y$ be left $G$-varieties with $h:X \to Y$ a $G$-equivariant morphism. Then there is an equivariant pullback functor
	\[
	h^{\ast}:D_G^b(Y;\overline{\Q}_{\ell}) \to D_G^b(X;\overline{\Q}_{\ell})
	\]
	and similarly for the non-$\ell$-adic derived categories.
\end{corollary}
\begin{proof}
	The categories $\SfResl_G(X)$ and $\SfResl_G(Y)$ are isomorphic via the maps $\Gamma \times X \mapsto \Gamma \times Y$ and $\Gamma \times Y \mapsto \Gamma \times X$ for all $\Sf(G)$-varieties $\Gamma$\footnote{Note that we are not saying that there is an isomorphism from $\Gamma \times X$ to $\Gamma \times Y$ in the category of $G$-varieties but instead simply that the assignment above induces an isomorphism of object collections $\SfResl_G(X)_0 \cong \SfResl_G(Y)_0$ which extends correspondingly to morphism classes.}. Furthermore, both admit $G$-quotient functors to $\Var_{/K}$ (one named $\quo_{X}:\SfResl_G(X) \to \Var_{/K}$ an the other named $\quo_{Y}:\SfResl_G(Y) \to \Var_{/K}$) of the forms 
	\[
	\quo_{X}(\Gamma \times X) := G \backslash (\Gamma \times X)
	\]
	and
	\[
	\quo_{Y}(\Gamma \times Y) := G \backslash (\Gamma \times Y),
	\]
	respectively (cf.\@ Proposition \ref{Prop: Section 2.1: The quotient prop}; note also the definition of each functor on morphisms is similar and omitted). There is a $2$-cell
	\[
	\begin{tikzcd}
		\SfResl_G(X) \ar[dr, swap]{}{\cong} \ar[rr, ""{name = U}]{}{\quo_{(-)\times X}} & & \Var_{/K} \\
		& \SfResl_G(Y) \ar[ur, swap]{}{\quo_{(-)\times Y}} \ar[from = U, to = 2-2, Rightarrow, shorten <= 4pt, shorten >= 4pt]{}{\varphi}
	\end{tikzcd}
	\]
	where the structure map
	\[
	\varphi_{\Gamma \times X}:G \backslash(\Gamma \times X) \to G \backslash (\Gamma \times Y)
	\]
	is given by $\varphi_{\Gamma \times X} := \overline{\id_{\Gamma} \times h}$, where this is the unique map making the diagram
	\[
	\begin{tikzcd}
		G \backslash(\Gamma \times X) \ar[d, swap]{}{\quo_{\Gamma \times X}} \ar[rr]{}{\id_{\Gamma} \times h} & & \Gamma \times Y \ar[d]{}{\quo_{\Gamma \times Y}} \\
		G \backslash(\Gamma \times X) \ar[rr, dashed, swap]{}{\overline{\id_{\Gamma} \times h}} & & G \backslash (\Gamma \times Y)
	\end{tikzcd}
	\]
	commute. That this is natural follows from the commutativity of the diagram
	\[
	\begin{tikzcd}
		& \Gamma \times Y \ar[dd, swap, near end]{}{\quo_{\Gamma \times Y}} \ar[rr]{}{f \times \id_Y} & & \Gamma^{\prime} \times Y \ar[dd]{}{\quo_{\Gamma^{\prime} \times Y}} \\
		\Gamma \times X \ar[dd, swap]{}{\quo_{\Gamma \times X}} \ar[rr, crossing over, near end]{}{f \times \id_X} \ar[ur]{}{\id_{\Gamma} \times h} & & \Gamma^{\prime} \times X \ar[ur,]{}{\id_{\Gamma^{\prime}}\times h} \\
		& G \backslash (\Gamma \times Y) \ar[rr, near end]{}{\overline{f \times \id_Y}} & & G \backslash (\Gamma^{\prime} \times Y) \\
		G \backslash (\Gamma \times X) \ar[ur]{}{\overline{\id_{\Gamma} \times h}} \ar[rr, swap]{}{\overline{f \times \id_X}} & & G \backslash (\Gamma^{\prime} \times X) \ar[ur, swap]{}{\overline{\id_{\Gamma^{\prime}} \times h}} \ar[from = 2-3, to = 4-3, crossing over, near end]{}{\quo_{\Gamma^{\prime} \times X}}
	\end{tikzcd}
	\]
	for any morphism $f:\Gamma \to \Gamma^{\prime}$ in $\Sf(G)$. The corollary now follows by applying Theorem \ref{Thm: Pseudocone Functors: Pullback induced by fibre functors in pseudofucntor} to the pseudofunctor $D^b_c(-;\overline{\Q}_{\ell}):\Var_{/K}^{\op} \to \fCat$ given by sending a variety to its (constructible) $\ell$-adic derived category and by sending a morphism to its corresponding pullback functor. The non-$\ell$-adic version follows mutatis mutandis.
\end{proof}
\begin{corollary}\label{Cor: Pseudocone Functors: Existence of equivariant pullback  for spaces}
	Let $G$ be a topological group and let $X, Y$ be left $G$-spaces with $h:X \to Y$ a $G$-equivariant morphism. Then there is an equivariant pullback functor
	\[
	h^{\ast}:D_G^b(Y) \to D_G^b(X).
	\]
\end{corollary}
\begin{proof}
	This follows mutatis mutandis to the proof of Corollary \ref{Cor: Pseudocone Functors: Existence of equivariant pullback  for schemes}.
\end{proof}

Much of what remains to discuss in this section is simply an application of the material in Section \ref{Subsection: Change of Fibre}; however, we present it explicitly because of its importance in applications to scheme-theoretic and topological situations. In particular, constructing the pushforward $Rh_{\ast}$, the proper pushforward $Rh_{!}$, and the exceptional pullback $h^{!}$ rely on these constructions while the adunctions they satisfy (namely that $h^{\ast} \dashv Rh_{\ast}$ and $Rh_{!} \dashv h^{!}$) follow from what we describe below. We will also see later in Section \ref{Subsection: The Equivariance Isos} that the existence of an ``equivariance'' isomorphism satisfying the GIT cocycle condtion (cf.\@ Defintion \ref{Defn: GIT Cocycle}) exists for pseudocone categories $\PC(\overline{F} \circ \quo^{\op})$ for pseudofunctors $\overline{F}:\Var_{/K}^{\op} \to \fCat$ and $\quo:\SfResl_G(X) \to \Var_{/K}$ by applying the formalism we develop here, so it is worth having some technical but precise results displayed explicitly. First we present a proposition which states when the functors $\ul{\alpha}$ preserve various limits, colimits, and monoidal structures and then describe adjoints between Change of Domain functors.

\begin{proposition}\label{Prop: Pseudocone Functors: Translations preserve stuff}
	Let $F:\Cscr^{\op} \to \fCat$ and $E:\Dscr^{\op} \to \fCat$ be pseudofunctors and let $\gamma:\Cscr^{\op} \to \Dscr^{\op}$. Assume either:
	\begin{enumerate}
		\item[Case 1:] $F$ admits pseudocone translations along $\gamma$ of shape $h = (\quot{h}{X}, \quot{h}{f})$ to $E$;
		\item[Case 2:] $E$ has an optranslation along $\gamma$ to $F$ of shape $k = (\quot{k}{X}, \quot{k}{f})$.
	\end{enumerate} 
	Then the following hold in Case 1:
	\begin{enumerate}
		\item If $I$ is an index category and each functor $\quot{h}{X}$ preserves all limits of shape $I$ then $\ul{h}$ does as well;
		\item If $I$ is an index category and each functor $\quot{h}{X}$ preserves all colimits of shape $I$ then $\ul{h}$ does as well;
		\item If $F$ and $E$ are (braided) monoidal and each of $\quot{h}{X}$ and $\quot{h}{f}_{\ast}$ are (braided) monoidal functors and transformations then $\ul{h}$ is (braided) monoidal as well;
	\end{enumerate}
	and analogously for Case 2.
\end{proposition}
\begin{proof}
	Claims (1) and (2) hold in both cases by virtue of Proposition \ref{Prop: Section 3: Change of fibre equivariant functor preserves limit of a specific shape}. The monoidal statement for Claim (3) follows in Cases 1 and 2 follows from Proposition \ref{Prop: Section 3: Monoidal pseudonats give rise to monoidal equivariant functors} while the braided monoidal statement follows by applying Corollary \ref{Cor: Pseudocone Functors: Braided monoidal pseudonats are braided monoidal functors}.
\end{proof}
As an immediate application of the above proposition we get the corollary below.
\begin{corollary}\label{Cor: Pseudocone Functors: When translations preserve all limits and colimits}
	Let $F:\Cscr^{\op} \to \fCat$ and $E:\Dscr^{\op} \to \fCat$ be pseudofunctors and let $\gamma:\Cscr^{\op} \to \Dscr^{\op}$. Assume either:
	\begin{enumerate}
		\item[Case 1:] $F$ admits pseudocone translations along $\gamma$ of shape $h = (\quot{h}{X}, \quot{h}{f})$ to $E$;
		\item[Case 2:] $E$ has an optranslation along $\gamma$ to $F$ of shape $k = (\quot{k}{X}, \quot{k}{f})$.
	\end{enumerate} 
	Then in both cases:
	\begin{itemize}
		\item If the pseudofunctors $F$ and $E$ are additive and every functor $\quot{h}{X}$ is additive, then the functor $\ul{h}$ is additive as well.
		\item If every functor $\quot{h}{X}$ is (left/right/two-sided) exact, so is $\ul{h}$.
		\item If every functor $\quot{h}{X}$ is continuous or cocontinuous, so is $\ul{h}$.
	\end{itemize}
\end{corollary}

We now discuss shortly adjunctions which arise between translation functors. As in Theorem \ref{Thm: Section 3: Gamma-wise adjoints lift to equivariant adjoints} it suffices in all cases to specify adjoint pairs of functors $\quot{h}{X} \dashv \quot{k}{X}$ which both vary pseudonaturally.

\begin{Theorem}\label{Thm: Pseudocone Functors: Adjoints between translations give adjoint functors}
	Let $F:\Cscr^{\op} \to \fCat$ and $E:\Dscr^{\op} \to \fCat$ be pseudofunctors and let $\gamma:\Cscr^{\op} \to \Dscr^{\op}$. Assume that $F$ admits pseudocone translations along $\gamma$ of shape $h = (\quot{h}{X}, \quot{h}{f})$ to $E$ and that $E$ has an optranslation along $\gamma$ to $F$ of shape $k = (\quot{k}{X}, \quot{k}{f})$. Then if for all objects $X \in \Cscr_0$ we have adjoints $\quot{h}{X} \dashv \quot{k}{X}$ there is an adjunction:
	\[
	\begin{tikzcd}
		\PC(E) \ar[rr, bend right = 20, swap, ""{name = D}]{}{\ul{k}} & & \PC(F) \ar[ll, bend right = 20, swap, ""{name = U}]{}{\ul{h}} \ar[from = U, to = D, symbol = \dashv]{}{}
	\end{tikzcd}
	\]
\end{Theorem}
\begin{proof}
	The assumptions in the statement of the theorem provide us with pseudonatural transformations
	\[
	\begin{tikzcd}
		\Cscr^{\op} \ar[rr, bend left = 30, ""{name = U}]{}{F} \ar[rr, bend right = 30, swap, ""{name = D}]{}{E \circ \gamma^{\op}} & & \fCat \ar[from = U, to = D, Rightarrow, shorten <= 4pt, shorten >= 4pt]{}{h}
	\end{tikzcd}
	\qquad
	\begin{tikzcd}
		\Cscr^{\op} \ar[rr, bend right = 30, swap, ""{name = U}]{}{F} \ar[rr, bend left = 30,""{name = D}]{}{E \circ \gamma^{\op}} & & \fCat \ar[from = D, to = U, Rightarrow, shorten <= 4pt, shorten >= 4pt]{}{k}
	\end{tikzcd}
	\]
	together with adjoints $\quot{h}{X} \dashv \quot{k}{X}$ for all $X \in \Cscr_0$. From here the result follows by applying Theorem \ref{Thm: Section 3: Gamma-wise adjoints lift to equivariant adjoints}.
\end{proof}

We now show some basic results regarding the use of these theorems and illustrate why we explicitly presented them. In particular, we'll be able to deduce the existence of equivariant pullback and pushforward adjunctions and tensor-hom adjunctions directly from these results. While the existence of the functors we describe is well-known (especially in the geometric and topological cases; cf.\@ \cite{BernLun} and \cite{LusztigCuspidal2}, for instance), we present the results more to illustrate how to apply this formalism to situations which need not be the equivariant derived category. For instance, these techniques may be adapted to give rise to an equivariant derived perverse push-pull adjunction
\[
\begin{tikzcd}
D_G^b(\Per(X)) \ar[rr, bend right = 20, swap, ""{name = D}]{}{{}^{p}Rh_{\ast}}  & & D_G^b(\Per(Y)) \ar[ll, bend right = 20, swap, ""{name = U}]{}{{}^{p}Lh^{\ast}} \ar[from = U, to = D, symbol = \dashv]
\end{tikzcd}
\]
although we do not spell this out explicitly here.

\begin{corollary}\label{Cor: Pseudocone Functors: Pushforward functors for equivariant maps for scheme sheaves}
	Let $G$ be a smooth algebraic group and let $h:X \to Y$ be a $G$-equivariant morphism of $G$-varieties. Then there is an equivariant pushforward functor $Rh_{\ast}:D_G^b(X; \overline{\Q}_{\ell}) \to D_G^b(Y; \overline{\Q}_{\ell})$ which satisfies the adjunction:
	\[
	\begin{tikzcd}
		D_G^b(X;\overline{\Q}_{\ell}) \ar[rr, bend right = 20, swap, ""{name = D}]{}{R\ul{h}_{\ast}} & & D_G^b(Y;\overline{\Q}_{\ell}) \ar[ll, bend right = 20, swap, ""{name = U}]{}{\ul{h}^{\ast}} \ar[from = U, to = D, symbol = \dashv]{}{}
	\end{tikzcd}
	\]
\end{corollary}
\begin{proof}
	Consider the $2$-cell
	\[
	\begin{tikzcd}
		\SfResl_G(X) \ar[dr, swap]{}{\cong} \ar[rr, ""{name = U}]{}{\quo_{(-)\times X}} & & \Var_{/K} \\
		& \SfResl_G(Y) \ar[ur, swap]{}{\quo_{(-)\times Y}} \ar[from = U, to = 2-2, Rightarrow, shorten <= 4pt, shorten >= 4pt]{}{\varphi}
	\end{tikzcd}
	\]
	described in Corollary \ref{Cor: Pseudocone Functors: Existence of equivariant pullback  for schemes}. Now define the psuedofunctor $F:\Var_{/K}^{\op} \to \fCat$ given by
	\[
	F(Z) := D^b_c(Z;\overline{\Q}_{\ell})
	\]
	for objects and on morphisms $f:Z \to W$ by
	\[
	F(f) := f^{\ast}:D^b_c(W;\overline{\Q}_{\ell}) \to D^b_c(Z;\overline{\Q}_{\ell}).
	\]
	Then if $f:\Gamma \to \Gamma^{\prime}$ is a morphism in $\Sf(G)$, the morphisms $\overline{f \times \id_X}$ and $\overline{f \times \id_Y}$ are smooth by Proposition \ref{Prop: Section 2.1: The quotient prop}. Thus, since the diagram
	\[
	\begin{tikzcd}
		\Gamma \times X \ar[rr]{}{f \times \id_X} \ar[d, swap]{}{\id_{\Gamma} \times h} & & \Gamma^{\prime} \times X \ar[d]{}{\id_{\Gamma^{\prime}} \times h} \\
		\Gamma \times Y \ar[rr, swap]{}{f \times \id_Y} & & \Gamma^{\prime} \times Y
	\end{tikzcd}
	\]
	is a pullback diagram which is preserved by taking $G$-quotients, we find that
	\[
	\begin{tikzcd}
		G\backslash(\Gamma \times X) \ar[rr]{}{\overline{f \times \id_X}} \ar[d, swap]{}{\overline{\id_{\Gamma} \times h}} & & G \backslash (\Gamma^{\prime} \times X) \ar[d]{}{\overline{\id_{\Gamma^{\prime}} \times h}} \\
		G\backslash(\Gamma \times Y) \ar[rr, swap]{}{\overline{f \times \id_Y}} & & 	G\backslash(\Gamma^{\prime} \times Y)
	\end{tikzcd}
	\]
	is a pullback. Consequently the diagram
	\[
	\begin{tikzcd}
		D^b_c\left(G \backslash (\Gamma^{\prime} \times X); \overline{\Q}_{\ell}\right) \ar[rrr, ""{name = U}]{}{R(\overline{\id_{\Gamma^{\prime}} \times h})_{\ast}} \ar[d, swap]{}{(\overline{f \times \id_{X}})^{\ast}} & & & D^b_c\left(G \backslash (\Gamma^{\prime} \times Y); \overline{\Q}_{\ell}\right) \ar[d]{}{(\overline{f \times \id_Y})^{\ast}} \\
		D^b_c\left(G \backslash(\Gamma \times X);\overline{\Q}_{\ell}\right) \ar[rrr, swap, ""{name = D}]{}{R(\overline{\id_{\Gamma} \times h})_{\ast}} & & & D^b_c\left(G \backslash(\Gamma \times Y)\right) \ar[from = U, to= D, Rightarrow, shorten <= 4pt, shorten >= 4pt]{}{\quot{h}{f}} \ar[from = U, to= D, Rightarrow, shorten <= 4pt, shorten >= 4pt, swap]{}{\cong}
	\end{tikzcd}
	\]
	commutes up to natural isomorphism by the Smooth Base Change Theorem (which may be applied because all the maps in sight are maps between varieties and hence quasi-compact, finitely presented, and separated, as well as because the maps $\overline{f \times \id_X}$ and $\overline{f \times \id_Y}$ are smooth), and this in turn varies pseudonaturally over $\SfResl_G(X)^{\op}$. This gives us the pseudonatural transformation
	\[
	\begin{tikzcd}
		\SfResl_{G}(X)^{\op} \ar[rrr, bend left = 40, ""{name = U}]{}{F \circ \quo_{(-) \times X}^{\op}} \ar[rrr, bend right = 40, swap, ""{name = D}]{}{F \circ \quo_{(-) \times Y}^{\op} \circ \gamma^{\op}} & & & \fCat \ar[from = U, to = D, Rightarrow, shorten <= 4pt, shorten >= 4pt]{}{h}
	\end{tikzcd}
	\]
	where $\quot{h}{\Gamma \times X} := R(\overline{\id_{\Gamma} \times h})_{\ast}$ for objects $\Gamma \times X$ in $\SfResl_G(X)$ and where $\quot{h}{f}$ is given as above. This gives us our desired equivariant pushforward $R\ul{h}_{\ast}:D_G^b(X;\overline{\Q}_{\ell}) \to D_G^b(Y;\overline{\Q}_{\ell})$ by simply applying Proposition \ref{Prop: Pseudocone Functors: TRanslation}. Finally that $R\ul{h}_{\ast}$ is right adjoint to $\ul{h}^{\ast}$ follows from the fact that since $\ul{h}^{\ast}$ is the functor induced by applying Corollary \ref{Cor: Pseudocone Functors: Existence of equivariant pullback  for schemes} to the psuedofunctor $F$ and because $(\overline{\id_{\Gamma} \times h})^{\ast} \dashv R(\overline{\id_{\Gamma} \times h})_{\ast}$ for all $\Gamma$, our result follows from Theorem \ref{Thm: Pseudocone Functors: Adjoints between translations give adjoint functors}.
\end{proof}
\begin{corollary}\label{Cor: Pseudocone Functors: Exceptional Pushforward functors for equivariant maps for scheme sheaves}
	Let $G$ be a smooth algebraic group and let $h:X \to Y$ be a $G$-equivariant morphism of $G$-varieties. Then there are exceptional equivariant pullback functors $h^{!}:D_G^b(Y;\overline{\Q}_{\ell}) \to D_G^b(X;\overline{\Q}_{\ell})$ and exceptional equivariant pushforward functors $Rh_{!}:D_G^b(X; \overline{\Q}_{\ell}) \to D_G^b(Y; \overline{\Q}_{\ell})$ which satisfy the adjunction:
	\[
	\begin{tikzcd}
		D_G^b(Y;\overline{\Q}_{\ell}) \ar[rr, bend right = 20, swap, ""{name = D}]{}{\ul{h}^{!}} & & D_G^b(X;\overline{\Q}_{\ell}) \ar[ll, bend right = 20, swap, ""{name = U}]{}{R\ul{h}_{!}} \ar[from = U, to = D, symbol = \dashv]{}{}
	\end{tikzcd}
	\]
\end{corollary}
\begin{proof}
	From the proof of Corollary \ref{Cor: Pseudocone Functors: Pushforward functors for equivariant maps for scheme sheaves} we need only show that the diagrams
	\[
	\begin{tikzcd}
		D^b_c\left(G \backslash (\Gamma^{\prime} \times X); \overline{\Q}_{\ell}\right) \ar[rrr, ""{name = U}]{}{R(\overline{\id_{\Gamma^{\prime}} \times h})_{!}} \ar[d, swap]{}{(\overline{f \times \id_{X}})^{\ast}} & & & D^b_c\left(G \backslash (\Gamma^{\prime} \times Y); \overline{\Q}_{\ell}\right) \ar[d]{}{(\overline{f \times \id_Y})^{\ast}} \\
		D^b_c\left(G \backslash(\Gamma \times X);\overline{\Q}_{\ell}\right) \ar[rrr, swap, ""{name = D}]{}{R(\overline{\id_{\Gamma} \times h})_{!}} & & & D^b_c\left(G \backslash(\Gamma \times Y)\right) \ar[from = U, to= D, Rightarrow, shorten <= 4pt, shorten >= 4pt]{}{\quot{h_!}{f}} \ar[from = U, to= D, Rightarrow, shorten <= 4pt, shorten >= 4pt, swap]{}{\cong}
	\end{tikzcd}
	\]
	and
	\[
	\begin{tikzcd}
		D^b_c\left(G \backslash (\Gamma^{\prime} \times Y); \overline{\Q}_{\ell}\right) \ar[rrr, ""{name = U}]{}{R(\overline{\id_{\Gamma^{\prime}} \times h})^{!}} \ar[d, swap]{}{(\overline{f \times \id_{X}})^{\ast}} & & & D^b_c\left(G \backslash (\Gamma^{\prime} \times X); \overline{\Q}_{\ell}\right) \ar[d]{}{(\overline{f \times \id_Y})^{\ast}} \\
		D^b_c\left(G \backslash(\Gamma \times Y);\overline{\Q}_{\ell}\right) \ar[rrr, swap, ""{name = D}]{}{R(\overline{\id_{\Gamma} \times h})^{!}} & & & D^b_c\left(G \backslash(\Gamma \times X)\right) \ar[from = U, to= D, Rightarrow, shorten <= 4pt, shorten >= 4pt]{}{\quot{h^!}{f}} \ar[from = U, to= D, Rightarrow, shorten <= 4pt, shorten >= 4pt, swap]{}{\cong}
	\end{tikzcd}
	\]
	form invertible $2$-cells. However, this follows on one hand from the Smooth Base Change with Compact Support and on the other hand from the fact that smooth pullback commutes with exceptional pullback. Consequently the corollary follows from Theorem \ref{Thm: Pseudocone Functors: Adjoints between translations give adjoint functors} and the adjunctions $R(\overline{\id_{\Gamma} \times h})_{!} \cong (\overline{\id_{\Gamma} \times h})^{!}$ which hold for all $\Gamma \in \Sf(G)_0$.
\end{proof}

\begin{corollary}\label{Cor: Pseudocone Functors: Topological EDC has equivariant pushpull}
	Let $G$ be a topological group such that for every $n \in \N$ there is an $n$-acyclic free $G$-space $M$ which is a manifold. If $h:X \to Y$ is a $G$-equivariant morphism of $G$-spaces then there is an equivariant pullback/pushforward adjunction.
\end{corollary}
\begin{proof}
	By the assumptions on $G$ it follows that for any $n \in \N$ if $M$ is a specified $n$-acyclic free $G$-space which is a manifold then the second projection $\pi_2:M \times X \to X$ is a smooth\index[terminology]{Smooth! Map of $G$-Spaces} (in the sense of \cite[Section 1.7]{BernLun}; cf.\@ also \cite[Section 1.2]{Landesman} --- the connection between the two notions follows from the existence of smooth homeomorphisms $(0,1)^n \cong \R^n$) $n$-acyclic resolution of $X$. From here taking the subcategory $\mathbf{SRes}_G(X)$\index[notation]{SRes@$\mathbf{SRes}_G(X)$} of $G$-resolutions of $X$ generated by the projections $M \times X \to X$, their pullbacks against each other, and the smooth maps between them, we get that there is an equivalence 
	\[
	D_G^b(X) \simeq \PC(D^b_c(G \backslash (-)):\mathbf{SRes}_G(X)^{\op} \to \fCat)
	\]
	by \cite[Section 2.4.4]{BernLun}.
	
	We now argue the foundations required in order to define the pullback/pushforward $h^{\ast}, h_{\ast}$ both hold. Observe that analogous to the scheme-theoretic situation, there is an isomorphism of categories
	\[
	\gamma:\mathbf{SRes}_G(X) \xrightarrow{\cong} \mathbf{SRes}_G(Y)
	\]
	induced on objects by the assignment $M \times X \mapsto M \times Y$. Moreover the existence of the map $h:X \to Y$, analogously to the scheme-theoretic case, gives rise to a commuting cube
	\[
	\begin{tikzcd}
		& M \times Y \ar[dd, swap, near end]{}{\quo_{M \times Y}} \ar[rr]{}{f \times \id_Y} & & N \times Y \ar[dd]{}{\quo_{N \times Y}} \\
		M \times X \ar[dd, swap]{}{\quo_{M \times X}} \ar[rr, crossing over, near end]{}{f \times \id_X} \ar[ur]{}{\id_{M} \times h} & & N \times X \ar[ur,]{}{\id_{N}\times h} \\
		& G \backslash (M \times Y) \ar[rr, near end]{}{\overline{f \times \id_Y}} & & G \backslash (N \times Y) \\
		G \backslash (M \times X) \ar[ur]{}{\overline{\id_{M} \times h}} \ar[rr, swap]{}{\overline{f \times \id_X}} & & G \backslash (N \times X) \ar[ur, swap]{}{\overline{\id_{N} \times h}} \ar[from = 2-3, to = 4-3, crossing over, near end]{}{\quo_{N \times X}}
	\end{tikzcd}
	\]
	in $\Top$; note here that $M, N \in \mathbf{SRes}_G(\ast)$. This allows us to define the $2$-cell
	\[
	\begin{tikzcd}
		\mathbf{SRes}_G(X) \ar[rr, ""{name=U}]{}{\quo_{(-) \times X}} \ar[dr, swap]{}{\gamma} & & \Top \\
		& \mathbf{SRes}_G(Y) \ar[ur, swap]{}{\quo_{(-) \times Y}} \ar[from = U, to = 2-2, Rightarrow, shorten <= 4pt, shorten >= 4pt]{}{\alpha}
	\end{tikzcd}
	\]
	where for all $M \times X \to X$ in $\mathbf{SRes}_G(X)_0$, we have $\alpha_{M \times X} := \overline{\id_{M}\times h}.$
	
	Now consider the\footnote{By which we really mean ``any equivalent model'' and leave the compositors unspecified.} derived category pseudofunctor $D^b:\Top^{\op} \to \fCat$. We then get the pullback functor $h^{\ast}:D_G^b(Y) \to D_G^b(X)$ by using Theorem \ref{Thm: Pseudocone Functors: Pullback induced by fibre functors in pseudofucntor} applied to the current situation. Alternatively, since we have that the square
	\[
	\begin{tikzcd}
		G \backslash (M \times X) \ar[rr]{}{\overline{\id_M \times h}} \ar[d, swap]{}{\overline{f \times \id_X}} & & G \backslash (M \times Y) \ar[d]{}{\overline{f \times \id_Y}} \\
		G \backslash (N \times X) \ar[rr, swap]{}{\overline{\id_N \times h}} & & G \backslash (N \times Y) 
	\end{tikzcd}
	\]
	is a pullback in $\Top$ for every map $f$ in $\mathbf{SRes}_G(\ast)$ (analgously to how it was in $\Var_{/K}$ in the scheme-theoretic case), it follows from the topological Smooth Base Change Theorem (cf.\@ \cite[Section 1.8, Appendix A]{BernLun}) that there is a natural isomorphism
	\[
	\begin{tikzcd}
		D^b(G \backslash (N \times X)) \ar[rr, ""{name = U}]{}{R(\overline{\id_N \times h})_{\ast}} \ar[d, swap]{}{(\overline{f \times \id_X})^{\ast}} & & D^b(G \backslash (N \times Y)) \ar[d]{}{(\overline{f \times \id_Y})^{\ast}} \\
		D^b(G \backslash (M \times X)) \ar[rr, swap, ""{name = D}]{}{R(\overline{\id_M \times h})_{\ast}} & & D^b(G \backslash (M \times Y)) \ar[from = U, to = D, Rightarrow, shorten <= 4pt, shorten >= 4pt]{}{\quot{h}{f}}
	\end{tikzcd}
	\]
	which varies pseudonaturally in $\mathbf{SRes}_G(X)^{\op}$. Assembling these natural isomorphisms together with the fibre functors $h_{M \times X} := R(\overline{\id_M \times h})_{\ast}$ gives rise to a psuedonatural transformation
	\[
	\begin{tikzcd}
		\mathbf{SRes}_G(X)^{\op} \ar[rr, bend left = 30, ""{name = U}]{}{D^b \circ \quo_{(-) \times X}^{\op}} \ar[rr, bend right = 30, swap, ""{name = D}]{}{D^b \circ \quo_{{-} \times Y}^{\op} \circ \gamma^{\op}} & & \fCat \ar[from = U, to = D, Rightarrow, shorten <= 4pt, shorten >= 4pt]{}{Rh_{\ast}}
	\end{tikzcd}
	\]
	and hence gives us our equivariant pushforward by Theorem \ref{Thm: Functor Section: Psuedonatural trans are pseudocone functors}. Finally using the adjoints $(\overline{\id_M \times h})^{\ast} \dashv R(\overline{\id_M \times h})_{\ast}$ for all objects $M \times X$ gives us the adjoint
	\[
	\begin{tikzcd}
		D^b_G(X) \ar[rr, swap, bend right = 20, ""{name = U}]{}{Rh_{\ast}} & & D^b_G(Y) \ar[ll, swap, bend right = 20, ""{name = D}]{}{h^{\ast}} \ar[from = D, to = U, symbol = \dashv]
	\end{tikzcd}
	\]
	by Theorem \ref{Thm: Pseudocone Functors: Adjoints between translations give adjoint functors}.
\end{proof}
\begin{corollary}\label{Cor: Pseudocone Functors: Topological EDC has equivariant exceptional pushpull}
	Let $G$ be a topological group such that for every $n \in \N$ there is an $n$-acyclic free $G$-space $M$ which is a manifold. If $h:X \to Y$ is a $G$-equivariant morphism of $G$-spaces then there is an equivariant exceptional pullback/pushforward adjunction.
\end{corollary}
\begin{proof}
	Follow the proof of Corollary \ref{Cor: Pseudocone Functors: Topological EDC has equivariant pushpull} mutatis mutandis to get the diagram
	\[
	\begin{tikzcd}
		\mathbf{SRes}_G(X) \ar[rr, ""{name=U}]{}{\quo_{(-) \times X}} \ar[dr, swap]{}{\gamma} & & \Top \\
		& \mathbf{SRes}_G(Y) \ar[ur, swap]{}{\quo_{(-) \times Y}} \ar[from = U, to = 2-2, Rightarrow, shorten <= 4pt, shorten >= 4pt]{}{\alpha}
	\end{tikzcd}
	\]
	of functors where $\alpha_{M \times X} = \overline{\id_{M} \times h}$ and consider the derived category pseudofunctor $D^b:\Top^{\op} \to \fCat$. We now claim that there are natural isomorphisms
	\[
	\begin{tikzcd}
		D^b(G \backslash (N \times X)) \ar[rr, ""{name = U}]{}{R(\overline{\id_N \times h})_{!}} \ar[d, swap]{}{(\overline{f \times \id_X})^{\ast}} & & D^b(G \backslash (N \times Y)) \ar[d]{}{(\overline{f \times \id_Y})^{\ast}} \\
		D^b(G \backslash (M \times X)) \ar[rr, swap, ""{name = D}]{}{R(\overline{\id_M \times h})_{!}} & & D^b(G \backslash (M \times Y)) \ar[from = U, to = D, Rightarrow, shorten <= 4pt, shorten >= 4pt]{}{\quot{Rh_!}{f}}
	\end{tikzcd}
	\]
	and
	\[
	\begin{tikzcd}
		D^b(G \backslash (N \times Y)) \ar[rr, ""{name = U}]{}{(\overline{\id_N \times h})^{!}} \ar[d, swap]{}{(\overline{f \times \id_X})^{\ast}} & & D^b(G \backslash (N \times X)) \ar[d]{}{(\overline{f \times \id_Y})^{\ast}} \\
		D^b(G \backslash (M \times Y)) \ar[rr, swap, ""{name = D}]{}{(\overline{\id_M \times h})^{!}} & & D^b(G \backslash (M \times X)) \ar[from = U, to = D, Rightarrow, shorten <= 4pt, shorten >= 4pt]{}{\quot{h^!}{f}}
	\end{tikzcd}
	\]
	which vary pseudonaturally in $\mathbf{SRes}_G(X)^{\op}$. To see this, we note  that the existence of $\quot{Rh_!}{f}$ foll follows from \cite[Section 1.8]{BernLun}, as does the existence of $\quot{h^{!}}{f}$. Finally applying \ref{Prop: Pseudocone Functors: TRanslation} twice gives the existence of the functors $Rh_!:D^b_G(X) \to D_G^b(Y)$ and $h^{!}:D^b_G(Y) \to D^b_G(X)$; using the adjoints $R(\overline{\id_M \times h})_! \dashv (\overline{\id_M \times h})^{!}$ together with Theorem \ref{Thm: Pseudocone Functors: Adjoints between translations give adjoint functors} gives the adjunction:
	\[
	\begin{tikzcd}
		D_G^b(Y) \ar[rr, bend right = 20, swap, ""{name = D}]{}{h^{!}} & & D_G^b(X) \ar[ll, bend right = 20, swap, ""{name = U}]{}{Rh_!} \ar[from = U, to = D, symbol = \dashv]
	\end{tikzcd}
	\]
\end{proof}
\newpage

\chapter{Pseudocone Homological Algebra}\label{Section: Pseudocone Algebra}

We now move to the homological algebraic chapter of the monograph: that of triangulated pseudocone categories and the $t$-structures which arise upon categories. Much of the structure that we describe and study in this section can be seen as an application to some of the results proved about pseudocone categories and their functors which involve additive and Abelian categorical structures. However, it will be important for applications in practice (such as when studying how the equivariant derived category relates to equivariant perverse sheaves, equivariant sheaves, equivariant cohomology, and the relations between each corresponding category and the functors which we can define; cf.\@ Theorems \ref{Thm: Section SFF: The big results on the EDC for schemes} and \ref{Thm: Section SFF: The big results on the EDC for spaces}) to have developed the language necessary that will allow us to discuss such homological algebraic situations at a reasonable level of generality. Of particular importance is the theorem below (cf.\@ Theorem \ref{Theorem: Heart of equivariant t-structure is equivariant heart of t-structure}) which says that if $F$ is a triangulated pseudofunctor on $X$ such that each functor $F(f)$ restricts to an appropriate functor $F(f)^{\heartsuit}$ between the hearts of the fibre categories $F(X)^{\heartsuit}$, then there is a canonical $t$-structure on $\PC(F)$ such that there is an \emph{equality}
\[
\PC(F)^{\heartsuit} = \PC(F^{\heartsuit}).
\]
I like to think of this as the ``Change of Heart'' theorem.
In the journey of describing the triangulations on equivariant categories, we will begin by showing how to take an additive pseudofunctor with suspensions on each fibre category and turn it into a suspended additive pseudocone category.

\section{Triangulated Pseudocone Categories}\label{Section: Triangles are go}
Our study of psuedocone homological algebra begins by describing when, given a pseudofunctor $F:\Cscr^{\op} \to \fCat$, the category $\PC(F)$ is triangulated with its triangulation induced by having each fibre category $F(X)$ be a triangulated category. Because of this we first need to recall what it means to be a category with suspension, deduce when $\PC(F)$ is a category with suspension, and then go on to define triangulated pseudofunctors.

\begin{definition}
	A category with suspension\index[terminology]{Category with Suspension} is a pair $(\Cscr, [1])$ where $\Cscr$ is an additive category and $[1]:\Cscr \to \Cscr$ is an autoequivalence. In such a case the functor $[1]$ is called the suspension functor.
\end{definition}
\begin{remark}
	We have gone with the algebraic-geometric naming convention for the suspension functor in our definition above; algebraic-topological naming coneventions tend to use $\Sigma$ for what we have called $[1]$. Occasionally we will use the topological notation, but generically we will use the geometric notation.
\end{remark}
\begin{proposition}\label{Prop: Section Triangle: Equivariant suspensions}
	Let $F:\Cscr^{\op} \to \fCat$ be an additive pseudofunctor such that each category $F(X)$ is equipped with an autoequivalence $\quot{\Sigma}{X}:F(X) \to F(X)$ and assume that these equivalences vary pseudonaturally in $\Cscr^{\op}$. Then $(\PC(F),\ul{\Sigma})$ is a category with suspension where $\Sigma$ is the autoequivalence $\Sigma:\PC(F) \to  \PC(F)$ induced by the $\quot{\Sigma}{X}$.
\end{proposition}
\begin{proof}
	By Corollary \ref{Cor: Section 2: Additive pseudofunctor gives additive equivariat cat} we have that $\PC(F)$ is additive and the assumptions in the statement of the proposition imply that the $\quot{\Sigma}{X}$ are the object functors in a pseudonatural equivalence:
	\[
	\begin{tikzcd}
		\Cscr^{\op} \ar[rr, bend left = 20, ""{name = U}]{}{F} \ar[rr, bend right = 20, swap, ""{name = D}]{}{F} & & \fCat \ar[from = U, to= D, Rightarrow, shorten <= 4pt, shorten >= 4pt]{}{\Sigma}
	\end{tikzcd}
	\]
	Let $\Omega:F \to F$ be an inverse pseudonatural transformation to $\Sigma$ and note that for all $X \in \Cscr_0$ we have an adjunction $\quot{\Omega}{X} \dashv \quot{\Sigma}{X}$. Applying Theorem \ref{Thm: Section 3: Gamma-wise adjoints lift to equivariant adjoints} gives the adjoint $\ul{\Omega} \dashv \ul{\Sigma}$ which is an equivalence by Corollary \ref{Cor: Section 3: Fibre-wise equivalences are equivalences of equivariant categories}. Thus $(\PC(F), \ul{\Sigma})$ is a category with suspension.
\end{proof}
\begin{corollary}
	For any smooth algebraic group $G$ and left $G$-variety $X$ the shift functors $[1]:D_G^b(X) \to D_G^b(X)$ and $[1]:\DbeqQl{X} \to \DbeqQl{X}$ given on objects by
	\[
	A[1] := \left\lbrace \quot{A}{\Gamma}[1]_{\Gamma} \; | \; \Gamma \times X\in \SfResl_G(X)_0, \quot{A}{\Gamma} \in A \right\rbrace,
	\] 
	where each $[1]_{\Gamma}$ is the corresponding shift functor, determine categories with suspension. 
\end{corollary}
\begin{proof}
	Checking that the functors $[1]$ are well-defined and induced by a pseudonatural transformation is immediate by noting that for any morphism $f \times \id_X:\Gamma \times X \to \Gamma^{\prime} \times X$, the corresponding map 
	\[
	\overline{f \times \id_X}:G \backslash (\Gamma \times X) \to G \backslash (\Gamma^{\prime} \times X)
	\]
	is smooth so the pullback $(\overline{f \times \id_X})^{\ast}$ commutes with shift functors up to isomorphism. The corollary now follows from an immediate application of Proposition \ref{Prop: Section Triangle: Equivariant suspensions}.
\end{proof}
\begin{remark}
	In the topological case, expressing the existence of $[1]:D_G^b(X) \to D_G^b(X)$ is surprisingly subtle in the non-pseudocone model due to the definition of objects as given in \cite{BernLun}; you need a notion of locally acyclic base change. This can be done directly when $X$ is some sort of manifold or pseudomanifold and $G$ is some sort of Lie group, but does hold more generally.
\end{remark}

The definitions we present below are the final remaining ingredients we need to be able to prove that if the pseudofunctor $F$ is triangulated then so is $\PC(F)$. However, when proving that $\PC(F)$ is triangulated we will follow the refinement of the triangulation axioms as developed by May and presented in \cite{MayTrace} as opposed to the original axiomatization of Verdier in \cite{VerdierDerivedCat}. This refinement has the benefit of having fewer technical conditions to check: May's refinement only asks us to check Verdier's axiom TR1, a weaker form of Verdier's axiom TR2 (we only need to be able to rotate triangles in one direction), and Verdier's axiom TR4 (the Octahedral Axiom); TR3 (the diagram filling axiom) is entirely redundant. We also follow the naming and numbering of the triangulation axioms of \cite{MayTrace} instead of \cite{VerdierDerivedCat}.

\begin{definition}
	Let $(\Ascr, [1])$ be a category with suspension with autoequivalence $[1]:\Ascr \to \Ascr$. A triangle\index[terminology]{Triangle} in $\Ascr$ is a triple of morphisms in $\Ascr$ of the form:
	\[
	\xymatrix{
		A \ar[r]^-{f} & B \ar[r]^-{g} & C \ar[r]^-{h} & A[1]
	}
	\]
	Two triangles are said to be isomorphic if there is a commuting diagram of the form
	\[
	\xymatrix{
		A \ar[r]^-{f} \ar[d]_{\rho} & B \ar[r]^-{g} \ar[d]_{\varphi} & C \ar[r]^-{h} \ar[d]^{\psi} & A[1] \ar[d]^{\rho[1]} \\
		X \ar[r]_-{f^{\prime}} & Y \ar[r]_-{g^{\prime}} & Z \ar[r]_-{h^{\prime}} & X[1]
	}
	\]
	where $\rho, \varphi, \psi$ are all isomorphisms. We will describe triangles either visually as above or as a sextuple $(A,B,C,f,g,h)$.
\end{definition}
\begin{definition}[{\cite{MayTrace}}]\label{Defn: Section Triangle: Traingulated Categories}
	A triangulated category\index[terminology]{Triangulated Category} is an additive category $\Ascr$ with an autoequivalence $[1]:\Ascr \to \Ascr$ equipped with a collection of triangles called distinguished triangles which satisfy the following axioms:
	\begin{itemize}
		\item[T1.] Let $A \in \Ascr_0$ and let $f:A \to B$ be in $\Ascr_1$. Then:
		\begin{itemize}
			\item The triangle $A \xrightarrow{\id} A \xrightarrow{0} 0 \to A[1]$ is distinguished;
			\item There exists a distinguished triangle which starts at $f$, i.e., there exist maps $g:B \to C, h:C \to A[1]$ such that the triangle 
			\[
			\xymatrix{
				A \ar[r]^-{f} & B \ar[r]^-{g} & C \ar[r]^-{h} & A[1]
			}
			\]
			is distinguished in $\Ascr$;
			\item Any triangle isomorphic to a distinguished triangle is itself distinguished.
		\end{itemize}
		\item[T2.] If the triangle
		\[
		\xymatrix{
			A \ar[r]^-{f} & B \ar[r]^-{g} & C \ar[r]^-{h} & A[1]
		}
		\]
		is distinguished, then so is:
		\[
		\xymatrix{
			B \ar[r]^-{g} & C \ar[r]^-{h} & A[1] \ar[r]^-{-f[1]} & B[1]
		}
		\]
		\item[T3.] Consider a diagram of the form
		
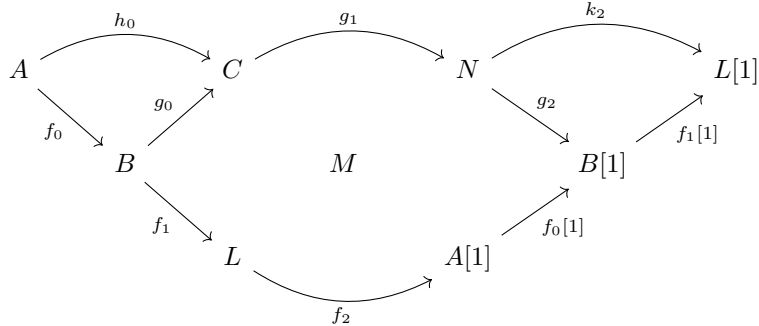
\[
		\begin{tikzcd}
			A \ar[dr, swap]{}{f_0} \ar[rr, bend left = 30]{}{h_0} & & C \ar[rr, bend left = 30]{}{g_1} & & N \ar[rr, bend left = 30]{}{k_2} \ar[dr]{}{g_2} & & L[1] \\
			& B \ar[ur]{}{g_0} \ar[dr, swap]{}{f_1} & & M & & B[1] \ar[ur, swap]{}{f_1[1]} \\
			& & L \ar[rr, bend right = 30, swap]{}{f_2} & & A[1] \ar[ur, swap]{}{f_0[1]}
		\end{tikzcd}
		\]
		where the triangles $(f_0,f_1,f_2)$ and $(g_0, g_1, g_2)$ are distinguished, $h_0 = g_0 \circ f_0$, and $k_2 =f_1[1] \circ  g_2$. If there are morphisms $h_1:C \to M$ and $h_2:M \to A[1]$ which make the triangle $(h_0, h_1, h_2)$ distinguished, then there exist maps $k_0:L \to M$ and $k_1:M \to N$ which make the triangle $(k_0, k_1, k_2)$ distinguished and the diagram
		\[
		\begin{tikzcd}
			A \ar[dr, swap]{}{f_0} \ar[rr, bend left = 30]{}{h_0} & & C \ar[dr]{}{h_1} \ar[rr, bend left = 30]{}{g_1} & & N \ar[rr, bend left = 30]{}{k_2} \ar[dr]{}{g_2} & & L[1] \\
			& B \ar[ur]{}{g_0} \ar[dr, swap]{}{f_1} & & M  \ar[ur]{}{k_1} \ar[dr]{}{h_2} & & B[1] \ar[ur, swap]{}{f_1[1]} \\
			& & L \ar[ur]{}{k_0} \ar[rr, bend right = 30, swap]{}{f_2} & & A[1] \ar[ur, swap]{}{f_0[1]}
		\end{tikzcd}
		\]
		commute in $\Ascr$. We call such a diagram a braid of triangles.\index[terminology]{Triangles! Braid of}
	\end{itemize}
\end{definition}
\begin{definition}\label{Defn: Section Triangle: Triangulated functor}
	Let $\Ascr, \Bscr$ be triangulated categories and let $F:\Ascr \to \Bscr$ be a functor. We say that $F$ is triangulated\index[terminology]{Triangulated Functor} if $F$ is additive and sends distinguished triangles in $\Ascr$ to distinguished triangles in $\Bscr$.
\end{definition}
\begin{definition}\label{Defn: Section Triangle: Triangulated Pseudofunctor}
	A pseudofunctor $F:\Cscr^{\op} \to \fCat$ is triangulated\index[terminology]{Triangulated Pseudofunctor} if for every object $X \in \Cscr_0$ the category $F(X)$ is triangulated and if for every morphism $f \in \Cscr_1$ the functor $F(f)$ is triangulated.
\end{definition}

\begin{Theorem}\label{Theorem: Section Triangle: Equivariant triangulation}
	Let $F:\Cscr^{\op} \to \fCat$ be a triangulated pseudofunctor. Then the category $\PC(F)$ is a triangulated category with triangulation induced by saying that a triangle
	\[
	\xymatrix{
		A \ar[r]^-{P} & B \ar[r]^-{\Sigma} & C \ar[r]^-{\Phi} & A[1]
	}
	\]
	is distinguished if and only if for all $X \in \Cscr_0$ the triangle
	\[
	\xymatrix{
		\quot{A}{X} \ar[r]^-{\quot{\rho}{X}} & \quot{B}{X} \ar[r]^-{\quot{\sigma}{X}} & \quot{C}{X} \ar[r]^-{\quot{\varphi}{X}} & \quot{A}{X}[1]
	}
	\]
	is distinguished in $F(\XGamma)$.
\end{Theorem}
\begin{proof}
	We begin by verifying Axiom T1 of \cite{MayTrace}. With this in mind let $A \in \PC(F)_0$. Then because for each $X \in \Cscr_0$ the triangle
	\[
	\xymatrix{
		\quot{A}{X} \ar@{=}[r] & \quot{A}{X} \ar[r] & \quot{0}{X} \ar[r] & \quot{A}{X}[1]
	}
	\]
	is distinguished in $F(\XGamma)$, it follows that the triangle
	\[
	\xymatrix{
		A \ar@{=}[r] & A \ar[r] & 0 \ar[r] & A[1]
	}
	\]
	is as well.
	
	We now verify the next aspect of Axiom T1. Assume that the triangle
	\[
	\xymatrix{
		A \ar[r]^-{P} & B \ar[r]^-{\Sigma} & C \ar[r]^-{\Phi} & A[1]
	}
	\]
	is distinguished in $\PC(F)$ and that the triangle
	\[
	\xymatrix{
		A^{\prime} \ar[r]^-{P^{\prime}} & B^{\prime} \ar[r]^-{\Sigma^{\prime}} & C^{\prime} \ar[r]^-{\Phi^{\prime}} & A^{\prime}[1]
	}
	\]
	is isomorphic to $(A,B,C,P,\Sigma,\Phi)$. Then for all $X \in \Cscr_0$ the triangle 
	\[
	(\quot{A}{X}, \quot{B}{X}, \quot{C}{X}, \quot{\rho}{X}, \quot{\sigma}{X}, \quot{\varphi}{X})
	\]
	is distinguished and $(\quot{A}{X}^{\prime},\quot{B}{X}^{\prime},\quot{C^{\prime}}{X},\quot{\rho^{\prime}}{X},\quot{\sigma^{\prime}}{X},\quot{\varphi^{\prime}}{X})$ is isomorphic to this triangle as well. Thus the triangles $(\quot{A}{X}^{\prime},\quot{B}{X}^{\prime},\quot{C^{\prime}}{X},\quot{\rho^{\prime}}{X},\quot{\sigma^{\prime}}{X},\quot{\varphi^{\prime}}{X})$ are distinguished for all $X \in \Cscr_0$ and so we conclude that the triangle
	\[
	\xymatrix{
		A^{\prime} \ar[r]^-{P^{\prime}} & B^{\prime} \ar[r]^-{\Sigma^{\prime}} & C^{\prime} \ar[r]^-{\Phi^{\prime}} & A^{\prime}[1]
	}
	\]
	is distinguished in $\PC(F)$.
		
	We now verify the last item of Axiom T1. Let $P:A \to B$ be a morphism in $\PC(F)$ and note that we now must find a distinguished triangle in $\PC(F)$ which begins at $P$. For this note that for all $X \in \Cscr_0$ we can find distinguished triangles
	\[
	\xymatrix{
		\quot{A}{X} \ar[r]^-{\quot{\rho}{X}} & \quot{B}{X} \ar[r]^-{\quot{\sigma}{X}} & \quot{C}{X} \ar[r]^-{\quot{\varphi}{X}} & \quot{A}{X}[1]
	}
	\]
	in each category $F(\XGamma)$. Now using that each functor $F(\overline{f})$ is triangulated allows us to produce the commuting diagram
	\[
	\xymatrix{
		F(f)\quot{A}{Y} \ar[rr]^-{F(f)\quot{\rho}{Y}} \ar[d]_{\tau_f^A} & & F(f)\quot{B}{Y} \ar[rr]^-{F(f)\quot{\sigma}{Y}} \ar[d]^{\tau_f^B} & & F(f)\quot{C}{Y} \ar[rr]^-{F(f)\quot{\varphi}{Y}} & & F(f)\quot{A}{Y} \ar[d]^{\tau_f^{A[1]}} \\
		\quot{A}{X} \ar[rr]_-{\quot{\rho}{X}} & & \quot{B}{X} \ar[rr]_-{\quot{\sigma}{X}} & & \quot{C}{X} \ar[rr]_-{\quot{\varphi}{X}} & & \quot{A}{X}[1]
	}
	\]
	in $F(X)$ where $f \in \Cscr_1$ with $f:X \to Y$. Applying \cite[Lemmas 2.2, 2.3]{MayTrace} then allows us to deduce the existence of an isomorphism $\alpha:F(f)\quot{C}{Y} \to \quot{C}{X}$ which make the diagram
	\[
	\xymatrix{
		F(f)\quot{A}{Y} \ar[rr]^-{F(f)\quot{\rho}{Y}} \ar[d]_{\tau_f^A} & & F(f)\quot{B}{Y} \ar[rr]^-{F(f)\quot{\sigma}{Y}} \ar[d]^{\tau_f^B} & & F(f)\quot{C}{Y} \ar[d]^{\alpha_f}_{\cong} \ar[rr]^-{F(f)\quot{\varphi}{Y}} & & F(f)\quot{A}{Y} \ar[d]^{\tau_f^{A[1]}} \\
		\quot{A}{X} \ar[rr]_-{\quot{\rho}{X}} & & \quot{B}{X} \ar[rr]_-{\quot{\sigma}{X}} & & \quot{C}{X} \ar[rr]_-{\quot{\varphi}{X}} & & \quot{A}{X}[1]
	}
	\]
	commute. Note that if we can prove that the isomorphisms $\alpha_f$ satisfy the cocycle condition they automatically make 
	\[
	\Sigma := \lbrace \quot{\sigma}{X} \; | \; X \in \Cscr_0\rbrace
	\] 
	and
	\[
	\Phi := \lbrace \quot{\varphi}{X} \; | \; X \in \Cscr_0 \rbrace
	\]
	morphisms in $\PC(F)$. These morphisms can be shown to make the diagrams, for any composable $f,g \in \Cscr_1$,
	\[
	\begin{tikzcd}
		F(f)F\quot{C}{Y} \ar[rr, ""{name = U}]{}{\alpha_f} & & \quot{C}{X} \\
		F(f)F(\overline{g})\quot{C}{Z} \ar[u]{}{F(f)\alpha_{g}} \ar[rr, swap, ""{name = L}]{}{\phi_{f,g}^{\quot{C}{Z}}} & & F(\overline{g}\circ f)\quot{C}{Z} \ar[u, swap]{}{\alpha_{g \circ f}} \ar[from = U, to = L, Rightarrow, shorten >= 4pt, shorten <= 4pt]{}{\cong}
	\end{tikzcd}
	\]
	commute up to a unique invertible $2$-cell in $F(X)$ because each morphism $\alpha_f$ is induced from a pseudocolimit in $F(\XGamma)$; this pseudocolimit exists because $F(X)$ is a triangulated category and $\alpha_f$ is induced by the corresponding pseudouniversal property. As such we invoke the Axiom of Choice to select our isomorphisms $\alpha_f$ such that for all composable arrows $X \xrightarrow{f} Y \xrightarrow{g} Z$ in $\Cscr$ the identity $\alpha_{g \circ f} \circ \phi_{f,g} = \alpha_f \circ F(f)\alpha_g$ hold. With this choice of isomorphisms we find that the pair
	\[
	(\lbrace \quot{C}{X} \; | \; X \in \Cscr_0 \rbrace, \lbrace \alpha_f \; | \; f \in \Cscr_1\rbrace)
	\]
	determines an object with morphisms $\Sigma:B \to C$ and $\Phi:C \to A[1]$ which make the triangle
	\[
	\xymatrix{
		A \ar[r]^-{P} & B \ar[r]^-{\Sigma} & C \ar[r]^-{\Phi} & A[1]
	}
	\]
	distinguished in $\PC(F)$. This completes the verification of Axiom T1.
	
	We now proceed to verify Axiom T2. Assume that the triangle
	\[
	\xymatrix{
		A \ar[r]^-{P} & B \ar[r]^-{\Sigma} & C \ar[r]^-{\Phi} & A[1]
	}
	\]
	is distinguished in $\PC(F)$. Then for each $X \in \Cscr_0$ the triangle
	\[
	\xymatrix{
		\quot{A}{X} \ar[r]^-{\quot{\rho}{X}} & \quot{B}{X} \ar[r]^-{\quot{\sigma}{X}} & \quot{C}{X} \ar[r]^-{\quot{\varphi}{X}} & \quot{A[1]}{X}
	}
	\]
	is distinguished in $F(\XGamma)$, so applying Axiom T2 in $F(\XGamma)$ tells us that the triangle
	\[
	\xymatrix{
		\quot{B}{X} \ar[r]^-{\quot{\sigma}{X}} & \quot{C}{X} \ar[r]^-{\quot{\varphi}{X}} & \quot{A}{X}[1] \ar[r]^-{-\quot{\rho[1]}{X}} & \quot{B}{X}[1]
	}
	\]
	is distinguished as well. Because this happens for each $X \in \Cscr_0$ and these give the corresponding $X$-local descriptions of the objects in $\PC(F)$, we conclude that the triangle
	\[
	\xymatrix{
		B \ar[r]^-{\Sigma} & C \ar[r]^-{\Phi} & A[1] \ar[r]^-{-P[1]} & B[1]
	}
	\]
	is distinguished in $\PC(F)$. This establishes Axiom T2.
	
	Finally we verify that Axiom T3 (Verdier's Axiom) holds in $\PC(F)$. Assume that we have a diagram in $\PC(F)$ of the form
	\[
	\begin{tikzcd}
		A \ar[dr, swap]{}{P_0} \ar[rr, bend left = 30]{}{\Phi_0} & & C \ar[rr, bend left = 30]{}{\Sigma_1} & & N \ar[rr, bend left = 30]{}{\Psi_2} \ar[dr]{}{\Sigma_2} & & L[1] \\
		& B \ar[ur]{}{\Sigma_0} \ar[dr, swap]{}{P_1} & & M & & B[1] \ar[ur, swap]{}{P_1[1]} \\
		& & L \ar[rr, bend right = 30, swap]{}{P_2} & & A[1] \ar[ur, swap]{}{P_0[1]}
	\end{tikzcd}
	\]
	in $\PC(F)$ in which the triangles $(P_0,P_1,P_2)$ and $(\Sigma_0,\Sigma_1,\Sigma_2)$ are both distinguished. Now assume that there are morphisms $\Phi_1:C \to M$ and $\Phi_2:M \to A[1]$ such that the triangle $(\Phi_0, \Phi_1, \Phi_2)$ is distinguished. We now need to construct morphisms $\Psi_0:L \to M$ and $\Psi_1:M \to N$ which make the triangle $(\Psi_0, \Psi_1, \Psi_2)$ and render the diagram commutative. To go about this, let $X \in \Cscr_0$ and invoke Axiom T3 for this $X$. We then produce morphisms $\quot{\psi_0}{X}:\quot{L}{X} \to \quot{M}{X}$ and $\quot{\psi_1}{X}:\quot{M}{X} \to \quot{N}{X}$ which make the diagram
	
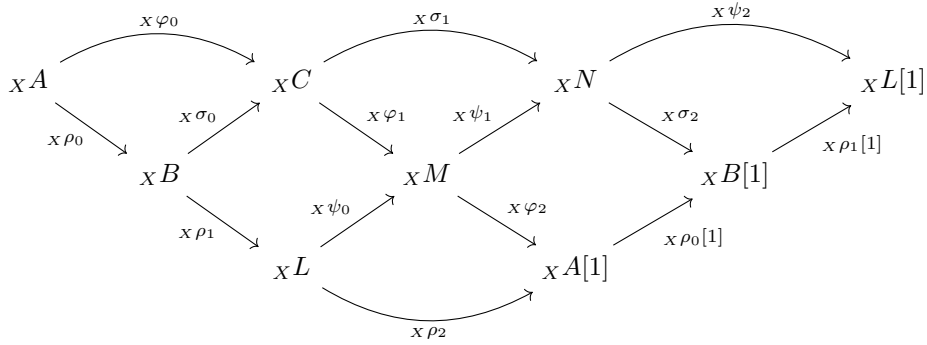
\[
	\begin{tikzcd}
		\quot{A}{X} \ar[dr, swap]{}{\quot{\rho_0}{X}} \ar[rr, bend left = 30]{}{\quot{\varphi_0}{X}} & & \quot{C}{X} \ar[dr]{}{\quot{\varphi_1}{X}} \ar[rr, bend left = 30]{}{\quot{\sigma_1}{X}} & & \quot{N}{X} \ar[rr, bend left = 30]{}{\quot{\psi_2}{X}} \ar[dr]{}{\quot{\sigma_2}{X}} & & \quot{L[1]}{X} \\
		& \quot{B}{X} \ar[ur]{}{\quot{\sigma_0}{X}} \ar[dr, swap]{}{\quot{\rho_1}{X}} & & \quot{M}{X} \ar[ur]{}{\quot{\psi_1}{X}} \ar[dr]{}{\quot{\varphi_2}{X}} & & \quot{B}{X}[1] \ar[ur, swap]{}{\quot{\rho_1[1]}{X}} \\
		& & \quot{L}{X} \ar[rr, bend right = 30, swap]{}{\quot{\rho_2}{X}} \ar[ur]{}{\quot{\psi_0}{X}} & & \quot{A}{X}[1] \ar[ur, swap]{}{\quot{\rho_0[1]}{X}}
	\end{tikzcd}
	\]
	commute with the triangle $(\quot{\psi_0}{X},\quot{\psi_1}{X}, \quot{\psi_2}{X})$ distinguished in addition to the others induced from the distinguished triangles given by assumption. Now using that for any $f \in \Cscr_1$ the functor $F(f)$ is triangulated, together with the fact that each of $P_1, P_2, P_0[1], \Sigma_1, \Sigma_2, \Phi_1,$ and $\Phi_2$ are $\PC(F)$-morphisms allows us to deduce via tedious diagram chases that the diagrams
	\[
	\begin{tikzcd}
		F(f)\quot{L}{Y} \ar[rr, ""{name = LU}]{}{F(f)\quot{\psi_0}{Y}} \ar[d, swap]{}{\tau_f^L} & & F(f)\quot{M}{Y} \ar[d]{}{\tau_f^M} & F(f)\quot{M}{Y} \ar[rr, ""{name = RU}]{}{F(f)\quot{\psi_1}{Y}} \ar[d, swap]{}{\tau_f^M} & & F(f)\quot{N}{Y} \ar[d]{}{\tau_f^N} \\
		\quot{L}{X} \ar[rr, swap, ""{name = LL}]{}{\quot{\psi_0}{X}} & & \quot{M}{X} & \quot{M}{X} \ar[rr, swap, ""{name = RL}]{}{\quot{\psi_1}{X}} & & \quot{N}{X}
		\ar[from = LU, to = LL, Rightarrow, shorten >= 4pt, shorten <= 4pt]{}{\cong} \ar[from = RU, to = RL, Rightarrow, shorten >= 4pt, shorten <= 4pt]{}{\cong}
	\end{tikzcd}
	\]
	both commute up to a natural isomorphism induced from the $2$-cells that compare any two morphisms which both provide braids of triangles. We once again invoke the Axiom of Choice to redefine the $\quot{\psi_i}{X}$ if any of these natural isomorphisms are not strict identities and replace them with their strict counterparts. This allows us to deduce that $\Psi_0 := \lbrace \quot{\psi_0}{X} \; | \; X \in \Cscr_0 \rbrace$ and $\Psi_1 := \lbrace \quot{\psi_1}{X} \; | \; X \in \Cscr_0 \rbrace$ are morphisms which makes the triangle $(\Psi_0,\Psi_1, \Psi_2)$ distinguished and renders the briad
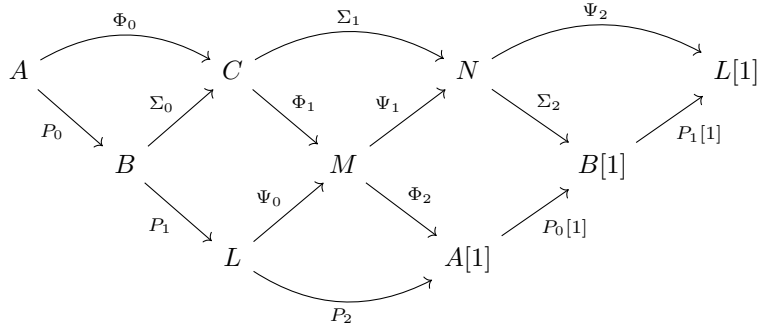
	\[
	\begin{tikzcd}
		A \ar[dr, swap]{}{P_0} \ar[rr, bend left = 30]{}{\Phi_0} & & C \ar[rr, bend left = 30]{}{\Sigma_1} \ar[dr]{}{\Phi_1} & & N \ar[rr, bend left = 30]{}{\Psi_2} \ar[dr]{}{\Sigma_2} & & L[1] \\
		& B \ar[ur]{}{\Sigma_0} \ar[dr, swap]{}{P_1} & & M \ar[ur]{}{\Psi_1} \ar[dr]{}{\Phi_2} & & B[1] \ar[ur, swap]{}{P_1[1]} \\
		& & L \ar[rr, bend right = 30, swap]{}{P_2} \ar[ur]{}{\Psi_0} & & A[1] \ar[ur, swap]{}{P_0[1]}
	\end{tikzcd}
	\]
	commutative. This completes the verification of the Axiom T3 and finishes the proof of the theorem.
\end{proof}
\begin{remark}
	In the above proof we have had to use the Axiom of Choice to artificially select morphisms to satisfy strictness conditions of our transition isomorphisms. This is largely because we are considering $1$-categorical information (strictly commuting diagrams) of what are inherently $2$-categorical (or $\infty$-categorical) constructions: the diagrams we construct only generically commute only up to coherent isomorphism provided by the pseudouniversal properties of pseudolimits and pseudocolimits. These issues disappear, however, in the $\infty$-categorical formalism, as the lack of strict equality becomes a feature rather than a bug; we record all such coherences and ways of moving between such isomorphisms as various higher morphisms and higher cells. In fact, by working with a notion of pseudocone $\infty$-category (the category defined by, for an $\infty$-pseudofunctor ``$F:\mathsf{Ner}(\Cscr^{\op}) \to \mathfrak{QCat}$''\index[notation]{QCat@$\mathfrak{QCat}$} for $\mathfrak{QCat}$ the $(\infty,2)$-category of quasicategories) into the $(\infty,2)$-category of quasi-categories (or ``$F:\SfResl_G(X)^{\op} \to \mathscr{K}$'' where $\mathscr{K}$ is an $\infty$-cosmos with a notion of stable $\infty$-category) whose fibre $\infty$-categories are all stable with triangulated $\infty$-fibre functors, we can sidestep these issues entirely. In such a case we should have that our $(\infty,1)$-category of pseudocones arises as the $(\infty,1)$-category
	\[
	\PC(F) = (\infty,2)\mathsf{Cat}(\mathsf{Ner}(\Cscr^{\op}),\mathfrak{QCat})(\cnst\One, F)
	\]
	and most, if not all, of the properties here will follow from taking the corresponding homotopy categories and homotopy $2$-categories. In particular, we expect a result analogous to Proposition \ref{Prop: Pseudocone Section: CSections are pseudolimits}: It should be the case that 
	\[
	\PC(F) \simeq \mathbf{CSect}_{B}(p:\El(F) \to \mathsf{Ner}(\Cscr^{\op}))
	\]
	for the corresponding $\infty$-categorical elements Cartesian fibration (cf.\@ \cite{GepnerEtAl}). We will not explore these notions here and leave them entirely speculative and as future work for the author; however, such a formalism will naturally lead to a development of an equivariant six functor formalism and pseudocone six functor formalism for pseudocone $\infty$-categories in the style of Mann as well as Scholze (cf.\@ \cite[Appendix A.5]{mann2022padic}, \cite{ScholzeSixFunctor}) and following the philosophy of \cite{gallauer2021introduction}.
\end{remark}
An immediate corollary of Theorem \ref{Theorem: Section Triangle: Equivariant triangulation} above is that it gives triangulated structures on the equivariant derived categories in the case where $G$ is a smooth algebraic group and when $X$ is a left $G$-variety.
\begin{corollary}\label{Cor: Section Triangle: EDC is triangle}
	Let $G$ be a smooth algebraic group and let $X$ be a left $G$-variety. Then	the categories $\DbeqQl{X}$ and $D_G^b(X)$ are triangulated.
\end{corollary}
When $G$ is a topological group and $X$ is a left $G$-space it is also true that $D_G^b(X)$ is a triangulated category, but we must instead use the proof of \cite[Section 2.5]{BernLun}. However, if we restrict to the following situation we can use Theorem \ref{Theorem: Section Triangle: Equivariant triangulation}.
\begin{corollary}\label{Cor: Section Triangle: EDC is triangle topology}
	Let $G$ be a topological group such that for every $n \in \N$ there is an $n$-acyclic free $G$-space $M$ which is a manifold. Then for any left $G$-space $X$ the category $D_G^b(X)$ is triangulated.
\end{corollary}
\begin{proof}
	Compare this to the proof of Corollary \ref{Cor: Pseudocone Functors: Monoidal Closed Derived Cat Topological}. By the assumptions on $G$ it follows that for any $n \in \N$ if $M$ is a specified $n$-acyclic free $G$-space which is a manifold then the second projection $\pi_2:M \times X \to X$ is a smooth $n$-acyclic resolution of $X$. From here taking the subcategory $\mathbf{SResl}_G(X)$ of $G$-resolutions of $X$ generated by the projections $M \times X \to X$, their pullbacks against each other, and the smooth maps between them, we get that there is an equivalence 
	\[
	D_G^b(X) \simeq \PC(D^b_c(G \backslash (-)):\mathbf{SResl}_G(X)^{\op} \to \fCat)
	\]
	by \cite[Section 2.4.4]{BernLun}. However, because every resolution $p^{\ast}:P \to X$ in $\mathbf{SResl}_G(X)$ is smooth and every morphism $f:P \to Q$ between such resolution is smooth, the pullback functors $p^{\ast}$ and $f^{\ast}$ all commute with suspension functors $[1]_{X}$. As such we can apply Theorem \ref{Theorem: Section Triangle: Equivariant triangulation} and translate the triangulation across the equivalence to deduce the fact that $D_G^b(X)$ is triangulated.
\end{proof}

We will need $t$-structures for what follows because arguably the crucial feature of the derived category is the presence of the universal property (as a localization) and the presence of both perverse and standard $t$-structures. As such, we recall their definition here. 
\begin{definition}[{\cite[Definition 1.3.1]{BBD}}]
	Let $\Tscr$ be a triangulated category. A $t$-structure on $\Tscr$\index[terminology]{Tstructure@$t$-structure} is a pair of full subcategories $(\Tscr^{\leq 0}, \Tscr^{\geq 0})$ which are both stable under equivalence such that the following are satisfied (where for any $n \in \Z$, we define $\Tscr^{\leq n} := \Tscr^{\leq 0}[-n]$ and $\Tscr^{\geq n} := \Tscr^{\geq 0}[-n]$):
	\begin{itemize}
		\item  If $A$ is an object of $\Tscr^{\geq 1}$ and $B$ is an object of $\Tscr^{\leq 0}$, then $\Tscr(A,B) = 0$;
		\item We have inclusions of categories $\Tscr^{\leq 0} \hookrightarrow \Tscr^{\leq 1}$ and $\Tscr^{\geq 1} \hookrightarrow \Tscr^{\geq 0}$;
		\item For any object $A \in \Tscr_0$ there exists a truncation distinguished triangle
		\[
		\xymatrix{
			X \ar[r] & A \ar[r] & Y \ar[r] & X[1]
		}
		\]
		where $X$ is an object of $\Tscr^{\leq 0}$ and $Y$ is an object of $\Tscr^{\geq 1}$.
	\end{itemize}
	Furthermore, given a $t$-structure $(\Tscr^{\leq 0}, \Tscr^{\geq 0})$, we define the heart of $\Tscr$ to be the category
	\[
	\Tscr^{\heartsuit} := \Tscr^{\leq 0} \cap \Tscr^{\geq 0}.
	\]
\end{definition}
\begin{example}\label{Example: Standard tstructure}
	Let $\Ascr$ be an Abelian category and let $D(\Ascr)$ be the derived category of chain complexes of $\Ascr$. Define the subcategory $D(\Ascr)^{\leq 0}$ (respectivly $D(\Ascr)^{\geq 0}$) to be the subcategory of $D(\Ascr)$ generated by the complexes $X^{\bullet}$ for which $H^{n}(X^{\bullet}) = 0$ for all $n > 0$ (respectively generated by the complexes $X^{\bullet}$ for which $H^n(X^{\bullet}) = 0$ for all $n < 0$). The pair $(D(\Ascr)^{\leq 0}, D(\Ascr)^{\geq 0})$ is a $t$-structure on $D(\Ascr)$ and its heart is the subcategory of objects whose cohomology is concentrated in degree $0$. In particular, there is an equivalence
	\[
	D(\Ascr)^{\heartsuit} \simeq \Ascr.
	\]
\end{example}
\begin{remark}
	For any $t$-structure on a triangulated category $\Tscr$, the heart $\Tscr^{\heartsuit}$ is an Abelian category by \cite[Th{\'e}or{\'e}me 1.3.6]{BBD}. If there are multiple $t$-structures on a category $\Tscr$, we will denote their corresponding subcategories as $({}^{t}\Tscr^{\leq 0}, {}^{t}\Tscr^{\geq 0})$ and $({}^{t^{\prime}}\Tscr^{\leq 0}, {}^{t^{\prime}}\Tscr^{\geq 0})$ while writing the corresponding hearts as $\Tscr^{\heartsuit_t}$ and $\Tscr^{\heartsuit_{t^{\prime}}}$, respectively.
\end{remark}
\begin{definition}\label{Defn: Standard tstructure}
	Let $\Dscr$ be a (bounded) derived category of some Abelian category or the bounded derived category of $\ell$-adic sheaves on some variety. We define the standard $t$-structure on $\Dscr$ to be the $t$-structure generated as in Example \ref{Example: Standard tstructure}.
\end{definition}

We now show how to build $t$-structures on $\PC(F)$ based on $t$-structures that come from those on $F(X)$. After proving the theorem (cf.\@ Theorem \ref{Theorem: Section Triangle: t-structure on our friend the ECat}) that says these are indeed $t$-structures on $\PC(F)$, we will then give some explicit examples of these $t$-structures as they may arise in nature (as much as perverse sheaves are natural, anyway). However, before proving this theorem, we will define truncated pseudofunctors, which essentially are triangulated pseudofunctors for which each fibre functor admits a $t$-structure and these $t$-structures vary pseudofunctorially through $F$. With this we will prove how to produce $t$-structures from pseudofunctors and also how a base category $\Tscr$ can give rise to $t$-structures on $\PC(F)$.

\begin{definition}\label{Defn: Truncated preq pseudofunctor}\index[terminology]{Pre-equivariant Pseudofunctor!Truncated}
	Let $F:\Cscr^{\op} \to \fCat$ be a triangulated pseudofunctor. Assume each fibre category has $t$-structure with truncation functors $(\quot{\tau^{\leq 0}}{X},\quot{\tau^{\geq 0}}{X})$ such that these give rise to pseudofunctors $F^{\leq 0}:\Cscr^{\op} \to \fCat$\index[notation]{Fleq0@$F^{\leq 0}$} and $F^{\geq 0}:\Cscr^{\op} \to \fCat$\index[notation]{Fgeq0@$F^{\geq 0}$} for which the functors $F^{\geq 0}(f)$ and $F^{\leq 0}(f)$ are $t$-exact functors which coincide on all common subcategories. Then if for all morphisms $f:X \to Y$ in $\Cscr$ there are natural isomorphisms
	\[
	\begin{tikzcd}
		F(Y) \ar[r, ""{name = U}]{}{F(f)} \ar[d, swap]{}{\quot{\tau^{\leq 0}}{Y}} & F(X) \ar[d]{}{\quot{\tau^{\leq 0}}{X}} \\
		F(Y)^{\leq 0} \ar[r, swap, ""{name = L}]{}{F^{\leq 0}(f)} & F(X)^{\leq 0} \ar[from = U, to = L, Rightarrow, shorten <= 4pt, shorten >= 4pt]{}{\theta_{f}^{\leq 0}}
	\end{tikzcd}
	\]
	and
	\[
	\begin{tikzcd}
		F(Y) \ar[r, ""{name = U}]{}{F(f)} \ar[d, swap]{}{\quot{\tau^{\geq 0}}{Y}} & F(X) \ar[d]{}{\quot{\tau^{\geq 0}}{X}} \\
		F(Y)^{\geq 0} \ar[r, swap, ""{name = L}]{}{F^{\geq 0}(f)} & F(X)^{\geq 0} \ar[from = U, to = L, Rightarrow, shorten <= 4pt, shorten >= 4pt]{}{\theta_{f}^{\geq 0}}
	\end{tikzcd}
	\]
	which satisfy the coherence diagrams presented below, we say that the pseudofunctor is truncated. The pasting diagram
	\[
	\begin{tikzcd}
		F(Z) \ar[rr, ""{name = UL}]{}{F(g)} \ar[d, swap]{}{\quot{\tau^{\leq 0}}{Z}} \ar[rrrr, bend left = 30, ""{name = U}]{}{F(g \circ f)} &  & F(Y) \ar[d]{}{\quot{\tau^{\leq 0}}{Y}} \ar[rr, ""{name = UR}]{}{F(f)} & & F(X) \ar[d]{}{\quot{\tau^{\leq 0}}{X}} \\
		F(Z)^{\leq 0} \ar[rr, swap, ""{name = LL}]{}{F^{\leq 0}(g)} \ar[rrrr, bend right = 30, swap, ""{name = L}]{}{F^{\leq 0}(g \circ f)} & & F(Y)^{\leq 0} \ar[rr, swap, ""{name = LR}]{}{F^{\leq 0}(f)} & & F(X)^{\leq 0} \ar[from = U, to = 1-3, Rightarrow, shorten <= 4pt, shorten <= 4pt]{}{\quot{\phi_{f,g}}{F}^{-1}} \ar[from = UR, to = LR, Rightarrow, shorten <= 4pt, shorten >= 4pt]{}{\theta_{f}^{\leq 0}} \ar[from = UL,, to = LL, Rightarrow, shorten <= 4pt, shorten >= 4pt]{}{\theta_g^{\leq 0}} \ar[from = 2-3, to = L, shorten <= 4pt, shorten >= 4pt, Rightarrow]{}{\quot{\phi_{f,g}}{F^{\leq 0}}}
	\end{tikzcd}
	\]
	is equal to the $2$-cell
	\[
	\begin{tikzcd}
		F(Z) \ar[rr, ""{name = U}]{}{F(g \circ f)} \ar[d, swap]{}{\quot{\tau^{\leq 0}}{Z}} & & F(X) \ar[d]{}{\quot{\tau^{\leq 0}}{X}} \\
		F(Z)^{\leq 0} \ar[rr, swap, ""{name = L}]{}{F^{\leq 0}( g \circ f)} & & F(X)^{\leq 0} \ar[from = U, to = L, Rightarrow, shorten <= 4pt, shorten >= 4pt]{}{\theta_{g \circ f}^{\leq 0}}
	\end{tikzcd}
	\]
	and dually the pasting diagram
	\[
	\begin{tikzcd}
		F(Z) \ar[rr, ""{name = UL}]{}{F(g)} \ar[d, swap]{}{\quot{\tau^{\geq 0}}{Z}} \ar[rrrr, bend left = 30, ""{name = U}]{}{F(g \circ f)} &  & F(Y) \ar[d]{}{\quot{\tau^{\geq 0}}{Y}} \ar[rr, ""{name = UR}]{}{F(f)} & & F(X) \ar[d]{}{\quot{\tau^{\geq 0}}{X}} \\
		F(Z)^{\geq 0} \ar[rr, swap, ""{name = LL}]{}{F^{\geq 0}(g)} \ar[rrrr, bend right = 30, swap, ""{name = L}]{}{F^{\geq 0}(g \circ f)} & & F(Y)^{\geq 0} \ar[rr, swap, ""{name = LR}]{}{F^{\geq 0}(f)} & & F(X)^{\geq 0} \ar[from = U, to = 1-3, Rightarrow, shorten <= 4pt, shorten <= 4pt]{}{\quot{\phi_{f,g}}{F}^{-1}} \ar[from = UR, to = LR, Rightarrow, shorten <= 4pt, shorten >= 4pt]{}{\theta_{f}^{\geq 0}} \ar[from = UL,, to = LL, Rightarrow, shorten <= 4pt, shorten >= 4pt]{}{\theta_g^{\geq 0}} \ar[from = 2-3, to = L, shorten <= 4pt, shorten >= 4pt, Rightarrow]{}{\quot{\phi_{f,g}}{F^{\geq 0}}}
	\end{tikzcd}
	\]
	is equal to the $2$-cell:
	\[
	\begin{tikzcd}
		F(Z) \ar[rr, ""{name = U}]{}{F(g \circ f)} \ar[d, swap]{}{\quot{\tau^{\geq 0}}{Z}} & & F(X) \ar[d]{}{\quot{\tau^{\geq 0}}{X}} \\
		F(Z)^{\geq 0} \ar[rr, swap, ""{name = L}]{}{F^{\geq 0}( g \circ f)} & & F(X)^{\geq 0} \ar[from = U, to = L, Rightarrow, shorten <= 4pt, shorten >= 4pt]{}{\theta_{g \circ f}^{\geq 0}}
	\end{tikzcd}
	\]
\end{definition}
\begin{proposition}\label{Prop: Section Triangle: Truncated preeq gives truncation functors}
	Let $F:\Cscr^{\op} \to \fCat$ be a truncated pseudofunctor. Then there are truncation functors
	\[
	\tau^{\leq 0}:\PC(F) \to \PC(F), \qquad \tau^{\geq 0}:\PC(F) \to \PC(F).
	\]
	The truncation $\tau^{\leq 0}$ is defined on objects by
	\[
	\tau^{\leq 0}\left(\left\lbrace \quot{A}{X}\; : \; X \in \Cscr_0 \right\rbrace\right) = \left\lbrace \quot{\tau^{\leq 0}}{X}(\quot{A}{X}) \; : \; X \in \Cscr_0 \right\rbrace,
	\]
	with transition isomorphisms $\tau_f^{A\leq 0}$ given by the diagram
	\[
	\xymatrix{
		F(f)\left(\quot{\tau^{\leq 0}}{Y}\quot{A}{Y}\right) \ar[rr]^-{\theta_f^{\leq 0}} \ar[drr]_{\tau_f^{A \leq 0}} & & \quot{\tau^{\leq 0}}{X}\left(F(f)\quot{A}{Y}\right) \ar[d]^{\quot{\tau^{\leq 0}}{X}\left(\tau_f^A\right)} \\
		& & \quot{\tau^{\leq 0}}{X}\left(\quot{A}{X}\right)
	}
	\]
	and defined on morphisms by
	\[
	\tau^{\leq 0}\left(\lbrace \quot{\rho}{X} \; | \; X \in \Cscr_0 \rbrace\right) = \left\lbrace \quot{\tau^{\leq 0}}{X}\left(\quot{\rho}{X}\right) \; | \; X \in \Cscr_0 \right\rbrace.
	\]
	Similarly, the other truncation $\tau^{\geq 0}$ is defined on objects by
	\[
	\tau^{\geq 0}\left(\left\lbrace \quot{A}{X}\; : \; X \in \Cscr_0 \right\rbrace\right) = \left\lbrace \quot{\tau^{\geq 0}}{X}\left(\quot{A}{X}\right) \; | \; X \in \Cscr_0 \right\rbrace,
	\]
	with transition isomorphisms $\tau_f^{A\geq 0}$ given by the diagram
	\[
	\xymatrix{
		F(f)\left(\quot{\tau^{\geq 0}}{Y}\quot{A}{Y}\right) \ar[rr]^-{\theta_f^{\geq 0}} \ar[drr]_{\tau_f^{A\leq 0}} & & \quot{\tau^{\geq 0}}{X}\left(F(f)\quot{A}{Y}\right) \ar[d]^{\quot{\tau^{\geq 0}}{X}\left(\tau_f^A\right)} \\
		& & \quot{\tau^{\geq 0}}{X}\left(\quot{A}{X}\right)
	}
	\]
	and defined on morphisms by
	\[
	\tau^{\geq 0}\left(\lbrace \quot{\rho}{X} \; | \; X \in \Cscr_0 \rbrace\right) = \left\lbrace \quot{\tau^{\geq 0}}{X}\left(\quot{\rho}{X}\right) \; | \; X \in \Cscr_0 \right\rbrace.
	\]
\end{proposition}
\begin{proof}
	The technical conditions of this proposition are formally similar to the tensor functors in Theorem \ref{Theorem: Section 2: Monoidal preequivariant pseudofunctor gives monoidal equivariant cat}. As such, proceeding mutatis mutandis as in the proof of Theorem \ref{Theorem: Section 2: Monoidal preequivariant pseudofunctor gives monoidal equivariant cat} proves the proposition.
\end{proof}
By taking appropriate shifts and setting $\PC(F)^{\leq n} := \PC(F)^{\leq 0}[-n]$ together with $\PC(F)^{\geq n} := \PC(F)^{\geq 0}[-n]$ for $n \in \Z$, we obtain the following corollary.
\begin{corollary}\label{Cor: Section Triangle: Truncation functors for all $n$}
	Let $F$ be a truncated pseudofunctor on $X$. Then for all $n \in \Z$ there are subcategories $\PC(F)^{\geq n}$ and $\PC(F)^{\leq n}$ of $\PC(F)$ equipped with truncation functors $\tau^{\geq n}:\PC(F) \to \PC(F)^{\geq n}$ and $\tau^{\leq n}:\PC(F) \to \PC(F)^{\leq n}$. As usual, these truncation functors satisfy the adjunctions:
	\[
	\begin{tikzcd}
		\PC(F)^{\geq n} \ar[rr, bend right = 20, swap, ""{name = D}]{}{\incl_{\geq n}} & & \PC(F) \ar[ll, bend right = 20, swap, ""{name = U}]{}{\tau^{\geq n}} \ar[from = U, to = D, symbol = \dashv]
	\end{tikzcd}\qquad
	\begin{tikzcd}
		\PC(F) \ar[rr, bend right = 20, swap, ""{name = D}]{}{\tau^{\leq n}} & & \PC(F)^{\leq n} \ar[ll, bend right = 20, swap, ""{name = U}]{}{\incl_{\leq n}} \ar[from = U, to = D, symbol = \dashv]
	\end{tikzcd}
	\]
\end{corollary}
\begin{proof}
	The existence of the subcategories $\PC(F)^{\leq n}$ and $\PC(F)^{\geq n}$ together with the inclusion functors $\incl_{\leq n}:\PC(F)^{\leq n} \to \PC(F)$ and $\incl_{\geq n}:\PC(F)^{\geq n}:\PC(F)^{\geq n} \to \PC(F)$ follow from Proposition \ref{Prop: Section Triangle: Truncated preeq gives truncation functors} by taking negative shifts of the subcategories and inclusions of $\PC(F)^{\leq 0} \to \PC(F)$ and $\PC(F)^{\geq 0} \to \PC(F)$. The existence of the truncations follow from the fact that $\tau^{\geq n}$ and $\tau^{\leq n}$ correspond, from the definition of a truncated pseudofunctor and the fact that each fibre functor $F(f)$ is $t$-exact, to degree $-n$ shifts of the pseudonatural transformations $\tau^{\geq 0}$ and $\tau^{\leq 0}$. Consequently they arise from the pseudonatural transformations
	\[
	\begin{tikzcd}
		\Cscr^{\op} \ar[rr, bend left = 30, ""{name = U}]{}{F} \ar[rr, bend right = 30, swap, ""{name = D}]{}{F^{\geq n}} & & \fCat \ar[from = U, to = D, Rightarrow, shorten <= 4pt, shorten >= 4pt]{}{\tau^{\geq n}}
	\end{tikzcd}
	\]
	and:
	\[
	\begin{tikzcd}
		\Cscr^{\op} \ar[rr, bend left = 30, ""{name = U}]{}{F} \ar[rr, bend right = 30, swap, ""{name = D}]{}{F^{\leq n}} & & \fCat  \ar[from = U, to = D, Rightarrow, shorten <= 4pt, shorten >= 4pt]{}{\tau^{\leq n}}
	\end{tikzcd}
	\]
	From here the adjunctions $\incl_{\geq n} \dashv \tau^{\geq n}$ and $\tau^{\leq n} \dashv \incl_{\leq n}$ follow from the $X$-local adjunctions $\quot{\incl_{\geq n}}{X} \dashv \quot{\tau^{\geq n}}{X}$ and $\quot{\tau^{\leq n}}{X} \dashv \quot{\incl_{\leq n}}{X}$ together with an application of Theorem \ref{Thm: Section 3: Gamma-wise adjoints lift to equivariant adjoints}.
\end{proof}

Using these truncation functors allows us to deduce that truncated psuedofunctors give rise to $t$-structures on $\PC(F)$.

\begin{Theorem}\label{Theorem: Section Triangle: t-structure on our friend the ECat}
	Let $F:\Cscr^{\op} \to \fCat$ be a truncated pseudofunctor. The truncation functors, when shifted for all $n \in \Z$, give rise to a $t$-structure on $\PC(F)$ for which the category of bounded above by degree $0$ objects of $\PC(F)$ is the $\PC(F^{\leq 0})$ while the category of objects of $\PC(F)$ bounded below by degree $0$ is the category $\PC(F^{\geq 0})$.
\end{Theorem}
\begin{remark}
	In the situation of Theorem \ref{Theorem: Section Triangle: t-structure on our friend the ECat}, we will freely interchange the notations
	\[
	\PC(F^{\leq 0}) = \PC(F)^{\leq 0},\qquad \PC(F^{\geq 0}) = \PC(F)^{\geq 0}.
	\]
	depending on context and how we want to think of each category at any given moment.
\end{remark}
\begin{proof}
	We first prove that we can realize $\PC(F^{\leq 0})$ and $\PC(F)$ as \emph{full} subcategories of $\PC(F)$; note that they are subcategories has already been shown. For this we observe first that since each category $F(X)^{\leq 0}$ and $F(X)^{\geq 0}$ is a strictly full subcategory of each fibre category $F(X)$. Now let $A \in \PC(F^{\geq 0})$; the case for $\PC(F^{\leq 0})$ follows dually and so is omitted here. Because $A \in \PC(F^{\geq 0})$, for all $X \in \Cscr_0$ we have that 
	\[
	\quot{A}{X} \in F^{\geq 0}(X)_0 = \left(F(X)^{\geq 0}\right)_{0}
	\] 
	so in turn each $A \in \PC(F)_0$. Now any morphism $P:A \to B$ in $\PC(F)$ with $A,B \in \PC(F^{\geq 0})_0$ has $X$-local components 
	\[
	\quot{\rho}{X}:\quot{A}{X} \to \quot{B}{X}. 
	\]
	Because $F(X)^{\geq 0}$ is full, we get that $\quot{\rho}{X} \in (F(X)^{\geq 0})_1$, and because this happens for each $X \in \Cscr_0$, it follows that $P \in \PC(F^{\geq 0})_1$ as well. This shows the fullness of $\PC(F^{\geq 0})$; that it is strictly full (i.e., the fact that the subcategory is closed with respect to isomorphisms) follows similarly. We thus conclude that $\PC(F^{\geq 0})$ and $\PC(F^{\leq 0})$ are both strictly full subcategories of $\PC(F)$.
	
	We now prove that $\PC(F)^{\geq 0}$ and $\PC(F)^{\leq 0}$ do indeed induce a $t$-structure on $\PC(F)$. For this we recall that the category $\PC(F)^{\geq n}$, for any $n \in \Z$, is induced by the construction on objects
	\[
	\left(\PC(F)^{\geq n}\right)_0 = \PC(F)^{\geq 0}[-n] = \left\lbrace A[-n] \; : \; A \in \PC(F^{\geq 0})_0 \right\rbrace.
	\]
	We now establish that for any $A \in \PC(F^{\leq 0})$ and any $B \in \PC(F)^{\geq 1}$, the hom-set $\PC(F)(A,B) = 0$. For this we note that since $A \in (\PC(F)^{\leq 0})_0$ and $B \in (\PC(F)^{\geq 1})_0$, for every $X \in \Cscr$ we have that $\quot{A}{X} \in (F(X)^{\leq 0})_0$ and $\quot{B}{X} \in (F(X)^{\geq 1})_0$. Thus we get that
	\[
	F(X)(\quot{A}{X},\quot{B}{X}) = \left\lbrace\quot{0}{X}\right\rbrace
	\]
	so we calculate that
	\begin{align*}
		0  &= \lbrace \quot{0}{X} \; | \; X \in \Cscr_0 \rbrace \subseteq \PC(F)(A,B) \\
		&\subseteq \bigcup_{X \in \Cscr_0}F(X)\left(\quot{A}{X},\quot{B}{X}\right) = \bigcup_{X \in \Cscr_0} \lbrace \quot{0}{X} \rbrace = \lbrace \quot{0}{X} \; | \; X \in \Cscr_0 \rbrace = 0,
	\end{align*}
	which establishes the claim.
	
	We now prove that the categories $\PC(F^{\leq 0})$ and $\PC(F^{\geq 0})$ are stable under shifts in the sense that $\PC(F^{\leq -1})$ is a subcategory of $\PC(F^{\leq 0})$ and $\PC(F^{\geq 0})$ is a subcategory of $\PC(F^{\geq -1})$. However, because these identities hold $X$-locally and the shift functors are defined by shifting $X$-locally, this is a trivial observation.
	
	We now prove that for any object $A$ that there is a distinguished triangle in $\PC(F)$ of the form
	\[
	\xymatrix{
		Z \ar[r] & A \ar[r] & Y \ar[r] & Z[1]
	}
	\]
	with $Z \in \PC(F^{\leq 0})_0$ and $Y \in \PC(F^{\geq 1})_0$. For this, however, we define $Z := \tau^{\leq 0}(A)$, $Y := \tau^{\geq 1}(A)$, and let the morphisms be induced from the $X$-local distinguished triangles
	\[
	\xymatrix{
		\quot{\tau^{\leq 0}}{X}\left(\quot{A}{X}\right) \ar[r] & \quot{A}{X} \ar[r] & \quot{\tau^{\geq 1}}{X}\left(\quot{A}{X}\right) \ar[r] & \quot{\tau^{\leq 0}\left(\quot{A}{X}\right)[1]}{X}
	}
	\]
	in each fibre category $F(X)$. However, because each of these triangles are distinguished we get that the triangle
	\[
	\xymatrix{
		\tau^{\leq 0}(A) \ar[r] & A \ar[r] & \tau^{\geq 1}(A) \ar[r] & \tau^{\leq  0}(A)[1]
	}
	\]
	is distinguished in $\PC(F)$, which in turn completes the proof of the claim and hence the theorem.
\end{proof}
\begin{example}
	Let $G$ be a smooth algebraic group and let $X$ be a left $G$-variety. We now give two different (families of) examples of truncated pseudofunctors $F:\SfResl_G(X)^{\op} \to \fCat$ by using the pseudofunctor $F = \DbQl{G\backslash(-)}$. Note that in in principle we can give four (families of) examples by using the same $t$-structures on $D_c^b(-)$ and proceeding mutatis mutandis. We can even get six (families of) examples in the topological case by letting $G$ be either a (closed) linear subgroup of $\GL_n(\R)$ for some $n \in \N$ or a Lie group with finitely many connected components and letting $X$ be a left $G$-space.
	
	The first truncated pseudofunctor we describe fully is the standard $t$-structure on $\DbQl{-}$ (cf.\@ Definition \ref{Defn: Standard tstructure}). For each $\Gamma \in \Sf(G)_0$ we say that
	\[
	F^{\leq 0}(\XGamma) := {}^{\text{stand}}\DbQl{\XGamma}^{\leq 0}
	\]
	and
	\[
	F^{\geq 0}(\XGamma) := {}^{\text{stand}}\DbQl{\XGamma}^{\geq 0}
	\]
	where stand denotes the standard $t$-structure on $\DbQl{-}$ and the truncations $\quot{\tau^{\leq 0}}{\Gamma}$ and $\quot{\tau^{\geq 0}}{\Gamma}$ are given from the standard $t$-structure on the derived $\ell$-adic category. Then because the morphism $\of$ is smooth, the functor $\of^{\ast}$ is $t$-exact for the standard $t$-structure and so we set $F^{\leq 0}(\of) = F^{\geq 0}(\of) = \of^{\ast}$. Then the diagrams
	\[
	\xymatrix{
		\DbQl{\XGammap} \ar[d]_{\quot{\tau^{\leq 0}}{\Gamma^{\prime}}} \ar[r]^-{\of^{\ast}} & \DbQl{\XGamma} \ar[d]^{\quot{\tau^{\leq 0}}{\Gamma}} \\
		{}^{\text{stand}}\DbQl{\XGamma}^{\leq 0} \ar[r]_-{\of^{\ast}} & {}^{\text{stand}}\DbQl{\XGamma}
	}
	\]
	\[
	\xymatrix{
		\DbQl{\XGammap} \ar[d]_{\quot{\tau^{\geq 0}}{\Gamma^{\prime}}} \ar[r]^-{\of^{\ast}} & \DbQl{\XGamma} \ar[d]^{\quot{\tau^{\geq 0}}{\Gamma}} \\
		{}^{\text{stand}}\DbQl{\XGammap}^{\geq 0} \ar[r]_-{\of^{\ast}} & {}^{\text{stand}}\DbQl{\XGamma}^{\geq 0}
	}
	\]
	commute strictly, which implies that in this case the trucantion coherences are satisfied automatically. This is the standard $t$-structure on $\DbeqQl{X}$ and it is routine (cf.\@ Corollary \ref{Cor: Section Triangle: Pseudofunctor of standard hearts is heart of pseudofunctor})to see that $\DbeqQl{X}^{\heartsuit_{\text{stand}}} \simeq \Shv_G(X;\overline{\Q}_{\ell})$.\index[notation]{DbGQLStand@${}^{\text{stand}}\DbeqQl{X}^{\heartsuit}$}\index[terminology]{Standard $t$-structure!On DbG@$\DbeqQl{X}$}
	
	The second truncated pre-equivariant pseudofunctor we construct on $\DbeqQl{X}$ is induced from the perverse $t$-structure. For each $\Gamma \in \Sf(G)_0$ we define
	\[
	F^{\leq 0}(\XGamma) = {}^{p}\DbQl{\XGamma}^{\leq 0}
	\]
	and
	\[
	F^{\geq 0}(\XGamma) = {}^{p}\DbQl{\XGamma}^{\geq 0}.
	\]
	Now for each morphism $f \in \Sf(G)_1$ let $d_f \in \Z$ denote the relative (pure) dimension of $\of$. Then because $\of$ is smooth, we have from \cite[Page 108]{BBD} that $\of^{\ast}[d_f]$ is a $t$-exact functor for the pervese $t$-structure, so we set for all $f \in \Sf(G)_1$
	\[
	F^{\geq 0}(\of) = F^{\leq 0}(\of) = \of^{\ast}[d_f].
	\]
	Let each truncation functor $\quot{\tau^{\geq 0}}{\Gamma}, \quot{\tau^{\leq 0}}{\Gamma}$ be truncations for the perverse $t$-structure; we will omit the ${(-)}^{p}$ left-handed superscript to reduce notational clutter. Then a routine calculation shows that for any $f \in \Sf(G)_1$ the diagram
	\[
	\begin{tikzcd}
		\DbQl{\XGammap} \ar[r, ""{name = LeftU}]{}{\of^{\ast}} \ar[d, swap]{}{\quot{\tau^{\leq 0}}{\Gamma^{\prime}}} &\DbQl{\XGamma} \ar[d]{}{\quot{\tau^{\leq 0}}{\Gamma}}  \\
		{}^{p}\DbQl{\XGammap}^{\leq 0} \ar[r, swap, ""{name= LeftL}]{}{\of^{\ast}[d_f]} & {}^{p}\DbQl{\XGamma}^{\leq 0} \ar[from = LeftU, to = LeftL, Rightarrow, shorten <= 4pt, shorten >= 4pt] 
	\end{tikzcd}
	\]
	\[
	\begin{tikzcd}
		\DbQl{\XGammap} \ar[r, ""{name = RightU}]{}{\of^{\ast}} \ar[d, swap]{}{\quot{\tau^{\geq 0}}{\Gamma^{\prime}}} & \DbQl{\XGamma} \ar[d]{}{\quot{\tau^{\geq 0}}{\Gamma}} \\
		{}^{p}\DbQl{\XGammap}^{\geq 0} \ar[r, swap, ""{name = RightL}]{}{\of^{\ast}[d_f]} & {}^{p}\DbQl{\XGamma}^{\geq 0} \ar[from = RightU, to = RightL, Rightarrow, shorten <= 4pt, shorten >= 4pt]
	\end{tikzcd}
	\]\index[notation]{DbGQLPer@$\DbeqQl{X}^{\heartsuit_{\text{per}}}$}
	commutes up to an invertible $2$-cell because of the $t$-exactness\footnote{A functor $F:\Ascr \to \Bscr$ of triangulated categories is $t$-exact if it is triangulated, $F(\Ascr^{\leq 0})$ is a subcategory of $\Bscr^{\leq 0}$, and $F(\Ascr^{\geq 0})$ is a subcategory of $\Bscr^{\geq 0}$.} of $\of^{\ast}[d_f]$, where the commutativity isomorphisms are induced based on how we identify and move between pullbacks. This defines the equivariant perverse $t$-structure on $\DbeqQl{X}$ and Theorem \ref{Theorem: Heart of equivariant t-structure is equivariant heart of t-structure} below shows that 
	\[
	(\DbeqQl{X})^{\heartsuit_{\text{per}}} = \Per_G(X).
	\]
\end{example}

We now observe that if $F:\Cscr^{\op} \to \fCat$ is a truncated pseudofunctor, because $\PC(F)$ admits a $t$-structure the full subcategory $\PC(F)^{\heartsuit}$ of $\PC(F)$ induced by object assignment 
\[
\PC(F)^{\heartsuit}_0 := \PC(F)^{\leq 0}_0 \cap \PC(F)^{\geq 0}_0
\]
is an Abelian subcategory of $\PC(F)$ with cohomological functor 
\[
\tau^{\geq 0} \circ \incl_{\leq 0} \circ \tau^{\leq 0}:\PC(F) \to \PC(F)^{\heartsuit}
\]
by \cite[Th{\'e}or{\`e}me 1.3.6]{BBD}\footnote{In \cite{BBD} the statement listed is that $\tau_{\geq 0} \circ \tau_{\leq 0}$ is the cohomological functor; however, because our functors $\tau_{\leq n}$ and $\tau_{\geq n}$ are typed to have domain $\PC(F)$, we require an inclusion to the category $\PC(F)$ for things to type correctly. This will not cause any significant issue for what we do, but is worth pointing out for the curious reader.}. A natural question then is whether or not the category $\PC(F)^{\heartsuit}$ is itself equivalent to a category of pseudocones for some pseudofunctor $E:\Cscr^{\op} \to \fCat$. We show that this is in fact the case below, and in a surprising way. The categories $F^{\heartsuit}(X)$ vary pseudofunctorially in $\Cscr^{\op}$ because the fibre functors $F(f)$ are $t$-exact and so form a pseudofunctor $F^{\heartsuit}:\Cscr^{\op} \to \fCat$. We show below that it is in fact the case that $\PC(F)^{\heartsuit} = \PC(F^{\heartsuit})$; note the strict equality!

\begin{Theorem}[Change of Heart]\label{Theorem: Heart of equivariant t-structure is equivariant heart of t-structure}
	Let $F:\Cscr^{\op} \to \fCat$ be a triangulated pseudofunctor. Define the Abelian pseudofunctor $F^{\heartsuit}:\Cscr^{\op} \to \fCat$ by setting
	\[
	F(X) := F^{\geq 0}(X) \cap F^{\leq 0}(X) = F(X)^{\heartsuit}
	\]
	and defining the fibre functors by 
	\[
	\left(F^{\geq 0}\right)|_{F^{\geq 0}(X) \cap F^{\leq 0}(X)}(f) =: F^{\heartsuit}(f).
	\]
	Then \index[notation]{FGXHeart@$F_G(X)^{\heartsuit}$}\index[notation]{FHeartGX@$(F^{\heartsuit})_G(X)$}
	\[
	\PC(F)^{\heartsuit} = \PC\left(F^{\heartsuit}\right).
	\]
\end{Theorem}
\begin{remark}
	By Definition \ref{Defn: Section Triangle: Triangulated Pseudofunctor} the definition of the fibre functor $F^{\heartsuit}(f)$
	\[
	\left(F^{\geq 0}\right)|_{F^{\geq 0}(X) \cap F^{\leq 0}(X)}(f) =: F^{\heartsuit}(f)
	\]
	does not depend seriously on $F^{\geq 0}$; we could have also chosen to use $F^{\leq 0}(f)$ instead as
	\[
	\left(F^{\geq 0}\right)|_{F^{\geq 0}(X) \cap F^{\leq 0}(X)}(f) = \left(F^{\leq 0}\right)|_{F^{\geq 0}(X) \cap F^{\leq 0}(X)}(f).
	\]
	The decision to use $F^{\geq 0}$ was simply to provide an actual explicit choice to ensure that $F^{\heartsuit}(f)$ exists.
\end{remark}
\begin{proof}
	Begin by observing that on one hand an object $A \in \PC(F^{\heartsuit})_0$ if and only if for all $X \in \Cscr_0$ the object $\quot{A}{X} \in F(X)^{\heartsuit}$ and if and only if each transition isomorphism $\tau_f^A:F(f)(\quot{A}{Y}) \xrightarrow{\cong} \quot{A}{X}$ is an isomorphism in $F(X)^{\heartsuit}$. On the other hand, by construction of the $t$-structure on $\PC(F)$ given in Theorem \ref{Theorem: Section Triangle: t-structure on our friend the ECat} an object $A$ in $\PC(F)$ is in $(\PC(F))^{\heartsuit}_0$ as well if and only if $A \in (\PC(F^{\geq 0}) \cap \PC(F^{\leq 0}))_0$. However, unwinding this definition we find that $A$ lies in the heart of this $t$-structure if and only if $X$-locally we have that
	\[
	\quot{A}{X} \in \left(F(X)^{\leq 0} \cap F(X)^{\geq 0}\right)_0.
	\]
	This in turn implies that for all $X \in \Cscr_0$ we must have that $\quot{A}{X} \in F(X)^{\heartsuit}_0$, so the only aspect of morphisms that may differ between $\PC(F)^{\heartsuit}$ and $\PC(F^{\heartsuit})$ lies in the transition isomorphisms. Because $\tau_f^{A}:F(f)(\quot{A}{Y}) \xrightarrow{\cong} \quot{A}{X}$ is a morphism in both $F(X)^{\geq 0}$ and $F(X)^{\leq 0}$, the isomorphism $\tau_f^A$ is an isomorphism in $F(X)^{\geq 0} \cap F(X)^{\leq 0}$ and hence in the heart of $F(X)$. As such, we find that $\PC(F^{\heartsuit})_0 = (\PC(F)^{\heartsuit})_0$. The verification that $(F^{\heartsuit})_{G}(X)_1 = (\PC(F))^{\heartsuit}_1$ follows a similar analysis as what we just performed using that the fibre functors all satisfy
	\[
	F^{\heartsuit}(f) = \left(F^{\geq 0}\right)|_{F^{\geq 0}(X) \cap F^{\leq 0}(X)}(f) = \left(F^{\leq 0}\right)|_{F^{\geq 0}(X) \cap F^{\leq 0}(X)}(f)
	\]
	is $t$-exact, allowing us to conclude that $\PC(F^{\heartsuit})_1 = (\PC(F)^{\heartsuit})_1$ as well. Finally, since composition, source, and target maps must also coincide by construction, it follows that $\PC(F)^{\heartsuit} = \PC(F^{\heartsuit})$.
\end{proof}
\begin{remark}
	The Change of Heart Theorem shows the heart of the standard $t$-structure on $\PC(F)$ coincides with the category of pseudocones with values in $F^{\heartsuit}$. 
\end{remark}

We now show that when combined with some specific knowledge of the equivariant derived category of a variety, our pseudocone formalism gives rise to equivalences of categories $D_G^b(X;\overline{\Q}_{\ell})^{\heartsuit_{\operatorname{per}}} \simeq \Per_G(X;\overline{\Q}_{\ell})$ and $D_G^b(X;\overline{\Q}_{\ell})^{\heartsuit_{\operatorname{stand}}} \simeq \Shv_G(X;\overline{\Q}_{\ell})$. In particular, this shows that the pseudocone formalism gives a language which allows for the simultaneous study for equivariant sheaves, equivariant perverse sheaves, and equivariant derived categories.

\begin{corollary}\label{Cor: Section Triangle: Pseudofunctor of standard hearts is heart of pseudofunctor}
	Let $G$ be a smooth algebraic group and let $X$ be a left $G$-variety. If $D:\SfResl_G(X)^{\op} \to \fCat$ is an $\ell$-adic derived category pseudofunctor
	\[
	D(\Gamma \times X) = D^b_c(G \backslash (\Gamma \times X);\overline{\Q}_{\ell}), \qquad D(f) = (\overline{f \times \id_X})^{\ast}
	\]
	then 
	\[
	\PC(D)^{\heartsuit} \simeq \PC(\Shv(G \backslash (-);\overline{\Q}_{\ell})) \simeq \Shv_G(X;\overline{\Q}_{\ell}).
	\] 
\end{corollary}
\begin{proof}
	The fact that for all $\Gamma \times X \in \SfResl_G(X)_0$ we have 
	\[
	D_c^b(G \backslash (\Gamma \times X); \overline{\Q}_{\ell})^{\heartsuit_{\operatorname{stand}}} \simeq \Shv(G \backslash (\Gamma \times X); \overline{\Q}_{\ell})
	\]
	is a classical fact regarding $\ell$-adic sheaves and that these equivalences vary pseudonaturally in $\SfResl_G(X)$ is routine. Regard $D$ as a truncated pseudofunctor with each truncation functor $\quot{\tau^{\geq n}}{\Gamma},$ $\quot{\tau^{\leq n}}{\Gamma}$ the standard truncation functors, i.e., equip $D_G^b(X;\overline{\Q}_{\ell})$ with the standard $t$-structure (cf.\@ the first example given in Example \ref{Example: Standard tstructure}). Applying Theorem \ref{Theorem: Heart of equivariant t-structure is equivariant heart of t-structure} then gives the equivalence
	\[
	\PC(D)^{\heartsuit} = \PC\left(D_c^b(G \backslash (-); \overline{\Q}_{\ell})^{\heartsuit_{\operatorname{stand}}}\right) \simeq \PC(\Shv(G \backslash (-); \overline{\Q}_{\ell}).
	\]
	For the remaining equivalence recall that by \cite[Lemma 6.4.8]{PramodBook} (see Remark \ref{Remark: Pseudocone Section: Using SFG instead of Resls} for where we used this) there is an equivalence of categories $\PC(\tilde{D}:\mathbf{Resl}_G(X)^{\op} \to \fCat) \simeq \PC(D)$ where $\mathbf{Resl}_G(X)$ is the category of all free $G$-resolutions of $X$ and $\tilde{D}$ is the extension of $D$ to $\mathbf{Resl}_G(X)^{\op}$. Using \cite[Exercise 6.4.3]{PramodBook} gives an equivalence $\PC(\tilde{D})^{\heartsuit} \simeq \Shv_G(X;\overline{\Q}_{\ell})$ and hence allows us to conclude
	\[
	\Shv_G(X;\overline{\Q}_{\ell}) \simeq \PC(D)^{\heartsuit} \simeq \PC(\Shv(G \backslash (-); \overline{\Q}_{\ell}).
	\]
\end{proof}
\begin{corollary}\label{Cor: Section Triangle: Psedofunctor of perv is equiv perv}
	Let $G$ be a smooth algebraic group and let $X$ be a left $G$-variety. If $D:\SfResl_G(X)^{\op} \to \fCat$ is an $\ell$-adic derived category pseudofunctor
	\[
	D(\Gamma \times X) = D^b_c(G \backslash (\Gamma \times X);\overline{\Q}_{\ell}), \qquad D(f) = (\overline{f \times \id_X})^{\ast}
	\]
	then 
	\[
	\PC(D)^{\heartsuit_{\operatorname{per}}} \simeq \PC(\Shv(G \backslash (-);\overline{\Q}_{\ell})) \simeq \Per_G(X;\overline{\Q}_{\ell}).
	\] 
\end{corollary}
\begin{proof}
	Regard $D$ as a truncated pseudofunctor with each truncation functor $\quot{\tau^{\geq n}}{\Gamma} = \quot{{}^{p}\tau^{\geq n}}{\Gamma}$ and $\quot{\tau^{\leq n}}{\Gamma} = \quot{{}^{p}\tau^{\leq n}}{\Gamma}$ the perverse truncation functors, i.e., equip $D_G^b(X;\overline{\Q}_{\ell})$ with the perverse $t$-structure (cf.\@ the second example given in Example \ref{Example: Standard tstructure}). The fact that for all $\Gamma \times X \in \SfResl_G(X)_0$ we have 
	\[
	D_c^b(G \backslash (\Gamma \times X); \overline{\Q}_{\ell})^{\heartsuit_{\operatorname{per}}} \simeq \Per(G \backslash (\Gamma \times X); \overline{\Q}_{\ell})
	\]
	is definition while verifying that these equivalences themselves vary pseudonaturally in $\SfResl_G(X)$ is routine. Applying Theorem \ref{Theorem: Heart of equivariant t-structure is equivariant heart of t-structure} then gives the equivalence
	\[
	\PC(D)^{\heartsuit_{\operatorname{per}}} = \PC\left(D_c^b(G \backslash (-); \overline{\Q}_{\ell})^{\heartsuit_{\operatorname{per}}}\right) \simeq \PC(\Per(G \backslash (-); \overline{\Q}_{\ell}).
	\]
	For the remaining equivalence recall that by \cite[Lemma 6.4.8]{PramodBook} (see Remark \ref{Remark: Pseudocone Section: Using SFG instead of Resls} for where we used this) there is an equivalence of categories $\PC(\tilde{D}:\mathbf{Resl}_G(X)^{\op} \to \fCat) \simeq \PC(D)$ where $\mathbf{Resl}_G(X)$ is the category of all free $G$-resolutions of $X$ and $\tilde{D}$ is the extension of $D$ to $\mathbf{Resl}_G(X)^{\op}$. Equip $\tilde{D}$ with a truncation structure induced by translating the peverse structure on $D$. Using \cite[Theorem 6.4.10]{PramodBook} gives an equivalence $\PC(\tilde{D})^{\heartsuit_{\operatorname{per}}} \simeq \Per_G(X;\overline{\Q}_{\ell})$ and hence allows us to conclude
	\[
	\Per_G(X;\overline{\Q}_{\ell}) \simeq D_G^b(X;\overline{\Q}_{\ell})^{\heartsuit_{\operatorname{per}}} \simeq \PC(\Per(G \backslash (-); \overline{\Q}_{\ell}).
	\]
\end{proof}
\newpage

\section{Pseudocone Categories and (Classical) Six Functor Formalisms}\label{Section: SFF for Pseudocones}

We now prove the last main purely homological algebraic result in this monograph and describe when and how (the classical) perspectives on six functor formalisms can arise on $\PC(F)$. We have already seen in Corollaries \ref{Cor: Pseudocone Functors: Monoidal Closed Derived Cat}, \ref{Cor: Pseudocone Functors: Existence of equivariant pullback  for schemes}, \ref{Cor: Pseudocone Functors: Pushforward functors for equivariant maps for scheme sheaves}, \ref{Cor: Pseudocone Functors: Exceptional Pushforward functors for equivariant maps for scheme sheaves} that for a smooth algebraic group $G$ and a $G$-equivariant morphism $h:X \to Y$ of left $G$-varieties, the categories $D_G^b(X;\overline{\Q}_{\ell})$ and their non-$\ell$-adic versions are symmetric monoidal and have adjunctions:
\[
\begin{tikzcd}
	D_G^b(X;\overline{\Q}_{\ell}) \ar[rr, bend right = 20, swap, ""{name = D}]{}{R\ul{h}_{\ast}} & & D_G^b(Y;\overline{\Q}_{\ell}) \ar[ll, bend right = 20, swap, ""{name = U}]{}{\ul{h}^{\ast}} \ar[from = U, to = D, symbol = \dashv]{}{}
\end{tikzcd}\qquad
\begin{tikzcd}
	D_G^b(Y;\overline{\Q}_{\ell}) \ar[rr, bend right = 20, swap, ""{name = D}]{}{\ul{h}^{!}} & & D_G^b(X;\overline{\Q}_{\ell}) \ar[ll, bend right = 20, swap, ""{name = U}]{}{R\ul{h}_{!}} \ar[from = U, to = D, symbol = \dashv]{}{}
\end{tikzcd}
\]
Similar results hold in the case that $G$ is a topological group which admits free $G$-spaces that are manifolds and $G$-equivariant maps $h:X \to Y$ by Corollaries \ref{Cor: Pseudocone Functors: Monoidal Closed Derived Cat Topological}, \ref{Cor: Pseudocone Functors: Existence of equivariant pullback  for spaces}, \ref{Cor: Pseudocone Functors: Topological EDC has equivariant pushpull}, and \ref{Cor: Pseudocone Functors: Topological EDC has equivariant exceptional pushpull}. Additionally the categories $D_G^b(X;\overline{\Q}_{\ell})$ are triangulated by Example \ref{Example: Standard tstructure} and we have already seen that the existence of the six functors (and the corresponding topological versions of the six functors on $D_G^b(X)$) depend primarily on properties of the pseudofunctors $D_G^b(G \backslash (-);\overline{\Q}_{\ell})$ and the various assemblies we can perform with these pseudofunctors and the myriad $2$-cells
\[
\begin{tikzcd}
	\SfResl_G(X) \ar[dr, swap]{}{\cong} \ar[rr, ""{name = U}]{}{\quo_{(-)\times X}} & & \Var_{/K} \\
	& \SfResl_G(Y) \ar[ur, swap]{}{\quo_{(-)\times Y}} \ar[from = U, to = 2-2, Rightarrow, shorten <= 4pt, shorten >= 4pt]{}{\alpha}
\end{tikzcd}
\]
induced by the morphism $h:X \to Y$. 

We now would like to study to what degree and in what sense the remaining properties  and relations of the equivariant six functor formalism depends only on the pseudofunctors, their local properties, and their categories of pseudocones. Additionally, we also have not seen how these equivariant six functor formalisms relate to the forgetful functors mapping to the standard six functor formalisms. While relations with the forgetful functors and the open/closed morphism distinguished triangles will be deferred to Section \ref{Section: SFF for Pseudocones}, as we need to have developed more of the actual geometry involved in equivariant geometry and topology, in this section we will examine the more algebraic and categorical relations we can study. In particular, we will show that the pseudocone analogues of the identities
\begin{align*}
	(Rh_{!}A) \otimes_{Y} B &\cong Rh_{!}\left(A \otimes_{X} h^{\ast}(B)\right)  \\
	\left[Rh_{!}(A), B\right]_Y &\cong Rh_{\ast}\left(\left[A, h^{!}(B)\right]_{X}\right) \\
	h^{!}\left(\left[B, C\right]_{Y}\right) &\cong \left[h^{\ast}(B), h^{!}(C)\right]_{X}
\end{align*}
which hold and are natural for all $A \in D_c^b(X;\overline{\Q}_{\ell})_0$ and for all $B, C \in D_c^b(Y;\overline{\Q}_{\ell})_0$.

As we proceed in this section, assume that we have a $2$-cell
\[
\begin{tikzcd}
	\Cscr \ar[dr, swap]{}{\varphi} \ar[rr, ""{name = U}]{}{\rho} & & \Escr \\
	& \Dscr \ar[ur, swap]{}{\psi} \ar[from = U, to = 2-2, Rightarrow, shorten <= 4pt, shorten >= 4pt]{}{\alpha}
\end{tikzcd}
\]
in $\fCat$ and a symmetric monoidal closed pseudofunctor $F:\Escr^{\op} \to \fCat$ so that we have the induced $2$-cells
\[
\begin{tikzcd}
	\Cscr^{\op} \ar[dr, swap]{}{\varphi^{\op}} \ar[rr, ""{name = U}]{}{\rho^{\op}} & & \Escr^{\op} \\
	& \Dscr^{\op} \ar[ur, swap]{}{\psi^{\op}} \ar[to = U, from = 2-2, Rightarrow, shorten <= 4pt, shorten >= 4pt]{}{\alpha^{\op}}
\end{tikzcd}
\]
and:
\[
\begin{tikzcd}
	\Cscr^{\op} \ar[rrr, bend left = 30, ""{name = U}]{}{F \circ \psi^{\op} \circ \varphi^{\op}} \ar[rrr, bend right = 30, swap, ""{name = D}]{}{F \circ \rho^{\op}} & & & \fCat \ar[from = U, to =D, Rightarrow, shorten <= 4pt, shorten >= 4pt, swap]{}{\iota_F \ast \alpha^{\op}}
\end{tikzcd}
\]
Write the functor induced by $\iota_F\ast \alpha^{\op}$ through Theorem \ref{Thm: Pseudocone Functors: Pullback induced by fibre functors in pseudofucntor} as $(F \ast \alpha)^{\square}$. Assume also that we have three additional pseudonatural transformations
\[
\begin{tikzcd}
	\Cscr^{\op} \ar[rrr, bend left = 30, ""{name = U}]{}{F \circ \rho^{\op}} \ar[rrr, bend right = 30, swap, ""{name = D}]{}{F \circ \psi^{\op} \circ \varphi^{\op}} & & & \fCat \ar[from = U, to = D, Rightarrow, shorten <= 4pt, shorten >= 4pt]{}{(F \ast \alpha)_{\square}}
\end{tikzcd}\qquad 
\begin{tikzcd}
	\Cscr^{\op} \ar[rrr, bend left = 30, ""{name = U}]{}{F \circ \rho^{\op}} \ar[rrr, bend right = 30, swap, ""{name = D}]{}{F \circ \psi^{\op} \circ \varphi^{\op}} & & & \fCat \ar[from = U, to = D, Rightarrow, shorten <= 4pt, shorten >= 4pt]{}{(F \ast \alpha)_{\dagger}}
\end{tikzcd}
\]
and:
\[
\begin{tikzcd}
	\Cscr^{\op} \ar[rrr, bend left = 30, ""{name = U}]{}{F \circ \psi^{\op} \circ \varphi^{\op}} \ar[rrr, bend right = 30, swap, ""{name = D}]{}{F \circ \rho^{\op}} & & & \fCat \ar[from = U, to = D, Rightarrow, shorten <= 4pt, shorten >= 4pt]{}{(F \ast \alpha)_{\square}}
\end{tikzcd}
\]
Assume also that the fibre functors all satisfy the adjunctions
\[
(F \ast \alpha)^{\square} \dashv (F \ast \alpha)_{\square}, \qquad (F \ast \alpha)_{\dagger} \dashv (F \ast \alpha)^{\dagger}.
\]
Note that we are writing $ (F \ast \alpha)^{\square}, (F \ast \alpha)_{\square}, (F \ast \alpha)_{\dagger},$ and $(F \ast \alpha)^{\dagger}$ suggestively as if they are pullback/pushforward functors and exceptional pullback/pushforward functors, but the notation is otherwise simply meant to represent how we want to argue with these functors; the notations otherwise carry no extra meaning. From Theorem \ref{Theorem: Section 2: Monoidal preequivariant pseudofunctor gives monoidal equivariant cat}, Corollary \ref{Cor: Pseudocone Functors: Braided monoidal pseudonats are braided monoidal functors}, Theorem \ref{Thm: Functor Section: Psuedonatural trans are pseudocone functors}, and Theorem \ref{Thm: Section 3: Gamma-wise adjoints lift to equivariant adjoints} we know that these give rise to adjoints
\[
\begin{tikzcd}
	\PC(F \circ \rho^{\op}) \ar[rr, bend right = 20, swap, ""{name = D}]{}{[A,-]} & & \PC(F \circ \rho^{\op}) \ar[ll, bend right = 20, swap, ""{name = U}]{}{(-) \otimes A} \ar[from = U, to = D, symbol = \dashv]
\end{tikzcd}
\]
\[
\begin{tikzcd}
	\PC(F \circ \psi^{\op} \circ \varphi^{\op}) \ar[rr, bend right = 20, swap, ""{name = D}]{}{[B,-]} & & \PC(F \circ \psi^{\op} \circ \varphi^{\op}) \ar[ll, bend right = 20, swap, ""{name = U}]{}{(-) \otimes B} \ar[from = U, to = D, symbol = \dashv]
\end{tikzcd}
\]
\[
\begin{tikzcd}
	\PC(F \circ \rho^{\op}) \ar[rr, bend right = 20, swap, ""{name = D}]{}{(F \ast \alpha)_{\square}} & & \PC(F \circ \psi^{\op} \circ \varphi^{\op}) \ar[ll, bend right = 20, swap, ""{name = U}]{}{(F \ast \alpha)^{\square}} \ar[from = U, to = D, symbol = \dashv]
\end{tikzcd}
\]
and
\[
\begin{tikzcd}
	\PC(F \circ \psi^{\op} \circ \varphi^{\op}) \ar[rr, bend right = 20, swap, ""{name = D}]{}{(F \ast \alpha)^{\dagger}} & & \PC(F \circ \rho^{\op}) \ar[ll, bend right = 20, swap, ""{name = U}]{}{(F \ast \alpha)_{\dagger}} \ar[from = U, to = D, symbol = \dashv]
\end{tikzcd}
\]
where $A \in \PC(F \circ \rho^{\op})_0$ and $B \in \PC(F \circ \psi^{\op} \circ \varphi^{\op})_0$ are arbitrary. When we have that $F$ is a truncated or triangulated pseudofunctor, we now have the basic set-up for a six functor formalism between $\PC(F \circ \rho^{\op})$ and $\PC(F \circ \psi^{\op} \circ \varphi^{\op})$. We now show that if these functors all $X$-locally in $\Cscr$ satisfy the natural isomorphisms required of our six functors, then it will be the case that all the functors above satisfy them as well.

\begin{Theorem}\label{Thm: Section Triangle: The Six Functor Relations}
	Assume that we have $\Cscr,$ $\Dscr,$ $\Escr,$ $F$, $\alpha$, $\varphi,$ $\psi,$ $\rho,$ $(F \ast \alpha)_{\square},$ $(F \ast \alpha)^{\square},$ $(F \ast \alpha)_{\dagger},$ and $(F \ast \alpha)^{\dagger}$ be given as in the setup prior to the theorem and assume that for all $X \in \Cscr_0$ and for all objects $A \in F(\rho X)_0$ and for all objects $B, C \in F((\psi \circ \varphi)X)_0$ we have natural isomorphisms
	\begin{align*}
		\quot{(F \ast \alpha)_{\dagger}}{X}\left(\quot{A}{X}\right) \us{F((\psi \circ \varphi)X)}{\otimes} \quot{B}{X} &\cong \quot{(F \ast \alpha)_{\dagger}}{X}\left(\quot{A}{X} \us{F(\rho X)}{\otimes} (F \ast \alpha)^{\square}(\quot{B}{X})\right)  \\
		\left[\quot{(F \ast \alpha)_{\dagger}}{X}(\quot{A}{X}), \quot{B}{X}\right]_{F((\psi \circ \varphi)(X))} &\cong \quot{(F \ast \alpha)_{\square}}{X}\left(\left[\quot{A}{X}, \quot{(F \ast \alpha)^{\dagger}}{X}(\quot{B}{X})\right]_{F(\rho X)}\right) \\
		\quot{(F \ast \alpha)^{\dagger}}{X}\left(\left[\quot{C}{X}, \quot{B}{X}\right]_{F((\psi \circ \varphi)X)}\right) &\cong \left[\quot{(F \ast \alpha)^{\square}}{X}(\quot{C}{X}), \quot{(F \ast \alpha)^{\dagger}}{X}(\quot{B}{X})\right]_{F(\rho X)}
	\end{align*}
	Then for any objects $A \in \PC(F \circ \rho^{\op})_0$ and $B, C \in \PC(F \circ \psi^{\op} \circ \varphi^{\op})_0$ there are natural isomorphisms:
	\begin{align*}
		(F \ast \alpha)_{\dagger}(A) \us{\PC(F \circ \psi^{\op} \circ \varphi^{\op})}{\otimes} B &\cong (F \ast \alpha)_{\dagger}\left(A \us{\PC(F \ast \rho^{\op})}{\otimes} (F \ast \alpha)^{\square}(B)\right) \\
		\left[(F \ast \alpha)_{\dagger}(A), B\right]_{\PC(F \circ \psi^{\op} \circ \varphi^{\op})} &\cong (F \ast \alpha)_{\square}\left(\left[A, (F \ast \alpha)^{\dagger}(B)\right]_{\PC(F \circ \rho^{\op})}\right) \\
		(F \ast \alpha)^{\dagger}\left(\left[C, B\right]_{\PC(F \circ \psi^{\op} \circ \varphi^{\op})}\right) &\cong \left[(F \ast \alpha)^{\square}(C), (F \ast \alpha)^{\dagger}(B)\right]_{\PC(F \circ \rho^{\op})}
	\end{align*}
\end{Theorem}
\begin{proof}
	Because each of the natural isomorphisms above hold at the $X$-local level for all $X \in \Cscr^{\op}$, we only need to check that they may be assembled into morphisms at the pseudocone category level to deduce that they exist. However, this is a routine but extremely tedious check using both the pseudonaturality of the various pseudonatural transformations coupled with the naturality of the isomorphisms applied to the transition morphisms $\tau_f^{A}:F(\rho^{\op}(f))(\quot{A}{Y}) \to \quot{A}{X}$, and similarly for $B$ and $C$.
\end{proof}
It now follows from the theorem above that for a smooth algebraic group $G$ and a $G$-equivariant morphism $h:X \to Y$ of left $G$-varieties, the categories $D_G^b(X;\overline{\Q}_{\ell})$ and $D_G^b(X)$ both admit six functor formalisms with the relations above with $\square = \ast$ and $\dagger = !$.
\begin{corollary}\label{Cor: Section Triangle: SFF Relations for EDC Scheme}
	Let $G$ be a smooth algebraic group and let $h:X \to Y$ be a $G$-equivariant morphism of left $G$-varieties. Consider the $2$-cell
	\[
	\begin{tikzcd}
		\SfResl_G(X) \ar[dr, swap]{}{\cong} \ar[rr, ""{name = U}]{}{\quo_{(-)\times X}} & & \Var_{/K} \ar[rr]{}{D^b_c(-; \overline{\Q}_{\ell})} & & \fCat \\
		& \SfResl_G(Y) \ar[ur, swap]{}{\quo_{(-)\times Y}} \ar[from = U, to = 2-2, Rightarrow, shorten <= 4pt, shorten >= 4pt]{}{\overline{\id_{-} \times h}}
	\end{tikzcd}
	\]
	and define the pseudonatural transformations/functors by assigning $(F \ast \alpha)^{\square} = h^{\ast}$, $(F \ast \alpha)_{\square} = Rh_{\ast},$ $(F \ast \alpha)^{\dagger} = h^{!}$, and $(F \ast \alpha)_{\dagger} = Rh_{!}$. Then there are natural isomorphisms
	\begin{align*}
		Rh_{!}(A) \us{Y}{\otimes} B &\cong Rh_{!}\left(A \us{X}{\otimes} h^{\ast}(B)\right), \\
		\left[Rh_{!}(A), B\right]_{Y} &\cong Rh_{\ast}\left(\left[A, h^{!}(B)\right]_{X}\right), \\
		h^{!}\left([C,B]_{Y}\right) &\cong \left[h^{\ast}(C), h^{!}(B)\right]_{X}
	\end{align*}
	where $A \in D_G^b(X; \overline{\Q}_{\ell})_0$ and $B, C \in D_G^b(Y;\overline{\Q}_{\ell})_0$. 
\end{corollary}
\begin{proof}
	By Corollaries \ref{Cor: Pseudocone Functors: Monoidal Closed Derived Cat}, \ref{Cor: Pseudocone Functors: Existence of equivariant pullback  for schemes}, \ref{Cor: Pseudocone Functors: Pushforward functors for equivariant maps for scheme sheaves}, \ref{Cor: Pseudocone Functors: Exceptional Pushforward functors for equivariant maps for scheme sheaves} we know that we satisfy the technical assumptions setting up Theorem \ref{Thm: Section Triangle: The Six Functor Relations}. Consequently, since each $\Gamma$-local fibre category and corresponding functors satisfy the natural isomorphisms
	\begin{align*}
		R(\overline{\id_{\Gamma} \times h})_{!}(\AGamma) \us{\otimes}{G \backslash (\Gamma \times Y)} \BGamma &\cong R(\overline{\id_{\Gamma} \times h})_{!}\left(\AGamma \us{\otimes}{G \backslash (\Gamma \times X)} (\overline{\id_{\Gamma} \times h})^{\ast}(\BGamma)\right), \\
		\left[R(\overline{\id_{\Gamma} \times h})_{!}(\AGamma), \BGamma\right]_{G \backslash (\Gamma \times Y)} &\cong R(\overline{\id_{\Gamma} \times h})_{\ast}\left(\left[\AGamma, (\overline{\id_{\Gamma} \times h})^{!}(\BGamma)\right]_{G \backslash (\Gamma \times X)}\right), \\
		(\overline{\id_{\Gamma} \times h})^{!}\left([\quot{C}{\Gamma},\BGamma]_{G \backslash (\Gamma \times Y)}\right) &\cong \left[(\overline{\id_{\Gamma} \times h})^{\ast}(\quot{C}{\Gamma}), (\overline{\id_{\Gamma} \times h})^{!}(\BGamma)\right]_{G \backslash (\Gamma \times X)}
	\end{align*}
	for all $\AGamma \in A, \BGamma \in B$, and $\quot{C}{\Gamma} \in C$ we simply apply Theorem \ref{Thm: Section Triangle: The Six Functor Relations} to give our result.
\end{proof}
\newpage

\chapter{Equivariant Descent and Categories of Pseudocones}\label{Section: Equivariant Descent}
We now begin to study equivariant algebraic geometry in earnest. While the first half of our work has shown that we can develop and study the theory of equivariant derived categories in topological and geometric contexts almost entirely from the perspective of pseudocone categories. In particular, the theory we have developed shows that there is a strong formal analogy that in a working sense equivariant descent comes in two moving components which meet transversally:
\begin{enumerate}
	\item The first component is the geometric/topological component. This is the part of the story that occupies the study of the category of $G$-resolutions, simplicial arguments and objects, which $G$-resolutions we should use or consider and why, properties of resolutions, etc.
	\item The second component is the algebraic/categorical component. This is the part of the story that focus on the various categories, pseudofunctors, and other gadgets we use and probe our objects with, as well as what algebraic properties these gadgets satisfy (or fail to satisfy) and what they tell us about the geometry/topology of the action $G \times X \to X$.
\end{enumerate}
Of course, these two components are transversal in nature and need not meet orthogonally: they need not move independently from each other and in fact can and do inform what techniques or settings we should use from either component in any given moment but also do not move parallel to each other. 

In this chapter we will begin to focus on the geometric and topological side of the story and start our in-depth study of equivariant descent. In particular, we will provide the framework that shows the category $\SfResl_G(X)$ has the properties we have claimed and used in Section \ref{Subsection: Equivariant Derived Cat of Var} (such as the existence of the quotient functor $\quo:\SfResl_G(X) \to \Var_{/K}$). Afterwards we will define what it means to be an equivariant category, what equivariant functors are, and introduce the ``equivariance'' isomorphism and cocycle condition attached to an equivariant category.

\section{Equivariant Categories and Smooth Free Varieties}\label{Section: Smooth Free Geometry}

In Section \ref{Subsection: Equivariant Derived Cat of Var} we saw some of the properties of the category $\SfResl_G(X)$ that allowed us to define $D_G^b(X; \overline{\Q}_{\ell})$. The most important ingredients in doing this were the fact that $\Gamma \times X$ admits a geometric quotient for any $G$-variety $X$ and for any locally isotrivial smooth free variety $\Gamma$ and these quotients extend to a functor $\quo:\SfResl_G(X) \to \Var_{/K}$; because $\Var_{/K}$ rarely admits coequalizers, this requires careful justification and set-up. In this and the next section (cf.\@ Section \ref{Subsection: The Equivariance Isos}) we will show how these facts arise as well as study some extra properties of the category $\Sf(G)$ as we vary groups and how these categories interact with the equivariant categories $F_G(X) = \PC(F)$ for a class of pseudofunctors which we call pre-equivariant.

Perhaps the best place to start this story is at the beginning; because of how important they are to what we do, we record various important properties of smooth morphisms and quotients of smooth morphisms both for the sake of completeness but also for the sake of non-algebraic geometers.. Recall that the definitions of smooth free $G$-varieties, the categories $\Sf(G)$ and $\SfResl_G(X)$, and other basics can be found in Section \ref{Subsection: Equivariant Derived Cat of Var} (cf.\@ Definitions \ref{Defn: Pseudocone Section: Smooth free varieties}, \ref{Defn: Pseudocone Section: Sf(G)}, \ref{Defn: Pseudocone Section: SfResl}).
\begin{proposition}\label{Prop: Smooth quotient from principal G var}
	If $\Gamma$ is a $G$-variety in $\Sf(G)_0$ then the quotient morphism $\quo_{\Gamma}:\Gamma \to G\backslash \Gamma$ is smooth and the diagram
	\[
	\xymatrix{
		G \times \Gamma \ar[r]^{\pi_2} \ar[d]_{\alpha_{\Gamma}} \pullbackcorner & \Gamma \ar[d]^{\quo_{\Gamma}} \\
		\Gamma \ar[r]_-{\quo_{\Gamma}} & G \backslash \Gamma
	}
	\]
	is a pullback square.
\end{proposition}
\begin{proof}
	This is \cite[Proposition 0.9]{GIT}.
\end{proof}
\begin{lemma}[{\cite[\href{https://stacks.math.columbia.edu/tag/02K5}{Tag 02K5}]{stacks-project}}]\label{Lemma: Section 2.1: Smooth morphism lemma}
	Consider the commutative diagram
	\[
	\xymatrix{
		X \ar[rr]^-{f} \ar[dr]_{h} & & Y \ar[dl]^{g} \\
		& Z
	}
	\]
	of schemes for which $f$ is surjective and smooth and $h$ is smooth. Then $g$ is smooth as well.
\end{lemma}
\begin{proposition}\label{Prop: Section 2.1: Smooth quotient map}
	Let $Y, Z$ be $G$-schemes with smooth quotient morphisms $\quo_Y:Y \to G \backslash Y$ and $\quo_Z:Z \to G \backslash Z$. Then if there is a smooth $G$-equivariant morphism $f:Y \to Z$ there exists a unique morphism $\of:G \backslash Y \to G \backslash Z$ making the diagram
	\[
	\xymatrix{
		Y \ar[r]^-{f} \ar[d]_{\quo_Y} & Z \ar[d]^{\quo_Z} \\
		G\backslash Y \ar[r]_-{\of} & G\backslash Z
	}
	\]
	commute. Moreover, $\of$ is smooth as well.
\end{proposition}
\begin{proof}
	The existence of the morphism $\of$ follows from the universal property of the quotient (as a coequalizer of $\alpha_Y, \pi_2:G \times Y \to Y$) applied to the morphism $(\quo_Z \circ f):Y \to G \backslash Z$. That $\of$ is smooth follows from Lemma \ref{Lemma: Section 2.1: Smooth morphism lemma} together with the observations that $\quo_{Y}$ is smooth and surjective while the composite $\quo_Z \circ f$ is smooth.
\end{proof}

The next proposition is a crucial technical result which allows us to deduce that for any $\Gamma \in \Sf(G)_0$ and any left $G$-variety $X$, the product $\Gamma \times X$ is free and admits a smooth geometric quotient. This allows us to resolve the action of $G$ on $X$ by looking at $G$-resolutions of the form $\pi_2^{\Gamma}:\Gamma \times X \to X$ as the $\Gamma$ vary through $\Sf(G)$. As we have seen, this is the fundamental technical tool which allows us to describe equivariant information by taking descent data indexed by a pseudofunctor defined on these resolutions.
\begin{proposition}\label{Prop: Section 2.1: The quotient prop}
	Let $\Gamma$ be a smooth free $G$-variety and let $Y$ be a $G$-variety. Then the variety $\Gamma \times Y$ is free with smooth quotient morphism 
	\[
	\quo_{\Gamma \times Y}:\Gamma \times Y \to G \backslash (\Gamma \times Y).
	\]
	Furthermore, if there is a $G$-equivariant morphism of smooth free $G$-varieties $f:\Gamma \to \Gamma^{\prime}$ then there is a smooth morphism making the diagram
	\[
	\xymatrix{
		\Gamma \times Y \ar[r]^-{f \times \id_Y} \ar[d]_{q_\Gamma} & \Gamma^{\prime} \times Y \ar[d]^{q_{\Gamma^{\prime}}} \\
		G\backslash(\Gamma \times Y) \ar@{-->}[r]_-{\exists!\overline{f}} & G \backslash(\Gamma^{\prime} \times Y)
	}
	\]
	commute in $\Var_{/\Spec K}$.
\end{proposition}
\begin{proof}
	Our argument for the existence of the quotient $G \backslash (\Gamma \times Y)$ follows the sketch given in \cite[Proposition 6.1.13]{PramodBook}. Because each quotient $\Gamma \to G \backslash \Gamma$ is locally isotrivial, by \cite[Proposition 4]{SerreAlgebraicFibreSpaces} it follows that if we can show for all finite subsets $F \subseteq \lvert Y \rvert$ there is an open affine $U$ of $Y$ for which $F \subseteq \lvert U \rvert$, then it follows that $\Gamma \times Y$ is free and admits a geometric quotient $\quo:\Gamma \times Y \to G \backslash (\Gamma \times Y)$. However, by \cite[Proposition 3.3.36]{LiuAGAC} because $Y$ is quasi-projective\footnote{Recall from Remark \ref{Remark: Notations and conventions} that we only consider quasi-projective varieties in this monograph.} over $\Spec K$, for any finite set $F \subseteq \lvert Y \rvert$ we can find an open affine subscheme $U$ of $Y$ for which $F \subseteq \lvert U \rvert$ so it follows that $\Gamma \times Y$ is a free $G$-variety. From this it follows that the morphism
	\[
	\quo_{\Gamma \times Y}:\Gamma \times Y \to G \backslash (\Gamma \times Y)
	\]
	is smooth, which establishes the first claim. For the second claim we simply apply Proposition \ref{Prop: Section 2.1: Smooth quotient map} to the $G$-equivariant morphism $f \times \id_Y$, which is smooth because $\id_Y$ is {\'e}tale, $f$ is smooth by assumption, and because the product of smooth morphisms is smooth.
\end{proof}
\begin{remark}\label{Remark: Section Descent: Affine SfG Vars}
By \cite[Proposition 6.1.13]{PramodBook} (see also \cite[Lemme XIV 1.4]{Raynaud}) it follows that when $G$ is a smooth affine algebraic group then asking for the quotient map $\Gamma \to G \backslash \Gamma$ to be locally isotrivial is equivalent to asking for the quotient map to be {\'e}tale locally trivializable.  In particular, when $G$ is affine smooth then the $\Sf(G)$-varieties are exactly the pure dimensional smooth free $G$-varieties which admit a geometric quotient.
\end{remark}
\begin{remark}
We suggest defining our smooth free $G$-varieties to be the (pure dimensional) locally isotrivial principal $G$-fibrations $\Gamma$ in $\Var_{/K}$ precisely because this is the class of $G$-varieties which allow us to have a quotient functor $\quo:\SfResl_G(X) \to \Var_{/K}$.  Because we have the induced equivalences of categories $D_G^b(\Gamma) \simeq D_c^b(G \backslash \Gamma)$ by virtue of Corollary \ref{Cor: Pseudocone Section: SfG variety quotient is edc} we know that this is a reasonable categorical language for performing equivariant descent and so suggests giving a comparison between $D_G^b(X)$ and $D_c^b([G \backslash X])$.

In the case where $G$ is an affine algebraic group we have equivalences of categories $D_G^b(X) \simeq D_c^b([G \backslash X])$ and $D_G^b(X;\overline{\Q}_{\ell}) \simeq D_c^b([G\backslash X];\overline{\Q}_{\ell})$; these are either definitional (as in \cite{PramodBook}) or the result of theorems (cf.\@ \cite{MyThesis}) depending on the perspective taken on the derived category of the quotient stack $[G \backslash X]$. When $G$ is nonaffine we still have the equivalences
\[
D_G^b(\Gamma) \simeq D_c^b(G \backslash \Gamma) \simeq D_c^b([G \backslash \Gamma])
\]
for all $\Sf(G)$-varieties $\Gamma$ (and similarly for the $\ell$-adic versions) precisely because in this case $[G \backslash \Gamma]$ is represented by the variety $G \backslash \Gamma$.  This suggests that by descending through the varieties $\Gamma \times X$ we get equivalences of categories
\[
D_G^b(\Gamma \times X) \simeq D_c^b\left(G \backslash (\Gamma \times X)\right) \simeq D_c^b\left(\left[G \backslash (\Gamma \times X)\right]\right)
\]
again because the stack $\left[G \backslash (\Gamma \times X)\right]$ is represented by the variety $G \backslash (\Gamma \times X)$. By using a simplicial presentation of $D_c^b([G \backslash X])$ as the category of groupoid sheaves $D_{\text{eq}}^b(\ul{G \backslash X}_{\bullet})$ (where the simplicial scheme $\ul{G \backslash X}_{\bullet}$ is the simplicial presentation of the translation groupoid of the group action $\alpha_X, \pi_2:G \times X \to X$ and has $n$-simplexes given by $\ul{G \backslash X}_n = G^n \times X$) we should be able to proceed with the techniques of \cite[Chapter 8]{MyThesis} in order to ultimately derive the equivalence $D_G^b(X) \simeq D_c^b([G \backslash X])$. This gives rise to the conjecture below which can be seen as justifying a ``model-agnostic'' perspective on the equivariant derived category.
\end{remark}
\begin{conjecture}\label{Conjecture: EDC of nonaffine group}
Let $G$ be a smooth algebraic group and let $X$ be a left $G$-variety. Then there are equivalences of categories
\[
\PC(D_c^b(G \backslash(-))) \simeq D_c^b([G \backslash X])
\]
and
\[
\PC(D_c^b(G \backslash (-);\overline{\Q}_{\ell})) \simeq D_c^b([G \backslash X]; \overline{\Q}_{\ell}).
\]
\end{conjecture}

Before proceeding to discuss the theory of equivariant categories, we record a crucial, but straightforward, scheme-theoretic lemma on algebraic groups. In particular we need to know that for any smooth algebraic group $G$, the group $G$ is itself an object in $\Sf(G)$. While this is obvious enough, it is fundamental to the constructions we will do later on (especially involving the change of group functors; cf.\@ Theorem \ref{Thm: Section 3: Change of groups functor}) that we make sure to record it.
\begin{lemma}\label{Lemma: Section 2: Smooth alg group G is in Sf(G)}
	Let $G$ be a smooth algebraic group over a field $K$. Then $G \in \Sf(G)_0$.
\end{lemma}
\begin{proof}
	The fact that $G$ admits a geometric quotient given by $G \backslash G \cong \Spec K$ is routine  by \cite[Theorem V.3.2]{SGA3} (cf.\@ also \cite{milneiAG}). That this is {\'e}tale locally trivializable and in particular locally isotrivial follows from the fact that the diagram
	\[
	\xymatrix{
		G \ar@{=}[r] \pullbackcorner \ar[d] & G \ar[d] \\
		\Spec K \ar@{=}[r] & \Spec K
	}
	\]
	is a pullback and $\id_{\Spec K}$ is a finite {\'e}tale cover of $\Spec K$. Finally, since $G$ is pure dimensional by \cite[Lemma VIB.1.5]{SGA3}, the lemma follows.
\end{proof}

We now define pre-equivariant pseudofunctors. These are the necessary ingredients we need in order to define both equivariant categories on varieties and equivariant categories on spaces. These pseudofunctors are particularly well-suited for equivariant descent and both capture and generalize the equivariant derived categories $D_G^b(X)$ and $D_G^b(X;\overline{\Q}_{\ell})$.

\begin{definition}\label{Defn: Equivariant Descent: Preequivariant Pseudofunctor}
	Let $G$ be a smooth algebraic group and let $X$ be a left $G$-variety. We define a pre-equivariant pseudofunctor\index[terminology]{Pre-Equivariant Pseudofunctor! On Varieties} on $X$ to be a pair of pseudofunctors $(F,\overline{F})$ where $F:\SfResl_G(X)^{\op} \to \fCat$ and $\overline{F}:\Var_{/K}^{\op} \to \fCat$ are pseudofunctors such that the diagram of (bi)categories
	\[
	\begin{tikzcd}
		\SfResl_G(X)^{\op} \ar[dr, swap]{}{F} \ar[r]{}{\quo^{\op}} & \Var_{/K}^{\op} \ar[d]{}{\overline{F}} \\
		& \fCat
	\end{tikzcd}
	\]
	commutes.
\end{definition}
\begin{remark}
	The definition given here for pre-equivariant pseudofunctors is slightly different than the one given in \cite{MyThesis} where they first appeared; it does, however, match with that given in \cite{DoretteMe}.
\end{remark}

We also now extend the definition of a pre-equivariant pseudofunctor to topological spaces.
\begin{definition}\label{Defn: Section Geometry: Preeq functor for spaces}
	Let $G$ be a topological group and let $X$ be a left $G$-space. A pre-equivariant pseudofunctor\index[terminology]{Pre-Equivariant Pseudofunctor! On Spaces} (for the space $X$) is a pair $(F, \overline{F})$ where $F:\mathbf{Resl}_G(X)^{\op} \to \fCat$ and $\overline{F}:\Top^{\op} \to \fCat$ are pseudofunctors for which there is a commuting diagram
	\[
	\begin{tikzcd}
		\mathbf{Resl}_G(X)^{\op} \ar[rr]{}{(P \mapsto G \backslash P)^{\op}} \ar[drr, swap]{}{F} & & \Top^{\op} \ar[d]{}{\overline{F}} \\
		& & \fCat
	\end{tikzcd}
	\]
	of (bi)categories.
\end{definition}

	To a pre-equivariant pseudofunctor $(F,\overline{F})$ there is an evident pseudocone category we can and should consider: the category
	\[
	\PC(F) = \PC(\overline{F} \circ \quo^{\op}).
	\]
	This is by far and away the most important class of categories to consider for our applications. When studying the categories $\PC(F)$, we will see in the remainder of this monograph that it is frequently important to leverage the properties of $\overline{F}$ in order to derive properties about $\PC(F)$ indirectly. For instance, it is relatively straightforward to show that if $\overline{F}$ is a triangulated pseudofunctor then so is $F$.

\begin{example}
	Here is a large list of pre-equivariant pseudofunctors for both schemes and spaces. Note that many of these examples overlap with Example \ref{Example: Pseudocone Section: Big list of examples for great justice}.
	\begin{enumerate}
		\item Let $G$ be a smooth algebraic group and let $X$ be a left $G$-variety. Then there is a pre-equivariant pseudofunctor pair $(F,\overline{F})$ for the category category of ($\ell$-adic) equivariant perverse sheaves. This pre-equivariant pseudofunctor is given by letting $F:\SfResl_G(X)^{\op} \to \fCat$ be a pseudofunctor with object and morphism assignments
		\[
		\Gamma \times X \mapsto \Per(G \backslash (\Gamma \times X); \overline{\Q}_{\ell})
		\]
		and
		\[
		f \times \id_X \mapsto {}^{p}\of^{\ast}:\Per(G \backslash (\Gamma^{\prime} \times X); \overline{\Q}_{\ell}) \to \Per(G \backslash (\Gamma \times X); \overline{\Q}_{\ell}),
		\]
		respectively. It's companion pseudofunctor $\overline{F}:\Var_{/K}^{\op} \to \fCat$ is then a pseudofunctor with $\overline{F}(Z) := \Per(Z;\overline{\Q}_{\ell})$ and with $\overline{F}(h) := {}^{p}h^{\ast}$ for morphisms of varieties $h$.
		\item Let $G$ be a smooth algebraic group over $K$ and let $X$ be a left $G$-variety. Then there is a local system pre-equivariant pseudofunctor on $X$ given as follows. The pseudofunctor $L:\SfResl_G(X)^{\op} \to \fCat$ is defined with object and morphism assignments
		\[
		\Gamma \times X \mapsto \Loc(G \backslash (\Gamma \times X); \overline{\Q}_{\ell})
		\]
		and
		\[
		f \times \id_X \mapsto \of^{\ast}:\Loc(G \backslash (\Gamma^{\prime} \times X)) \to \Loc(G \backslash (\Gamma \times X)),
		\]
		respectively. The companion pseudofunctor $\overline{L}:\Var_{/K}^{\op} \to \fCat$ is given by $\overline{L}(Z) = \Loc(Z;\overline{\Q}_{\ell})$ on objects and $\overline{L}(f) = f^{\ast}$ for morphisms.
		\item Let $G$ be a smooth algebraic group over $K$ and let $X$ be a left $G$-variety. Then the pre-equivariant pseudofunctor of the derived category of perverse sheaves on $X$, is determined by a pseudofunctor $F:\SfResl_G(X)^{\op} \to \fCat$ where $F$ has object and morphism assignments
		\[
		\Gamma \times X \mapsto D^b(\Per(G \backslash (\Gamma \times X)))
		\]
		and
		\[
		f \times \id_X \mapsto {}^{p}L\of^{\ast}:D^b(\Per(G \backslash (\Gamma^{\prime} \times X))) \to D^b(\Per(G \backslash (\Gamma \times X))),
		\]
		respectively. The companion pseudofunctor is a pseudofunctor $\overline{F}:\Var_{/K}^{\op} \to \fCat$ with $\overline{F}(Z) = D^b(\Per(Z))$ and $\overline{F}(h) = {}^{p}Lh^{\ast}$ on morphisms. If $G$ is a complex Lie group and $X$ is a complex analytic variety then there is a corresponding topological version as well defined mutatis mutandis.
		\item Let $G$ be a topological group and let $X$ be a left $G$-space. The pre-equivariant pseudofuncotr of perverse sheaves on $X$ is given as follows. Define the pseudofunctor $F:\FResl_G(X)^{\op} \to \fCat$ where $F$ has object and morphism assignments
		\[
		P \mapsto \Per(G \backslash P)
		\]
		and
		\[
		f \mapsto {}^{p}\of^{\ast}:\Per(G \backslash Q) \to \Per(G \backslash P)
		\]
		respectively. The companion $\overline{F}:\mathbf{TopStrat}^{\op} \to \fCat$ is given by sending a stratified space $Z$ to its category of perverse sheaves $\Per(Z)$ and sending a morphism $h$ of stratified spaces to the corresponding category of perverse sheaves.
		\item Let $G$ be a smooth algebraic group and let $X$ be a left $G$-variety. We can then describe the pre-equivariant pseudofunctor of $G$-equivariant $D$-modules as follows. Let $D:\SfResl_G(X)^{\op} \to \fCat$ is a pseudofunctor with object and morphism assignments
		\[
		D(\Gamma \times X) := \DMod(G \backslash (\Gamma \times X))
		\]
		and
		\[
		D(f \times \id_X) := \of^{\ast}:\DMod(G \backslash (\Gamma^{\prime} \times X)) \to \DMod(G \backslash (\Gamma \times X))
		\]
		respectively. Note that in this case the functor $\of^{\ast}$ exists by \cite[Section 4.1]{benzvi2015character} because each morphism $\of$ is smooth by \@ Proposition \ref{Prop: Section 2.1: The quotient prop}. We similarly define $\overline{D}:\Var_{K,\text{sm}}^{\op}$ to send a $K$-variety $Z$ to its category $\DMod(Z)$ of $\Dcal$-modules and send a morphism to the corresponding $\Dcal$-module pullback (which also exists by \cite[Section 4.1]{benzvi2015character} and the smoothness of morphisms in $\Var_{/K,\text{sm}}$ given by assumption).
		\item Let $G$ be a Lie group, regarded as a topological group, and let $X$ be a smooth manifold (regarded as a topological space) with a smooth $G$-action. Then the pre-equivariant pseudofunctor modeling $G$-equivariant sheaves on $X$ is given as follows. Define $F:\FResl_G(X)^{\op} \to \fCat$ as a pseudofunctor with object and morphism assignments
		\[
		P \mapsto \Shv(G \backslash P)
		\]
		and
		\[
		f \mapsto \of^{\ast}:\Shv(G \backslash Q) \to \Shv(G \backslash P),
		\]
		respectively. The companion $\overline{F}:\Var_{/K}^{\op} \to \fCat$ then sends a $K$-variety to its corresponding category of sheaves and sends a morphism of varieties to the corresponding pullback between sheaf categories. There is a topological version of this example given mutatis mutandis as well.
	\end{enumerate}
\end{example}

The definitions of pre-equivariant pseudofunctors for both schemes and spaces allow us to define equivariant categories as well as a convention of convenience we will use throughout the remainder of this work. As we work with these categories (especially in the cases of working with the categories $\SfResl_G(X)$ when $G$ is a smooth algebraic group and $X$ is a left $G$-variety or when $\mathbf{SRes}_G(X)$ when $G$ is a topological group which has $n$-acyclic free spaces which are manifolds and $X$ is a left $G$-space) we will need to be indexing our objects and morphisms in the pseudocone categories $\PC(F)$ by the quotient spaces $G \backslash (\Gamma \times X)$ or $G \backslash (M \times X)$; however, as we have seen already (cf.\@ Corollaries \ref{Cor: Pseudocone Functors: Monoidal Closed Derived Cat}, \ref{Cor: Pseudocone Functors: Monoidal Closed Derived Cat Topological}, \ref{Cor: Pseudocone Functors: Existence of equivariant pullback  for schemes}, \ref{Cor: Pseudocone Functors: Existence of equivariant pullback  for spaces}, \ref{Cor: Pseudocone Functors: Exceptional Pushforward functors for equivariant maps for scheme sheaves}, \ref{Cor: Pseudocone Functors: Pushforward functors for equivariant maps for scheme sheaves}, \ref{Cor: Pseudocone Functors: Topological EDC has equivariant exceptional pushpull}, \ref{Cor: Pseudocone Functors: Topological EDC has equivariant pushpull}, \ref{Cor: Section Triangle: SFF Relations for EDC Scheme}, for instance) this can be not only notationally cumbersome but also hides an important aspect to the story. Because each of the categories $\SfResl_G(X)$ carries an isomorphism $\SfResl_G(X) \cong \Sf(G)$ in the geometric case (and similarly in the topological case each category $\mathbf{SRes}_G(X)$ carries an isomorphism $\mathbf{SRes}_G(X) \cong \mathbf{FAMan}(G)$\index[notation]{FAMan@$\mathbf{FAMan}(G)$} where $\mathbf{FAMan}(G)$ denotes the category whose objects are free  $G$-spaces which are manifolds $n$-acyclic for some $n \in \N$ and whose morphisms are smooth $G$-equivariant morphisms between such resolutions). As such to a large degree what is important is that we are indexing our objects all by the smooth free $G$-spaces $\Gamma$, and baking this into our notation will not only reduce the notational clutter in our discussion and formalism but also will introduce a convenient formalism for keeping track of objects which may require multiple indexes.

\begin{definition}
	Let $G$ be a smooth algebraic group and let $X$ be a left $G$-variety. For any variety $\Gamma \in \Sf(G)_0$ we define the shorthand\index[notation]{XGamma@$\XGamma$}
	\[
	\XGamma := G \backslash (\Gamma \times X)
	\]
	Similarly, if $G$ is a topological group admitting free $n$-acyclic $G$-spaces which are manifolds, $X$ is a left $G$-space, and $M \in \mathbf{FAMan}(G)_0$, we define\index[notation]{XMamma@$\quot{X}{M}$}
	\[
	\quot{X}{M} := G \backslash (M \times X).
	\]
\end{definition}
\begin{definition}
	Let $G$ be a smooth algebraic group and let $X$ be a left $G$-variety. For any morphism $f:\Gamma \to \Gamma^{\prime}$ in $\Sf(G)$ we define the map
	\[
	\of:\XGamma \to \XGammap
	\]
	by taking $\of = \overline{f \times \id_X} = G \backslash(f \times \id_X)$. Similarly, if $G$ is if $G$ is a topological group admitting free $n$-acyclic $G$-spaces which are manifolds and $X$ is a left $G$-space, for any map $f:M \to N$ in $\mathbf{FAMan}(G)$ we define
	\[
	\of:\quot{X}{M} \to \quot{X}{N}
	\]
	by taking $\of = \overline{f \times \id_X} = G \backslash (f \times \id_X)$.
\end{definition}

We now define equivariant categories and show that there is a forgetful functor which ``forgets'' the $G$-equivariance recorded by the equivariant descent contained in $\SfResl_G(X)$ and the composite pseudofunctor $\overline{F}\circ \quo^{\op}$. In what follows we give the definition of equivariant categories for both topological and geometric situations, but will mainly focus on the geometric situation. Note that the two definitions differ functionally only in the domain of the pseudofunctors $F$.

\begin{definition}\label{Defn: Equivariant Category}
	Let $G$ be a smooth algebraic group and let $X$ be a left $G$-variety. Let $(F,\overline{F})$ be a pre-equivariant pseudofunctor on $X$. We define the equivariant category\index[terminology]{Equivariant Category} $F_G(X)$ on $X$ to be the category
	\[
	F_G(X) := \PC(F) = \PC(F \circ \quo^{\op}).
	\]
\end{definition}
\begin{definition}
	Let $G$ be a topological space and let $X$ be a left $G$-space. Let $(F,\overline{F})$ be a pre-equivariant pseudofunctor on $X$. We define the equivariant category $F_G(X)$ on $X$ to be the category
	\[
	F_G(X) := \PC(F) = \PC(F \circ \quo^{\op}).
	\]
\end{definition}

Let us do a quick analysis of the situation that surrounds an equivariant category and a pre-equivariant pseudofunctor. For any smooth algebraic group $G$, there is a non-equivariant isomorphism of varieties
\[
\psi:X \xrightarrow{\cong} \quot{X}{G} = G \backslash (G \times X)
\]
which is induced, in the set-theoretic language, by the map $x \mapsto [1,x]$. There is a similar non-equivariant isomorphism when $G$ is a topological group and $X$ is a left $G$-space as well. From Theorem \ref{Thm: Section Pseudocones: PCF is the pseudolimit of F} we know that
\[
\PC(F) = \PC(\overline{F} \circ \quo^{\op})
\]
and so there is a functor
\[
p_G:\PC(F) \to \overline{F}(\quot{X}{G})
\]
induced by the pseudouniversal property. Post-composing this with the equivalence $\overline{F}(\psi)$ then gives rise to a functor
\[
F_G(X) \xrightarrow{p_G} \overline{F}(\quot{X}{G}) \xrightarrow{\overline{F}(\psi)} \overline{F}(X)
\]
which serves as our forgetful functor\footnote{And in fact does so in a strong way by Proposition \ref{Prop: Section 3.3: Forgetful functor is de-equivariantification functor}.} (and in fact we will call this functor $\Forget$ because it forgot the ``$G$'' information).

\begin{proposition}\label{Prop: Section ECSV: Forgetful functor}
	Let $G$ and $X$ be given as in either case:
	\begin{enumerate}
		\item $G$ is a smooth algebraic group and $X$ is a left $G$-variety;
		\item $G$ is a topological group and $X$ is a left $G$-space;
	\end{enumerate} 
	and let $(F,\overline{F})$ be a pre-equivariant pseudofunctor on $X$. There is a forgetful functor
	\[
	\Forget:F_G(X) \to \overline{F}(X).
	\] 
	Moreover, there is also a functor $\Eq:\overline{F}(X) \to F_G(X)$ which gives rise to an invertible $2$-cell:
	\[
	\begin{tikzcd}
		& F_G(X) \ar[dr]{}{\Forget} \\
		\overline{F}(X) \ar[ur]{}{\Eq} \ar[rr, swap, equals, ""{name = D}]{}{} & & \overline{F}(X) \ar[from = 1-2, to = D, Rightarrow, shorten <= 4pt, shorten >= 4pt]{}{\nu}
		\ar[from = 1-2, to = D, Rightarrow, shorten <= 4pt, shorten >= 4pt, swap]{}{\cong}
	\end{tikzcd}
	\]
\end{proposition}
\begin{proof}
	The existence of $\Forget:F_G(X) \to \overline{F}(X)$ as $\Forget = \overline{F}(\psi) \circ p_G$ follows from the description prior to the statement of the proposition. We prove the existence of the $\Eq$ in the scheme-theoretic case; the topological case follows mutatis mutandis by replacing every $\Gamma$ with an $M$. 
	
	To see the existence of the equivariantification functor we note that for any $\Gamma \in \Sf(G)_0$, the morphism
	\[
	\begin{tikzcd}
		\Gamma \times X \ar[r]{}{\pi_2} & X \ar[r]{}{\psi^{-1}} & \quot{X}{G}
	\end{tikzcd}
	\]
	is $G$-trivial by virtue of ending in the object $\quot{X}{G}$. It then follows from the universal property of the quotient object $\XGamma$ that there is a unique morphism $\overline{\Gamma}$ making the diagram
	\[
	\begin{tikzcd}
		\Gamma \times X \ar[r]{}{\pi_2} \ar[d, swap]{}{\quo_{\Gamma}} & X \ar[d]{}{\psi^{-1}} \\
		\XGamma \ar[r, swap, dashed]{}{\exists!\overline{\Gamma}} & \quot{X}{G}
	\end{tikzcd}
	\]
	commute. The various morphisms $\overline{\Gamma}$ satisfy the property that for any $f:\Gamma \to \Gamma^{\prime}$ in $\Sf(G)$, the diagram
	\[
	\begin{tikzcd}
		& & & \quot{X}{G}\\
		\\
		& X \ar[uurr]{}{\psi} & \XGamma \ar[rr, swap]{}{\of} \ar[uur, swap]{}{\overline{\Gamma}} & & \XGammap \ar[uul]{}{\overline{\Gamma}^{\prime}}\\
		\\
		\Gamma \times X  \ar[rr, swap]{}{f \times \id_X} \ar[uurr, swap, near start]{}{\quo_{\Gamma}} & & \Gamma^{\prime} \times X \ar[uurr, swap]{}{\quo_{\Gamma^{\prime}}}
		\ar[from = 5-1, to = 3-2, crossing over]{}{\pi_2^{\Gamma}} \ar[from = 5-3, to = 3-2, swap, crossing over, near start]{}{\pi_2^{\Gamma^{\prime}}}
	\end{tikzcd}
	\]
	commutes in $\Var_{/K}$. From this we define, for any $A \in \overline{F}(X)_0$ the object
	\[
	\Eq A := \left\lbrace \overline{F}(\overline{\Gamma})\left(F(\psi^{-1})(A)\right) \; | \; \Gamma \in \Sf(G)_0 \right\rbrace
	\]
	with transition isomorphisms induced by the diagram
	\[
	\begin{tikzcd}
		\overline{F}(\of)\left(\overline{F}(\overline{\Gamma}^{\prime})\left(\overline{F}(\psi^{-1})A\right)\right) \ar[dd, equals] \ar[r]{}{\phi_{\of,\overline{\Gamma}}^{\overline{F}(\psi^{-1})A}} & \overline{F}(\overline{\Gamma}^{\prime} \circ \of)\left(\overline{F}(\psi^{-1})A\right) \ar[d, equals] \\
		& \left(\overline{F}(\overline{\Gamma}) \circ \overline{F}(\psi^{-1})\right)(A) \ar[d, equals]\\
		\overline{F}(\of)\left(\quot{\Eq A}{\Gamma^{\prime}}\right) \ar[r, swap]{}{\tau_f^{\Eq A}} & \quot{\Eq A}{\Gamma}
	\end{tikzcd}
	\]
	for any $f:\Gamma \to \Gamma^{\prime}$ in $\Sf(G)_1$. That these transition isomorphisms satisfy the cocycle condition is immediate from the pseudofuncoriality of $\overline{F}$, so we do indeed have an object of$F_G(X)$. The assignment on morphisms is given similarly: For any $\alpha \in \overline{F}(X)_1,$
	\[
	\Eq \alpha := \left\lbrace \overline{F}(\overline{\Gamma})\left(\overline{F}(\psi^{-1})\alpha\right) \; | \; \Gamma \in \Sf(G)_0 \right\rbrace
	\]
	and that these do indeed define morphisms follows from the fact that the compositors $\phi$ are natural isomorphisms. Thus we get our desired functor $\Eq$.
	
	For the claim on the natural isomorphism we simply calculate that
	\[
	(\Forget \circ \Eq)(A) = \overline{F}(\psi)\left(\overline{F}(\psi^{-1})A\right)
	\]
	so setting $\nu := \phi_{\psi,\psi^{-1}}$ proves the proposition.
\end{proof}

\newpage

%
%

\section{The Equivariance Isomorphisms}\label{Subsection: The Equivariance Isos}
We now study the geometric information carried by the objects of an equivariant category and precisely in what sense it is that the objects ``respect the $G$-action on $X$\footnote{Analogously to how each of the objects in $\Shv_G(X)$, $\Per_G(X)$, etc.\@ contain sheaves that respect the $G$-action on $X$.}.'' This ultimately means prove first that equivariant categories $F_G(X)$ are formally na{\"i}evely equivariant in the sense that there are functors $\pi_2^{\ast}, \alpha^{\ast}:F_G(X) \to F_G(G \times X)$ where $\pi_2$ is induced by the projection map $\pi_2:G \times X \to X$ and $\alpha_X^{\ast}:F_G(X) \to F_G(G \times X)$ is induced by the action map $\alpha_X:G \times X \to X$. We then use these functors to construct a natural isomorphism
\[
\Theta:\alpha_X^{\ast}A \xRightarrow{\cong} \pi_2^{\ast}A
\]
which satisfies a pseudofunctorial generalization of the cocycle condition of \cite[Definition 1.6]{GIT} (which we spell out later; cf.\@ Definition \ref{Defn: GIT Cocycle}). Our first step here, consequently, is to introduce some language and then construct the functors $\pi_2^{\ast}$ and $\alpha_X^{\ast}$.
\begin{remark}
	For the impatient, the GIT cocycle condition is the condition described in \cite[Definition 1.6]{GIT} as well as \cite[Chapter 0]{BernLun} and often stated as writing $\Theta$ as an isomorphism of sheaves $\Theta:d_0^{\ast}\Fscr \Rightarrow d_1^{\ast}\Fscr$ which satisfies the equation $d_0^{\ast}\Theta \circ d_2^{\ast}\Theta = d_1^{\ast}\Theta$. In this case the maps $d_0, d_1, d_2$ come from the simplicial presentation
	\[
	\begin{tikzcd}
	G \times G \times X \ar[r, shift left = 1.5 ex]{}{d_0} \ar[r]{}[description]{d_1} \ar[r, shift right = 1.5 ex, swap]{}{d_2} & G \times X \ar[r, shift left = 0.5ex]{}{d_0} \ar[r, shift right = 0.5ex]{}{d_1} & X
	\end{tikzcd}
	\]
	of the translation groupoid of the action $\alpha_X:G \times X \to X$.
	We explain and describe this condition in more detail in Definition \ref{Defn: GIT Cocycle} below.
\end{remark}

\begin{definition}
	Let $G$ be a smooth algebraic group and let $X$ and $Y$ be left $G$-varieties or let $G$ be a topological group and let $X$ and $Y$ be left $G$-spaces. We say that a simultaneous pre-equivariant equivariant\index[terminology]{Pre-equivariant Pseudofunctor! Simultaneous} pseudofunctor over $X$ and $Y$ is a triple $\left((F,\overline{F}),(F^{\prime},\overline{F}^{\prime}),e:\overline{F} \to \overline{F}^{\prime}\right)$ where $(F,\overline{F})$ is a pre-equivariant pseudofunctor on $X$, $(F^{\prime},\overline{F}^{\prime})$ is a pre-equivariant pseudofunctor on $Y$, and $e:\overline{F}\Rightarrow \overline{F}^{\prime}$ is an equivalence $\overline{F} \simeq \overline{F}^{\prime}$  in the bicategory $\Bicat(\Var_{/K}^{\op},\fCat)$.
\end{definition}
\begin{remark}
	It follows immediately from the pseudoconification $2$-functor described in Lemma \ref{Lemma: Section 3: The embedding of pre-equivariant functors to equivariant categories is a strict 2-functor} that there is an equivalence of categories $\PC(\overline{F}) \simeq \PC(\overline{F}^{\prime})$.
\end{remark}
\begin{lemma}\label{Lemma: Section Equivaraince: simultaneous preeq}
	For any smooth algebraic group $G$ and any left $G$-variety $X$ or any topological group $G$ and any left $G$-space $X$, there is an inclusion functor $\SfResl_G(G \times X) \to \SfResl_G(X)$ or $\mathbf{SRes}_G(G \times X) \to \mathbf{SRes}_G(X)$. In particular, any pre-equivariant pseudofunctor $(F, \overline{F})$ on $X$ induces a simultaneous pre-equivariant pseudofunctor over $X$ and $G \times X$.
\end{lemma}
\begin{proof}
	We show this for the scheme-theoretic situation; the topological case follows mutatis mutandis. The first claim in the lemma follows because each variety $\Gamma \times (G \times X) \in \SfResl_G(G \times X)_0$ also satisfies
	\[
	\Gamma \times (G \times X) \cong (\Gamma \times G) \times X \in \SfResl_G(X)_0
	\]
	because $\Gamma \times G \in \Sf(G)_0$. Call the corresponding inclusion functor 
	\[
	\incl:\SfResl_G(G \times X) \to \SfResl_G(X)
	\]
	and note that the corresponding quotient functor \[
	\quo:\SfResl_G(G \times X) \to \Var_{/K}
	\] 
	factors through $\incl$. As such, we get a commuting diagram:
	\[
	\begin{tikzcd}
		& \SfResl_G(X)^{\op} \ar[dr]{}{\quo^{\op}} & \\
		\SfResl_G(G \times X)^{\op} \ar[rr] \ar[dr, swap]{}{F \circ \incl^{\op}} \ar[ur]{}{\incl^{\op}} & & \Var_{/K}^{\op} \ar[dl, swap]{}{\overline{F}} \\
		& \fCat \ar[from = 1-2, to = 3-2, crossing over, near end]{}{F}
	\end{tikzcd}
	\]
	Consequently we define the pre-equivariant pseudofunctors comprising the simultaneous pre-equivariant pseudofunctor as $(F,\overline{F})$ on $X$ and $(F \circ \operatorname{incl}^{\op}, \overline{F})$ on $G \times X$ and then taking the equivalence $e$ to be the identity transformation on $\overline{F}$.
\end{proof}
\begin{remark}
Reading the proof above we see that it is even unnecessary to start with a pre-equivariant pseudofunctor on $X$. If we instead simply have a pseudofunctor $\overline{F}:\Var_{/K}^{\op} \to \fCat$ or $\overline{F}:\Top^{\op} \to \fCat$ then we get a pre-equivariant pseudofunctor on $X$ by simply setting $F$ to be the pre-composition of $\overline{F}$ by $\quo^{\op}$.
\end{remark}
\begin{lemma}\label{Lemma: Section Change of Space: pi2 pullback}
	For any pre-equivariant pseudofunctor $F$ over $X$, there is a functor
	\[
	\pi_2^{\ast}:F_G(X) \to F_G(G \times X).
	\]
\end{lemma}
\begin{proof}
	Since $\pi_2:G \times X \to X$ is $G$-equivariant, we apply Theorem \ref{Thm: Pseudocone Functors: Pullback induced by fibre functors in pseudofucntor} to produce $\pi_2^{\ast}$ as the functor induced by the map $\pi_2$.
\end{proof}
Constructing the functor $\alpha_X^{\ast}:F_G(X) \to F_G(G \times X)$ is more difficult than constructing $\pi_2^{\ast}$. The action morphism $\alpha_X:G \times X \to X$ is rarely $G$-equivariant; the corresponding diagram takes the form (where $\alpha_{G \times X}$ uses the diagonal $G$-action on $G \times X$)
\[
\xymatrix{
	G \times G \times X \ar[rr]^-{\alpha_{G \times X}} \ar[d]_{\id_G\times \alpha_X} & & G \times X \ar[d]^{\alpha_X} \\
	G \times X \ar[rr]_-{\alpha_X} & & X
}
\]
which commutes if and only if 
\[
\alpha_{X}(\alpha_{G \times X}(h,g,x)) = hghx = hgx = \alpha_X((\id_G \times \alpha_X)(h,g,x))
\]
for every $x \in X$ and $g,h \in G$. Because this implies the action on $X$ is trivial, we cannot assume this. As such, to construct the pullback of the action $\alpha_X^{\ast}:F_G(X) \to F_G(G \times X)$ we need to use the language of Change of Domain functors in order to overcome this obstruction. We will give two distinct arguments for this proof: one, which is presented as the main proof goes through the exact categorical structure used and illustrates how Theorem \ref{Thm: Pseudocone Functors: Pullback induced by fibre functors in pseudofucntor} determines the pullback against the action map. The other (more direct) argument given in Remark \ref{Remark: Second Proof of the action pullback} shows in a quick-and-dirty fashion how to avoid the categorical arguments and jump straight to the result.

\begin{lemma}\label{Lemma: Section Change of Space: Action map pullback}
	For any pre-equivariant pseudofunctor there is a functor
	\[
	\alpha_X^{\ast}:F_G(X) \to F_G(G \times X).
	\]
\end{lemma}
\begin{proof}
	We argue the scheme-theoretic case; the topological case follows similarly. Our strategy is to prove that $F$ admits pseudocone translations along $\alpha_X$. Begin letting $f:\Gamma \to \Gamma^{\prime}$ be a morphism in  $\Sf(G)$. Then because $f$ is $G$-equivariant we can produce a commuting cube:
	\[
	\begin{tikzcd}
		& \Gamma^{\prime} \times G \times X \ar[rr]{}{\id_{\Gamma^{\prime}} \times \alpha_{X}} \ar[dd, swap, near end]{}{\quo_{\Gamma^{\prime} \times G}} & & \Gamma^{\prime} \times X  \ar[dd]{}{\quo_{\Gamma^{\prime}}} \\
		\Gamma \times G \times X \ar[dd, swap]{}{\quo_{\Gamma \times G}} \ar[ur]{}{f \times \id_{G \times X}} & & \Gamma \times X \ar[ur]{}{f \times \id_X} \\
		& \quot{X}{\Gamma^{\prime} \times G} \ar[rr, near start, swap]{}{\overline{\id_{\Gamma^{\prime}} \times \alpha_X}} & & \XGammap \\
		\quot{X}{\Gamma \times G} \ar[ur]{}{\overline{f \times \id_{G \times X}}} \ar[rr, swap]{}{\overline{\id_{\Gamma} \times \alpha_X}} & & \XGamma \ar[ur, swap]{}{\overline{f \times \id_X}} \ar[near end, from = 2-1, to = 2-3, crossing over]{}{\id_{\Gamma} \times \alpha_X} \ar[from = 2-3, to = 4-3, near start, crossing over]{}{\quo_{\Gamma}}
	\end{tikzcd}
	\]
	Note that the equivariance of $f$ is used to deduce the commutativity of the vertically oriented faces of the cube while the top face commutes by construction. This in turn provides the commuting square of quotient varieties:
	\[
	\begin{tikzcd}
		\quot{(G \times X)}{\Gamma} \ar[rr, equals]{}{} \ar[d, swap]{}{\overline{f \times \id_{G \times X}}} & & \quot{X}{\Gamma \times G} \ar[rr]{}{\overline{\id_{\Gamma} \times \alpha_X}} \ar[d, swap]{}{\overline{f \times \id_{G \times X}}} & & \XGamma \ar[d]{}{\overline{f \times \id_X}} \\
		\quot{(G \times X)}{\Gamma^{\prime}} \ar[rr, equals] & & \quot{X}{\Gamma^{\prime} \times G} \ar[rr, swap]{}{\overline{\id_{\Gamma^{\prime}} \times \alpha_X}} & & \XGammap
	\end{tikzcd}
	\]
	Let $\ul{X}:\Sf(G) \to \SfResl_G(X)$ and $\ul{G \times X}:\Sf(G) \to \SfResl_G(G \times X)$ be the isomorphisms of Remark \ref{Remark: Pseudocone Section: SfResl is iso to Sf}. As the above diagrams of quotient varieties vary naturally in $\Gamma$ we induce a $2$-cell
	\[
	\begin{tikzcd}
		\SfResl_G(G \times X) \ar[dr, swap]{}{\incl} \ar[rr, ""{name = U}]{}{\quo} & & \Var_{/K} \\
		& \SfResl_G(X) \ar[ur, swap]{}{\quo} \ar[from = U, to = 2-2, Rightarrow, shorten <= 4pt, shorten >= 4pt]{}{\tilde{\alpha}_X}
	\end{tikzcd}
	\]
	where for all $\Gamma \in \Sf(G)_0$, $\quot{\tilde{\alpha}_X}{\Gamma} := \overline{\id_{\Gamma} \times \alpha_X}$. Post-composing by these isomorphisms gives rise to the $2$-cell
	\[
	\begin{tikzcd}
		\Sf(G)^{\op} \ar[rrr, ""{name = U}]{}{\ul{X}^{\op}} \ar[dd, swap]{}{(\ul{G\times X})^{\op}} & & & \SfResl_G(X)^{\op} \ar[dd]{}{\quo^{\op}} \\
		\\
		\SfResl_G(G \times X)^{\op} \ar[rrr, swap, ""{name = D}]{}{\quo^{\op}} & & & \Var_{/K}^{\op} \ar[from = U, to = D, Rightarrow, shorten <= 4pt, shorten >= 4pt]{}{\ul{\alpha}_X^{\op}} 
	\end{tikzcd}
	\]
	which factors as the pasting diagram
	\[
	\begin{tikzcd}
		\Sf(G)^{\op} \ar[rrrr, ""{name = U}]{}{\ul{X}^{\op}} \ar[ddd, swap]{}{(\ul{G\times X})^{\op}} & & & &\SfResl_G(X)^{\op} \ar[ddd]{}{\quo^{\op}} \\
		\\
		\\
		\SfResl_G(G \times X)^{\op} \ar[uuurrrr, ""{name = M}]{}[description]{\incl^{\op}}\ar[rrrr, swap, ""{name = D}]{}{\quo^{\op}} & & & & \Var_{/K}^{\op} \ar[from = U, to = M, Rightarrow, shorten <= 4pt, shorten >= 4pt]{}[description]{(\id_{-} \times \alpha_X)^{\op}} 
	\end{tikzcd}
	\]
	where the empty cell is used as a shorthand for the identity $2$-cell. Post-composing by the pseudofunctor $\overline{F}:\Var^{\op} \to \fCat$ then provides us with the induced $2$-cell
	\[
	\begin{tikzcd}
		\Sf(G)^{\op} \ar[rrr, bend left = 30, ""{name = U}]{}{\overline{F} \circ \quo^{\op} \circ \ul{X}^{\op}} \ar[rrr, bend right = 30, swap, ""{name = D}]{}{\overline{F} \circ \quo^{\op} \circ \incl^{\op} \circ (\ul{G \times X})^{\op}} & & & \fCat \ar[from = U, to = D, Rightarrow, shorten <= 4pt, shorten >= 4pt]{}{\ul{\alpha}_X}
	\end{tikzcd}
	\]
	where $\ul{\alpha}_X = \iota_{\overline{F}} \ast \iota_{\quo^{\op} \circ \incl^{\op}} \ast (\id_{-} \times \alpha_X)^{\op}$. Applying Theorem \ref{Thm: Pseudocone Functors: Pullback induced by fibre functors in pseudofucntor} then gives rise to a functor
	\[
	\hat{\alpha}_X^{\square}:\PC\left(\overline{F} \circ \quo^{\op} \circ \ul{X}^{\op}\right) \longrightarrow \PC\left(\overline{F} \circ \quo^{\op} \circ \incl^{\op} \circ (\ul{G \times X})^{\op}\right)
	\]
	as induced by $\overline{F}$. However, observe that since $\ul{X}$ and $\ul{G \times X}$ are isomorphisms, there are induced equivalences of categories
	\[
	\PC\left(\overline{F} \circ \quo^{\op} \circ \ul{X}^{\op}\right) \simeq \PC\left(\overline{F}\circ \quo^{\op}\right)
	\]
	and
	\[
	\PC\left(\overline{F} \circ \quo^{\op} \circ \incl^{\op} \circ (\ul{G \times X})^{\op}\right) \simeq \PC\left(\overline{F} \circ \quo^{\op} \circ \incl^{\op}\right).
	\]
	We now use Lemma \ref{Lemma: Section Equivaraince: simultaneous preeq} to write
	\[
	\PC\left(\overline{F} \circ \quo^{\op}\right) = F_G(X), \qquad \PC\left(\overline{F} \circ \quo^{\op} \circ \incl^{\op}\right) = F_G(G \times X)
	\]
	and hence produce our desired functor
	\[
	\alpha_X^{\ast}:F_G(X) \to F_G(G \times X).
	\]
\end{proof}
\begin{remark}[A Second Proof of Lemma \ref{Lemma: Section Change of Space: Action map pullback}]\label{Remark: Second Proof of the action pullback}
	Proceed as in the first proof above to produce the diagram of varieties
	\[
	\begin{tikzcd}
		\quot{(G \times X)}{\Gamma} \ar[rr, equals]{}{} \ar[d, swap]{}{\overline{f \times \id_{G \times X}}} & & \quot{X}{\Gamma \times G} \ar[rr]{}{\overline{\id_{\Gamma} \times \alpha_X}} \ar[d, swap]{}{\overline{f \times \id_{G \times X}}} & & \XGamma \ar[d]{}{\overline{f \times \id_X}} \\
		\quot{(G \times X)}{\Gamma^{\prime}} \ar[rr, equals] & & \quot{X}{\Gamma^{\prime} \times G} \ar[rr, swap]{}{\overline{\id_{\Gamma^{\prime}} \times \alpha_X}} & & \XGammap
	\end{tikzcd}
	\]
	and now apply the functor $\overline{F}$ to produce the pasting diagram:
	\[
	\begin{tikzcd}
		\overline{F}\left(\XGammap\right) \ar[d, swap]{}{\overline{F}\left(\overline{f \times \id_X}\right)} \ar[rr, ""{name = U}]{}{\overline{F}\left(\overline{\id_{\Gamma^{\prime}} \times \alpha_X}\right)} & & \overline{F}\left(\quot{X}{\Gamma^{\prime} \times G}\right) \ar[rr, equals] \ar[d]{}{\overline{F}\left(\overline{f \times \id_{G \times X}}\right)} & & \overline{F}\left(\quot{(G \times X)}{\Gamma^{\prime}}\right) \ar[d]{}{\overline{F}\left(\overline{f \times \id_{G \times X}}\right)}\\
		\overline{F}\left(\XGamma\right)  \ar[rr, swap, ""{name = D}]{}{\overline{F}\left(\overline{\id_{\Gamma} \times \alpha_X}\right)}  & & \overline{F}\left(\quot{X}{\Gamma \times G}\right)  \ar[rr, equals]{}{} & & 	\overline{F}\left(\quot{(G \times X)}{\Gamma}\right) \ar[from = U, to = D, Rightarrow, shorten <= 4pt, shorten >= 4pt]{}{\phi}
	\end{tikzcd}
	\]
	where 
	\[
	\phi = \phi^{-1}_{\overline{\id_{\Gamma} \times \alpha_X},\overline{f \times \id_X}}  \circ \phi_{\overline{f \times \id_{G \times X}}, \overline{\id_{\Gamma^{\prime}} \times \alpha_X}}.
	\]
	It then follows from a verification mutatis mutandis to Theorems \ref{Thm: Functor Section: Psuedonatural trans are pseudocone functors} and \ref{Thm: Pseudocone Functors: Pullback induced by fibre functors in pseudofucntor} that this gives rise to a functor $F_G(X) \to F_G(G \times X)$ induced by the assignments, for $A \in F_G(X)_0$,
	\[
	\alpha_X^{\ast}A = \left\lbrace F(\quot{\overline{\alpha}_X}{\Gamma})\left(\AGamma\right) \; | \; \Gamma \in \Sf(G)_0, \AGamma \in A\right\rbrace
	\]
	on object collecations with the transition isomorphisms given by
	\[
	T_{\alpha_X^{\ast}A} = \left\lbrace F(\quot{\overline{\alpha}_X^{\ast}}{\Gamma})(\tau_f^A) \circ \phi^{-1}_{\overline{\id_{\Gamma} \times \alpha_X},\overline{f \times \id_X}}  \circ \phi_{\overline{f \times \id_{G \times X}}, \overline{\id_{\Gamma^{\prime}} \times \alpha_X}} \; : \; f \in \Sf(G)_1 \right\rbrace.
	\]
\end{remark}
\begin{remark}
	In what follows, we will write
	\[
	\quot{\overline{\alpha}_X}{\Gamma} := \overline{\id_{\Gamma} \times \alpha_X}:\quot{X}{\Gamma \times G} \to \XGamma
	\]\index[notation]{alphaoverline@$\quot{\overline{\alpha}_X}{\Gamma}$}
	in order to reduce the already significant notational bloat.
\end{remark}
\begin{remark}\label{Remark: Section Change of Space: Explict pullback functors}
	In the interest of being explicit, we note that the functor $\pi_2^{\ast}:F_G(X) \to F_G(G \times X)$ is given as follows on objects. Let $(A,T_A) \in F_G(X)_0$. By Theorem \ref{Thm: Pseudocone Functors: Pullback induced by fibre functors in pseudofucntor}
	\[
	\pi_2^{\ast}A = \left\lbrace F(\quot{\opi_2}{\Gamma})\left(\AGamma\right) \; | \; \Gamma \in \Sf(G)_0, \AGamma \in A \right\rbrace,
	\]
	where $\quot{\pi_2}{\Gamma}$ is the projection map $\Gamma \times G \times X \to \Gamma \times X$, and the transition isomorphisms are given as
	\[
	T_{\pi_2^{\ast}A} = \left\lbrace \overline{F}(\quot{\opi_2}{\Gamma})\tau_f^A \circ \left(\phi_{\quot{\pi_2}{\Gamma}, \quot{f}{X}}^{\AGammap}\right)^{-1} \circ \phi_{\quot{f}{G \times X}, \quot{\pi_2}{\Gamma^{\prime}}}^{\AGammap}  \; : \; f \in \Sf(G)_1 \right\rbrace.
	\]
	The topological version of the functor is similar save with every $\Gamma$ replaced with an $M$.
\end{remark}
\begin{remark}
	In order to reduce the (yet even further) notational clutter creeping its way into this section, if $f:\Gamma \to \Gamma^{\prime}$ is a morphism in $\Sf(G)$ and if $h:X \to Y$ is a morphism of $G$-varieties, we write
	\[
	\quot{\of}{X} := \overline{f \times \id_X}, \quot{\of}{Y} := \overline{f \times \id_{Y}}, \quot{\overline{h}}{\Gamma} := \overline{\id_{\Gamma} \times h},
	\]
	and $\quot{\overline{h}}{\Gamma^{\prime}} := \overline{\id_{\Gamma^{\prime}} \times h}$.
\end{remark}

In order to produce the isomorphism $\Theta$, we follow a process discovered by Clifton Cunningham during the BIRS focused research group ``the Voganish Project'' in 2019. The main idea here is to find, for each variety $\Gamma \in \Sf(G)_0$, an $\Sf(G)$-variety $\Gamma_c$ which induces an isomorphism of quotient varieties
\[
\quot{X}{\Gamma_c} \cong \quot{(G \times X)}{\Gamma}
\]
that is natural over $\Sf(G)$ with respect to the actions of $\quot{\overline{\alpha}_X}{\Gamma}$ and $\quot{\opi_2}{\Gamma}$ so that we can use the transition isomorphisms of an object of $F_G(X)$ to build an isomorphism with the objects $\pi_2^{\ast}A$ and $\alpha_X^{\ast}A$. We present this as the lemma below, which is the main ingredient in cooking up our desired isomorphism $\Theta$.
\begin{lemma}\label{Lemma: Clifton's Naivification Variety}\index[notation]{Gammac@$\Gamma_c$}
	For all $\Gamma \in \Sf(G)_0$ there exists a variety $\Gamma_c \in \Sf(G)_0$ such that there is an isomorphism of quotient varieties
	\[
	\quot{(G \times X)}{\Gamma} \xrightarrow{\mu_{\Gamma}} \quot{X}{\Gamma_c}.
	\]
	Moreover, if $f \in \Sf(G)_1$ with $f \in \Sf(G)(\Gamma, \Gamma^{\prime})$ there is an $\Sf(G)$-morphism 
	\[
	f_c:\Gamma_c \to \Gamma^{\prime}_c
	\] 
	for which the diagram
	\[
	\xymatrix{
		\quot{(G \times X)}{\Gamma} \ar[rr]^-{\mu_{\Gamma}} \ar[d]_{\quot{\of}{G \times X}} & & \quot{X}{\Gamma_c} \ar[d]^{\overline{f}_c} \\
		\quot{(G \times X)}{\Gamma^{\prime}} \ar[rr]_-{\mu_{\Gamma^{\prime}}} & & \quot{X}{\Gamma^{\prime}_c}
	}
	\]
	commutes. There are also $\Sf(G)$-morphisms $a_{\Gamma}, p_{\Gamma}:\Gamma_c \to \Gamma$ for which the diagrams
	\[
	\xymatrix{
		\quot{(G \times X)}{\Gamma} \ar[r]^-{\mu_{\Gamma}} \ar[dr]_{\quot{\overline{\alpha}_X}{\Gamma}} & \quot{X}{\Gamma_c} \ar[d]^{\overline{a}_{\Gamma}} & \quot{(G \times X)}{\Gamma} \ar[r]^-{\mu_{\Gamma}} \ar[dr]_{\quot{\opi_2}{\Gamma}} & \quot{X}{\Gamma_c} \ar[d]^{\overline{p}_{\Gamma}}\\
		& \XGamma & & \XGamma
	}
	\]
	commute. Finally, for any $f \in \Sf(G)_1$ the diagrams
	\[
	\xymatrix{
		\Gamma_c \ar[d]_{f_c} \ar[r]^-{a_{\Gamma}} & \Gamma \ar[d]^{f} & \Gamma_c \ar[r]^-{p_{\Gamma}} \ar[d]_{f_c} & \Gamma \ar[d]^{f} \\
		\Gamma_c^{\prime} \ar[r]_-{a_{\Gamma^{\prime}}} & \Gamma^{\prime} & \Gamma_c^{\prime} \ar[r]_-{p_{\Gamma^{\prime}}} & \Gamma^{\prime}
	}
	\]
	commute as well.
\end{lemma}
\begin{proof}
	We give this proof by doing a sketch using the points of the quotient varieties for readability. Let $\Gamma \in \Sf(G)_0$. Then define $\Gamma_c$ by setting
	\[
	\Gamma_c := \Gamma \times G
	\]
	with action map
	\[
	\alpha_{\Gamma_c} = \alpha_{\Gamma} \times \Ad
	\]
	where $\Ad$ is the adjoint action of $G$ on itself. More explicitly, $\Ad:G \times G \to G$ is the morphism given on sets by $\Ad(h,g) = hgh^{-1}$ and in categorical terms by
	\[
	\xymatrix{
		G \times G \ar[d]_{\Delta \times \id_G} \ar[rrr]^-{\Ad} & & & G \\
		G \times G \times G \ar[rr]_-{\id_G \times \inv \times \id_G} & & G \times G \times G \ar[r]_-{\id_G \times s} & G \times G \times G \ar[u]_{\mu \circ (\id_G \times \mu)}
	}
	\]
	where $\inv:G \to G$ is the inversion map and $s$ is the switching isomorphism of the product. Now because $\Sf(G)$ admits products and both $\Gamma$ and $G$ are $\Sf(G)$-varieties, $\Gamma_c$ is a pure dimensional smooth variety. That it is a principal $G$-variety follows from the fact that $\Gamma \in \Sf(G)_0$ and Proposition \ref{Prop: Smooth quotient from principal G var}.
	
	We define the morphism of quotient varieties $\mu_{\Gamma}:\quot{(G \times X)}{\Gamma} \to \quot{X}{\Gamma_c}$ by (in set-theoretic terms)
	\[
	[\gamma, (g,x)]_G \mapsto [(\gamma, g),x]_G;
	\]
	that this is an isomorphism is trivial to verify. Define $f_c:\Gamma_c \to \Gamma^{\prime}_c$ by $f_c := f \times \id_G$. Then we calculate that
	\begin{align*}
		(\of_c \circ \mu_{\Gamma})[\gamma, (g,x)]_G &= \of_c[(\gamma, g),x]_G = [(f(\gamma),g),x]_G = \mu_{\Gamma^{\prime}}[f(\gamma),(g,x)]_G \\
		&= (\mu_{\Gamma^{\prime}} \circ \quot{\of}{G \times X})[\gamma, (g,x)]_G,
	\end{align*}
	which shows that the first diagram commutes.
	
	We define the morphisms $a_{\Gamma}, p_{\Gamma}:\Gamma_c \to \Gamma$.  The map $a_{\Gamma}$ is defined by the diagram
	\[
	\xymatrix{
		\Gamma_c \ar[d]_{a_{\Gamma}} \ar@{=}[r] & \Gamma \times G \ar[r]_-{\id_{\Gamma} \times \inv} & \Gamma \times G \ar[d]_{s}\\
		\Gamma & & G \times \Gamma \ar[ll]^-{\alpha_{\Gamma}}
	}
	\]
	(which says in set-theoretic terms that $a_{\Gamma}(\gamma,g) = g^{-1}\gamma$) while we define $p_{\Gamma}$ $=$ $\pi_2^{\Gamma G}:\Gamma \times G \to G$ to be the second projection. Because $\pi_2^{\Gamma G}$ is trivially verified to be an $\Sf(G)$-morphism, we only need to do this for $a_{\Gamma}$. For $a_{\Gamma}$ we note that it is smooth and of constant fibre dimension, so we only need to show that it is $G$-equivariant. We calculate that on one hand
	\[
	(a_{\Gamma} \circ \alpha_{\Gamma_c})(h,\gamma,g) = a_{\Gamma}(h\gamma, hgh^{-1}) = (hgh^{-1})^{-1}h\gamma = hg^{-1}hh^{-1}\gamma = hg^{-1}\gamma
	\]
	while on the other hand
	\[
	\big(\alpha_{\Gamma} \circ (\id_{G} \times a_{\Gamma})\big)(h,\gamma,g) = (h,g^{-1}\gamma) = hg^{-1}\gamma.
	\]
	Thus $a_{\Gamma}$ is $G$-equivariant and hence an $\Sf(G)$-morphism.
	
	We now verify the commutativity of the triangles. For this note that the triangle involving $p_{\Gamma}$ commutes trivially, so we only verify the triangle for $a_{\Gamma}$. In this case we find that
	\begin{align*}
		(\overline{a}_{\Gamma} \circ \mu_{\Gamma})[\gamma, (g,x)]_G &= \overline{a}_{\Gamma}[(\gamma,g),x]_G = [g^{-1}\gamma,x]_G = [\gamma,gx]_G \\
		&=\quot{\overline{\alpha}_X}{\Gamma}[\gamma, (g, x)]_G,
	\end{align*}
	as was desired. Finally if $f \in \Sf(G)_1$, the last two squares commute trivially in the cases of $p_{\Gamma}, p_{\Gamma^{\prime}}$ and exactly because $f$ is $G$-equivariant in the cases of $a_{\Gamma}, a_{\Gamma^{\prime}}$.
\end{proof}
Because of the set-theoretic nature of the sketch given above we get the following topological version of Cunningham's na{\"i}vification variety.
\begin{lemma}\label{Lemma: Clifton's Naivification Space}
For all $M \in \mathbf{FAMan}(G)_0$ there exists a free $G$-acyclic manifold $M_c \in \mathbf{FAMan}(G)_0$ such that there is an isomorphism of quotient varieties
\[
\quot{(G \times X)}{M} \xrightarrow{\mu_{M}} \quot{X}{M_c}.
\]
Moreover, if $f \in \mathbf{FAMan}(G)_1$ with $f \in \mathbf{FAMan}(G)(M, N)$ there is a $\mathbf{FAMan}(G)$-morphism 
\[
f_c:M_c \to N_c
\] 
for which the diagram
\[
\xymatrix{
	\quot{(G \times X)}{M} \ar[rr]^-{\mu_{M}} \ar[d]_{\quot{\of}{G \times X}} & & \quot{X}{M_c} \ar[d]^{\overline{f}_c} \\
	\quot{(G \times X)}{N} \ar[rr]_-{\mu_{N}} & & \quot{X}{N_c}
}
\]
commutes. There are also $\mathbf{FAMan}(G)$-morphisms $a_{M}, p_{M}:M_c \to M$ for which the diagrams
\[
\xymatrix{
	\quot{(G \times X)}{M} \ar[r]^-{\mu_{M}} \ar[dr]_{\quot{\overline{\alpha}_X}{M}} & \quot{X}{M_c} \ar[d]^{\overline{a}_{M}} & \quot{(G \times X)}{M} \ar[r]^-{\mu_{M}} \ar[dr]_{\quot{\opi_2}{M}} & \quot{X}{M_c} \ar[d]^{\overline{p}_{M}}\\
	& \quot{X}{M} & & \quot{X}{M}
}
\]
commute. Finally, for any $f \in \mathbf{FAMan}(G)_1$ the diagrams
\[
\xymatrix{
	M_c \ar[d]_{f_c} \ar[r]^-{a_{\Gamma}} & M \ar[d]^{f} & M_c \ar[r]^-{p_{M}} \ar[d]_{f_c} & M \ar[d]^{f} \\
	N_c \ar[r]_-{a_{N}} & N & N_c \ar[r]_-{p_{N}} & N
}
\]
commute as well.
\end{lemma}

We now finally have the tools to provide our claimed isomorphism $\Theta:\alpha_X^{\ast}A \to \pi_2^{\ast}A$ for any object of $F_G(X)$. The construction of this isomorphism, in view of Lemma \ref{Lemma: Clifton's Naivification Variety} is straightforward, but the proof that it is an $F_G(X)$-morphism is quite difficult and involved. The techniques we use are essentially a generalization of an unpublished proof of C.\@ Cunningham and myself for the equivariant derived category, as the proof we gave establishes (after generalizing) the case when $F$ is a strict pseudofunctor.

\begin{Theorem}\label{Theorem: Formal naive equivariance}
Let $G$ be a smooth algebraic group and $X$ a left $G$-variety or let $G$ be a topological group and $X$ a left $G$-space. For any pre-equivariant pseudofunctor $F = (F, \overline{F})$ over $X$ and any $A \in F_G(X)_0$, there is an isomorphism
	\[
	\Theta_A:\alpha_X^{\ast}A \to \pi_2^{\ast}A
	\]
	such that for any $P:A \to B$ in $F_G(X)$ the diagram
	\[
	\xymatrix{
		\alpha_X^{\ast}A \ar[d]_{\alpha_X^{\ast}P} \ar[r]^-{\Theta_A} & \pi_2^{\ast}A \ar[d]^{\pi_2^{\ast}P} \\
		\alpha_X^{\ast}B \ar[r]_-{\Theta_B} & \pi_2^{\ast}B
	}
	\]
	commutes in $F_G(G \times X)$.
\end{Theorem}
\begin{proof}
We prove this in the variety situation, as the topological case follows similarly. Begin by fixing $A \in F_G(X)_0$ and defining the isomorphism $\quot{\theta}{\Gamma}$ as in the diagram below:
	\[
	\begin{tikzcd}
		\overline{F}(\quot{\overline{\alpha}_X}{\Gamma})(\AGamma) \ar[rrr]{}{\quot{\theta}{\Gamma}} \ar[d, equals]  & & & \overline{F}(\quot{\opi_2}{\Gamma}) \\
		\overline{F}(\aGamma \circ \muGamma)(\AGamma) \ar[d, swap]{}{\left(\phi_{\muGamma,a_{\Gamma}}^{\AGamma}\right)^{-1}} & & & \overline{F}(\pGamma \circ \muGamma)(\AGamma) \ar[u, equals] \\
		\overline{F}(\muGamma)\left(\overline{F}(\aGamma)(\AGamma)\right) \ar[r, swap]{}{\overline{F}(\muGamma)\tau_{a_{\Gamma}}^{A}} & \overline{F}(\muGamma)\quot{A}{\Gamma_c} \ar[rr, swap]{}{\overline{F}(\muGamma)(\tau_{p_{\Gamma}}^{A})^{-1}} & & \overline{F}(\muGamma)\left(\overline{F}(\pGamma)(\AGamma)\right) \ar[u, swap]{}{\phi_{\muGamma, p_{\Gamma}}}
	\end{tikzcd}
	\]
	We now begin our arduous task of showing that the equation 
	\[
	\tau_f^{\pi_2^{\ast}A} \circ \overline{F}(\fGX)\quot{\theta}{\Gamma^{\prime}} = \quot{\theta}{\Gamma} \circ \tau_f^{\alpha_X^{\ast}A}
	\] 
	holds. For this we first find
	\begin{align*}
		&\tau_f^{\pi_2^{\ast}A} \circ \overline{F}(\fGX)\quot{\theta}{\Gamma^{\prime}} \\
		&= \tau_f^{\pi_2^{\ast}A} \circ \overline{F}(\fGX)\left(\phi_{\muGammap,p_{\Gamma^{\prime}}}^{\AGammap} \circ \overline{F}(\muGammap) \left(\left(\tau_{p_{\Gamma^{\prime}}}^{A}\right)^{-1} \circ \tau_{a_{\Gamma^{\prime}}}^{A}\right) \circ \left(\phi_{\muGammap, a_{\Gamma^{\prime}}}^{\AGammap}\right)^{-1}\right) \\
		&=\tau_f^{\pi_2^{\ast}A} \circ \overline{F}(\fGX)\left(\phi_{\muGammap,p_{\Gamma^{\prime}}}^{\AGammap}\right) \circ  \big(\overline{F}(\fGX) \circ \overline{F}(\muGammap)\big)\left(\left(\tau_{p_{\Gamma^{\prime}}}^{A}\right)^{-1} \circ \tau_{a_{\Gamma^{\prime}}}^{A}\right) \\
		&\circ  \overline{F}(\fGX)\left(\phi_{\muGammap, a_{\Gamma^{\prime}}}^{\AGammap}\right)^{-1}.
	\end{align*} 
	Now using the naturality of the compositor isomorphisms we calculate that 
	\begin{align*}
		& \big(\overline{F}(\fGX) \circ \overline{F}(\muGammap)\big)\left(\left(\tau_{p_{\Gamma^{\prime}}}^{A}\right)^{-1} \circ \tau_{a_{\Gamma^{\prime}}}^{A}\right) \\
		&= \left(\phi_{\fGX, \muGammap}^{\overline{F}(\pGammap)\AGammap}\right)^{-1} \circ  \overline{F}(\muGammap \circ \fGX)\left(\left(\tau_{p_{\Gamma^{\prime}}}^{A}\right)^{-1} \circ\tau_{a_{\Gamma^{\prime}}}^{A}\right) \circ \phi_{\fGX, \muGammap}^{\overline{F}(\aGammap)\AGammap}.
	\end{align*}
	Using that there is an equality of functors
	\[
	\overline{F}(\muGammap \circ \fGX) = \overline{F}(\of_c \circ \muGamma)
	\]
	we have
	\begin{align*}
		&\overline{F}(\muGammap \circ \fGX)\left(\left(\tau_{p_{\Gamma^{\prime}}}^{A}\right)^{-1} \circ\tau_{a_{\Gamma^{\prime}}}^{A}\right) =\overline{F}(\of_c \circ \muGamma)\left(\left(\tau_{p_{\Gamma^{\prime}}}^{A}\right)^{-1} \circ\tau_{a_{\Gamma^{\prime}}}^{A}\right)
	\end{align*}
	so using the same trick we write
	\begin{align*}
		&\overline{F}(\of_c \circ \muGamma)\left(\left(\tau_{p_{\Gamma^{\prime}}}^{A}\right)^{-1} \circ\tau_{a_{\Gamma^{\prime}}}^{A}\right) \\
		&= \phi_{\muGamma, \of_c}^{\overline{F}(\aGammap)\AGammap} \circ \big(\overline{F}(\muGamma) \circ \overline{F}(\of_c)\big)\left(\left(\tau_{p_{\Gamma^{\prime}}}^{A}\right)^{-1} \circ\tau_{a_{\Gamma^{\prime}}}^{A}\right) \circ \left(\phi_{\muGamma, \of_c}^{\overline{F}(\aGammap)\AGammap}\right)^{-1}
	\end{align*}
	and get that
	\begin{align*}
		&\tau_f^{\pi_2^{\ast}A} \circ \overline{F}(\fGX)\quot{\theta}{\Gamma^{\prime}} \\
		&= \tau_f^{\pi_2^{\ast}A} \circ \overline{F}(\fGX)\left(\phi_{\muGammap,p_{\Gamma^{\prime}}}^{\AGammap} \circ \overline{F}(\muGammap) \left(\left(\tau_{p_{\Gamma^{\prime}}}^{A}\right)^{-1} \circ \tau_{a_{\Gamma^{\prime}}}^{A}\right) \circ \left(\phi_{\muGammap, a_{\Gamma^{\prime}}}^{\AGammap}\right)^{-1}\right) \\
		&=\tau_f^{\pi_2^{\ast}A} \circ \overline{F}(\fGX)\left(\phi_{\muGammap,p_{\Gamma^{\prime}}}^{\AGammap}\right) \circ  \big(\overline{F}(\fGX) \circ \overline{F}(\muGammap)\big)\left(\left(\tau_{p_{\Gamma^{\prime}}}^{A}\right)^{-1} \circ \tau_{a_{\Gamma^{\prime}}}^{A}\right) \\
		&\circ  \overline{F}(\fGX)\left(\phi_{\muGammap, a_{\Gamma^{\prime}}}^{\AGammap}\right)^{-1} \\
		&=\tau_f^{\pi_2^{\ast}A} \circ \overline{F}(\fGX)\left(\phi_{\muGammap,p_{\Gamma^{\prime}}}^{\AGammap}\right) \circ  \left(\phi_{\fGX, \muGammap}^{\overline{F}(\pGammap)\AGammap}\right)^{-1} \circ  \phi_{\muGamma, \of_c}^{\overline{F}(\pGammap)\AGammap} \\
		&\circ \big(\overline{F}(\muGamma) \circ \overline{F}(\of_c)\big)\left(\left(\tau_{p_{\Gamma^{\prime}}}^{A}\right)^{-1} \circ\tau_{a_{\Gamma^{\prime}}}^{A}\right) \circ \left(\phi_{\muGamma, \of_c}^{\overline{F}(\aGammap)\AGammap}\right)^{-1} \circ \phi_{\fGX, \muGammap}^{\overline{F}(\aGammap)\AGammap} \\
		&\circ  \overline{F}(\fGX)\left(\phi_{\muGammap, a_{\Gamma^{\prime}}}^{\AGammap}\right)^{-1}
	\end{align*}
	while
	\begin{align*}
		\tau_f^{\pi_2^{\ast}A} &= \overline{F}(\opiGamma)\tau_f^A \circ \left(\phi_{\quot{\pi_2}{\Gamma}, \fX}^{\AGammap}\right)^{-1} \circ \phi_{\fGX, \quot{\pi_2}{\Gamma}}^{\AGammap} \\
		&= \overline{F}(\pGamma \circ \muGamma)\tau_f^A \circ \left(\phi_{\pGamma \circ \muGamma, \fX}^{\AGammap}\right)^{-1} \circ \phi_{\fGX, \pGammap \circ \muGammap}^{\AGammap} \\
		&= \phi_{\muGamma, \pGamma}^{\AGamma} \circ \big(\overline{F}(\muGamma) \circ \overline{F}(\pGamma)\big)\tau_f^A \circ \left(\phi_{\muGamma, \pGamma}^{\overline{F}(\fX)\AGammap}\right)^{-1} \circ \left(\phi_{\pGamma \circ \muGamma, \fX}^{\AGammap}\right)^{-1} \\
		&\circ \phi_{\fGX, \pGammap \circ \muGammap}^{\AGammap}.
	\end{align*}
	In order to manipulate the (implicit) monster composition, we need to derive some identities with the compositor natural isomorphisms to simplify the expression. We will show that we can simplify the morphism
	\begin{align*}
		&\left(\phi_{\muGamma, \pGamma}^{\overline{F}(\fX)\AGammap}\right)^{-1} \circ \left(\phi_{\pGamma \circ \muGamma, \fX}^{\AGammap}\right)^{-1} \circ \phi_{\fGX, \pGammap \circ \muGammap}^{\AGammap} \circ \overline{F}(\fGX)\left(\phi_{\muGammap,p_{\Gamma^{\prime}}}^{\AGammap}\right) \\
		&\circ  \left(\phi_{\fGX, \muGammap}^{\overline{F}(\pGammap)\AGammap}\right)^{-1}
	\end{align*}
	by manipulating the coherences of the compositors.
	
	For this we first calculate, using the pasting diagrams of pseudofunctors, that the diagram
	\[
	\begin{tikzcd}
		\overline{F}(\XGammap) \ar[r]{}{\overline{F}(\fX)} \ar[rr, bend right = 60, swap, ""{name = L}]{}{\overline{F}(\fX \circ \pGamma \circ \muGamma)} & \overline{F}(\Gamma) \ar[r, bend left = 30, ""{name = U}]{}{\overline{F}(\muGamma) \circ \overline{F}(\pGamma)} \ar[r, bend right = 30, swap, ""{name = M}]{}{\overline{F}(\pGamma \circ \muGamma)} & \overline{F}(\quot{(G \times X)}{\Gamma}) \ar[from = U, to = M, Rightarrow, shorten <= 4pt, shorten >= 4pt] \ar[from = 1-2, to = L, Rightarrow, shorten <= 4pt, shorten >= 4pt]{}{}
	\end{tikzcd}
	\]
	is equivalent to the pasting diagram:
	\[
	\begin{tikzcd}
		\overline{F}(\XGammap) \ar[r, bend left = 30, ""{name = U}]{}{\overline{F}(\pGamma) \circ \overline{F}(\fX)} \ar[r, bend right = 30, swap, ""{name = M}]{}{\overline{F}(\fX \circ \pGamma)} \ar[rr, bend right = 70, swap, ""{name = L}]{}{\overline{F}(\fX \circ \pGamma \circ \muGamma)} & \overline{F}(\quot{X}{\Gamma_c}) \ar[r]{}{\overline{F}(\muGamma)} & \overline{F}(\quot{(G \times X)}{\Gamma}) \ar[from = U, to = M, Rightarrow, shorten <= 4pt, shorten >= 4pt] \ar[from = 1-2, to = L, Rightarrow, shorten <= 4pt, shorten >= 4pt]{}{}
	\end{tikzcd}
	\]
	Note that in both cases, the $2$-cells are given by the corresponding compositors. This gives us that
	\[
	\left(\phi_{\mu_{\Gamma},\pGamma}^{\overline{F}(\fX)\AGammap}\right)^{-1} \circ \left(\phi_{\pGamma \circ \muGamma, \fX}^{\AGammap}\right)^{-1} = \overline{F}(\muGamma)\left(\phi_{\pGamma, \fX}^{\AGammap}\right)^{-1} \circ \left(\phi_{\muGamma, \fX \circ \pGamma}^{\AGammap}\right)^{-1}.
	\]
	Similarly, the pasting diagram
	\[
	\begin{tikzcd}
		\overline{F}(\XGammap) \ar[r, bend left = 30, ""{name = U}]{}{\overline{F}(\muGammap) \circ \overline{F}(\pGammap)} \ar[r, bend right = 30, swap, ""{name = M}]{}{\overline{F}(\opiGammap)} \ar[rr, bend right = 60, swap, ""{name = L}]{}{\overline{F}(\pGammap \circ \muGammap \circ \fGX)} & \overline{F}(\quot{(G \times X)}{\Gamma^{\prime}}) \ar[r]{}{\overline{F}(\fGX)} & \overline{F}(\quot{(G \times X)}{\Gamma}) \ar[from = U, to = M, Rightarrow, shorten <= 4pt, shorten >= 4pt]{}{} \ar[from = 1-2, to = L, Rightarrow, shorten <= 4pt, shorten >= 4pt]{}{}
	\end{tikzcd}
	\]
	is equivalent to the pasting diagram:
	\[
	\begin{tikzcd}
		\overline{F}(\XGammap) \ar[r]{}{\overline{F}(\pGammap)} \ar[rr, bend right = 70, swap, ""{name = L}]{}{\overline{F}(\pGammap \circ \muGammap \circ \fGX)} & \overline{F}(\quot{X}{\Gamma_c^{\prime}}) \ar[r, bend left = 30, ""{name = U}]{}{\overline{F}(\fGX) \circ \overline{F}(\muGammap)} \ar[r, swap, bend right = 30, ""{name = M}]{}{\overline{F}(\muGammap \circ \fGX)} &  \overline{F}(\quot{(G \times X)}{\Gamma}) \ar[from = U, to = M, Rightarrow, shorten <= 4pt, shorten >= 4pt]{}{} \ar[from = 1-2, to = L, Rightarrow, shorten <= 4pt, shorten >= 4pt]{}{}
	\end{tikzcd}
	\]
	It follows that the $2$-cells are given by the corresponding compositors. This allows us to deduce that
	\[
	\phi_{\fGX, \pGammap \circ \muGammap}^{\AGammap} \circ \overline{F}(\fGX)\phi_{\muGammap, \pGammap}^{\AGammap} = \phi_{\muGammap \circ \fGX, \pGammap}^{\AGammap} \circ \phi_{\fGX, \muGammap}^{\overline{F}(\pGammap)\AGammap}.
	\]
	This allows us to compute that
	\begin{align*}
		&\left(\phi_{\muGamma, \pGamma}^{\overline{F}(\fX)\AGammap}\right)^{-1} \circ \left(\phi_{\pGamma \circ \muGamma, \fX}^{\AGammap}\right)^{-1} \circ \phi_{\fGX, \pGammap \circ \muGammap}^{\AGammap} \circ \overline{F}(\fGX)\left(\phi_{\muGammap,p_{\Gamma^{\prime}}}^{\AGammap}\right) \\
		&\circ  \left(\phi_{\fGX, \muGammap}^{\overline{F}(\pGammap)\AGammap}\right)^{-1} \\
		&= \overline{F}(\muGamma)\left(\phi_{\pGamma, \fX}^{\AGammap}\right)^{-1} \circ \left(\phi_{\muGamma, \fX \circ \pGamma}^{\AGammap}\right)^{-1} \circ \phi_{\muGammap \circ \fGX, \pGammap}^{\AGammap} \circ \phi_{\fGX, \muGammap}^{\overline{F}(\pGammap)\AGammap} \\
		&\circ \left(\phi_{\fGX, \muGammap}^{\overline{F}(\pGammap)\AGammap}\right)^{-1} \\
		&=\overline{F}(\muGamma)\left(\phi_{\pGamma, \fX}^{\AGammap}\right)^{-1} \circ \left(\phi_{\muGamma, \fX \circ \pGamma}^{\AGammap}\right)^{-1} \circ \phi_{\muGammap \circ \fGX, \pGammap}^{\AGammap}.
	\end{align*}
	We now use the same tricks to simplify the composition
	\[
	\overline{F}(\muGamma)\left(\phi_{\pGamma, \fX}^{\AGammap}\right)^{-1} \circ \left(\phi_{\muGamma, \fX \circ \pGamma}^{\AGammap}\right)^{-1} \circ \phi_{\muGammap \circ \fGX, \pGammap}^{\AGammap} \circ  \phi_{\muGamma, \of_c}^{\overline{F}(\pGammap)\AGammap}.
	\]
	First note that $\muGammap \circ \fGX = \of_c \circ \muGamma$ so
	\[
	\phi_{\muGammap \circ \fGX, \pGammap}^{\AGammap} \circ  \phi_{\muGamma, \of_c}^{\overline{F}(\pGammap)\AGammap} = \phi_{\of_c \circ \muGamma, \pGammap}^{\AGammap} \circ  \phi_{\muGamma, \of_c}^{\overline{F}(\pGammap)\AGammap}.
	\]
	Now because the pasting diagram
	\[
	\begin{tikzcd}
		\overline{F}(\XGammap) \ar[r]{}{\overline{F}(\pGammap)} \ar[rr, bend right = 70, swap, ""{name = L}]{}{\overline{F}(\pGammap \circ \of_c \circ \muGamma)} & \overline{F}(\quot{X}{\Gamma_c^{\prime}}) \ar[r, bend left = 30, ""{name = U}]{}{\overline{F}(\muGamma) \circ \overline{F}(\of_c)} \ar[r, swap, bend right = 30, ""{name = M}]{}{\overline{F}(\of_c \circ \muGamma)} &  \overline{F}(\quot{(G \times X)}{\Gamma}) \ar[from = U, to = M, Rightarrow, shorten <= 4pt, shorten >= 4pt]{}{} \ar[from = 1-2, to = L, Rightarrow, shorten <= 4pt, shorten >= 4pt]{}{}
	\end{tikzcd}
	\]
	is equivalent to
	\[
	\begin{tikzcd}
		\overline{F}(\XGammap) \ar[r, bend left = 30, ""{name = U}]{}{\overline{F}(\of_c) \circ \overline{F}(\pGammap)} \ar[r, bend right = 30, swap, ""{name = M}]{}{\overline{F}(\pGammap \circ \of_c)} \ar[rr, bend right = 60, swap, ""{name = L}]{}{\overline{F}(\pGammap \circ \of_c \circ \muGamma)} & \overline{F}(\quot{(G \times X)}{\Gamma^{\prime}}) \ar[r]{}{\overline{F}(\muGamma)} & \overline{F}(\quot{(G \times X)}{\Gamma}) \ar[from = U, to = M, Rightarrow, shorten <= 4pt, shorten >= 4pt]{}{} \ar[from = 1-2, to = L, Rightarrow, shorten <= 4pt, shorten >= 4pt]{}{}
	\end{tikzcd}
	\]
	we get that
	\begin{align*}
		\phi_{\muGammap \circ \fGX, \pGammap}^{\AGammap} \circ  \phi_{\muGamma, \of_c}^{\overline{F}(\pGammap)\AGammap} &= \phi_{\of_c \circ \muGamma, \pGammap}^{\AGammap} \circ  \phi_{\muGamma, \of_c}^{\overline{F}(\pGammap)\AGammap} \\
		&=\phi_{\muGamma, \pGammap \circ \of_c}^{\AGammap} \circ \overline{F}(\muGamma)\phi_{\of_c, \pGammap} \\
		&=\phi_{\muGamma, \fX \circ \pGamma}^{\AGammap} \circ \overline{F}(\muGamma)\phi_{\of_c, \pGammap}^{\AGammap}.
	\end{align*}
	Using this we derive that
	\begin{align*}
		&\overline{F}(\muGamma)\left(\phi_{\pGamma, \fX}^{\AGammap}\right)^{-1} \circ \left(\phi_{\muGamma, \fX \circ \pGamma}^{\AGammap}\right)^{-1} \circ \phi_{\muGammap \circ \fGX, \pGammap}^{\AGammap} \circ  \phi_{\muGamma, \of_c}^{\overline{F}(\pGammap)\AGammap} \\
		&= \overline{F}(\muGamma)\left(\phi_{\pGamma, \fX}^{\AGammap}\right)^{-1} \circ \left(\phi_{\muGamma, \fX \circ \pGamma}^{\AGammap}\right)^{-1} \circ \phi_{\muGamma, \fX \circ \pGamma}^{\AGammap} \circ \overline{F}(\muGamma)\phi_{\of_c, \pGammap}^{\AGammap} \\
		&= \overline{F}(\muGamma)\left(\phi_{\pGamma, \fX}^{\AGammap}\right)^{-1} \circ \overline{F}(\muGamma)\phi_{\of_c, \pGammap}^{\AGammap}.
	\end{align*}
	In summary, we have shown so far that
	\begin{align*}
		&\left(\phi_{\muGamma,\pGamma}^{\overline{F}(\fX)\AGammap}\right)^{-1} \circ \left(\phi_{\quot{\pi_2}{\Gamma},\fX}^{\AGammap}\right)^{-1} \circ \phi_{\fGX, \quot{\pi_2}{\Gamma^{\prime}}}^{\AGammap} 
		\circ \overline{F}(\fGX)\phi_{\muGammap, \pGammap}^{\AGammap} \\
		&\circ \left(\phi_{\fGX, \muGammap}^{\overline{F}(\pGammap)\AGammap}\right)^{-1} \circ  \phi_{\muGamma, \of_c}^{\overline{F}(\pGammap)\AGammap} \\
		&= \overline{F}(\muGamma)\left(\phi_{\pGamma, \fX}^{\AGammap}\right)^{-1} \circ \overline{F}(\muGamma)\phi_{\of_c, \pGammap}^{\AGammap}.
	\end{align*}
	Altogether this allows us to deduce that
	\begin{align*}
		&\tau^{\pi_2^{\ast}A} \circ \overline{F}(\fGX)\quot{\theta}{\Gamma^{\prime}} \\
		&= \phi_{\muGamma, \pGamma}^{\AGamma} \circ \big(\overline{F}(\muGamma) \circ \overline{F}(\pGamma)\big)\tau_f^A \circ \left(\phi_{\muGamma,\pGamma}^{\overline{F}(\fX)\AGammap}\right)^{-1} \circ \left(\phi_{\quot{\pi_2}{\Gamma},\fX}^{\AGammap}\right)^{-1} \circ \phi_{\fGX, \quot{\pi_2}{\Gamma^{\prime}}}^{\AGammap} \\
		&\circ \overline{F}(\fGX)\phi_{\muGammap, \pGammap}^{\AGammap} \circ \left(\phi_{\fGX, \muGammap}^{\overline{F}(\pGammap)\AGammap}\right)^{-1} \circ  \phi_{\muGamma, \of_c}^{\overline{F}(\pGammap)\AGammap} \\
		&\circ \big(\overline{F}(\muGamma) \circ \overline{F}(\of_c)\big)\left(\left(\tau_{p_{\Gamma^{\prime}}}^{A}\right)^{-1} \circ\tau_{a_{\Gamma^{\prime}}}^{A}\right) \circ \left(\phi_{\muGamma, \of_c}^{\overline{F}(\aGammap)\AGammap}\right)^{-1} \circ \phi_{\fGX, \muGammap}^{\overline{F}(\aGammap)\AGammap} \\
		&\circ  \overline{F}(\fGX)\left(\phi_{\muGammap, a_{\Gamma^{\prime}}}^{\AGammap}\right)^{-1} \\
		&= \phi_{\muGamma, \pGamma}^{\AGamma} \circ \big(\overline{F}(\muGamma) \circ \overline{F}(\pGamma)\big)\tau_f^A \circ \overline{F}(\muGamma)\left(\phi_{\pGamma, \fX}^{\AGammap}\right)^{-1} \circ \overline{F}(\muGamma)\phi_{\of_c, \pGammap}^{\AGammap} \\
		&\circ \big(\overline{F}(\muGamma) \circ \overline{F}(\of_c)\big)\left(\left(\tau_{p_{\Gamma^{\prime}}}^{A}\right)^{-1} \circ\tau_{a_{\Gamma^{\prime}}}^{A}\right) \circ \left(\phi_{\muGamma, \of_c}^{\overline{F}(\aGammap)\AGammap}\right)^{-1} \circ \phi_{\fGX, \muGammap}^{\overline{F}(\aGammap)\AGammap} \\
		&\circ  \overline{F}(\fGX)\left(\phi_{\muGammap, a_{\Gamma^{\prime}}}^{\AGammap}\right)^{-1}.
	\end{align*}
	
	To proceed, we first observe that
	\begin{align*}
		&\left(\phi_{\muGamma, \pGamma}^{\AGamma}\right)^{-1} \circ \tau_{f}^{\pi_2^{\ast}A} \circ \overline{F}(\fGX)\quot{\theta}{\Gamma^{\prime}} \circ  \overline{F}(\fGX)\left(\phi_{\muGammap, a_{\Gamma^{\prime}}}^{\AGammap}\right) \circ \left(\phi_{\fGX, \muGammap}^{\overline{F}(\aGammap)\AGammap}\right)^{-1} \\
		&\circ \left(\phi_{\muGamma, \of_c}^{\overline{F}(\aGammap)\AGammap}\right) \\
		&=\big(\overline{F}(\muGamma) \circ \overline{F}(\pGamma)\big)\tau_f^A \circ \overline{F}(\muGamma)\left(\phi_{\pGamma, \fX}^{\AGammap}\right)^{-1} \circ \overline{F}(\muGamma)\phi_{\of_c, \pGammap}^{\AGammap} \\
		&\circ \big(\overline{F}(\muGamma) \circ \overline{F}(\of_c)\big)\left(\left(\tau_{p_{\Gamma^{\prime}}}^{A}\right)^{-1} \circ\tau_{a_{\Gamma^{\prime}}}^{A}\right)
	\end{align*} 
	We will now study what happens when we post-compose this expression by $\overline{F}(\muGamma)\tau_{\pGamma}^A$. Explicitly,
	\begin{align*}
		&\overline{F}(\muGamma)\tau_{\pGamma}^{A} \circ \big(\overline{F}(\muGamma) \circ \overline{F}(\pGamma)\big)\tau_f^A \circ \overline{F}(\muGamma)\left(\phi_{\pGamma, \fX}^{\AGammap}\right)^{-1}\circ \overline{F}(\muGamma)\phi_{\of_c, \pGammap}^{\AGammap} \\
		& \circ \big(\overline{F}(\muGamma) \circ \overline{F}(\of_c)\big)\left(\left(\tau_{p_{\Gamma^{\prime}}}^{A}\right)^{-1} \circ\tau_{a_{\Gamma^{\prime}}}^{A}\right) \\
		&= \overline{F}(\muGamma)\left(\tau_{\pGamma}^{A} \circ \overline{F}(\pGamma)\tau_f^A \circ \left(\phi_{\pGamma, \fX}^{\AGammap}\right)^{-1} \circ \cdots \circ \overline{F}(\of_c)\left(\left(\tau_{p_{\Gamma^{\prime}}}^{A}\right)^{-1} \circ\tau_{a_{\Gamma^{\prime}}}^{A}\right)\right) \\
		&= \overline{F}(\muGamma)\left(\tau_{\fX \circ \pGamma}^{A} \circ \phi_{\pGamma, \fX}^{\AGammap} \circ \left(\phi_{\pGamma, \fX}^{\AGammap}\right)^{-1} \circ \cdots \circ \overline{F}(\of_c)\left(\left(\tau_{p_{\Gamma^{\prime}}}^{A}\right)^{-1} \circ\tau_{a_{\Gamma^{\prime}}}^{A}\right)\right) \\
		&=\overline{F}(\muGamma)\left(\tau_{\pGammap \circ \of_c}^{A} \circ \phi_{\pGamma, \of_c}^{\AGammap} \circ  \overline{F}(\of_c)\left(\left(\tau_{p_{\Gamma^{\prime}}}^{A}\right)^{-1} \circ\tau_{a_{\Gamma^{\prime}}}^{A}\right)\right) \\
		&=\overline{F}(\muGamma)\left(\tau_{\of_c}^{A} \circ \overline{F}(\of_c)\tau_{\pGammap}^A \circ \left(\phi_{\pGammap, \of_c}^{\AGammap}\right)^{-1} \circ \phi_{\pGamma, \of_c}^{\AGammap} \overline{F}(\of_c)\left(\left(\tau_{p_{\Gamma^{\prime}}}^{A}\right)^{-1} \circ\tau_{a_{\Gamma^{\prime}}}^{A}\right)\right) \\
		&=\overline{F}(\muGamma)\left(\tau_{\of_c}^{A} \circ \overline{F}(\of_c)\tau_{\pGammap}^A \circ  \overline{F}(\of_c)\left(\tau_{p_{\Gamma^{\prime}}}^{A}\right)^{-1} \circ \overline{F}(\of_c)\tau_{a_{\Gamma^{\prime}}}^{A}\right) \\
		&= \overline{F}(\muGamma)\left(\tau_{f_c}^{A} \circ \overline{F}(\of_c)\tau_{a_{\Gamma^{\prime}}}^{A}\right) = \overline{F}(\muGamma)\left(\tau_{\aGammap \circ \of_c}^{A} \circ \phi_{f_c, \aGammap}^{\AGammap}\right) \\
		&= \overline{F}(\muGamma)\left(\tau_{\fX \circ \aGamma}^{A} \circ \phi_{f_c, \aGammap}^{\AGammap}\right) \\
		&=\overline{F}(\muGamma)\left(\tau_{\aGamma}^{A} \circ \overline{F}(\aGamma)\tau_{f}^{A} \circ \left(\phi_{\aGamma, \fX}^{\AGammap}\right)^{-1}\circ\phi_{f_c, \aGammap}^{\AGammap}\right).
	\end{align*}
	This implies in turn that
	\begin{align*}
		&\overline{F}(\muGamma)\tau_{\pGamma}^{A} \circ \left(\phi_{\muGamma, \pGamma}^{\AGamma}\right)^{-1} \circ \tau_{f}^{\pi_2^{\ast}A} \circ \overline{F}(\fGX)\quot{\theta}{\Gamma^{\prime}} \circ  \overline{F}(\fGX)\left(\phi_{\muGammap, a_{\Gamma^{\prime}}}^{\AGammap}\right) \\
		&\circ \left(\phi_{\fGX, \muGammap}^{\overline{F}(\aGammap)\AGammap}\right)^{-1} \\ 
		&=\overline{F}(\muGamma)\left(\tau_{\aGamma}^{A} \circ \overline{F}(\aGamma)\tau_{f}^{A} \circ \left(\phi_{\aGamma, \fX}^{\AGammap}\right)^{-1}\circ\phi_{f_c, \aGammap}^{\AGammap}\right)
	\end{align*}
	which allows us to deduce that
	\begin{align*}
		&\tau_f^{\pi_2^{\ast}A} \circ \overline{F}(\fGX)\quot{\theta}{\Gamma^{\prime}}\\
		&= \phi_{\muGamma, \pGamma}^{\AGamma} \circ \overline{F}(\muGamma)\left(\left(\tau_{\pGamma}^{A}\right)^{-1} \circ \tau_{\aGamma}^{A}\right) \circ \big(\overline{F}(\muGamma) \circ \overline{F}(\aGamma)\big)\tau_f^A \circ \overline{F}(\muGamma)\left(\phi_{\aGamma, \fX}^{\AGammap}\right)^{-1}\\
		&\circ \overline{F}(\muGamma)\phi_{f_c, \aGammap}^{\AGammap} \circ \left(\phi_{\muGamma, \of_c}^{\overline{F}(\aGammap)\AGammap}\right)^{-1} \circ \phi_{\fGX, \muGammap}^{\overline{F}(\aGammap)\AGammap} \circ  \overline{F}(\fGX)\left(\phi_{\muGammap, a_{\Gamma^{\prime}}}^{\AGammap}\right)^{-1}.
	\end{align*}
	To complete the proof that $\Theta$ is an $F_G(X)$-morphism, it suffices at this point to show that the equality
	\begin{align*}
		&\big(\overline{F}(\muGamma) \circ \overline{F}(\aGamma)\big)\tau_f^A \circ \overline{F}(\muGamma)\left(\phi_{\aGamma, \fX}^{\AGammap}\right)^{-1} \circ \overline{F}(\muGamma)\phi_{f_c, \aGammap}^{\AGammap}  \\
		&\circ \left(\phi_{\muGamma, \of_c}^{\overline{F}(\aGammap)\AGammap}\right)^{-1} \circ \phi_{\fGX, \muGammap}^{\overline{F}(\aGammap)\AGammap} \circ  \overline{F}(\fGX)\left(\phi_{\muGammap, a_{\Gamma^{\prime}}}^{\AGammap}\right)^{-1} \\
		&= \left(\phi_{\muGamma, \aGamma}^{\AGamma}\right)^{-1} \circ \overline{F}(\oalphaGamma)\tau_f^A \circ \left(\phi_{\oalphaGamma, \fX}^{\AGammap}\right)^{-1} \circ \phi_{\fGX, \oalphaGammap}^{\AGammap}
	\end{align*}
	holds. However, this follows from an extremely tedious manipulation of the naturality of the compositors and the relations on the compositors as induced by the pseudofunctoriality of $\overline{F}$ and Lemma \ref{Lemma: Clifton's Naivification Variety}. From this it follows that
	\begin{align*}
		&\tau_f^{\pi_2^{\ast}A} \circ \overline{F}(\fGX)\quot{\theta}{\Gamma^{\prime}}\\
		&= \phi_{\muGamma, \pGamma}^{\AGamma} \circ \overline{F}(\muGamma)\left(\left(\tau_{\pGamma}^{A}\right)^{-1} \circ \tau_{\aGamma}^{A}\right) \circ \big(\overline{F}(\muGamma) \circ \overline{F}(\aGamma)\big)\tau_f^A \circ \overline{F}(\muGamma)\left(\phi_{\aGamma, \fX}^{\AGammap}\right)^{-1}\\
		&\circ \overline{F}(\muGamma)\phi_{f_c, \aGammap}^{\AGammap} \circ \left(\phi_{\muGamma, \of_c}^{\overline{F}(\aGammap)\AGammap}\right)^{-1} \circ \phi_{\fGX, \muGammap}^{\overline{F}(\aGammap)\AGammap} \circ  \overline{F}(\fGX)\left(\phi_{\muGammap, a_{\Gamma^{\prime}}}^{\AGammap}\right)^{-1} \\
		&=\phi_{\muGamma, \pGamma}^{\AGamma} \circ \overline{F}(\muGamma)\left(\left(\tau_{\pGamma}^{A}\right)^{-1} \circ \tau_{\aGamma}^{A}\right) \circ \left(\phi_{\muGamma, \aGamma}^{\AGamma}\right)^{-1} \circ \overline{F}(\oalphaGamma)\tau_f^A \circ \left(\phi_{\oalphaGamma, \fX}^{\AGammap}\right)^{-1} \\
		&\circ \phi_{\fGX, \oalphaGammap}^{\AGammap} \\
		&= \quot{\theta}{\Gamma} \circ \tau_f^{\alpha_X^{\ast}A}.
	\end{align*}
	This proves, finally, that $\tau_f^{\pi_2^{\ast}A} \circ \overline{F}(\fGX)\quot{\theta}{\Gamma^{\prime}} = \quot{\theta}{\Gamma} \circ \tau_f^{\alpha_X^{\ast}A}$ and hence shows that $\Theta \in F_G(X)_1$. That $\Theta$ is an isomorphism is immediate from each $\quot{\theta}{\Gamma}$ being an isomorphism, which proves that there is an isomorphism
	\[
	\alpha_X^{\ast}A \xrightarrow[\Theta]{\cong} \pi_2^{\ast}A
	\]
	in $F_G(G \times X)$ for any $A \in F_G(X)_0$.
	
	Having finally proved that $\Theta$ is an $F_G(X)$-isomorphism, we now show that the $\Theta$ are natural in $A$. For this fix a $P:A \to B$ in $F_G(X)$. Fix a $\Gamma \in \Sf(G)_0$ and note that we must check that the diagram
	\[
	\xymatrix{
		\alpha_X^{\ast}A \ar[d]_{\alpha_X^{\ast}P} \ar[r]^-{\Theta_A} & \pi_2^{\ast}A \ar[d]^{\pi_2^{\ast}P} \\
		\alpha_X^{\ast}B \ar[r]_-{\Theta_B} & \pi_2^{\ast}B
	}
	\]
	commutes. We calculate
	\begin{align*}
		&\overline{F}(\opiGamma)\rhoGamma \circ \quot{\theta_A}{\Gamma} \\
		&= \overline{F}(\opiGamma)\rhoGamma \circ \phi_{\muGamma, p_{\Gamma}}^{\AGamma} \circ \overline{F}(\muGamma)\left(\left(\tau_{p_{\Gamma}}^{A}\right)^{-1} \circ \tau_{a_{\Gamma}}^{A}\right) \circ \left(\phi_{\muGamma, a_{\Gamma}}^{\AGamma}\right)^{-1} \\
		&= \phi_{\muGamma, p_{\Gamma}}^{\BGamma} \circ \overline{F}(\muGamma)\left(\overline{F}(\pGamma)\rhoGamma \circ \left(\tau_{p_{\Gamma}}^{A}\right)^{-1} \circ \tau_{a_{\Gamma}}^{A}\right) \circ \left(\phi_{\muGamma, a_{\Gamma}}^{\AGamma}\right)^{-1}  \\
		&=\phi_{\muGamma, p_{\Gamma}}^{\BGamma} \circ \overline{F}(\muGamma)\left(\left(\tau_{p_{\Gamma}}^{B}\right)^{-1} \circ \rhoGamma \circ \tau_{a_{\Gamma}}^{A}\right) \circ \left(\phi_{\muGamma, a_{\Gamma}}^{\AGamma}\right)^{-1}  \\
		&= \phi_{\muGamma, p_{\Gamma}}^{\BGamma} \circ \overline{F}(\muGamma)\left(\left(\tau_{p_{\Gamma}}^{B}\right)^{-1} \circ \tau_{a_{\Gamma}}^{B}\right) \circ \overline{F}(\muGamma)\left(\overline{F}(\aGamma)\rhoGamma\right) \circ \left(\phi_{\muGamma, a_{\Gamma}}^{\AGamma}\right)^{-1} \\
		&= \phi_{\muGamma, p_{\Gamma}}^{\BGamma} \circ \overline{F}(\muGamma)\left(\left(\tau_{p_{\Gamma}}^{B}\right)^{-1} \circ \tau_{a_{\Gamma}}^{B}\right)  \circ \left(\phi_{\muGamma, a_{\Gamma}}^{\BGamma}\right)^{-1}\circ \overline{F}(\aGamma \circ \muGamma)\rhoGamma \\
		&= \phi_{\muGamma, p_{\Gamma}}^{\BGamma} \circ \overline{F}(\muGamma)\left(\left(\tau_{p_{\Gamma}}^{B}\right)^{-1} \circ \tau_{a_{\Gamma}}^{B}\right)  \circ \left(\phi_{\muGamma, a_{\Gamma}}^{\BGamma}\right)^{-1}\circ \overline{F}(\oalphaGamma)\rhoGamma \\
		&= \quot{\theta_B}{\Gamma} \circ \overline{F}(\oalphaGamma)\rhoGamma,
	\end{align*}
	which verifies the commutativity of the diagram.
\end{proof}
\begin{corollary}
	For any left $G$-variety $X$ or any left $G$-space $X$ and any pre-equivariant pseudofunctor $F$ on $X$, there is a natural isomorphism of functors $\pi_2^{\ast} \cong \alpha_X^{\ast}:F_G(X) \to F_G(G \times X)$.
\end{corollary}
Before proving the cocycle condition that $\Theta$ satisfies, we actually define it so that we know what we must \emph{actually} prove.
\begin{definition}\label{Defn: GIT Cocycle}\index[terminology]{GIT Cocycle Condition}
	Let $F$ be a pre-equivariant pseudofunctor on $X$, let $A \in F_G(X)_0$, and let $\Theta:\alpha_X^{\ast}A \to \pi_2^{\ast}A$ be an isomorphism in $F_G(G \times X)$. We say that $\Theta$ satisfies the GIT-cocycle condition if, for morphisms
	\begin{align*}
		d_0^1 &:= \pi_{23}:G \times G \times X \to G \times X; \\
		d_1^1 &:= \mu_G \times \id_X:G \times G \times X \to G \times X; \\
		d_2^1 &:= \id_G \times \alpha_X:G \times G \times X \to G \times X
	\end{align*}\index[notation]{dFaceMaps@$d_0^1, d_1^1, d_2^1$}
	for all $\Gamma \in \Sf(G)_0$ the composite
	\begin{align*}
		\phi_{\quot{d_1^1}{\Gamma}, \opiGamma}^{\AGamma} \circ F(\quot{d_1^1}{\Gamma})\quot{\theta}{\Gamma} \circ \left(\phi_{\quot{d_1^1}{\Gamma}, \oalphaGamma}^{\AGamma}\right)
	\end{align*}
	is equal to the composite
	\begin{align*}
		&\phi_{\quot{d_0^1}{\Gamma}, \opiGamma}^{\AGamma} \circ F(\quot{d_0^1}{\Gamma})\quot{\theta}{\Gamma} \circ \left(\phi_{\quot{d_0^1}{\Gamma}, \oalphaGamma}^{\AGamma}\right)^{-1} \circ \phi_{\quot{d_2^1}{\Gamma},\opiGamma}^{\AGamma} \circ F(\quot{d_2^1}{\Gamma})\quot{\theta}{\Gamma} \circ \left(\phi_{\quot{d_2^1}{\Gamma},\oalphaGamma}^{\AGamma}\right)^{-1}.
	\end{align*}
\end{definition}

Let us now motivate the GIT cocycle condition by studying how it arose. In \cite{GIT} the authors described an equivariance isomorphism on a sheaf $\Fscr$ on $X$ as an isomorphism
\[
\theta:\alpha_X^{\ast}\Fscr \xrightarrow{\cong} \pi_2^{\ast} \Fscr
\]
in $\Shv(G \times X)$. This isomorphism $\theta$ was used to illustrate the how the action of $G$ on $X$ affects the space $X$ and the sheaf $\Fscr$, as well as how these conditions need to be suitably compatible in order to define an equivariant sheaf. Let us describe this more explicitly. The action $\alpha_X$ needs to intertwine the multiplication of the group $G$ on itself together with the action on $X$ in the sense that the diagram
\[
\xymatrix{
	G \times G \times X \ar[rr]^-{\id_G \times \alpha_X} \ar[d]_{\mu_G \times \id_X} & & G \times X \ar[d]^{\alpha_X} \\
	G \times X \ar[rr]_-{\alpha_X} & & X
}
\]
commutes. When looking at this at the level of stalks of $\Fscr$, this requires an isomorphism
\[
\Fscr_x \cong \Fscr_{(gh)x} \cong \Fscr_{g(hx)}
\]
for any stalk $(g,h,x) \in G \times G \times X$. However, these isomorphisms should come from the equivariance $\theta$ lifted to $G \times G \times X$ by either changing the order of multiplication,  action, and projection more-or-less freely (provided we do so in a way that makes sense). There are three maps from $G \times G \times X$ to $G \times X$ that we use to lift $\theta$; these three maps coincide with the simplicial presentation of both the quotient variety $\underline{G\backslash X}_{\bullet}$ and the simplicial information encoded by the internal translation groupoid $\mathbb{G\times X}$ (cf.\@  \cite[Appendix B]{MyThesis}. These three maps, together with the two maps from $G \times X \to X$ in the diagram
\[
\begin{tikzcd}
	G \times G \times X \ar[rr, shift left = 3]{}{\pi_{23}} \ar[rr]{}[description]{\mu_G\times \id_X} \ar[rr, shift right = 3, swap]{}{\id_G \times \alpha_X} & & G \times X \ar[rr, shift left]{}{\pi_2} \ar[rr, shift right, swap]{}{\alpha_X} & & X
\end{tikzcd}
\]
satisfy the equations
\begin{align*}
	\alpha_X \circ (\mu_G\times \id_X) &= \alpha_X \circ (\id_G \times \alpha_X); \\
	\alpha_X \circ \pi_{23} &= \pi_2 \circ (\id_G \times \alpha_X); \\
	\pi_2 \circ \pi_{23} &= \pi_2 \circ (\mu_G \times \id_X).
\end{align*}
Define, in order to save space and illustrate the simplicial connection, the morphisms
\begin{align*}
	\pi_{23} &:= d_0^1; \\
	\mu_G \times \id_X &:= d_1^1; \\
	\id_G \times \alpha_X &:= d_2^1.
\end{align*}
Using these maps and relations we then can produce a diagram (assuming that the underlying pseudofunctor $\Shv(-)$ is strict in this case)
\[
\begin{tikzcd}
	(d_2^1)^{\ast}\big(\alpha_X^{\ast}\Fscr\big) \ar[r]{}{(d_2^1)^{\ast}\theta} \ar[d, equals] & (d_2^1)^{\ast}\big(\pi_{2}^{\ast}\Fscr\big) \ar[r, equals] & \left( \pi_2 \circ d_2^1\right)^{\ast}\Fscr \ar[d, equals] \\
	(\alpha_X \circ d_2^1)^{\ast}\Fscr \ar[d, equals] & & (\alpha_X \circ d_0^1)^{\ast}\Fscr \ar[d, equals] \\
	(\alpha_X \circ d_1^1)^{\ast}\Fscr \ar[d, equals] & & (d_0^1)^{\ast}\big(\alpha_X^{\ast}\Fscr\big) \ar[d]{}{(d_0^1)^{\ast}\theta} \\
	(d_1^1)^{\ast}\alpha_X^{\ast}\Fscr \ar[dd, swap]{}{(d_1^1)^{\ast}\theta} &  & (d_0^1)^{\ast}\big(\pi_2^{\ast}\Fscr\big) \ar[d, equals] \\
	&  & (\pi_2 \circ d_0^1)^{\ast}\Fscr \ar[d, equals] \\
	(d_1^1)^{\ast}\big(\pi_2^{\ast}\Fscr\big)& & \ar[ll, equals] (\pi_2 \circ d_1^1)^{\ast}\Fscr
\end{tikzcd}
\]
The cocycle condition of \cite{GIT} asks exactly that the above diagram commute, i.e.,
\[
(d_1^1)^{\ast}\theta = (d_0^1)^{\ast}\theta \circ (d_2^1)^{\ast}\theta.
\] 
\begin{remark}
	In the case when $\Fscr$ is a sheaf on a topological space $X$, the GIT cocycle condition is equivalent to asking that if $\operatorname{esp}(\Fscr)$ is the espace {\'e}tal{\'e} of $\Fscr$, then $\operatorname{esp}(\Fscr)$ carries a $G$-action for which the {\'e}tale map of spaces $p:\operatorname{esp}(\Fscr) \to X$ is $G$-equivariant.
\end{remark}
Our pseudofunctorial GIT-cocycle condition then asks for a pre-equivariant pseudofunctor $F$ over $X$ that for any $A \in F_G(X)_0$ and $\Gamma \in \Sf(G)_0$, the diagram
\[
\begin{tikzcd}
	F(\oalphaGamma \circ \odGamma{2}{1})(\AGamma) \ar[rr]{}{\left(\phi_{\odGamma{2}{1}, \oalphaGamma}^{\AGamma}\right)^{-1}} \ar[d, equals] & & F(\odGamma{2}{1})\big(F(\oalphaGamma)(\AGamma)\big) \ar[d]{}{F(\odGamma{2}{1})\quot{\theta}{\Gamma}}  \\
	F(\oalphaGamma \circ \odGamma{1}{1})(\AGamma) \ar[d, swap]{}{\left(\phi_{\odGamma{1}{1}, \oalphaGamma}^{\AGamma}\right)^{-1}} & & F(\odGamma{2}{1})\big(F(\opiGamma)(\AGamma)\big) \ar[d]{}{\phi_{\odGamma{2}{1}, \opiGamma}^{\AGamma}} \\
	F(\odGamma{1}{1})\big(F(\oalphaGamma)(\AGamma)\big) \ar[d, swap]{}{F(\odGamma{1}{1})\quot{\theta}{\Gamma}} & & F(\opiGamma \circ \odGamma{2}{1}) \ar[d, equals] \\
	F(\odGamma{1}{1})\big(F(\opiGamma)(\AGamma)\big) \ar[d, swap]{}{\phi_{\odGamma{1}{1}, \opiGamma}^{\AGamma}} & &  F(\oalphaGamma \circ \odGamma{0}{1})(\AGamma) \ar[d]{}{\left(\phi_{\odGamma{0}{1}, \oalphaGamma}^{\AGamma}\right)^{-1}} \\
	F(\opiGamma \circ \odGamma{1}{1})(\AGamma) & & F(\odGamma{0}{1})\big(F(\oalphaGamma)(\AGamma)\big) \ar[d]{}{F(\odGamma{0}{1})\quot{\theta}{\Gamma}} \\
	F(\opiGamma \circ \odGamma{0}{1})(\AGamma) \ar[u, equals] & & F(\odGamma{0}{1})\big(F(\opiGamma)(\AGamma)\big) \ar[ll]{}{\phi_{\odGamma{0}{1},\opiGamma}^{\AGamma}}
\end{tikzcd}
\]
commutes. Note that while we have not proved it explicitly, it follows mutatis mutandis to Lemmas \ref{Lemma: Section Change of Space: pi2 pullback} and \ref{Lemma: Section Change of Space: Action map pullback} that $F$ admits pseudocone translations all of $d_0^1,$ $d_1^1,$ and $d_2^1$. This gives rise to induce the corresponding pullback functors
\[
(d_0^1)^{\ast}, (d_1^1)^{\ast}, (d_2^1)^{\ast}:F_G(G \times X) \to F_G(G \times G \times X).
\]

As we move to show that the isomorphism $\Theta$ of Theorem \ref{Theorem: Formal naive equivariance} satisfies the GIT cocycle condtion (cf.\@ Theorem \ref{Thm: GIT Cocycle Condition Formal Equivariance} below), we will need a short lemma on the interaction of the maps $\odGamma{2}{1}$ and $\odGamma{0}{1}$ with $\muGamma$.
\begin{lemma}\label{Lemma: MuGamma maps same with odgammas}
	For any $\Gamma \in \Sf(G)_0$ when $G$ is a smooth algebraic group or $M \in \mathbf{FAMan}(G)_0$ when $G$ is a topological group, if $\muGamma$ is the isomorphism of Lemma \ref{Lemma: Clifton's Naivification Variety} and $\mu_M$ is the isomorphism of Lemma \ref{Lemma: Clifton's Naivification Space}, then
	\[
	\muGamma \circ \odGamma{0}{1} = \muGamma \circ \odGamma{2}{1}
	\]
	and
	\[
	\mu_M \circ \quot{\overline{d}_0^1}{M} = \mu_{M} \circ \quot{\overline{d}_2^1}{M}.
	\]
\end{lemma}
\begin{proof}
	Observe that since
	\[
	\quot{d_0^1}{\Gamma} = \pi_{134}:\Gamma \times G \times G \times X \to \Gamma \times G \times X
	\]
	and
	\[
	\quot{d_2^1}{\Gamma} = \id_{\Gamma \times G} \times \alpha_{X}:\Gamma \times G \times G \times X \to \Gamma \times G \times X
	\]
	there are unique morphisms making the diagrams
	\[
	\xymatrix{
		\Gamma \times G \times G \times X \ar@<.5ex>[r]^-{\quot{d_2^1}{\Gamma}} \ar@<-.5ex>[r]_-{\quot{d_0^1}{\Gamma}} \ar[d]_{\quo_{\Gamma\times G \times G}} & \Gamma \times G \times X \ar[d]^{\quo_{G \times \Gamma}} \\
		\quot{(G \times G \times X)}{\Gamma} \ar@<.5ex>[r]^-{\odGamma{2}{1}} \ar@<-.5ex>[r]_-{\odGamma{0}{1}} & \quot{(G \times X)}{\Gamma}
	}
	\]
	\[
	\xymatrix{
		\Gamma \times G \times G \times X \ar@<.5ex>[r]^-{\quot{d_2^1}{\Gamma}} \ar@<-.5ex>[r]_-{\quot{d_0^1}{\Gamma}} \ar[d]_{\quo_{(\Gamma \times G)_{c}}} & \Gamma \times G \times X \ar[d]^{\quo_{\Gamma_c}} \\
		\quot{X}{(\Gamma \times G)_c} \ar@<.5ex>[r]^-{(\odGamma{2}{1})_c} \ar@<-.5ex>[r]_-{(\odGamma{0}{1})_c} & \quot{X}{\Gamma_c}
	}
	\]
	commute. However, since the second diagram factors through the diagrams
	\[
	\xymatrix{
		\Gamma \times G \times G \times X \ar@<.5ex>[r]^-{\quot{d_2^1}{\Gamma}} \ar@<-.5ex>[r]_-{\quot{d_0^1}{\Gamma}} \ar[d]_{\quo_{(\Gamma \times G)_{c}}} & \Gamma \times G \times X \ar[d]^{\quo_{\Gamma_c}} \\
		\quot{(G \times G \times X)}{\Gamma} \ar@<.5ex>[r]^-{\odGamma{2}{1}} \ar@<-.5ex>[r]_-{\odGamma{0}{1}} \ar[d]_{\mu_{\Gamma \times G}} & \quot{(G \times X)}{\Gamma} \ar[d]^{\muGamma} \\
		\quot{X}{(\Gamma \times G)_c} \ar@<.5ex>[r]^-{(\odGamma{2}{1})_c} \ar@<-.5ex>[r]_-{(\odGamma{0}{1})_c} & \quot{X}{\Gamma_c}
	}
	\]
	the result follows.
\end{proof}
\begin{corollary}\label{Cor: F of muG circ d01 is F of muG circ d21}
	For any pre-equivariant pseudofunctor $F$ over $X$,
	\[
	\overline{F}\left(\muGamma \circ \odGamma{0}{1}\right) = \overline{F}\left(\muGamma \circ \odGamma{2}{1}\right)
	\]
	and
	\[
		\overline{F}\left(\mu_M \circ \quot{\overline{d}_0^1}{M}\right) = \overline{F}\left(\mu_{M} \circ \quot{\overline{d}_2^1}{M}\right)
	\]
\end{corollary}
\begin{Theorem}\label{Thm: GIT Cocycle Condition Formal Equivariance}
	Let $F$ be a pre-equivariant pseudofunctor over $X$ and let $A \in F_G(X)_0$. Then the equivariance $\Theta$ of Theorem \ref{Theorem: Formal naive equivariance} satisfies the GIT cocycle condition.
\end{Theorem}
\begin{proof}
As usual, we prove the variety case and omit the topological case, as the topological case follows mutatis mutandis. Begin by noting that from Theorem \ref{Theorem: Formal naive equivariance} we must prove that
	\begin{align*}
		&\phi_{\quot{d_1^1}{\Gamma}, \opiGamma}^{\AGamma} \circ \overline{F}(\quot{d_1^1}{\Gamma})\quot{\theta}{\Gamma} \circ \left(\phi_{\quot{d_1^1}{\Gamma}, \oalphaGamma}^{\AGamma}\right) \\
		&= \phi_{\quot{d_1^1}{\Gamma}, \opiGamma}^{\AGamma} \circ \overline{F}(\quot{d_1^1}{\Gamma})\left(\phi_{\muGamma, \pGamma}^{\AGamma}\circ \overline{F}(\muGamma)\big(\big(\tau_{\pGamma}^{A}\big)^{-1} \circ \tau_{\aGamma}^{A}\big) \circ \left(\phi_{\muGamma, \aGamma}^{\AGamma}\right)^{-1}\right) \\
		&\circ \left(\phi_{\quot{d_1^1}{\Gamma}, \oalphaGamma}^{\AGamma}\right)
	\end{align*}
	is equal to the composite
	\begin{align*}
		&\phi_{\quot{d_0^1}{\Gamma}, \opiGamma}^{\AGamma} \circ \overline{F}(\quot{d_0^1}{\Gamma})\quot{\theta}{\Gamma} \circ \left(\phi_{\quot{d_0^1}{\Gamma}, \oalphaGamma}^{\AGamma}\right)^{-1} \circ \phi_{\quot{d_2^1}{\Gamma},\opiGamma}^{\AGamma} \circ \overline{F}(\quot{d_2^1}{\Gamma})\quot{\theta}{\Gamma} \circ \left(\phi_{\quot{d_2^1}{\Gamma},\oalphaGamma}^{\AGamma}\right)^{-1} \\
		&= \phi_{\quot{d_0^1}{\Gamma}, \opiGamma}^{\AGamma} \circ \overline{F}(\odGamma{2}{1})\phi_{\muGamma, \pGamma}^{\AGamma} \circ \big(\overline{F}(\odGamma{2}{1})\circ \overline{F}(\muGamma)\big)\left(\big(\tau_{\pGamma}^{A}\big)^{-1} \circ \tau_{\aGamma}^A\right) \\
		&\circ \overline{F}(\odGamma{0}{1})\left(\phi_{\muGamma, \aGamma}^{\AGamma}\right)^{-1} \circ \left(\phi_{\quot{d_0^1}{\Gamma}, \oalphaGamma}^{\AGamma}\right)^{-1} \circ \phi_{\quot{d_2^1}{\Gamma},\opiGamma}^{\AGamma} \circ \overline{F}(\odGamma{2}{1})\phi_{\muGamma, \pGamma}^{\AGamma} \\
		&\circ \big(\overline{F}(\odGamma{2}{1}) \circ \overline{F}(\muGamma)\big)\left(\big(\tau_{\pGamma}^{A}\big)^{-1} \circ \tau_{\aGamma}^{A}\right) \circ \overline{F}(\odGamma{2}{1})\left(\phi_{\muGamma, \aGamma}^{\AGamma}\right)^{-1} \circ \left(\phi_{\quot{d_2^1}{\Gamma},\oalphaGamma}^{\AGamma}\right)^{-1}.
	\end{align*}
	We will prove this by systematically attacking and simplifying the second expression. Begin by recalling that either by Lemma \ref{Lemma: Clifton's Naivification Variety} or the exposition prior to the Lemma (respectively Lemma \ref{Lemma: Clifton's Naivification Space}), the identities
	\begin{align*}
		\opiGamma &= \pGamma \circ \muGamma; \\
		\oalphaGamma &= \aGamma \circ \muGamma.
	\end{align*}
	both hold. Now, using these and the identities from before the statment of the theorem, we observe as in the proof of Thoerem \ref{Theorem: Formal naive equivariance} that the pasting diagram
	\[
	\begin{tikzcd}
		\overline{F}(\XGamma) \ar[r, bend left = 30, ""{name = U}]{}{\overline{F}(\muGamma) \circ \overline{F}(\pGammap)} \ar[r, bend right = 30, swap, ""{name = M}]{}{\overline{F}(\opiGamma)} \ar[rr, bend right = 60, swap, ""{name = L}]{}{\overline{F}(\pGamma \circ \muGamma \circ \odGamma{0}{1})} & \overline{F}(\quot{(G \times X)}{\Gamma}) \ar[r]{}{\overline{F}(\odGamma{0}{1})} & \overline{F}(\quot{(G \times G \times X)}{\Gamma}) \ar[from = U, to = M, Rightarrow, shorten <= 4pt, shorten >= 4pt]{}{} \ar[from = 1-2, to = L, Rightarrow, shorten <= 4pt, shorten >= 4pt]{}{}
	\end{tikzcd}
	\]
	is equal to the pasting diagram
	\[
	\begin{tikzcd}
		\overline{F}(\XGamma) \ar[r]{}{\overline{F}(\pGamma)} \ar[rr, bend right = 60, swap, ""{name = L}]{}{\overline{F}(\pGamma \circ \muGamma \circ \odGamma{0}{1})} & \overline{F}(\quot{X}{\Gamma_c}) \ar[r, bend left = 30, ""{name = U}]{}{\overline{F}(\odGamma{0}{1}) \circ \overline{F}(\muGamma)} \ar[r, bend right = 30, swap, ""{name = M}]{}{\overline{F}(\muGamma \circ \odGamma{0}{1})} & \overline{F}(\quot{(G \times G \times X)}{\Gamma}) \ar[from = U, to = M, Rightarrow, shorten <= 4pt, shorten >= 4pt] \ar[from = 1-2, to = L, Rightarrow, shorten <= 4pt, shorten >= 4pt]{}{}
	\end{tikzcd}
	\]
	where the the $2$-cells are again the compositor isomorphisms. This allows us to deduce that
	\[
	\phi_{\quot{d_0^1}{\Gamma}, \opiGamma}^{\AGamma} \circ \overline{F}(\odGamma{0}{1})\phi_{\muGamma, \pGamma}^{\AGamma} = \phi_{\muGamma \circ \odGamma{0}{1}, \pGamma}^{\AGamma} \circ \phi_{\odGamma{0}{1}, \muGamma}^{\overline{F}(\pGamma)\AGamma}.
	\]
	Similarly, using that
	\[
	\overline{F}(\odGamma{0}{1})\left(\phi_{\muGamma, \aGamma}^{\AGamma}\right)^{-1} \circ \left(\phi_{\quot{d_0^1}{\Gamma}, \oalphaGamma}^{\AGamma}\right)^{-1} = \left(\phi_{\odGamma{0}{1}, \oalphaGamma}^{\AGamma} \circ \overline{F}(\odGamma{0}{1})\phi_{\muGamma, \aGamma}^{\AGamma} \circ \right)^{-1}
	\]
	and proceeding to analyze the pasting diagrams as before, we get that
	\[
	\overline{F}(\odGamma{0}{1})\left(\phi_{\muGamma, \aGamma}^{\AGamma}\right)^{-1} \circ \left(\phi_{\quot{d_0^1}{\Gamma}, \oalphaGamma}^{\AGamma}\right)^{-1} = \left(\phi_{\muGamma,\odGamma{0}{1}}^{\overline{F}(\aGamma)}\right)^{-1} \circ \left(\phi_{\muGamma \circ \odGamma{0}{1}, \aGamma}^{\AGamma}\right)^{-1}.
	\]
	We also derive that
	\[
	\phi_{\quot{d_2^1}{\Gamma}, \opiGamma}^{\AGamma} \circ \overline{F}(\odGamma{2}{1})\phi_{\muGamma, \pGamma}^{\AGamma} = \phi_{\muGamma \circ \odGamma{2}{1}, \pGamma}^{\AGamma} \circ \phi_{\odGamma{2}{1}, \muGamma}^{\overline{F}(\pGamma)\AGamma}
	\]
	and
	\[
	\overline{F}(\odGamma{2}{1})\left(\phi_{\muGamma, \aGamma}^{\AGamma}\right)^{-1} \circ \left(\phi_{\quot{d_2^1}{\Gamma}, \oalphaGamma}^{\AGamma}\right)^{-1} = \left(\phi_{\muGamma,\odGamma{2}{1}}^{\overline{F}(\aGamma)}\right)^{-1} \circ \left(\phi_{\muGamma \circ \odGamma{2}{1}, \aGamma}^{\AGamma}\right)^{-1}
	\]
	through the same technique. Putting these together we get that
	\begin{align*}
		&\phi_{\quot{d_0^1}{\Gamma}, \opiGamma}^{\AGamma} \circ \overline{F}(\odGamma{2}{1})\phi_{\muGamma, \pGamma}^{\AGamma} \circ \big(\overline{F}(\odGamma{2}{1})\circ \overline{F}(\muGamma)\big)\left(\big(\tau_{\pGamma}^{A}\big)^{-1} \circ \tau_{\aGamma}^A\right) \\
		&\circ \overline{F}(\odGamma{0}{1})\left(\phi_{\muGamma, \aGamma}^{\AGamma}\right)^{-1} \circ \left(\phi_{\quot{d_0^1}{\Gamma}, \oalphaGamma}^{\AGamma}\right)^{-1} \circ \phi_{\quot{d_2^1}{\Gamma},\opiGamma}^{\AGamma} \circ \overline{F}(\odGamma{2}{1})\phi_{\muGamma, \pGamma}^{\AGamma} \\
		&\circ \big(\overline{F}(\odGamma{2}{1}) \circ \overline{F}(\muGamma)\big)\left(\big(\tau_{\pGamma}^{A}\big)^{-1} \circ \tau_{\aGamma}^{A}\right) \circ \overline{F}(\odGamma{2}{1})\left(\phi_{\muGamma, \aGamma}^{\AGamma}\right)^{-1} \circ \left(\phi_{\quot{d_2^1}{\Gamma},\oalphaGamma}^{\AGamma}\right)^{-1} \\
		&= \phi_{\muGamma \circ \odGamma{0}{1}, \pGamma}^{\AGamma} \circ \phi_{\odGamma{0}{1}, \muGamma}^{\overline{F}(\pGamma)\AGamma} \circ \big(\overline{F}(\odGamma{2}{1})\circ \overline{F}(\muGamma)\big)\left(\big(\tau_{\pGamma}^{A}\big)^{-1} \circ \tau_{\aGamma}^A\right) \\
		& \circ \left(\phi_{\muGamma,\odGamma{0}{1}}^{\overline{F}(\aGamma)}\right)^{-1} \circ \left(\phi_{\muGamma \circ \odGamma{0}{1}, \aGamma}^{\AGamma}\right)^{-1} \circ \phi_{\muGamma \circ \odGamma{2}{1}, \pGamma}^{\AGamma} \circ \phi_{\odGamma{2}{1}, \muGamma}^{\overline{F}(\pGamma)\AGamma} \\
		&\circ \big(\overline{F}(\odGamma{2}{1}) \circ \overline{F}(\muGamma)\big)\left(\big(\tau_{\pGamma}^{A}\big)^{-1} \circ \tau_{\aGamma}^{A}\right) \circ \left(\phi_{\muGamma,\odGamma{2}{1}}^{\overline{F}(\aGamma)}\right)^{-1} \circ \left(\phi_{\muGamma \circ \odGamma{2}{1}, \aGamma}^{\AGamma}\right)^{-1}.
	\end{align*}
	
	To simplify the expression above, we now focus on manipulating the transition isomorphisms and their interactions with the compositor natural isomorphisms. To begin this process note that since the diagram
	\[
	\begin{tikzcd}
		\big(\overline{F}(\odGamma{0}{1})\circ \overline{F}(\muGamma)\big)\big(\overline{F}(\aGamma)(\AGamma)\big) \ar[rrr]{}{\big(\overline{F}(\odGamma{0}{1}) \circ \overline{F}(\muGamma)\big)\tau_{\aGamma}^{A}} & & & \big(\overline{F}(\odGamma{0}{1}) \circ \overline{F}(\muGamma)\big)(\quot{A}{\Gamma_c}) \\ \overline{F}(\muGamma \circ \odGamma{0}{1})\big(\overline{F}(\aGamma)(\AGamma)\big) \ar[u]{}{\left(\phi_{\odGamma{0}{1}, \muGamma}^{\AGamma}\right)^{-1}} \ar[rrr, swap]{}{\overline{F}(\muGamma \circ \odGamma{0}{1})\tau_{\aGamma}^{A}} & & & \overline{F}(\muGamma \circ \odGamma{0}{1})(\quot{A}{\Gamma_c}) \ar[u, swap]{}{\left(\phi_{\odGamma{0}{1}, \muGamma}^{\quot{A}{\Gamma_c}}\right)^{-1}}
	\end{tikzcd}
	\]
	the equation
	\[
	\big(\overline{F}(\odGamma{0}{1}) \circ \overline{F}(\muGamma)\big)\tau_{\aGamma}^{A} \circ \left(\phi_{\odGamma{0}{1}, \muGamma}^{\AGamma}\right)^{-1} = \left(\phi_{\odGamma{0}{1}, \muGamma}^{\quot{A}{\Gamma_c}}\right)^{-1} \circ \overline{F}(\muGamma \circ \odGamma{0}{1})\tau_{\aGamma}^{A}
	\]
	holds. Similarly, the naturality of the compositors and their inverses allows us to derive that
	\[
	\big(\overline{F}(\odGamma{0}{1}) \circ \overline{F}(\muGamma)\big)(\tau_{\pGamma}^{A})^{-1} \circ \left(\phi_{\odGamma{0}{1}, \muGamma}^{\quot{A}{\Gamma_c}}\right)^{-1} = \left(\phi_{\odGamma{0}{1}, \muGamma}^{\overline{F}(\pGamma)\AGamma}\right)^{-1} \circ \overline{F}(\muGamma \circ \odGamma{0}{1})(\tau_{\pGamma}^{A})^{-1}
	\]
	which allows us to further calculate that
	\begin{align*}
		&\phi_{\muGamma \circ \odGamma{0}{1}, \pGamma}^{\AGamma} \circ \phi_{\odGamma{0}{1}, \muGamma}^{\overline{F}(\pGamma)\AGamma} \circ \big(\overline{F}(\odGamma{2}{1})\circ \overline{F}(\muGamma)\big)\left(\big(\tau_{\pGamma}^{A}\big)^{-1} \circ \tau_{\aGamma}^A\right)\circ \left(\phi_{\muGamma,\odGamma{0}{1}}^{\overline{F}(\aGamma)}\right)^{-1} \\
		&= \phi_{\muGamma \circ \odGamma{0}{1}, \pGamma}^{\AGamma} \circ \phi_{\odGamma{0}{1}, \muGamma}^{\overline{F}(\pGamma)\AGamma} \circ \left(\phi_{\odGamma{0}{1}, \muGamma}^{\overline{F}(\pGamma)\AGamma}\right)^{-1} \circ \overline{F}(\muGamma \circ \odGamma{0}{1})\big((\tau_{\pGamma}^{A})^{-1} \circ \tau_{\aGamma}^{A}\big) \\
		&= \phi_{\muGamma \circ \odGamma{0}{1}, \pGamma}^{\AGamma} \circ \overline{F}(\muGamma \circ \odGamma{0}{1})\big((\tau_{\pGamma}^{A})^{-1} \circ \tau_{\aGamma}^{A}\big).
	\end{align*}
	Now by using Corollary \ref{Cor: F of muG circ d01 is F of muG circ d21} we get that 
	\[
	\overline{F}(\muGamma \circ \odGamma{0}{1})(\quot{A}{\Gamma_c}) = \overline{F}(\muGamma \circ \odGamma{2}{1})(\quot{A}{\Gamma_c}).
	\]
	Using this we derive that the diagram
	\[
	\begin{tikzcd}
		\overline{F}(\aGamma \circ \muGamma \circ \odGamma{0}{1})(\AGamma) \ar[rr, equals] \ar[d, swap]{}{\left(\phi_{\muGamma \circ \odGamma{0}{1}}^{\AGamma}\right)^{-1}} & & \overline{F}(\pGamma \circ \muGamma \circ \odGamma{2}{1})(\AGamma) \ar[d]{}{\left(\phi_{\muGamma \circ \odGamma{2}{1}}^{\AGamma}\right)^{-1}} \\
		\overline{F}(\muGamma \circ \odGamma{0}{1})\big(\overline{F}(\aGamma)(\AGamma)\big) \ar[d, swap]{}{\overline{F}(\muGamma \circ \odGamma{0}{1})} & & \overline{F}(\muGamma \circ \odGamma{2}{1})\big(\overline{F}(\pGamma)(\AGamma)\big) \ar[d]{}{\overline{F}(\muGamma \circ \odGamma{2}{1})\tau_{\pGamma}^{A}} \ar[ll]{}{\left(\phi_{\muGamma \circ \odGamma{0}{1}}^{\AGamma}\right)^{-1} \circ \phi_{\muGamma \circ \odGamma{2}{1}}^{\AGamma}} \\
		\overline{F}(\muGamma \circ \odGamma{0}{1})(\quot{A}{\Gamma_c}) \ar[rr, equals] & & \overline{F}(\muGamma \circ \odGamma{2}{1})(\quot{A}{\Gamma_c})
	\end{tikzcd}
	\]
	commutes which tells us that
	\[
	\overline{F}(\muGamma \circ \odGamma{2}{1})\tau_{\pGamma}^{A} = \overline{F}(\muGamma \circ \odGamma{0}{1})\tau_{\aGamma}^{A} \circ \left(\phi_{\muGamma \circ \odGamma{0}{1}}^{\AGamma}\right)^{-1} \circ \phi_{\muGamma \circ \odGamma{2}{1}}^{\AGamma}.
	\]
	This allows us to compute that
	\begin{align*}
		&\overline{F}(\muGamma \circ \odGamma{0}{1})\tau_{\aGamma}^{A} \circ \left(\phi_{\muGamma \circ \odGamma{0}{1}}^{\AGamma}\right)^{-1} \circ \phi_{\muGamma \circ \odGamma{2}{1}}^{\AGamma} \circ \phi_{\odGamma{2}{1}, \muGamma}^{\overline{F}(\pGamma)\AGamma} \\
		&\circ \big(\overline{F}(\odGamma{2}{1}) \circ \overline{F}(\muGamma)\big)(\tau_{\pGamma}^{A})^{-1} \\
		&=\overline{F}(\muGamma \circ \odGamma{2}{1})\tau_{\pGamma}^{A} \circ \phi_{\muGamma \circ \odGamma{2}{1}}^{\AGamma} \circ \big(\overline{F}(\odGamma{2}{1}) \circ \overline{F}(\muGamma)\big)(\tau_{\pGamma}^{A})^{-1} \\
		&=\overline{F}(\muGamma \circ \odGamma{2}{1})\tau_{\pGamma}^{A} \overline{F}(\muGamma \circ \odGamma{2}{1})(\tau_{\pGamma}^{A})^{-1} \circ \phi_{\muGamma \circ \odGamma{2}{1}}^{\quot{A}{\Gamma_c}} \\
		&= \phi_{\muGamma \circ \odGamma{2}{1}}^{\quot{A}{\Gamma_c}}.
	\end{align*}
	Using these various calculations we find that
	\begin{align*}
		& \phi_{\muGamma \circ \odGamma{0}{1}, \pGamma}^{\AGamma} \circ \phi_{\odGamma{0}{1}, \muGamma}^{\overline{F}(\pGamma)\AGamma} \circ \big(\overline{F}(\odGamma{2}{1})\circ \overline{F}(\muGamma)\big)\left(\big(\tau_{\pGamma}^{A}\big)^{-1} \circ \tau_{\aGamma}^A\right) \\
		& \circ \left(\phi_{\muGamma,\odGamma{0}{1}}^{\overline{F}(\aGamma)}\right)^{-1} \circ \left(\phi_{\muGamma \circ \odGamma{0}{1}, \aGamma}^{\AGamma}\right)^{-1} \circ \phi_{\muGamma \circ \odGamma{2}{1}, \pGamma}^{\AGamma} \circ \phi_{\odGamma{2}{1}, \muGamma}^{\overline{F}(\pGamma)\AGamma} \\
		&\circ \big(\overline{F}(\odGamma{2}{1}) \circ \overline{F}(\muGamma)\big)\left(\big(\tau_{\pGamma}^{A}\big)^{-1} \circ \tau_{\aGamma}^{A}\right) \circ \left(\phi_{\muGamma,\odGamma{2}{1}}^{\overline{F}(\aGamma)}\right)^{-1} \circ \left(\phi_{\muGamma \circ \odGamma{2}{1}, \aGamma}^{\AGamma}\right)^{-1} \\
		&= \phi_{\muGamma \circ \odGamma{0}{1}, \pGamma}^{\AGamma} \circ \overline{F}(\muGamma \circ \odGamma{0}{1})\big((\tau_{\pGamma}^{A})^{-1} \circ \tau_{\aGamma}^{A}\big)  \circ \left(\phi_{\muGamma \circ \odGamma{0}{1}, \aGamma}^{\AGamma}\right)^{-1} \circ \phi_{\muGamma \circ \odGamma{2}{1}, \pGamma}^{\AGamma}  \\
		&\circ \phi_{\odGamma{2}{1}, \muGamma}^{\overline{F}(\pGamma)\AGamma} \circ \big(\overline{F}(\odGamma{2}{1}) \circ \overline{F}(\muGamma)\big)\left(\big(\tau_{\pGamma}^{A}\big)^{-1} \circ \tau_{\aGamma}^{A}\right) \circ \left(\phi_{\muGamma,\odGamma{2}{1}}^{\overline{F}(\aGamma)}\right)^{-1} \\
		&\circ \left(\phi_{\muGamma \circ \odGamma{2}{1}, \aGamma}^{\AGamma}\right)^{-1} \\
		&= \phi_{\muGamma \circ \odGamma{0}{1}, \pGamma}^{\AGamma} \circ \overline{F}(\muGamma \circ \odGamma{0}{1})(\tau_{\pGamma}^{A})^{-1} \circ \phi_{\odGamma{2}{1}, \muGamma}^{\quot{A}{\Gamma_c}} \circ \big(\overline{F}(\odGamma{2}{1}) \circ \overline{F}(\muGamma)\big)\tau_{\aGamma}^{A}  \\
		&\circ \left(\phi_{\muGamma,\odGamma{2}{1}}^{\overline{F}(\aGamma)}\right)^{-1}\circ \left(\phi_{\muGamma \circ \odGamma{2}{1}, \aGamma}^{\AGamma}\right)^{-1}  \\
		&=\phi_{\muGamma \circ \odGamma{0}{1}, \pGamma}^{\AGamma} \circ \overline{F}(\muGamma \circ \odGamma{0}{1})(\tau_{\pGamma}^{A})^{-1} \circ \overline{F}(\muGamma \circ \odGamma{2}{1})\tau_{\aGamma}^{A} \circ \phi_{\odGamma{2}{1}, \muGamma}^{\overline{F}(\aGamma)\AGamma} \\
		&\circ \left(\phi_{\muGamma,\odGamma{2}{1}}^{\overline{F}(\aGamma)}\right)^{-1}\circ \left(\phi_{\muGamma \circ \odGamma{2}{1}, \aGamma}^{\AGamma}\right)^{-1}  \\
		&= \phi_{\muGamma \circ \odGamma{0}{1}, \pGamma}^{\AGamma} \circ \overline{F}(\muGamma \circ \odGamma{0}{1})(\tau_{\pGamma}^{A})^{-1} \circ \overline{F}(\muGamma \circ \odGamma{2}{1})\tau_{\aGamma}^{A} \circ \left(\phi_{\muGamma \circ \odGamma{2}{1}, \aGamma}^{\AGamma}\right)^{-1}.
	\end{align*}
	
	To further manipulate the expression above into our desired format, we will need some calculations to translate between the various compositions of $\overline{F}(\pGamma \circ \muGamma \circ \odGamma{0}{1})$ and $\overline{F}(\pGamma \circ \muGamma \circ \odGamma{1}{1})$ (and similarly with $\overline{F}(\aGamma \circ \muGamma \circ \odGamma{2}{1})$ and $\overline{F}(\aGamma \circ \muGamma \circ \odGamma{1}{1})$). Note that because
	\[
	\overline{F}(\pGamma \circ \muGamma \circ \odGamma{0}{1}) = \overline{F}(\pGamma \circ \muGamma \circ \odGamma{1}{1})
	\]
	we have the $2$-commutative diagram:
	\[
	\begin{tikzcd}
		&\overline{F}(\quot{X}{\Gamma_c}) \ar[ddr, near start, ""{name = UpDiag}]{}{\overline{F}(\muGamma \circ \odGamma{0}{1})} \ar[rr, ""{name = UpRight}]{}{\overline{F}(\muGamma)} & & \overline{F}(\quot{(G \times X)}{\Gamma}) \ar[ddl, ""{name = UpUp}]{}{\overline{F}(\odGamma{0}{1})} \\
		\\
		\overline{F}(\XGamma) \ar[ddr, swap]{}{\overline{F}(\pGamma)} \ar[uur]{}{\overline{F}(\pGamma)} \ar[rr, bend left = 20, ""{name = MidUp}]{}{\overline{F}(\pGamma \circ \muGamma \circ \odGamma{0}{1})} \ar[rr, bend right = 20, swap, ""{name = MidLow}]{}{\overline{F}(\pGamma \circ \muGamma \circ \odGamma{1}{1})} & & \overline{F}(\quot{(G \times G \times X)}{\Gamma}) \\
		\\
		& \overline{F}(\quot{X}{\Gamma_c}) \ar[uur, swap, near start, ""{name = DownDiag}]{}{\overline{F}(\muGamma \circ \odGamma{1}{1})} \ar[rr, swap, ""{name = DownRight}]{}{\overline{F}(\muGamma)} & & \overline{F}(\quot{(G \times X)}{\Gamma}) \ar[uul, swap]{}{\overline{F}(\odGamma{1}{1})}
		\ar[from = MidUp, to= MidLow, equals, shorten <= 4pt, shorten >= 4pt]
		\ar[from = 1-2, to = MidUp, Rightarrow, shorten <= 4pt, shorten >= 8pt]
		\ar[from = UpRight, to = 3-3, Rightarrow, shorten <= 4pt, shorten >= 4pt]
		\ar[from = 5-2, to = MidLow, Rightarrow, shorten <= 4pt, shorten >= 8pt]
		\ar[from = DownRight, to = 3-3, Rightarrow, shorten <= 4pt, shorten >= 4pt]
	\end{tikzcd}
	\]
	These allow us to deduce that the (equivalent) pasting diagrams
	\[
	\begin{tikzcd}
		\overline{F}(\XGamma) \ar[r, bend left = 30, ""{name = U}]{}{\overline{F}(\muGamma) \circ \overline{F}(\pGammap)} \ar[r, bend right = 30, swap, ""{name = M}]{}{\overline{F}(\opiGamma)} \ar[rr, bend right = 60, swap, ""{name = L}]{}{\overline{F}(\pGamma \circ \muGamma \circ \odGamma{0}{1})} & \overline{F}(\quot{(G \times X)}{\Gamma}) \ar[r]{}{\overline{F}(\odGamma{0}{1})} & \overline{F}(\quot{(G \times G \times X)}{\Gamma}) \ar[from = U, to = M, Rightarrow, shorten <= 4pt, shorten >= 4pt]{}{} \ar[from = 1-2, to = L, Rightarrow, shorten <= 4pt, shorten >= 4pt]{}{}
	\end{tikzcd}
	\]
	\[
	\begin{tikzcd}
		\overline{F}(\XGamma) \ar[r]{}{\overline{F}(\pGamma)} \ar[rr, bend right = 60, swap, ""{name = L}]{}{\overline{F}(\pGamma \circ \muGamma \circ \odGamma{0}{1})} & \overline{F}(\quot{X}{\Gamma_c}) \ar[r, bend left = 30, ""{name = U}]{}{\overline{F}(\odGamma{0}{1}) \circ \overline{F}(\muGamma)} \ar[r, bend right = 30, swap, ""{name = M}]{}{\overline{F}(\muGamma \circ \odGamma{0}{1})} & \overline{F}(\quot{(G \times G \times X)}{\Gamma}) \ar[from = U, to = M, Rightarrow, shorten <= 4pt, shorten >= 4pt] \ar[from = 1-2, to = L, Rightarrow, shorten <= 4pt, shorten >= 4pt]{}{}
	\end{tikzcd}
	\]
	are equal to the (also equivalent) pasting diagrams
	\[
	\begin{tikzcd}
		\overline{F}(\XGamma) \ar[r, bend left = 30, ""{name = U}]{}{\overline{F}(\muGamma) \circ \overline{F}(\pGammap)} \ar[r, bend right = 30, swap, ""{name = M}]{}{\overline{F}(\opiGamma)} \ar[rr, bend right = 60, swap, ""{name = L}]{}{\overline{F}(\pGamma \circ \muGamma \circ \odGamma{1}{1})} & \overline{F}(\quot{(G \times X)}{\Gamma}) \ar[r]{}{\overline{F}(\odGamma{1}{1})} & \overline{F}(\quot{(G \times G \times X)}{\Gamma}) \ar[from = U, to = M, Rightarrow, shorten <= 4pt, shorten >= 4pt]{}{} \ar[from = 1-2, to = L, Rightarrow, shorten <= 4pt, shorten >= 4pt]{}{}
	\end{tikzcd}
	\]
	\[
	\begin{tikzcd}
		\overline{F}(\XGamma) \ar[r]{}{\overline{F}(\pGamma)} \ar[rr, bend right = 60, swap, ""{name = L}]{}{\overline{F}(\pGamma \circ \muGamma \circ \odGamma{1}{1})} & \overline{F}(\quot{X}{\Gamma_c}) \ar[r, bend left = 30, ""{name = U}]{}{\overline{F}(\odGamma{1}{1}) \circ \overline{F}(\muGamma)} \ar[r, bend right = 30, swap, ""{name = M}]{}{\overline{F}(\muGamma \circ \odGamma{1}{1})} & \overline{F}(\quot{(G \times G \times X)}{\Gamma}) \ar[from = U, to = M, Rightarrow, shorten <= 4pt, shorten >= 4pt] \ar[from = 1-2, to = L, Rightarrow, shorten <= 4pt, shorten >= 4pt]{}{}
	\end{tikzcd}
	\]
	and hence deduce that
	\begin{align*}
		\phi_{\opiGamma, \odGamma{1}{1}}^{\AGamma} \circ \overline{F}(\odGamma{1}{1})\phi_{\muGamma, \pGamma}^{\AGamma} &= \phi_{\muGamma \circ \odGamma{1}{1}, \pGamma}^{\AGamma} \circ \phi_{\odGamma{1}{1}, \muGamma}^{\overline{F}(\pGamma)\AGamma} = \phi_{\muGamma \circ \odGamma{0}{1}, \pGamma}^{\AGamma} \circ \phi_{\odGamma{0}{1}, \muGamma}^{\overline{F}(\pGamma)\AGamma}.
	\end{align*}
	Similarly, from the equality
	\[
	\overline{F}(\aGamma \circ \muGamma \circ \odGamma{2}{1}) = \overline{F}(\aGamma \circ \muGamma \circ \odGamma{1}{1})
	\]
	we get the $2$-cell
	\[
	\begin{tikzcd}
		&\overline{F}(\quot{X}{\Gamma_c}) \ar[ddr, near start, ""{name = UpDiag}]{}{\overline{F}(\muGamma \circ \odGamma{2}{1})} \ar[rr, ""{name = UpRight}]{}{\overline{F}(\muGamma)} & & \overline{F}(\quot{(G \times X)}{\Gamma}) \ar[ddl, ""{name = UpUp}]{}{\overline{F}(\odGamma{2}{1})} \\
		\\
		\overline{F}(\XGamma) \ar[ddr, swap]{}{\overline{F}(\aGamma)} \ar[uur]{}{\overline{F}(\aGamma)} \ar[rr, bend left = 20, ""{name = MidUp}]{}{\overline{F}(\aGamma \circ \muGamma \circ \odGamma{2}{1})} \ar[rr, bend right = 20, swap, ""{name = MidLow}]{}{\overline{F}(\aGamma \circ \muGamma \circ \odGamma{1}{1})} & & \overline{F}(\quot{(G \times G \times X)}{\Gamma}) \\
		\\
		& \overline{F}(\quot{X}{\Gamma_c}) \ar[uur, swap, near start, ""{name = DownDiag}]{}{\overline{F}(\muGamma \circ \odGamma{1}{1})} \ar[rr, swap, ""{name = DownRight}]{}{\overline{F}(\muGamma)} & & \overline{F}(\quot{(G \times X)}{\Gamma}) \ar[uul, swap]{}{\overline{F}(\odGamma{1}{1})}
		\ar[from = MidUp, to= MidLow, equals, shorten <= 4pt, shorten >= 4pt]
		\ar[from = 1-2, to = MidUp, Rightarrow, shorten <= 4pt, shorten >= 8pt]
		\ar[from = UpRight, to = 3-3, Rightarrow, shorten <= 4pt, shorten >= 4pt]
		\ar[from = 5-2, to = MidLow, Rightarrow, shorten <= 4pt, shorten >= 8pt]
		\ar[from = DownRight, to = 3-3, Rightarrow, shorten <= 4pt, shorten >= 4pt]
	\end{tikzcd}
	\]
	which allows us to deduce that
	\begin{align*}
		\left(\phi_{\odGamma{2}{1},\muGamma}^{\overline{F}(\aGamma)\AGamma}\right)^{-1} \circ \left(\phi_{\muGamma \circ \odGamma{2}{1}, \aGamma}^{\AGamma}\right)^{-1} &= \left(\phi_{\odGamma{1}{1}, \muGamma}^{\overline{F}(\aGamma)\AGamma}\right)^{-1} \circ \left(\phi_{\muGamma \circ \odGamma{1}{1}, \aGamma}^{\AGamma}\right)^{-1} \\
		&= \left(\overline{F}(\odGamma{1}{1})\phi_{\muGamma, \pGamma}^{\AGamma}\right)^{-1} \circ \left(\phi_{\oalphaGamma, \odGamma{1}{1}}^{\AGamma}\right)^{-1}.
	\end{align*}
	
	Before we can complete the proof of the theorem, we need three last observations. First note that by using the compositor isomorphisms we have that
	\[
	\overline{F}(\muGamma \circ \odGamma{2}{1})\tau_{\aGamma}^{A} = \phi_{\odGamma{2}{1}, \muGamma}^{\quot{A}{\Gamma_c}} \circ \big(\overline{F}(\odGamma{2}{1}) \circ \overline{F}(\muGamma)\big)\tau_{\aGamma}^{A} \circ \left(\phi_{\odGamma{2}{1}, \muGamma}^{\overline{F}(\aGamma)\AGamma}\right)^{-1}
	\]
	and similarly that
	\[
	\overline{F}(\muGamma \circ \odGamma{0}{1})(\tau_{\pGamma}^{A})^{-1} = \phi_{\odGamma{0}{1}, \muGamma}^{\overline{F}(\pGamma)\AGamma} \circ \big(\overline{F}(\odGamma{0}{1}) \circ \overline{F}(\muGamma)\big)(\tau_{\pGamma}^{A})^{-1} \circ \left(\phi_{\odGamma{0}{1}, \muGamma}^{\quot{A}{\Gamma_c}}\right)^{-1}.
	\]
	However, using that
	\begin{align*}
		&\big(\overline{F}(\odGamma{0}{1}) \circ \overline{F}(\muGamma)\big)(\tau_{\pGamma}^{A})^{-1} \circ \left(\phi_{\odGamma{0}{1}, \muGamma}^{\quot{A}{\Gamma_c}}\right)^{-1} \circ \phi_{\odGamma{2}{1}, \muGamma}^{\quot{A}{\Gamma_c}} \circ \big(\overline{F}(\odGamma{2}{1}) \circ \overline{F}(\muGamma)\big)\tau_{\aGamma}^{A} \\
		&= \big(\overline{F}(\odGamma{1}{1}) \circ \overline{F}(\muGamma)\big)(\tau_{\pGamma}^{A})^{-1} \circ \left(\phi_{\odGamma{1}{1}, \muGamma}^{\quot{A}{\Gamma_c}}\right)^{-1} \circ \phi_{\odGamma{1}{1}, \muGamma}^{\quot{A}{\Gamma_c}} \circ \big(\overline{F}(\odGamma{1}{1}) \circ \overline{F}(\muGamma)\big)\tau_{\aGamma}^{A} \\
		&= \big(\overline{F}(\odGamma{1}{1}) \circ \overline{F}(\muGamma)\big)(\tau_{\pGamma}^{A})^{-1} \circ \big(\overline{F}(\odGamma{1}{1}) \circ \overline{F}(\muGamma)\big)\tau_{\aGamma}^{A} \\
		&= \big(\overline{F}(\odGamma{1}{1}) \circ \overline{F}(\muGamma)\big)\left((\tau_{\pGamma}^{A})^{-1} \circ \tau_{\aGamma}^{A}\right)
	\end{align*}
	we derive that
	\begin{align*}
		&\phi_{\quot{d_0^1}{\Gamma}, \opiGamma}^{\AGamma} \circ \overline{F}(\quot{d_0^1}{\Gamma})\quot{\theta}{\Gamma} \circ \left(\phi_{\quot{d_0^1}{\Gamma}, \oalphaGamma}^{\AGamma}\right)^{-1} \circ \phi_{\quot{d_2^1}{\Gamma},\opiGamma}^{\AGamma} \circ \overline{F}(\quot{d_2^1}{\Gamma})\quot{\theta}{\Gamma} \circ \left(\phi_{\quot{d_2^1}{\Gamma},\oalphaGamma}^{\AGamma}\right)^{-1} \\
		&=\phi_{\muGamma \circ \odGamma{0}{1}, \pGamma}^{\AGamma} \circ \overline{F}(\muGamma \circ \odGamma{0}{1})(\tau_{\pGamma}^{A})^{-1} \circ \overline{F}(\muGamma \circ \odGamma{2}{1})\tau_{\aGamma}^{A} \circ \left(\phi_{\muGamma \circ \odGamma{2}{1}, \aGamma}^{\AGamma}\right)^{-1} \\
		&= \phi_{\muGamma \circ \odGamma{0}{1}, \pGamma}^{\AGamma} \circ \phi_{\odGamma{0}{1}, \muGamma}^{\overline{F}(\pGamma)\AGamma} \circ \big(\overline{F}(\odGamma{1}{1}) \circ \overline{F}(\muGamma)\big)\left((\tau_{\pGamma}^{A})^{-1} \circ \tau_{\aGamma}^{A}\right) \\
		&\circ \left(\phi_{\odGamma{2}{1}, \muGamma}^{\overline{F}(\aGamma)\AGamma}\right)^{-1} \circ \left(\phi_{\muGamma \circ \odGamma{2}{1}, \aGamma}^{\AGamma}\right)^{-1} \\
		&= \phi_{\opiGamma, \odGamma{1}{1}}^{\AGamma} \circ \overline{F}(\odGamma{1}{1})\phi_{\muGamma, \pGamma}^{\AGamma} \circ \big(\overline{F}(\odGamma{1}{1}) \circ \overline{F}(\muGamma)\big)\left((\tau_{\pGamma}^{A})^{-1} \circ \tau_{\aGamma}^{A}\right) \\ 
		&\circ \left(\overline{F}(\odGamma{1}{1})\phi_{\muGamma, \pGamma}^{\AGamma}\right)^{-1} \circ \left(\phi_{\oalphaGamma, \odGamma{1}{1}}^{\AGamma}\right)^{-1} \\
		&= \phi_{\opiGamma, \odGamma{1}{1}}^{\AGamma} \circ \overline{F}(\odGamma{1}{1})\quot{\theta}{\Gamma} \circ \left(\phi_{\oalphaGamma, \odGamma{1}{1}}^{\AGamma}\right)^{-1}.
	\end{align*}
	This establishes the GIT cocycle condition.
\end{proof}
\begin{definition}\label{Defn: Section Equivariance: Equivariance}
	Let $G$ be a smooth algebraic group and let $X$ be a left $G$-variety or let $G$ be a topological group and let $X$ be a left $G$-space. If $(F,\overline{F})$ is a pre-equivariant pseudofunctor on $X$ and $A \in F_G(X)_0$ then we say that the natural isomorphism $\Theta^A$ is the equivariance\index[terminology]{Equivariance} on $A$.
\end{definition}

We now close this chapter with a short comparison between how the equivariance $\Theta$ attached to a pre-equivariant pseudofunctor $(F,\overline{F})$ interacts with a $G$-equivariant morphism $h:X \to Y$. Because a pre-equivariant pseudofunctor $(F,\overline{F})$ attached to a variety $X$ comes equipped with the pseudofunctor $\overline{F}:\Var_{/K}^{\op} \to \fCat$, there is another pre-equivariant pseudofunctor $(F^{\prime}, \overline{F})$ on $Y$ with $F^{\prime} = \overline{F} \circ \quo_{G \backslash (-)}:\SfResl_G(Y)^{\op} \to \fCat$ (and similarly for the topological situation). The fact that the morphism $h:X \to Y$ is equivariant implies that the constructions in Lemma \ref{Lemma: Clifton's Naivification Variety} are natural in the sense that for any $\Gamma$ in $\Sf(G)_0$ the diagrams
\[
\begin{tikzcd}
	\quot{(G \times X)}{\Gamma} \ar[r]{}{\mu_{\Gamma}^{X}} \ar[d, swap]{}{G \backslash (\id_{\Gamma \times G} \times f)} & \quot{X}{\Gamma_c} \ar[d]{}[description]{G \backslash (\id_{\Gamma_c} \times f)} \ar[r]{}{a_{\Gamma}^{X}} & \quot{X}{\Gamma} \ar[d]{}{G \backslash (\id_{\Gamma} \times  f)} \\
	\quot{(G \times Y)}{\Gamma} \ar[r, swap]{}{\mu_{\Gamma}^{Y}} & \quot{Y}{\Gamma_c} \ar[r, swap]{}{a_{\Gamma}^{Y}} & \quot{Y}{\Gamma}
\end{tikzcd}
\]
and
\[
\begin{tikzcd}
	\quot{(G \times X)}{\Gamma} \ar[r]{}{\mu_{\Gamma}^{X}} \ar[d, swap]{}{G \backslash (\id_{\Gamma \times G} \times f)} & \quot{X}{\Gamma_c} \ar[d]{}[description]{G \backslash (\id_{\Gamma_c} \times f)} \ar[r]{}{p_{\Gamma}^{X}} & \quot{X}{\Gamma} \ar[d]{}{G \backslash (\id_{\Gamma} \times  f)} \\
	\quot{(G \times Y)}{\Gamma} \ar[r, swap]{}{\mu_{\Gamma}^{Y}} & \quot{Y}{\Gamma_c} \ar[r, swap]{}{p_{\Gamma}^{Y}} & \quot{Y}{\Gamma}
\end{tikzcd}
\]
both commute. Because of this we can apply the formalism of Theorem \ref{Thm: Pseudocone Functors: Pullback induced by fibre functors in pseudofucntor} to give rise to functors $h^{\ast}:F_G(Y) \to F_G(X)$ and $(\id_G \times h)^{\ast}:F_G(G \times Y) \to F_G(G \times X)$. We can then consider the two distinct equivariances
\[
\begin{tikzcd}
	F_G(X) \ar[rr, bend left = 20, ""{name = U}]{}{\alpha_X^{\ast}} \ar[rr, bend right = 20, swap, ""{name = D}]{}{(\pi_2^X)^{\ast}} & & F_G(G \times X) \ar[from = U, to = D, Rightarrow, shorten <= 4pt, shorten >= 4pt]{}{\Theta_X}
\end{tikzcd}
\]
and:
\[
\begin{tikzcd}
	F_G(Y) \ar[rr, bend left = 20, ""{name = U}]{}{\alpha_X^{\ast}} \ar[rr, bend right = 20, swap, ""{name = D}]{}{(\pi_2^Y)^{\ast}} & & F_G(G \times Y) \ar[from = U, to = D, Rightarrow, shorten <= 4pt, shorten >= 4pt]{}{\Theta_Y}
\end{tikzcd}
\]
By pre-and-post whiskering by pullbacks along $h$ in different ways we get natural isomorphisms
\[
\Theta_X \ast h^{\ast}:\alpha_X^{\ast} \circ h^{\ast} \xRightarrow{\cong} \left(\pi_2^X\right)^{\ast} \circ h^{\ast}
\]
and 
\[
(\id_G \times h)^{\ast} \ast \Theta_Y: (\id_G \times h)^{\ast} \circ \alpha_Y^{\ast} \xRightarrow{\cong} (\id_G \times h)^{\ast} \circ \left(\pi_2^Y\right)^{\ast}.
\]
However, because every pullback functor in sight is induced by $F$ and $\Gamma$-locally has the form
\[
\quot{\left(\alpha_X^{\ast} \circ h^{\ast}\right)\left(A\right)}{\Gamma} = \left(\overline{F}\left(\quot{\alpha_X}{\Gamma}\right) \circ \overline{F}\left(\quot{h}{\Gamma}\right)\right)\left(\quot{A}{\Gamma}\right),
\]
\[
\quot{\left(\left(\id_G \times h\right)^{\ast} \circ \alpha_Y^{\ast}\right)\left(A\right)}{\Gamma} = \left(\overline{F}\left(\quot{\overline{h}}{\Gamma \times G}\right) \circ \overline{F}\left(\quot{\alpha_Y}{\Gamma}\right)\right)\left(\quot{A}{\Gamma}\right)
\]
and similarly for the $\pi_2$ versions. However, because $h$ is $G$-equivariant, we have that $\alpha_Y \circ (\id_G \times h) = h \circ \alpha_X$ and so there are natural isomorphisms
\[
\alpha_X^{\ast} \circ h^{\ast} \cong \left(\id_G \times h\right)^{\ast} \circ \alpha_Y^{\ast}
\]
and similarly
\[
\left(\pi_2^X\right)^{\ast} \circ h^{\ast} \cong (\id_G \times h)^{\ast} \circ \left(\pi_2^Y\right)^{\ast}.
\]
This leads to a natural question: how different are the distinct pullbacks of the equivariances? More explicitly, what is the difference, if any, between the natural isomorphisms
$
(\id_G \times h)^{\ast} \ast \Theta_Y
$
and
$
\Theta_X \ast h^{\ast}
$
when they are manipulated through the compositors to the pullback functor $(h \circ \alpha_X)^{\ast} = (\alpha_Y \circ (\id_G \times h))^{\ast}$? We give an explicit answer in Proposition \ref{Prop: Section Equivariance: Pullbacks and Equivariances} below and show there is no difference in the strict pseudofunctor case in Corollary \ref{Cor: Section Equivariance: Pullbacks and Equivariances for Strict pseudofunctors}.

\begin{proposition}\label{Prop: Section Equivariance: Pullbacks and Equivariances}
	For any smooth algebraic group $G$ and $G$-equivariant morphism of varieties $h:X \to Y$, if $(\overline{F} \circ \quo_X^{\op}, \overline{F})$ and $(\overline{F} \circ \quo_{Y}^{\op}, \overline{F})$ are pre-equivariant morphisms on $X$ and $Y$, respectively, then the natural isomorphisms $(\id_G \times h)^{\ast}(\Theta_Y^A)$ and $\Theta_X^{h^{\ast}A}$ satisfy the equation
	\[
	\phi_{\quot{h}{\Gamma \times G}, \quot{\pi_2^{Y}}{\Gamma}}^{A} \circ \quot{(\id_G \times h)^{\ast}\left(\Theta_{Y}^{A}\right)}{\Gamma} \circ \left(\phi_{\quot{h}{\Gamma \times G}, \quot{\alpha_Y}{\Gamma}}^{A}\right)^{-1} = \phi_{\quot{\pi_2^X}{\Gamma}, \quot{h}{\Gamma}}^{A} \circ \quot{\Theta_X^{h^{\ast}A}}{\Gamma} \circ \left(\phi_{\quot{\alpha_X}{\Gamma}, \quot{h}{\Gamma}}^{A}\right)^{-1}
	\]
	for all $\Gamma \in \Sf(G)_0$ and for every $A \in F_G(Y)_0$. The same holds for topological groups $G$, left $G$-spaces $X$, and all $M \in \mathbf{FAMan}(G)_0$.
\end{proposition}
\begin{proof}
	We prove this for the geometric case, as the topological case follows mutatis mutandis. The pseudofunctoriality of the pre-equivariant pseudofunctor $(F,\overline{F})$ implies that for any object $A$ of $F_G(Y)$ and for any $\Gamma \in \Sf(G)_0$ the diagrams
	\[
	\begin{tikzcd}
		\overline{F}\left(\mu_{\Gamma}^{X}\right)\left(\overline{F}(\quot{h}{\Gamma_c})\left(\overline{F}(p_{\Gamma}^{Y})\left(\quot{A}{\Gamma}\right)\right)\right) \ar[d, swap]{}{\overline{F}\left(\mu_{\Gamma}^{X}\right)\left(\phi_{\quot{h}{\Gamma_c}, p_{\Gamma}^{Y}}^{A}\right)} \ar[r]{}{\phi_{\mu_{\Gamma}^{X},\quot{h}{\Gamma_c}}^{\overline{F}(p_{\Gamma}^{Y})A}} & \overline{F}\left(\quot{h}{\Gamma_c} \circ \mu_{\Gamma}^{X}\right)\left(\overline{F}\left(p_{\Gamma}^{Y}\right)\left(\quot{A}{\Gamma}\right)\right) \ar[d, equals] \\
		\overline{F}\left(\mu_{\Gamma}^{X}\right)\left(\overline{F}\left(p_{\Gamma}^{Y} \circ \quot{h}{\Gamma_c}\right)\left(\quot{A}{\Gamma}\right)\right) \ar[d, swap]{}{\overline{F}\left(\mu_{\Gamma}^{X}\right)\left(\phi_{p_{\Gamma}^{X},\quot{h}{\Gamma}}^{A}\right)^{-1}} & \overline{F}\left(\mu_{\Gamma}^{Y} \circ \quot{\overline{h}}{\Gamma_c}\right)\left(\overline{F}\left(p_{\Gamma}^{Y}\right)\left(\quot{A}{\Gamma}\right)\right) \ar[d]{}{\left(\phi_{\mu_{\Gamma}^{Y},\quot{h}{\Gamma \times G}}^{\overline{F}(p_{\Gamma}^{Y})A}\right)^{-1}} \\
		\left(\overline{F}\left(\mu_{\Gamma}^{X}\right) \circ \overline{F}\left(p_{\Gamma}^{X}\right) \circ \overline{F}\left(\quot{h}{\Gamma}\right)\right)\left(\quot{A}{\Gamma}\right) \ar[d, swap]{}{\phi_{\mu_{\Gamma}^{X},p_{\Gamma}^{X}}^{\overline{F}(\quot{h}{\Gamma})A}} &  \left(\overline{F}\left(\quot{\overline{h}}{\Gamma \times G}\right)\circ \overline{F}\left(\mu_{\Gamma}^{Y}\right)\circ \overline{F}\left(p_{\Gamma}^{Y}\right)\right)\left(\quot{A}{\Gamma}\right) \ar[d]{}{\overline{F}\left(\quot{h}{\Gamma \times G}\right)\left(\phi_{\mu_{\Gamma}^{Y},p_{\Gamma}^{Y}}^A\right)} \\
		\overline{F}\left(\quot{\pi_{2}^{X}}{\Gamma}\right)\left(\overline{F}\left(\quot{h}{\Gamma}\right)\left(\quot{A}{\Gamma}\right)\right) \ar[d, swap]{}{\phi_{\quot{\pi_{2}^{X}}{\Gamma}, \quot{h}{\Gamma}}^{A}} & \overline{F}\left(\quot{h}{\Gamma \times G}\right)\left(\overline{F}\left(\quot{\overline{\pi}_2^Y}{\Gamma}\right)\left(\quot{A}{\Gamma}\right)\right) \\
		\overline{F}\left(\quot{h}{\Gamma} \circ \quot{\pi_2^{X}}{\Gamma}\right)\left(\quot{A}{\Gamma}\right) \ar[ur, swap]{}{\left(\phi_{\quot{h}{\Gamma \times G}, \quot{\pi_2^Y}{\Gamma}}^{A}\right)^{-1}}
	\end{tikzcd}
	\]
	and
	\[
	\begin{tikzcd}
		\left(\overline{F}\left(\quot{h}{\Gamma \times G}\right) \circ \overline{F}\left(a_{\Gamma}^{Y} \circ \mu_{\Gamma}^{Y}\right)\right)\left(\quot{A}{\Gamma}\right) \ar[r]{}{\overline{F}\left(\quot{h}{\Gamma \times G}\right)} \ar[d, swap]{}{\phi_{\quot{h}{\Gamma \times G}, \quot{\alpha_Y}{\Gamma}}^{A}} & \left(\overline{F}\left(\quot{h}{\Gamma \times G}\right) \circ \overline{F}\left(\mu_{\Gamma}^{Y}\right) \circ \overline{F}\left(a_{\Gamma}^{Y}\right)\right) \ar[d]{}{\phi_{\quot{h}{\Gamma\times G},\mu_{\Gamma}^{Y}}^{\overline{F}(a_{\Gamma}^{Y})A}} \\
		\overline{F}\left(\quot{\alpha_Y}{\Gamma} \circ \quot{h}{\Gamma \times G}\right)\left(\quot{A}{\Gamma}\right) \ar[d, swap]{}{\left(\phi_{\quot{\alpha_X}{\Gamma},\quot{h}{\Gamma \times G}}^{A}\right)^{-1}} & \left(\overline{F}\left(\mu_{\Gamma}^{X}\right) \circ \overline{F}\left(\quot{h}{\Gamma_c}\right) \circ \overline{F}\left(a_{\Gamma}^{Y}\right)\right)\left(\quot{A}{\Gamma}\right) \ar[d]{}{\left(\phi_{\mu_{\Gamma}^{X},\quot{h}{\Gamma_c}}^{\overline{F}(a_{\Gamma}^{Y})A}\right)^{-1}} \\
		\left(\overline{F}\left(a_{\Gamma}^{X} \circ \mu_{\Gamma}^{X}\right) \circ \overline{F}\left(\quot{h}{\Gamma}\right)\right)\left(\quot{A}{\Gamma}\right) \ar[d, swap]{}{\left(\phi_{\mu_{\Gamma}^{X},a_{\Gamma}^{X}}^{\overline{F}(\quot{h}{\Gamma})A}\right)^{-1}} & \left(\overline{F}\left(\mu_{\Gamma}^{X}\right) \circ \overline{F}\left(\quot{h}{\Gamma_c}\right) \circ \overline{F}\left(a_{\Gamma}^{Y}\right)\right)\left(\quot{A}{\Gamma}\right) \\
		\left(\overline{F}\left(\mu_{\Gamma}^{X}\right) \circ \overline{F}\left(a_{\Gamma}^{X}\right) \circ \overline{F}\left(\quot{h}{\Gamma}\right)\right)\left(\quot{A}{\Gamma}\right) \ar[r,swap]{}{\overline{F}\left(\mu_{\Gamma}^{X}\right)\left(\phi_{a_{\Gamma}^{X},\quot{h}{\Gamma}}\right)} & \left(\overline{F}\left(\mu_{\Gamma}^{X}\right) \circ \overline{F}\left(\quot{h}{\Gamma} \circ a_{\Gamma}^{X}\right)\right)\left(\quot{A}{\Gamma}\right) \ar[u, swap]{}{\overline{F}\left(\mu_{\Gamma}^{X}\right)\left(\phi_{\quot{h}{\Gamma_c},a_{\Gamma}^{Y}}^{A}\right)^{-1}}
	\end{tikzcd}
	\]
	both commute. We then calculate on one hand that
	\begin{align*}
		&\quot{\Theta^{h^{\ast}A}_{X}}{\Gamma} = \quot{\theta^{\overline{F}(\quot{h}{\Gamma})A}_{X}}{\Gamma} \\
		&=\phi_{\mu_{\Gamma}^{X},p_{\Gamma}^{X}}^{\overline{F}(\quot{h}{\Gamma})A} \circ \overline{F}\left(\mu_{\Gamma}^{X}\right)\left(\left(\tau_{p_{\Gamma}^{X}}^{h^{\ast}A}\right)^{-1} \circ \tau_{a_{\Gamma}^{X}}^{h^{\ast}A}\right) \circ \left(\phi_{\mu_{\Gamma}^{X}, a_{\Gamma}^{X}}^{\overline{F}(\quot{h}{\Gamma})A}\right)^{-1} \\
		&= \phi_{\mu_{\Gamma}^{X},p_{\Gamma}^{X}}^{\overline{F}(\quot{h}{\Gamma})A} \circ \overline{F}\left(\mu_{\Gamma}^{X}\right)\left(\left(\overline{F}\left(\quot{h}{\Gamma_c}\right)\left(\tau_{p_{\Gamma}^{Y}}^{A}\right) \circ \left(\phi_{\quot{h}{\Gamma_c}, p_{\Gamma}^{Y}}^{A}\right) \circ \phi_{p_{\Gamma}^{X}, \quot{h}{\Gamma}}^{A}\right)^{-1}\right) \\
		& \circ \overline{F}\left(\mu_{\Gamma}^{X}\right)\left(\overline{F}\left(\quot{h}{\Gamma_c}\right)\left(\tau_{a_{\Gamma}^{Y}}^{A}\right) \circ \left(\phi_{\quot{h}{\Gamma_c},a_{\Gamma}^{Y}}^{A}\right)^{-1} \circ \phi_{a_{\Gamma}^{X}, \quot{h}{\Gamma}}^{A}\right) \circ \left(\phi_{\mu_{\Gamma}^{X}, a_{\Gamma}^{X}}^{\overline{F}(\quot{h}{\Gamma})A}\right)^{-1} \\
		&= \phi_{\mu_{\Gamma}^{X},p_{\Gamma}^{X}}^{\overline{F}(\quot{h}{\Gamma})A} \circ \overline{F}\left(\mu_{\Gamma}^{X}\right)\left(\left(\phi_{p_{\Gamma}^{X},\quot{h}{\Gamma}}^{A}\right)^{-1} \circ \phi_{\quot{h}{\Gamma_c}, p_{\Gamma}^{Y}}^{A}\right) \circ \left(\overline{F}\left(\mu_{\Gamma}^{X}\right) \circ  \overline{F}\left(\quot{h}{\Gamma_c}\right)\right)\left(\tau_{p_{\Gamma}^{Y}}^{A}\right)^{-1} \\
		&\circ \left(\overline{F}\left(\mu_{\Gamma}^{X}\right) \circ \overline{F}\left(\quot{h}{\Gamma_c}\right)\right)\left(\tau_{a_{\Gamma}^{Y}}^{A}\right) \circ \overline{F}\left(\mu_{\Gamma}^{X}\right)\left(\left(\phi_{\quot{h}{\Gamma_c},a_{\Gamma}^{Y}}^{A}\right)^{-1} \circ \phi_{a_{\Gamma}^{X}, \quot{h}{\Gamma}}^{A}\right) \circ \left(\phi_{\mu_{\Gamma}^{X}, a_{\Gamma}^{X}}^{\overline{F}(\quot{h}{\Gamma})A}\right)^{-1}
	\end{align*}
	which we suggestively display as
	\begin{align*}
		\quot{\Theta^{h^{\ast}A}_{X}}{\Gamma} &= \phi_{\mu_{\Gamma}^{X},p_{\Gamma}^{X}}^{\overline{F}(\quot{h}{\Gamma})A} \circ \overline{F}\left(\mu_{\Gamma}^{X}\right)\left(\left(\phi_{p_{\Gamma}^{X},\quot{h}{\Gamma}}^{A}\right)^{-1} \circ \phi_{\quot{h}{\Gamma_c}, p_{\Gamma}^{Y}}^{A}\right)  \\
		&\circ \left(\overline{F}\left(\mu_{\Gamma}^{X}\right)\circ  \overline{F}\left(\quot{h}{\Gamma_c}\right)\right)\left(\left(\tau_{p_{\Gamma}^{Y}}^{A}\right)^{-1}\circ \tau_{a_{\Gamma}^{Y}}^{A}\right) \\
		&\circ \overline{F}\left(\mu_{\Gamma}^{X}\right)\left(\left(\phi_{\quot{h}{\Gamma_c},a_{\Gamma}^{Y}}^{A}\right)^{-1} \circ \phi_{a_{\Gamma}^{X}, \quot{h}{\Gamma}}^{A}\right) \circ \left(\phi_{\mu_{\Gamma}^{X}, a_{\Gamma}^{X}}^{\overline{F}(\quot{h}{\Gamma})A}\right)^{-1}
	\end{align*}
	On the other hand we calculate
	\begin{align*}
		&\quot{(\id_G\times h)^{\ast}\left(\theta_Y^A\right)}{\Gamma} = \overline{F}\left(\quot{h}{\Gamma \times G}\right)\left(\quot{\theta_Y^{A}}{\Gamma}\right) \\
		&=\overline{F}\left(\quot{h}{\Gamma \times G}\right)\left(\phi_{\mu_{\Gamma}^{Y}, p_{\Gamma}^{Y}}^{A} \circ \overline{F}\left(\mu_{\Gamma}^{Y}\right)\left(\left(\tau_{p_{\Gamma}^{Y}}^{A}\right)^{-1} \circ \tau_{a_{\Gamma}^{Y}}^{A}\right) \circ \left(\phi_{\mu_{\Gamma}^{Y},a_{\Gamma}^{Y}}^{A}\right)^{-1}\right) \\
		&= \overline{F}\left(\quot{h}{\Gamma \times G}\right)\left(\phi_{\mu_{\Gamma}^{Y}, p_{\Gamma}^{Y}}^{A}\right) \circ \left(\overline{F}\left(\quot{h}{\Gamma \times G}\right) \circ \overline{F}\left(\mu_{\Gamma}^{Y}\right)\right)\left(\left(\tau_{p_{\Gamma}^{Y}}^{A}\right)^{-1} \circ \tau_{a_{\Gamma}^{Y}}^{A}\right) \\
		&\circ  \overline{F}\left(\quot{h}{\Gamma \times G}\right)\left(\phi_{\mu_{\Gamma}^{Y},a_{\Gamma}^{Y}}^{A}\right)^{-1} \\
		&= \overline{F}\left(\quot{h}{\Gamma \times G}\right)\left(\phi_{\mu_{\Gamma}^{Y}, p_{\Gamma}^{Y}}^{A}\right) \circ \left(\phi_{\quot{h}{\Gamma \times G}, \mu_{\Gamma}^{Y}}^{\overline{F}(p_{\Gamma}^{Y}A)}\right)^{-1} \circ  \overline{F}\left(\mu_{\Gamma}^{Y} \circ \quot{h}{\Gamma \times G}\right)\left(\left(\tau_{p_{\Gamma}^{Y}}^{A}\right)^{-1} \circ  \tau_{a_{\Gamma}^{Y}}^{A}\right) \\
		&\circ \phi_{\quot{h}{\Gamma \times G}, \mu_{\Gamma}^{Y}}^{\overline{F}(a_{\Gamma}^{Y})A} \circ \overline{F}\left(\quot{h}{\Gamma \times G}\right)\left(\phi_{\mu_{\Gamma}^{Y},a_{\Gamma}^{Y}}^{A}\right)^{-1} \\
		&= \overline{F}\left(\quot{h}{\Gamma \times G}\right)\left(\phi_{\mu_{\Gamma}^{Y}, p_{\Gamma}^{Y}}^{A}\right) \circ \left(\phi_{\quot{h}{\Gamma \times G}, \mu_{\Gamma}^{Y}}^{\overline{F}(p_{\Gamma}^{Y}A)}\right)^{-1} \circ  \overline{F}\left(\quot{h}{\Gamma_c} \circ \mu_{\Gamma}^{X}\right)\left(\left(\tau_{p_{\Gamma}^{Y}}^{A}\right)^{-1} \circ  \tau_{a_{\Gamma}^{Y}}^{A}\right) \\
		&\circ \phi_{\quot{h}{\Gamma \times G}, \mu_{\Gamma}^{Y}}^{\overline{F}(a_{\Gamma}^{Y})A} \circ \overline{F}\left(\quot{h}{\Gamma \times G}\right)\left(\phi_{\mu_{\Gamma}^{Y},a_{\Gamma}^{Y}}^{A}\right)^{-1} \\
		&= \overline{F}\left(\quot{h}{\Gamma \times G}\right)\left(\phi_{\mu_{\Gamma}^{Y}, p_{\Gamma}^{Y}}^{A}\right) \circ \left(\phi_{\quot{h}{\Gamma \times G}, \mu_{\Gamma}^{Y}}^{\overline{F}(p_{\Gamma}^{Y}A)}\right)^{-1} \circ \phi_{\mu_{\Gamma}^{X}, \quot{h}{\Gamma_c}}^{\overline{F}(p_{\Gamma}^{Y})A} \circ \left(\overline{F}\left(\mu_{\Gamma}^{X}\right) \circ \overline{F}\left(\mu_{\Gamma}^{X}\right)\right)\left(\tau_{p_{\Gamma}^{Y}}^{A}\right)^{-1} \\
		&\circ \left(\overline{F}\left(\mu_{\Gamma}^{X}\right) \circ \overline{F}\left(\mu_{\Gamma}^{X}\right)\right)\left(\tau_{a_{\Gamma}^{Y}}^{A} \right)
		\circ \left(\phi_{\mu_{\Gamma}^{X},\quot{h}{\Gamma \times G}}^{\overline{F}(a_{\Gamma}^{Y})A}\right)^{-1}\circ \phi_{\quot{h}{\Gamma \times G}, \mu_{\Gamma}^{Y}}^{\overline{F}(a_{\Gamma}^{Y})A} \circ \overline{F}\left(\quot{h}{\Gamma \times G}\right)\left(\phi_{\mu_{\Gamma}^{Y},a_{\Gamma}^{Y}}^{A}\right)^{-1}.
	\end{align*}
	Using the first commutative diagram tells us that
	\begin{align*}
		&\overline{F}\left(\quot{h}{\Gamma \times G}\right)\left(\phi_{\mu_{\Gamma}^{Y}, p_{\Gamma}^{Y}}^{A}\right) \circ \left(\phi_{\quot{h}{\Gamma \times G}, \mu_{\Gamma}^{Y}}^{\overline{F}(p_{\Gamma}^{Y}A)}\right)^{-1} \circ \phi_{\mu_{\Gamma}^{X}, \quot{h}{\Gamma_c}}^{\overline{F}(p_{\Gamma}^{Y})A} \\
		&= \left(\phi_{\quot{h}{\Gamma \times G}, \quot{\pi_2^Y}{\Gamma}}^{A}\right)^{-1} \circ \phi_{\quot{\pi_{2}^{X}}{\Gamma}, \quot{h}{\Gamma}}^{A} \circ \phi_{\mu_{\Gamma}^{X},p_{\Gamma}^{X}}^{\overline{F}(\quot{h}{\Gamma})A} \circ \overline{F}\left(\mu_{\Gamma}^{X}\right)\left(\left(\phi_{p_{\Gamma}^{X},\quot{h}{\Gamma}}^{A}\right)^{-1} \circ \left(\phi_{\quot{h}{\Gamma_c}, p_{\Gamma}^{Y}}^{A}\right)\right)
	\end{align*}
	while the second gives that
	\begin{align*}
		&\left(\phi_{\mu_{\Gamma}^{X},\quot{h}{\Gamma \times G}}^{\overline{F}(a_{\Gamma}^{Y})A}\right)^{-1}\circ \phi_{\quot{h}{\Gamma \times G}, \mu_{\Gamma}^{Y}}^{\overline{F}(a_{\Gamma}^{Y})A} \circ \overline{F}\left(\quot{h}{\Gamma \times G}\right)\left(\phi_{\mu_{\Gamma}^{Y},a_{\Gamma}^{Y}}^{A}\right)^{-1} \\
		&=\overline{F}\left(\mu_{\Gamma}^{X}\right)\left(\left(\phi_{\quot{h}{\Gamma_c}, a_{\Gamma}^{Y}}^{A}\right)^{-1} \circ \phi_{a_{\Gamma}^{X}, \quot{h}{\Gamma}}^{A}\right) \circ \phi_{\mu_{\Gamma}^{X}, a_{\Gamma}^{X}}^{\overline{F}(\quot{h}{\Gamma})} \circ \left(\phi_{\quot{\alpha_X}{\Gamma}, \quot{h}{\Gamma}}^{A}\right)^{-1} \circ \phi_{\quot{h}{\Gamma \times G}, \quot{\alpha_Y}{\Gamma}}^{A}.
	\end{align*} 
	Substituting these equations into our description of $\quot{\Theta_X^{h^{\ast}A}}{\Gamma}$ gives that
	\begin{align*}
		&\left(\phi_{\quot{h}{\Gamma \times G}, \quot{\pi_2^Y}{\Gamma}}^{A}\right)^{-1} \circ \phi_{\quot{\pi_{2}^{X}}{\Gamma}, \quot{h}{\Gamma}}^{A} \circ \quot{\Theta_X^{h^{\ast}A}}{\Gamma} \circ \left(\phi_{\quot{\alpha_X}{\Gamma}, \quot{h}{\Gamma}}^{A}\right)^{-1} \circ \phi_{\quot{h}{\Gamma \times G}, \quot{\alpha_Y}{\Gamma}}^{A}\\
		&= \left(\phi_{\quot{h}{\Gamma \times G}, \quot{\pi_2^Y}{\Gamma}}^{A}\right)^{-1} \circ \phi_{\quot{\pi_{2}^{X}}{\Gamma}, \quot{h}{\Gamma}}^{A} \circ \phi_{\mu_{\Gamma}^{X},p_{\Gamma}^{X}}^{\overline{F}(\quot{h}{\Gamma})A} \circ \overline{F}\left(\mu_{\Gamma}^{X}\right)\left(\left(\phi_{p_{\Gamma}^{X},\quot{h}{\Gamma}}^{A}\right)^{-1} \circ \left(\phi_{\quot{h}{\Gamma_c}, p_{\Gamma}^{Y}}^{A}\right)\right) \\
		&\circ \left(\overline{F}\left(\mu_{\Gamma}^{X}\right)\circ  \overline{F}\left(\quot{h}{\Gamma_c}\right)\right)\left(\left(\tau_{p_{\Gamma}^{Y}}^{A}\right)^{-1}\circ \tau_{a_{\Gamma}^{Y}}^{A}\right) \\
		& \circ \overline{F}\left(\mu_{\Gamma}^{X}\right)\left(\left(\phi_{\quot{h}{\Gamma_c}, a_{\Gamma}^{Y}}^{A}\right)^{-1} \circ \phi_{a_{\Gamma}^{X}, \quot{h}{\Gamma}}^{A}\right) \circ \phi_{\mu_{\Gamma}^{X}, a_{\Gamma}^{X}}^{\overline{F}(\quot{h}{\Gamma})} \circ \left(\phi_{\quot{\alpha_X}{\Gamma}, \quot{h}{\Gamma}}^{A}\right)^{-1} \circ \phi_{\quot{h}{\Gamma \times G}, \quot{\alpha_Y}{\Gamma}}^{A} \\
		&= \overline{F}\left(\quot{h}{\Gamma \times G}\right)\left(\phi_{\mu_{\Gamma}^{Y}, p_{\Gamma}^{Y}}^{A}\right) \circ \left(\phi_{\quot{h}{\Gamma \times G}, \mu_{\Gamma}^{Y}}^{\overline{F}(p_{\Gamma}^{Y}A)}\right)^{-1} \circ \phi_{\mu_{\Gamma}^{X}, \quot{h}{\Gamma_c}}^{\overline{F}(p_{\Gamma}^{Y})A} \\
		&\circ \left(\overline{F}\left(\mu_{\Gamma}^{X}\right)\circ  \overline{F}\left(\quot{h}{\Gamma_c}\right)\right)\left(\left(\tau_{p_{\Gamma}^{Y}}^{A}\right)^{-1}\circ \tau_{a_{\Gamma}^{Y}}^{A}\right) \\
		&\circ \left(\phi_{\mu_{\Gamma}^{X},\quot{h}{\Gamma \times G}}^{\overline{F}(a_{\Gamma}^{Y})A}\right)^{-1}\circ \phi_{\quot{h}{\Gamma \times G}, \mu_{\Gamma}^{Y}}^{\overline{F}(a_{\Gamma}^{Y})A} \circ \overline{F}\left(\quot{h}{\Gamma \times G}\right)\left(\phi_{\mu_{\Gamma}^{Y},a_{\Gamma}^{Y}}^{A}\right)^{-1} \\
		&= \quot{(\id_{G} \times h)^{\ast}\Theta_{Y}^{A}}{\Gamma}
	\end{align*}
	and hence proves that
	\[
	\quot{(\id_{G} \times h)^{\ast}\Theta_{Y}^{A}}{\Gamma} = \left(\phi_{\quot{h}{\Gamma \times G}, \quot{\pi_2^Y}{\Gamma}}^{A}\right)^{-1} \circ \phi_{\quot{\pi_{2}^{X}}{\Gamma}, \quot{h}{\Gamma}}^{A} \circ \quot{\Theta_X^{h^{\ast}A}}{\Gamma} \circ \left(\phi_{\quot{\alpha_X}{\Gamma}, \quot{h}{\Gamma}}^{A}\right)^{-1} \circ \phi_{\quot{h}{\Gamma \times G}, \quot{\alpha_Y}{\Gamma}}^{A}
	\]
	for every $\Gamma \in \Sf(G)_0$. Rearranging the equation to have the form
	\begin{equation}\label{Eqn: Discrepancy between equivariances and pullbacks}
		\phi_{\quot{h}{\Gamma \times G}, \quot{\pi_2^{Y}}{\Gamma}}^{A} \circ \quot{(\id_G \times h)^{\ast}\left(\Theta_{Y}^{A}\right)}{\Gamma} \circ \left(\phi_{\quot{h}{\Gamma \times G}, \quot{\alpha_Y}{\Gamma}}^{A}\right)^{-1} = \phi_{\quot{\pi_2^X}{\Gamma}, \quot{h}{\Gamma}}^{A} \circ \quot{\Theta_X^{h^{\ast}A}}{\Gamma} \circ \left(\phi_{\quot{\alpha_X}{\Gamma}, \quot{h}{\Gamma}}^{A}\right)^{-1}
	\end{equation}
	gives the equation in the statement of the proposition.
\end{proof}
\begin{corollary}\label{Cor: Section Equivariance: Pullbacks and Equivariances for Strict pseudofunctors}
	Let $G$ be a smooth algebraic group and let $h:X \to Y$ be a $G$-equivariant morphism of varieties. Then if $F$ is a strict pre-equviariant pseudofunctor $(\id_G \times h)^{\ast}(\Theta_Y^A) = \Theta_X^{h^{\ast}A}$ for every object $A$ of $F_G(Y)$.
\end{corollary}
\begin{proof}
	Equation \ref{Eqn: Discrepancy between equivariances and pullbacks} gives a measure between the difference between the isomorphisms $(\id_G \times h)^{\ast}\Theta_Y^A$ and $\Theta_X^{h^{\ast}A}$. However, for a strict pseudofunctor the compositors become identity maps and the corollary follows.
\end{proof}

\newpage

\chapter{Change of Group Functors}\label{Chapter: Chofg}
Let $\varphi:G \to H$ be a morphism of algebraic groups and let $X$ be a left $H$-variety or let $\varphi: G \to H$ be a morphism of topological groups and let $X$ be a left $H$-space. In both cases, by translating the action of $G$ to $H$ along $\varphi$, we can equip $X$ with a $G$-action. This technique is frequently important in representation theory, algebraic geometry, and algebraic topology and especially of interest when $G$ is a subgroup of $H$. Some examples of these arise, of course, in the theory of endoscopy in the Langlands Programme. Some explicit examples for the case of algebraic groups are when:
\begin{itemize}
	\item $G = \GL_2$ or $G = \GL_1 \times \GL_1$ and $H = G_2$ over a $p$-adic field $F$;
	\item $G = G_2$ and $H = \operatorname{SO}_7$ over a $p$-adic field $F$;
\end{itemize} 
see \cite{G2Cubics} for more details regarding the $p$-adic algebraic group $G_2$. Because of this we want to show how we can take an $H$-equivariant category, $F_H(X)$, and produce a functor $\varphi^{\sharp}:F_H(X) \to F_G(X)$ which provides a algebraic/categorical mirror that reflects the geometric situation described by translating the $H$-action of $X$ to a $G$-action through the map $\varphi$. 

While the Change of Groups functors have been in principle known in the equivariant derived categorical case since \cite{LusztigCuspidal2} and \cite{BernLun}, there are significant technical issues with the presentation therein. In \cite{BernLun} Change of Groups functors are only defined for subgroup inclusions $i:G \to H$. The issues in \cite{LusztigCuspidal2} are more significant, however, and there are significant details of the construction which are either missing or not even mentioned as a technical issue which needs verification. Consequently much of what we do is provide the equivariant descent-theoretic background needed in order to do the Change of Groups fully, carefully, and at a full level of detail as well as extend it to the higher-categorical pseudocone expression of equivariant descent.

\section{Equivariant Functors: Change of Groups}\label{Section: Section 3: Change of Groups}

In what follows in the geometric situation we fix smooth algebraic groups $G$ and $H$ over the field $K$ and let $\varphi \in \AlgGrp(G,H)$ be a morphism where $\AlgGrp$\index[notation]{SAlgGrp@$\AlgGrp$} is the category of smooth algebraic groups over $\Spec K$. Given the algebraic group $G$, we also write $\mu_G:G \times G \to G$\index[notation]{Mutipliciation@$\mu_G$} for the multiplication map, $\inv_G:G \to G$\index[notation]{InverseMapGroup@$\inv_G$} for the inversion map, and $1_G:\Spec K \to G$\index[notation]{Unit@$1_G$} for the unit morphism of $G$ (and write similarly $H$). 

In the topological situation we fix topological groups $G$ and $H$ and let $\varphi:G \to H$ be a morphism of topological groups. As in the geometric case we write $\mu_G:G \times G \to G$ for the multiplication, $1_G:\lbrace \ast \rbrace \to G$ for the unit map, and $\inv_G:G \to G$ for the inversion map (and similarly for $H$).

\begin{lemma}\label{Lemma: Section 3: Obvious Pullback functor for restriction of action}
	There is a pullback functor $\varphi^{\ast}:\HVar \to \GVar$ (respectively $\varphi^{\ast}:H$-$\Top \to G$-$\Top$) given by sending an $H$-variety $(X,\alpha_X:H \times X \to X)$ to the $G$-variety $(X,\alpha_X\circ (\varphi \times \id_X):G \times X \to X)$ and sending an $H$-equivariant morphism $f:X \to Y$ to itself (and similarly for the topological case).
\end{lemma}
\begin{proof}
	The verification that the morphism
	\[
	\xymatrix{
		G \times X \ar[r]^-{\varphi \times \id_X} & H \times X \ar[r]^-{\alpha_X} & X
	}
	\]
	makes $X$ into a $G$-variety is a trivial application of the fact that $X$ is an $H$-variety and the fact that $\varphi$ is a morphism of algebraic groups. Similarly, that $f$ remains as a $G$-equivariant morphism follows because the diagram
	\[
	\xymatrix{
		G \times X \ar[rrr]^-{\alpha_X \circ (\varphi \times \id_X)} \ar[d]_{\id_G \times f} & & & X \ar[d]^{f} \\
		G \times Y \ar[rrr]_-{\alpha_Y \circ (\varphi \times \id_Y)} & & & Y
	}
	\]
	factors as
	\[
	\xymatrix{
		G \times X \ar[r]^-{\varphi \times \id_X} \ar[d]_{\id_G \times f} & H \times X \ar[r]^-{\alpha_X} \ar[d]^{\id_H \times f} & X \ar[d]^{f} \\
		G \times Y \ar[r]_-{\varphi \times \id_Y} & H \times Y \ar[r]_-{\alpha_Y} & Y
	}
	\]
	and both squares commute.
\end{proof}
\begin{remark}
	We will write the action of $G$ on $\varphi^{\ast}X$ as $\varphi^{\ast}\alpha_X$ or $\alpha_{\varphi^{\ast}X}$ (depending on the situation).
\end{remark}
\begin{definition}\label{Remark: GH pre-equiv pseudofunctors}
	Let $G, H, \varphi,$ and $X$ be given in either two cases:
	\begin{enumerate}
		\item $G$ and $H$ are smooth algebraic groups, $\varphi:G \to H$ is a morphism of algebraic groups, and $X$ is a left $H$-variety;
		\item $G$ and $H$ are topological groups, $\varphi:G \to H$ is a morphism of topological groups, and $X$ is a left $H$-space.
	\end{enumerate}
	In the geometric case a $(G,H)$-pre-equivariant pseudofunctor on $X$ is comprised of a triple $\left((F,\overline{F}),(F^{\prime},\overline{F}^{\prime}), e:\overline{F} \Rightarrow \overline{F}^{\prime}\right)$ where $(F,\overline{F})$ is a pre-equivariant pseudofunctor
	\[
	\begin{tikzcd}
		\SfResl_G(\varphi^{\ast}X)^{\op} \ar[dr, swap]{}{F} \ar[r]{}{\quo_G^{\op}} & \Var_{/K}^{\op} \ar[d]{}{\overline{F}} \\
		& \fCat
	\end{tikzcd}
	\]
	on $\varphi^{\ast}X$, $(F^{\prime},\overline{F}^{\prime})$ is a pre-equivariant pseudofunctor
	\[
	\begin{tikzcd}
		\SfResl_H(X)^{\op} \ar[r]{}{\quo_{H}^{\op}} \ar[dr, swap]{}{F^{\prime}} & \Var_{/K}^{\op} \ar[d]{}{\overline{F}^{\prime}} \\
		& \fCat
	\end{tikzcd}
	\]
	on $X$, and $e:\overline{F} \to \overline{F}^{\prime}$ is an equivalence in the bicategory $\Bicat(\Var_{/K}^{\op},\fCat)$. In the topological case a $(G,H)$-pre-equivariant pseudofunctor on $X$ is comprised of a triple $\left((F,\overline{F}),(F^{\prime},\overline{F}^{\prime}), e:\overline{F} \Rightarrow \overline{F}^{\prime}\right)$ where $(F,\overline{F})$ is a pre-equivariant pseudofunctor
	\[
	\begin{tikzcd}
		\FResl_G(\varphi^{\ast}X)^{\op} \ar[dr, swap]{}{F} \ar[r]{}{\quo_G^{\op}} & \Top^{\op} \ar[d]{}{\overline{F}} \\
		& \fCat
	\end{tikzcd}
	\]
	on $\varphi^{\ast}X$, $(F^{\prime},\overline{F}^{\prime})$ is a pre-equivariant pseudofunctor
	\[
	\begin{tikzcd}
		\FResl_H(X)^{\op} \ar[r]{}{\quo_{H}^{\op}} \ar[dr, swap]{}{F^{\prime}} & \Top^{\op} \ar[d]{}{\overline{F}^{\prime}} \\
		& \fCat
	\end{tikzcd}
	\]
	on $X$, and $e:\overline{F} \to \overline{F}^{\prime}$ is an equivalence in the bicategory $\Bicat(\Top^{\op},\fCat)$.
	
	By abuse of notation we will call $\overline{F}$ a $(G,H)$-pre-equivariant\index[terminology]{Pre-equivariant Pseudofunctor! GH@ $(G,H)$-pre-equivariant Psuedofunctor} pseudofunctor on $X$ and leave all the above data implicit.
\end{definition}

We now give our pseudofunctorial generalization of the change of groups functor presented in \cite{LusztigCuspidal2}. In \cite{LusztigCuspidal2} this functor is given at heart by restricting an object collection $A = \lbrace \AGamma \; | \; \Gamma \in \Sf(G)_0 \rbrace$ to a specific sub-collection which records the information coming from a morphism $\varphi:G \to H$ of algebraic groups; our construction will give the functor in essentially the same way, save that we distinguish when varieties are isomorphic and bang-on equal to each other.

For what follows we will need the induction space of Bernstein as described in \cite{PramodBook}, \cite{BernLun}, \cite{Bien}, and \cite{MirkovicVilonen}. While we present this for varieties, note that the definition given extends freely to topological spaces; in fact, the categorical perspective adapts mutatis mutandis and the set-theoretic sketches readily apply to the topological case in a less informal fashion. To describe the induction space, we first let $Z$ be a $G$-variety and consider the variety $H \times Z$. We now define a $G$-action $\alpha_{H \times Z}$ on $H \times Z$ via the diagram below:
\begin{equation}\label{Eqn: The action on H times X}
	\begin{tikzcd}
		G \times H \times Z \ar[dd, swap]{}{\alpha_{H \times Z}} \ar[rr]{}{\Delta_G \times \id_{H\times Z}} & & G \times G \times H \times Z \ar[r]{}{\cong} & H \times G \times G \times Z \ar[d]{}{\id_H \times \inv_G \times \id_{G \times Z}}\\
		& & & H \times G \times G \times Z \ar[d]{}{\id_H\times \varphi \times \id_{G \times X}} \\
		H \times Z & & & H \times H \times G \times Z \ar[lll]{}{\mu_H\times \alpha_Z}
	\end{tikzcd}
\end{equation}
Note that in set-theoretic terms, this is the action given by 
\[
g\cdot(h,z) := (h\varphi(g)^{-1},gz).
\] 
Denote the quotient variety $G \backslash (H \times X)$ by $H \times^{G} X$; note that it exists because the action of $G$ on $H$ is image-isomorphic to the action of an algebraic subgroup of $H$ on the whole group, and because $H$ is a smooth free $H$-variety it follows that $H \times X$ admits a quotient by actions of algebraic subgroups as well. Note that in the topological case we need not take extra pains to argue for the existence of the quotient defining $H \times^{G}X$, as $\Top$ is cocomplete.

We now define a further action of $H$ on $H \times^{G} X$. Write $[h,x]_G$ for a generic ``element'' of the variety $H \times^{G} X$ and define the $H$-action on $H \times^{G} X$ by
\[
h\cdot[h^{\prime},x]_G := [hh^{\prime},x]_G.
\]
Note that this action is well-defined by a scalar restriction-style argument and the fact that any ``element'' $h^{\prime}$ coming from the algebraic subgroup of $H$ isomorphic to $G/\Ker \varphi$ gets absorbed/moved by the quotient action onto the $X$ component.
\begin{definition}\label{Defn: Induction Space}\index[terminology]{Induction Space}
Let $G, H, \varphi,$ and $X$ be given in either of the two cases:
\begin{itemize}
	\item $G$ and $H$ are smooth algebraic groups, $\varphi:G \to H$ is a morphism of algebraic groups, and $X$ is a left $H$-variety;
	\item $G$ and $H$ are topological groups, $\varphi:G \to H$ is a morphism of topological groups, and $X$ is a left $H$-space.
\end{itemize}
In either case the object $H \times^G X$ together with the $H$-action
	\[
	h\cdot[h^{\prime},x]_G := [hh^{\prime},x]_G
	\]
	is called the induction space\footnote{Perhaps it would be better to call the variety-theoretic version of $H \times^G X$ the induction variety and the topological version the induction space, but because such terminology is nonstandard we will not use it save for perhaps informally.} of $X$. 
\end{definition}
\begin{remark}
	It is common in representation-theoretic literature (cf.\@ \cite{BernLun}, \cite{Bien}, \cite{LusztigCuspidal2}, \cite{MirkovicVilonen}, for instance) to see the above variety written as $H \times_G X$. We will steer far away from that notation, as it conflicts with the pullback notation we have used throughout the monograph up to this point. As a generic rule there should also be \textit{less} notational overlap and clash in mathematics whenever possible, so we follow \cite{PramodBook} and write $H \times^G X$ instead.
\end{remark}

We proceed with a key structural lemma that allows us to see that for any $\Gamma \in \Sf(G)_0$, the scheme $H \times^{G} \Gamma$ is an object in $\Sf(H)$. This works essentially because the actions that arise in this construction and quotients are all smooth, while pure dimension is automatic by construction and the fact that before taking quotients the action has fibres of constant pure dimension as well. To see that this is a principal $H$-variety on the surface seems less trivial than the pure dimension or smoothness of $H \times^{G} \Gamma$. However, because both $\Gamma$ and $H$ admit {\'e}tale locally isotrivial quotients by $G$ and $H$, respectively, it follows that $H \times^{G} \Gamma$ admits an {\'e}tale locally isotrivial quotient by $H$ as well. Thus the actual only nontrivial point is to check that the fibres take the correct form, i.e., for such a trivializing {\'e}tale cover $\lbrace U_i \to H \times^{G} \Gamma \; | \; i \in I \rbrace$, we have $\alpha^{-1}(U_i) \cong H \times U_i$.

\begin{proposition}\label{Prop: Section  3: Bernstein Quotient is Sf(H)}
	Let $G$ and $H$ be smooth algebraic groups over $\Spec K$ with $\varphi \in \AlgGrp(G,H)$ and $\Gamma \in \Sf(G)$. Then $H \times^{G} \Gamma \in \Sf(H)_0$.
\end{proposition}
\begin{proof}
	Let $\alpha:H \times (H \times^{G} \Gamma) \to H \times^{G} \Gamma$ be the action morphism and write $q:H \times^{G} \Gamma \to H\backslash(H \times^{G} \Gamma)$ for the quotient morphism. The discussion prior to the proposition argues why $q$ is finite {\'e}tale locally trivializable so let $\lbrace \varphi_i:U_i \to H\backslash (H \times^{G}\Gamma) \; | \; i \in I \rbrace$ be a finite {\'e}tale trivialization of $q$. Note that because the $U_i$ are a local trivialization of $q$, we have that for each $i \in I$ there is a $Z_i$ for which the diagram
	\[
	\begin{tikzcd}
	q^{-1}(U_i) \ar[r]{}{\cong} \ar[dr, swap]{}{} & Z_i \times U_i \ar[d]{}{\pi_2} \\
	 & U_i
	\end{tikzcd}
	\] 
	commutes. Because $\Gamma$ is a principal $G$-variety (and hence locally isotrivial) and the action $\alpha$ takes the form
	\[
	\alpha(h, [h^{\prime},\gamma]_G) = [hh^{\prime},\gamma]_G
	\]
	it follows that $Z_i \times U_i$ is a subobject of $H \times U_i$. Consequently we need to argue that $H \times U_i$ comes with a map to $Z_i \times U_i \cong q^{-1}(U_i)$ as the universal property of the pullback will give us our isomorphism. For this we give a set-theoretic sketch. Fix an element $u$ of $U_i$ and write $u$ as $[[h,\gamma]_G]_H$. Now note that this implies that every term of the form $(k, [[h,\gamma]_G]_H)$ in $H \times U_i$ sits in the fibre over $u$ as the assignment
	\[
	\left(k,\left[[h,\gamma]_G\right]_H\right) \mapsto \left[kh,\gamma\right]_G
	\] 
	is well-defined by the fact that $\Gamma$ is a free $G$-variety and sits over $u$ because $[[kh,\gamma]_G]_H = [[h,\gamma]_G]_H = u$. But then this gives our desired map 
	\[
	H \times U_i \to q^{-1}(U_i) \xrightarrow{\cong} Z_i \times U_i
	\] 
	and hence allows us to deduce that $H \times U_i \cong Z_i \times U_i$. 
\end{proof}
\begin{remark}
	The argument sketched above also shows that if $G$ and $H$ are topological groups and if $P$ is a free $G$-space (cf.\@ Definition \ref{Defn: Free Gspace}) then the induction space $H \times^{G} P$ is a locally trivializable $H$-space as well. However, mimicking the same argument for the covering families by replacing the finite {\'e}tale covers with open covers allows us to deduce that $H \times^G P$ is a free $H$-space.
\end{remark}
\begin{proposition}\label{Prop: Section  3: Bernstein Quotient is FSf(H)}
	Let $G$ and $H$ be topological groups with $\varphi:G \to H$ a map of topological groups and $P$ a free $G$-space. Then $H \times^{G} P$ is a free $H$-space.
\end{proposition}

We now show that this produces a functor over $\Sf(G)$, which is a necessary step in concluding that we get a functor $F_H(X) \to F_G(X)$. Essentially this just says that the constructions above are functorial in free $G$-objects, which is a routine argument we present below. The crucial ingredient is that if we have a map $f \in \Sf(G)(\Gamma,\Gamma^{\prime})$ then the induced morphism
\[
\id_H \times f:H \times \Gamma \to H \times \Gamma^{\prime}
\]
descends to a map between the induction spaces. This holds because the actions $\alpha_{H \times \Gamma}$ and $\alpha_{H \times \Gamma^{\prime}}$ as described in Diagram \ref{Eqn: The action on H times X} do not intermingle the spaces $H$ and $\Gamma$ (and also $H$ and $\Gamma^{\prime}$) in nontrivial ways. In particular, this says that we get the commuting diagram
\[
\xymatrix{
	H \times X \ar[rr]^-{\id_H\times f} \ar[d]_{\quo_{H \times^G X}} & & H \times Y \ar[d]^{\quo_{H\times Y}} \\
	G\backslash{(H \times X)} \ar[rr]_-{G \backslash(\id_H \times f)} & & G\backslash (H \times Y)
}
\]
We claim the (suggestively written) morphism $\id_H \times^G f := G \backslash (\id_H\times f)$ is also $H$-equivariant, i.e., is a morphism $H \times^G X \to H \times^G Y$ of $H$-varieties. However, this is a routine argument; in set-theoretic terms, this morphism is given via
\[
[h,x]_G \mapsto [h,f(x)]_G
\]
and since the action is given by $h[h^{\prime},x]_G = [hh^{\prime},x]_G$, we get that
\[
(\id_H \times^G f)(h[h^{\prime},x]_G) = [hh^{\prime},f(x)]_G = h[h^{\prime},f(x)]_G = h\left((\id_H \times^G f)[h^{\prime},x]_G\right).
\]
It is also routine to verify using the argument above that this construction is functorial; identities get preserved immediately, while the fact that composition is preserved follows from the uniqueness of the induced maps between quotient schemes. Putting these observations together with Proposition \ref{Prop: Section  3: Bernstein Quotient is Sf(H)} then gives us the corollary below, which is crucial  our change of groups functors. 

\begin{corollary}\label{Cor: Section 3: Functor from Sf(G) to Sf(H) by taking Bernstein quotients}
	There is a functor $\Sf(G) \to \Sf(H)$ which sends a variety $\Gamma$ to the $H$-variety $H \times^{G} \Gamma$ and which sends a morphism $f:\Gamma \to \Gamma^{\prime}$ to the map $(\id_{H} \times^G f):H \times^{G} \Gamma \to H \times^{G} \Gamma^{\prime}$.
\end{corollary}
\begin{proof}
	Given the remarks prior to the statement of this corollary and Proposition \ref{Prop: Section  3: Bernstein Quotient is Sf(H)}, all that we need do is check that the fibres of the smooth morphism $\id_H \times^G f$ are of constant pure dimension for all $f \in \Sf(G)_1$. However, this is immediate from the fact that $\id_H \times f$ has fibres of constant pure dimension ($\id_H$ trivially so and $f$ by assumption) and from the fact that the $G$-action and its quotient do not change fibre dimensions in a non-uniform fashion.
\end{proof}
Note that this translates to the topological case mutatis mutandis.
\begin{corollary}\label{Cor: Section 3: Top Functor from Sf(G) to Sf(H) by taking Bernstein quotients}
	There is a functor $\mathbf{Free}(G) \to \mathbf{Free}(H)$ which sends a free $G$-space $P$ to the free $G$-space $H \times^{G} P$ and which sends a morphism $f:P \to Q$ to $(\id_{H} \times^{G} f):H \times^{G} P \to H \times^{G} Q$.
\end{corollary}

We need one final lemma before we can present the change of groups functors. This lemma realizes an isomorphism of varieties
\[
G \backslash (\Gamma \times X) = \quot{X}{\Gamma} \cong \quot{X}{H \times^G \Gamma} = H \backslash \left(\left(H \times^G \Gamma\right) \times X\right),
\]
and respectively of topological spaces
\[
G \backslash (M \times X) = \quot{X}{M} \cong \quot{X}{H \times^G M} = H \backslash\left(\left(H \times^G M\right) \times X\right), 
\]
which allows us to define an object $A \in F_G(X)$ by restricting an object of $F_H(X)$ to those pieces whose fibres are equivalent to fibres which come from the varieties $\XGamma$ for $\Gamma \in \Sf(G)_0$ in the geometric case and which come from the spaces $\quot{X}{M}$ for $M \in \mathbf{Free}(G)$. It is worth remarking that the geometric version of this is an unproven side remark in \cite{LusztigCuspidal2} which quickly becomes a massive technical issue when working directly with these objects at this pseudocone level of abstraction. Because it is absolutely fundamental to getting the change of groups functors running, we sketch a proof of this fact here.
\begin{lemma}\label{Lemma: Section 3: Induction space quotient is iso to the Sf(G) quotient}
	For any $\Gamma \in \Sf(G)_0$, there is an isomorphism of $K$-varieties
	\[
	G\backslash(\Gamma \times X)= \XGamma \cong \quot{X}{H \times^G \Gamma} = H \backslash \big((H \times^G \Gamma) \times X\big).
	\]
\end{lemma}
\begin{proof}[Sketch]
	We prove this by way of a set-theoretical group action sketch.
	
	Begin by fixing a variety $\Gamma \in \Sf(G)_0$ and recall from Proposition \ref{Prop: Section  3: Bernstein Quotient is Sf(H)} that $H \times^{G} \Gamma$ is an object in $\Sf(H)$. Thus the quotient 
	\[
	H\backslash\big((H \times^{G} \Gamma) \times X\big) =: \quot{X}{H \times^{G} \Gamma}
	\] 
	exists for all $H$-varieties $X$. Fix such an $X$. By taking the $H$ quotient of $H \times^{G} \Gamma \times X$ we get an $H$-variety which, upon pulling back to $\GVar$ as in Lemma \ref{Lemma: Section 3: Obvious Pullback functor for restriction of action} realizes $\quot{X}{H \times^{G} \Gamma}$ as a $G$-variety with a trivial $G$-action. In particular, it is routine to check that there is a trivial $G$-morphism $\Gamma \times X \to \quot{X}{H \times^{G} \Gamma}$ induced by the set-theoretic argument
	\[
	(\gamma,x) \mapsto \big[[1_H,\gamma]_G, x\big]_H;
	\]
	note that because the $H$-action on $H \times^{G} \Gamma$ is given by $h[h^{\prime},\gamma]_G = [hh^{\prime},\gamma]_G$ and $[h,g\gamma] = [h\varphi(g)^{-1},\gamma]_G$, we have that this well-defines a morphism. In particular, because $\XGamma$ is a categorical quotient, this implies that there is a unique morphism $\zeta:\XGamma \to \quot{X}{\quot{(H \times \Gamma)}{G}}$ making
	\[
	\xymatrix{
		\Gamma \times X \ar[rr] \ar[dr]_{q} & & \quot{X}{H \times^{G} \Gamma} \\
		& \XGamma \ar@{-->}[ur]_{\exists!\zeta}
	}
	\]
	commute.
	
	We now construct a $\theta:\quot{X}{H \times^{G} \Gamma} \to \XGamma$ factoring the diagram:
	\[
	\xymatrix{
		\Gamma \times X \ar[rr]^-{q} \ar[dr] & & \XGamma \\
		& \quot{X}{H \times^{G} \Gamma} \ar@{-->}[ur]_{\exists\,\theta}
	}
	\]
	To construct this morphism, observe that the assignment
	\[
	\big[[h,\gamma]_G,x\big]_H \mapsto [\gamma,x]_G
	\]
	well-defines a morphism between $G$-varieties. Moreover, it is straightforward to check that the longer path of the diagram sends $(\gamma,x)$ to
	\[
	\theta\big[[1_H,\gamma]_G,x\big]_H = [\gamma,x]_G = q(\gamma,x)
	\]
	so the diagram does indeed commute. Moreover, it is also routine to check from this construction that from the universal property of $\XGamma$, we have that
	\[
	\theta \circ \zeta = \id_{\XGamma}.
	\]
	Thus to prove that $\XGamma \cong \quot{X}{H \times^{G} \Gamma},$ we only need to verify the other direction of the isomorphism.
	
	To verify that $\zeta \circ \theta = \id_{\quot{X}{H \times^G \Gamma}}$, we consider that the functor $\varphi^{\ast}$ of Lemma \ref{Lemma: Section 3: Obvious Pullback functor for restriction of action} induces the commuting diagram
	\[
	\xymatrix{
		\Gamma \times X \ar[dr] \ar[d] \\
		\left(H \times^{G} \Gamma\right) \times X \ar[r] & \quot{X}{H \times^{G} \Gamma}
	}
	\]
	as the map $(\gamma, x) \mapsto ([1_H,\gamma]_G,x)$ is $G$-equivariant from the action
	\[
	g ([1_H,\gamma]_G,x) = (g[1_H,\gamma]_G,gx) = ([1_H,\gamma]_G,\varphi(g)x).
	\]
	This allows us to note that upon taking $G$-orbits of the above action that $\zeta[\gamma,x] = \big[[h,\gamma]_G,x\big]_H$ is well-defined (as permuting the $g$-action and using how the quotient $[h,\gamma]_G$ is defined allows us manipulate the shift by an appropriate $G$-multiplication $x$ and then cancel via a quotient). Thus we have that
	\[
	(\zeta \circ \theta)\big[[h,\gamma]_G,x\big] = \zeta[\gamma,x]_G = \big[[h,\gamma]_G,x\big]
	\]
	so we get that $\zeta \circ \theta = \id_{\quot{X}{H \times^{G} \Gamma}}$. Thus $\XGamma \cong \quot{X}{H \times^{G} \Gamma}$, as desired.
\end{proof}
\begin{lemma}\label{Lemma: Section 3: Induction space quotient is iso to the Free(G) quotient}
	For any $M \in \Free(G)_M$ and any left $H$-space $X$, there is an isomorphism of quotient spaces
	\[
	G\backslash(M \times X)= \quot{X}{M} \cong \quot{X}{H \times^G M} = H \backslash \big((H \times^G M) \times X\big).
	\]
\end{lemma}

We finally have the tools at hand to give the change of group functors. Let us describe our technique in the geometric case; the topological situation follows by the same techniques by way of replacing every variety with a free $G$-space. When given a morphism $\varphi:G \to H$ of smooth algebraic groups, a left $H$-variety $X$, and an $F_H(X)$-object $A$, we can restrict the collection $\lbrace \AGamma \; | \; \Gamma \in \Sf(H)_0 \rbrace$ to the subcollection 
\[
\lbrace \quot{A}{H \times^{G} \Gamma^{\prime}} \; | \; \Gamma^{\prime} \in \Sf(G)_0 \rbrace
\] 
varying through $\Sf(G)$ taking values in the quotient varieties $H \times^{G} \Gamma^{\prime}$. However, to make this into an honest functor to $F_G(X)$, we need each of the objects to take values in $\overline{F}(\XGamma)$ and not the induction spaces $\quot{X}{H \times^{G} \Gamma}$. For this, however, we use Lemma \ref{Lemma: Section 3: Induction space quotient is iso to the Sf(G) quotient} and apply the equivalence of categories
\[
\overline{F}(h_{\Gamma}):\overline{F}(\quot{X}{H \times^{G} \Gamma}) \xrightarrow{\simeq} \overline{F}(\XGamma)
\] 
induced by the isomorphism $h_{\Gamma}:\XGamma \xrightarrow{\cong} \quot{X}{H \times^{G} \Gamma}$. That this does indeed define a functor is essentially due to Corollary \ref{Cor: Section 3: Functor from Sf(G) to Sf(H) by taking Bernstein quotients} and Lemma \ref{Lemma: Section 3: Induction space quotient is iso to the Sf(G) quotient}, as the two of these together say that we can pass through the equivalences 
\[
\overline{F}(\XGamma) \simeq \overline{F}(\quot{X}{H \times^G \Gamma})
\] 
induced by the isomorphisms $\XGamma \cong \quot{X}{H \times^{G} \Gamma}$ in $\Var_{/K}$ to get our functor.
\begin{Theorem}\label{Thm: Section 3: Change of groups functor}
	Let $G, H,$ $\varphi:G \to H$, $X$, and $F$ be given as in either of the cases below:
	\begin{itemize}
		\item $G$ and $H$ are smooth algebraic groups, $\varphi$ is a morphism of algebraic groups, $X$ is a left $H$-variety, and $F$ is an $H$-pre-equivariant pseudofunctor;
		\item $G$ and $H$ are topological groups, $\varphi$ is a moprhism of topological groups, $X$ is a left $H$-space, and $F$ is an $H$-pre-equivariant pseudofunctor.
	\end{itemize} 
	Then in both cases there is a functor $\varphi^{\sharp}:F_H(X) \to F_G(X)$\index[notation]{PhiSharp@$\varphi^{\sharp}$} given by sending an object $(A,T_A)$ to the pair $(A^{\prime},T_{A^{\prime}})$ and a morphism $P$ to $P^{\prime}$, where, in the geometric case, the pair $(A^{\prime}, T_{A^{\prime}})$ is induced by the equation
	\[
	A^{\prime} := \left\lbrace \overline{F}(h_{\Gamma})\left(\quot{A}{H \times^{G} \Gamma^{\prime}}\right) \; : \; \Gamma^{\prime} \in \Sf(G)_0, \quot{A}{H \times^{G} \Gamma^{\prime}} \in A \right\rbrace
	\]
	and on $P^{\prime}$ is induced by
	\[
	P^{\prime} := \left\lbrace \overline{F}(h_{\Gamma})\left(\quot{\rho}{H \times^{G} \Gamma^{\prime}}\right) \; : \; \Gamma^{\prime} \in \Sf(G)_0, \quot{\rho}{H \times^{G} \Gamma^{\prime}} \in P  \right\rbrace,
	\]
	where $h_{\Gamma}$ is the isomorphism $\XGamma \xrightarrow{\cong} \quot{X}{H \times^{G} \Gamma}$ of Lemma \ref{Lemma: Section 3: Induction space quotient is iso to the Sf(G) quotient}. The topological case is given mutatis mutandis.
\end{Theorem}
\begin{proof}
	By Corollary \ref{Cor: Section 3: Functor from Sf(G) to Sf(H) by taking Bernstein quotients} and Lemma \ref{Lemma: Section 3: Induction space quotient is iso to the Sf(G) quotient} (respectively Corollary \ref{Cor: Section 3: Top Functor from Sf(G) to Sf(H) by taking Bernstein quotients} and Lemma \ref{Lemma: Section 3: Induction space quotient is iso to the Free(G) quotient} in the topological case) there are functorial isomorphisms
	\[
	\XGamma \to \quot{X}{H \times^{G} \Gamma}
	\]
	for all $\Gamma \in \Sf(G)_0$. Let $h_{\Gamma}:\XGamma \to \quot{X}{H\times^G\Gamma}$ be the induced isomorphism of varieties and consider the induced equivalence of categories
	\[
	\overline{F}(h_{\Gamma}):\overline{F}(\quot{X}{H \times^{G} \Gamma}) \xrightarrow{\simeq} \overline{F}(\XGamma).
	\]
	Now observe that the diagram of schemes
	\[
	\xymatrix{
		\XGamma \ar[dr]^{k_f} \ar[r]^-{h_{\Gamma}} \ar[d]_{\overline{f}} & \quot{X}{H \times^{G} \Gamma} \ar[d]^{\overline{\id_H \times^G f}} \\
		\XGammap \ar[r]_-{h_{\Gamma^{\prime}}} & \quot{X}{{(H \times \Gamma^{\prime})}_{G}}
	}
	\]
	commutes for all $f:\Gamma \to \Gamma^{\prime}$ in $\Sf(G)$. Set
	\[
	k_f := h_{\Gamma^{\prime}} \circ \overline{f}
	\]
	so that
	\[
	h_{\Gamma^{\prime}} \circ \overline{f} = k_f = \overline{\id_H \times^G f} \circ h_{\Gamma}.
	\]
	This in turn implies that if $(A,T_A) \in F_H(X)_0$ we must have that
	\[
	\tau_f^{\varphi^{\sharp}} = \overline{F}(h_{\Gamma})\left(\tau_{\overline{\id_H \times^G f}}^{A}\right) \circ \left(\phi_{h_{\Gamma},\overline{\id_H \times^{G} f}}\right)^{-1} \circ \phi_{\overline{f},h_{\Gamma^{\prime}}}.
	\]
	Checking that $\varphi^{\sharp}$ defines a functor is tedious but follows mutatis mutandis to the proof of Theorem \ref{Thm: Functor Section: Psuedonatural trans are pseudocone functors} with the role of the transformation $\quot{\alpha}{f}^{-1}$ played by $\phi_{h_{\Gamma},\overline{\id_H \times^{G} f}}^{-1} \circ \phi_{\overline{f},h_{\Gamma^{\prime}}}$ instead.
\end{proof}
\begin{remark}\label{Remark: Cat theory of change of groups}
	In category theoretic-terms the functor $\varphi^{\sharp}:F_H(X) \to F_G(X)$ gives a functor of pseudocones
	\[
	F_H(X) = \PC(\overline{F} \circ \quo_H^{\op}) \to \PC(\overline{F} \circ \quo_G^{\op}) = F_G(X)
	\]
	by identifying a subcategory of $\Var_{/K}$ in the image of $\quo_H$ which is equivalent to the image of $\quo_G$ and then passing $\overline{F}$ through this equivalence. We note also that by keeping track of the pseudonatural equivalence $e$ in the definition of a $(G,H)$-pre-equivariant pseudofunctor we can produce a Change of Groups functor
	\[
	\quot{\varphi^{\sharp}_F}{F^{\prime}}:F_H^{\prime}(X) \to F_G(X)
	\]
	by using the factorization:
	\[
	\begin{tikzcd}
		F^{\prime}_H(X) \ar[d, swap]{}{\quot{\varphi^{\sharp}_{F}}{F^{\prime}}} \ar[rrr, equals] & & & \PC\left(\overline{F}^{\prime} \circ \quo_{H}^{\op}\right) \ar[d]{}{\varphi^{\sharp}_{F^{\prime}}} \\
		F_G(X) & \PC\left(\overline{F} \circ \quo_G^{\op}\right) \ar[l, equals] & & \PC\left(\overline{F}^{\prime} \circ \quo_G^{\op}\right) \ar[ll]{}{e^{-1} \ast \quo_G^{\op}}
	\end{tikzcd}
	\]
\end{remark}
\begin{definition}\index[terminology]{Equivariant Functor! Change of Group}
	A functor $R:F^{\prime}_H(X) \to F_G(X)$ is a Change of Groups functor if it arises as in the diagram expressed in Remark \ref{Remark: Cat theory of change of groups}.
\end{definition}

Let us now proceed to study how the Change of Groups functors interact when we change fibres and when we change spaces. As we have seen from Theorem \ref{Thm: Functor Section: Psuedonatural trans are pseudocone functors} and its various consequences (cf.\@ Proposition \ref{Prop: Pseudocone Functors: TRanslation}, Theorem \ref{Thm: Pseudocone Functors: Pullback induced by fibre functors in pseudofucntor}, Corollaries \ref{Cor: Pseudocone Functors: Existence of equivariant pullback  for schemes}, \ref{Cor: Pseudocone Functors: Existence of equivariant pullback  for spaces}, \ref{Cor: Pseudocone Functors: Exceptional Pushforward functors for equivariant maps for scheme sheaves}, etc.\@) a functor $F_G(X) \to F_G(Y)$ is  encoded as a Change of Domain functor and hence is not-so-secretly a Change of Fibre functor in the sense of Definition \ref{Defn: Functors Section: Fibre Functor}. As such, to study how Change of Groups functors interact with both changing fibres and spaces (in the sense of studying interactions between functors $F_H(X) \to F_G(X) \to E_G(X)$ as opposed to $F_H(X) \to E_H(X) \to E_G(X)$ or $F_H(X) \to F_H(Y) \to F_G(Y)$ as opposed to $F_H(X) \to F_G(X) \to F_G(Y)$), it suffices to give a careful study of how Change of Groups functors interact with pseudonatural transformations of $(G,H)$-pre-equivariant pseudofunctors. Note that as we proceed we describe the geometric situation; the topological case is analogous and omitted for space reasons.

To set the stage, let $G$ and $H$ be smooth algebraic groups, $X$ a left $H$-variety, and let $\varphi:G \to H$ is a morphism of algebraic groups. Fix $(G,H)$-pre-equivariant pseudofunctors $\left((F,\overline{F}), (F^{\prime}, \overline{F}^{\prime}), e:\overline{F}\Rightarrow \overline{F}^{\prime}\right)$ and $\left((E,\overline{E}), (E^{\prime}, \overline{E}^{\prime}), e^{\prime}:\overline{E}\Rightarrow\overline{E}^{\prime}\right)$. To have functors $\ul{\alpha}:F_G(X) \to E_G(X)$ and $F_H(X) \to E_H(X)$ defined simultaneously it suffices to have pseudonatural transformations
\[
\begin{tikzcd}
	\Var_{/K}^{\op} \ar[rr, bend left = 30, ""{name = R}]{}{\overline{F}} \ar[rr, bend right = 30, swap, ""{name = L}]{}{\overline{E}} & & \fCat \ar[from = R, to = L, Rightarrow, shorten <= 4pt, shorten >= 4pt]{}{\alpha}
\end{tikzcd}
\]
and
\[
\begin{tikzcd}
	\Var_{/K}^{\op} \ar[rr, bend left = 30, ""{name = R}]{}{\overline{F}^{\prime}} \ar[rr, bend right = 30, swap, ""{name = L}]{}{\overline{E}^{\prime}} & & \fCat \ar[from = R, to = L, Rightarrow, shorten <= 4pt, shorten >= 4pt]{}{\alpha^{\prime}}
\end{tikzcd}
\]
which are compatible with the equivalences $e$ and $e^{\prime}$  in the sense that the equation
\[
e^{\prime} \circ \alpha = \alpha^{\prime} \circ e
\]
and its inverse version
\[
\alpha \circ e^{-1} = (e^{\prime})^{-1} \circ \alpha^{\prime}
\]
both hold. Note that we have made the convention that $\alpha:\overline{F} \Rightarrow \overline{E}$ and $\alpha^{\prime}:\overline{F}^{\prime} \Rightarrow \overline{E}^{\prime}$ for notational convenience; if we instead have a transformation $\alpha:\overline{F} \Rightarrow \overline{E}^{\prime}$ we can simply use the inverse equivalence $(e^{\prime})^{-1}$ of $e^{\prime}$ and proceed as we did above using the inverse instead of $e^{\prime}$. Now observe that the assumptions above, together with the pseudonaturality of the equivalences $e, e^{\prime},$ and their inverses, give rise to the commuting diagram of pseudocone categories
\[
\begin{tikzcd}
	F^{\prime}_H(X) \ar[d, equals] \ar[rr]{}{\ul{\alpha}^{\prime}} & & E^{\prime}_H(X) \ar[d, equals]{}{} \\
	\PC\left(\overline{F}^{\prime} \circ \quo_H^{\op}\right)\ar[rr]{}{\ul{\alpha^{\prime} \ast \quo_{H}^{\op}}} \ar[d, swap]{}{\varphi^{\sharp}_{F^{\prime}}} \ar[dd, bend right = 80, swap]{}{\quot{\varphi^{\sharp}_F}{F^{\prime}}} & & \PC\left(\overline{E}^{\prime} \circ \quo_H^{\op}\right) \ar[d]{}{\varphi_{E^{\prime}}^{\sharp}} \ar[dd, bend left = 80]{}{\quot{\varphi^{\sharp}_E}{E^{\prime}}} \\
	\PC\left(\overline{F}^{\prime} \circ \quo_G^{\op}\right) \ar[d, swap]{}{\ul{e^{-1} \ast \quo_G^{\op}}}  & &  \PC\left(\overline{E}^{\prime} \circ \quo_G^{\op}\right) \ar[d]{}{\ul{(e^{\prime})^{-1} \ast \quo_G^{\op}}} \\
	\PC\left(\overline{F} \circ \quo_G^{\op}\right) \ar[d, equals] \ar[rr]{}{\ul{\alpha \ast \quo_G^{\op}}} & & \PC\left(\overline{E} \circ \quo_G^{\op}\right) \ar[d, equals] \\
	F_G(X) \ar[rr, swap]{}{\ul{\alpha}} & &E_G(X)
\end{tikzcd}
\]
by repeated use of Theorem \ref{Thm: Functor Section: Psuedonatural trans are pseudocone functors}. However, this implies that the diagram
\[
\begin{tikzcd}
	F^{\prime}_G(X) \ar[r]{}{\ul{\alpha}^{\prime}} \ar[d, swap]{}{\quot{\varphi^{\sharp}_{F}}{F^{\prime}}} \ar[r]{}{\ul{\alpha}^{\prime}} & E^{\prime}_H(X) \ar[d]{}{\quot{\varphi^{\sharp}_{E}}{E^{\prime}}} \\
	F_G(X) \ar[r, swap]{}{\ul{\alpha}} & E_G(X)
\end{tikzcd}
\]
commutes, which in shows that all we need to ensure that we can change groups functorially over a $(G,H)$-pre-equivariant pseudofunctor is to have pseudonatural transformations compatible with the equivalences in the $(G,H)$-pre-equviariant pseudofunctors. We record this as the proposition below.
\begin{proposition}\label{Prop: Section Chofg: Change of groups together with pullbacks}
	To set the stage, let $G$ and $H$ be smooth algebraic groups, $X$ a left $H$-variety, and let $\varphi:G \to H$ is a morphism of algebraic groups. Let $\ul{F} = \left((F,\overline{F}), (F^{\prime}, \overline{F}^{\prime}), e:\overline{F}\Rightarrow \overline{F}^{\prime}\right)$ and $\ul{E} = \left((E,\overline{E}), (E^{\prime}, \overline{E}^{\prime}), e^{\prime}:\overline{E}\Rightarrow\overline{E}^{\prime}\right)$ be $(G,H)$-pre-equivariant pseudofunctors  and assume that there are pseudonatural transformations $\alpha:\overline{F} \Rightarrow \overline{E}$ and $\alpha^{\prime}:\overline{F}^{\prime} \Rightarrow \overline{E}^{\prime}$ which commute with the equivalences $e, e^{\prime}$, and their inverses. Then the diagram
	\[
	\begin{tikzcd}
		F^{\prime}_G(X) \ar[r]{}{\ul{\alpha}^{\prime}} \ar[d, swap]{}{\quot{\varphi^{\sharp}_{F}}{F^{\prime}}} \ar[r]{}{\ul{\alpha}^{\prime}} & E^{\prime}_H(X) \ar[d]{}{\quot{\varphi^{\sharp}_{E}}{E^{\prime}}} \\
		F_G(X) \ar[r, swap]{}{\ul{\alpha}} & E_G(X)
	\end{tikzcd}
	\]
	commutes.
\end{proposition}
\begin{remark}
	This proposition captures each of \cite[Propositions 4.3.13, 4.3.14, 4.3.15]{MyThesis} by virtue of expressing Change of Domain functors in terms of Change of Fibre functors and because the $(G,H)$-pre-equivariant pseudofunctors which appear in \cite{MyThesis} only handle the case when both equivalences for $\ul{F}$ and $\ul{E}$ satisfy $e = \id$ and $e^{\prime} = \id$.
\end{remark}

We now discuss how Change of Group functors interact with other Change of Group functors \footnote{While we motivate and express our setting in the geometric language, we can translate to the topological case freely by virtue of the same techniques we have used so far. More explicitly, while we motivate the Change of Groups interaction in terms of varieties and algebraic groups, by simply replacing algebraic groups with topological groups and smooth free varieties with free $G$ (respectively $H$)-spaces, we get the topological argument with what is ultimately a ``find and replace'' style of argument.}. In particular, if we have algebraic group morphisms $\varphi\in\AlgGrp(G_0,G_1)$ and $\psi \in \AlgGrp(G_1, G_2)$ for smooth algebraic groups $G_0,$ $G_1,$ $G_2$, we would like to study how the composites $\varphi^{\sharp} \circ \psi^{\sharp}$ and $(\psi \circ \varphi)^{\sharp}$ compare. For this we will need a pre-equivariant pseudofunctor on a $G_2$-variety $X$ for both pairs $(G_0, G_1)$ and $(G_1, G_2)$; let us denote this by a $(G_0, G_1, G_2)$-pre-equivariant pseudofunctor on $X$. There are then two distinct Change of Groups functors from $F_{G_2}(X) \to F_{G_0}(X)$: first, the total composite
\[
\begin{tikzcd}
	F_{G_2}(X) \ar[r]{}{\psi^{\sharp}} & F_{G_1}(X) \ar[r]{}{\varphi^{\sharp}} & F_{G_0}(X)
\end{tikzcd}
\]
and second the morphism induced in one step by $\psi \circ \varphi$:
\[
\begin{tikzcd}
	F_{G_2}(X) \ar[r]{}{(\psi \circ \varphi)^{\sharp}} & F_{G_0}(X)
\end{tikzcd}
\]
While in general these two functors should not be expected to be the same, in the proposition below we will prove that the are isomorphic to each other by a natural isomorphism induced simultaneously from the transition isomorphisms in $F_{G_2}(X)$ that an object must satisfy and the compositor isomorphisms that the pseudofunctor $F$ has at hand. We start by giving a very concrete and explicit proof of the case when our $(G_0,G_1,G_2)$-pre-equivariant pseudofunctors have witness equivalence equal to the identity equivalence.

\begin{proposition}\label{Prop: Section 3.3: Chagnge of groups interacting with change of grups}
	Let $G_0, G_1, G_2, X, F, \varphi,$ and $\psi$ be given as in either case below:
	\begin{itemize}
		\item $G_0, G_1, G_2$ are smooth algebraic groups\footnote{It is unfortunate that there is a clash of notation, but in this case $G_2$ is meant only as the third entry in a list of algebraic groups and \emph{not} as the exceptional Lie group over $\Spec K$.}, $X$ is a left $G_2$-variety, $F$ is a $G_2$-pre-equivariant pseudofunctor on $X$, and both $\varphi:G_0 \to G_1$ and $\psi:G_1 \to G_2$ are morphisms of algebraic groups.
		\item $G_0, G_1, G_2$ are topological groups, $X$ is a left $G_2$-space, $F$ is a $G_2$-pre-equivariant pseudofunctor on $X$, and both $\varphi:G_0 \to G_1$ and $\psi:G_1 \to G_2$ are morphisms of topological groups.
	\end{itemize}
	Then there is an invertible $2$-cell
	\[
	\begin{tikzcd}
		F_{G_2}(X) \ar[dr, swap]{}{(\psi \circ \varphi)^{\sharp}} \ar[rr]{}{\psi^{\sharp}} & {} &F_{G_1}(X) \ar[dl, swap]{}{\varphi^{\sharp}} \\
		& F_{G_0}(X) \ar[from = 2-2, to = 1-2, Rightarrow, shorten >= 4pt, shorten <= 4pt]{}{\alpha}
	\end{tikzcd}
	\]
	in $\fCat$.
\end{proposition}
\begin{proof}
As usual we give the argument for the geometric case; the topological case follows similarly. Begin by observing that for any variety $\Gamma$ in $\Sf(G_0)$, it is straightforward but tedious to prove that there is a $G_2$-equivariant isomorphism of varieties
	\[
	G_2 \times^{G_1} \left(G_1 \times^{G_0} \Gamma\right) \xrightarrow[\rho]{\cong} G_2 \times^{G_0} \Gamma;
	\]
	in a set-theoretic sketch, this isomorphism is given by $[g_2, [g_1,\gamma]_{G_0}]_{G_1} \mapsto [g_2,\gamma]_{G_0}$ and that this is $G_2$-equivariant is immediate because the $G_2$-action occurs only in the $G_2$-component. Now, because each of the above varieties are $\Sf(G_2)$-objects, by three applications of Proposition \ref{Prop: Section  3: Bernstein Quotient is Sf(H)} (one for $G_2 \times^{G_0}\Gamma$ directly, one for realizing $G_1 \times^{G_0} \Gamma$ as an $\Sf(G_1)$-variety, and then one for the case of realizing $(G_2 \times^{G_1} (G_1 \times^{G_0} \Gamma))$ as an $\Sf(G_2)$-variety) it follows from the fact that isomorphisms are smooth (and in fact zero-dimensional {\'e}tale with constant fibres) that $\rho$ is a morphism in $\Sf(G_2)$.
	
	We calculate our two composites so as to begin our comparison. First fix an object $A$ in $F_{G_2}(X)$. Now observe that on one hand
	\begin{align*}
		\varphi^{\sharp}(\psi^{\sharp}A) &= \varphi^{\sharp}\left(\left\lbrace \overline{F}\left(h^{G_1G_2}_{\Gamma^{\prime}}\right)\left(\quot{A}{(G_2 \times \Gamma^{\prime})_{G_1}}\right) \; : \; \Gamma^{\prime} \in \Sf(G_1)_0 \right\rbrace\right) \\
		&= \left\lbrace \overline{F}\left(h_{\Gamma}^{G_0G_1}\right)\left(\overline{F}\left(h^{G_1G_2}_{G_1 \times^{G_0} \Gamma}\right)\left(\quot{A}{G_2 \times^{G_1} \left(G_1 \times^{G_1} \Gamma\right)}\right)\right) \; : \; \Gamma \in \Sf(G_0)_0 \right\rbrace
	\end{align*}
	while on the other hand
	\[
	(\psi \circ \varphi)^{\sharp}A = \left\lbrace \overline{F}\left(h_{\Gamma}^{G_0G_2}\right)\left(\quot{A}{G_2 \times^{G_0} \Gamma}\right) \; : \; \Gamma \in \Sf(G)_0 \right\rbrace.
	\]
	To compare these two constructions, we note that it is routine to verify that the diagram
	\[
	\xymatrix{
		\quot{X}{\Gamma} \ar[d]_{h_{\Gamma}^{G_0G_2}} \ar[rr]^-{h_{\Gamma}^{G_0G_1}} & & \quot{X}{G_1 \times^{G_0} \Gamma} \ar[d]^{h_{G_1 \times^{G_0} \Gamma}^{G_1G_2}} \\
		\quot{X}{G_2 \times^{G_0} \Gamma} & & \quot{X}{G_2 \times^{G_1} \left(G_1 \times^{G_1} \Gamma\right)} \ar[ll]^-{\overline{\rho}}
	}
	\]
	commutes by using the universal properties of quotients. With this we produce the pasting diagram:
	\[
	\begin{tikzcd}
		& \overline{F}\left(\quot{X}{G_1 \times^{G_0} \Gamma}\right) \ar[dr, bend left = 30]{}{\overline{F}(h_{\Gamma}^{G_0G_1})} & \\
		\overline{F}\left(\quot{X}{G_2 \times^{G_1} \left(G_1 \times^{G_1} \Gamma\right)}\right) \ar[rr]{}[description]{\overline{F}(h_{G_1 \times^{G_0} \Gamma} \circ h_{\Gamma})} \ar[ur, bend left = 30]{}{\overline{F}\left(h_{G_1 \times^{G_0} \Gamma}^{G_1G_2}\right)} & {} & \overline{F}\left(\quot{X}{\Gamma}\right) \\
		& \overline{F}\left(\quot{X}{G_2 \times^{G_0} \Gamma}\right) \ar[ur, swap, bend right = 30]{}{\overline{F}\left(h_{\Gamma}^{G_0G_2}\right)} \ar[ul, bend left = 30]{}{\overline{F}(\overline{\rho})} & \ar[from = 2-2, to = 3-2, swap, Rightarrow, shorten >= 4pt, shorten <= 4pt]{}{\phi_{h^{G_1G_2} \circ h^{G_0G_1}, \rho}} \ar[from = 1-2, to = 2-2, Rightarrow, shorten <= 4pt, shorten >=4pt, swap]{}{\phi_{h^{G_0G_1}, h^{G_1G_2}}}
	\end{tikzcd}
	\]
	From the fact that $\rho$ is a morphism in $\Sf(G_2)$, we have a transition isomorphism
	\[
	\tau_{\rho}^{A}:\overline{F}(\overline{\rho})\left(\quot{A}{G_2 \times^{G_0} \Gamma}\right) \xrightarrow{\cong} \quot{A}{G_2 \times^{G_1} \left(G_1 \times^{G_1} \Gamma\right)}.
	\]
	Applying the functor $\overline{F}(h_{\Gamma}^{G_0G_1}) \circ \overline{F}(h_{G_1 \times^{G_0} \Gamma}^{G_1G_2})$ to $\tau_{\rho}^{A}$ above gives an isomorphism:
	\[
	\begin{tikzcd}
		\overline{F}\left(h_{\Gamma}^{G_0G_1}\right)\left(\overline{F}\left(h_{G_1 \times^{G_0} \Gamma}^{G_1G_2}\right)\left(\overline{F}(\overline{\rho})\left(\quot{A}{G_2 \times^{G_0} \Gamma}\right)\right)\right) \ar[d]{}{\cong} \\
		\overline{F}\left(h_{\Gamma}^{G_0G_1}\right)\left(\overline{F}\left(h_{G_1 \times^{G_0} \Gamma}^{G_1G_2}\right)\left(\quot{A}{G_2 \times^{G_1} \left(G_1 \times^{G_1} \Gamma\right)}\right)\right)
	\end{tikzcd}
	\]
	Pre-composing this with the inverses of the two compositors in the pasting diagram then gives an isomorphism
	\[
	\quot{\alpha_A}{\Gamma} := \overline{F}\left(h_{\Gamma}^{G_0G_1}\right)\left(\overline{F}\left(h_{G_1 \times^{G_0} \Gamma}^{G_1G_2}\right)(\tau_{\rho}^{A})\right) \circ \left(\phi_{h^{G_0G_1}, h^{G_1G_2}}^{\overline{F}(\overline{\rho})A}\right)^{-1} \circ \left(\phi_{h^{G_1G_2} \circ h^{G_0G_1}, \rho}^{A}\right)^{-1}
	\]
	between the objects:
	\[
	\begin{tikzcd}
		\quot{(\psi \circ \varphi)^{\sharp}(A)}{\Gamma} \ar[r, equals] \ar[d, swap]{}{\cong} & \overline{F}\left(h_{\Gamma}^{G_0G_2}\right)\left(\quot{A}{G_2 \times^{G_0} \Gamma}\right) \ar[d]{}{\cong} \\ \quot{(\varphi^{\sharp}(\psi^{\sharp}A))}{\Gamma} \ar[r, equals] & \overline{F}\left(h_{\Gamma}^{G_0G_1}\right)\left(\overline{F}\left(h_{G_1 \times^{G_0} \Gamma}^{G_1G_2}\right)\left(\quot{A}{G_2 \times^{G_1} \left(G_1 \times^{G_1} \Gamma\right)}\right)\right)
	\end{tikzcd}
	\]
	Define the map $\alpha_A: (\psi \circ \varphi)^{\sharp}A \to (\varphi^{\sharp} \circ \psi^{\sharp})A$ by setting the $\Gamma$-local components as the morphisms $\quot{\alpha}{\Gamma}$ given above.
	
	We now show that $\alpha$, as defined, is a morphism in $F_{G_0}(X)$ from $(\psi \circ \varphi)^{\sharp}A$ to $(\varphi^{\sharp} \circ \psi^{\sharp})A$. For this and for the rest of the proof we set some notation to reduce the already significant notational clutter. Make the following definitions of convenience:
	\begin{align*}
		\Gamma_1 &:= G_1 \times^{G_0} \Gamma & \Gamma_2 &:= G_2 \times^{G_1} \Gamma_1 & \Gamma_3 &:= G_2 \times^{G_0} \Gamma \\
		\Gamma^{\prime}_{1} &:= G_1 \times^{G_0} \Gamma^{\prime} & \Gamma^{\prime}_2 &:= G_2 \times^{G_1} \Gamma_{1}^{\prime} & \Gamma_{3}^{\prime} &:= G_2 \times^{G_0} \Gamma \\
		f_1 &:= \id_{G_1} \times^{G_0} f & f_2 & := \id_{G_2} \times^{G_1} f_1 & f_3 &:= \id_{G_2} \times^{G_0} f \\
		h_{\Gamma}^{01} &:= h_{\Gamma}^{G_0G_1} & h_{\Gamma_1}^{12} &:= h_{\Gamma_1}^{G_1G_2} & h_{\Gamma}^{02} &:= h_{\Gamma}^{G_0G_2} \\
		h_{\Gamma^{\prime}}^{01} &:= h_{\Gamma^{\prime}}^{G_0G_1} & h_{\Gamma_1^{\prime}}^{12} &:= h_{\Gamma_1^{\prime}}^{G_1G_2} & h_{\Gamma^{\prime}}^{02} &:= h_{\Gamma^{\prime}}^{G_0G_2}
	\end{align*}
	To prove that $\alpha_A$ is a morphism in $F_{G_0}(X)$, we must show that the diagram
	\[
	\begin{tikzcd}
		\overline{F}(\of)\left(\quot{(\psi \circ \varphi)^{\sharp}A}{\Gamma^{\prime}}\right) \ar[rr]{}{\overline{F}(\of)\left(\quot{\alpha_A}{\Gamma^{\prime}}\right)} \ar[d, swap]{}{\tau_f^{(\psi \circ \varphi)^{\sharp}A}} & & \overline{F}(\of)\left(\quot{\varphi^{\sharp}(\psi^{\sharp}A)}{\Gamma^{\prime}}\right) \ar[d]{}{\tau_f^{\varphi^{\sharp}(\psi^{\sharp}A)}} \\
		\quot{(\psi \circ \varphi)^{\sharp}A}{\Gamma} \ar[rr, swap]{}{\quot{\alpha}{\Gamma}} & & \quot{\varphi^{\sharp}(\psi^{\sharp}A)}{\Gamma}
	\end{tikzcd}
	\]
	commutes for any morphism $f:\Gamma \to \Gamma^{\prime}$ in $\Sf(G_0)$. To establish this we first calculate using Theorem \ref{Thm: Section 3: Change of groups functor} that on one hand
	\[
	\tau_f^{(\psi \circ \varphi)^{\sharp}A} = \overline{F}(h_{\Gamma}^{02})\tau_{f_3}^{A} \circ \left(\phi_{h_{\Gamma}^{02}, f_3}^{A}\right)^{-1} \circ \phi_{f,h_{\Gamma^{\prime}}^{02}}^{A}
	\]
	while on the other hand
	\begin{align*}
		\tau^{\varphi^{\sharp}(\psi^{\sharp}A)}_f &= \overline{F}(h_{\Gamma}^{01})(\tau_{f_1}^{\psi^{\sharp}A}) \circ \left(\phi_{h_{\Gamma}^{01},f_1}^{\psi^{\sharp}A}\right)^{-1} \circ \phi_{f,h_{\Gamma^{\prime}}^{01}}^{\psi^{\sharp}A} \\
		&= \overline{F}(h_{\Gamma}^{01})\left(\overline{F}(h_{\Gamma_1}^{12})(\tau_{f_2}^{A}) \circ \left(\phi_{h_{\Gamma_1},f_2}^{A}\right)^{-1} \circ \phi_{f_1, h_{\Gamma^{\prime}_1}}^{A}\right) \circ \left(\phi_{h_{\Gamma}^{01},f_1}^{\psi^{\sharp}A}\right)^{-1} \circ \phi_{f,h_{\Gamma^{\prime}}^{01}}^{\psi^{\sharp}A}.
	\end{align*}
	In order to prove the commutativity of the above diagram, we will break these morphisms into pieces and manipulate them by starting with $\tau_f^{\varphi^{\sharp}(\psi^{\sharp}A)} \circ \overline{F}(\of)\quot{\alpha_A}{\Gamma^{\prime}}$ (or a piece of it, anyway). For this it will be particularly helpful to keep the commuting diagram of varieties
	\[
	\begin{tikzcd}
		\quot{X}{\Gamma} \ar[rrr, bend left = 30]{}{h_{\Gamma}^{02}} \ar[r]{}{h_{\Gamma}^{01}} \ar[d, swap]{}{\of} & \quot{X}{\Gamma_1} \ar[d, swap]{}{\of_1} \ar[r]{}{h_{\Gamma_1}^{12}} & \quot{X}{\Gamma_2} \ar[d]{}{\of_2} \ar[r]{}{\rho} & \quot{X}{\Gamma_3} \ar[d]{}{\of_3} \\
		\quot{X}{\Gamma^{\prime}} \ar[rrr, bend right = 30, swap]{}{h_{\Gamma^{\prime}}^{02}} \ar[r, swap]{}{h_{\Gamma^{\prime}}^{01}} & \quot{X}{\Gamma_1^{\prime}} \ar[r, swap]{}{h_{\Gamma_1^{\prime}}} & \quot{X}{\Gamma_2^{\prime}} \ar[r, swap]{}{\rho^{\prime}} & \quot{X}{\Gamma_3^{\prime}}
	\end{tikzcd}
	\] 
	in mind. Begin by observing that
	\begin{align*}
		&\left(\phi_{h_{\Gamma}^{01},f_1}^{\psi^{\sharp}A}\right)^{-1} \circ \phi_{f,h_{\Gamma^{\prime}}^{01}}^{\psi^{\sharp}A} \circ \overline{F}(\of)\quot{\alpha_A}{\Gamma^{\prime}} \\
		&= \left(\phi_{h_{\Gamma}^{01},f_1}^{\overline{F}(h_{\Gamma_1})\quot{A}{\Gamma_2}}\right)^{-1} \circ \phi_{f,h_{\Gamma^{\prime}}^{01}}^{\overline{F}(h_{\Gamma_1})\quot{A}{\Gamma_2}}  \\
		&\circ \overline{F}(\of)\left(\left(\overline{F}(h_{\Gamma^{\prime}}^{01}) \circ \overline{F}(h_{\Gamma_1^{\prime}}^{12})\right)\tau_{\rho^{\prime}}^{A} \circ \left(\phi_{h_{\Gamma^{\prime}}, h_{\Gamma^{\prime}_1}}^{\overline{F}(\overline{\rho}^{\prime})A}\right)^{-1} \circ \left(\phi_{h_{\Gamma_1^{\prime}} \circ h_{\Gamma^{\prime}}, \rho^{\prime}}^{A}\right)^{-1}\right) \\
		&= \left(\phi_{h_{\Gamma}^{01},f_1}^{\overline{F}(h_{\Gamma_1})\quot{A}{\Gamma_2}}\right)^{-1} \circ \phi_{f,h_{\Gamma^{\prime}}^{01}}^{\overline{F}(h_{\Gamma_1})\quot{A}{\Gamma_2}} \circ \left(\overline{F}(\of)\circ \overline{F}(h_{\Gamma^{\prime}}^{01}) \circ \overline{F}(h_{\Gamma_1^{\prime}}^{12})\right)\tau_{\rho^{\prime}}^{A} \\
		&\circ \overline{F}(\of)\left(\left(\phi_{h_{\Gamma^{\prime}}, h_{\Gamma^{\prime}_1}}^{\overline{F}(\overline{\rho}^{\prime})A}\right)^{-1} \circ \left(\phi_{h_{\Gamma_1^{\prime}} \circ h_{\Gamma^{\prime}}, \rho^{\prime}}^{A}\right)^{-1}\right) \\
		&= \left(\phi_{h_{\Gamma}^{01},f_1}^{\overline{F}(h_{\Gamma_1})\quot{A}{\Gamma_2}}\right)^{-1} \circ \left(\overline{F}(h_{\Gamma^{\prime}} \circ \of)\circ \overline{F}(h_{\Gamma^{\prime}_1}^{12})\right)\tau_{{\rho}^{\prime}}^{A} \circ \phi_{f, h_{\Gamma^{\prime}}}^{(\overline{F}(h_{\Gamma_1^{\prime}}) \circ \overline{F}(\overline{\rho}^{\prime}))A} \\
		&\circ \overline{F}(\of)\left(\left(\phi_{h_{\Gamma^{\prime}}, h_{\Gamma^{\prime}_1}}^{\overline{F}(\overline{\rho}^{\prime})A}\right)^{-1} \circ \left(\phi_{h_{\Gamma_1^{\prime}} \circ h_{\Gamma^{\prime}}, \rho^{\prime}}^{A}\right)^{-1}\right)  \\
		&= \left(\overline{F}(h_{\Gamma}^{01})\circ \overline{F}(\of_1)\circ \overline{F}(h_{\Gamma_{1}^{\prime}}^{12})\right)\tau_{\rho^{\prime}}^{A} \circ \left(\phi_{h_{\Gamma},f_1}^{\overline{F}(h_{\Gamma_1^{\prime}})\overline{F}(\overline{\rho}^{\prime})A}\right)^{-1} \circ \phi_{f, h_{\Gamma^{\prime}}}^{\overline{F}(h_{\Gamma_1^{\prime}})\overline{F}(\overline{\rho}^{\prime})A} \\
		&\circ \overline{F}(\of)\left(\left(\phi_{h_{\Gamma^{\prime}}, h_{\Gamma^{\prime}_1}}^{\overline{F}(\overline{\rho}^{\prime})A}\right)^{-1} \circ \left(\phi_{h_{\Gamma_1^{\prime}} \circ h_{\Gamma^{\prime}}, \rho^{\prime}}^{A}\right)^{-1}\right).
	\end{align*}
	We now calculate that
	\begin{align*}
		&\overline{F}(h_{\Gamma}^{01})\left(\overline{F}(h_{\Gamma_1}^{12})\tau_{f_2}^{A} \circ \left(\phi_{h_{\Gamma_1},f_2}^{A}\right)^{-1} \circ \phi_{f_1, h_{\Gamma_1}^{\prime}}^{A}\right) \circ \left(\overline{F}(h_{\Gamma}^{01})\circ \overline{F}(\of_1)\circ \overline{F}(h_{\Gamma_1^{\prime}}^{12})\right)\tau_{\rho^{\prime}}^{A} \\
		&= \overline{F}(h_{\Gamma}^{01})\left(\overline{F}(h_{\Gamma_1}^{12})\tau_{f_2}^{A} \circ \left(\phi_{h_{\Gamma_1},f_2}^{A}\right)^{-1} \circ \phi_{f_1, h_{\Gamma_1}^{\prime}}^{A} \circ \left(\overline{F}(\of_1)\circ \overline{F}(h_{\Gamma_1^{\prime}}^{12})\right)\tau_{\rho^{\prime}}^{A} \right) \\
		&= \overline{F}(h_{\Gamma}^{01})\left(\overline{F}(h_{\Gamma_1}^{12})\tau_{f_2}^{A} \circ \left(\phi_{h_{\Gamma_1},f_2}^{A}\right)^{-1} \circ \overline{F}(h_{\Gamma_1^{\prime}}^{12} \circ \of_1)\tau_{\rho^{\prime}}^{A} \circ \phi_{f_1, h_{\Gamma_1^{\prime}}}^{\overline{F}(\overline{\rho}^{\prime})A} \right) \\
		&= \overline{F}(h_{\Gamma}^{01})\left(\overline{F}(h_{\Gamma_1}^{12})\tau_{f_2}^{A} \circ \left(\overline{F}(h_{\Gamma_1})\circ \overline{F}(\of_2)\right)\tau_{\rho^{\prime}}^{A} \circ \left(\phi_{h_{\Gamma_1}, f_2}^{\overline{F}(\overline{\rho}^{\prime})A}\right)^{-1} \circ \phi_{f_1, h_{\Gamma_1^{\prime}}}^{\overline{F}(\overline{\rho}^{\prime})A}\right) \\
		&= \left(\overline{F}(h_{\Gamma}^{01}) \circ \overline{F}(h_{\Gamma_{1}}^{12})\right)\left(\tau_{f_2}^{A} \circ \overline{F}(\of_2)\tau_{\rho^{\prime}}^{A}\right) \circ \overline{F}(h_{\Gamma}^{01})\left(\left(\phi_{h_{\Gamma_1}, f_2}^{\overline{F}(\overline{\rho}^{\prime})A}\right)^{-1} \circ \phi_{f_1, h_{\Gamma_1^{\prime}}}^{\overline{F}(\overline{\rho}^{\prime})A}\right) \\
		&= \left(\overline{F}(h_{\Gamma}^{01})\circ \overline{F}(h_{\Gamma_1}^{12})\right)\left(\tau_{\rho^{\prime} \circ f_2}^{A} \circ \phi_{f_2,\rho^{\prime}}^{A}\right) \circ \overline{F}(h_{\Gamma}^{01})\left(\left(\phi_{h_{\Gamma_1}, f_2}^{\overline{F}(\overline{\rho}^{\prime})A}\right)^{-1} \circ \phi_{f_1, h_{\Gamma_1^{\prime}}}^{\overline{F}(\overline{\rho}^{\prime})A}\right) \\
		&= \left(\overline{F}(h_{\Gamma}^{01})\circ \overline{F}(h_{\Gamma_1}^{12})\right)\left(\tau_{f_3 \circ \rho}^{A} \circ \phi_{f_2,\rho^{\prime}}^{A}\right) \circ \overline{F}(h_{\Gamma}^{01})\left(\left(\phi_{h_{\Gamma_1}, f_2}^{\overline{F}(\overline{\rho}^{\prime})A}\right)^{-1} \circ \phi_{f_1, h_{\Gamma_1^{\prime}}}^{\overline{F}(\overline{\rho}^{\prime})A}\right) \\
		&= \left(\overline{F}(h_{\Gamma}^{01})\circ \overline{F}(h_{\Gamma_1}^{12})\right)\left(\tau_{\rho}^{A} \circ \overline{F}(\rho)\tau_{f_3}^{A} \circ \left(\phi_{\rho,f_3}^{A}\right)^{-1} \circ \phi_{f_2,\rho^{\prime}}^{A}\right) \\
		&\circ \overline{F}(h_{\Gamma}^{01})\left(\left(\phi_{h_{\Gamma_1}, f_2}^{\overline{F}(\overline{\rho}^{\prime})A}\right)^{-1} \circ \phi_{f_1, h_{\Gamma_1^{\prime}}}^{\overline{F}(\overline{\rho}^{\prime})A}\right) \\
		&= \left(\overline{F}(h_{\Gamma}^{01})\circ \overline{F}(h_{\Gamma_1}^{12})\right)\tau_{\rho}^{A} \circ  \left(\overline{F}(h_{\Gamma}^{01})\circ \overline{F}(h_{\Gamma_1}^{12})\right)\left(\overline{F}(\rho)\tau_{f_3}^{A} \circ \left(\phi_{\rho,f_3}^{A}\right)^{-1} \circ \phi_{f_2,\rho^{\prime}}^{A}\right) \\
		&\circ \overline{F}(h_{\Gamma}^{01})\left(\left(\phi_{h_{\Gamma_1}, f_2}^{\overline{F}(\overline{\rho}^{\prime})A}\right)^{-1} \circ \phi_{f_1, h_{\Gamma_1^{\prime}}}^{\overline{F}(\overline{\rho}^{\prime})A}\right).
	\end{align*}
	Computing that
	\begin{align*}
		&\quot{\alpha_A}{\Gamma} \circ \tau_f^{(\psi \circ \varphi)^{\sharp}A} \\
		&= \left(\overline{F}(h_{\Gamma}^{01})\circ \overline{F}(h_{\Gamma_1}^{12})\right)\tau_{\rho}^{A} \circ \left(\phi_{h^{01}_{\Gamma}, h_{\Gamma_1}}^{\overline{F}(\overline{\rho})A}\right)^{-1} \circ \left(\phi_{h^{12}_{\Gamma_1} \circ h^{01}_{\Gamma}, \rho}\right)^{-1} \circ  \overline{F}(h_{\Gamma}^{02})\tau_{f_3}^{A} \\
		&\circ \left(\phi_{h_{\Gamma}^{02}, f_3}^{A}\right)^{-1} \circ \phi_{f,h_{\Gamma^{\prime}}^{02}}^{A}
	\end{align*}
	and realizing that by combining our derivations above we have calculated and manipulated the morphism $\tau_f^{\varphi^{\sharp}(\psi^{\sharp}A)} \circ \overline{F}(\of)\quot{\alpha_A}{\Gamma^{\prime}}$, to prove that the two composites coincide it suffices to verify that the morphism
	\begin{align*}
		&\left(\overline{F}(h_{\Gamma}^{01})\circ \overline{F}(h_{\Gamma_1}^{12})\right)\left(\overline{F}(\rho)\tau_{f_3}^{A} \circ \left(\phi_{\rho,f_3}^{A}\right)^{-1} \circ \phi_{f_2,\rho^{\prime}}^{A}\right) \circ \overline{F}(h_{\Gamma}^{01})\left(\left(\phi_{h_{\Gamma_1}, f_2}^{\overline{F}(\overline{\rho}^{\prime})A}\right)^{-1} \circ \phi_{f_1, h_{\Gamma_1^{\prime}}}^{\overline{F}(\overline{\rho}^{\prime})A}\right) \\
		&\circ \overline{F}(\of)\left(\left(\phi_{h_{\Gamma^{\prime}}, h_{\Gamma^{\prime}_1}}^{\overline{F}(\overline{\rho}^{\prime})A}\right)^{-1} \circ \left(\phi_{h_{\Gamma_1^{\prime}} \circ h_{\Gamma^{\prime}}, \rho^{\prime}}^{A}\right)^{-1}\right)
	\end{align*}
	is equal to the morphism
	\[
	\left(\phi_{h^{01}_{\Gamma}, h_{\Gamma_1}}^{\overline{F}(\overline{\rho})A}\right)^{-1} \circ \left(\phi_{h^{12}_{\Gamma_1} \circ h^{01}_{\Gamma}, \rho}\right)^{-1} \circ  \overline{F}(h_{\Gamma}^{02})\tau_{f_3}^{A} \circ \left(\phi_{h_{\Gamma}^{02}, f_3}^{A}\right)^{-1} \circ \phi_{f,h_{\Gamma^{\prime}}^{02}}^{A}.
	\]
	Proving this is an extremely tedious but routine verification that follows from the pseudofunctoriality of $\overline{F}$ and the induced relations provided by the pasting diagram:
	\[
	\begin{tikzcd}
		\overline{F}(\quot{X}{\Gamma_3^{\prime}}) \ar[dddrr, ""{name = LeftMid}] \ar[rrrrrr, bend left = 50, ""{name = UU}]{}{\overline{F}(h_{\Gamma^{\prime}}^{02})} \ar[rr, ""{name = UpLeft}]{}[description]{\overline{F}(\overline{\rho}^{\prime})} \ar[ddd, swap]{}{\overline{F}(\of_3)} & & \overline{F}(\quot{X}{\Gamma_2^{\prime}}) \ar[dddrr, ""{name = MM}]\ar[rr, ""{name = UL}]{}[description]{\overline{F}(h_{\Gamma^{\prime}}^{12})} \ar[ddd, swap]{}{\overline{F}(\of_2)} & & \overline{F}(\quot{X}{\Gamma_1^{\prime}}) \ar[dddrr, ""{name = RightMid}]{}{} \ar[ddd, swap]{}{\overline{F}(\of_1)} \ar[rr, ""{name = UpRight}]{}[description]{\overline{F}(h_{\Gamma^{\prime}}^{01})} & & \overline{F}(\quot{X}{\Gamma}) \ar[ddd]{}{\overline{F}(\of)} \\
		\\
		\\
		\overline{F}(\quot{X}{\Gamma_3}) \ar[rrrrrr, swap, bend right = 50, ""{name = LL}]{}{\overline{F}(h_{\Gamma}^{02})} \ar[rr, ""{name = LowLeft}]{}[description]{\overline{F}(\overline{\rho})} & & \overline{F}(\quot{X}{\Gamma_2}) \ar[rr, swap, ""{name = LU}]{}[description]{\overline{F}(h_{\Gamma_1}^{12})} & & \overline{F}(\quot{X}{\Gamma_1}) \ar[rr, ""{name = LowRight}]{}[description]{\overline{F}(h_{\Gamma}^{01})} & & \overline{F}(\quot{X}{\Gamma}) 
		\ar[from = UU, to = UL,  Rightarrow, shorten <= 4pt, shorten >= 4pt]{}{(\phi_{h_{\Gamma_1^{\prime} \circ h_{\Gamma^{\prime}, \rho}}} \ast \phi_{h_{\Gamma^{\prime}}, h_{\Gamma_1^{\prime}}})^{-1}}
		\ar[from = LU, to = LL, Rightarrow, shorten <= 4pt, shorten >= 4pt]{}{\phi_{h_{\Gamma_1} \circ h_{\Gamma}, \rho} \ast \phi_{h_{\Gamma}, h_{\Gamma_1}}}
		\ar[from = UpLeft, to = LeftMid, Rightarrow, shorten >= 4pt, shorten <= 4pt, near end]{}{\phi_{ f_2,\rho^{\prime}}}
		\ar[from = LeftMid, to = LowLeft, Rightarrow, shorten >= 4pt, shorten <= 4pt, swap]{}{(\phi_{ \rho,f_3})^{-1}}
		\ar[from = UL, to = MM, Rightarrow, shorten <= 4pt, shorten >= 4pt, near end]{}{\phi_{f_1,h_{\Gamma^{\prime}}}}
		\ar[from = MM, to = LU, Rightarrow, swap, shorten >=4pt, shorten <= 4pt]{}{(\phi_{h_{\Gamma_1}, f_2})^{-1}}
		\ar[from = UpRight, to = RightMid, Rightarrow, near end, shorten <= 4pt, shorten >= 4pt]{}{\phi_{f,h_{\Gamma^{\prime}}}}
		\ar[from = RightMid, to = LowRight, Rightarrow, swap, shorten <= 4pt, shorten >= 4pt]{}{(\phi_{h_{\Gamma},f_1})^{-1}}
	\end{tikzcd}
	\]
	and the fact that it is equal to the invertible $2$-cell:
	\[
	\begin{tikzcd}
		\overline{F}(\quot{X}{\Gamma_3^{\prime}}) \ar[rrr, ""{name = U}]{}{\overline{F}(h_{\Gamma^{\prime}}^{02})} \ar[d, swap]{}{\overline{F}(\of_3)} & & & \overline{F}(\XGammap) \ar[d]{}{\overline{F}(\of)} \\
		\overline{F}(\quot{X}{\Gamma_3}) \ar[rrr, swap, ""{name = L}]{}{\overline{F}(h_{\Gamma}^{02})} & & & \overline{F}(\XGamma) \ar[from = U, to = L, Rightarrow, shorten <= 4pt, shorten >= 4pt]{}{\phi_{h_{\Gamma}^{02},f_3}^{-1} \circ \phi_{f,h_{\Gamma^{\prime}}^{02}}}
	\end{tikzcd}
	\]
	Doing this shows that the two morphisms are indeed equivalent and establishes that
	\[
	\tau_{f}^{\varphi^{\sharp}(\psi^{\sharp}A)} \circ \overline{F}(\overline{f})\quot{\alpha_A}{\Gamma^{\prime}} = \quot{\alpha_A}{\Gamma} \circ \tau_f^{(\psi \circ \varphi)^{\sharp}A},
	\]
	which in turn completes the verification that $\alpha_A:(\psi \circ \varphi)^{\sharp}A \to \varphi^{\sharp}(\psi^{\sharp}A)$ is a morphism in $\overline{F}_{G_0}(X)$.
	
	We finish by verifying the naturality of $\alpha$, i.e., that for any morphism
	\[
	\Sigma = \lbrace \quot{\sigma}{\Gamma} \; | \; \Gamma \in \Sf(G_2)_0 \rbrace \in \overline{F}_{G_2}(X)(A,B)
	\]
	the diagram
	\[
	\xymatrix{
		(\psi \circ \varphi)^{\sharp}A \ar[r]^-{\alpha_A} \ar[d]_{(\psi \circ \varphi)^{\sharp}\Sigma} & \varphi^{\sharp}(\psi^{\sharp}A) \ar[d]^{\varphi^{\sharp}(\psi^{\sharp}\Sigma)} \\
		(\psi \circ \varphi)^{\sharp}B \ar[r]_-{\alpha_B} & \varphi^{\sharp}(\psi^{\sharp}B)
	}
	\]
	commutes. Now observe that since $\Sigma$ is an $\Sf(G_2)$-morphism, for any object $\Gamma$ of $\Sf(G_0)$ the diagram
	\[
	\xymatrix{
		\overline{F}(\overline{\rho})\quot{A}{\Gamma_3} \ar[rr]^-{\overline{F}(\overline{\rho})\quot{\sigma}{\Gamma_3}} \ar[d]_{\tau_{\rho}^{A}} & & \overline{F}(\overline{\rho})\quot{B}{\Gamma_3} \ar[d]^{\tau_{\rho}^{B}} \\
		\quot{A}{\Gamma_2} \ar[rr]_-{\sigma_{\Gamma_2}} & & \quot{B}{\Gamma_2}
	}
	\]
	commutes. Using this we verify that
	\begin{align*}
		&\quot{\alpha_B}{\Gamma} \circ \quot{(\psi \circ \varphi)^{\sharp}\Sigma}{\Gamma} \\
		&= \left(\overline{F}(h_{\Gamma}^{01})\circ \overline{F}(h_{\Gamma_1}^{12})\right)\tau_{\rho}^{B} \circ \left(\phi_{h_{\Gamma},h_{\Gamma_1}}^{\overline{F}(\overline{\rho})B}\right)^{-1} \circ \left(\phi_{h_{\Gamma_1} \circ h_{\Gamma}, \rho}^{B}\right)^{-1} \circ \overline{F}(h_{\Gamma}^{02})\quot{\sigma}{\Gamma_3} \\
		&=\left(\overline{F}(h_{\Gamma}^{01})\circ \overline{F}(h_{\Gamma_1}^{12})\right)\tau_{\rho}^{B} \circ \left(\phi_{h_{\Gamma},h_{\Gamma_1}}^{\overline{F}(\overline{\rho})B}\right)^{-1} \circ \left(\phi_{h_{\Gamma_1} \circ h_{\Gamma}, \rho}^{B}\right)^{-1} \circ \overline{F}(\overline{\rho} \circ h_{\Gamma_1}^{12} \circ h_{\Gamma}^{01})\quot{\sigma}{\Gamma_3} \\
		&= \left(\overline{F}(h_{\Gamma}^{01})\circ \overline{F}(h_{\Gamma_1}^{12})\right)\tau_{\rho}^{B} \circ \left(\phi_{h_{\Gamma},h_{\Gamma_1}}^{\overline{F}(\overline{\rho})B}\right)^{-1} \circ \left(\overline{F}(h_{\Gamma_1}^{12} \circ h_{\Gamma}^{01})\circ \overline{F}(\overline{\rho})\right)\quot{\sigma}{\Gamma_3} \\
		& \circ \left(\phi_{h_{\Gamma_1} 	\circ h_{\Gamma}, \rho}^{A}\right)^{-1} \\
		&= \left(\overline{F}(h_{\Gamma}^{01})\circ \overline{F}(h_{\Gamma_1}^{12})\right)\tau_{\rho}^{B} \circ \left(\overline{F}(h_{\Gamma}^{01})\circ \overline{F}(h_{\Gamma_1}^{12})\circ \overline{F}(\overline{\rho})\right)\quot{\sigma}{\Gamma_3} \circ \left(\phi_{h_{\Gamma},h_{\Gamma_1}}^{\overline{F}(\overline{\rho})A}\right)^{-1}\\
		& \circ \left(\phi_{h_{\Gamma_1} \circ h_{\Gamma}, \rho}^{A}\right)^{-1} \\
		&= \left(\overline{F}(h_{\Gamma}^{01})\circ \overline{F}(h_{\Gamma_1}^{12})\right)\left(\tau_{\rho}^{B} \circ \overline{F}(\overline{\rho})\quot{\sigma}{\Gamma_3}\right) \circ \left(\phi_{h_{\Gamma},h_{\Gamma_1}}^{\overline{F}(\overline{\rho})A}\right)^{-1} \circ \left(\phi_{h_{\Gamma_1} \circ h_{\Gamma}, \rho}^{A}\right)^{-1} \\
		&= \left(\overline{F}(h_{\Gamma}^{01})\circ \overline{F}(h_{\Gamma_1}^{12})\right)\left(\quot{\sigma}{\Gamma_2} \circ \tau_{\rho}^{A}\right) \circ \left(\phi_{h_{\Gamma},h_{\Gamma_1}}^{\overline{F}(\overline{\rho})A}\right)^{-1} \circ \left(\phi_{h_{\Gamma_1} \circ h_{\Gamma}, \rho}^{A}\right)^{-1} \\
		&= \left(\overline{F}(h_{\Gamma}^{01})\circ \overline{F}(h_{\Gamma_1}^{12})\right)\quot{\sigma}{\Gamma_2} \circ \left(\overline{F}(h_{\Gamma}^{01})\circ \overline{F}(h_{\Gamma_1}^{12})\right)\tau_{\rho}^{A} \circ \left(\phi_{h_{\Gamma},h_{\Gamma_1}}^{\overline{F}(\overline{\rho})A}\right)^{-1} \circ \left(\phi_{h_{\Gamma_1} \circ h_{\Gamma}, \rho}^{A}\right)^{-1} \\
		&= \quot{\varphi^{\sharp}(\psi^{\sharp}\Sigma)}{\Gamma} \circ \quot{\alpha_A}{\Gamma},
	\end{align*}
	which establishes the naturality of $\alpha$.
\end{proof}

We now discuss the de-equivariantification functor, which is a Change of Groups analog of the forgetful functor of Proposition \ref{Prop: Section ECSV: Forgetful functor}, and then give an explicit comparison between the functor $1_G^{\sharp}:F_G(X) \to F_{\Spec K}(X)$ (respectively $1_G^{\sharp}:F_G(X) \to F_{\lbrace \ast \rbrace}(X)$ in the topological case) and $\Forget:F_G(X) \to \overline{F}(X)$. This comparison is a useful consequence of the change of groups formalism; it gives rise to a functor $F_G(X) \to F(X)$ for any group $G$ and any space $X$. Since we always have a morphism of algebraic groups $1_G:\Spec K \to G$ for any algebraic group $G$, from Theorem \ref{Thm: Section 3: Change of groups functor} we get a functor $1_G^{\sharp}:F_G(X) \to F_{\Spec K}(X)$ that restricts along this inclusion of the identity of $G$. Recall that there is an equivalence of categories $F_{\Spec K}(X) \simeq \overline{F}(X)$ (cf.\@ Proposition \ref{Prop: Pseudocone Section: Terminal object in C gives us pseudocones as global sections}); let $\overline{E}:F_{\Spec K}(X) \to \overline{F}(X)$ be one such functor witnessing the equivalence. We then obtain a functor
\[
F_G(X) \xrightarrow{1_G^{\sharp}}  F_{\Spec K}(X) \xrightarrow{\overline{E}} \overline{F}(X)
\]
which suitably ``de-equivariantifies'' the category $F_G(X)$; note that in the topological case this is the functor:
\[
\begin{tikzcd}
F_G(X) \ar[r]{}{1_G^{\sharp}} & F_{\lbrace \ast \rbrace}(X) \ar[r]{}{\overline{E}} & \overline{F}(X) 
\end{tikzcd}
\] 
We call the corresponding functor $\overline{E} \circ 1_G^{\sharp}:F_G(X) \to \overline{F}(X)$  or its topological companion a de-equivariantification functor. In fact, this can be seen as an explanation using the Change of Groups formalism exactly what the forgetful functor $\Forget:F_G(X) \to F(X)$ is doing in a group-theoretic sense. While the proof we present below is functionally very similar to that of Proposition \ref{Prop: Section 3.3: Chagnge of groups interacting with change of grups} above, we present it in detail below to illustrate how the forgetful functor may be seen as a Change of Groups in disguise.
\begin{proposition}\label{Prop: Section 3.3: Forgetful functor is de-equivariantification functor}
	Let $G$ be a smooth algebraic group and let $X$ be a left $G$-variety with unit map $1_G:\Spec K \to G$. Then for any pre-equivariant pseudofunctor $F$ on $X$ there is an invertible $2$-cell:
	\[
	\begin{tikzcd}
		F_G(X) \ar[rr]{}{1_G^{\sharp}} \ar[dr, swap]{}{\Forget} & {} & F_{\Spec K}(X) \ar[dl]{}{\overline{E}} \\
		& \overline{F}(X) \ar[from = 2-2, to = 1-2, Rightarrow, shorten >= 4pt, shorten <= 4pt]{}{\alpha}
	\end{tikzcd}
	\]
	Similarly, if $G$ is a topological group with unit map $1_G:\lbrace \ast \rbrace \to G$ and $X$ is a left $G$-space then for any pre-equivariant pseudofunctor on $X$ there is an invertible $2$-cell:
	\[
		\begin{tikzcd}
		F_G(X) \ar[rr]{}{1_G^{\sharp}} \ar[dr, swap]{}{\Forget} & {} & F_{\lbrace \ast \rbrace}(X) \ar[dl]{}{\overline{E}} \\
		& \overline{F}(X) \ar[from = 2-2, to = 1-2, Rightarrow, shorten >= 4pt, shorten <= 4pt]{}{\alpha}
	\end{tikzcd}
	\]
\end{proposition}
\begin{proof}
	We begin by observing that there is a $G$-equivariant isomorphism of varieites in $\Sf(G)$ of the form
	\[
	G \times^{\Spec K} \Spec K \cong G
	\] 
	(this can be seen in set-theoretic terms via the mappings $[g,1]_{1} \mapsto g, g \mapsto [g,1]_{1}$), and this isomorphism induces the further isomorphism
	\[
	\quot{X}{G \times^{\Spec K} \Spec K} \xrightarrow[\rho]{\cong} \quot{X}{G}.
	\]
	Now consider the isomorphism
	\[
	h_{\Spec K}:\quot{X}{\Spec K} \xrightarrow{\cong} \quot{X}{G \times^{\Spec K} \Spec K}
	\] 
	and write the isomorphism from $X$ to $\quot{X}{\Spec K}$as
	\[
	X \xrightarrow[\varphi]{\cong} \quot{X}{\Spec K}.
	\]
	Finally, let $\psi:X \to \quot{X}{G}$ be the standard (non-$G$-equivariant) isomorphism of varieties. Let us now observe that on one hand we have that for any object $A \in F_G(X)_0$
	\begin{align*}
		\overline{E}(1_G^{\sharp}(A)) &= \overline{E}\left(\left\lbrace \overline{F}(h_{\Gamma})\left(\quot{A}{G \times^{\Spec K} \Gamma}\right) \; : \; \Gamma \in \Sf(G)_0 \right\rbrace\right) \\
		&= \left(\overline{F}(\varphi) \circ \overline{F}(h_{\Spec K})\right)\left(\quot{A}{G \times^{\Spec K} \Spec K}\right)
	\end{align*}
	and similarly for morphisms. On the other hand we have that
	\[
	\Forget(A) = \overline{F}(\psi)(\quot{A}{G}).
	\]
	Now since the diagram
	\[
	\xymatrix{
		\quot{X}{G \times^{\Spec K} \Spec K} \ar[d]_{\overline{\rho}} & & \quot{X}{\Spec K} \ar[ll]_-{h_{\Spec K}} \\
		\quot{X}{G} & & X \ar[ll]^-{\psi} \ar[u]_{\varphi}
	}
	\]
	of varieties commutes by the uniqueness of isomorphisms between quotient schemes, we have the following pasting diagram (where we implicitly use that on quotient varieties $F = \overline{F}$)
	\[
	\begin{tikzcd}
		& \overline{F}(\quot{X}{\Spec K}) \ar[dr]{}{\overline{F}(\varphi)} & \\
		\overline{F}\left(\quot{X}{G \times^{\Spec K} \Spec K}\right) \ar[rr]{}[description]{\overline{F}(h_{\Spec K} \circ \varphi)} \ar[ur]{}{\overline{F}(h_{\Spec K})} & {} & \overline{F}(X) \\
		& \overline{F}\left(\quot{X}{G}\right) \ar[ur, swap]{}{\overline{F}(\psi)} \ar[ul]{}{\overline{F}(\overline{\rho})} & \ar[from = 2-2, to = 3-2, swap, Rightarrow, shorten >= 4pt, shorten <= 4pt]{}{\phi_{h \circ \varphi, \rho}} \ar[from = 1-2, to = 2-2, Rightarrow, shorten <= 4pt, shorten >=4pt, swap]{}{\phi_{\varphi, h}}
	\end{tikzcd}
	\]
	in $\fCat$. Notice that because $\rho \in \Sf(G)_1$, we have an isomorphism
	\[
	\tau_{\rho}^{A}:\overline{F}(\overline{\rho})\quot{A}{G} \xrightarrow{\cong} \quot{A}{(G \times \Spec K)_{\Spec K}}
	\]
	in $T_A$. Applying $\overline{F}(\varphi) \circ \overline{F}(h_{\Spec K})$ to this isomorphism gives us an isomorphism
	\[
	\begin{tikzcd}
		& & \left(\overline{F}(\varphi)\circ \overline{F}(h_{\Spec K})\right)\left(\quot{A}{G \times^{\Spec K} \Spec K}\right) \ar[dd, equals] \\
		\left(F(\varphi)\circ F(h_{\Spec K})\circ F(\overline{\rho})\right)\left(\quot{A}{G}\right) \ar[drr] \ar[urr]{}{\cong} \ar[urr, swap]{}{\overline{F}(\varphi)(F(h_{\Spec K})\tau_{\rho}^{A})}  \\
		& & (\overline{E}\circ 1_G^{\sharp})(A)
	\end{tikzcd}
	\]
	To get our isomorphism to have the desired domain of 
	\[
	\Forget(A) = \overline{F}(\psi)\left(\quot{A}{G}\right) = \overline{F}(\overline{\rho} \circ h_{\Spec K} \circ \varphi)\left(\quot{A}{G}\right),
	\] 
	we pre-compose by the inverses of the two natural isomorphisms to produce the isomorphism
	\[
	\alpha_A:= \left(\overline{F}(\varphi) \circ \overline{F}(h_{\Spec K})\right)\tau_{\rho}^{A} \circ \left(\phi_{\varphi,h}^{F(\overline{\rho})\quot{A}{G}}\right)^{-1} \circ \left(\phi_{h \circ \varphi, \rho}^{\quot{A}{G}}\right)^{-1}.
	\]
	Then $\alpha_A:\Forget(A) \xrightarrow{\cong} (\overline{F} \circ 1_G^{\sharp})(A)$, as desired.
	
	We end by showing that $\alpha$ is natural in $F_G(X)$, i.e.,that for a morphism
	\[
	\Sigma = \lbrace \quot{\sigma}{\Gamma} \; | \; \Gamma \in \Sf(G)_0 \rbrace \in F_G(X)(A,B)
	\]
	that the diagram
	\[
	\xymatrix{
		\Forget(A) \ar[r]^-{\alpha_A} \ar[d]_{\Forget(\Sigma)} & (\overline{F} \circ 1_G^{\sharp})(A) \ar[d]^{(\overline{F}\circ 1_G^{\sharp})(\Sigma)} \\
		\Forget(B) \ar[r]_-{\alpha_B} & (\overline{F} \circ 1_G^{\sharp})(B)
	}
	\]
	commutes. First recall that since $\Sigma$ is an $F_G(X)$-morphisms, the diagram
	\[
	\begin{tikzcd}
		F(\overline{\rho})(\quot{A}{G}) \ar[d, swap]{}{\tau_{\rho}^{A}} \ar[rrr]{}{F(\overline{\rho})\quot{\sigma}{G}} & & & F(\overline{\rho})\quot{B}{G} \ar[d]{}{\tau_{\rho}^{B}} \\
		\quot{A}{G \times^{\Spec K} \Spec K} \ar[rrr, swap]{}{\quot{\sigma}{G \times^{\Spec K} \Spec K}} & & & \quot{B}{G \times^{\Spec K} \Spec K}
	\end{tikzcd}
	\]
	commutes in $\overline{F}(\quot{X}{G \times^{\Spec K} \Spec K}).$ We then conclude the proof by calculating
	\begin{align*}
		&(\overline{E}\circ 1_G^{\sharp})(\Sigma) \circ \alpha_A \\
		&= \left(\overline{F}(\varphi)\circ \overline{F}(h_{\Spec K})\right)\left(\quot{\sigma}{G \times^{\Spec K} \Spec K}\right) \circ \left(\overline{F}(\varphi)\circ\overline{F}(h_{\Spec K})\right)\tau_{\rho}^{A} \\ 
		&\circ \left(\phi_{\varphi,h}^{F(\overline{\rho})\quot{A}{G}}\right)^{-1} \circ \left(\phi_{h \circ \varphi, \rho}^{\quot{A}{G}}\right)^{-1} \\
		&= \left(\overline{F}(\varphi)\circ\overline{F}(h_{\Spec K})\right)\left(\quot{\sigma}{G \times^{\Spec K} \Spec K} \circ \tau_{\rho}^{A}\right) \circ \left(\phi_{\varphi,h}^{F(\overline{\rho})\quot{A}{G}}\right)^{-1} \circ \left(\phi_{h \circ \varphi, \rho}^{\quot{A}{G}}\right)^{-1} \\
		&= \left(\overline{F}(\varphi)\circ\overline{F}(h_{\Spec K})\right)\left(\tau_{\rho}^{B} \circ F(\overline{\rho})\quot{\sigma}{G}\right) \circ \left(\phi_{\varphi,h}^{F(\overline{\rho})\quot{A}{G}}\right)^{-1} \circ \left(\phi_{h \circ \varphi, \rho}^{\quot{A}{G}}\right)^{-1} \\
		&= \left(\overline{F}(\varphi)\circ\overline{F}(h_{\Spec K})\right)(\tau_{\rho}^{B}) \circ \left(\phi_{\varphi,h}^{F(\overline{\rho})\quot{B}{G}}\right)^{-1} \circ \left(\phi_{h \circ \varphi, \rho}^{\quot{B}{G}}\right)^{-1} \circ \overline{F}(\psi)\quot{\sigma}{G} \\
		&= \alpha_B \circ \Forget(\Sigma).
	\end{align*}
\end{proof}
\newpage

\section{A High-Level Perspective On Change of Groups Functors}\label{Section: High Level Perspective Change of Groups}
In this section we give a high-level perspective on the functors we have just developed and studied. Because of the high degree of technical information involved in the construction Change of Group functors, it can be easy to lose ourselves in the weeds; as such, we should have a $2$-category and a pseudofunctor or two running around to help us record exactly what those details are describing and recording. With this in mind, we start this section with a goal of introducing a $2$-category $\GHPreq{G}{H}$ of $(G,H)$-pre-equivariant pseudofunctors on a variety $X$ for fixed algebraic groups $G, H$ with a morphism $\varphi:G \to H$ (respectively of $(G,H)$-pre-equivariant pseudofunctors on an $H$-space $X$ for fixed topological groups $G, H$) before describing a more general $2$-category $\GHPreq{-}{H}$ which leaves only $H$ fixed. Consequently, we need to have a notion of morphism between $(G,H)$-pre-equivariant pseudofunctors as well as transformations between such morphisms.

\begin{definition}\label{Defn: Morphism of GHPreq}
	Let $G$, $H$, $\varphi:G \to H$, and $X$ be as in either of the cases below:\index[terminology]{Morphism! Of $(G,H)$-pre-equivariant Pseudofunctors}
	\begin{enumerate}
		\item $G$ and $H$ are smooth algebraic groups over a field $K$, $X$ is a left $H$-variety, and $\varphi$ is an algebraic group morphism $\varphi:G \to H$;
		\item $G$ and $H$ are topological groups, $X$ is a left $H$-space, and $\varphi$ is a topological group morphism $\varphi:G \to H$.
	\end{enumerate} 
	If $\ul{F}_1` = \left((F_1,\overline{F}_1),(F_1^{\prime},\overline{F}_1^{\prime}),e_1\right)$ and $\ul{F}_2 = \left((F_2,\overline{F}_2),(F_2^{\prime},\overline{F}_2^{\prime}), e_2\right)$ are $(G,H)$-pre-equivariant pseudofunctors on $X$ then a morphism $\alpha:\ul{F}_1 \to \ul{F}_2$ of $(G,H)$-pre-equivariant pseudofunctors is a quadruple $\alpha = (\alpha, \overline{\alpha},\alpha^{\prime}, \overline{\alpha}^{\prime})$ where:
	\begin{itemize}
		\item $\alpha:F_1 \Rightarrow F_2, \overline{\alpha}:\overline{F}_1 \Rightarrow \overline{F}_2, \alpha^{\prime}:F_1^{\prime} \Rightarrow F_2^{\prime}$, and $\overline{\alpha}^{\prime}:\overline{F}_1^{\prime} \Rightarrow \overline{F}_2^{\prime}$ are pseudonatural transformations.
		\item The whiskering identities $\alpha = \overline{\alpha} \ast \quo_H^{\op}$ and $\alpha^{\prime} = \overline{\alpha}^{\prime} \ast \quo_G^{\op}$ hold.
		\item The diagram
		\[
		\begin{tikzcd}
			\overline{F}_1 \ar[r]{}{\overline{\alpha}} \ar[d, swap]{}{e_1} & \overline{F}_2 \ar[d]{}{e_2} \\
			\overline{F}^{\prime}_1 \ar[r, swap]{}{\overline{\alpha}^{\prime}} & \overline{F}^{\prime}_2
		\end{tikzcd}
		\]
		of pseudofunctors and pseudonatural transformations commutes.
	\end{itemize}
\end{definition}
Immediate from the definition and Theorem \ref{Thm: Functor Section: Psuedonatural trans are pseudocone functors} is the following lemma.
\begin{lemma}\label{Lemma: GHPreq morphisms give four functors and a cohernece}
	A morphism $\alpha$ of $(G,H)$-pre-equivariant pseudofunctors gives rise to functors $(F_1)_G(X) \to (F_2)_G(X)$, $\PC(\overline{F}_1) \to \PC(\overline{F}_2)$, $(F_1^{\prime})_H(X) \to (F_2^{\prime})_H(X)$, and $\PC(\overline{F}^{\prime}_1) \to \PC(\overline{F}_2^{\prime})$ for which the diagram
	\[
	\begin{tikzcd}
		\PC(\overline{F}_1) \ar[r]{}{\PC(\overline{\alpha})} \ar[d, swap]{}{\PC(e_1)} & \PC(\overline{F}_2) \ar[d]{}{\PC(e_2)} \\
		\PC(\overline{F}_1^{\prime}) \ar[r, swap]{}{\PC(\overline{\alpha}^{\prime})} & \PC(\overline{F}^{\prime}_2)
	\end{tikzcd}
	\]
	of categories and functors commutes. 
\end{lemma}

\begin{definition}\label{Defn: Transformation of GHPreq}\index[terminology]{Transformation! Of $(G,H)$-pre-equivariant Pseudofunctors}
	Let $G$, $H$, $X$, and $\varphi$ all be given as in Cases 1 or 2 of Definition \ref{Defn: Morphism of GHPreq}. Let $\ul{F}_1` = \left((F_1,\overline{F}_1),(F_1^{\prime},\overline{F}_1^{\prime}),e_1\right)$ and $\ul{F}_2 = \left((F_2,\overline{F}_2),(F_2^{\prime},\overline{F}_2^{\prime}), e_2\right)$ be $(G,H)$-pre-equivariant pseudofunctors on $X$ with morphisms $\alpha,\beta:\ul{F}_1 \to \ul{F}_2$ of $(G,H)$-pre-equivariant pseudofunctors. A transformation $\rho:\alpha \Rightarrow \beta$ is a quadruple $\rho = (\rho, \overline{\rho}, \rho^{\prime}, \overline{\rho})$ where:
	\begin{itemize}
		\item $\rho:\alpha \Rrightarrow \beta, \overline{\rho}:\overline{\alpha} \Rrightarrow \overline{\beta}, \rho^{\prime}:\alpha^{\prime} \Rrightarrow \beta^{\prime}$, and $\overline{\rho}^{\prime}:\overline{\alpha}^{\prime} \Rrightarrow \overline{\beta}^{\prime}$ are modifications.
		\item The whiskering identity $e_2 \ast \overline{\rho} = \overline{\rho}^{\prime} \ast e_1$ holds.
		\item The whiskering identities $\rho = \overline{\rho} \ast \iota_{\quo_G^{\op}}$ and $\rho^{\prime} = \overline{\rho}^{\prime} \ast \iota_{\quo_G^{\op}}$ both hold.
	\end{itemize}
\end{definition}
We, similar to Lemma \ref{Lemma: GHPreq morphisms give four functors and a cohernece} above, derive the following lemma immediately from applying Lemma \ref{Lemma: Modifications give equivariant natural transformations} to the definition above.
\begin{lemma}\label{Lemma: GHPreq transf give nat transforms}
	A transformation $\rho:\alpha \Rightarrow \beta$ between $(G,H)$-pre-equivariant pseudofunctor morphisms $\alpha,\beta:\ul{F}_1 \to \ul{F}_2$ give rise to natural transformations
	\[
	\begin{tikzcd}
		(F_1)_G(X) \ar[rr, bend left = 30, ""{name = U}]{}{\alpha} \ar[rr, bend right = 30, swap, ""{name = D}]{}{\beta} & & (F_2)_G(X) \ar[from = U, to = D, Rightarrow, shorten <= 4pt, shorten >= 4pt]{}{\rho}
	\end{tikzcd}\qquad
	\begin{tikzcd}
		(F_1^{\prime})_H(X) \ar[rr, bend left = 30, ""{name = U}]{}{\alpha^{\prime}} \ar[rr, bend right = 30, swap, ""{name = D}]{}{\beta^{\prime}} & & (F_2^{\prime})_H(X) \ar[from = U, to = D, Rightarrow, shorten <= 4pt, shorten >= 4pt]{}{\rho^{\prime}}
	\end{tikzcd}
	\]
	\[
	\begin{tikzcd}
		\PC(\overline{F}_1) \ar[rr, bend left = 30, ""{name = U}]{}{\PC(\overline{\alpha})} \ar[rr, bend right = 30, swap, ""{name = D}]{}{\PC(\overline{\beta})} & & \PC(\overline{F}_2) \ar[from = U, to = D, Rightarrow, shorten <= 4pt, shorten >= 4pt]{}[description]{\PC(\overline{\rho})}
	\end{tikzcd}\qquad
	\begin{tikzcd}
		\PC(\overline{F}_1^{\prime}) \ar[rr, bend left = 30, ""{name = U}]{}{\PC(\overline{\alpha}^{\prime})} \ar[rr, bend right = 30, swap, ""{name = D}]{}{\PC(\overline{\beta}^{\prime})} & & \PC(\overline{F}_2) \ar[from = U, to = D, Rightarrow, shorten <= 4pt, shorten >= 4pt]{}[description]{\PC(\overline{\rho}^{\prime})}
	\end{tikzcd}
	\]
	for which the identity $\PC(\overline{\rho}^{\prime}) \ast \PC(e_1) = \PC(e_2) \ast \PC(\overline{\rho})$ holds.
\end{lemma}

A straightforward consequence of these constructions is that it allows us to organize all three pieces of information into a $2$-category $\GHPreq{G}{H}$ of $(G,H)$-pre-equivariant pseudofunctors on $X$.
\begin{proposition}\label{Prop: Section Chofg: GHPreq twocat}
	Let $K$ be a field, $G$ and $H$ be smooth algebraic groups over $K$ with $X$ a left $H$-variety and $\varphi:G \to H$ a morphism of smooth algebraic groups. Then there is a $2$-category $\GHPreq{G}{H}$ of $(G,H)$-pre-equivariant pseudofunctors on $X$ where:
	\begin{itemize}
		\item The $0$-cells are $(G,H)$-pre-equivariant pseudofunctors;
		\item The $1$-cells are morphisms of $(G,H)$-pre-equivariant pseudofunctors;
		\item The $2$-cells are transformations of morphisms of $(G,H)$-pre-equivariant pseudofunctors;
		\item Composition of $1$-cells is induced by the vertical composition of pseudonatural transformations;
		\item Composition of $2$-cells is induced by the composition of $2$-cells in the $2$-categories $\Bicat(-,\fCat)$. 
	\end{itemize}
\end{proposition}
\begin{proof}
	Fix two $(G,H)$-pre-equivariant pseudofunctors $\ul{F}_1$ and $\ul{F}_2$. Define the hom-category $\GHPreq{G}{H}(\ul{F}_1,\ul{F}_2)$ as follows:
	\begin{itemize}
		\item Objects: morphisms $\alpha:\ul{F}_1 \to \ul{F}_2$.
		\item Morphisms: transformations $\rho:\alpha \to \beta$.
		\item Composition: Given morphisms $\alpha,\beta,\gamma:\ul{F}_1 \to \ul{F}_2$ with transformations $\varphi:\alpha \to \beta$ and $\psi:\beta \to \gamma$, we define
		\[
		\psi \circ \varphi := \left(\psi \circ \varphi, \overline{\psi} \circ \overline{\varphi}, \psi^{\prime} \circ \varphi^{\prime}, \overline{\psi}^{\prime} \circ \overline{\varphi}^{\prime}\right).
		\]
		\item Identities: Given a morphism $\alpha:\ul{F}_1 \to \ul{F}_2$, the identity transformation is $\iota_{\alpha} = (\iota_{\alpha},\iota_{\overline{\alpha}},\iota_{\alpha^{\prime}}, \iota_{\overline{\alpha}^{\prime}})$.
	\end{itemize}
	That this is a category is trivial to verify, as the fact that composition in the assignment $\GHPreq{G}{H}(\ul{F}_1,\ul{F}_2)$ is strictly associative in the geometric case follows from the fact that each of the three bicategories $\Bicat(\SfResl_G(X)^{\op},\fCat)$, $\Bicat(\SfResl_H(X)^{\op},\fCat)$, and $\Bicat(\Var_{/K}^{\op},\fCat)$ are actually $2$-categories; in the topological case that composition is strictly associative follows from the fact that the three bicategories $\Bicat(\mathbf{FRes}_G(X)^{\op},\fCat)$, $\Bicat(\mathbf{FRes}_H(X)^{\op},\fCat)$, and $\Bicat(\Top^{\op},\fCat)$ are all $2$-categories.
	
	Alternatively, if we have $(G,H)$-pre-equivariant pseudofunctors $\ul{F}_1, \ul{F}_2$, and $\ul{F}_3$ on $X$ then the horizontal composition functors
	\[
	\begin{tikzcd}
		\GHPreq{G}{H}(\ul{F}_2, \ul{F}_3) \times \GHPreq{G}{H}(\ul{F}_1, \ul{F}_2) \ar[d]{}{\ast} \\
		\GHPreq{G}{H}(\ul{F}_1, \ul{F}_3)
	\end{tikzcd}
	\]
	are defined by taking, for each pair of morphisms $\alpha:\ul{F}_1 \to \ul{F}_2$ and $\beta:\ul{F}_2 \to \ul{F}_3$,
	\[
	\ast(\beta, \alpha) = \left(\beta \ast \alpha, \overline{\beta}\ast\overline{\alpha}, \beta^{\prime} \ast \alpha^{\prime}, \overline{\beta}^{\prime} \ast \overline{\alpha}^{\prime}\right)
	\]
	as the object assignment (where the horizontal composition is defined as at the beginning of Section \ref{Subsection: Change of Fibre}) and
	\[
	\ast(\psi,\varphi) = \left(\psi \ast \varphi, \overline{\psi} \ast \overline{\varphi}, \psi^{\prime} \ast \varphi^{\prime}, \overline{\psi}^{\prime} \ast \overline{\varphi}^{\prime}\right)
	\]
	for transformations $\varphi:\alpha \Rightarrow \beta:\ul{F}_1 \to \ul{F}_2$ and $\psi:\gamma \Rightarrow \delta:\ul{F}_2 \to \ul{F}_3$; note that the horizontal composition of modifications is also defined at the beginning of Section \ref{Subsection: Change of Fibre}. Then because each bicategory $\Bicat(\SfResl_G(X)^{\op},\fCat)$, \\$\Bicat(\SfResl_H(X)^{\op},\fCat)$, and $\Bicat(\Var_{/K}^{\op},\fCat)$ is a $2$-category (and similarly for the topological case), it follows that $\ast$ is a functor. In particular, this allows us to deduce that $\GHPreq{G}{H}$ is a category enriched in $\Cat$ and hence is a $2$-category.
\end{proof}
\begin{remark}
	Following Remark \ref{Remark: Notations and conventions}, we will write $(G,H)\PreEq(X)$ for the $1$-category of $(G,H)$-pre-equivariant pseudofunctors and their morphisms.
\end{remark}

A straightforward translation of the proof above allows us to construct two further $2$-categories useful in the study of Change of Groups functors. The first of these $2$-categories is in essence the $2$-category of $(G,H)$-pre-equivariant pseudofunctors on $X$ where we allow $G$ to vary through all the algebraic groups (respectively topological groups) and structure maps $\varphi:G \to H$ appearing in the category $\mathbf{SAlgGrp}_{/K} \downarrow H$ (respectively also the category $\mathbf{TopGrp}\downarrow H$; cf.\@ Remark \ref{Remark: Notations and conventions} for the definition of the categories $\mathbf{SAlgGrp}_{/K}$ and $\mathbf{TopGrp}\downarrow H$) while the second is the $2$-category generated by all  $(G,G^{\prime})$-pre-equivariant pseudofunctors as we vary through morphisms $\varphi:G \to G^{\prime}$ in $\mathbf{SAlgGrp} \downarrow H$ (respectively in $\mathbf{TopGrp}\downarrow H$).

\begin{definition}
Let $G, H$, and $X$ be given as in either cases:
\begin{itemize}
	\item $G$ is an arbitrary smooth algebraic group over $K$ and for which there is a morphism $\varphi:G \to H$ of smooth algebraic groups and $X$ is a left $H$-variety;
	\item $G$ is an arbitrary topological group for which there is a morphism $\varphi:G \to H$ of topological groups and $X$ is a left $H$-space.
\end{itemize} 
The $2$-category $\GHPreq{-}{H}$\index[notation]{GHPreqdash@$\GHPreq{-}{H}$} is defined as follows:
	\begin{itemize}
		\item The $0$-cells of $\GHPreq{-}{H}$ are $(G,H)$-pre-equivariant pseudofunctors on $X$ in either case:

		More succinctly,
		\[
		\GHPreq{-}{H}_0 = \bigcup_{G \in (\mathbf{SAlgGrp}_{/K} \downarrow H)_0} \GHPreq{G}{H}.
		\]
		\item The $1$-cells of $\GHPreq{-}{H}$ are generated by
		\[
		\GHPreq{-}{H}_1 = \left\langle \bigcup_{G \in (\mathbf{SAlgGrp}_{/K}\downarrow H)_0} \GHPreq{G}{H}_1 \right\rangle.
		\]
		\item The $2$-cells are generated by
		\[
		\GHPreq{-}{H}_2 = \left\langle\bigcup_{G \in (\mathbf{SAlgGrp}_{/K}\downarrow H)_0} \GHPreq{G}{H}_2\right\rangle.
		\]
	\end{itemize}
\end{definition}
\begin{remark}\label{Remark: It makes sense to do stuff chofg}
	For any groups $G, G^{\prime}$ over $H$ with a morphism $\varphi:G \to G^{\prime}$, we can equip the $H$-object $X$ with the pullback action against the structure maps $\nu_G:G \to H$ and $\nu_{G^{\prime}}:G^{\prime} \to H$ in order to produce the category $\GHPreq{G}{G^{\prime}}$; note that the subtlety here lies in equipping $X$ with $G^{\prime}$ and $G$ actions which are suitably compatible with each other through the pullback action of $G$ on $X$ as induced through the map $\varphi:G \to G^{\prime}$ (which, in this case, works out exactly because $\varphi$ is a morphism in the category of groups over $H$, giving that the equation $\nu_{G^{\prime}} \circ \varphi = \nu_G$ holds).
\end{remark}

Our main goal for the moment is to show that the Change of Groups functors are pseudofunctorial by using $\GHPreq{-}{H}$ as a tool to mediate and describe exactly what it is that Change of Groups functors are doing. We first record that the $2$-category $\GHPreq{-}{H}$ arises as the elements of the pseudofunctor $(\mathbf{SAlgGrp}_{/K} \downarrow H)^{\op} \to \fCat$ which sends an algebraic group to its $1$-category of $(G,H)$-pre-equivariant pseudofunctors.

\begin{lemma}\label{Lemma: Is this a pseudofunctor}
Let $H$ and $X$ be given as in either of the two cases below:
\begin{enumerate}
	\item $H$ is a smooth algebraic group and $X$ is a left $H$-variety.
	\item $H$ is a topological group and $X$ is a left $H$-space.
\end{enumerate}
There are pseudofunctors
	\[
	\Ebb:(\mathbf{SAlgGrp}_{/K} \downarrow H)^{\op} \to \fCat
	\]
	and
	\[
	\Fbb:(\mathbf{TopGrp}_{/K} \downarrow H)^{\op} \to \fCat
	\]
	defined by
	\[
	\Ebb(G) := (G,H)\PreEq(X)
	\]
	and
	\[
	\Fbb(G) := (G,H)\PreEq(X)
	\]
	on objects. For morphisms $\varphi:G \to G^{\prime}$, the fibre functors
	\[
	\Ebb(\varphi) := \varphi^{\natural}:(G^{\prime},H)\PreEq(X) \to (G,H)\PreEq(X)
	\]
	and
	\[
	\Fbb(\varphi) := \varphi^{\natural}:(G^{\prime},H)\PreEq(X) \to (G,H)\PreEq(X)
	\]
	are defined by sending a $(G^{\prime},H)$-pre-equivariant pseudofunctor $\left((F,\overline{F}),(F^{\prime},\overline{F}^{\prime}),e\right)$ to the $(G,H)$-pre-equivariant pseudofunctor 
	\[
	\left(\left(F \circ\left(G^{\prime} \times^{G}(-)\right)^{\op}, \overline{F} \right), (F^{\prime},\overline{F}^{\prime}),e\right)
	\]
	and similarly on morphisms.
\end{lemma}
\begin{proof}
	As usual, we present the geometric argument as the topological argument is carried out mutatis mutandis to the geometric one. We first need to type-check $\varphi^{\natural}$. However, note that since $G^{\prime} \times^{G} (-):\SfResl_G(X) \to \SfResl_{G^{\prime}}(X)$ is a functor, by construction $(F \circ (G \times^{G} (-))^{\op}, \overline{F})$ is a $G$-pre-equivariant pseudofunctor on $X$; that it fits into the $(G,H)$-pre-equivariant pseudofunctor claimed follows by construction as well. Finally, that this is pseudofunctorial in $\mathbf{SAlgGrp}_{/K} \downarrow H$ follows from the natural isomorphisms
	\[
	G^{\prime\prime} \times^{G^{\prime}} \left(G^{\prime} \times^{G} (-)\right) \cong G^{\prime\prime} \times^{G} (-)
	\]
	of $G^{\prime\prime}$-varieties for any composable morphisms $G \xrightarrow{\varphi} G^{\prime} \xrightarrow{\psi} G^{\prime\prime}$ in $\mathbf{SAlgGrp} \downarrow H$.
\end{proof}

In what follows we will examine to what degree the Change of Groups functors are associative. The main task in front of us is to verify that if we have algebraic group (respectively topological group) morphisms
\[
G_0 \xrightarrow{\varphi_{01}} G_1 \xrightarrow{\varphi_{12}} G_2 \xrightarrow{\varphi_{23}} G_3
\]
then for any left $G_3$-object $X$ and any $G_3$-pre-equivariant pseudofunctor $F = (F,\overline{F})$ on $X$, the pasting diagram
\begin{equation}\label{Eqn: Pasto 1 for chofg associative}
	\begin{tikzcd}
		& F_{G_2}(X) \ar[dr]{}{\varphi_{12}^{\sharp}} \\
		F_{G_3}(X) \ar[ur]{}{\varphi_{23}^{\sharp}} \ar[dr, swap]{}{(\varphi_{23} \circ \varphi_{12} \circ \varphi_{01})^{\sharp}} \ar[rr, ""{name = M}]{}[description]{(\varphi_{23} \circ \varphi_{12})^{\sharp}} & & F_{G_1}(X) \ar[dl]{}{\varphi_{01}^{\sharp}} 	 \\
		& F_{G_0}(X) \ar[from = 3-2, to = M, Rightarrow, shorten <= 4pt, shorten >= 4pt]{}{\alpha_{0,123}} \ar[from = M, to = 1-2, Rightarrow, shorten <= 4pt, shorten >= 4pt]{}{\alpha_{123}}
	\end{tikzcd}
\end{equation}
is equal to the pasting diagram:
\begin{equation}\label{Eqn: Pasto 2 for chofg associative}
	\begin{tikzcd}
		& F_{G_2}(X) \ar[dr]{}{\varphi_{12}^{\sharp}} \ar[dd, ""{name = M}]{}{}  \ar[dd, crossing over, near start]{}[description]{(\varphi_{12} \circ \varphi_{01})^{\sharp}} \\
		F_{G_3}(X) \ar[ur]{}{\varphi_{23}^{\sharp}} \ar[dr, swap]{}{(\varphi_{23} \circ \varphi_{12} \circ \varphi_{01})^{\sharp}}  & & F_{G_1}(X) \ar[dl]{}{\varphi_{01}^{\sharp}} 	 \\
		& F_{G_0}(X)  \ar[from = M, to = 2-3, Rightarrow, swap, shorten <= 4pt, shorten <= 4pt]{}{\alpha_{012}} \ar[from = 2-1, to = M, swap, Rightarrow, shorten <= 4pt, shorten >= 4pt]{}{\alpha_{012,3}}
	\end{tikzcd}
\end{equation}
However, because of the nature of each of these functors and constructions involving successive induction space constructions and cancelations (such as the examples
\begin{align*}
	&G_3 \times^{G_2}\left(G_2 \times^{G_1}\left(G_1 \times^{G_0} \Gamma \right)\right) &\qquad G_3 \times^{G_0} \Gamma  \\
	&G_2 \times^{G_1}\left(G_1 \times^{G_0} \Gamma\right) &\qquad G_3 \times^{G_1}\left(G_1 \times^{G_0}\Gamma\right)\\
	&G_3 \times^{G_2}\left(G_2 \times^{G_0} \Gamma\right) &\qquad G_1 \times^{G_0} \Gamma
\end{align*}
just listed, among others), we will require some notational shorthand to reduce simultaneously the notational clutter but also to expedite the process of weaving through the induction spaces\footnote{A secondary benefit of this is that the diagrams which involve the induction spaces, morphisms between them, and their quotients now actually fit within the page margins.}. As such we introduce the following notation for any object $\Gamma$ of $\Sf(G_0)$ (respectively if $\Gamma$ is a free $G$-space).

\begin{definition}
	Let $G_0, G_1, G_2, G_3, \Gamma, \hat{\Gamma}, \tilde{\Gamma}$ be given as in either of the following cases:
	\begin{enumerate}
		\item $G_0, G_1, G_2,$ and $G_3,$ are smooth algebraic groups; $\Gamma$ is an $\Sf(G_0)$-variety; $\hat{\Gamma}$ is an $\Sf(G_1)$-variety; and $\tilde{\Gamma}$ is an $\Sf(G_2)$-variety;
		\item $G_0, G_1, G_2,$ and $G_3,$ are topological groups; $\Gamma$ is a free $G_0$-space; $\hat{\Gamma}$ is a free $G_1$-space; and $\tilde{\Gamma}$ is a free $G_2$-space.
	\end{enumerate}
	Assume that we have morphisms
	\[
	G_0 \xrightarrow{\varphi_{01}} G_1 \xrightarrow{\varphi_{12}} G_2 \xrightarrow{\varphi_{23}} G_3
	\]
	of group objects. Then we define:
	\begin{align*}
		\Gamma_{01} &:= G_1 \times^{G_0} \Gamma & \Gamma_{012} &:= G_2 \times^{G_1} \Gamma_{01} & \Gamma_{02} &:= G_2 \times^{G_0} \Gamma \\
		\Gamma_{03} &:= G_3 \times^{G_0} \Gamma & \Gamma_{0123} &:= G_3 \times^{G_2} \Gamma_{012} &  \tilde{\Gamma}_{23}&:= G_3 \times^{G_2} \tilde{\Gamma} \\
		\hat{\Gamma}_{12} &:= G_2 \times^{G_1} \hat{\Gamma}  & \hat{\Gamma}_{123} & := G_3 \times^{G_2}\tilde{\Gamma}_{12} & \Gamma_{013} &:= G_3 \times^{G_1} \Gamma_{01} \\
		\Gamma_{023} &:= G_3 \times^{G_2} \Gamma_{02} & & & \hat{\Gamma}_{13} &:= G^3 \times^{G_1} \hat{\Gamma}
	\end{align*}
	Furthermore, we write the following for the isomorphisms of Lemma \ref{Lemma: Section 3: Induction space quotient is iso to the Sf(G) quotient} applied to the various situations above:
	\begin{align*}
		h_{\Gamma}^{01}:G_0 \backslash(\Gamma\times X) &\xrightarrow{\cong} G_1 \backslash\left(G_1 \times^{G_0}\left(\Gamma \times X\right)\right) =\quot{X}{\Gamma_{01}} \\
		h_{\Gamma}^{02}:G_0 \backslash(\Gamma \times \Gamma) &\xrightarrow{\cong} G_2\backslash \left(G_2 \times^{G_0}\left(\Gamma \times X\right)\right)=\quot{X}{\Gamma_{02}} \\
		h_{\Gamma}^{03}:G_0 \backslash(\Gamma \times \Gamma) &\xrightarrow{\cong} G_3\backslash \left(G_3 \times^{G_0}\left(\Gamma \times X\right)\right)=\quot{X}{\Gamma_{03}} \\
		h_{\Gamma}^{12}:G_1\backslash\left(\hat{\Gamma} \times X\right) &\xrightarrow{\cong} G_2 \backslash \left(G_2 \times^{G_1}\left(\hat{\Gamma} \times X\right)\right) = \quot{X}{\hat{\Gamma}_{12}} \\
		h_{\hat{\Gamma}}^{13}:G_1 \backslash\left(\hat{\Gamma} \times X\right) &\xrightarrow{\cong} G_3 \backslash\left(G_3 \times^{G_1} \left(\hat{\Gamma} \times X\right)\right) = \quot{X}{\hat{\Gamma}_{13}} \\
		h^{23}_{\tilde{\Gamma}}:G_2 \backslash\left(\tilde{\Gamma} \times X\right) &\xrightarrow{\cong} G_3\backslash\left(G_3 \times^{G_2}\left(\tilde{\Gamma} \times X\right)\right) = \quot{X}{\tilde{\Gamma}_{23}}
	\end{align*}
	Similarly, we define the following analogues of the maps $\rho$ of Proposition \ref{Prop: Section 3.3: Chagnge of groups interacting with change of grups}:
	\begin{align*}
		\rho_{012}:\Gamma_{012} &\xrightarrow{\cong} \Gamma_{02} \\
		\rho_{013}:\Gamma_{013} &\xrightarrow{\cong} \Gamma_{03} \\
		\rho_{023}:\Gamma_{023} &\xrightarrow{\cong} \Gamma_{03} \\
		\rho_{123}:\hat{\Gamma}_{123} &\xrightarrow{\cong} \hat{\Gamma}_{13}
	\end{align*}
\end{definition}

As a result of the definition of the $\rho_{ijk}$ maps above, we also have two induced isomorphisms $\Gamma_{0123} \to \Gamma_{30}$ for any $\Sf(G_0)$-variety $\Gamma$ based on whether we do our ``cancellation'' inside-out (by first contracting the $G_1$ and then $G_2$-components) or outside-in (by first contracting the $G_2$ and then $G_1$-components):
\[
\begin{tikzcd}
	& \Gamma_{013} \ar[dr]{}{\rho_{013}} \\
	\Gamma_{0123} \ar[rr, swap]{}{\rho_{0,123}} \ar[ur]{}{\rho_{123} \times^{G_0} \id_{\Gamma}} & & \Gamma_{03}
\end{tikzcd}\qquad
\begin{tikzcd}
	& \Gamma_{023} \ar[dr]{}{\rho_{023}} \\
	\Gamma_{0123} \ar[rr, swap]{}{\rho_{012,3}} \ar[ur]{}{\id_{G_3} \times^{G_2} \rho_{012}} & & \Gamma_{03}
\end{tikzcd}
\]
Upon taking $G_3$-quotients, we see that these maps coincide by the universal property of the quotient object and hence give rise to a commuting diagram:
\[
\begin{tikzcd}
	& \quot{X}{\Gamma_{013}} \ar[dr]{}{\overline{\rho}_{013}} \\
	\quot{X}{\Gamma_{0123}} \ar[dr, swap]{}{\overline{\id_{G_3} \times^{G_2} \rho_{012}}} \ar[ur]{}{\overline{\rho_{123} \times^{G_0}\Gamma}} \ar[rr, bend left = 15, ""{name = U}]{}{\overline{\rho}_{0,123}} \ar[rr, bend right = 15, swap, ""{name = D}]{}{\overline{\rho}_{012,3}} & & \quot{X}{\Gamma_{03}} \ar[from = U, to = D, equals, shorten <= 4pt, shorten >= 4pt]{}{} \\
	& \quot{X}{\Gamma_{023}} \ar[ur, swap]{}{\overline{\rho}_{023}}
\end{tikzcd}
\]
Our next two lemmas show that all these morphisms induce the same isomorphisms $\XGamma \xrightarrow{\cong} \quot{X}{\Gamma_{0123}}$ and $h_{\Gamma}^{03}:\XGamma \to \quot{X}{\Gamma_{03}}$. In this sense, they are ``cancellation associative'' and provide a way of weaving between perspectives on our quotients.

\begin{lemma}\label{Lemma: Cancellation Associative}
	Let $G_0 \xrightarrow{\varphi} G_1 \xrightarrow{\varphi} G_2 \xrightarrow{\rho} G_3$ be morphisms of group objects and let $X$ be a left $G_2$-object. Then the natural isomorphisms $h_{\Gamma}^{01},$ $h_{\Gamma}^{02},$ $h_{\Gamma_{01}}^{13}$, and $h_{\Gamma_{02}}^{23}$	
	are cancellation associative in the sense that for every smooth free $G_0$-variety $\Gamma$, the diagrams
	\[
	\begin{tikzcd}
		& \quot{X}{\Gamma_{013}} \ar[dl, swap]{}{\left(\overline{\rho_{123} \times^{G_0} \id_{\Gamma}}\right)^{-1}} \\
		\quot{X}{\Gamma_{0123}} & & \quot{X}{\Gamma_{01}} \ar[ul, swap]{}{h_{\Gamma_{01}}^{13}} \\
		& \XGamma \ar[ur, swap]{}{h_{\Gamma}^{01}} \ar[ul]{}{}
	\end{tikzcd}
	\]
	and
	\[
	\begin{tikzcd}
		& \quot{X}{\Gamma_{023}} \ar[dl, swap]{}{\left(\overline{\id_{G_3} \times^{G_2} \rho_{012}}\right)^{-1}} \\
		\quot{X}{\Gamma_{0123}} & & \quot{X}{\Gamma_{02}} \ar[ul, swap]{}{h_{\Gamma_{02}}^{23}} \\
		& \XGamma \ar[ur, swap]{}{h_{\Gamma}^{02}} \ar[ul]
	\end{tikzcd}
	\]
	both commute and coincide.
\end{lemma}
\begin{proof}
	This is routine by using the universal property of the quotients which appear in the diagrams above.
\end{proof}
Because of the lack of ambiguity, we define the map $h_{\Gamma}^{0123}:\XGamma \to \quot{X}{\Gamma_{0123}}$ by 
\begin{equation}\label{Eqn: HGamma0123}
	h_{\Gamma}^{0123} := \left(\overline{\rho_{123} \times^{G_0} \id_{\Gamma}}\right)^{-1} \circ h_{\Gamma_{01}}^{13} \circ h_{\Gamma}^{01}.
\end{equation}
\begin{lemma}\label{Lemma: The two factorizations of h03Gamma}
	Let $\Gamma$ be a smooth free $G_0$-variety, let $\hat{\Gamma}$ be a smooth free $G_1$-variety, and let $X$ be a $G_3$-variety. Then the diagrams
	\[
	\begin{tikzcd}
		\quot{X}{\hat{\Gamma}} \ar[d, swap]{}{h_{\hat{\Gamma}}^{13}} \ar[r]{}{h_{\hat{\Gamma}}^{12}} & \quot{X}{\hat{\Gamma}_{12}} \ar[d]{}{h^{23}_{\hat{\Gamma}_{12}}} \\
		\quot{X}{\hat{\Gamma}_{13}} & \quot{X}{\hat{\Gamma}_{123}} \ar[l]{}{\overline{\rho}_{123}}
	\end{tikzcd}
	\]
	\[
	\begin{tikzcd}
		\XGamma \ar[r]{}{h_{\Gamma}^{01}} \ar[d, swap]{}{h_{\Gamma}^{03}} & \quot{X}{\Gamma_{01}} \ar[d]{}{h_{\Gamma_{01}}^{13}} \\
		\quot{X}{\Gamma_{03}} & \quot{X}{\Gamma_{013}} \ar[l]{}{\overline{\rho}_{013}}
	\end{tikzcd}
	\]
	and
	\[
	\begin{tikzcd}
		\XGamma \ar[r]{}{h_{\Gamma}^{02}} \ar[d, swap]{}{h_{\Gamma}^{03}} & \quot{X}{\Gamma_{02}} \ar[d]{}{h_{\Gamma_{02}}^{23}} \\
		\quot{X}{\Gamma_{03}} & \quot{X}{\Gamma_{023}} \ar[l]{}{\overline{\rho}_{023}}
	\end{tikzcd}
	\]
	all commute.
\end{lemma}
\begin{proof}
	This is a routine proof using the universal property of each quotient variety in sight.
\end{proof}

Combining these lemmas together allows us to produce the commuting cube:
\begin{equation}\label{Eqn: Commuting cube for chofg assoc}
	\begin{tikzcd}
		& \quot{X}{\Gamma_{02}} \ar[rr]{}{h_{\Gamma_{02}}^{23}} \ar[dd, near start]{}{\overline{\rho}_{012}^{-1}} & & \quot{X}{\Gamma_{023}} \ar[dd]{}{\left(\overline{\id_{G_3} \times^{G_2} \rho_{012}}\right)^{-1}} \\
		\XGamma \ar[ur]{}{h_{\Gamma}^{02}} \ar[dd, swap]{}{h_{\Gamma}^{01}} & & \quot{X}{\Gamma_{03}} \ar[ur, swap]{}{\overline{\rho}_{023}^{-1}} \\
		& \quot{X}{\Gamma_{012}} \ar[rr, near start, swap]{}{h_{\Gamma_{012}}^{23}} & & \quot{X}{\Gamma_{0123}} \\
		\quot{X}{\Gamma_{01}} \ar[rr, swap]{}{h_{\Gamma_{01}}^{13}} \ar[ur]{}{h_{\Gamma_{01}}^{12}} & & \quot{X}{\Gamma_{013}} \ar[ur, swap]{}{\left(\overline{\rho_{123} \times^{G_0} \id_{\Gamma}}\right)^{-1}}
		\ar[from = 2-1, to = 2-3, crossing over, near start]{}{h_{\Gamma}^{03}}
		\ar[from = 2-3, to = 4-3, crossing over, near start]{}{\overline{\rho}_{013}^{-1}}
	\end{tikzcd}
\end{equation}
An immediate consequence of the lemmas above, the commuting cube, and the definition of $h_{\Gamma}^{0123}$ are the following myriad alternate factorizations of $h_{\Gamma}^{0123}$.
\begin{lemma}\label{Lemma: The Many factorizations of h0123}
	For any left $G_3$-object $X$ and any smooth free $G_0$-variety (respectively free $G_0$-space) $\Gamma$ we have the following factorizations of the morphism $h_{\Gamma}^{0123}$:
	\begin{align*}
		h_{\Gamma}^{0123} &= \left(\overline{\rho_{123} \times^{G_0} \id_{\Gamma}}\right)^{-1} \circ h_{\Gamma_{01}}^{13} \circ h_{\Gamma}^{01} = \left(\overline{\id_{G_3}\times^{G_2}\rho_{123} }\right)^{-1} \circ h_{\Gamma_{02}}^{23} \circ h_{\Gamma}^{02} \\
		&= h^{23}_{\Gamma_{012}} \circ \left(\overline{\rho}_{012}\right)^{-1} \circ h_{\Gamma}^{02} = h_{\Gamma_{012}}^{23} \circ h_{\Gamma_{01}} \circ h_{\Gamma}^{01} \\
		&= \left(\overline{\rho_{123} \times^{G_0} \id_{\Gamma}}\right)^{-1} \circ \left(\overline{\rho}_{013}\right)^{-1} \circ h_{\Gamma}^{03} = \left(\overline{\id_{G_3}\times^{G_2}\rho_{123}}\right)^{-1} \circ \left(\overline{\rho}_{023}\right)^{-1} \circ h_{\Gamma}^{03}.
	\end{align*}
\end{lemma}
\begin{proof}
	These different factorizations are simply the total paths starting at $\XGamma$ and ending at $\quot{X}{\Gamma_{0123}}$ in the cube presented in Diagram \ref{Eqn: Commuting cube for chofg assoc}.
\end{proof}

Returning explicitly to the task of verifying the (pseudofunctorial) associativity of the Change of Groups functors, we recall that we must show that the pasting diagram presented in Diagram \ref{Eqn: Pasto 1 for chofg associative} is equal to the pasting diagram presented in Diagram \ref{Eqn: Pasto 2 for chofg associative}. However, the first pasting diagram is ultimately given as the composite
\[
\left(\varphi_{01}^{\sharp} \ast \alpha_{123}\right) \circ \alpha_{0,123}:(\varphi_{23}\circ \varphi_{12} \circ \varphi_{01})^{\sharp} \Rightarrow \left(\varphi_{01}^{\sharp} \circ (\varphi_{23} \circ \varphi_{12})^{\sharp}\right) \Rightarrow \varphi_{01}^{\sharp} \circ \varphi_{12}^{\sharp} \circ \varphi_{23}^{\sharp}
\]
while the second is given by
\[
\left(\alpha_{012} \ast \varphi_{23}^{\sharp}\right) \circ \alpha_{012,3}:(\varphi_{23}\circ \varphi_{12} \circ \varphi_{01})^{\sharp} \Rightarrow  \left((\varphi_{12} \circ \varphi_{01})^{\sharp}  \circ\varphi_{23}^{\sharp}\right) \Rightarrow \varphi_{01}^{\sharp} \circ \varphi_{12}^{\sharp} \circ \varphi_{23}^{\sharp}.
\]
It thus suffices to prove that 
\begin{equation}\label{Eqn: The goal}
	(\varphi_{01}^{\sharp} \ast \alpha_{123}) \circ \alpha_{0,123} = (\alpha_{012} \ast \varphi_{23}^{\sharp}) \circ \alpha_{012,3}
\end{equation}
in order to conclude that the Change of Groups functors vary pseudofunctorially over the category of smooth algebraic groups. To this end we compute using the proof of Proposition \ref{Prop: Section 3.3: Chagnge of groups interacting with change of grups} that on one hand for any object $\Gamma$ in $\Sf(G_0)$ (respectively any object $\Gamma$ in $\mathbf{Free}(G_0)$)
\begin{align*}
	&\quot{\alpha_{0,123}}{\Gamma} \\
	&= \overline{F}\left(h_{\Gamma}^{01}\right)\left(\overline{F}\left(h_{\Gamma_{01}}^{13}\right)\left(\tau_{\overline{\rho}_{013}}^{A}\right)\right) \circ \quot{\left(\phi_{h^{01},h^{13}} \ast \overline{F}(\overline{\rho}_{013})\right)^{-1}}{\Gamma} \circ \quot{\left(\phi_{h^{13} \circ h^{01},\overline{\rho}_{013}}\right)^{-1}}{\Gamma},
\end{align*}
while on the other hand
\begin{align*}
	&\quot{\alpha_{012,3}}{\Gamma} \\
	&= \overline{F}\left(h_{\Gamma}^{02}\right)\left(\overline{F}\left(h^{23}_{\Gamma_{02}}\right)\left(\tau_{\overline{\rho}_{023}}^{A}\right)\right) \circ \quot{\left(\phi_{h^{02}, h^{23}} \ast \overline{F}(\overline{\rho}_{023})\right)^{-1}}{\Gamma} \circ \quot{\left(\phi_{h^{23} \circ h^{02}, \rho_{023}}\right)^{-1}}{\Gamma}.
\end{align*}
Computing also using that $\overline{F}$ is a pseudofunctor gives that
\begin{align*}
	&\quot{\left(\phi_{h^{01},h^{13}} \ast \overline{F}(\overline{\rho}_{013})\right)^{-1}}{\Gamma} \circ \quot{\left(\phi_{h^{13} \circ h^{01},\overline{\rho}_{013}}\right)^{-1}}{\Gamma} \\
	&= \quot{\left(\phi_{h^{01} \circ h^{13},\rho_{013}} \circ \left(\phi_{h^{01},h^{13}} \ast \overline{F}(\overline{\rho}_{013})\right)\right)^{-1}}{\Gamma} \\
	&= \quot{\left(\phi_{h^{01},\overline{\rho}_{013} \circ h^{13}} \circ \left(\overline{F}\left(h^{01}_{\Gamma}\right) \ast \phi_{h^{13},\rho_{013}}\right)\right)^{-1}}{\Gamma}
	\\
	&= \quot{\left(\overline{F}\left(h^{01}_{\Gamma}\right) \ast \phi_{h^{13},\rho_{013}}\right)^{-1}}{\Gamma}\circ\quot{\left(\phi_{h^{01},\overline{\rho}_{013} \circ h^{13}}\right)^{-1}}{\Gamma} 
\end{align*}
and similarly
\begin{align*}
	&\quot{\left(\phi_{h^{02},h^{23}} \ast \overline{F}(\overline{\rho}_{023})\right)^{-1}}{\Gamma} \circ \quot{\left(\phi_{h^{23} \circ h^{02},\overline{\rho}_{023}}\right)^{-1}}{\Gamma} \\
	&= \quot{\left(\phi_{h^{02} \circ h^{23},\rho_{023}} \circ \left(\phi_{h^{02},h^{23}} \ast \overline{F}(\overline{\rho}_{023})\right)\right)^{-1}}{\Gamma} \\
	&= \quot{\left(\phi_{h^{02},\overline{\rho}_{023} \circ h^{23}} \circ \left(\overline{F}\left(h^{02}_{\Gamma}\right) \ast \phi_{h^{23},\rho_{023}}\right)\right)^{-1}}{\Gamma}
	\\
	&= \quot{\left(\overline{F}\left(h^{02}_{\Gamma}\right) \ast \phi_{h^{23},\rho_{023}}\right)^{-1}}{\Gamma}\circ\quot{\left(\phi_{h^{02},\overline{\rho}_{023} \circ h^{23}}\right)^{-1}}{\Gamma}.
\end{align*}
Putting these observations together gives the calculation
\begin{align*}
	&\quot{\alpha_{0,123}}{\Gamma} \\
	&=\overline{F}\left(h_{\Gamma}^{01}\right)\left(\overline{F}\left(h_{\Gamma_{01}}^{13}\right)\left(\tau_{\overline{\rho}_{013}}^{A}\right)\right)\circ \quot{\left(\overline{F}\left(h^{01}_{\Gamma}\right) \ast \phi_{h^{13},\rho_{013}}\right)^{-1}}{\Gamma}\circ\quot{\left(\phi_{h^{01},\overline{\rho}_{013} \circ h^{13}}\right)^{-1}}{\Gamma} \\
	&=\overline{F}\left(h_{\Gamma}^{01}\right)\left(\overline{F}\left(h_{\Gamma_{01}}^{13}\right)\left(\tau_{\overline{\rho}_{013}}^{A}\right) \circ\left(\phi_{h^{13},\rho_{013}}\right)^{-1}\right) \circ \quot{\left(\phi_{h^{01},\overline{\rho}_{013} \circ h^{13}}\right)^{-1}}{\Gamma}.
\end{align*}
and similarly
\begin{align*}
	&\quot{\alpha_{012,3}}{\Gamma} \\
	&= \overline{F}\left(h_{\Gamma}^{02}\right)\left(\overline{F}\left(h_{\Gamma_{02}}^{23}\right)\left(\tau_{\overline{\rho}_{023}}^{A}\right) \circ\left(\phi_{h^{23},\rho_{023}}\right)^{-1}\right) \circ \quot{\left(\phi_{h^{02},\overline{\rho}_{023} \circ h^{23}}\right)^{-1}}{\Gamma}.
\end{align*}

To complete the computation of the pasting diagrams we note
\begin{align*}
	&\quot{\left(\alpha_{012} \ast \varphi_{23}^{\sharp}\right)}{\Gamma} \\
	&=\overline{F}\left(h_{\Gamma}^{01}\right)\left(\overline{F}\left(h_{\Gamma_{01}}^{12}\right)\left(\tau_{\overline{\rho}_{012}}^{A}\right)\right) \circ \left(\phi_{h^{01},h^{12}}^{\varphi_{23}^{\sharp}A} \ast \overline{F}(\overline{\rho}_{012})\right)^{-1} \circ \left(\phi_{h^{12} \circ h^{01}, \overline{\rho}_{012}}^{\varphi_{23}^{\sharp}A}\right)^{-1} \\
	&= \overline{F}\left(h_{\Gamma}^{01}\right)\left(\overline{F}\left(h_{\Gamma_{01}}^{12}\right)\left(\tau_{\overline{\rho}_{012}}^{A}\right)\right) \circ \left(\phi_{h^{12} \circ h^{01}, \overline{\rho}_{012}}^{\varphi_{23}^{\sharp}A} \circ \left(\phi_{h^{01},h^{12}}^{\varphi_{23}^{\sharp}A} \ast \overline{F}(\overline{\rho}_{012})\right)\right)^{-1} \\
	&=\overline{F}\left(h_{\Gamma}^{01}\right)\left(\overline{F}\left(h_{\Gamma_{01}}^{12}\right)\left(\tau_{\overline{\rho}_{012}}^{A}\right)\right) \circ \left(\phi_{h^{01},\overline{\rho}_{012} \circ h^{12}}^{\varphi_{23}^{\sharp}A} \circ \left(\overline{F}\left(h^{01}_{\Gamma}\right) \ast \phi_{h^{23},\overline{\rho}_{023}}^{\varphi_{23}^{\sharp}A}\right)\right)^{-1} \\
	&=\overline{F}\left(h_{\Gamma}^{01}\right)\left(\overline{F}\left(h_{\Gamma_{01}}^{12}\right)\left(\tau_{\overline{\rho}_{012}}^{A}\right)\right) \circ \left(\overline{F}\left(h^{01}_{\Gamma}\right) \ast \phi_{h^{23},\overline{\rho}_{023}}^{\varphi_{23}^{\sharp}A}\right)^{-1} \circ \left(\phi_{h^{01},\overline{\rho}_{012} \circ h^{12}}^{\varphi_{23}^{\sharp}A}\right)^{-1} \\
	&=\overline{F}\left(h_{\Gamma}^{01}\right)\left(\overline{F}\left(h_{\Gamma_{01}}^{12}\right)\left(\tau_{\overline{\rho}_{012}}^{A}\right) \circ \phi_{h^{23},\overline{\rho}_{023}}^{\varphi_{23}^{\sharp}A}\right)^{-1} \circ \left(\phi_{h^{01},\overline{\rho}_{012} \circ h^{12}}^{\varphi_{23}^{\sharp}A}\right)^{-1}
\end{align*}
and, alternatively,
\begin{align*}
	&\quot{\left(\varphi_{01}^{\sharp} \ast \alpha_{123}\right)}{\Gamma} =\overline{F}\left(h_{\Gamma}^{01}\right)\left(\quot{\alpha_{123}}{\Gamma_{01}}\right) \\
	&= \overline{F}\left(h_{\Gamma}^{01}\right)\left(\overline{F}\left(h_{\Gamma_{01}}^{12}\right)\left(\overline{F}\left(h_{\Gamma_{012}}^{23}\right)\left(\tau_{\overline{\rho}_{123}}^{A}\right)\right) \circ \left(\phi_{h^{23} \circ h^{12},\overline{\rho}_{123}} \circ \left(\phi_{h^{12},h^{23}} \ast \overline{F}(\overline{\rho}_{123})\right)\right)^{-1}\right) \\
	&= \overline{F}\left(h_{\Gamma}^{01}\right)\left(\overline{F}\left(h_{\Gamma_{01}}^{12}\right)\left(\overline{F}\left(h_{\Gamma_{012}}^{23}\right)\left(\tau_{\overline{\rho}_{123}}^{A}\right)\right) \circ \left(\phi_{h^{23} \circ h^{12},\overline{\rho}_{123}} \circ \left(\overline{F}\left(h_{\Gamma_{01}}^{12}\right) \ast \phi_{h^{23}, \overline{\rho}_{123}}\right)\right)^{-1}\right) \\
	&= \overline{F}\left(h_{\Gamma}^{01}\right)\left(\overline{F}\left(h_{\Gamma_{01}}^{12}\right)\left(\overline{F}\left(h_{\Gamma_{012}}^{23}\right)\left(\tau_{\overline{\rho}_{123}}^{A}\right) \circ \left(\phi_{h^{23}, \overline{\rho}_{123}}\right)^{-1}\right) \circ \left(\phi_{h^{12}, \overline{\rho}_{123} \circ h^{23}}\right)^{-1} \right)
\end{align*}
Observe that in the last line we used the identities
\begin{align*}
	&\left(\phi_{h^{23} \circ h^{12},\overline{\rho}_{123}} \circ \left(\phi_{h^{12},h^{23}} \ast \overline{F}(\overline{\rho}_{123})\right)\right)^{-1} \\
	&=\left(\phi_{h^{12},\overline{\rho}_{123} \circ h^{23}} \circ \left(\overline{F}\left(h_{\Gamma_{01}}^{12}\right) \ast \phi_{h^{23}, \overline{\rho}_{123}}\right)\right)^{-1} \\
	&= \overline{F}\left(h_{\Gamma_{01}}^{12}\right)\left(\phi_{h^{23}, \overline{\rho}_{123}}\right)^{-1} \circ \left(\phi_{h^{12},\overline{\rho}_{123} \circ h^{23}}\right)^{-1}.
\end{align*}
Combining these above calculations give the lemma below.
\begin{lemma}\label{Lemma: The explicit Description of the pasting diagrams}
	For any smooth free $G_0$-variety $\Gamma$ (respectively any free $G_0$-space $\Gamma$), the equations
	\begin{align*}
		&\quot{\left(\varphi_{01}^{\sharp} \ast \alpha_{123}\right)}{\Gamma} \circ \quot{\alpha_{0,123}}{\Gamma} \\
		&= \overline{F}\left(h_{\Gamma}^{01}\right)\left(\overline{F}\left(h_{\Gamma_{01}}^{12}\right)\left(\overline{F}\left(h_{\Gamma_{012}}^{23}\right)\left(\tau_{\overline{\rho}_{123}}^{A}\right) \circ \left(\phi_{h^{23}, \overline{\rho}_{123}}\right)^{-1}\right) \circ \left(\phi_{h^{12}, \overline{\rho}_{123} \circ h^{23}}\right)^{-1} \right) \\
		&\circ \overline{F}\left(h_{\Gamma}^{01}\right)\left(\overline{F}\left(h_{\Gamma_{01}}^{13}\right)\left(\tau_{\overline{\rho}_{013}}^{A}\right) \circ\left(\phi_{h^{13},\rho_{013}}\right)^{-1}\right) \circ \quot{\left(\phi_{h^{01},\overline{\rho}_{013} \circ h^{13}}\right)^{-1}}{\Gamma}
	\end{align*}
	and
	\begin{align*}
		&\quot{\left(\alpha_{012} \ast \varphi_{23}^{\sharp}\right)}{\Gamma} \circ \quot{\alpha_{012,3}}{\Gamma} \\
		&= \overline{F}\left(h_{\Gamma}^{01}\right)\left(\overline{F}\left(h_{\Gamma_{01}}^{12}\right)\left(\tau_{\overline{\rho}_{012}}^{A}\right) \circ \phi_{h^{23},\overline{\rho}_{023}}^{\varphi_{23}^{\sharp}A}\right)^{-1} \circ \left(\phi_{h^{01},\overline{\rho}_{012} \circ h^{12}}^{\varphi_{23}^{\sharp}A}\right)^{-1}	\\
		&\circ \overline{F}\left(h_{\Gamma}^{02}\right)\left(\overline{F}\left(h^{23}_{\Gamma_{02}}\right)\left(\tau_{\overline{\rho}_{023}}^{A}\right)\right) \circ \quot{\left(\phi_{h^{02}, h^{23}} \ast \overline{F}(\overline{\rho}_{023})\right)^{-1}}{\Gamma} \circ \quot{\left(\phi_{h^{23} \circ h^{02}, \rho_{023}}\right)^{-1}}{\Gamma}
	\end{align*}
	hold.
\end{lemma}

\begin{Theorem}\label{Thm: Section Change of Groups: Associativity of chofg compositors}
	The pasting diagram
	\begin{equation*}
		\begin{tikzcd}
			& F_{G_2}(X) \ar[dr]{}{\varphi_{12}^{\sharp}} \\
			F_{G_3}(X) \ar[ur]{}{\varphi_{23}^{\sharp}} \ar[dr, swap]{}{(\varphi_{23} \circ \varphi_{12} \circ \varphi_{01})^{\sharp}} \ar[rr, ""{name = M}]{}[description]{(\varphi_{23} \circ \varphi_{12})^{\sharp}} & & F_{G_1}(X) \ar[dl]{}{\varphi_{01}^{\sharp}} 	 \\
			& F_{G_0}(X) \ar[from = 3-2, to = M, Rightarrow, shorten <= 4pt, shorten >= 4pt]{}{\alpha_{0,123}} \ar[from = M, to = 1-2, Rightarrow, shorten <= 4pt, shorten >= 4pt]{}{\alpha_{123}}
		\end{tikzcd}
	\end{equation*}
	is equal to the pasting diagram:
	\begin{equation*}
		\begin{tikzcd}
			& F_{G_2}(X) \ar[dr]{}{\varphi_{12}^{\sharp}} \ar[dd, ""{name = M}]{}{}  \ar[dd, crossing over, near start]{}[description]{(\varphi_{12} \circ \varphi_{01})^{\sharp}} \\
			F_{G_3}(X) \ar[ur]{}{\varphi_{23}^{\sharp}} \ar[dr, swap]{}{(\varphi_{23} \circ \varphi_{12} \circ \varphi_{01})^{\sharp}}  & & F_{G_1}(X) \ar[dl]{}{\varphi_{01}^{\sharp}} 	 \\
			& F_{G_0}(X)  \ar[from = M, to = 2-3, Rightarrow, swap, shorten <= 4pt, shorten <= 4pt]{}{\alpha_{012}} \ar[from = 2-1, to = M, swap, Rightarrow, shorten <= 4pt, shorten >= 4pt]{}{\alpha_{012,3}}
		\end{tikzcd}
	\end{equation*}
	coincide.
\end{Theorem}
\begin{proof}
	As explained earlier (cf.\@ Equation \ref{Eqn: The goal}), it suffices to prove that the equation
	\[
	(\varphi_{01}^{\sharp} \ast \alpha_{123}) \circ \alpha_{0,123} = (\alpha_{012} \ast \varphi_{23}^{\sharp}) \circ \alpha_{012,3}
	\]
	holds. However, begin by noting that since
	\[
	\left(\overline{\rho_{123} \times^{G_0} \id_{\Gamma}}\right)^{-1} \circ h_{\Gamma_{01}}^{13} \circ h_{\Gamma}^{01} = h_{\Gamma}^{0123} = \left(\overline{\id_{G_3} \times^{G_2} \rho_{012}}\right)^{-1} \circ h_{\Gamma_{02}}^{23} \circ h_{\Gamma}^{02}
	\]
	by Lemma \ref{Lemma: The Many factorizations of h0123} we have
	\[
	\left(\overline{\id_{G_3} \times^{G_2} \rho_{012}}\right) \circ \left(\overline{\rho_{123} \times^{G_0} \id_{\Gamma}}\right)^{-1} \circ h_{\Gamma_{01}}^{13} \circ h_{\Gamma}^{01} = h_{\Gamma_{02}}^{23} \circ h_{\Gamma}^{02}
	\]
	and, dually,
	\[
	h_{\Gamma_{01}}^{13} \circ h_{\Gamma}^{01} =  \left(\overline{\rho_{123} \times^{G_0} \id_{\Gamma}}\right) \circ\left(\overline{\id_{G_3} \times^{G_2} \rho_{012}}\right)^{-1} \circ h_{\Gamma_{02}}^{23} \circ h_{\Gamma}^{02}
	\]
	Starting from here we also derive that
	\begin{align*}
		&\left(\overline{\rho_{123} \times^{G_0} \id_{\Gamma}}\right) \circ\left(\overline{\id_{G_3} \times^{G_2} \rho_{012}}\right)^{-1} \circ h_{\Gamma_{02}}^{23} \circ h_{\Gamma}^{02} \\
		&= \left(\overline{\rho_{123} \times^{G_0} \id_{\Gamma}}\right) \circ h_{\Gamma_{012}}^{23} \circ \left(\overline{\rho}_{012}\right)^{-1} \circ h_{\Gamma}^{02} \\
		&= \left(\overline{\rho_{123} \times^{G_0} \id_{\Gamma}}\right) \circ h_{\Gamma_{012}}^{23} \circ h_{\Gamma_{01}}^{12} \circ h_{\Gamma}^{01}.
	\end{align*}
	Using the psuedofunctoriality of $\overline{F}$, the various rewritings
	\[
	\phi_{b \circ a, c} \circ (\phi_{a,b} \ast \overline{F}(c)) = \phi_{a, c \circ b} \circ (\overline{F}(a) \ast \phi_{b,c})
	\] 
	of the compositors when necessary, myriad identities implied by Diagram \ref{Eqn: Commuting cube for chofg assoc}, and performing an extremely tedious algebraic chase allows us to derive that
	\begin{align*}
		&\overline{F}\left(h_{\Gamma}^{01}\right)\left(\overline{F}\left(h_{\Gamma_{01}}^{12}\right)\left(\overline{F}\left(h_{\Gamma_{012}}^{23}\right)\left(\tau_{\overline{\rho}_{123}}^{A}\right) \circ \left(\phi_{h^{23}, \overline{\rho}_{123}}\right)^{-1}\right) \circ \left(\phi_{h^{12}, \overline{\rho}_{123} \circ h^{23}}\right)^{-1} \right) \\
		&\circ \overline{F}\left(h_{\Gamma}^{01}\right)\left(\overline{F}\left(h_{\Gamma_{01}}^{13}\right)\left(\tau_{\overline{\rho}_{013}}^{A}\right) \circ\left(\phi_{h^{13},\rho_{013}}\right)^{-1}\right) \circ \quot{\left(\phi_{h^{01},\overline{\rho}_{013} \circ h^{13}}\right)^{-1}}{\Gamma} \\
		&= \overline{F}\left(h_{\Gamma}^{01}\right)\left(\overline{F}\left(h_{\Gamma_{01}}^{12}\right)\left(\tau_{\overline{\rho}_{012}}^{A}\right) \circ \phi_{h^{23},\overline{\rho}_{023}}^{\varphi_{23}^{\sharp}A}\right)^{-1} \circ \left(\phi_{h^{01},\overline{\rho}_{012} \circ h^{12}}^{\varphi_{23}^{\sharp}A}\right)^{-1}	\\
		&\circ \overline{F}\left(h_{\Gamma}^{02}\right)\left(\overline{F}\left(h^{23}_{\Gamma_{02}}\right)\left(\tau_{\overline{\rho}_{023}}^{A}\right)\right) \circ \quot{\left(\phi_{h^{02}, h^{23}} \ast \overline{F}(\overline{\rho}_{023})\right)^{-1}}{\Gamma} \circ \quot{\left(\phi_{h^{23} \circ h^{02}, \rho_{023}}\right)^{-1}}{\Gamma}.
	\end{align*}
	However, this exactly allows us to deduce the identity
	\[
	(\varphi_{01}^{\sharp} \ast \alpha_{123}) \circ \alpha_{0,123} = (\alpha_{012} \ast \varphi_{23}^{\sharp}) \circ \alpha_{012,3},
	\]
	as was desired.
\end{proof}

As an immediate corollary to this theorem we can regard the Change of Groups functors as the fibre functors for a pseudofunctor with domain $\mathbf{SAlgGrp} \downarrow H$ in the geometric case and $\mathbf{TopGrp} \downarrow H$ in the topological case.
\begin{corollary}\label{Cor: Finally Chofg is pseudofunctorial}
	Let $H$ be a smooth algebraic group (respectively a topological group), let $X$ be a left $H$-variety (respectively a left $H$-space), and let $(F,\overline{F})$ be an $H$-pre-equivariant psuedofunctor on $X$. Then there is a pseudofunctor
	\[
	F_{-}(X):\left(\mathbf{SAlgGrp} \downarrow H\right)^{\op} \to \fCat
	\]
	where:
	\begin{itemize}
		\item For each smooth algebraic group morphism $\varphi:G \to H$, $F_{-}(X)(\varphi) := F_G(X)$ for the $(G,H)$-pre-equivariant pseudofunctor
		\[
		\left(\left(\overline{F}\left(G \backslash (H \times^{G} ((-) \times X))\right), \overline{F}\right), (F,\overline{F}), \id_{\overline{F}}\right);
		\]
		\item For each morphism $\varphi:G \to G^{\prime}$ in $\mathbf{SAlgGrp} \downarrow H$, the fibre functor is given by the Change of Groups functor
		\[
		\varphi^{\sharp}:F_{G^{\prime}}(X) \to F_G(X).
		\]
		\item For each composable pair of morphisms $\varphi_{01}:G_0 \to G_1$ and $\varphi_{12}:G_1 \to G_2$ in $\mathbf{SAlgGrp} \downarrow H$, the compositor $\phi_{\varphi_{01}, \varphi_{12}}$ is defined by
		\[
		\phi_{\varphi_{01}, \varphi_{12}} := \alpha_{012}^{-1},
		\]
		where $\alpha_{012}$ is the natural isomorphism of Proposition \ref{Prop: Section 3.3: Chagnge of groups interacting with change of grups}.
	\end{itemize}
	Respectively, there is also a pseudofunctor $F_{-}(X):\left(\mathbf{TopGrp} \downarrow H\right)^{\op} \to \fCat$ whose assignment is given similarly.
\end{corollary}
\begin{proof}
	The only difficult aspect of this corollary lies in verifying that the compositors vary pseudofunctorially in the category $(\mathbf{SAlgGrp} \downarrow H)^{\op}$ for the geometric case and in the category $(\mathbf{TopGrp} \downarrow H)^{\op}$ in the topological case; however, this is exactly the content of Theorem \ref{Thm: Section Change of Groups: Associativity of chofg compositors}.
\end{proof}

This allows us to formulate two different perspectives on Change of Groups: first a more general version valid for any $(G,H)$-pre-equivariant pseudofunctor and second a perspective given in terms of the $2$-category $\GHPreq{-}{H}$. 

In the first case assume that we are in either of the situations below:
\begin{enumerate}
	\item smooth algebraic groups $G$ and $H$, a group morphism $\varphi:G \to H$, and also a left $H$-variety $X$;
\item topological groups $G$ and $H$, a group morphism $\varphi:G \to H$, and also a left $H$-space $X$;
\end{enumerate}
Now consider a $(G,H)$-pre-equivariant pseudofunctor $\left((F,\overline{F}),(E,\overline{E}),e:\overline{F}\Rightarrow \overline{E}\right).$ We then have two distinct ways to define the Change of Groups functors $\varphi^{\sharp}:F_H(X) \to E_G(X)$: the first is via the composition
\[
\begin{tikzcd}
	F_H(X) \ar[r]{}{\quot{\varphi^{\sharp}}{F}} & F_G(X) \ar[r, equals] & \PC(\overline{F}\circ\quo_G^{\op}) \ar[d]{}{\PC(e \ast \quo_G^{\op})}\\ 
	& E_G(X) & \PC(\overline{E} \circ \quo_G^{\op}) \ar[l, equals]{}{} &
\end{tikzcd}
\]
while the second is via
\[
\begin{tikzcd}
	F_H(X) \ar[r, equals]{}{} & \PC(\overline{F}\circ \quo_H^{\op}) \ar[rr]{}{\PC(e \ast \quo_H^{\op})} & & \PC(\overline{E} \circ \quo_H^{\op}) \ar[d]{}{\quot{\varphi^{\sharp}}{E}} \\
	& & E_G(X) & \PC(\overline{E} \circ \quo_G^{\op}) \ar[l, equals]{}{} & 
\end{tikzcd}
\]
where we use the left-handed subscripts $\quot{\varphi^{\sharp}}{F}$ and $\quot{\varphi^{\sharp}}{E}$ to denote the pseudofunctor with which we change groups. The first composite acts on objects $A$ by defining, for each $\Gamma \in \Sf(G)_0$ (respectively for each $\Gamma \in \mathbf{Free}(G)_0$),
\[
\quot{\left(\PC(e \ast \quo_G^{\op}) \circ \quot{\varphi^{\sharp}}{F}\right)(A)}{\Gamma} = \quot{e}{\Gamma}\left(\overline{F}\left(h_{\Gamma}\right)\left(\quot{A}{H \times^G \Gamma}\right)\right)
\]
while the second acts by setting
\[
\quot{\left(\quot{\varphi^{\sharp}}{E} \circ \PC(e \ast \quo_H^{\op})\right)}{\Gamma} = \overline{E}\left(h_{\Gamma}\right)\left(\quot{e}{\Gamma}\left(\quot{A}{H \times^G \Gamma}\right)\right).
\]
These objects are then equal by the pseudofunctoriality of the Change of Groups functors together with the fact that by definition $e$ is a pseudonatural equivalence from $\overline{F}$ to $\overline{E}$. In particular, for a general $(G,H)$-pre-equivariant pseudofunctor we get an unambiguous definition for our Change of Groups functors.
\begin{lemma}
Assume we are in either of the cases below.
\begin{enumerate}
	\item smooth algebraic groups $G$ and $H$, a group morphism $\varphi:G \to H$, and also a left $H$-variety $X$;
	\item topological groups $G$ and $H$, a group morphism $\varphi:G \to H$, and also a left $H$-space $X$.
\end{enumerate}
If $\left((F,\overline{F}),(E,\overline{E}),e:\overline{F} \Rightarrow \overline{E}\right)$ is a $(G,H)$-pre-equivariant pseudofunctor then
	\[
	\PC(e \ast \quo_{G}^{\op}) \circ \quot{\varphi^{\sharp}}{F} = \quot{\varphi^{\sharp}}{E} \circ \PC(e \ast \quo_H^{\op}).
	\]
\end{lemma}
\begin{definition}\label{Defn: Change of Groups, General}
Assume we are in either of the cases below.
\begin{enumerate}
	\item smooth algebraic groups $G$ and $H$, a group morphism $\varphi:G \to H$, and also a left $H$-variety $X$;
	\item topological groups $G$ and $H$, a group morphism $\varphi:G \to H$, and also a left $H$-space $X$.
\end{enumerate} 
Then for any $(G,H)$-pre-equivariant pseudofunctor $\left((F,\overline{F}),(E,\overline{E}),e:\overline{F} \Rightarrow \overline{E}\right)$ on $X$, the Change of Groups functor 
	\[
	\varphi^{\sharp}:F_H(X) \to E_G(X)
	\]
	is given by $\varphi^{\sharp}:= \PC(e \ast \quo_{G}^{\op}) \circ \quot{\varphi^{\sharp}}{F}.$
\end{definition}

Using this we can now define our primary four $2$-functors $\GHPreq{-}{H} \to \fCat$. Define the functors as follows:
\begin{itemize}
	\item The functor $\PC_1:\GHPreq{-}{H} \to \fCat$ is induced by the assignment
	\[
	\left((F,\overline{F}),(E,\overline{E}),e\right) \mapsto \PC(\overline{F})
	\]
	on objects.
	\item The functor $\PC_2:\GHPreq{-}{H} \to \fCat$ is induced by the assignment
	\[
	\left((F,\overline{F}),(E,\overline{E}),e\right) \mapsto \PC(\overline{E})
	\]
	on objects.
	\item The functor $(-)_{H}(X):\GHPreq{-}{H} \to \fCat$ is induced by the assignment
	\[
	\left((F,\overline{F}),(E,\overline{E}),e\right) \mapsto F_H(X)
	\]
	on objects.
	\item The functor $(-)_{\square}(X) \to \fCat$ is induced by the assignment
	\[
	\left((F,\overline{F}),(E,\overline{E}),e\right) \mapsto E_G(X)
	\]
	where $G$ is the group for which the object quintuple is a $(G,H)$-pre-equivariant pseudofunctor.
\end{itemize}
By construction there is a (pseudo)natural equivalence of $2$-functors
\[
\begin{tikzcd}
	\GHPreq{-}{H} \ar[rrrr, bend left = 20, ""{name = U}]{}{\PC_1} \ar[rrrr, bend right = 20, swap, ""{name = D}]{}{\PC_2} & & & & \fCat \ar[from = U, to = D, Rightarrow, shorten <= 4pt, shorten >= 4pt]{}{\simeq}
\end{tikzcd}
\]
which is induced by $e$. Our high-level perspective on the Change of Groups functors is that they give rise to a pseudonatural transformation
\[
\begin{tikzcd}
	\GHPreq{-}{H} \ar[rrrr, bend left = 20, ""{name = U}]{}{(-)_{H}(X)} \ar[rrrr, bend right = 20, swap, ""{name = D}]{}{(-)_{\square}(X)} & & & & \fCat \ar[from = U, to = D,  Rightarrow, shorten <= 4pt, shorten >=4pt]{}{(-)^{\sharp}}
\end{tikzcd}
\]
which we describe in the next proposition.

\begin{proposition}\label{Prop: Pseudonat Chofg perspective}
Assume we are in either of the cases below.
\begin{enumerate}
	\item smooth algebraic groups $G$ and $H$, a group morphism $\varphi:G \to H$, and also a left $H$-variety $X$;
	\item topological groups $G$ and $H$, a group morphism $\varphi:G \to H$, and also a left $H$-space $X$.
\end{enumerate}
There is then a pseudonatural transformation
	\[
	\begin{tikzcd}
		\GHPreq{-}{H} \ar[rrrr, bend left = 20, ""{name = U}]{}{(-)_{H}(X)} \ar[rrrr, bend right = 20, swap, ""{name = D}]{}{(-)_{\square}(X)} & & & & \fCat \ar[from = U, to = D, Rightarrow, shorten <= 4pt, shorten >=4pt]{}{(-)^{\sharp}}
	\end{tikzcd}
	\]
	induced by declaring the object functors $\quot{(-)^{\sharp}}{\ul{F}}$ for a $(G,H)$-pre-equivariant pseudofunctor
	\[
	\ul{F} = \left((F,\overline{F}),(E,\overline{E}),e:\overline{F} \xRightarrow{\simeq} \overline{E}\right)
	\]
	to be given by the Change of Groups functor $\varphi^{\sharp}:F_H(X) \to E_G(X)$.
\end{proposition}
\begin{proof}
	The witness natural isomorphisms and transformations are induced by the various pseudofunctor witness isomorphisms associated to the Change of Groups. That this determines a pseudonatural transformation is a straightforward and routine verification and omitted.
\end{proof}

\newpage

\section{The Induction and Quotient Equivalences}\label{Section: Quot and Ind Equiv}

We now present a short section giving sanity check results by extending what are known as the induction and quotient equivalences to our pseudocone formalism in both the topological and geometric cases. It has been known since \cite{BernLun} in the topological case that for any topological group $G$ with topological subgroup $H \leq G$, for any $G$-space $X$ there are equivalences of categories
\[
D_H^b(X) \simeq D_G^b\left(G \times^{H} X\right)
\]
and, whenever $X$ is a free $G$-space,
\[
D_G^b(X) \simeq D_{G/H}^b\left(H \backslash X\right)
\]
which are called the induction and quotient equivalences, respectively. The proofs of these results in \cite{BernLun} are given entirely in terms of the resolution categories (cf.\@ Lemmas \ref{Lemma: Section IndQuot Equiv: GroupResls for quotient space is Resls for space} and \ref{Lemma: Section IndQuot Equiv: GroupResls for induction space is SubgroupResls for restriction action}), and so carry over to the pseudocone language freely; consequently, our main goal is to derive these for the variety-theoretic case and also give suggestions for how to further generalize these results to other situations.

\begin{lemma}[\cite{BernLun}]\label{Lemma: Section IndQuot Equiv: GroupResls for quotient space is Resls for space}
Let $G$ be a topological group and let $H \trianglelefteq G$ be a normal topological subgroup with $X$ a free $G$-space, there is an equivalence of categories
\[
\Resl_G(X) \simeq \Resl_{G/H}(H\backslash X).
\]
\end{lemma}
\begin{proof}[Sketch]
This is argued in the proof of the unnumbered theorem in \cite[Section 2.6.2]{BernLun}.
\end{proof}
\begin{lemma}[\cite{BernLun}]\label{Lemma: Section IndQuot Equiv: GroupResls for induction space is SubgroupResls for restriction action}
Let $H$ be a topological subgroup of a topological group $G$ with inclusion $i:H \to G$ and let $X$ be a left $G$-space. Then there is an equivalence of categories
\[
\Resl_G\left(G \times^H X\right) \simeq \Resl_H(X).
\]
\end{lemma}
\begin{proof}[Sketch]
This is sketched in the proof of the unnumbered theorem in \cite[Section 2.6.3]{BernLun}.
\end{proof}

We will give the geometric analogues of these Lemmas below and illustrate more-or-less how they are proved; for the moment we simply need to know that they hold. From these lemmas we get the propositions below that give the topological versions of the quotient and induction equivalences.

\begin{proposition}\label{Prop: Section IndQuot Equiv: Quotient equiv for top}
Let $G$ be a topological group, $X$ a free $G$-space, and $H \trianglelefteq G$ a normal topological subgroup with inclusion $i:H \to G$. Then for any pseudofunctor $\overline{F}:\Top^{\op} \to \fCat$ there is an equivalence of categories
\[
F_G(X) \simeq F_{G/H}(H \backslash X)
\]
for the induced pre-equivariant pseudofunctors on $X$ and $H \backslash X$, respectively.
\end{proposition}
\begin{proof}
Because of the pseudofunctor $\overline{F}:\Top^{\op} \to \fCat$, we get corresponding induced pre-equivariant pseudofunctors $(\overline{F} \circ \quo_G^{\op}, \overline{F})$ and $(\overline{F} \circ \quo_{G/H}^{\op})$ on $X$ and $H \backslash X$, respectively. Then if we write
\[
e:\Resl_G(X) \xrightarrow{\simeq} \Resl_{G/H}(H \backslash X)
\]
for one direction of the equivalence of Lemma \ref{Lemma: Section IndQuot Equiv: GroupResls for quotient space is Resls for space}, we get that the functor 
\[
\PC(\overline{F} \ast e):\PC(\overline{F}\circ \quo_{G}^{\op}) \xrightarrow{\simeq} \PC(\overline{F} \circ \quo_{G/H}^{\op})
\]
is an equivalence. Because $\PC(\overline{F}\circ \quo_G^{\op}) = F_G(X)$ and $\PC(\overline{F} \ast \quo_{G/H}^{\op}) = F_{G/H}(H \backslash X)$, the result follows.
\end{proof}
\begin{proposition}\label{Prop: Section IndQuot Equiv: Ind equiv for top}
	Let $G$ be a topological group, $X$ a left $G$-space, and $H \leq G$ a topological subgroup with inclusion $i:H \to G$. Then for any pseudofunctor $\overline{F}:\Top^{\op} \to \fCat$ there is an equivalence of categories
	\[
	F_H(X) \simeq F_{G}\left(G \times^{H} X\right)
	\]
	for the induced pre-equivariant pseudofunctors on $X$ and $G \times^H X$, respectively.
\end{proposition}
\begin{proof}
	Because of the pseudofunctor $\overline{F}:\Top^{\op} \to \fCat$, we get corresponding induced pre-equivariant pseudofunctors $(\overline{F} \circ \quo_H^{\op}, \overline{F})$ and $(\overline{F} \circ \quo_{G}^{\op})$ on $X$ and $G \times^{H} X$, respectively. Then if we write
	\[
	e:\Resl_H(X) \xrightarrow{\simeq} \Resl_{G}\left(G \times^{H} X\right)
	\]
	for one direction of the equivalence of Lemma \ref{Lemma: Section IndQuot Equiv: GroupResls for induction space is SubgroupResls for restriction action}, we get that the functor 
	\[
	\PC(\overline{F} \ast e):\PC(\overline{F}\circ \quo_{H}^{\op}) \xrightarrow{\simeq} \PC(\overline{F} \circ \quo_{G}^{\op})
	\]
	is an equivalence. Because $\PC(\overline{F}\circ \quo_G^{\op}) = F_G(X)$ and $\PC(\overline{F} \ast \quo_{G/H}^{\op}) = F_{G/H}(H \backslash X)$, the result follows.
\end{proof}

We now prove the geometric analogues of the results above by following the same basic techniques as before: providing equivalences between resolution categories. Our first requirement is to show the induction resolution equivalences. Fix a smooth algebraic group $G$ with a left $G$-variety $X$ and with an algebraic subgroup $i:H \to G$. Upon restricting the $G$-action on $X$ to an $H$-action, we obtain a morphism $X \to G \times^{H} X$ given by $x \mapsto [1_G,x]_H$ on points. 

\begin{lemma}\label{Lemma: Section IndQuot Equiv: GroupResls for Induction Variety}
Let $G$ be a smooth algebraic group, let $X$ be a left $G$-variety, and let $I:H \to G$ be a smooth algebraic subgroup. Then there is an equivalence of categories
\[
\Resl_H(X) \simeq \Resl_G\left(G \times^H X\right).
\] 
\end{lemma}
\begin{proof}
Because $H$ is a group subvariety of $G$, we find that the induced functor $\Resl_H(X) \to \Resl_G(G \times^H X)$ given by $P \mapsto G \times^H P$ is fully faithful. To see essential surjectivity, let $Q \to G \times^H X$ be a locally isotrivial $G$-fibration and write the corresponding quotient maps as $q:Q \to G\backslash Q$. Because $Q$ is a locally isotrivial $G$-fibration, we can find a cover $\lbrace \varphi_i:U_i \to G \backslash (G \times^H X) \; | \; i \in I \rbrace$ of finite {\'e}tale opens for which there are commuting diagrams
\[
\begin{tikzcd}
q^{-1}(U_i) \ar[rr]{}{\cong} \ar[dr] & & G \times U_i \ar[dl]{}{\pi_2} \\
 & U_i
\end{tikzcd}
\]
for each $i \in I$. But now observe that $G \backslash Q \cong H \backslash i^{\ast}Q$ where $i^{\ast}Q$ takes the variety $Q$ and regards it as an $H$-variety via the restriction action. Then the composition
\[
U_i \xrightarrow{\varphi_i} G \backslash Q \xrightarrow{\cong} H \backslash i^{\ast}Q
\]
is a finite {\'e}tale open and the maps $\lbrace U_i \to H \backslash i^{\ast}Q \; | \; i \in I \rbrace$ constitute a cover by virtue of the facts that covers compose and isomorphisms are covers. Write $q^{\prime}:i^{\ast}Q \to H \backslash i^{\ast}Q$ for the quotient map. It is then routine to check that in the pullback diagram
\[
\begin{tikzcd}
(q^{\prime})^{-1}(U_i) \ar[r] \ar[d] & i^{\ast}Q \ar[d]{}{q^{\prime}} \\
U_i \ar[r] & H \backslash i^{\ast}Q
\end{tikzcd}
\]
we have an isomorphism $(q^{\prime})^{-1}(U_i) \cong H \times U_i$ which fits into the diagram:
\[
\begin{tikzcd}
	(q^{\prime})^{-1}(U_i) \ar[rr]{}{\cong} \ar[dr] & & H \times U_i \ar[dl]{}{\pi_2} \\
	& U_i
\end{tikzcd}
\]
It then follows that $i^{\ast}Q$ is a locally isotrival $H$-variety. Moreover, by construction, we also have that the map $Q \to G \times^H i^{\ast}Q$ is an isomorphism $G \times^H i^{\ast}Q \cong Q$ given by $q \mapsto [1_G,q]_H$. Thus the functor $\Resl_H(X) \to \Resl_G(G \times^H X)$ given by $P \mapsto G \times^H X$ is essentially surjective and hence an equivalence.
\end{proof}
\begin{proposition}\label{Prop: Section IndQuot Equiv: Ind equiv for var}
Let $G$ be a smooth algebraic group, $X$ a left $G$-variety, and $H \leq G$ a smooth algebraic subgroup with inclusion $i:H \to G$. Then for any pseudofunctor $\overline{F}:\Var_{/K}^{\op} \to \fCat$ there is an equivalence of categories
\[
F_H(X) \simeq F_{G}\left(G \times^{H} X\right)
\]
for the induced pre-equivariant pseudofunctors on $X$ and $G \times^H X$, respectively.
\end{proposition}
\begin{proof}
This follows mutatis mutandis to the proof of Proposition \ref{Prop: Section IndQuot Equiv: Ind equiv for top} by replacing the use of Lemma \ref{Lemma: Section IndQuot Equiv: GroupResls for induction space is SubgroupResls for restriction action} with Lemma \ref{Lemma: Section IndQuot Equiv: GroupResls for Induction Variety} instead.
\end{proof}

We now show the geometric analogue of the quotient equivalence. This is slightly more involved than the induction case, as we will need 

\begin{lemma}\label{Lemma: Section IndQuot Equiv: Iso for quotients}
Let $G$ be a smooth algebraic group, $H$ a  smooth algebraic normal subgroup of $G$ with inclusion $i:H \to G$, let $X$ be a free $G$-variety, and let $P \to X$ be a smooth free $G$-resolution. Then $H \backslash P$ is a smooth free $G/H$-resolution of $H\backslash X$. In particular there is a functor $\Resl_G(X) \to \Resl_{G \backslash H}(H \backslash X)$ given by $P \mapsto H \backslash P$ on objects.
\end{lemma}
\begin{proof}
The existence of the quotients $H \backslash P$ and $H \backslash X$ follow from the fact that $G$ is a smooth free $H$-variety. Now proving that $H \backslash P$ is a smooth free $(G/H)$-space follows from the isomorphism
\[
G \backslash P \cong \left(G \backslash H\right) \backslash \left(H \backslash P\right)
\] 
and techniques used in the proof of Lemma \ref{Lemma: Section IndQuot Equiv: GroupResls for Induction Variety} by passing a trivialization $U_i \to G \backslash P$ through the isomorphism above give the desired local isotrivalization. The statement regarding functoriality is also immediate.
\end{proof}

We now need a geometric analogue of \cite[Lemma 2.1.1]{BernLun} which says that if $X$ is a free $G$-space then the only equivariant morphisms $P \to X$ are necessarily $G$-resolutions and satisfy the isomorphism
\[
P \cong X \times_{G\backslash X} G \backslash P.
\]
This isomorphism is necessary to build the other half of our equivalence of categories $\Resl_{G/H}(H \backslash X) \to \Resl_G(X)$.
\begin{lemma}[{\cite[Lemma 2.1.1]{BernLun}}]\label{Lemma: Section IndQuot Equiv: Geometric BernLun Lemma 211}
Let $G$ be a smooth algebraic group, $X$ a locally isotrivial $G$-space, and let $P \to X$ be a $G$-resolution. Then there is an isomorphism of $G$-varieties
\[
P \cong X \times_{G \backslash X} (G \backslash P)
\]
\end{lemma}
\begin{proof}[Sketch]
We geometrize the proof of \cite[Lemma 2.1.1]{BernLun}, which is only given for topological spaces. Because $P \to X$ is a $G$-resolution, by construction we have that there is unique map $\overline{p}$ making the diagram
\[
\begin{tikzcd}
	P \ar[r]{}{p} \ar[d, swap]{}{\quo_{P}} & X \ar[d]{}{\quo_X} \\
	G \backslash P \ar[r, dashed, swap]{}{\overline{p}} & G \backslash X
\end{tikzcd}
\]
commute so there is a natural morphism $\theta:X \times_{G \backslash X} G \backslash P$ of cones over the diagram $X \to G \backslash X \leftarrow G \backslash P$. Because every quotient in sight is a geometric quotient and $\quo_X:X \to G \backslash X$ is locally isotrivial, in order to prove that $P$ is isomorphic to the claimed pullback it suffices to work over a finite {\'e}tale trivialization of the quotient $\quo_{X}:X \to G \backslash X$. 

With this in mind let $\lbrace U_i \to G \backslash X \; | \; i \in I \rbrace$ be such a trivialization and fix some $U = U_i$. Define the finite {\'e}tale open $V \to G \backslash P$ via the pullback:
\[
\xymatrix{
V \ar[r] \ar[d]_{\varphi} \pullbackcorner & G \backslash P \ar[d]^{\overline{p}} \\
U \ar[r] & G \backslash X
}
\]
Since $U$ is a trivializing {\'e}tale open of $\quo_X$ by construction we have a commuting diagram:
\[
\begin{tikzcd}
\quo_{P}^{-1}(V) \ar[r] \ar[d] & G \times U \ar[d]{}{\pi_1^{GU}} \\
\quo_X^{-1}(U) \ar[ur]{}{\cong} \ar[r] & G
\end{tikzcd}
\]
Use these to define the map $\psi_V:\quo_{P}^{-1}(V) \to G$ via:
\[
\begin{tikzcd}
\quo_P^{-1}(V) \ar[r] & \quo_{X}^{-1}(U) \ar[r]{}{\cong} & G \times U \ar[r]{}{\pi_1^{GU}} & G
\end{tikzcd}
\]
Because the $\lbrace \quo_{P}^{-1}(V_i) \to P \; | \; i \in I \rbrace$ are a $G$-stable {\'e}tale cover of $P$ and the quotient $P \to G \backslash P$ is geometric, we can glue the maps $\psi_V$ to produce a morphism $\tau:P \to P$ defined locally by:
\[
\begin{tikzcd}
\quo_{P}^{-1}(V) \ar[drr, swap]{}{\tau} \ar[r]{}{\Delta} & \quo_P^{-1}(V) \times \quo_P^{-1}(V) \ar[r]{}{\psi_V \times \id} & G \times \quo_P^{-1}(V) \ar[d] \\
 & & \quo_P^{-1}(V)
\end{tikzcd}
\]
Because each of these actions are orbit-stable, it follows that $\tau$ coequalizes $\alpha_P$ and $\pi_2^{GP}$, i.e., the diagram
\[
\begin{tikzcd}
G \times P \ar[r, shift left = 1]{}{\alpha_P} \ar[r, swap, shift right = 1]{}{\pi_2^{GP}} & P \ar[r]{}{\tau} & P
\end{tikzcd}
\]
commutes. Consequently, by the universal property of the quotient, there is a unique morphism $\overline{\tau}$ making the triangle
\[
\begin{tikzcd}
P \ar[rr]{}{\tau} \ar[dr, swap]{}{\quo_P} & & P \\
 & G \backslash P \ar[ur, swap]{}{\overline{\tau}}
\end{tikzcd}
\]
commute. By intertwining $\overline{\tau}$ with the $G$-action we locally induce a morphism $\sigma_V:G \times V \to \quo_P^{-1}(V)$ which makes 
\[
\begin{tikzcd}
G \times V \ar[d, swap]{}{\sigma_V} \ar[dr]{}{} \\
\quo_{P}^{-1}(V) \ar[r] & V
\end{tikzcd}
\]
commute. However, since
\[
G \times V \cong (G \times U) \times_U V \cong \quo_{X}^{-1}(U) \times_U V 
\]
the induced diagram
\[
\begin{tikzcd}
G \times V \ar[dr] \ar[d, swap]{}{\sigma_V} \\
\quo_P^{-1}(V) \ar[r] & \quo_{X}^{-1}(U)
\end{tikzcd}
\]
also commutes. Thus the $\sigma_V$ glue to a map 
\[
\sigma:X \times_{G \backslash X} (G \backslash P) \to P
\] 
fitting into the commuting diagram
\[
\begin{tikzcd}
X \times_{G \backslash X} (G \backslash P) \ar[drr, bend left = 20]{}{} \ar[ddr, bend right = 20, swap]{}{} \ar[dr]{}[description]{\sigma} & & \\
 & P \ar[r]{}{p} \ar[d, swap]{}{\quo_P} & X \ar[d]{}{\quo_X} \\
 & G \backslash P \ar[r, swap]{}{\overline{p}} & G \backslash X
\end{tikzcd}
\]
which implies that $\theta$ and $\sigma$ are necessarily inverse to each other. Furthermore, that these isomorphisms are natural follows from the fact that each one is induced by the universal property of the pullback.
\end{proof}
\begin{lemma}\label{Lemma: Section IndQuot Equiv: The Remaining Iso}
Let $G$ be a smooth algebraic group, $H$ a smooth algebraic subgroup of $G$ with inclusion $i:H \to G$, and $X$ a locally isotrivial $G$-space. Then for any $G$-resolution $P$ of $X$, there is a natural isomorphism
\[
P \cong X \times_{H \backslash X} \left(H \backslash P\right).
\]
In particular, there is an invertible $2$-cell:
\[
\begin{tikzcd}
 & \Resl_{G/H}\left(H \backslash X\right) \ar[dr]{}{Q \mapsto X \times_{H \backslash X} Q} \\
\Resl_G(X) \ar[rr, equals] \ar[ur]{}{P \mapsto H \backslash P} & {} & \Resl_G(X) \ar[from = 1-2, to = 2-2, Rightarrow, shorten <= 4pt, shorten >= 4pt]{}{\cong}
\end{tikzcd}
\]
\end{lemma}
\begin{proof}
Because $X$ is a free $H$-variety by virtue of $G$ being a free $H$-variety, the lemma follows by an application of Lemma \ref{Lemma: Section IndQuot Equiv: Geometric BernLun Lemma 211} once we equip $X \times_{H \backslash X} Q$ with the structure of a $G$-resolution of $X$ via the first projection.
\end{proof}
\begin{proposition}\label{Prop: Section IndQuot Equiv: Quotient equiv for variety}
Let $G$ be a smooth algebraic group, $X$ a free left $G$-variety, and $H \leq G$ a smooth algebraic normal subgroup with inclusion $i:H \to G$. Then for any pseudofunctor $\overline{F}:\Var_{/K}^{\op} \to \fCat$ there is an equivalence of categories
\[
F_G(X) \simeq F_{G/H}\left(H \backslash X\right)
\]
for the induced pre-equivariant pseudofunctors on $X$ and $G \times^H X$, respectively.
\end{proposition}
\begin{proof}
By Lemma \ref{Lemma: Section IndQuot Equiv: The Remaining Iso} we have a natural isomorphism $\id_{\Resl_G(X)} \cong X \times_{H \backslash X} (H \backslash (-))$. A straightforward calculation also shows that for any $G/H$-resolution $Q \to H \backslash X$ we have that
\[
H \backslash \left(X \times_{H \backslash X} Q\right) \cong H \backslash X \times_{H \backslash X} Q \cong Q
\]
so $\id_{\Resl_{G/H}(H \backslash X)} \cong H \backslash (X \times_{H \backslash X} (-))$. Consequently it follows that we have an equivalence of categories
\[
\Resl_G(X) \simeq \Resl_{G/H}(H \backslash X).
\]
From here the proof of the proposition follows mutatis mutandis to that of Proposition \ref{Prop: Section IndQuot Equiv: Quotient equiv for top}.
\end{proof}

\newpage

\section{Averaging Functors}\label{Subsection: Averaging Functors}
In this short section we describe averaging functors. While we do not use them in any major fashion in this monograph, we expect them to be of use in future work using the pseudocone perspective on equivariant geometry and topology. We can see some evidence of this already in \cite{PramodBook} where the averaging functors are shown to exist for complex algebraic groups $G$, their closed algebraic subgroups $H$, and left $G$-varieties $X$. Our brief discussion follows the basic construction of \cite[Section 6.6]{PramodBook} as a point of departure for generalization. As a first step, however, we need to do some setup and discuss the quotient/induction equivalences. 

As we proceed, assume we are in either of the cases below:
\begin{enumerate}
	\item Let $G \leq H$ be smooth algebraic groups $G$ and $H$ with subgroup inclusion morphism $i:G \to H$;
	\item Let $G \leq H$ be topological groups with $G$ a closed subgroup and inclusion morphism $i:G \to H$.
\end{enumerate}
If $X$ is a left $H$-object, it is frequently of interest in representation theory to forget the $H$-action on $X$ and regard $X$ as a $G$-variety. In doing so, we also want to keep track of what the change of actions on $X$ does to the equivariant information by tracking the various Change of Groups functors $\varphi^{\sharp}:F_H(X) \to F_G(X)$. However, a natural question is whether or not we can return from the category $F_G(X)$ to $F_H(X)$ by enriching the $G$-equivariant information recorded by $F_G(X)$ to the $H$-equivariant information recorded by $F_H(X)$ in a way compatible with the Change of Groups functors (so in particular we should not enrich our $G$-equivariant information to $H$-equivariant infromation by fundamentally altering the $G$-equivariant information we started with). This leads to the notion of what we call left and right averaging functors\footnote{In \cite{PramodBook} these functors are simply called averaging functors and only differentiated by notation. Because the scope of \cite{PramodBook} is focused on the equivariant derived category of a complex variety, this is reasonable as both averaging functors exist; however, because of the more general $2$-categorical nature of what we are doing here it makes sense to make a more explicit distinction between left and right averaging.} (cf.\@ \cite[Definition 6.6.2]{PramodBook}).

\begin{definition}\label{Defn: Averaging Functors}
If
	\[
	\ul{F} = \left((F,\overline{F}),(E,\overline{E}),e\right)
	\] 
	is a $(G,H)$-pre-equivariant pseudofunctor we say that $\ul{F}$ has a left averaging functor $\LAvg:E_G(X) \to F_H(X)$ if $\LAvg$ is a left adjoint of $\varphi^{\sharp}$
	\[
	\begin{tikzcd}
		F_H(X) \ar[rr, bend right = 20, swap, ""{name = D}]{}{\varphi^{\sharp}} & & E_G(X) \ar[ll, bend right = 20, swap, ""{name = U}]{}{\LAvg} \ar[from = U, to = D, symbol = \dashv]
	\end{tikzcd}
	\]
	and dually we say that $\ul{F}$ has a right averaging functor $\RAvg:E_G(X) \to F_H(X)$ if there $\RAvg$ is right adjoint to $\varphi^{\ast}$:
	\[
	\begin{tikzcd}
		E_G(X) \ar[rr, bend right = 20, swap, ""{name = D}]{}{\RAvg} & & F_H(X) \ar[ll, bend right = 20, swap, ""{name = U}]{}{\varphi^{\sharp}} \ar[from = U, to = D, symbol = \dashv]	
	\end{tikzcd}
	\]
\end{definition}
\begin{example}\label{Example: Left/Right averaging}
	By \cite[Theorem 6.6.1]{PramodBook} when $H$ is a complex algebraic group and $G$ is a closed algebraic subgroup of $H$ then for the bounded derived constructible category 
	\[
	\ul{D} = \left((D^b_c(-),D_c^b(-)),(D_c^b(-),D_c^b(-)),\iota_{D_c^b(-)}\right)
	\]
	the left and right averaging functors $\LAvg$ and $\RAvg$ both exist for $\ul{D}$. Note that that $\RAvg$ is constructed from equivariant pullback and pushforward while $\LAvg$ is constructed from exceptional equivariant pullback and pushforward, this example also shows that in general the left averaging functor differs from the right averaging functors.
\end{example}
We now follow the general idea of the proof of \cite[Theorem 6.6.1]{PramodBook} to give conditions on when the averaging functors exist and some corollaries which suggest strategies for using the averaging functors. Recall that there is an $H$-equivariant morphism $i:X \to H \times^G X$ given, in point-set form, by $x \mapsto [1_H,x]_G$. Similarly, the action map $\alpha_X:H \times X \to X$ induces a morphism
\[
\overline{\alpha}_X:H \times^G X \to X
\]
which factors as
\[
\begin{tikzcd}
	H \times^{G} X \ar[rr]{}{[h,x]_G \mapsto [hG,hx]_G} \ar[dr, swap]{}{\overline{\alpha}_X} & & H/G \times X \ar[dl]{}{\pi_2} \\
	& X
\end{tikzcd}
\]
with the top row an $H$-equivariant isomorphism (cf.\@ \cite[Section 6.6]{PramodBook}). By construction these maps satisfy the commuting diagram
\[
\begin{tikzcd}
	X \ar[dr, equals]{}{} \ar[r]{}{i} & H \times^{G} X \ar[d]{}{\overline{\alpha}_X} \\
	& X
\end{tikzcd}
\]
and hence, for any pseudofunctor $\overline{F}$, induce a natural isomorphism
\[
\begin{tikzcd}
	& \overline{F}\left(H \times^G X\right) \ar[dr]{}{\overline{F}(i)} \\
	\overline{F}(X) \ar[ur]{}{\overline{F}(\overline{\alpha}_X)} \ar[rr, equals] & {} & \overline{F}(X) \ar[from = 1-2, to = 2-2, Rightarrow, shorten <= 4pt, shorten >= 4pt]{}{\phi_{i,\overline{\alpha}}} \ar[from = 1-2, to = 2-2, Rightarrow, shorten <= 4pt, shorten >= 4pt, swap]{}{\cong}
\end{tikzcd}
\]
By applying the Induction Equivalence of Propositions \ref{Prop: Section IndQuot Equiv: Ind equiv for top} and \ref{Prop: Section IndQuot Equiv: Ind equiv for var} it follows that for any $(G,H)$-pre-equivariant psuedofunctor
\[
\ul{F} = \left((F,\overline{F}),(E,\overline{E}),e\right)
\]
we can factor the Change of Groups functor as
\[
\begin{tikzcd}
	F_H(X) \ar[rr]{}{\overline{\alpha}_X^{\ast}} \ar[d, swap, ""{name = L}]{}{\varphi^{\sharp}} & & F_H\left(H \times^G X\right) \ar[d, ""{name = U}]{}{\simeq} \\
	E_G(X) & &  F_G(X) \ar[ll]{}{\PC(e \ast \quo_G^{\op})}  \ar[from = U, to = L, Rightarrow, shorten <= 4pt, shorten >= 4pt]{}{\cong}
\end{tikzcd}
\]
where $\overline{\alpha}_X^{\ast}$ denotes the pullback functor as induced by Theorem \ref{Thm: Pseudocone Functors: Pullback induced by fibre functors in pseudofucntor}. Because the functors $F_H(H \times^G X) \to F_G(X)$ and $F_G(X) \to E_G(X)$ are equivalences, assuming that we have promoted the Induction equivalence to an adjoint equivalence\footnote{Because $e$ comes from an adjoint equivalence in a bicategory $\PC(e \ast \quo_G^{\op})$ and $\PC(e^{-1} \ast \quo_G^{\op})$ are already adjoint equivalences.}, we always have adjoints
\[
\begin{tikzcd}
F_H\left(H \times^{G} X\right) \ar[rr, bend right = 20, ""{name = D}, swap]{}{P} & & F_G(X) \ar[ll, bend right = 20, ""{name = U}, swap]{}{\Lambda} \ar[from = U, to = D, symbol = \dashv]
\end{tikzcd}
\]
and
\[
\begin{tikzcd}
	F_G\left(X\right) \ar[rr, bend right = 20, ""{name = D}, swap]{}{\PC(e \ast \quo_G^{\op})} & & F_G(X) \ar[ll, bend right = 20, ""{name = U}, swap]{}{\PC(e^{-1}\ast\quo_G^{\op})} \ar[from = U, to = D, symbol = \dashv]
\end{tikzcd}
\]
where the corresponding functors are both left and right adjoint to each other. Thus in determining if $\varphi^{\sharp}:F_H(X) \to E_G(X)$ has a left or right adjoint we can reduce to asking if the functor $\overline{\alpha}_X^{\ast}:F_H(X) \to F_H(H \times^{G} X)$ admits a left or right adjoint.
\begin{proposition}\label{Prop: Averaging functors for pullbacks}
Assume we are in either of the cases below:
\begin{enumerate}
	\item Let $G \leq H$ be smooth algebraic groups $G$ and $H$ with subgroup inclusion morphism $i:G \to H$;
	\item Let $G \leq H$ be topological groups with $G$ a closed subgroup and inclusion morphism $i:G \to H$.
\end{enumerate}
If the functor $\overline{\alpha}_X:F_H(X) \to F_H(H \times^G X)$ induced by Theorem \ref{Thm: Pseudocone Functors: Pullback induced by fibre functors in pseudofucntor} has a left (respectively right) adjoint $L:F_H(H \times^G X) \to F_H(X)$ then $\varphi^{\sharp}:F_H(X) \to E_G(X)$ admits a left (respecitvely right) averaging functor.
\end{proposition}
\begin{proof}
We show only the left averaging case as the right case is dual. Assume that $\overline{\alpha}_X^{\ast}$ has a left adjoint $L$. Because the existence of adjoints is stable under isomorphism, we can replace $\varphi^{\sharp}$ with the composite:
\[
F_H(X) \xrightarrow{\overline{\alpha}_X^{\ast}} F_H\left(H \times^{G} X\right) \xrightarrow{P} F_G(X) \xrightarrow{\PC(e \ast \quo_G^{\op})} E_G(X)
\]
Writing the adjoint diagram
\[
\begin{tikzcd}
F_H(X) \ar[rr, bend right = 20, swap, ""{name = DL}]{}{\overline{\alpha}_X^{\ast}} & & F_H\left(H \times^{G}X\right) \ar[rr, bend right = 20, swap, ""{name = DR}]{}{\PC(e \ast \quo_G^{\op}) \circ P} \ar[ll, swap, bend right = 20, ""{name = UL}]{}{L} & & E_G(X) \ar[ll, swap, bend right = 20, ""{name = UR}]{}{\Lambda \circ \PC(e^{-1} \ast \quo_G^{\op})} \ar[from = UR, to = DR, symbol = \dashv] \ar[from = UL, to = DL, symbol = \dashv]
\end{tikzcd}
\]
shows that $L \circ \Lambda \circ \PC(e^{-1} \ast \quo_G^{\op}) \dashv \PC(e \ast \quo_G^{\op}) \circ P \circ \overline{\alpha}_X$ and hence that
$L \circ \Lambda \circ \PC(e^{-1} \ast \quo_G^{\op}) \dashv \varphi^{\sharp}$. Setting $\LAvg := L \circ \Lambda \circ \PC(e^{-1} \ast \quo_G^{\op})$ thus completes the proof.
\end{proof}

While the pullback construction above is very useful, Example \ref{Example: Left/Right averaging} shows that the Change of Groups functor $\varphi^{\sharp}$ can admit left and right averaging functors where at least one of the averaging functors is not induced by pullback functors given by Theorem \ref{Thm: Pseudocone Functors: Pullback induced by fibre functors in pseudofucntor}; instead we need to invoke the more general language of psuedocone translations to give rise to the technology which allows us to admit both left and right averaging functors. In essence this discussion simply proceeds mutatis mutandis to the discussion prior to the statmenet of Proposition \ref{Prop: Averaging functors for pullbacks}. If we have a pseudocone translation functor $\overline{\alpha}^{\square}:F_H(X) \to F_H(H \times^G X)$ for which we can factor the Change of Groups functor as
\[
\begin{tikzcd}
	F_H(X) \ar[rr]{}{\overline{\alpha}_X^{\square}} \ar[d, swap, ""{name = L}]{}{\varphi^{\sharp}} & & F_H\left(H \times^G X\right) \ar[d, ""{name = U}]{}{\simeq} \\
	E_G(X) & &  F_G(X) \ar[ll]{}{\PC(e \ast \quo_G^{\op})}  \ar[from = U, to = L, Rightarrow, shorten <= 4pt, shorten >= 4pt]{}{\cong}
\end{tikzcd}
\]
then we can proceed exactly as before and use adjoints to $\overline{\alpha}_X^{\square}$ to give the existence of Averaging Functors. Namely, the essential ingredients towards this technique for proving the existence of averaging functors comes down to factoring the Change of Groups functor itself by a pseudocone translation against $\overline{\alpha}_X$ followed by the Induction Equivalence and then using the existence of left/right adjoints to $\overline{\alpha}_X^{\square}$ to give left/right averaging functors.

\begin{proposition}\label{Prop: Section Chofg: General form of averaging}
Assume we are in either of the cases below:
\begin{enumerate}
	\item Let $G \leq H$ be smooth algebraic groups $G$ and $H$ with subgroup inclusion morphism $i:G \to H$;
	\item Let $G \leq H$ be topological groups with $G$ a closed subgroup and inclusion morphism $i:G \to H$.
\end{enumerate}
Let $\ul{F}$ be a $(G,H)$-pre-equivariant pseudofunctor on $X$. Assume that we have pseudocone a translation
	\[
	\begin{tikzcd}
		\SfResl_H(X)^{\op} \ar[rr, ""{name = U}]{}{\overline{F} \circ \quo_H^{\op}} \ar[dr, swap]{}{\gamma^{\op}} & & \fCat \\
		& \SfResl_H\left(H \times^G X\right)^{\op} \ar[ur, swap]{}{\overline{F}\circ\quo_{H}^{\op}} \ar[from = U, to = 2-2, Rightarrow, shorten <= 4pt, shorten >= 4pt]{}{\overline{\alpha}^{\square}}
	\end{tikzcd}
	\]
	for which there is an invertible $2$-cell:
	\[
	\begin{tikzcd}
		F_H(X) \ar[rr]{}{\overline{\alpha}_X^{\square}} \ar[d, swap, ""{name = L}]{}{\varphi^{\sharp}} & & F_H\left(H \times^G X\right) \ar[d, ""{name = U}]{}{\simeq} \\
		E_G(X) & &  F_G(X) \ar[ll]{}{\PC(e \ast \quo_G^{\op})}  \ar[from = U, to = L, Rightarrow, shorten <= 4pt, shorten >= 4pt]{}{\cong}
	\end{tikzcd}
	\]
	Then if $\overline{\alpha}_X^{\square}$ has a right adjoint (respectively a left adjoint), $\ul{F}$ has a right averaging functor (respectively a left averaging functor).
\end{proposition}
\begin{proof}
	The proof follows mutatis mutandis to the proof of Proposition \ref{Prop: Averaging functors for pullbacks}.
\end{proof}

We close this section with an example illustrating the existence of averaging functors for the $\ell$-adic equivariant derived category.
\begin{example}[{\cite[Theorem 6.6.1]{PramodBook}}]
	Let $H$ be a smooth algebraic group with a closed smooth algebraic subgroup $G$ and inclusion morphism $\varphi:G \to H$. Let $X$ be a left $H$-variety and consider the $(G,H)$-pre-equivariant pseudofunctor $\ul{D}$ induced by the constructible bounded derived pseudofunctor
	\[
	D^b_c(-; \overline{\Q}_{\ell}):\Var_{/K}^{\op} \to \fCat.
	\]
	Then because $\overline{\alpha}^{\ast} \dashv R\overline{\alpha}_{\ast}$, the right averaging functor exists by Proposition \ref{Prop: Averaging functors for pullbacks}. Similarly, using the exceptional pullbacks and pushforwards of Corollary \ref{Cor: Pseudocone Functors: Exceptional Pushforward functors for equivariant maps for scheme sheaves} also allows us to factor the Change of Groups functor as
	\[
	\begin{tikzcd}
		D^b_H(X;\overline{\Q}_{\ell}) \ar[rr]{}{\overline{\alpha}^{!}} \ar[d, swap, ""{name = L}]{}{\varphi^{\sharp}} & & D_H^b\left(H \times^G X; \overline{\Q}_{\ell}\right) \ar[d, ""{name = U}]{}{\simeq} \\
		D_G^b(X;\overline{\Q}_{\ell}) & &  D_G^b(X;\overline{\Q}_{\ell}) \ar[ll, equals]{}{}  \ar[from = U, to = L, Rightarrow, shorten <= 4pt, shorten >= 4pt]{}{\cong}
	\end{tikzcd}
	\] 
	This, together with the adjunction $R\overline{\alpha}_! \dashv \overline{\alpha}^{!}$ gives the existence of the of the left averaging functor $\LAvg$ by Proposition \ref{Prop: Section Chofg: General form of averaging}.
\end{example}

\chapter{Equivariant Six Functor Formalism and Equivariant Trace}\label{Section: Six Functor Formalism}

In this chapter we give a culmination of our categorical and geometric/topological techniques by illustrating how they can be used to study both the six functor formalisms on the categories $\DbeqQl{X}$ in what amount to essentially purely categorical terms as well as traces for equivariant derived categories. 

While the six functor formalism for $\DbeqQl{X}$ is well-known, our goal here is to illustrate how to use the pseudocone and descent-theoretic techniques can be used to make developing and using the six functor formalism relatively straightforward. For instance, we will show below how the commutativity of the cohomology functors $H_G:\DbeqQl{X} \to \Shv_G(X;\overline{\Q}_{\ell})$ and ${}^{p}H_G:\DbeqQl{X} \to \Per_G(X;\overline{\Q}_{\ell})$ with the forgetful functors $\Forget:\DbeqQl{G}{X} \to \DbQl{X}$ and $\Forget:\Shv_G(X;\overline{\Q}_{\ell}) \to \Shv(X;\overline{\Q}_{\ell})$ are essentially consequences of the pseudocone formalism developed in the first half of this monograph. We will also see that the open/closed distinguished triangles which appear in the six functor formalism for the equivariant derived category $\Dbeq{G}{X}$ arise from pseudocone techniques combined with some equivariant topology and geometry we will develop. 

Changing gears from studying the various extra properties of the six functor formalism on $\DbeqQl{X}$, we move to establish and study some basic properties of a notion of trace for the equivariant derived category $\DbeqQl{X}$. The notion of a trace associated to an object is very important in representation theory, algebraic geometry, and number theory; in fact, Deligne's first proof of the third Weil conjecture (cf.\@ \cite{DeligneWeil1} or \cite{FreitagKiehl} for an expository textbook account), following a trick of Grothendieck, used the trace of the Frobenius operators acting on the $\ell$-adic cohomology $H^{\bullet}_{\text{{\'e}t}}(V;\overline{\Q}_{\ell})$ of an $\Fbb_{q}$-variety $V$ as a crucial step. Extending these traces to objects (as is important in homological algebra, algebraic topology, and the determination of dualizable objects; cf.\@, for instance, \cite{MayTrace}) directly by associating the trace of a canonical endomorphism allows us to use these techniques in a general level. While this has been done in a stack-theoretic way for the equivariant derived category (using, implicitly, an equivalence $D_G^b(X;\overline{\Q}_{\ell}) \simeq D^b_c([G \backslash X];\overline{\Q}_{\ell})$ which is either definitional or a theorem) we show that the existence of traces for the equivariant derived category arises, in essence, as a consequence of the pseudocone construction of $\DbeqQl{X}$. Additionally we establish various functoriality properties of the trace with an eye towards future applications.

\section{A Categorical Approach towards Six Functor Formalisms}\label{Section: SFF}
We begin this chaptering by proving the existence of open/closed exact triangles in the equivariant derived category $\DbeqQl{X}$. For this, however, we will need some basic results regarding the existence these triangles at the $\Gamma$-local level and how they interact with immersions. Consequently, our main topological/geometric technique will be to show that immersions descend through the category of resolutions, i.e., that if we have an (open/closed) immersion $i:U \to V$ then the induced map $\quot{\overline{\imath}}{\Gamma}:\quot{U}{\Gamma} \to \quot{V}{\Gamma}$ is an (open/closed) immersion as well. We start this by considering closed immersions. As we proceed through this section, we will also establish the topological analogue of the geometric results as well.

\begin{lemma}\label{Lemma: Section: Yet More Stuff: Closed immersions descend}
	Let $i:V \to X$ be an equivariant closed immersion of varieties. Then for any $\Gamma \in \Sf(G)_0$ the morphism
	\[
	\overline{(\id_{\Gamma} \times i)}:\quot{V}{\Gamma} \to \XGamma
	\]
	is a closed immersion.
\end{lemma}
\begin{proof}
	Asking for the map $i$ to be equivariant is simply to ensure that the product $\Gamma \times V$ is a $G$-variety so $\quot{V}{\Gamma}$ is defined. Recall that a map $\zeta:Y \to Z$ determines a closed immersion if and only if $\lvert \zeta \rvert(\lvert Y \rvert)$ is closed in $\lvert Z  \rvert$ and there is a short exact sequence of quasi-coherent sheaves
	\[
	\begin{tikzcd}
		0 \ar[r] & \Iscr \ar[r] & \CalO_Z \ar[r]{}{\zeta^{\sharp}} &  \zeta_{\ast}\CalO_Y \ar[r] & 0
	\end{tikzcd}
	\]
	where $\Iscr$ is a sheaf of ideals on $Z$. We will first prove that we obtain the desired short exact sequence of quasi-coherent sheaves on $\XGamma$ and then show the set
	\[
	\left\lvert \overline{(\id_{\Gamma} \times i)} \right\rvert \left(\lvert \quot{V}{\Gamma} \rvert\right)
	\]
	is closed in $\lvert \XGamma \rvert$. Consider the commuting diagram
	\[
	\begin{tikzcd}
		\Gamma \times V \ar[r]{}{\id_{\Gamma} \times i} \ar[d, swap]{}{\quo_{\Gamma \times V}} & \Gamma \times X \ar[d]{}{\quo_{\Gamma \times X}} \\
		\quot{V}{\Gamma} \ar[r, swap]{}{\overline{(\id_{\Gamma} \times i)}} & \XGamma
	\end{tikzcd}
	\]
	of varieties. Because the identity map is closed and $i$ is closed, $\id_{\Gamma} \times i$ is closed and so there is a short exact sequence of quasi-coherent sheaves
	\[
	\begin{tikzcd}
		0 \ar[r] & \Kscr \ar[r] & \CalO_{\Gamma \times X} \ar[r] & (\id_{\Gamma} \times i)_{\ast}\CalO_{\Gamma \times V} \ar[r] & 0
	\end{tikzcd}
	\]
	on $\Gamma \times X$ where $\Kscr$ is a sheaf of ideals on $\Gamma \times X$. 
	
	We now apply the functor $(\quo_{\Gamma \times X})_{\ast}$ and use the fact that since $\quo_{\Gamma \times X}:\Gamma \times X \to \XGamma$ is a smooth morphism of varieties (and hence is a smooth separated quasi-compact map between Noetherian schemes), the pushforward gives an exact functor
	\[
	(\quo_{\Gamma \times X})_{\ast}:\QCoh(\Gamma \times X) \to \QCoh(\XGamma).
	\]
	As such the sequence of quasi-coherent sheaves
	\[
	\begin{tikzcd}
		0 \ar[r] & (\quo_{\Gamma \times X})_{\ast}\Iscr \ar[r] & (\quo_{\Gamma \times X})_{\ast}\CalO_{\Gamma \times X} \ar[d]\\
		& 0 &  (\quo_{\Gamma \times X})_{\ast}\left((\id_{\Gamma} \times i)_{\ast}\CalO_{\Gamma \times V}\right) \ar[l]
	\end{tikzcd}
	\]
	is short exact (and hence the sheaf $(\quo_{\Gamma \times X})_{\ast}\Iscr$ is a sheaf of ideals on $\XGamma$). Furthermore, because the quotient $\XGamma$ is a geometric quotient, the Zariski topology on $\XGamma$ coincides with the quotient topology. Thus there is an isomorphism
	\[
	(\quo_{\Gamma \times X})_{\ast}\CalO_{\Gamma \times X} \cong \CalO_{\XGamma}
	\]
	so we get that our short exact sequence above is isomorphic to the short exact sequence:
	\[
	\begin{tikzcd}
		0 \ar[r] & (\quo_{\Gamma \times X})_{\ast}\Iscr \ar[r] & \CalO_{\XGamma} \ar[d] \\ 
		& 0 & \ar[l] (\quo_{\Gamma \times X})_{\ast}\left((\id_{\Gamma} \times i)_{\ast}\CalO_{\Gamma \times V}\right) 
	\end{tikzcd}
	\]
	Now from
	\begin{align*}
		(\quo_{\Gamma \times X})_{\ast}\left((\id_{\Gamma} \times i)_{\ast}\CalO_{\Gamma \times V}\right) &\cong \left(\quo_{\Gamma \times X} \circ (\id_\Gamma \times i)\right)_{\ast}\CalO_{\Gamma \times V} \\
		& \cong \left(\overline{(\id_{\Gamma} \times i)} \circ \quo_{\Gamma \times V}\right)_{\ast}\CalO_{\Gamma \times V} \\ 
		& \cong \left(\overline{(\id_{\Gamma} \times i)}\right)_{\ast}\left((\quo_{\Gamma \times V})_{\ast}\CalO_{\Gamma \times V}\right) \\
		&\cong \left(\overline{(\id_{\Gamma} \times i)}\right)_{\ast}\CalO_{\quot{V}{\Gamma}},
	\end{align*}
	where the last isomorphism follows mutatis mutandis to the isomorphism \[
	(\quo_{\Gamma \times X})_{\ast}\CalO_{\Gamma \times X} \cong \CalO_{\XGamma},
	\] 
	we find that the short exact sequence is isomorphic to the short exact sequence
	\[
	\begin{tikzcd}
		0 \ar[r] & (\quo_{\Gamma \times X})_{\ast}\Iscr \ar[r] & \CalO_{\XGamma} \ar[r] & \left(\overline{(\id_{\Gamma} \times i)}\right)_{\ast}\CalO_{\quot{V}{\Gamma}} \ar[r] & 0
	\end{tikzcd}
	\]
	of quasi-coherent sheaves on $\XGamma$.
	
	We now argue that $\lvert \overline{(\id_{\Gamma} \times i)} \rvert \left(\lvert \quot{V}{\Gamma} \rvert\right)$ is closed in $\lvert \XGamma\rvert$. For this we recall that since the topology on $\XGamma$ is a quotient topology, for any set $W \subseteq \lvert \XGamma \rvert,$ $W$ is closed (open, respectively) if and only if $\lvert \quo_{\Gamma \times X}\rvert^{-1}(W)$ is closed (open, respectively) in $\lvert \Gamma \times X \rvert$. As such consider that
	\[
	\left\lvert \overline{(\id_{\Gamma} \times i)} \right\rvert \left(\lvert \quot{V}{\Gamma} \rvert\right) = \lbrace [\gamma, i(v)]_G \; : \; \gamma \in \lvert\Gamma\rvert, v \in \lvert V \rvert \rbrace
	\]
	and note
	\begin{align*}
		&\left\lvert \quo_{\Gamma \times X}\right\rvert^{-1}\left(\left\lvert \overline{(\id_{\Gamma} \times i)} \right\rvert \left(\lvert \quot{V}{\Gamma} \rvert\right)\right) \\
		&= \left\lbrace (\gamma, x) \in \lvert \Gamma \times X \rvert \; : \; (\gamma, x) \in \lvert \overline{(\id_{\Gamma} \times i)} \rvert \left(\lvert \quot{V}{\Gamma} \rvert\right)\right\rbrace \\
		&= \left\lbrace (\gamma, x) \; : \; \exists\, v\in \lvert V \rvert.\, [\gamma, x]_G = [\gamma,i(v)]_G \right\rbrace = \lvert \Gamma \times V \rvert \\
		&= \lvert\id_{\Gamma} \times i\rvert \left(\lvert\Gamma \times V \rvert\right).
	\end{align*}
	Because $\id_{\Gamma} \times i$ is closed, the above preimage is closed as well. This shows that $\lvert \overline{(\id_{\Gamma} \times i)}\rvert(\lvert \quot{V}{\Gamma}\rvert)$ is closed in $\lvert \XGamma\rvert$. Thus $\overline{(\id_{\Gamma} \times i)}$ is a closed immersion.
\end{proof}

Because the argument here adapts mutatis mutandis to topological spaces (by simply taking the closed component of the argument), we get the following topological analogue of Lemma \ref{Lemma: Section: Yet More Stuff: Closed immersions descend}.
\begin{lemma}\label{Lemma: Section SFF: The Inclusion}
	Let $G$ be a topological group which admits $n$-acyclic resolutions which are manifolds for all $n \in \N$ and let $X$ be a left $G$-space. If $i:V \to X$ is an equivariant closed inclusion then for any smooth $G$-resolution $\Gamma$, the inclusion
	\[
	G \backslash(\id_{\Gamma} \times i):G \backslash (\Gamma \times V) \to G \backslash (\Gamma \times X)
	\]
	is closed.
\end{lemma}

\begin{lemma}\label{Lemma: Section: Yet More Stuff: Open immerions descend}
	Let $j:U \to X$ be an equivariant open immersion of varieties and let $\Gamma \in \Sf(G)_0$. The map
	\[
	\overline{(\id_{\Gamma} \times j)}:\quot{U}{\Gamma} \to \XGamma
	\]
	is an open immersion of varieties.
\end{lemma}
\begin{proof}
	That the subset $\lvert\overline{(\id_{\Gamma} \times j)}\rvert(\lvert \quot{U}{\Gamma}\rvert)$ is open in $\lvert \XGamma \rvert$ is dual to the fact that the image $\lvert\overline{(\id_{\Gamma} \times i)}\rvert\left(\lvert \quot{V}{\Gamma}\rvert\right)$ in Lemma \ref{Lemma: Section: Yet More Stuff: Closed immersions descend} is closed and hence is omitted. By \cite[Proposition 4.2.2.a]{EGA1} we only need to show that for any point $[\gamma, u]_G\in \lvert \quot{U}{\Gamma} \rvert$ there is an isomorphism of local rings
	\[
	(\overline{(\id_{\Gamma} \times j)})^{\ast}_{[\gamma,u]_G}:\CalO_{\XGamma,[\gamma,j(u)]_G} \xrightarrow{\cong} \CalO_{\quot{U}{\Gamma}, [\gamma,u]_G}.
	\]
	For this we note that because $\id_{\Gamma} \times j$ is itself an open immersion, for each $(\gamma,u) \in \lvert \Gamma \times U\rvert$ there is an isomorphism of local rings
	\[
	(\id_{\Gamma} \times j)^{\ast}_{(\gamma,u)}:\CalO_{\Gamma \times X, (\gamma, j(u))} \xrightarrow{\cong} \CalO_{\Gamma \times U, (\gamma,u)}.
	\]
	Using the smooth surjectivity of the quotient morphisms, taking stalks, and the commuting diagram
	\[
	\begin{tikzcd}
		\Gamma \times U \ar[d, swap]{}{\quo_{\Gamma \times U}} \ar[r]{}{\id_{\Gamma} \times j} & \Gamma \times X \ar[d]{}{\quo_{\Gamma \times X}} \\
		\quot{U}{\Gamma} \ar[r, swap]{}{\overline{(\id_{\Gamma} \times j)}} & \XGamma
	\end{tikzcd}
	\]
	we find that the isomorphism $(\id_{\Gamma} \times j)^{\ast}_{(j,u)}$ descends to an isomorphism of local rings based induced by the stalks at the point $\quo_{\Gamma \times U}(\gamma,u) = [\gamma,u]_G$:
	\[
	\rho_{[\gamma,u]_G}:(\quo_{\Gamma \times X})_{\ast}\left(\CalO_{\Gamma \times X, (\gamma,j(u))}\right) \xrightarrow{\cong} (\quo_{\Gamma \times U})_{\ast}\left(\CalO_{\Gamma \times U, (\gamma, u)}\right)
	\]
	which is equal to
	\[
	\rho_{[\gamma, u]_{G}}:\left((\quo_{\Gamma \times X})_{\ast}\CalO_{\Gamma \times X}\right)_{[\gamma, j(u)]_G} \xrightarrow{\cong} \left((\quo_{\Gamma \times U})_{\ast}\CalO_{\Gamma \times U}\right)_{[\gamma,u]_G}.
	\]
	Using the isomorphisms $(\quo_{\Gamma \times X})_{\ast}\CalO_{\Gamma \times X} \cong \CalO_{\XGamma}$ and $(\quo_{\Gamma \times U})_{\ast}\CalO_{\Gamma \times U} \cong \CalO_{\quot{U}{\Gamma}}$ we then find $\rho_{[\gamma, u]_G}$ fits into a commuting diagram 
	\[
	\begin{tikzcd}
		\left((\quo_{\Gamma \times X})_{\ast}\CalO_{\Gamma \times X}\right)_{[\gamma, j(u)]_G} \ar[rr]{}{\rho_{[\gamma,u]_G}} \ar[d, swap]{}{\cong} & & \left((\quo_{\Gamma \times U})_{\ast}\CalO_{\Gamma \times U}\right)_{[\gamma,u]_G} \ar[d]{}{\cong} \\
		\CalO_{\XGamma, [\gamma, j(u)]_G} \ar[rr, swap]{}{\overline{(\id_{\Gamma} \times j)}^{\ast}_{[\gamma, j(u)]_G}} &  & \CalO_{\quot{U}{\Gamma}, [\gamma, u]_G}
	\end{tikzcd}
	\]
	where the bottom map is $\overline{(\id_{\Gamma} \times j)}_{[\gamma,u]_G}^{\ast}$ by construction of taking these particular stalks and passing through the various isomorphisms. From a diagram chase it follows that $\overline{(\id_{\Gamma} \times j)}_{[\gamma,u]_G}^{\ast}$ is an isomorphism and hence that $\overline{(\id_{\Gamma} \times j)}$ is an open immersion.
\end{proof}
\begin{lemma}\label{Lemma: Section SFF: Open inclusions descend}
	Let $j:U \to X$ be an equivariant open immersion of spaces and let $\Gamma$ be a smooth $G$-resolution. The map
	\[
	\overline{(\id_{\Gamma} \times j)}:\quot{U}{\Gamma} \to \XGamma
	\]
	is open.
\end{lemma}
\begin{proof}
	In this case the result follows immediately from the fact that quotients by group actions are open maps and the diagram
	\[
	\begin{tikzcd}
		\Gamma \times U \ar[rr]{}{\id_{\Gamma} \times j} \ar[d, swap]{}{\quo_{\Gamma \times U}} & & \Gamma \times X \ar[d]{}{\quo_{\Gamma \times X}} \\
		G \backslash(\Gamma \times U) \ar[rr, swap]{}{\overline{\id_{\Gamma} \times j}} & & G \backslash (\Gamma \times X)
	\end{tikzcd}
	\]
	commutes with each of $\id_{\Gamma} \times j,$ $\quo_{\Gamma \times U}$, and $\quo_{\Gamma \times X}$ open morphisms.
\end{proof}

\begin{lemma}\label{Lemma: Yet More Stuff: Immersions}
	Let $h:W \to X$ be an equivariant immersion of varieties and let $\Gamma \in \Sf(G)_0$. Then the map
	\[
	\overline{(\id_{\Gamma} \times h)}:\quot{W}{\Gamma} \to \XGamma
	\]
	is an immersion.
\end{lemma}
\begin{proof}
	If $h = i \circ j$ for $i$ closed and $j$ open, applying Lemmas \ref{Lemma: Section: Yet More Stuff: Open immerions descend} and \ref{Lemma: Section: Yet More Stuff: Closed immersions descend} together with the observation
	\[
	\overline{(\id_{\Gamma} \times h)} = \overline{(\id_{\Gamma} \times (i \circ j))} = \overline{(\id_{\Gamma} \times i)} \circ \overline{(\id_{\Gamma} \times j)}
	\]
	gives the result. In the other case, where $h = j \circ i$ for $j$ open and $i$ closed, we simply apply Lemmas \ref{Lemma: Section: Yet More Stuff: Closed immersions descend} and \ref{Lemma: Section: Yet More Stuff: Open immerions descend} in the opposite order of the prior case.
\end{proof}
Following the same basic structure as the proof above we get the corresponding topological version.
\begin{lemma}\label{Lemma: Section SFF: Topolgocial locally closed immersions descend}
	Let $h:W \to X$ be a locally closed inclusion of topolgoical spaces and let $\Gamma$ be a smooth $G$-resolution. Then the map
	\[
	G \backslash (\id_{\Gamma} \times h): G\backslash (\Gamma \times W) \to G \backslash (\Gamma \times X)
	\]
	is a locally closed inclusion as well provided that $W$ is equipped with a $G$-action for which $h$ is $G$-equivariant.
\end{lemma}

These various immersion lemmas allow us to deduce the existence of the open and closed immersion distinguished triangles in the equivariant derived category. However, because the existence of these triangles holds in principle for a larger class of equivariant categories than just the equivariant derived categories, we introduce some terminology to describe these.
\begin{definition}\label{Defn: Section SFF: Triangulated Preeq}\index[terminology]{Pre-equivariant Pseudofunctor! Triangulated}
Let $\ul{F} = (F,\overline{F})$ be a pre-equivariant pseudofunctor. If $\overline{F}$ is a triangulated pseudofunctor, then we say that $\ul{F}$ is a triangulated pre-equivariant pseudofunctor.
\end{definition}
\begin{definition}\label{Defn: Section SFF: Truncated Preeq}\index[terminology]{Pre-equivariant Pseudofunctor! Truncated}
	Let $\ul{F} = (F,\overline{F})$ be a pre-equivariant pseudofunctor. If $\overline{F}$ is a truncated pseudofunctor, then we say that $\ul{F}$ is a truncated pre-equivariant pseudofunctor.
\end{definition}
\begin{proposition}\label{Prop: Exact open/closed triangle general version}
Assume we are in either of the cases below:
\begin{enumerate}
	\item Let $G$ be a smooth algebraic group, $X$ a left $G$-variety, $j:U \to X$ a $G$-equivariant open immersion, and $i:V \to X$ the corresponding ($G$-equivariant) complementary closed immersion;
	\item Let $G$ be a topological group, $X$ a left $G$-space, $j:U \to X$ a $G$-equivariant open inclusion, and $i:V \to X$ the corresponding ($G$-equivariant) complementary closed inclusion.
\end{enumerate}
Assume $\ul{F} = (F,\overline{F})$ is a triangulated pre-equivariant pseudofunctor on $X$ such that:
\begin{itemize}
	\item The pullback pseudofunctors $j^{\ast}:F_G(X) \to F_G(U)$ and $i^{\ast}:F_G(X) \to F_G(V)$ induced by Theorem \ref{Thm: Pseudocone Functors: Pullback induced by fibre functors in pseudofucntor} have right adjoints
	\[
	\begin{tikzcd}
	F_G(U) \ar[rr, bend right = 20, swap, ""{name = D}]{}{j_{\ast}} & & F_G(X) \ar[ll, bend right = 20, swap, ""{name = U}]{}{j^{\ast}} \ar[from = U, to = D, symbol = \dashv]
	\end{tikzcd}
	\]
	and
	\[
		\begin{tikzcd}
		F_G(V) \ar[rr, bend right = 20, swap, ""{name = D}]{}{i_{\ast}} & & F_G(X) \ar[ll, bend right = 20, swap, ""{name = U}]{}{i^{\ast}} \ar[from = U, to = D, symbol = \dashv]
	\end{tikzcd}
	\]
	induced by Theorem \ref{Thm: Pseudocone Functors: Adjoints between translations give adjoint functors}.
	\item There are pseudocone translations and optranslations which give rise to functors
	\[
	j^{\square}:F_G(X) \to F_G(U)\qquad i^{\square}:F_G(X) \to F_G(V)
	\]
	and
	\[
	j_{\square}:F_G(U) \to F_G(X)\qquad i_{\square}:F_G(V) \to F_G(X)
	\]
	which have adjoints
		\[
	\begin{tikzcd}
		F_G(X) \ar[rr, bend right = 20, swap, ""{name = D}]{}{j^{\square}} & & F_G(U) \ar[ll, bend right = 20, swap, ""{name = U}]{}{j_{\square}} \ar[from = U, to = D, symbol = \dashv]
	\end{tikzcd}
	\]
	and:
	\[
	\begin{tikzcd}
		F_G(X) \ar[rr, bend right = 20, swap, ""{name = D}]{}{i^{\square}} & & F_G(V) \ar[ll, bend right = 20, swap, ""{name = U}]{}{i_{\square}} \ar[from = U, to = D, symbol = \dashv]
	\end{tikzcd}
	\]
	\item For any free $G$-resolution $\Gamma$ there is a natural open/closed distinguished triangle
	\[
	\begin{tikzcd}
	(\quot{j}{\Gamma}_{\square} \circ \quot{j}{\Gamma}^{\square})(-) \ar[r]{}{\quot{\epsilon^{j}}{\Gamma}} & \id_{\overline{F}(\XGamma)} \ar[r]{}{\quot{\eta^{i}}{\Gamma}} & (\quot{i}{\Gamma}_{\ast} \circ \quot{i}{\Gamma}^{\ast})(-) \ar[r] & \left(\left(\quot{j}{\Gamma}_{\square} \circ \quot{j}{\Gamma}^{\square}\right)[1]\right)(-)
	\end{tikzcd}
	\]
	in $\overline{F}(\XGamma)$ where $\quot{\epsilon}{\Gamma}^{j}:\quot{j}{\Gamma}_{\square} \circ \quot{j}{\Gamma}^{\square} \Rightarrow \id_{\overline{F}(\XGamma)}$ and $\quot{\eta}{\Gamma}^i:\id_{\overline{F}(\XGamma)} \Rightarrow \quot{i}{\Gamma}_{\ast} \circ \quot{i}{\Gamma}^{\ast}$ are the counits and units of adjunction, respectively.
\end{itemize}
Then for every object $A$ in $F_G(X)$ there is an open/closed distinguished triangle:
\[
\begin{tikzcd}
(j_{\square} \circ j^{\square})(A) \ar[r]{}{\epsilon^j_A} & A \ar[r]{}{\eta_A^i} & (i_{\ast} \circ i^{\ast})(A) \ar[r] & \left(\left(j_{\square} \circ j^{\square}\right)(A)\right)[1]
\end{tikzcd}
\]
\end{proposition}
\begin{proof}
The assumptions in the statement of the proposition imply by Theorem \ref{Theorem: Section Triangle: Equivariant triangulation} that if we can show that $\quot{U}{\Gamma}$ and $\quot{V}{\Gamma}$ are complementary open and closed subobjects of $\XGamma$ for all $\Gamma$ then the corresponding triangle in $F_G(X)$ is distinguished. However, this is immediate by the fact that $i$ and $j$ are $G$-equivariant: fix a term $[\gamma,x]_G$ in $\quot{X}{\Gamma}$. Note that because $i$ and $j$ are equivariant, $[\gamma,x]_G$ is in the space $\quot{V}{\Gamma} \cap \quot{U}{\Gamma}$ if and only if $x \in U \cap V$, which is impossible because $U \cap V = \emptyset$.
\end{proof}

We now can use the proposition above to prove the existence of open/closed distinguished triangles for the equivariant derived categories on varieties and on spaces.
\begin{corollary}\label{Cor: Exact open/closed triangle}
	Let $j:U \to X$ and be a $G$-equivariant open immersion and let $i:V \to X$ be the corresponding closed immersion. Then for any object $A \in \DbeqQl{X}_0$ there is a distinguished triangle of the form
	\[
	\begin{tikzcd}
		Rj_!\left(j^{!}A\right) \ar[r]{}{\epsilon_A^j} & A \ar[r]{}{\eta_A^i} & Ri_{\ast}\left(i^{\ast}A\right) \ar[r] & Rj_!\left(j^!A\right)[1]
	\end{tikzcd}
	\]
	where $\epsilon^j:Rj_{!} \circ j^{!} \Rightarrow \id_{\DbeqQl{X}}$ is the counit of adjunction and $\eta^i:\id_{\DbeqQl{X}} \Rightarrow Ri_{\ast} \circ i^{\ast}$ is the unit of adjunction.
\end{corollary}
\begin{proof}
Because the derived categories $D_c^b(Z)$ and $D_c^b(Z;\overline{\Q}_{\ell})$ have open/closed distinguished triangles with respect to the pushforward/pullback and exceptional pullback/exceptional pushforward adjunctions, the result follows by applying Proposition \ref{Prop: Exact open/closed triangle general version}.
\end{proof}
\begin{corollary}\label{Cor: Section SFF: Openclosed Dist Triangles}
	Let $G$ be a topological group for which $G$ admits $n$-acylic resolutions which are manifolds for all $n \in \N$ and let $X$ be a left $G$-space. Let $j:U \to X$ be a $G$-equivariant open inclusion and let $i:V \to X$ be the corresponding ($G$-equivariant) closed inclusion. Then for any object $A$ of $D^b_G(X)$ there are distinguished triangles of the form
	\[
	\begin{tikzcd}
		Rj_!\left(j^{!}A\right) \ar[r]{}{\epsilon^j_A} & A \ar[r]{}{\eta^i_A} & Ri_{\ast}\left(i^{\ast}A\right) \ar[r] & Rj_!\left(j^!A\right)[1]
	\end{tikzcd}
	\]
	and
	\[
	\begin{tikzcd}
		Ri_!\left(i^!A\right) \ar[r]{}{\epsilon^i_A} & A \ar[r]{}{\eta^j_A} & Rj_{\ast}\left(j^{\ast} A\right) \ar[r] & \left(Ri_!\left(i^!A\right)\right)[1]
	\end{tikzcd}
	\]
	where $\epsilon^j:Rj_{!} \circ j^{!} \Rightarrow \id_{D_G^b(X)}$ and $\epsilon^i:Ri_{!} \circ i^{!} \Rightarrow \id_{D_G^b(X)}$ are the counits of adjunction while $\eta^i:\id_{D_G^b({X})} \Rightarrow Ri_{\ast} \circ i^{\ast}$ and $\eta^j:\id_{D^b_G(X)} \Rightarrow Rj_{\ast} \circ j^{\ast}$ are the units of adjunction.
\end{corollary}
\begin{proof}
	The basic structure of the argument follows mutatis mutandis to the proof of Corollary \ref{Cor: Exact open/closed triangle} in the case of the first (family of) distinguished triangle(s). From here the existence of the second (family of) distinguished triangle(s) follows from the fact that there are distinguished triangles
	\[
	\begin{tikzcd}
	\left(R(\quot{\overline{\imath}}{M})_{!} \circ \quot{\overline{\imath}}{M}^{!}\right)(-) \ar[r]{}{\quot{\epsilon^{i}}{M}} & \id_{D^b(\quot{X}{M})}(-) \ar[r]{}{\quot{\eta^{j}}{M}} & \left(R(\quot{\overline{\jmath}}{M})_{\ast} \circ \quot{\overline{\jmath}}{M}^{\ast}\right)(-) \ar[d] \\
	 & & \left(\left(R(\quot{\overline{\imath}}{M})_{!} \circ \quot{\overline{\imath}}{M}^{!}\right)(-)\right)[1]
	\end{tikzcd}
	\]
	for each $M \in \mathbf{FAMan}(G)_0$ and then proceeding as in the proof of Proposition \ref{Prop: Exact open/closed triangle general version}.
\end{proof}

We now prove that the forgetful functor of Proposition \ref{Prop: Section ECSV: Forgetful functor} (and also Proposition \ref{Prop: Section 3.3: Forgetful functor is de-equivariantification functor}), applied to the equivariant derived category pre-equvariant pseudofunctors and the standard and perverse cohomology functors $H^0_G:\DbeqQl{X} \to \Shv_G(X;\overline{\Q}_{\ell})$, ${}^{p}H_G^0:\DbeqQl{X} \to \Per_G(X;\overline{\Q}_{\ell})$, commute up to isomorphism for smooth algebraic groups $G$ and left $G$-varieties $X$ and in both the variety-theoretic and topological cases for the equivariant derived categories.

For what follows we will need to recall a technical result from algebraic geometry. Additionally, we will also need to define a special notation to denote the conventional perspective on equivariant sheaves and perverse sheaves which is distinct from the pseudocone language we have been using so far. This leads us to the following lemma and definition of convenience.
\begin{lemma}[{\cite[\href{https://stacks.math.columbia.edu/tag/0F7L}{Tag 0F7L}]{stacks-project}}]\label{LEmma: Proper base change commutes iwth pullback}
	Consider a pullback diagram of quasi-compact and quasi-separated schemes
	\[
	\begin{tikzcd}
		W \ar[r]{}{g} \ar[d, swap]{}{k} & Z \ar[d]{}{f} \\
		Y \ar[r, swap]{}{h} & X 
	\end{tikzcd}
	\]
	where $f$ is a separated morphism of finite type. Then there is a natural isomorphism
	\[
	h^{\ast} \circ Rf_! \cong Rk_! \circ g^{\ast}.
	\]
\end{lemma}
\begin{definition}
	Let $G$ be either a smooth algebraic group or a topological space which admits $n$-acylcic resolutions which are manifolds for all $n \in \N$ and let $X$ be a left $G$-variety or a left $G$-space, respectively. We then write the following when $G$ is an algebraic group and $X$ is a left $G$-variety\index[notation]{ShvGnaive@$\Shv_G(X)^{\textrm{na{\"i}ve}}$}\index[notation]{ShvGnaiveQl@$\Shv_G^{\textrm{na{\"i}ve}}(X;\overline{\Q}_{\ell})$}\index[notation]{PerGnaive@$\Per_G^{\textrm{na{\"i}ve}}(X)$}\index[notation]{PerGnaiveQl@$\Per_G^{\textrm{na{\"i}ve}}(X;\overline{\Q}_{\ell})$}
	\begin{align*}
		\Shv_G(X) &:= \PC\left(\Shv\left(G \backslash (-)\right):\SfResl_G(X)^{\op} \to \fCat\right)	\\
		\Shv_G^{\textrm{na{\"i}ve}}(X) &:= \text{traditional definition of equivariant sheaves} \\
		\Shv_G(X; \overline{\Q}_{\ell}) &:= \PC\left(\Shv\left(G \backslash (-); \overline{\Q}_{\ell}\right):\SfResl_G(X)^{\op} \to \fCat\right) \\
		\Shv_G^{\textrm{na{\"i}ve}}(X;\overline{\Q}_{\ell}) &:= \text{traditional definition of equivariant}\,\ell\text{-adic sheaves} \\
		\Per_G(X) &:= \PC\left(\Per\left(G \backslash (-)\right):\SfResl_G(X)^{\op} \to \fCat\right)	\\
		\Per_G^{\textrm{na{\"i}ve}}(X) &:= \text{traditional definition of equivariant perverse sheaves} \\
		\Per_G(X; \overline{\Q}_{\ell}) &:= \PC\left(\Per\left(G \backslash (-); \overline{\Q}_{\ell}\right):\SfResl_G(X)^{\op} \to \fCat\right) \\
		\Per_G^{\textrm{na{\"i}ve}}(X;\overline{\Q}_{\ell}) &:= \text{traditional definition of equivariant}\,\ell\text{-adic perverse sheaves}
	\end{align*}
	and
	\begin{align*}
		\Shv_G(X) &:= \PC\left(\Shv\left(G \backslash (-)\right):\SfResl_G(X)^{\op} \to \fCat\right)	\\
		\Shv_G^{\textrm{na{\"i}ve}}(X) &:= \text{traditional definition of equivariant sheaves} \\
		\Per_G(X) &:= \PC\left(\Per\left(G \backslash (-)\right):\SfResl_G(X)^{\op} \to \fCat\right)	\\
		\Per_G^{\textrm{na{\"i}ve}}(X) &:= \text{traditional definition of equivariant perverse sheaves} 
	\end{align*}
	when $G$ is a topological group and $X$ is a left $G$-space.
\end{definition}

Observe that because of the equivalences of Corollaries \ref{Cor: Section Triangle: Pseudofunctor of standard hearts is heart of pseudofunctor}  and Corollary \ref{Cor: Section Triangle: Psedofunctor of perv is equiv perv} we have invertible $2$-cells
\[
\begin{tikzcd}
	\DbeqQl{X} \ar[rr, bend left = 30, ""{name = Up}]{}{\Forget} \ar[d, swap]{}{H_G^0} \ar[r, ""{name = LU}]{}{\DbQl{1_G^{\sharp}}} & D^b_{\Spec K}(X;\overline{\Q}_{\ell}) \ar[d]{}{H_{\Spec K}^{0}} \ar[r, ""{name = RU}]{}{\overline{D}^b} & \DbQl{X} \ar[d]{}{H^0} \\
	\Shv_G(X;\overline{\Q}_{\ell}) \ar[r, swap, ""{name = LL}]{}{\Shv(1_G^{\sharp})} \ar[rr, swap, bend right = 30, ""{name = Down}]{}{\Forget} & \Shv_{\Spec K}(X;\overline{\Q}_{\ell}) \ar[r, swap, ""{name = RL}]{}{\overline{\Shv}} & \Shv(X;\overline{\Q}_{\ell}) \ar[from = LU, to = LL, Rightarrow, shorten <= 4pt, shorten >= 4pt]{}{\cong} \ar[from = RU, to = RL, Rightarrow, shorten <= 4pt, shorten >= 4pt]{}{\cong} \ar[from = Up, to = 1-2, Rightarrow, shorten >= 4pt, shorten <= 4pt]{}{\cong} \ar[from = 2-2, to = Down, Rightarrow, shorten <= 4pt, shorten >= 4pt]{}{\cong}
\end{tikzcd}
\]
where $\DbQl{1_G^{\sharp}}$ is the Change of Groups functor along the identity inclusion $1_G:\Spec K \to G$ against the pre-equivariant pseudofunctor $\DbQl{-}$ (and similarly for $\Shv(1_G^{\sharp})$; cf.\@ Theorem \ref{Thm: Section 3: Change of groups functor}) and where $\overline{D}^b$ and $\overline{\Shv}$ are the equivalence witnesses for
\[
D^b_{\Spec K}(X;\overline{\Q}_{\ell}) \simeq \DbQl{X}
\]
and
\[
\Shv_{\Spec K}(X;\overline{\Q}_{\ell}) \simeq \Shv(X;\overline{\Q}_{\ell})
\]
furnished by Proposition \ref{Prop: Pseudocone Section: Terminal object in C gives us pseudocones as global sections}. Note that the inner squares form invertible $2$-cells via a routine use of Proposition \ref{Prop: Section Chofg: Change of groups together with pullbacks}. Pasting all the cells together we get the invertible $2$-cell:
\[
\begin{tikzcd}
	\DbeqQl{X} \ar[rr, ""{name = Up}]{}{\Forget} \ar[d, swap]{}{H_G^0} & & \DbQl{X} \ar[d]{}{H^0} \\
	\Shv_G(X;\overline{\Q}_{\ell}) \ar[rr, swap, ""{name = Down}]{}{\Forget} & & \Shv(X;\overline{\Q}_{\ell}) \ar[from = Up, to = Down, Rightarrow, shorten >= 4pt, shorten <= 4pt]{}{\cong}
\end{tikzcd}
\]
Adding the invertible $2$-cell induced via the equivalence of categories 
\[
\Shv_G(X;\overline{\Q}_{\ell}) \simeq \Shv_G^{\textrm{na{\"i}ve}}(X;\overline{\Q}_{\ell})
\] 
of Corollary \ref{Cor: Section Triangle: Pseudofunctor of standard hearts is heart of pseudofunctor} fitting into the pasting diagram
\[
\begin{tikzcd}
	\DbeqQl{X} \ar[rr, ""{name = Up}]{}{\Forget} \ar[d, swap]{}{H_G^0} & & \DbQl{X} \ar[d]{}{H^0} \\
	\Shv_G(X;\overline{\Q}_{\ell}) \ar[rr, swap, ""{name = Down}]{}{\Forget} \ar[dr, swap]{}{\simeq} & & \Shv(X;\overline{\Q}_{\ell}) \ar[from = Up, to = Down, Rightarrow, shorten >= 4pt, shorten <= 4pt]{}{\cong} \\
	& \Shv_G^{\textrm{na{\"i}ve}}(X;\overline{\Q}_{\ell}) \ar[ur, swap]{}{\Forget} \ar[from = Down, to = 3-2, Rightarrow, shorten <= 6pt, shorten >= 4pt]{}{\cong}
\end{tikzcd}
\]
produces the total diagram:
\[
\begin{tikzcd}
	\DbeqQl{X} \ar[rr, ""{name = Up}]{}{\Forget} \ar[d, swap]{}{H_G^0} & & \DbQl{X} \ar[d]{}{H^0} \\
	\Shv_G^{\textrm{na{\"i}ve}}(X;\overline{\Q}_{\ell}) \ar[rr, swap, ""{name = Down}]{}{\Forget} & & \Shv(X;\overline{\Q}_{\ell}) \ar[from = Up, to = Down, Rightarrow, shorten >= 4pt, shorten <= 4pt]{}{\cong}
\end{tikzcd}
\]
Note that the existence of the invertible $2$-cell
\[
\begin{tikzcd}
	\DbeqQl{X} \ar[rr, ""{name = Up}]{}{\Forget} \ar[d, swap]{}{{}^{p}H_G^0} & & \DbQl{X} \ar[d]{}{{}^{p}H^0} \\
	\Per_G^{\textrm{na{\"i}ve}}(X;\overline{\Q}_{\ell}) \ar[rr, swap, ""{name = Down}]{}{\Forget} & & \Per(X;\overline{\Q}_{\ell}) \ar[from = Up, to = Down, Rightarrow, shorten >= 4pt, shorten <= 4pt]{}{\cong}
\end{tikzcd}
\]
follows mutatis mutandis to the diagram above, save from using Corollary \ref{Cor: Section Triangle: Psedofunctor of perv is equiv perv} instead of  Corollary \ref{Cor: Section Triangle: Pseudofunctor of standard hearts is heart of pseudofunctor}. Putting these together gives the proposition below.

\begin{proposition}\label{Prop: Forget restricts to naive forget}
	For any smooth algebraic group $G$ and any left $G$-variety $X$ there are invertible $2$-cells
	\[
	\begin{tikzcd}
		\DbeqQl{X} \ar[rr, ""{name = Up}]{}{\Forget} \ar[d, swap]{}{H_G^0} & & \DbQl{X} \ar[d]{}{H^0} \\
		\Shv_G^{\operatorname{na{\ddot{\imath}}ve}}(X;\overline{\Q}_{\ell}) \ar[rr, swap, ""{name = Down}]{}{\Forget} & & \Shv(X;\overline{\Q}_{\ell}) \ar[from = Up, to = Down, Rightarrow, shorten >= 4pt, shorten <= 4pt]{}{\cong}
	\end{tikzcd}
	\]
	and:
	\[
	\begin{tikzcd}
		\DbeqQl{X} \ar[rr, ""{name = Up}]{}{\Forget} \ar[d, swap]{}{{}^{p}H_G^0} & & \DbQl{X} \ar[d]{}{{}^{p}H^0} \\
		\Per_G^{\operatorname{na{\ddot{\imath}}ve}}(X;\overline{\Q}_{\ell}) \ar[rr, swap, ""{name = Down}]{}{\Forget} & & \Per(X;\overline{\Q}_{\ell}) \ar[from = Up, to = Down, Rightarrow, shorten >= 4pt, shorten <= 4pt]{}{\cong}
	\end{tikzcd}
	\]
	That is, the forgetful functors on the equivariant derived category restrict to the forgetful functors on equivariant perverse sheaves and equivarant $\ell$-adic sheaves. Moreover, we also have invertible $2$-cells
	\[
	\begin{tikzcd}
		D_G^b({X}) \ar[rr, ""{name = Up}]{}{\Forget} \ar[d, swap]{}{H_G^0} & & D_c^b(X) \ar[d]{}{H^0} \\
		\Shv_G^{\operatorname{na\ddot{\imath}ve}}(X) \ar[rr, swap, ""{name = Down}]{}{\Forget} & & \Shv(X) \ar[from = Up, to = Down, Rightarrow, shorten >= 4pt, shorten <= 4pt]{}{\cong}
	\end{tikzcd}
	\]
	and:
	\[
	\begin{tikzcd}
		D_G^b({X}) \ar[rr, ""{name = Up}]{}{\Forget} \ar[d, swap]{}{{}^{p}H_G^0} & & D_c^b({X}) \ar[d]{}{{}^{p}H^0} \\
		\Per_G^{\operatorname{na{\ddot{\imath}}ve}}(X) \ar[rr, swap, ""{name = Down}]{}{\Forget} & & \Per(X) \ar[from = Up, to = Down, Rightarrow, shorten >= 4pt, shorten <= 4pt]{}{\cong}
	\end{tikzcd}
	\]
	whenever $G$ is a smooth algebraic group and $X$ is a left $G$-variety or whenever $G$ is a topological group which admits free $n$-acyclic $G$-spaces which are manifolds for all $n \in \N$ and $X$ is a left $G$-space.
\end{proposition}

We now show that the each of the $\ell$-adic equivariant six functors $h^{\ast},$ $Rh_{\ast},$ $Rh_{!},$ $h^{!},$ $(-)\otimes(-),$ and $R[-,-]$ commute up to isomorphism with the forgetful functors when $G$ is a smooth algebraic group and $X$ is a left $G$-variety. Because the proof we give below is essentially entirely formal in nature and mainly organizes various coherence results and commutativity witnesses together, we will also be able to adapt it to the topological and non-$\ell$-adic cases more-or-less freely.

\begin{proposition}\label{Prop: Six functor forgetful commute}
	Let $h:X \to Y$ be a $G$-equivariant morphism for a smooth algebraic group $G$ and left $G$-varieties $X$ and $Y$. Then each of the six functors of Corollaries \ref{Cor: Pseudocone Functors: Monoidal Closed Derived Cat}, \ref{Cor: Pseudocone Functors: Existence of equivariant pullback  for schemes}, \ref{Cor: Pseudocone Functors: Pushforward functors for equivariant maps for scheme sheaves}, and \ref{Cor: Pseudocone Functors: Exceptional Pushforward functors for equivariant maps for scheme sheaves}  commute with the forgetful functors up to isomorphism for both the $\ell$-adic equivariant derived category and standard equivariant derived category.
\end{proposition}
\begin{proof}
	We begin by writing the non-equivariant isomorphisms $\psi_X:X \to \quot{X}{G}$ and $\psi_{Y}:Y\to \quot{Y}{G}$ of Proposition \ref{Prop: Section ECSV: Forgetful functor}. Note that these are natural in the sense that the diagram of varieties
	\[
	\begin{tikzcd}
		X \ar[r]{}{h} \ar[d, swap]{}{\psi_{X}} & Y \ar[d]{}{\psi_Y} \\
		\quot{X}{G} \ar[r, swap]{}{\quot{\overline{h}}{G}} & \quot{Y}{G}
	\end{tikzcd}
	\]
	commutes. We also write $\Forget_X:\DbeqQl{X} \to \DbQl{X}$ for the corresponding forgetful functors of Proposition \ref{Prop: Section ECSV: Forgetful functor} associated to the equivariant derived category$D_G^b(X;\overline{\Q}_{\ell})$ and similarly for $\Forget_Y:\DbeqQl{Y} \to \DbQl{Y}$.
	
	Let us now show that the natural isomorphism $\Forget_X \circ h^{\ast} \cong h^{\ast} \circ \Forget_Y$
	exists. We calculate that on one hand, for any $A \in \DbeqQl{Y}_0$,
	\begin{align*}
		\Forget_X\left(h^{\ast}A\right) &= \Forget_X\left(\left\lbrace \hGamma^{\ast}\left(\AGamma\right) \; | \; \Gamma \in \Sf(G)_0, \AGamma \in A \right\rbrace\right) \\
		&= \psi_X^{\ast}\left(\quot{\overline{h}}{G}^{\ast}\left(\quot{A}{G}\right)\right) = \left(\psi_X^{\ast} \circ \quot{\overline{h}}{G}^{\ast}\right)\quot{A}{G}.
	\end{align*}
	On the other hand
	\[
	h^{\ast}\left(\Forget_Y(A)\right) = h^{\ast}\left(\psi_Y^{\ast}\left(\quot{A}{G}\right)\right) = \left(h^{\ast} \circ \psi_Y^{\ast}\right)\left(\quot{A}{G}\right);
	\]
	using the natural isomorphism of functors $h^{\ast} \circ \psi_Y^{\ast} \cong \psi_X^{\ast} \circ \quot{\overline{h}}{G}^{\ast}$ thus gives the desired invertible $2$-cell:
	\[
	\begin{tikzcd}
		\DbeqQl{Y} \ar[r, ""{name = Up}]{}{h^{\ast}} \ar[d, swap]{}{\Forget_Y} & \DbeqQl{X} \ar[d]{}{\Forget_X} \\
		\DbQl{Y} \ar[r, swap, ""{name = Down}]{}{h^{\ast}} & \DbQl{X} \ar[from = Up, to = Down, Rightarrow, shorten >= 4pt, shorten <= 4pt]{}{\cong}
	\end{tikzcd}
	\]
	
	We now show the existence of the natural isomorphism
	\[
	\Forget_X \circ h^{!} \cong h^{!} \circ \Forget_Y.
	\]
	Because the diagram
	\[
	\begin{tikzcd}
		X \ar[r]{}{h} \ar[d, swap]{}{\psi_{X}} & Y \ar[d]{}{\psi_Y} \\
		\quot{X}{G} \ar[r, swap]{}{\quot{\overline{h}}{G}} & \quot{Y}{G}
	\end{tikzcd}
	\]
	is a pullback diagram with both morphisms $\psi_X$ and $\psi_Y$ smooth (by virtue of being isomorphisms --- isomorphisms are {\'e}tale and hence smooth of relative dimension zero) we have a natural isomorphism of functors $h^{!}\circ\psi_Y^{\ast}  \cong \psi_X^{\ast} \circ \quot{\overline{h}}{G}^{!}$. Using this we calculate for any $A \in \DbeqQl{Y}_0$
	\begin{align*}
		\Forget_X\left(h^{!}A\right) &= \Forget_X\left(\left\lbrace \hGamma^{!}\left(\AGamma\right) \; | \; \Gamma \in \Sf(G)_0, \AGamma \in A \right\rbrace\right) \\
		&= \psi_X^{\ast}\left(\quot{\overline{h}}{G}^{!}\left(\quot{A}{G}\right)\right) = \left(\psi_X^{\ast} \circ \quot{\overline{h}}{G}^{!}\right)\quot{A}{G} \\
		&\cong \left(h^{!}\circ\psi_Y^{\ast}\right)\quot{A}{G} = h^{!}\left(\psi_Y^{\ast}(\quot{A}{G})\right) \\
		&= h^{!}\left(\Forget_Y\left\lbrace \AGamma \; | \; \Gamma \in \Sf(G)_0\right\rbrace\right) = h^{!}\left(\Forget_Y(A)\right)
	\end{align*}
	This establishes the invertible $2$-cell:
	\[
	\begin{tikzcd}
		\DbeqQl{Y} \ar[r, ""{name = Up}]{}{h^{!}} \ar[d, swap]{}{\Forget_Y} & \DbeqQl{X} \ar[d]{}{\Forget_X} \\
		\DbQl{Y} \ar[r, swap, ""{name = Down}]{}{h^{!}} & \DbQl{X} \ar[from = Up, to = Down, Rightarrow, shorten >= 4pt, shorten <= 4pt]{}{\cong}
	\end{tikzcd}
	\]
	
	We now prove the existence of the natural isomorphism
	\[
	Rh_{\ast} \circ \Forget_X \cong \Forget_Y \circ Rh_{\ast}.
	\]
	As in the case above, we use that the diagram
	\[
	\begin{tikzcd}
		X \ar[r]{}{h} \ar[d, swap]{}{\psi_{X}} & Y \ar[d]{}{\psi_Y} \\
		\quot{X}{G} \ar[r, swap]{}{\quot{\overline{h}}{G}} & \quot{Y}{G}
	\end{tikzcd}
	\]
	is a pullback diagram of varieties with both $\psi_X$ and $\psi_Y$ smooth morphisms. Thus we get isomorphisms of functors $Rh_{\ast} \circ \psi_X^{\ast} \cong \psi_Y^{\ast} \circ R(\quot{\overline{h}}{G})_{\ast}$ by the Smooth Base Change Theorem. Proceeding as in the two prior cases we derive that for any $A \in \DbeqQl{X}_0$,
	\begin{align*}
		(Rh_{\ast} \circ \Forget_X)A &= \left(Rh_{\ast} \circ \psi_X^{\ast}\right)\quot{A}{G} \cong \left(\psi_Y^{\ast} \circ R(\quot{\overline{h}}{G})_{\ast}\right)\quot{A}{G} \\
		&= (\Forget_Y \circ Rh_{\ast})A
	\end{align*}
	which gives the invertible $2$-cell:
	\[
	\begin{tikzcd}
		\DbeqQl{X} \ar[r, ""{name = Up}]{}{Rh_{\ast}} \ar[d, swap]{}{\Forget_Y} & \DbeqQl{Y} \ar[d]{}{\Forget_X} \\
		\DbQl{X} \ar[r, swap, ""{name = Down}]{}{Rh_{\ast}} & \DbQl{Y} \ar[from = Up, to = Down, Rightarrow, shorten >= 4pt, shorten <= 4pt]{}{\cong}
	\end{tikzcd}
	\]
	In fact, using that the diagram
	\[
	\begin{tikzcd}
		X \ar[r]{}{h} \ar[d, swap]{}{\psi_{X}} & Y \ar[d]{}{\psi_Y} \\
		\quot{X}{G} \ar[r, swap]{}{\quot{\overline{h}}{G}} & \quot{Y}{G}
	\end{tikzcd}
	\]
	satisfies the hypotheses of Lemma \ref{LEmma: Proper base change commutes iwth pullback} (because every morphism and scheme in sight is finite type and separated) we get that $ Rh_{!} \circ \psi_X^{\ast}  \cong \psi_Y^{\ast} \circ R(\quot{\overline{h}}{G})_{!}$. Proceeding mutatis mutandis to the pusforward case gives the invertible $2$-cell:
	\[
	\begin{tikzcd}
		\DbeqQl{X} \ar[r, ""{name = Up}]{}{Rh_{!}} \ar[d, swap]{}{\Forget_Y} & \DbeqQl{Y} \ar[d]{}{\Forget_X} \\
		\DbQl{X} \ar[r, swap, ""{name = Down}]{}{Rh_{!}} & \DbQl{Y} \ar[from = Up, to = Down, Rightarrow, shorten >= 4pt, shorten <= 4pt]{}{\cong}
	\end{tikzcd}
	\]
	
	We now show that the tensor products commute with the forgetful functor, i.e., there is a natural isomorphism
	\[
	\Forget_X\left(A \os{L}{\otimes} B\right) \cong \Forget_X(A) \os{L}{\otimes} \Forget_X(B)
	\]
	for any $A,B \in \DbeqQl{X}$. For this we calculate that on one hand
	\begin{align*}
		\Forget_X\left(A \os{L}{\otimes} B\right) &= \psi_X^{\ast}\left(\quot{A}{G} \ds{L}{\XGamma}{\otimes} \quot{B}{G}\right)
	\end{align*}
	while on the other hand
	\[
	\Forget_X(A) \os{L}{\otimes} \Forget_X(B) = \psi_X^{\ast}\left(\quot{A}{G}\right) \os{L}{\otimes} \psi_X^{\ast}\left(\quot{B}{G}\right).
	\]
	Because $\psi:X \to \quot{X}{G}$ is smooth, there is a a natural isomorphism
	\[
	\psi_X^{\ast}\left(\quot{A}{G} \ds{L}{\XGamma}{\otimes} \quot{B}{G}\right) \cong \psi_X^{\ast}\left(\quot{A}{G}\right) \os{L}{\otimes} \psi_X^{\ast}\left(\quot{B}{G}\right)
	\]
	which gives an invertible $2$-cell:
	\[
	\begin{tikzcd}
		\DbeqQl{X} \times \DbeqQl{X} \ar[rr, ""{name = Up}]{}{\os{L}{\otimes}} \ar[d, swap]{}{\Forget_X \times \Forget_X} & & \DbeqQl{X} \ar[d]{}{\Forget_X} \\
		\DbQl{X} \times \DbQl{X} \ar[rr, swap, ""{name = Down}]{}{\os{L}{\otimes}} & & \DbQl{X} \ar[from = Up, to = Down, Rightarrow, shorten >= 4pt, shorten <= 4pt]{}{\cong}
	\end{tikzcd}
	\]
	Following this argument mutatis mutandis for the functors $R[-,-]$ gives an invertible $2$-cell
	\[
	\begin{tikzcd}
		\DbeqQl{X}^{\op} \times \DbeqQl{X} \ar[rr, ""{name = Up}]{}{R[-,-]} \ar[d, swap]{}{\Forget_X^{\op} \times \Forget_X} & & \DbeqQl{X} \ar[d]{}{\Forget_X} \\
		\DbQl{X}^{\op} \times \DbQl{X} \ar[rr, swap, ""{name = Down}]{}{R[-,-]} & & \DbQl{X} \ar[from = Up, to = Down, Rightarrow, shorten >= 4pt, shorten <= 4pt]{}{\cong}
	\end{tikzcd}
	\]
	and completes the proof of the proposition for the $\ell$-adic case. Proceeding mutatis mutandis also establishes the same results for the equivariant derived categories $D_G^b(X)$ and $D_G^b(Y)$.
\end{proof}

Following the same proof mutatis mutandis for the topological case gives the following topological results.
\begin{proposition}\label{Prop: Six functor forgetful commute topology}
	Let $h:X \to Y$ be a $G$-equivariant morphism for a topology group $G$ which admits free $n$-acyclic $G$-spaces which are manifolds and left $G$-varieties $X$ and $Y$. Then each of the six functors of Corollaries \ref{Cor: Pseudocone Functors: Monoidal Closed Derived Cat Topological}, \ref{Cor: Pseudocone Functors: Existence of equivariant pullback  for spaces}, \ref{Cor: Pseudocone Functors: Topological EDC has equivariant pushpull}, and \ref{Cor: Pseudocone Functors: Topological EDC has equivariant exceptional pushpull} between $D^b_G(X)$ and $D_G^b(Y)$ commute with the forgetful functors up to isomorphism.
\end{proposition}

Putting everything together in this section gives the theorems below which state that in both the topological and geometric cases the equivariant derived categories are good homes for equivariant cohomology. The good news, however, is that we have already done most of the work required to show this and now just have to collect the various corollaries, propositions, and such to record this explicitly. In particular, this tells us quite practically that the existence and structure of the equivariant six functor formalism follows from the properties of the $\PC(-)$ construction and how it interacts with the corresponding equivariant geometry and equivariant topology.

\begin{Theorem}\label{Thm: Section SFF: The big results on the EDC for schemes}
	Let $G$ be a smooth algebraic group and let $h:X \to Y$ be a $G$-equivariant morphism of left $G$-varieties. Then the following hold for both pseudofunctors $D_c^b(-;\overline{\Q}_{\ell}):\Var_{/K}^{\op} \to \fCat$ and $D_c^b(-):\Var_{/K}^{\op} \to \fCat$.
	\begin{enumerate}
		\item The categories $D_G^b(X)$ and $D^b_G(X;\overline{\Q}_{\ell})$ are symmetric monoidal closed categories.
		\item The categories $D_G^b(X)$ and $D_G^b(X;\overline{\Q}_{\ell})$ are triangulated categories which have standard and perverse $t$-structures for which
		\[
		D_G^b(X;\overline{\Q}_{\ell})^{\heartsuit_{\operatorname{stand}}} \simeq \Shv_G^{\operatorname{na{\ddot{\imath}}ve}}(X;\overline{\Q}_{\ell})
		\]
		\[
		D_G^b(X;\overline{\Q}_{\ell})^{\heartsuit_{\operatorname{per}}} \simeq \Per_G^{\operatorname{na{\ddot{\imath}}ve}}(X;\overline{\Q}_{\ell})
		\]
		and similarly for $D_G^b(X)$.
		\item There are equivariant pullback and equivariant exceptional pullback functors $h^{\ast}:D_G^b(Y;\overline{\Q}_{\ell}) \to D_G^b(X;\overline{\Q}_{\ell})$ and $h^{!}:D_G^b(Y;\overline{\Q}_{\ell}) \to D_G^b(X;\overline{\Q}_{\ell})$ as well as equivariant pushforward and equivariant proper pushforward functors $Rh_{\ast}:D_G^b(X;\overline{\Q}_{\ell}) \to D_G^b(Y;\overline{\Q}_{\ell})$ and $Rh_{!}:D_G^b(X;\overline{\Q}_{\ell}) \to D_G^b(Y;\overline{\Q}_{\ell})$. The same also holds for $D_G^b(X)$ and $D_G^b(Y)$.
		\item There are adjunctions $h^{\ast} \dashv Rh_{\ast}$ and $Rh_{!} \dashv h^{!}$.
		\item For any $G$-equivariant open immersion $j:U \to X$ and corresponding $G$-equivariant closed immersion $i:V \to X$ there is an open/closed distinguished triangle
		\[
		\begin{tikzcd}
			Rj_!\left(j^!(-)\right) \ar[r]{}{\epsilon^{j}_{-}} & \id_{D_G^b(X;\overline{\Q}_{\ell})}(-) \ar[r]{}{\eta_{-}^{i}} & Ri_{\ast}\left(i^{\ast}(-)\right) \ar[r] & \left(Rj_!\left(j^!(-)\right)\right)[1]
		\end{tikzcd}
		\]
		in $D_G^b(X;\overline{\Q}_{\ell})$ and the same is true of $D_G^b(X)$.
		\item The standard and perverse $t$-structure cohomology functors on $D_G^b(X;\overline{\Q}_{\ell})$ commute up to isomorphism with the forgetful functors and the same is true of the standard and perverse $t$-structures on $D_G^b(X)$.
		\item The six functors $(-)\otimes(-),$ $R[-,-]$, $h^{\ast}$, $Rh_{\ast}$, $Rh_{!}$, and $h^{!}$ all commute up to isomorphism with the corresponding forgetful functors.
		\item There are natural isomorphisms
		\begin{align*}
			Rh_{!}(A) \us{Y}{\otimes} B &\cong Rh_{!}\left(A \us{X}{\otimes} h^{\ast}(B)\right), \\
			\left[Rh_{!}(A), B\right]_{Y} &\cong Rh_{\ast}\left(\left[A, h^{!}(B)\right]_{X}\right), \\
			h^{!}\left([C,B]_{Y}\right) &\cong \left[h^{\ast}(C), h^{!}(B)\right]_{X}
		\end{align*}
		where $A \in D_G^b(X; \overline{\Q}_{\ell})_0$ and $B, C \in D_G^b(Y;\overline{\Q}_{\ell})_0$. The same also holds for $D_G^b(X)$ and $D_G^b(Y)$.
	\end{enumerate}
\end{Theorem}
\begin{proof}
	We prove the different aspects of the theorem in the order listed in the statement of the theorem.
	\begin{enumerate}
		\item The fact that both categories $D_G^b(X)$ and $D_G^b(X;\overline{\Q}_{\ell})$ are symmetric monoidal closed is Corollary \ref{Cor: Pseudocone Functors: Monoidal Closed Derived Cat}.
		\item That both $D_G^b(X;\overline{\Q}_{\ell})$ and $D_G^b(X)$ are triangulated follows from Corollary \ref{Cor: Section Triangle: EDC is triangle} while the existence of the standard and perverse $t$-structures follow from Corollaries \ref{Cor: Section Triangle: Pseudofunctor of standard hearts is heart of pseudofunctor} and \ref{Cor: Section Triangle: Psedofunctor of perv is equiv perv}.
		\item The existence of the functors listed follows from Corollaries \ref{Cor: Pseudocone Functors: Existence of equivariant pullback  for schemes}, \ref{Cor: Pseudocone Functors: Pushforward functors for equivariant maps for scheme sheaves}, and \ref{Cor: Pseudocone Functors: Exceptional Pushforward functors for equivariant maps for scheme sheaves}.
		\item The adjunctions listed are given in Corollaries \ref{Cor: Pseudocone Functors: Topological EDC has equivariant pushpull} and \ref{Cor: Pseudocone Functors: Exceptional Pushforward functors for equivariant maps for scheme sheaves}.
		\item The existence of the open/closed distinguished triangle is Corollary \ref{Cor: Exact open/closed triangle}.
		\item The commutativity of the cohomology functors with the forgetful functors, up to isomorphism, is the content of Proposition \ref{Prop: Forget restricts to naive forget}.
		\item The commutativity of the six functors with the forgetful functors is Proposition \ref{Prop: Six functor forgetful commute}.
		\item The given natural isomorphisms follow from Corollary \ref{Cor: Section Triangle: SFF Relations for EDC Scheme}.
	\end{enumerate}
\end{proof}
\begin{Theorem}\label{Thm: Section SFF: The big results on the EDC for spaces}
	Let $G$ be a topological group which admits $n$-acyclic free $G$-spaces which are manifolds and let $h:X \to Y$ be a $G$-equivariant morphism of left $G$-spaces. Then the following hold for the pseudofunctor $D_c^b(-):\Top^{\op} \to \fCat$.
	\begin{enumerate}
		\item The category $D_G^b(X)$ is symmetric monoidal closed.
		\item The category $D_G^b(X)$ is triangulated and has standard and perverse $t$-structures for which
		\[
		D_G^b(X;\overline{\Q}_{\ell})^{\heartsuit_{\operatorname{stand}}} \simeq \Shv_G^{\operatorname{na{\ddot{\imath}}ve}}(X;\overline{\Q}_{\ell}),
		\]
		and
		\[
		D_G^b(X;\overline{\Q}_{\ell})^{\heartsuit_{\operatorname{per}}} \simeq \Per_G^{\operatorname{na{\ddot{\imath}}ve}}(X;\overline{\Q}_{\ell}).
		\]
		\item There are equivariant pullback and equivariant exceptional pullback functors $h^{\ast}:D_G^b(Y) \to D_G^b(X)$ and $h^{!}:D_G^b(Y) \to D_G^b(X)$ as well as equivariant pushforward and equivariant proper pushforward functors $Rh_{\ast}:D_G^b(X) \to D_G^b(Y)$ and $Rh_{!}:D_G^b(X) \to D_G^b(Y)$.
		\item There are adjunctions $h^{\ast} \dashv Rh_{\ast}$ and $Rh_{!} \dashv h^{!}$.
		\item For any $G$-equivariant open inclusion $j:U \to X$ and corresponding $G$-equivariant closed inclusion $i:V \to X$ there are open/closed distinguished triangles
		\[
		\begin{tikzcd}
			Rj_!\left(j^!(-)\right) \ar[r]{}{\epsilon^{j}_{-}} & \id_{D_G^b(X;\overline{\Q}_{\ell})}(-) \ar[r]{}{\eta_{-}^{i}} & Ri_{\ast}\left(i^{\ast}(-)\right) \ar[r] & \left(Rj_!\left(j^!(-)\right)\right)[1]
		\end{tikzcd}
		\]
		and
		\[
		\begin{tikzcd}
			Ri_!\left(i^!(-)\right) \ar[r]{}{\epsilon^i_{-}} & \id_{D^b_G(X)}(-) \ar[r]{}{\eta^j_{-}} & Rj_{\ast}\left(j^{\ast}(-)\right) \ar[r] & \left(Ri_!\left(i^!(-)\right)\right)[1]
		\end{tikzcd}
		\]
		in $D_G^b(X)$.
		\item The standard and perverse $t$-structure cohomology functors on $D_G^b(X)$ commute up to isomorphism with the forgetful functors.
		\item The six functors $(-)\otimes(-),$ $R[-,-]$, $h^{\ast}$, $Rh_{\ast}$, $Rh_{!}$, and $h^{!}$ all commute up to isomorphism with the corresponding forgetful functors.
	\end{enumerate}
\end{Theorem}
\begin{proof}
	As before, we prove the different aspects of the theorem in the order listed in the statement of the theorem.
	\begin{enumerate}
		\item This is Corollary \ref{Cor: Pseudocone Functors: Monoidal Closed Derived Cat Topological}.
		\item The triangulation of $D_G^b(X)$ is Corollary \ref{Cor: Section Triangle: EDC is triangle topology}. The standard $t$-structure is given in \cite[Proposition 2.5.3]{BernLun} and the perverse $t$-structure is given similarly to Corollary \ref{Cor: Section Triangle: Psedofunctor of perv is equiv perv}.
		\item The existence of the functors listed follows from Corollaries \ref{Cor: Pseudocone Functors: Existence of equivariant pullback  for spaces}, \ref{Cor: Pseudocone Functors: Topological EDC has equivariant pushpull}, and \ref{Cor: Pseudocone Functors: Topological EDC has equivariant exceptional pushpull}.
		\item The adjucntions listed are given in Corollaries \ref{Cor: Pseudocone Functors: Topological EDC has equivariant pushpull} and \ref{Cor: Pseudocone Functors: Topological EDC has equivariant exceptional pushpull}.
		\item The existence of the open/closed distinguished triangle is Corollary \ref{Cor: Section SFF: Openclosed Dist Triangles}.
		\item The commutativity of the cohomology functors with the forgetful functors, up to isomorphism, is described in Proposition \ref{Prop: Forget restricts to naive forget}.
		\item The commutatvity of the six functors with the forgetful functor is given in Proposition \ref{Prop: Six functor forgetful commute topology}.
	\end{enumerate}
\end{proof}
\newpage

\section{Equivariant Stalks and Trace}\label{Section: Trace and Stalks}
We now change gears to the final topic of the monograph: that of equivariant traces and stalks. While we have discussed some of the general motivation for having a theory of traces and stalks already, let us say some words regarding motivation from the Langlands Programme for $p$-adic groups.

In \cite{CFMMX} the authors studied $p$-adic groups using microlocal geometry and the equivariant derived category of $\ell$-adic sheaves by extending the microlocal analytic approach of \cite{ABV} (which was used in the Archimedian setting) to $p$-adic groups. In doing so they developed and gave a systematic study of what they called ABV packets and how those packets gave a geometric approach towards studying the representations (and in particular the Arthur packets) of a $p$-adic group $G$. 

Recall that a $p$-adic group $G$ is an algebraic group $G$ defined over a $p$-adic field $F$, i.e., $G$ is a group object in $\Var_{/F}$ for $F/\Q_p$ a finite degree (algebraic) extension of some field $\Q_p$ of $p$-adic numbers where $p$ is an integer prime. Fix an integer prime $\ell$ coprime to $p$. In \cite[Section 8.3, Conjecture 1]{CFMMX} the authors give a conjecture which relates the geometry of either a special orthogonal $p$-adic group $G$ or a quasisiplit symplectic $p$-adic group to the representation theoretic aspects of $G$ (in particular by comparing the pure Arthur packets to the ABV packets) in a precise and careful way. This conjecture, following \cite[Equation 8.8]{CFMMX}, is equivalent to verifying the ``fixed point formula'' for every perverse sheaf $\Pcal \in \Per_{H_{\lambda}}(V_{\lambda};\overline{\Q}_{\ell})$
\[
\langle \eta_{\psi,s}, \Pcal\rangle = (-1)^{\dim C_{\psi}}\trace_{a_s}\left(\mathsf{NEvs}_{\psi}\Pcal\right)
\]
where the pairing on the left takes place in the Grothendieck group $\mathsf{K}\Per_{H_{\lambda}}(V_{\lambda};\overline{\Q}_{\ell})$ and the trace operations are calculated in the category $D^b_{H_{\lambda}}(V_{\lambda};\overline{\Q}_{\ell})$. Note that we have left the meaning of $\lambda,$ $\psi$, $H_{\lambda}$, $V_{\lambda}$, $s$, $a_s$, $C_{\psi}$, $\eta_{\psi,s}$, and $\mathsf{NEvs}$ purposely vague; the important thing here is that these representation-theoretic conjectures ultimately reduce to performing a calculation in the equivariant derived category $D_{H_{\lambda}}^b(V_{\lambda};\overline{\Q}_{\ell})$ involving traces and functorial constructions. Consequently, in order to study this conjecture in detail from a categorical geometric perspective requires a more in depth understanding of trace operators in the categories $D_G^b(X;\overline{\Q}_{\ell})$ and, more generally, in the pseudocone categories $\PC(F)$ for a pre-equivariant pseudofunctor $(F,\overline{F})$ on $X$.

We begin our stalk and trace journey by discussing how to define equivariant stalks. While our first focus will be on extending the definition of stalks of geometric points to the equivariant situation, we will at least sketch how to define the stalks of topological points to the equivariant situation as this will be relevant for both the geometric and topological cases. First, however, let us recall the definition of a geometric point of a $K$-variety $X$.

\begin{definition}
	Let $K$ be a field with a fixed separable closure $K^{\operatorname{sep}}$ and let $X$ be a $K$-variety. A geometric point\index[terminology]{Geometric point} of $X$ is a morphism $x:\Spec K^{\operatorname{sep}} \to X$.
\end{definition}
\begin{definition}
	Let $K$ be a field and let $X$ be a $K$-variety with a geometric point $x:\Spec K^{\operatorname{sep}} \to X$. An {\'e}tale neighbourhood\index[terminology]{Neighbourhood! {\'E}tale} of $x$ is a pair $(U\xrightarrow{\varphi}X,u)$ where $\varphi$ is an {\'e}tale morphism and $u$ is a geometric point $u:\Spec K^{\operatorname{sep}} \to U$ making the diagram
	\[
	\begin{tikzcd}
		& \Spec K^{\operatorname{sep}} \ar[dr]{}{x} \ar[dl, swap]{}{u} \\
		U \ar[rr, swap]{}{\varphi} & & X
	\end{tikzcd}
	\]
	commute.
\end{definition}
\begin{remark}
	By abuse of notation we will often write an {\'e}tale neighbourhood of a geometric point $x$ as $(U,u)$ instead of $(U \xrightarrow{\varphi}X,u)$.
\end{remark}
\begin{lemma}\label{Lemma: Section SFF: Translations along equivariant etale}
	Let $x:\Spec K^{\operatorname{sep}} \to X$ be a geometric point of $X$ with an {\'e}tale neighbourhood  $(U,u)$ of $x$. Then if $G$ is a smooth algebraic group and $U$ and $X$ are left $G$-varieties for which $\varphi$ is $G$-equivariant then for any pseudofunctor $\overline{F}:\Var_{/K}^{\op} \to \fCat$ there is a pseudocone translation
	\[
	\begin{tikzcd}
		\SfResl_G(X)^{\op} \ar[rr, ""{name = U}]{}{\overline{F} \circ \quo_{G\backslash(-)}^{\op}} \ar[dr, swap]{}{\gamma} & & \fCat \\
		& \SfResl_G(U)^{\op} \ar[ur, swap]{}{\overline{F} \circ \quo_{G \backslash(-)}^{\op}} \ar[from = U, to = 2-2, Rightarrow, shorten <= 4pt, shorten >= 4pt]{}{\varphi^{\square}}
	\end{tikzcd}	
	\]
	where $\gamma$ is the isomorphism of Remark \ref{Remark: Pseudocone Section: SfResl is iso to Sf}.
\end{lemma}
\begin{proof}
	The assumptions on $U$ being equivariant mean that the quotient varieties $G \backslash (\Gamma \times U) = \quo_{G \backslash (\Gamma \times U)}$, and hence the quotient functor, are defined for every $\Gamma \in \Sf(G)_0$ and make the cube
	\[
	\begin{tikzcd}
		& \Gamma^{\prime} \times U \ar[rr]{}{\id_{\Gamma^{\prime}} \times \varphi} \ar[dd, near start]{}{\quo_{\Gamma^{\prime} \times U}} & & \Gamma^{\prime} \times X \ar[dd]{}{\quo_{\Gamma^{\prime} \times X}} \\
		\Gamma \times U \ar[dd, swap]{}{\quo_{\Gamma \times U}} \ar[ur]{}{f \times \id_U} & & \Gamma \times X \ar[ur]{}{f \times \id_X} \\
		& \quot{U}{\Gamma^{\prime}} \ar[rr, near start]{}{G \backslash (\id_{\Gamma^{\prime}} \times \varphi)} & & \quot{X}{\Gamma^{\prime}} \\
		\quot{U}{\Gamma} \ar[rr, swap]{}{G \backslash (\id_{\Gamma} \times \varphi)} \ar[ur, swap]{}{G \backslash (f \times \id_{U})} & & \quot{X}{\Gamma} \ar[ur, swap]{}{G \backslash (f \times \id_X)}
		\ar[from = 2-1, to = 2-3, crossing over, swap, near end]{}{\id_{\Gamma} \times \varphi} \ar[from = 2-3, to = 4-3, crossing over, near start]{}{\quo_{\Gamma \times X}}
	\end{tikzcd}
	\]
	commute. The result now follows from a routine application of the proof of Theorem \ref{Thm: Pseudocone Functors: Pullback induced by fibre functors in pseudofucntor}.
\end{proof}
Before we proceed to discuss the stalks of at geometric points we need to verify that if $\varphi:U \to X$ is {\'e}tale then so is $G \backslash (\id_{\Gamma} \times \varphi):\quot{U}{\Gamma} \to \quot{X}{\Gamma}$ for every $\Gamma \in \Sf(G)_0$, as otherwise the evaluation $\quot{A}{\Gamma}(\quot{U}{\Gamma})$ we anticipate need not be valid.
\begin{lemma}
	if $\varphi:U \to X$ is {\'e}tale and a $G$-equivariant morphism then $G \backslash (\id_{\Gamma} \times \varphi):\quot{U}{\Gamma} \to \quot{X}{\Gamma}$ is {\'e}tale for every $\Gamma \in \Sf(G)_0$.
\end{lemma}
\begin{proof}
	Recall that by construction we have a commuting diagram
	\[
	\begin{tikzcd}
		\Gamma \times U \ar[rr]{}{\id_{\Gamma} \times \varphi} \ar[d, swap]{}{\quo_{\Gamma \times U}} & & \Gamma \times X \ar[d]{}{\quo_{\Gamma \times X}} \\
		G \backslash (\Gamma \times U) \ar[rr, swap]{}{G \backslash(\id_{\Gamma} \times \varphi)} & & G \backslash (\Gamma \times X)
	\end{tikzcd}
	\]
	for every $\Gamma$ in $\Sf(G)$ with $\quo_{\Gamma \times U}$ smooth and surjective and $\quo_{\Gamma \times X} \circ (\id_{\Gamma} \times \varphi)$ smooth. It then follows by Lemma \ref{Lemma: Section 2.1: Smooth morphism lemma}. Since $G \backslash (\id_{\Gamma} \times \varphi)$ is smooth, it suffices to verify that it is of relative dimension zero to conclude that it is {\'e}tale. However, this follows from the fact that $\Gamma$ is an {\'e}tale locally trivializable $G$-fibration, the identity is zero-dimensional, and the fact that $\varphi$ is of relative dimension zero.
\end{proof}

These pseudonatural transformations may then be constructed for every morphism between {\'e}tale neighbourhoods of a geometric point. This allows us us to consider a filtered colimit of the form
\[
\colim_{\substack{(U,u)\,\text{{\'e}tale neighbourhood} \\ \text{of the geometric point}\,x}} \quot{A}{\Gamma}(\quot{U}{\Gamma})
\]
for every $\Gamma \in \Sf(G)_0$ and every object $A$ in some equivariant category $F_G(X)$ provided that the object $\quot{A}{\Gamma}$ can be ``evaluated'' at the scheme $\quot{U}{\Gamma}$ and the colimit displayed exists; this makes sense, for instance, if $\quot{A}{\Gamma}$ is an {\'e}tale sheaf or a (derived) complex of {\'e}tale sheaves or even a (derived) complex of (pre)sheaves on $\Sch_{/\XGamma}$, if $\quot{A}{\Gamma}$ is a pseudofunctor, or some other such similar situation; cf.\@ Definitions \ref{Defn: Section SFF: Evaluate at opens} and \ref{Defn: Section SFF: Evaluate at etale} for a more precise notion of what we mean by evaluating at an ({\'e}tale) open and Remark \ref{Remark: Section SFF: Evlatuation not precise} for a discussion of some inherent imprecision in this notion.

This consideration highlights the two main issues which arise in our search for stalks of an equivariant object:
\begin{enumerate}
	\item Whatever categories the objects $\quot{A}{\Gamma}$ arise in, we at least need to ensure that the colimits
	\[
	\colim_{\substack{(U,u)\,\text{{\'e}tale neighbourhood} \\ \text{of geometric point}\,x}} \quot{A}{\Gamma}(\quot{U}{\Gamma})
	\]
	make sense to discuss in order to hope to have a notion of a stalk at $x$.
	\item We have a large amount of choice as to \emph{which} smooth free $G$-fibration $\Gamma$ we use when considering the stalks; is there one fibration $\Gamma$ we \emph{should} choose over the others?
\end{enumerate}
The first issue is a technical issue. Stalks only exist, even in the classical sheaf-theoretic sense, when the underlying sheaf category considers sheaves with values in a categroy $\Ascr$ which admits certain types of colimits\footnote{We leave this purposefully vague; the type of colimits we ask $\Ascr$ to satisfy changes based on site-theoretic considerations on the domain of the sheaves. For instance, if we are considering sheaves on a topological space it suffices to ask for $\Ascr$ to admit filtered colimits in order to have a notion of stalks.}. We can axiomatize the properties required of $\Ascr$ in order to have a notion of stalks which is sufficient for our site-theoretic needs; in this case, in terms of the {\'e}tale topology or, later on, in terms of the open topology. 

The second issue is more subtle and requires some choices. There are a few different ways we can handle the issue of having to choose the $G$-fibration $\Gamma$ which is ``most natural'' for defining and considering stalks:
\begin{enumerate}
	\item Avoid the choice altogether and work with the family of stalks
	\[
	A_{x} = \left\lbrace \colim_{\substack{(U,u)\,\text{{\'e}tale neighbourhood} \\ \text{of geometric point}\,x}} \quot{A}{\Gamma}\left(\quot{U}{\Gamma}\right) \; : \; \Gamma \in \Sf(G)_0 \right\rbrace
	\]
	all at the same time.
	\item Wrap up the choices into one object by taking a colimit
	\[
	A_x = \colim_{\Gamma \in \Sf(G)_0}\left(\colim_{\substack{(U,u)\,\text{{\'e}tale neighbourhood} \\ \text{of geometric point}\,x}} \quot{A}{\Gamma}\left(\quot{U}{\Gamma}\right)\right).
	\]
	\item Work with a specific $\Gamma_0$ chosen for a reason and define the stalk to be given as
	\[
	A_{x} = \colim_{\substack{(U,u)\,\text{{\'e}tale neighbourhood} \\ \text{of geometric point}\,x}} \quot{A}{\Gamma_0}\left(\quot{U}{\Gamma_0}\right).
	\]
\end{enumerate}
Each of these approaches has its merits and downsides. The first is an approach is consistent with our approach towards equivariance in the sense that it gives us a stalk for each resolution $\Gamma \times U$ and then collects the stalks together. As a downside, however, the first approach makes it difficult to produce just one object and forces us to still choose a way of producing a stalk at the end. The second approach is a good way of producing a single object because it gives a canonical ``one stalk'' out of all potential stalks and smooths out the geometry of the resolutions and how the object $A$ records said geometry. However, the way the second definition relates to conventional stalks and the forgetful functor is unclear and depending on the categories in which the objects $\quot{A}{\Gamma}$ reside, it is unclear that the colimit indexed by $\Sf(G)_0$ would always exist. The third approach has the benefit that it compares readily with the forgetful functor and also does not ask much formally in order to be defined. However, there had better be a good reason or justification in choosing just one variety $\Gamma_0$ in order to define \emph{the} stalk at $x$ to be given as just the $\quot{A}{\Gamma_0}$ stalk at $x$.

We choose the third approach for defining equivariant stalks and leave examination and study of the two remaining approaches for future work. Our justification comes from the following desideratum:
\begin{desideratum}\label{Desideratum for stals}
	Let $F$ be a pre-equivariant pseudofunctor on $X$ and let $x$ be a geometric point of $X$. Any equivariant stalk of an object $A$ in $F_G(X)$, if it exists, should satisfy
	\[
	\left(\Forget(A)\right)_x \cong \Forget(A_x).
	\] 
\end{desideratum}
Consequently we are obliged to choose $\Gamma_0 := G$ as our variety and then define what we mean by evaluating the objects of $F$ at {\'e}tale opens before we can define the equivariant stalk at $x$. We also must do the same for topological points (both in the geometric and topological cases), but this is deferred to Definitions \ref{Defn: Top point for space}, \ref{Defn: Top point for scheme}, and \ref{Defn: Equiv stalks for topological spaces}.
\begin{definition}\label{Defn: Section SFF: Evaluate at opens}
	Let $X$, $G$, and $F$ be given as in either of the two cases below:
	\begin{enumerate}
		\item $G$ is a topological group, $X$ is a left $G$-space, and $(F,\overline{F}) = F$ is a pre-equivariant pseudofunctor on $X$;
		\item $G$ is a smooth algebraic group, $X$ is a left $G$-variety, and $(F,\overline{F}) = F$ is a pre-equivariant pseudofunctor on $X$.
	\end{enumerate}
	Then in either case we say that the objects of $F$ can be evaluated at open subspaces of $X$\index[terminology]{Pre-equivariant Pseudofucntor! Evalaution at Opens} if for every object $A$ of $F(\quot{X}{G})$ and for every open subobject $U$ of $\quot{X}{G}$, there is an assignment $A(U)$ taking values (pseudo)functorially in $\Open(\quot{X}{G})^{\op}$ in some category $\Cscr$.
\end{definition}
\begin{definition}\label{Defn: Section SFF: Evaluate at etale}
	Let $G$ be a smooth algebraic group, $X$ a left $G$-variety, and $(F,\overline{F}) = F$ is a pre-equivariant pseudofunctor on $X$. Then we say that the objects of $F$ can be evaluated at {\'e}tale opens of $X$\index[terminology]{Pre-equivariant Pseudofucntor! Evalaution at {\'E}tale Opens} if for every object $A$ of $F(\quot{X}{G})$ and for every scheme $U$ of $\operatorname{\acute{E}t}(\quot{X}{G})$, there is an assignment $A(U)$ taking values functorially in $\operatorname{\acute{E}t}(X)^{\op}$ in some category $\Cscr$.
\end{definition}
\begin{remark}\label{Remark: Section SFF: Evlatuation not precise}
	We have been somewhat imprecise in our definition of what it means to evaluate an object at an open of $X$. The definitions given are sufficiently general that we can capture the primary examples of interest to us: sheaves on a space/scheme $X$ in the open topology, quasi-coherent sheaves on a scheme, sheaves on a scheme in the {\'e}tale topology, $\ell$-adic sheaves on a scheme, and each of the corresponding chain complex and derived categories. However, because these categories need not be a sheaf category, presheaf category, or even a derived category of a (pre)sheaf category\footnote{Most distressingly this is seen with both the categories of $\ell$-adic sheaves and the $\ell$-adic derived category; infamously $D_c^b(X;\overline{\Q}_{\ell})$ arises as a pseudolimit of a construction which itself is based on the Ore localization of a subcategory of $[\Nbb^{\op},\mathbf{TorShv}(X)]$; cf.\@ \cite[Expos{\'e} V, VI]{SGA5} or \cite{FreitagKiehl} for explicit details.} we have left the definition of what it means to evaluate at this level so as to apply to at least the examples listed above.
\end{remark}
\begin{remark}
	Perhaps a better name for the structure described in Definitions \ref{Defn: Section SFF: Evaluate at opens} and \ref{Defn: Section SFF: Evaluate at etale} would be to say that $F$ may be globally evaluated at ({\'e}tale) opens of $X$ and instead leave the definition of evaluation at ({\'e}tale) opens without a modifier to mean that for all $\Gamma \in \Sf(G)_0$ (or for all $\Gamma \in \Resl_G(X)_0$ in the topological space case) that the evaluations $A(U)$ are defined for all objects $A$ of $F(\XGamma)$ and for all opens (respectively {\'e}tale opens) $U$ of $\XGamma$ (or for $G \backslash \Gamma$ in the topological case) of the relevant type. In either case, both situations can be readily adapted to other Grothendieck topologies over each object $\XGamma$, but for the scope of this work it suffices to focus on both the open and {\'e}tale topologies as well as the case where the objects $A$ of $F(\quot{X}{G})$ may be evaluated at the opens (in whichever relevant topology) of $\quot{X}{G}$.
\end{remark}
\begin{definition}\label{Defn: Equiv stalks, geometric point}
	Let $G$ be a smooth algebraic group and let $X$ be a left $G$-variety with a geometric point $x$. If $F$ is a pre-equivariant pseudofunctor on $X$ for which $\overline{F}(\quot{X}{G})$ admits filtered colimits and for which each object $A$ of $\overline{F}(\quot{X}{G})$ may be evaluated at {\'e}tale schemes over $\quot{X}{G}$. We then define the equivariant stalk\index[terminology]{Equivariant Stalk! {\'E}tale} of an object $A$ of $F_G(X)$ at $x$ by
	\[
	A_x := \colim_{\substack{(U,u)\,\text{{\'e}tale neighbourhood} \\ \text{of geometric point}\,x}} \quot{A}{G}\left(\quot{U}{G}\right).
	\]
\end{definition}

We now give two short sanity check results. The first shows that our choice to index the colimit defining the stalk by ({\'e}tale) opens of $X$ does not pose a problem while the second that this definition indeed satisfies Desideratum \ref{Desideratum for stals}.
\begin{proposition}
	For any smooth algebraic group $G$ and any left $G$-variety $X$ there is an equivalence of categories $\operatorname{\acute{E}t}(X) \simeq \operatorname{\acute{E}t}(\quot{X}{G})$.
\end{proposition}
\begin{proof}
	This is immediate from the isomorphism of varieties $X \cong \quot{X}{G}$.
\end{proof}
\begin{lemma}\label{Lemma: Section SFF: Stalk satisfy desiderat}
	Let $G$ be a smooth algebraic group and let $X$ be a left $G$-variety with a geometric point $x$. If $F$ is a pre-equivariant pseudofunctor on $X$ for which $F(\XGamma)$ admits filtered colimits for all objects $\Gamma \times X$ in $\SfResl_G(X)$ and for which each object $\quot{A}{\Gamma}$ may be evaluated at {\'e}tale schemes over $\XGamma$. Then for any object $A$ of $F_G(X)$, $A$ satisfies Desideratum \ref{Desideratum for stals}:
	\[
	\Forget(A_x) \cong \left(\Forget(A)\right)_x.
	\]
\end{lemma}
\begin{proof}
	Recall that $\psi:\quot{X}{G} \to X$ is the non-equivariant isomorphism of varieties furnishing the existence of the forgetful functor $\Forget:F_G(X) \to \overline{F}(X)$. Then we calculate on one hand
	\[
	\Forget(A)_X = \overline{F}(\psi)\left(\quot{A}{G}\right)_x = \colim_{\substack{(U,u)\,\text{{\'e}tale neighbourhood} \\ \text{of geometric point}\,x}}\left(\overline{F}(\psi)\left(\quot{A}{G}\right)\right)(U)
	\]
	while on the other hand
	\begin{align*}
		\Forget(A_x) &= \Forget\left(\colim_{\substack{(U,u)\,\text{{\'e}tale neighbourhood} \\ \text{of geometric point}\,x}}\quot{A}{G}(\quot{U}{G})\right) \\
		&= \overline{F}(\psi)\left(\colim_{\substack{(U,u)\,\text{{\'e}tale neighbourhood} \\ \text{of geometric point}\,x}}\quot{A}{G}(\quot{U}{G})\right)
	\end{align*}
	Now observe that since $\psi$ is an isomorphism, $\overline{F}(\psi)$ is an equivalence of categories and hence is cocontinuous. Consequnetly there is an isomorphism
	\[
	\overline{F}(\psi)\left(\colim_{\substack{(U,u)\,\text{{\'e}tale neighbourhood} \\ \text{of geometric point}\,x}}\quot{A}{G}(\quot{U}{G})\right) \cong \colim_{\substack{(U,u)\,\text{{\'e}tale neighbourhood} \\ \text{of geometric point}\,x}} \overline{F}(\psi)\left(\quot{A}{G}(\quot{U}{G})\right).
	\]
	Note that there are non-equivariant isomorphisms $\psi_U:\quot{U}{G} \to U$ such that
	\[
	\begin{tikzcd}
	\quot{U}{G} \ar[r]{}{\psi_U} \ar[d, swap]{}{\varphi} & U \ar[d]{}{\varphi} \\
	\quot{X}{G} \ar[r, swap]{}{\psi_X} & X
	\end{tikzcd}
	\] 
	commutes for every {\'e}tale map $\varphi:U \to X$. This, combined with the pseudofunctoriality of $\overline{F}$, ensures that there is a natural isomorphism
	\[
	\overline{F}(\psi)\left(\quot{A}{G}(\quot{U}{G})\right) \cong \left(\overline{F}(\psi)\left(\quot{A}{G}\right)(\right)U)
	\]
	in $\text{{\'E}t}(X)$ and hence gives rise to an isomorphism
	\[
	\colim_{\substack{(U,u)\,\text{{\'e}tale neighbourhood} \\ \text{of geometric point}\,x}} \overline{F}(\psi)\left(\quot{A}{G}(\quot{U}{G})\right) \cong \colim_{\substack{(U,u)\,\text{{\'e}tale neighbourhood} \\ \text{of geometric point}\,x}}\left(\overline{F}(\psi)\left(\quot{A}{G}\right)\right)(U).
	\]
	Putting these together gives the isomorphism
	\[
	\Forget(A_x) \cong \Forget(A)_x
	\]
	and hence proves the lemma.
\end{proof}

The benefit of this perspective is that we can give the definition of a topological stalk of $X$ mutatis mutandis and also derive that this notion of an equivariant topological stalk also commutes with the forgetful functor. We record both these below after recalling the definitions of topological points both for topological spaces and for schemes.

\begin{definition}\label{Defn: Top point for space}
	Let $X$ be a topological space. A (topological) point of $X$ is a map $x:\lbrace \ast \rbrace \to X$.\index[terminology]{Topological point! Of a Topological Space}
\end{definition}
\begin{definition}\label{Defn: Top point for scheme}
	Let $X$ be a scheme. A topological point\index[terminology]{Topological point! Of a Scheme} of $X$ is a map $x:\Spec k \to X$ where $k \cong \kappa(x) = \Ocal_{X,x}/\mfrak_{X,x}$ for some $x \in \lvert X \rvert$.
\end{definition}
\begin{remark}
	While Definition \ref{Defn: Top point for space} seems like it would be best phrased in terms of saying a point is an element $x \in X$, we want to think of points in terms of maps. Because the category of topological spaces is reasonably well-behaved the underlying set of a space $X$ may be recovered by $X \cong \Top(\lbrace \ast \rbrace, X)$, i.e., the terminal object of $X$ co-represents the identity functor on $\Top$.
	
	Alternatively, writing the elements $x \in \lvert X \rvert$ when $X$ is a scheme cannot be done directly. As such, we need to exploit the fact that schemes $X$ are locally ringed spaces. Because each local ring $\Ocal_{X,x}$ of $X$ admits a residue field $\kappa(x) := \Ocal_{X,x}/\mfrak_{X,x}$ we get that there is a scheme map $\Spec \kappa(x) \to X$ given by $x \mapsto x$ locally and defining the sheaf map $\Ocal_X \to i^{\ast}\Ocal_{\kappa(x)}$ by noting for every open $U \subseteq \lvert X \rvert$, if $x \in U$ then there is a canonical map given by the composition
	\[
	\Ocal_X(U) \to \colim_{\substack{U \subseteq X\,\text{open} \\x \in U}} \Ocal_{X}(U) \to \kappa(x)
	\]
	and if $x \notin U$ then $i_{\ast}\Ocal_{\kappa(x)}(U) = \Ocal_{\kappa(x)}(i^{-1}(U)) = 0$ so we simply use the zero map in this case.
\end{remark}

\begin{definition}\label{Defn: Equiv stalks for topological spaces}
	Let $G$ be a topological group and let $X$ be a topological space. If $F$ is a pre-equivariant pseudofunctor for which each category $\overline{F}(G \backslash (G \times X))$ admits filtered colimits objects of $F(G \backslash (G \times X))$ may be evaluated on opens of $X$, then if $x \in X$ is any point of $X$ we define the equivariant stalk\index[terminology]{Equivariant Stalk! Of a Topological Point of a Topological Space} of an object $A$ of $F_G(X)$ to be given by
	\[
	A_x := \colim_{\substack{U \subseteq X\,\text{open} \\x \in U}}\quot{A}{G}\left(\quot{U}{G}\right).
	\]
\end{definition} 
\begin{definition}\label{Defn: Equiv topolgical stalks for schemes}
	Let $G$ be a smooth algebraic group and let $X$ be a left $G$-variety with a topological point $x$.\index[terminology]{Equivariant Stalk! Of a Topological Point of a Variety} If $F$ is a pre-equivariant pseudofunctor on $X$ for which $F(\quot{X}{G})$ admits filtered colimits and objects can be evaluated at topological points of $X$ then we define the equivariant stalk of an object $A$ of $F_G(X)$ at $x$ by
	\[
	A_x := \colim_{\substack{U \subseteq X\,\text{open} \\x \in U}} \quot{A}{G}\left(\quot{U}{G}\right).
	\]
\end{definition}
\begin{remark}
	In both cases above note that the definitions type because $\quot{U}{G}$ is open by Lemmas \ref{Lemma: Section: Yet More Stuff: Open immerions descend} and \ref{Lemma: Section SFF: Open inclusions descend}.
\end{remark}
\begin{lemma}\label{Lemma: Section SFF: Topological stalks sat desideratum}
	Let $G, X,$ $x$, and $F$ be either of the following:
	\begin{enumerate}
		\item $G$ is a topological group, $X$ is a left $G$-space, $x:\lbrace \ast \rbrace \to X$ is a point, and $F$ is a pre-equivariant pseudofunctor such that the category $\overline{F}(G \backslash (G \times X))$ admits filtered colimits and for which each object $A$ of $\overline{F}(G \backslash (G \times X))$ can be evaluated at opens of $G \backslash (G \times X)$.
		\item $G$ is a smooth algebraic group, $X$ is a left $G$-variety, $x$ is a topological point of $X$, and $F$ is a pre-equivariant pseudofunctor such that the category $\overline{F}(\quot{X}{G})$ admits filtered colimits and for which each object $A$ of $\overline{F}(\quot{X}{G})$ can be evaluated at the opens $\quot{U}{\Gamma}$ for $U$ open in $X$.
	\end{enumerate} 
	Then in either case for any object $A$ of $F_G(X)$, $A$ satisfies Desideratum \ref{Desideratum for stals}:
	\[
	\Forget(A_x) \cong \left(\Forget(A)\right)_x.
	\]
\end{lemma}
\begin{proof}
	This follows mutatis mutandis to the proof of Lemma \ref{Lemma: Section SFF: Stalk satisfy desiderat} and uses the fact that the non-equivariant isomorphisms $\psi$ are isomorphisms in the exact same way.
\end{proof}

We now introduce some terminology to describe the pre-equivariant pseudofunctors for which we can take equivariant stalks as well as those which give rise to stalks which behave functorially over the category of $G$-objects.

\begin{definition}\label{Defn: Section SFF: Herbacesous Preeqs}
	Let $G$, $X$, and $(F,\overline{F})$ satisfy any of the cases below:
	\begin{enumerate}
		\item $G$ is a smooth algebraic group, $X$ is a left $G$-variety, and $(F,\overline{F})$ is a pre-equivariant pseudofunctor on $X$ such that the category $\overline{F}(\quot{X}{G})$ admits filtered colimits and each object $A$ of $\overline{F}(\quot{X}{G})$ can be evaluated on {\'e}tale schemes $V \to \quot{X}{G}$.
		\item $G$ is a smooth algebraic group, $X$ is a left $G$-variety, and $(F,\overline{F})$ is a pre-equivariant pseudofunctor on $X$ such that the category $\overline{F}(\quot{X}{G})$ admits filtered colimits and each object $A$ of $\overline{F}(\quot{X}{G})$ can be evaluated on open subschemes $U \to \quot{X}{G}$.
		\item $G$ is a topological group, $X$ is a left $G$-space, and $(F,\overline{F})$ is a pre-equivariant pseudofunctor on $X$ such that the category $\overline{F}(\quot{X}{G})$ admits filtered colimits and each object $A$ of $\overline{F}(\quot{X}{G})$ can be evaluated on open subspaces $U \subseteq \quot{X}{G}$.
	\end{enumerate}
	Then in Case (1) we say that $(F,\overline{F})$ is an {\'e}tale weakly herbaceous pre-equivariant pseudofunctor\index[terminology]{Pre-equivariant Pseudofunctor! Weakly Herbaceous} and in Cases (2) and (3) $(F,\overline{F})$ is a weakly topological herbaceous pre-equivariant pseudofunctor. 
	
	Finally, if $h:X \to Y$ is any $G$-equivariant morphism, if $(F,\overline{F})$ is a weakly herbaceous pre-equivariant pseudofunctor, and if $(\quot{F}{Y},\overline{F})$ is the induced pre-equivariant pseudofunctor on $Y$ with $h^{\ast}:F_G(Y) \to F_G(X)$ the pullback functor induced by Theorem \ref{Thm: Pseudocone Functors: Pullback induced by fibre functors in pseudofucntor} satisfies the fact that for all relevant points $x$ of $X$ there is an isomorphism
	\[
	\left(h^{\ast}A\right)_{x} \cong A_{h(x)},
	\]
	we say that $(F,\overline{F})$ is an herbaceous\index[terminology]{Pre-equivariant Pseudofunctor! Herbaceous} pre-equivariant pseudofunctor\footnote{The notation $h(x)$ means, in the geometric point case, the map $\Spec k^{\operatorname{sep}} \xrightarrow{x} X \xrightarrow{h} Y$.}.
\end{definition}
\begin{remark}
	The conditions given for a weakly ({\'e}tale or topological) herbaceous pre-equivariant pseudofunctor mean that for every object $A$ of $F_G(X)$ and for any relevant type of point $x$ of $X$, the stalk $A_x$ is well-defined. The additional criterion required for $(F,\overline{F})$ to be a herbaceous pre-equivariant pseudofunctor asks that the stalks vary continuously in the varieties $X$ and $Y$.
\end{remark}
\begin{example}
	Some standard examples of {\'e}tale herbaceous pre-equivariant pseudofunctors are the equivariant derived $\ell$-adic category pre-equivariant pseudofunctor, the {\'e}tale sheaf pre-equivariant pseudofunctor, and the {\'e}tale $D$-module pre-equivariant pseudofunctor on $X$. Some topological herbaceous pre-equivariant pseudofunctors may be defined similarly save for making sure that we take topological sheaves instead of {\'e}tale sheaves in each case.
\end{example}

Now that we have a class of pre-equivariant pseudofunctors with which we can take stalks, we are led to ask what of the equivariant trace? This requires us to tackle understanding two distinct but related aspects of what goes into defining our desired trace object
\[
\trace_g(A,x)
\]
for an element $g \in G$, some specific point $x \in X$, and an object $A \in F_G(X)$. 

To this end we will first explore the aspect of defining the trace which can be done for \emph{any} herbaceous pre-equivariant pseudofunctor. Recall that Theorem \ref{Theorem: Formal naive equivariance} gives the existence of an equivariance (cf.\@ Definition \ref{Defn: Section Equivariance: Equivariance})
\[
\Theta:\alpha_X^{\ast} \xRightarrow{\cong} \pi_2^{\ast}
\]
where $\alpha_X^{\ast}$ and $\pi_2^{\ast}$ are the pullbacks against the action map and the projection given by an application of Theorem \ref{Thm: Pseudocone Functors: Pullback induced by fibre functors in pseudofucntor} in the case of $\pi_2$ and Lemma \ref{Lemma: Section Change of Space: Action map pullback} in the case of $\alpha_X$. To each object $A$ we can thus associate the equivariance $\Theta^A$ which, we will see, forms the backbone of our trace formalism. First, however, in order to have a well-behaved notion of trace we need to understand the interaction between stalks and pullback functors.

\begin{proposition}\label{Prop: Section SFF: Pullback of stalk is stalk of pullbakcs}
	Let $(F,\overline{F})$ be a weakly herbaceous pre-equivariant pseudofunctor on $X$ for which $\overline{F}$ is strict. Then for any point $g$ of $G$ and any point $x$ of $X$ we have that
	\[
	\left((\id_G \times h)^{\ast}\Theta^A_Y\right)_{(g,x)} = \left(\Theta_{X}^{h^{\ast}A}\right)_{(g,x)}.
	\]
\end{proposition}
\begin{proof}
Since $\overline{F}$ is a strict pseudofunctor, we have an equality
\[
(\id_G \times h)^{\ast}\Theta_{Y}^{A} = \Theta_{X}^{h^{\ast}A}
\]
for every object $A$ of $F_G(Y)$ by Corollary \ref{Cor: Section Equivariance: Pullbacks and Equivariances for Strict pseudofunctors}.
\end{proof}

With this above technical detail established we can now define and describe the equivariant trace. We simply need to define the kinds of pre-equivariant pseudofunctors for which we can even hope to define equivariant trace before establishing the basic properties which it satisfies. To do this carefully we simply recall what it means to be fully dualizable in a symmetric monoidal category before defining trace-class pre-equivariant pseudofunctors.

\begin{definition}\label{Defn: Section SFF: dualizable smc object}
	Let $(\Cscr, \otimes, I, \lambda, \rho, \alpha, \gamma)$ be a symmetric monoidal category. We say that an object $X$ in $\Cscr$ is dualizable if $X$ admits a right adjoint $X^{\ast}$ in the delooping bicategory $\underline{B}\Cscr$. More explicitly, $X$ is dualizable\index[terminology]{Dualizable Object} if there is an object $X^{\ast}$ and morphisms $\eta_X:I \to X^{\ast} \otimes X$ and $\epsilon_X:X \otimes X^{\ast} \to I$ for which the diagrams
	\[
	\begin{tikzcd}
		X \otimes I \ar[rr]{}{\id_X \otimes \eta_X} \ar[d, swap]{}{\lambda_X \circ \rho_X^{-1}} & & X \otimes (X^{\ast} \otimes X) \ar[d]{}{\alpha_{X,X^{\ast},X}^{-1}} \\
		I \otimes X & & (X \otimes X^{\ast}) \otimes X \ar[ll]{}{\epsilon_X \otimes \id_X}
	\end{tikzcd}
	\]
	and
	\[
	\begin{tikzcd}
		I \otimes X^{\ast} \ar[rr]{}{\eta_X \otimes \id_{X^{\ast}}} \ar[d, swap]{}{\rho_{X^{\ast}} \circ \lambda_{X^{\ast}}^{-1}} & & (X^{\ast} \otimes X) \otimes X^{\ast} \ar[d]{}{\alpha_{X^{\ast},X,X^{\ast}}} \\
		X^{\ast} \otimes I & & X^{\ast} \otimes (X \otimes X^{\ast}) \ar[ll]{}{\id_{X^{\ast}} \circ \epsilon_X}
	\end{tikzcd}
	\]
	commute. The subcategory of $\Cscr$ generated by the dualizable objects is denoted $\Cscr_{\operatorname{f.d.}}$.\index[notation]{Cfd@$\Cscr_{\operatorname{f.d.}}$}
\end{definition}
There is, of course, a well-established and well-studied notion of fully dualizable objects which is important in representation theory and the theory of TQFTs. We will not use this more general higher-categorical notion in this monograph but instead refer the reader to \cite[Appendix A]{Chris2Noah} and \cite{LurieCobord}. The connection to the $1$-categorical notion used in this monograph lies in the fact that every $1$-morphism of a symmetric monoidal $(\infty,1)$-category $\mathsf{C}$ is dualizable in the sense that every $1$-morphism in the homotopy bicategory $h_2^{(0)}\mathsf{C}$ admits a left/right adjoint and a $1$-morphism in the 
homotopy bicategory $h_2^{(-1)}\mathsf{C}$ has left/right adjoints if and only if the corresponding object in $h\mathsf{C}$ is dualizable in the sense of Definition \ref{Defn: Section SFF: dualizable smc object}.

\begin{definition}\label{Defn: Section SFF: Trace class preeq}
	Let $F = (F,\overline{F})$ be a pre-equivariant pseudofunctor on a left $G$-variety $X$. We say that $F$ is a trace-class pre-equivariant\index[terminology]{Pre-equivariant Pseudofunctor! Trace-class} pseudofunctor if:
	\begin{itemize}
		\item $F$ is herbaceous;
		\item The categories $\overline{F}(X)$ and $\overline{F}(\quot{X}{G})$ have objects which may be evaluated on opens of $X$ for which each object $A$ takes values in a category $\Cscr_{\operatorname{f.d.}}$, where $\Cscr_{\operatorname{f.d.}}$ is the subcategory of fully dualizable objects in a symmetric monoidal category $\Cscr$.
	\end{itemize}
\end{definition}
\begin{example}
The pre-equivariant pseudofunctors $\Shv(-;\Z_{\ell}\mathbf{PrMod}_{\operatorname{f.g.}})$ of finitely generated projective $\ell$-adic sheaves and $\DbQl{-},$ the bounded derived category of $\overline{\Q}_{\ell}$-sheaves are trace-class. Additional examples can be found in the categories $\Shv(-;\mathbf{Vect}_{\operatorname{f.d.}})$ of sheaves of finite-dimensional vector spaces and bounded derived complexes of perfect chain complexes $D_c^b(-)_{\operatorname{perf}}$. A more silly example can be found by taking the pre-equivariant pseudofunctor $F$ induced by taking $\overline{F}$ to be the constant pseudofunctor on the terminal category $\One$.
\end{example}
\begin{example}
	For a non-example of a trace-class pre-equivariant pseudofunctor let $F$ be the pre-equivariant pseudofunctor induced by the pseudofunctor $\overline{F}$ defined by object assignments $\overline{F}(Z) := \KVect$ for a field $K$ and taking translation functors $\overline{F}(f)$ to be given by the identity functor on $\KVect$. Then $F$ is not trace-class because the category of fully dualizable objects in $\KVect$, $\KVect_{\operatorname{f.d}}$, is the category of finite dimensional vectors
\end{example}

\begin{proposition}
	If $F$ is a trace-class pre-equivariant pseudofunctor then for any object $A$ of $F_G(X)$, any endomorphism $f:A \to A$, and any point $x$ of $X$ the stalks $A_x$ and $f_x$ admit traces.
\end{proposition}
\begin{proof}
	This is immediate from the definitions with $\trace(A_x)$ being given by
	\[
	\begin{tikzcd}
		\Ibb \ar[r]{}{\eta_{A_x}} & A_x^{\ast} \otimes A_x \ar[r]{}{\gamma_{A_x^{\ast},A_x}} & A_x \otimes A_x^{\ast} \ar[r]{}{\epsilon_{A_x}} & \Ibb
	\end{tikzcd}
	\]
	and $\trace(f_x)$ given by
	\[
	\begin{tikzcd}
		\Ibb \ar[r]{}{\eta_{A_x}} & A_x^{\ast} \otimes A_x \ar[rr]{}{\id_{A_x^{\ast}} \otimes f_x} & & A_x^{\ast} \otimes A_x \ar[r]{}{\gamma_{A_x^{\ast},A_x}} & A_x \otimes A_x^{\ast} \ar[r]{}{\epsilon_{A_x}} & \Ibb
	\end{tikzcd}
	\]
	where $\Ibb$ is the monoidal unit of $\overline{F}(X)$ and $\gamma_{B,C}:B \otimes C \to C \otimes B$ is the braiding map.
\end{proof}

We now are in position to define the trace of an equivariant category. Let $F$ be a trace-class pre-equivariant pseudofunctor and let $g$ be a point of $G$ and $x$ be a point of $X$ for which $x$ is fixed by $g$, i.e., $gx = x$. Then on one hand by virtue of $F$ being herbaceous there is an isomorphism of stalks
\[
\pi_2^{\ast}(A)_{(g,x)} \cong A_{\pi_2(g,x)} = A_{x}
\]
while on the other hand
\[
\alpha_X^{\ast}(A)_{(g,x)} \cong A_{\alpha_X(g,x)} = A_{gx} = A_x.
\]
Consequently the map 
\[
\quot{\theta}{G}_{(g,x)}^{A}:A_x \to A_x
\] 
gives an automorphism of $A_x$ which we use to define the equivariant trace of $A$ at the point $(g,x)$.

\begin{definition}\label{Defn: Section SFF: Equivariant Trace}
	Let $F = (F,\overline{F})$ be a trace-class pre-equivariant pseudofunctor on $X$ and let $g$ be a point of $G$ and $x$ a point of $X$ fixed by $g$. Then if $A \in F_G(X)_0$ we define the trace\index[terminology]{Pre-equivariant Pseudofunctor! Trace} at $g$ and $x$ of $A$ to be given by
	\[
	\trace_g(A,x) := \trace\left(\quot{\theta_{(g,x)}^{A}}{G}\right).
	\]
\end{definition}

Our first goal is to prove that when $F$ is a triangulated trace-class pre-equivariant pseudofunctor with right exact translation functors, trace is additive on distinguished triangles in the sense that if
\[
\begin{tikzcd}
	A \ar[r]{}{P} & B \ar[r]{}{\Phi} & C \ar[r]{}{\Psi} & A[1]
\end{tikzcd}
\]
is a distinguished triangle in $F_G(X)$ then
\[
\trace_g(B,x) = \trace_g(A,x) + \trace_g(C,x).
\]
To go about this proof we will need one technical lemma that allows us to see, in the presence of a regular epimorphism $e:B \to C$, the equivariance $\Theta^C$ is functionally determined by the ``quotient'' of $\Theta^B$ by $e$. We will also collect some folkloric facts about triangulated categories which are well-known to experts solely for the sake of completeness.

\begin{lemma}\label{Lemma: Section SFF: Equivariances and regular epis}
	Let $F$ be a pre-equivariant pseudofunctor for which the functor $\alpha_X^{\ast}:F_G(X) \to F_G(G \times X)$ is right exact. Then for any regular epimorphism $e:B \to C$ in $F_G(X)$, the equivariance $\Theta^C$ is a quotient of $\Theta^B$ in the sense that it is induced by $\Theta^B$ and the universal property of an epimorphism.
\end{lemma}
\begin{proof}
	Because we have that $\alpha_X^{\ast}e$ is a regular epimorphism by assumption and because the equivariances $\Theta$ are natural in $F_G(X)$, the diagram
	\[
	\begin{tikzcd}
		\alpha_X^{\ast}B \ar[d, swap]{}{\Theta^B} \ar[r]{}{\alpha_X^{\ast}e} & \alpha_X^{\ast}C \ar[d]{}{\Theta^C} \\
		\pi_2^{\ast}B \ar[r, swap]{}{\pi_2^{\ast}e} & \pi_2^{\ast}C
	\end{tikzcd}
	\]
	commutes with $\alpha_X^{\ast}e$ a regular epimorphism. Then it follows that $\Theta^C$ is the unique map making
	\[
	\begin{tikzcd}
		& \alpha_X^{\ast}C \ar[dr]{}{\Theta^C} \\
		\alpha_X^{\ast}B \ar[ur]{}{\alpha_X^{\ast}C} \ar[rr, swap]{}{\pi_2^{\ast}e \circ \Theta^B} & & \pi_2^{\ast}C
	\end{tikzcd}
	\]
	commute, which in turn proves the lemma.
\end{proof}
\begin{lemma}[{\cite{GregStevensonSplitMonic}}]\label{Lemma: Section SFF: Epics and Monics split in trinagulated cats}
	Let $\Cscr$ be a triangulated category (cf.\@ Definition \ref{Defn: Section Triangle: Traingulated Categories}). Then $\Cscr$ is a semisimple category, i.e., every epic and monic split in $\Cscr$. 
\end{lemma}
\begin{proof}
	We present the argument for monics, as the example for epics follows dually. Fix a monic $\mu:A \to B$ and complete it to a distinguished triangle:
	\[
	\begin{tikzcd}
		A \ar[r]{}{\mu} & B \ar[r]{}{\varphi} & C \ar[r]{}{\psi} & A[1]
	\end{tikzcd}
	\]
	Rotating the triangle to the left gives a distinguished triangle
	\[
	\begin{tikzcd}
		C[-1] \ar[r]{}{\psi[-1]} & A \ar[r]{}{\mu} & B \ar[r]{}{\varphi} & C 
	\end{tikzcd}
	\]
	which, because adjacent morphisms in a triangle compose to zero give that $\mu \circ \psi[-1] = 0 = \mu \circ 0$ and hence show that $\psi = 0[1] = 0$. Thus it follows that the triangle is isomorphic to
	\[
	\begin{tikzcd}
		A \ar[r] & A \oplus C \ar[r] & C \ar[r]{}{0} & A[1]
	\end{tikzcd}
	\]
	by \cite[Corollary 1.2.7]{NeemanTriangulated}, giving that $\mu$ splits.
\end{proof}
\begin{lemma}\label{Lemma: Section SFF: Additive functors between triangulated cats}
	If $\Cscr$ and $\Dscr$ are triangulated categories then any additive functor $F:\Cscr \to \Dscr$ is exact.
\end{lemma}
\begin{proof}
	Because additive functors of additive categories being exact is equivalent to preserving monics and epics, and because every epic and monic in $\Cscr$ splits by Lemma \ref{Lemma: Section SFF: Epics and Monics split in trinagulated cats}, the result follows because every functor preserves split epics and split monics.
\end{proof}
\begin{Theorem}\label{Thm: Section SFF: Additive trace}
Let $F$ be a triangulated trace-class pre-equivariant pseudofunctor on $X$ and $G \times X$ for which:
\begin{itemize}
	\item The category $\Cscr$ of values the objects $A$ of $\overline{F}(\quot{X}{G})$, $\overline{F}(X)$, $\overline{F}(\quot{(G \times X)}{G})$ and $\overline{F}(G \times X)$ take is triangulated and arises as the ``derived category'' of a ``homotopy category'' in the sense of \cite[Section 5]{MayTrace}.
	\item The operation of taking stalks gives triangulated exact functors $(-)_{(g,x)}:F_G(G \times X) \to \Cscr$ and $(-)_{x}:F_G(X) \to \Cscr$.
	\item Each functor 
	\[
	\overline{F}\left(\overline{\id_{\Gamma} \times \alpha_X}\right):\overline{F}(\XGamma) \to \overline{F}\left(\quot{(G \times X)}{\Gamma}\right)
	\]
	(cf.\@ the proof of Lemma \ref{Lemma: Section Change of Space: Action map pullback}) is additive. 
\end{itemize}
Then equivariant trace is additive on distinguished triangles in the sense that if
\[
\begin{tikzcd}
	A \ar[r]{}{P} & B \ar[r]{}{\Phi} & C \ar[r]{}{\Psi} & A[1]
\end{tikzcd}
\]
is a distinguished triangle for objects $A, B, C$ in $F_G(X)$ then
\[
\trace_g(B,x) = \trace_g(A,x) + \trace_g(C,x).
\]
\end{Theorem}
\begin{proof}
	Note that the assumption that each functor $\overline{F}\left(\overline{\id_{\Gamma} \times \alpha_X}\right)$ is exact implies, because $\overline{F}(\XGamma)$ and $\overline{F}(\quot{(G \times X)}{\Gamma})$ are triangulated, that each corresponding functor is exact by Lemma \ref{Lemma: Section SFF: Additive functors between triangulated cats} and hence so is $\alpha_X^{\ast}:F_G(X) \to F_G(G \times X)$ by Corollary \ref{Cor: Pseudocone Functors: When translations preserve all limits and colimits}. The functor $\pi_2^{\ast}:F_G(X) \to F_G(G \times X)$ and the functors $\quot{\overline{\pi}_2}{\Gamma}:\overline{F}(\XGamma) \to \overline{F}\left(\quot{(G \times X)}{\Gamma}\right)$ are exact because $F$ is a triangulated pre-equivariant pseudofunctor and hence are additive and exact for similar reasons. 
	
	Assume that
	\[
	\begin{tikzcd}
		A \ar[r]{}{P} & B \ar[r]{}{\Phi} & C \ar[r]{}{\Psi} & A[1]
	\end{tikzcd}
	\]
	is a distinguished triangle in $F_G(X)$. Applying the functors $\pi_2^{\ast}$ and $\alpha_X^{\ast}$ gives us that the triangles
	\[
	\begin{tikzcd}
		\alpha_X^{\ast}A \ar[r]{}{\alpha_X^{\ast}P} & \alpha_X^{\ast}B \ar[r]{}{\alpha_X^{\ast}\Phi} & \alpha_X^{\ast}C \ar[r]{}{\alpha_X^{\ast}\Psi} & \alpha_X^{\ast}(A)[1]
	\end{tikzcd}
	\]
	and
	\[
	\begin{tikzcd}
		\pi_2^{\ast}A \ar[r]{}{\pi_2^{\ast}P} & \pi_2^{\ast}B \ar[r]{}{\pi_2^{\ast}\Phi} & \pi_2^{\ast}C \ar[r]{}{\pi_2^{\ast}\Psi} & \pi_2^{\ast}(A)[1]
	\end{tikzcd}
	\]
	are distinguished in $F_G(G \times X)$. The naturality of $\Theta$ proved in Theorem \ref{Theorem: Formal naive equivariance} gives rise to a commuting diagram
	\[
	\begin{tikzcd}
		\alpha_X^{\ast}A \ar[r]{}{\alpha_X^{\ast}P} \ar[d, swap]{}{\Theta^A} & \alpha_X^{\ast}B \ar[r]{}{\alpha_X^{\ast}\Phi} \ar[d]{}[description]{\Theta^B} & \alpha_X^{\ast}C \ar[d]{}[description]{\Theta^C} \ar[r]{}{\alpha_X^{\ast}\Psi} & \alpha_X^{\ast}(A)[1] \ar[d]{}{\Theta^{A[1]}} \\
		\pi_2^{\ast}A \ar[r, swap]{}{\pi_2^{\ast}P} & \pi_2^{\ast}B \ar[r, swap]{}{\pi_2^{\ast}\Phi} & \pi_2^{\ast}C \ar[r, swap]{}{\pi_2^{\ast}\Psi} & \pi_2^{\ast}(A)[1]
	\end{tikzcd}
	\]
	in $F_G(G \times X)$. Taking the $(G \times X)$-component of the diagram above and using the definition of being distinguished in $F_G(G \times X)$ (cf.\@ Theorem \ref{Theorem: Section Triangle: Equivariant triangulation}) produces the commuting diagram
	\[
	\begin{tikzcd}
		\overline{F}\left(\quot{\overline{\alpha}_X}{G}\right)\left(\quot{A}{G }\right) \ar[d, swap]{}{\overline{F}\left(\quot{\overline{\alpha}_X}{G}\right)\left(\quot{\rho}{G}\right)} \ar[rr]{}{\quot{\theta^A}{G}} & & \overline{F}\left(\quot{\overline{\pi}_2}{G}\right)\left(\quot{A}{G}\right) \ar[d]{}{\overline{F}\left(\quot{\overline{\pi}_2}{G}\right)\left(\quot{\rho}{G}\right)} \\
		\overline{F}\left(\quot{\overline{\alpha}_X}{G}\right)\left(\quot{B}{G}\right) \ar[d, swap]{}{\overline{F}\left(\quot{\overline{\alpha}_X}{G}\right)\left(\quot{\varphi}{G}\right)}	\ar[rr]{}[description]{\quot{\theta^B}{G}} & & \overline{F}\left(\quot{\overline{\pi}_2}{G}\right)\left(\quot{B}{G}\right) \ar[d]{}{\overline{F}\left(\quot{\overline{\pi}_2}{G}\right)\left(\quot{\varphi}{G}\right)} \\
		\overline{F}\left(\quot{\overline{\alpha}_X}{G}\right)\left(\quot{C}{G}\right) \ar[rr]{}[description]{\quot{\theta^C}{G}} \ar[d, swap]{}{\overline{F}\left(\quot{\overline{\alpha}_X}{G}\right)\left(\quot{\psi}{G}\right)} & & \overline{F}\left(\quot{\overline{\pi}_2}{G}\right)\left(\quot{C}{G}\right) \ar[d]{}{\overline{F}\left(\quot{\overline{\pi}_2}{G}\right)\left(\quot{\psi}{G}\right)} \\
		\overline{F}\left(\quot{\overline{\alpha}_X}{G}\right)\left(\quot{A}{G}\right)[1] \ar[rr, swap]{}{\quot{\theta^{A[1]}}{G}} & & \overline{F}\left(\quot{\overline{\pi}_2}{G}\right)\left(\quot{A}{G}\right)[1]
	\end{tikzcd}
	\]
	for which both columns are distinguished triangles. Now because $\overline{F}$ is triangulated, the equivalence of categories $\overline{F}(\psi_{G \times X}):\overline{F}(\quot{(G \times X)}{G}) \to \overline{F}(G \times X)$ is triangulated and hence sends distinguished triangles to distinguished triangles. As such, applying the functor $\overline{F}\left(\quot{\psi}{G \times G}\right)$ gives the commuting diagram
	\[
	\begin{tikzcd}
		\overline{F}\left(\quot{\psi}{G \times G}\right)\left(\overline{F}\left(\quot{\overline{\alpha}_X}{G}\right)\left(\quot{A}{G}\right)\right) \ar[dd]{}[description]{\overline{F}\left(\quot{\psi}{G \times G}\right)\left(\overline{F}\left(\quot{\overline{\alpha}_X}{G}\right)\left(\quot{\rho}{G}\right)\right)} \ar[rr]{}{\overline{F}\left(\quot{\psi}{G \times G}\right)\left(\quot{\theta^A}{G}\right)} & & \overline{F}\left(\quot{\psi}{G \times G}\right)\left(\overline{F}\left(\quot{\overline{\pi}_2}{G}\right)\left(\quot{A}{G}\right)\right) \ar[dd]{}[description]{\overline{F}\left(\quot{\psi}{G \times G}\right)\left(\overline{F}\left(\quot{\overline{\pi}_2}{G}\right)\left(\quot{\rho}{G}\right)\right)} \\
		\\
		\overline{F}\left(\quot{\psi}{G \times G}\right)\left(\overline{F}\left(\quot{\overline{\alpha}_X}{G}\right)\left(\quot{B}{G}\right)\right) \ar[dd]{}[description]{\overline{F}\left(\quot{\psi}{G \times G}\right)\left(\overline{F}\left(\quot{\overline{\alpha}_X}{G}\right)\left(\quot{\varphi}{G}\right)\right)}	\ar[rr]{}{\overline{F}\left(\quot{\psi}{G \times G}\right)\left(\quot{\theta^B}{G}\right)} & & \overline{F}\left(\quot{\psi}{G \times G}\right)\left(\overline{F}\left(\quot{\overline{\pi}_2}{G}\right)\left(\quot{B}{G}\right)\right) \ar[dd]{}[description]{\overline{F}\left(\quot{\psi}{G \times G}\right)\left(\overline{F}\left(\quot{\overline{\pi}_2}{G}\right)\left(\quot{\varphi}{G}\right)\right)} \\
		\\
		\overline{F}\left(\quot{\psi}{G \times G}\right)\left(\overline{F}\left(\quot{\overline{\alpha}_X}{G}\right)\left(\quot{C}{G}\right)\right) \ar[rr]{}{\overline{F}\left(\quot{\psi}{G \times G}\right)\left(\quot{\theta^C}{G}\right)} \ar[dd]{}[description]{\overline{F}\left(\quot{\psi}{G \times G}\right)\left(\overline{F}\left(\quot{\overline{\alpha}_X}{G}\right)\left(\quot{\psi}{G}\right)\right)} & & \overline{F}\left(\quot{\psi}{G \times G}\right)\left(\overline{F}\left(\quot{\overline{\pi}_2}{G}\right)\left(\quot{C}{G}\right)\right) \ar[dd]{}[description]{\overline{F}\left(\quot{\psi}{G \times G}\right)\left(\overline{F}\left(\quot{\overline{\pi}_2}{G}\right)\left(\quot{\psi}{G}\right)\right)} \\
		\\
		\overline{F}\left(\quot{\psi}{G \times G}\right)\left(\overline{F}\left(\quot{\overline{\alpha}_X}{G}\right)\left(\quot{A}{G}\right)\right)[1] \ar[rr, swap]{}{\overline{F}\left(\quot{\psi}{G \times G}\right)\left(\quot{\theta^{A[1]}}{G}\right)} & & \overline{F}\left(\quot{\psi}{G \times G}\right)\left(\overline{F}\left(\quot{\overline{\pi}_2}{G}\right)\left(\quot{A}{G}\right)\right)[1]
	\end{tikzcd}
	\]
	with each column a distinguished triangle. Taking stalks of the diagram at the pair $(g,x)$ where $g$ is a point of $G$ and $x$ is a $g$-fixed point of $X$ then gives a commuting diagram which is equal to
	\[
	\begin{tikzcd}
		\alpha_X^{\ast}(A)_{(g,x)} \ar[d, swap]{}{\Theta_{(g,x)}^{A}} \ar[r]{}{\alpha_X^{\ast}(P)_{(g,x)}} & \alpha_X^{\ast}(B)_{(g,x)} \ar[d, swap]{}{\Theta_{(g,x)}^{B}} \ar[r]{}{\alpha_X^{\ast}(\Phi)_{(g,x)}} & \alpha_X^{\ast}(C)_{(g,x)} \ar[d]{}{\Theta_{(g,x)}^{C}} \ar[r]{}{\alpha_X^{\ast}(\Psi)_{(g,x)}} & \alpha_X^{\ast}(A)[1]_{(g,x)} \ar[d]{}{\Theta_{(g,x)}^{A[1]}} \\
		\pi_2^{\ast}(A)_{(g,x)} \ar[r, swap]{}{\pi_2^{\ast}(P)_{(g,x)}} & \pi_2^{\ast}(B)_{(g,x)} \ar[r, swap]{}{\pi_2^{\ast}(\Phi)_{(g,x)}} & \pi_2^{\ast}(C)_{(g,x)} \ar[r, swap]{}{\pi_2^{\ast}(\Psi)_{(g,x)}} & \pi_2^{\ast}(A)[1]_{(g,x)}
	\end{tikzcd}
	\]
	where each row is a distinguished triangle. Now by virtue of $F$ being herbaceous we can rewrite this to the diagram
	\[
	\begin{tikzcd}
		A_x \ar[d, swap]{}{\Theta_{(g,x)}^{A}} \ar[r]{}{\rho_x} & B_x \ar[d, swap]{}{\Theta_{(g,x)}^{B}} \ar[r]{}{\varphi_x} & C_x \ar[d]{}{\Theta_{(g,x)}^{C}} \ar[r]{}{\psi_x} & A[1]_x \ar[d]{}{\Theta_{(g,x)}^{A[1]}} \\
		A_x \ar[r, swap]{}{\rho_x} & B_x \ar[r, swap]{}{\varphi_x} & C_x \ar[r, swap]{}{\psi_x} & A[1]_x
	\end{tikzcd}
	\]
	in which each row is a distinguished triangle. Using now that $\Cscr$ is the category of fully dualizable objects in a triangulated symmetric monoidal category allow us to invoke \cite[Theorem 1.9]{MayTrace} and produce a map $\omega$ for which the diagram
	\[
	\begin{tikzcd}
		A_x \ar[d, swap]{}{\Theta_{(g,x)}^{A}} \ar[r]{}{\rho_x} & B_x \ar[d, swap]{}{\Theta_{(g,x)}^{B}} \ar[r]{}{\varphi_x} & C_x \ar[d]{}{\Theta_{(g,x)}^{C}} \ar[r]{}{\psi_x} & A[1]_x \ar[d]{}{\omega} \\
		A_x \ar[r, swap]{}{\rho_x} & B_x \ar[r, swap]{}{\varphi_x} & C_x \ar[r, swap]{}{\psi_x} & A[1]_x
	\end{tikzcd}
	\]
	commutes and the additivity relation
	\[
	\trace\left(\Theta^B_{(g,x)}\right) = \trace\left(\Theta^A_{(g,x)}\right) + \trace(\omega)
	\]
	holds. However, a reading of the proof of \cite[Theorem 1.9]{MayTrace} (explicitly \cite[Section 8]{MayTrace}) shows that $\omega$ is given in essence as a homotopy class of quotients. Lemma \ref{Lemma: Section SFF: Equivariances and regular epis} shows that the equivariances $\Theta^{-}$ respect these quotient constructions and hence allows us to choose $\omega = \Theta^C_{(g,x)}$. In this way we get that
	\begin{align*}
		\trace_g(B,x) &= \trace\left(\Theta_{(g,x)}^{B}\right) = \trace\left(\Theta_{(g,x)}^{A}\right) + \trace\left(\Theta_{(g,x)}^{C}\right) \\
		&= \trace_g(A,x) + \trace_g(C,x),
	\end{align*}
	as was desired.
\end{proof}
\begin{corollary}
	Equivariant trace is additive in any of the cases below:
	\begin{itemize}
		\item $G$ is a smooth algebraic group, $X$ is a left $G$-variety, and $F$ is either the equivariant derived pre-equivariant pseudofunctor or the $\ell$-adic equivariant derived pre-equviariant pseudofunctor;
		\item $G$ is a topological group which admits $n$-acyclic free $G$-spaces which are manifolds, $X$ is a left $G$-space, and $F$ is the equivariant derived category pre-equivariant pseudofunctor.
	\end{itemize}
\end{corollary}
\begin{proof}
	The proof ultimately comes down to the fact that, in any of the cases listed, because pullback functors are right exact\footnote{As pullback functors are adjoints $f^{\ast} \dashv f_{\ast}$, giving that $f^{\ast}$ preserves colimits because it is a left adjoint.}, each functor $(\overline{\id_{\Gamma} \times \alpha_X})^{\ast}$ is right exact; consequently $\alpha_X^{\ast}$ is right exact by Corollary \ref{Cor: Pseudocone Functors: When translations preserve all limits and colimits} and so is exact by Lemma \ref{Lemma: Section SFF: Additive functors between triangulated cats}. From here the corollary follows by a direct application of Theorem \ref{Thm: Section SFF: Additive trace}.
\end{proof}

We now conclude the monograph by showing that there are two ``continuity'' properties of equivariant trace which apply to strict pre-equivariant pseudofunctors: first, if $h:X \to Y$ is an equivariant map and if $g$ is a point of $G$ with $x$ a point of $X$ fixed by $g$, then for any object $A$ of $F_G(Y)$,
\[
\trace_g\left(h^{\ast}A,x\right) = \trace_g\left(A,h(x)\right).
\]
Second, we show that if $i:X \to Y$ is a $G$-equivariant immersion and if $F$ admits a pushforward functor along $i$ (right adjoint to pullback against $i$) with evaluation on opens calculated by
\[
h_{\ast}A(\quot{U}{G}) = A\left(X \times_Y \quot{U}{G}\right)
\]
then for any object $A$ of $F_G(X)$ there is an equality of traces
\[
\trace_g(A,x) = \trace_g\left(i_{\ast}A,i(x)\right).
\]
The assumptions regarding the strictness are not too harsh in practice. Many practical models of the equivariant derived category which appear in the literature (cf.\@ those used in \cite{PramodBook}, \cite{BernLun}, \cite{CFMMX}, \cite{LusztigCuspidal2}, \cite{MirkVilDuality}, \cite[Chapters 6 -- 9]{MyThesis}, \cite{Geordie}) are those induced from the corresponding strict derived category pseudofunctor (or its $\ell$-adic counterpart).

\begin{proposition}
Let $F = (F, \overline{F})$ be a strict trace-class pre-equivariant pseudofunctor on $X$ and let $h:X \to Y$ be a $G$-equivariant morphism. Then for any object $A$ of $F_G(Y)$ and any point $g$ of $G$ with $x$ a point of $X$ fixed by $g$, we have that
\[
\trace_g(h^{\ast}A, x) = \trace_g\left(A,h(x)\right)
\]
where $h^{\ast}$ is induced from $\overline{F}$ by Theorem \ref{Thm: Pseudocone Functors: Pullback induced by fibre functors in pseudofucntor}.
\end{proposition}
\begin{proof}
Recall from Corollary \ref{Cor: Section Equivariance: Pullbacks and Equivariances for Strict pseudofunctors} we have that
\[
\Theta_{X}^{h^{\ast}A} = (\id_G \times h)^{\ast}\Theta_Y^{A}
\]
so it follows that in this case we have that
\[
\left(\quot{\theta_X^{h^{\ast}A}}{G}\right)_{(g,x)} = \left(\overline{F}\left(\quot{\overline{h}}{G \times G}\right)\left(\quot{\theta_Y^{A}}{G}\right)\right)_{(g,x)}.
\]
Now, writing $(g,x) \to U$ for opens $U$ to mean that $U$ is an open neighbourhood in the relevant (Grothendieck) topology around the point $(g,x)$, we find that
\begin{align*}
&\trace_g(h^{\ast}A,x) \\
&= \trace\left(\quot{\theta_{(g,x)}^{h^{\ast}A}}{G}\right) = \trace\left(\quot{(h^{\ast}\theta^A_Y)_{(g,x)}}{G}\right) = \trace\left(\left(\quot{\overline{h}}{G \times G}^{\ast}(\quot{\theta^{A}_Y}{G})\right)_{(g,x)}\right) \\
&= \trace\left(\colim_{\substack{(g,x) \to W \\ W \to G \times X\, \text{open}}}\left(\overline{F}\left(\quot{\overline{h}}{G \times G}\right)\left(\quot{\theta^A}{G}\right)\right)(\quot{W}{G})\right) \\
&=\trace\left(\colim_{\substack{g \to V, x \to W^{\prime} \\ V \to G,  W^{\prime} \to X\, \text{open}}}\left(\overline{F}\left(\quot{\overline{h}}{G \times G}\right)\left(\quot{\theta^A}{G}\right)\right)\left(\quot{(V \times W^{\prime})}{G}\right)\right) \\
&=\trace\left(\colim_{\substack{g \to V, x \to W^{\prime} \\ V \to G,  W^{\prime} \to X\, \text{open}}}\left(\colim_{\substack{\overline{F}\left(\quot{\overline{h}}{G \times G}\right)(\quot{(V \times W^{\prime})}{G}) \to \quot{U^{\prime\prime}}{G} \\ U^{\prime\prime} \to G \times Y\, \text{open}}}\left(\quot{\theta^A}{G}\right)(\quot{ U^{\prime\prime}}{G})\right)\right) \\
&=\trace\left(\colim_{\substack{g \to V, x \to W^{\prime} \\ V \to G,  W^{\prime} \to X\, \text{open}}}\left(\colim_{\substack{\overline{F}\left(\quot{\overline{h}}{G \times G}\right)(\quot{(V \times W^{\prime})}{G}) \to \quot{(V \times U^{\prime})}{G} \\ U^{\prime} \to Y\, \text{open}}}\left(\quot{\theta^A}{G}\right)(\quot{(V \times U^{\prime})}{G})\right)\right) \\
&= \trace\left(\colim_{\substack{g \to V, h(x) \to U^{\prime}\\ V \to G,  U^{\prime} \to Y\, \text{open}}}\left(\quot{\theta^A}{G}\right)(\quot{(V \times U^{\prime})}{G})\right) \\
&= \trace\left(\colim_{\substack{(g,h(x)) \to U \\ U \to G \times X\, \text{open}}} \quot{\theta^A}{G}(\quot{U}{G})\right) = \trace\left(\quot{\theta^A_{(g,h(x))}}{G}\right) = \trace_g(A,h(x))
\end{align*}
which was what was desired.
\end{proof}

Our goal now is to show that trace preserves base points for immersions $i:X \to Y$ in the sense that
\[
\trace_g(A,x) = \trace_g\left(i_{\ast}A,i(x)\right)
\]
in mild conditions. To go about this, we now assume that there is a right adjoint
\[
\begin{tikzcd}
F_G(X) \ar[rr, bend right = 20, swap, ""{name = D}]{}{i_{\ast}} & & F_G(Y) \ar[ll, bend right = 20, swap, ""{name = U}]{}{i^{\ast}} \ar[from = U, to = D, symbol = \dashv]
\end{tikzcd}
\]
induced by adjunctions
\[
\begin{tikzcd}
	\overline{F}(\quot{X}{\Gamma}) \ar[rr, bend right = 20, swap, ""{name = D}]{}{\quot{\overline{\imath}_{\ast}}{\Gamma}} & & \overline{F}(\quot{Y}{\Gamma}) \ar[ll, bend right = 20, swap, ""{name = U}]{}{\overline{F}(\quot{\overline{\imath}}{\Gamma})} \ar[from = U, to = D, symbol = \dashv]
\end{tikzcd}
\]
for all $\Gamma \in \Sf(G)_0$. We also assume that evaluation of the objects of $i_{\ast}A$ at opens of $Y$ takes the form
\[
(i_{\ast}A)(U) = A(U \times_Y X).
\]

As a first sanity check we show that for immersions $i:X \to Y$, stalks are stable under pushforward. 
In this case we first note that if $A$ is an object of $F_G(X)$ then for any point $x$ of $X$,
\[
(i_{\ast}A)_{i(x)} = \colim_{\substack{i(x) \to U \\ U \to Y\,\text{open} }} \left(\left(\quot{\overline{\imath}_{\ast}}{G}\right)\left(\quot{A}{G}\right)\right)\left(\quot{U}{G}\right) = \colim_{\substack{i(x) \to U \\ U \to Y\,\text{open} }}\quot{A}{G}\left(\quot{X}{G} \times_{\quot{Y}{G}} \quot{U}{G}\right).
\]
Now asking that we have a map $i(x) \to U$ is equivalent to asking for a map $x \to U \times_Y X$, so we can rewrite the colimit as
\[
\colim_{\substack{i(x) \to U \\ U \to Y\,\text{open} }}\quot{A}{G}\left(\quot{X}{G} \times_{\quot{Y}{G}} \quot{U}{G}\right) = \colim_{\substack{x \to X \times_Y U \\ U \to Y\,\text{open} }}\quot{A}{G}\left(\quot{X}{G} \times_{\quot{Y}{G}} \quot{U}{G}\right)
\]
However, since the system of opens $\quot{X}{G} \times_{\quot{Y}{G}} \quot{U}{G}$ is cofinal in the class of opens of $\quot{X}{G}$ because $i$ is an immersion, we can further rewrite the colimit as
\[
\colim_{\substack{x \to X \times_Y U \\ U \to Y\,\text{open} }}\quot{A}{G}\left(\quot{X}{G} \times_{\quot{Y}{G}} \quot{U}{G}\right) = \colim_{\substack{x \to V \\ V \to X\,\text{open} }}\quot{A}{G}\left(\quot{V}{G}\right) = A_x.
\]
Putting this all together gives that
\begin{align*}
(i_{\ast}A)_{i(x)} &= \colim_{\substack{i(x) \to U \\ U \to Y\,\text{open} }} \left(\left(\quot{\overline{\imath}_{\ast}}{G}\right)\left(\quot{A}{G}\right)\right)\left(\quot{U}{G}\right) \\
&= \colim_{\substack{x \to V \\ V \to X\,\text{open} }}\quot{A}{G}\left(\quot{V}{G}\right) = A_x
\end{align*}
as well as our next proposition.
\begin{proposition}
Let $i:X \to Y$ be a $G$-equivariant immersion and let $F = (F,\overline{F})$ be a trace-class pre-equivariant pseudofunctor on $X$ for which $\overline{F}$ is strict and for which pushforwards $i_{\ast}A$ are evaluated by the equation
\[
\left(\left(\quot{\overline{\imath}}{G}\right)\left(\quot{A}{G}\right)\right)\left(U\right) = \quot{A}{G}\left(\quot{X}{G} \times_{\quot{Y}{G}} U\right).
\]
Then for any point $x$ of $X$ and any object $A$ of $F_G(X)$,
\[
(i_{\ast}A)_{i(x)} = A_{x}.
\]
\end{proposition}
 
We can now use the techniques that appear in this derivation to establish that immersions preserve traces.
\begin{proposition}\label{Prop: Section SFF: Traces of pushforwards}
Let $i:X \to Y$ be a $G$-equivariant immersion and let $F = (F,\overline{F})$ be a trace-class pre-equivariant pseudofunctor on $X$ for which $\overline{F}$ is strict and for which pushforwards $i_{\ast}A$ are evaluated by the equation
\[
\left(\left(\quot{\overline{\imath}}{G}\right)\left(\quot{A}{G}\right)\right)\left(U\right) = \quot{A}{G}\left(\quot{X}{G} \times_{\quot{Y}{G}} U\right).
\] 
Then for any point $g$ of $G$ and any $g$-fixed point $x$ of $X$, for any object $A$ of $F_G(X)$ we have that
\[
\quot{\theta^{i_{\ast}A}_{(g,i(x))}}{G} = \quot{\theta^A_{(g,x)}}{G}
\]
and
\[
\trace_g(i_{\ast}A,i(x)) = \trace_g(A,x).
\]
\end{proposition}
\begin{proof}
Note that if we can establish the first claim the second follows immediately because
\[
\trace_g(i_{\ast}A,i(x)) = \trace\left(\quot{\theta^{i_{\ast}A}_{(g,i(x))}}{G}\right) = \trace\left(\quot{\theta^{A}_{(g,x)}}{G}\right) = \trace_g(A,x).
\]
To establish the first claim we make the following calculation:
\begin{align*}
\quot{\theta^{i_{\ast}A}_{(g,i(x))}}{G} &=\colim_{\substack{(g,i(x)) \to U \\ U \to G \times Y\,\text{open}}} \quot{\theta_Y^{i_{\ast}A}}{G}\left(\quot{U}{G}\right) = \colim_{\substack{g \to V, i(x) \to W \\ V \to G, W \to Y\,\text{open}}} \quot{\theta_Y^{i_{\ast}A}}{G}\left(\quot{(V \times W)}{G}\right) \\
&= \colim_{\substack{g \to V, i(x) \to W \\ V \to G, W \to Y\,\text{open}}}\left(\left(\quot{\overline{\imath}}{G}\right)_{\ast}\left(\quot{\theta^{A}_{X}}{G}\right)\right)\left(\quot{(V \times W)}{G}\right) \\
&= \colim_{\substack{g \to V, i(x) \to W \\ V \to G, W \to Y\,\text{open}}}\quot{\theta^{A}_{X}}{G}\left(\quot{\left(V \times \left(X \times_Y W\right)\right)}{G}\right) \\
&=\colim_{\substack{g \to V, x \to W^{\prime} \\ V \to G, W^{\prime} \to X\,\text{open}}}\quot{\theta^{A}_{X}}{G}\left(\quot{\left(V \times W^{\prime}\right)}{G}\right) \\
&\colim_{\substack{(g,x) \to U^{\prime} \\ U^{\prime} \to G \times X\,\text{open}}}\quot{\theta^{A}_{X}}{G}\left(\quot{U^{\prime}}{G}\right) = \quot{\theta_{(g,x)}^{A}}{G},
\end{align*}
which proves the first claim and hence the proposition. 
\end{proof}

\bibliography{ConesBib}

\providecommand{\bysame}{\leavevmode\hbox to3em{\hrulefill}\thinspace}
\providecommand{\MR}{\relax\ifhmode\unskip\space\fi MR }
\providecommand{\MRhref}[2]{%
  \href{http://www.ams.org/mathscinet-getitem?mr=#1}{#2}
}
\providecommand{\href}[2]{#2}
\begin{thebibliography}{10}

\bibitem{PramodBook}
P.\@~N.\@ Achar, \emph{Perverse sheaves and applications to representation
  theory}, Mathematical Surveys and Monographs, vol. 258, American Mathematical
  Society, 2021.

\bibitem{ABV}
J.\@ Adams, D.\@ Barbasch, and D.\@~A.\@ Vogan, Jr.\@, \emph{The {L}anglands
  classification and irreducible characters for real reductive groups},
  Progress in Mathematics, vol. 104, Birkh\"{a}user Boston, Inc., Boston, MA,
  1992. \MR{1162533}

\bibitem{AguiMahajan.}
M.\@ Aguiar and S.\@ Mahajan, \emph{Monoidal functors, species and {H}opf
  algebras}, CRM Monograph Series, vol.~29, American Mathematical Society,
  Providence, RI, 2010, With forewords by Kenneth Brown and Stephen Chase and
  Andr\'{e} Joyal.

\bibitem{BBD}
A.\@ Belinson, J.\@ Bernstein, P.\@ Deligne, and O.\@ Gabber, \emph{Faisceaux
  pevers}, analyse et topologie sur les espace singuliers, Ast{\'e}risque, vol.
  100, Soc. Math. France, 1982, pp.~5 -- 171.

\bibitem{benzvi2015character}
D.\@ Ben-Zvi and D.\@ Nadler, \emph{The character theory of a complex group},
  Available at \url{https://arxiv.org/abs/0904.1247}.

\bibitem{BernLun}
J.\@ Bernstein and V.\@ Lunts, \emph{Equivariant sheaves and functors}, Lecture
  Notes in Mathematics, vol. 1578, Springer-Verlag, Berlin, 1994.

\bibitem{Bien}
F.\@~V.\@ Bien, \emph{Spherical {$\mathcal{D}$}-modules and representations of
  reductive {L}ie groups}, ProQuest LLC, Ann Arbor, MI, 1986, Thesis
  (Ph.D.)--Massachusetts Institute of Technology. \MR{2941104}

\bibitem{BorceuxCatAlg1}
F.\@ Borceux, \emph{Handbook of categorical algebra}, Encyclopedia of
  Mathematics and its Applications, vol.~1, Cambridge University Press, 1994.

\bibitem{brandenburg2021large}
M.\@ Brandenburg, \emph{Large limit sketches and topological space objects},
  2021, Available at \url{https://arxiv.org/abs/2106.11115}.

\bibitem{G2Cubics}
C.\@ Cunnigham, A.\@ Fiori, and Q.\@ Zhang, \emph{Arthur packets for {$G_2$}
  and perverse sheaves on cubics}, Adv. Math. \textbf{395} (2022).

\bibitem{CFMMX}
C.\@ Cunningham, A.\@ Fiori, A.\@ Moussaoui, J.\@ Mracek, and B.\@ Xu,
  \emph{Arthur packets for {$p$}-adic groups by way of microlocal vanishing
  cycles of perverse sheaves, with examples}, Memoirs of the AMS (2022).

\bibitem{Arun}
A.\@ Debray, \emph{The low-energy {TQFT} of the generalized double semion
  model}, Comm. Math. Phys. \textbf{375} (2020), no.~2, 1079--1115.
  \MR{4083885}

\bibitem{Delaney2019fusion}
C.\@ Delaney, \emph{Fusion rules for permutation extensions of modular tensor
  categories}, 2019, Available at \url{https://arxiv.org/abs/1604.06429}.

\bibitem{DeligneWeil1}
P.\@ Deligne, \emph{La conjecture de {W}eil. {I}}, Inst. Hautes \'{E}tudes Sci.
  Publ. Math. (1974), no.~43, 273--307. \MR{340258}

\bibitem{DeligneTanak}
P.~Deligne, \emph{Cat\'{e}gories {T}annakiennes}, The {G}rothendieck
  {F}estschrift, {V}ol. {II}, Progr. Math., vol.~87, Birkh\"{a}user Boston,
  Boston, MA, 1990, pp.~111--195. \MR{1106898}

\bibitem{SGA3}
M.\@ Demazure, A.\@ Grothendieck, and M.\@ Artin, \emph{Sch{\'e}mas en groupes
  {SGA} 3: Propri{\'e}t{\'e}s g{\'e}n{\'e}rales des sch{\'e}mas en groupes},
  Documents math{\'e}matiques, Soci{\'e}t{\'e} mathematique de France, 2011.

\bibitem{Chris2Noah}
C.~L. Douglas, C.~Schommer-Pries, and N.~Snyder, \emph{Dualizable tensor
  categories}, Mem. Amer. Math. Soc. \textbf{268} (2020), no.~1308, vii+88.
  \MR{4254952}

\bibitem{Torsten}
T.\@ Ekedahl, \emph{Is this a definition of equivariant derived category?},
  MathOverflow, URL:https://mathoverflow.net/q/27662 (version: 2010-06-10).

\bibitem{FreitagKiehl}
E.\@ Freitag and R.\@ Kiehl, \emph{\'{E}tale cohomology and the {W}eil
  conjecture}, Ergebnisse der Mathematik und ihrer Grenzgebiete (3) [Results in
  Mathematics and Related Areas (3)], vol.~13, Springer-Verlag, Berlin, 1988,
  Translated from the German by Betty S. Waterhouse and William C. Waterhouse,
  With an historical introduction by J. A. Dieudonn\'{e}.

\bibitem{GabrielZisman}
P.\@ Gabriel and M.\@ Zisman, \emph{Calculus of fractions and homotopy theory},
  Ergebnisse der Mathematik und ihrer Grenzgebiete, Band 35, Springer-Verlag
  New York, Inc., New York, 1967.

\bibitem{gallauer2021introduction}
M.\@ Gallauer, \emph{An introduction to six-functor formalisms}, 2021.

\bibitem{GepnerEtAl}
D.\@ Gepner, R.\@ Haugseng, and T.\@ Nikolaus, \emph{Lax colimits and free
  fibrations in {$\infty$}-categories}, Doc. Math. \textbf{22} (2017),
  1225--1266. \MR{3690268}

\bibitem{GiraudCN}
J.\@ Giraud, \emph{Cohomologie non ab\'{e}lienne}, Die Grundlehren der
  mathematischen Wissenschaften, Band 179, Springer-Verlag, Berlin-New York,
  1971. \MR{344253}

\bibitem{EGA1}
A.\@ Grothendieck, \emph{{\'E}l\'ements de g\'eom\'etrie alg\'ebrique : I. le
  langage des sch\'emas}, Publications Math\'ematiques de l'IH\'ES \textbf{4}
  (1960), 5--228 (fr).

\bibitem{SGA1}
\bysame, \emph{Rev\^etements \'etales et groupe fondamental ({SGA} 1)}, Lecture
  notes in mathematics, vol. 224, Springer-Verlag, 1971.

\bibitem{SGA7}
A.\@ Grothendieck, P.\@ Deligne, and N.\@ Katz, \emph{Groupes de monodromie en
  g{\'e}om{\'e}trie alg{\'e}brique: Sga 7}, Groupes de Monodromie en
  G{\'e}om{\'e}trie Alg{\'e}brique: SGA 7, no. v. 1, Springer-Verlag, 1972.

\bibitem{SGA5}
A.\@ Grothendieck and L.\@ Illusie, \emph{Cohomologie l-adique et fonctions l:
  Seminaire de geometrie algebrique du bois-marie 1965-66, sga 5}, Lecture
  Notes in Mathematics, Springer Berlin Heidelberg, 1977.

\bibitem{HiltonStammHA}
P.\@~J.\@ Hilton and U.\@ Stammbach, \emph{A course in homological algebra},
  second ed., Graduate Texts in Mathematics, vol.~4, Springer-Verlag, New York,
  1997. \MR{1438546}

\bibitem{TwoDimCat}
N.\@ Johnson and D.\@ Yau, \emph{2-dimensional categories}, Oxford University
  Press, Oxford, 2021. \MR{4261588}

\bibitem{TheoClaudia}
T.\@ Johnson-Freyd and C.\@ Scheimbauer, \emph{({O}p)lax natural
  transformations, twisted quantum field theories, and ``even higher'' {M}orita
  categories}, Adv. Math. \textbf{307} (2017), 147--223. \MR{3590516}

\bibitem{Keller}
B.\@ Keller, \emph{On differential graded categories}, International {C}ongress
  of {M}athematicians. {V}ol. {II}, Eur. Math. Soc., Z\"{u}rich, 2006,
  pp.~151--190.

\bibitem{KellyStreet}
G.\@~M.\@ Kelly and R.\@ Street, \emph{Review of the elements of
  {$2$}-categories}, Category {S}eminar ({P}roc. {S}em., {S}ydney, 1972/1973),
  1974, pp.~75--103. Lecture Notes in Math., Vol. 420.

\bibitem{RepGradedHecke}
C.\@ Kriloff and A.\@ Ram, \emph{Representations of graded {H}ecke algebras},
  Represent. Theory \textbf{6} (2002), 31--69.

\bibitem{Landesman}
A.\@ Landesman, \emph{The smooth base change theorem}, Available at
  \url{https://virtualmath1.stanford.edu/~conrad/Weil2seminar/Notes/L7-8.pdf}.
  Accessed 2023 September 24.

\bibitem{Leinster}
T.~Leinster, \emph{Higher operads, higher categories}, London Mathematical
  Society Lecture Note Series, vol. 298, Cambridge University Press, Cambridge,
  2004.

\bibitem{LiuAGAC}
Q.\@ Liu, \emph{Algebraic geometry and arithmetic curves}, Oxford Graduate
  Texts in Mathematics, vol.~6, Oxford University Press, Oxford, 2002,
  Translated from the French by R.\@ Ern\'{e}, Oxford Science Publications.

\bibitem{LurieCobord}
J.\@ Lurie, \emph{On the classification of topological field theories}, Current
  developments in mathematics, 2008, Int. Press, Somerville, MA, 2009,
  pp.~129--280. \MR{2555928}

\bibitem{LusztigCuspidal2}
G.\@ Lusztig, \emph{Cuspidal local systems and graded {H}ecke algebras. {II}},
  Representations of groups ({B}anff, {AB}, 1994), CMS Conf. Proc., vol.~16,
  Amer. Math. Soc., Providence, RI, 1995, With errata for Part I [Inst. Hautes
  \'{E}tudes Sci. Publ. Math. No. 67 (1988), 145--202; MR0972345 (90e:22029)],
  pp.~217--275.

\bibitem{MacLaneCWM}
S.\@ Mac~Lane, \emph{Categories for the working mathematician}, second ed.,
  Graduate Texts in Mathematics, vol.~5, Springer-Verlag, New York, 1998.

\bibitem{MacLaneMoerdijk}
S.\@ Mac~Lane and I.\@ Moerdijk, \emph{Sheaves in geometry and logic},
  Universitext, Springer-Verlag, New York, 1994, A first introduction to topos
  theory, Corrected reprint of the 1992 edition.

\bibitem{BenThesis}
B.\@ MacAdam, \emph{The functorial semantics of {L}ie theory}, 2022, PhD
  Thesis. Availabe at
  \url{https://prism.ucalgary.ca/items/05830be5-fbdf-4be2-a9a5-cc1e73762d7f/full}.

\bibitem{mann2022padic}
L.\@ Mann, \emph{A $p$-adic 6-functor formalism in rigid-analytic geometry},
  2022, Available at \url{https://arxiv.org/abs/2206.02022}.

\bibitem{MayTrace}
J.\@~P.\@ May, \emph{The additivity of traces in triangulated categories}, Adv.
  Math. \textbf{163} (2001), no.~1, 34--73.

\bibitem{milneiAG}
J.\@~S.\@ Milne, \emph{Algebraic groups (v2.00)}, 2015, Available at
  www.jmilne.org/math/, p.~528.

\bibitem{MirkovicVilonen}
I.\@ Mirkovi\'{c} and K.\@ Vilonen, \emph{Characteristic varieties of character
  sheaves}, Invent. Math. \textbf{93} (1988), no.~2, 405--418.

\bibitem{MirkVilDuality}
\bysame, \emph{Geometric {L}anglands duality and representations of algebraic
  groups over commutative rings}, Ann. of Math. (2) \textbf{166} (2007), no.~1,
  95--143.

\bibitem{GIT}
D.\@ Mumford and J.\@ Fogarty, \emph{Geometric invariant theory}, second ed.,
  Ergebnisse der Mathematik und ihrer Grenzgebiete [Results in Mathematics and
  Related Areas], vol.~34, Springer-Verlag, Berlin, 1982.

\bibitem{NeemanTriangulated}
A.\@ Neeman, \emph{Triangulated categories.}, Ann.Math.Studies 148, Princeton
  University Press, 2001.

\bibitem{Orlov_2018}
D.\@~O.\@ Orlov, \emph{Derived noncommutative schemes, geometric realizations,
  and finite dimensional algebras}, Russian Mathematical Surveys \textbf{73}
  (2018), no.~5, 865–918.

\bibitem{DoretteMe}
D.\@ Pronk and G.\@ Vooys, \emph{Equivariant tangent categories}, 2023,
  Available at \url{https://arxiv.org/abs/2308.11753}.

\bibitem{Raynaud}
M.~Raynaud, \emph{Faisceaux amples sur les schemas en groupes et les espaces
  homogenes}, Lecture Notes in Mathematics, Springer Berlin Heidelberg, 2006.

\bibitem{RiehlVertityElementss}
E.\@ Riehl and D.\@ Verity, \emph{Elements of {$\infty$}-category theory},
  2022.

\bibitem{FreeAdj}
S.\@ Schanuel and R.\@ Street, \emph{The free adjunction}, Cahiers Topologie
  G\'{e}om. Diff\'{e}rentielle Cat\'{e}g. \textbf{27} (1986), no.~1, 81--83.

\bibitem{ScholzeSixFunctor}
P.\@ Scholze, \emph{Six functor formalisms}, 2023, Available at
  \url{https://people.mpim-bonn.mpg.de/scholze/SixFunctors.pdf}.

\bibitem{SerreAlgebraicFibreSpaces}
J.-P.\@ Serre, \emph{Espaces fibr\'{e}s alg\'{e}briques (d'apr\`es {A}ndr\'{e}
  {W}eil)}, S\'{e}minaire {B}ourbaki, {V}ol. 2, Soc. Math. France, Paris, 1995,
  pp.~Exp. No. 82, 305--311.

\bibitem{stacks-project}
The {Stacks Project Authors}, \emph{\textit{Stacks Project}},
  \url{https://stacks.math.columbia.edu}, 2020.

\bibitem{GregStevensonSplitMonic}
G.\@ Stevenson, \emph{How do {I} know the derived category is {NOT} abelian?},
  MathOverflow, Available at \url{https://mathoverflow.net/q/15662}.

\bibitem{toën2014derived}
B.\@ To{\"e}n, \emph{Derived algebraic geometry}, 2014, Available at
  \url{https://arxiv.org/abs/1401.1044}.

\bibitem{ToenDAGDQ}
\bysame, \emph{Derived algebraic geometry and deformation quantization},
  Proceedings of the {I}nternational {C}ongress of {M}athematicians---{S}eoul
  2014. {V}ol. {II}, Kyung Moon Sa, Seoul, 2014, pp.~769--792.

\bibitem{Vakil}
R.\@ Vakil, \emph{The rising sea foundations of algebraic geometry}, 2022,
  Online text available at
  \url{https://math.stanford.edu/~vakil/216blog/FOAGaug2922public.pdf}; note
  the 2022 August 29 version was used. Accessed 2024 January 17.

\bibitem{VerdierDerivedCat}
J.-L.\@ Verdier, \emph{Cat\'{e}gories d\'{e}riv\'{e}es: quelques r\'{e}sultats
  (\'{e}tat 0)}, Cohomologie \'{e}tale, Lecture Notes in Math., vol. 569,
  Springer, Berlin, 1977, pp.~262--311.

\bibitem{Vistoli}
A.\@ Vistoli, \emph{Notes on {G}rothendieck topologies, fibered categories and
  descent theory}, \url{arXiv:math/0412512v4}, 5 2007.

\bibitem{MyThesis}
G.\@ Vooys, \emph{Equivariant functors and sheaves}, 2021, PhD.\@ thesis.
  Available at \url{https://arxiv.org/abs/2110.01130}.

\bibitem{Geordie}
G.\@ Williamson, \emph{Bernstein and {L}unts' fundamental example}, Essay
  available at
  \url{http://people.mpim-bonn.mpg.de/geordie/FundamentalExample.pdf}.

\bibitem{Yetter}
D.\@~N.\@ Yetter, \emph{Functorial knot theory}, Series on Knots and
  Everything, vol.~26, World Scientific Publishing Co., Inc., River Edge, NJ,
  2001, Categories of tangles, coherence, categorical deformations, and
  topological invariants.

\end{thebibliography}
\bibliographystyle{amsplain}

\printindex[terminology]
\printindex[notation]

\end{document}